\documentclass[graybox,envcountchap,12pt]{svmono}
\usepackage{natbib}
\usepackage{bibentry}
\nobibliography*
\usepackage{makecell} 
\usepackage{mathptmx}
\usepackage{helvet}
\usepackage{courier}
\usepackage{amsmath}
\usepackage{amssymb}
\usepackage{mathrsfs}
\usepackage{wrapfig}
\usepackage{pdfpages}
\usepackage{geometry}
\usepackage{marginnote}
\usepackage{siunitx}
\usepackage[nobiblatex]{xurl}
\usepackage{xcolor}
\usepackage{color}
\usepackage{transparent}
\usepackage{ulem}
\usepackage{imakeidx}
\newcommand{\indexme}[1]{\index{#1}}
\usepackage{mathtools}
\usepackage{fancyvrb}
\usepackage{caption}
\usepackage{enumitem} 
\newcounter{saveenumi}
\setenumerate{itemsep=4pt,topsep=4pt}
\newcounter{saveenum}
\newcommand{\saveenumerate}{%
\setcounter{saveenum}{\value{enumi}}}
\newcommand{\restoreenumerate}{%
  \setcounter{enumi}{\value{saveenum}}}
\usepackage{quotchap}
\usepackage{epigraph}
\usepackage[bottom]{footmisc} 
\usepackage{graphicx,makeidx,multicol}
\usepackage{tikz}
\newcommand*\numcircledtikz[1]{\tikz[baseline=(char.base)]{
            \node[shape=circle,draw,inner sep=1.2pt] (char) {\upshape{#1}};}} %
\makeindex             

\usepackage{amsthm}
\newcommand{\B}[1]{\mathbf{#1}}
\renewcommand{\qedsymbol}{\ensuremath{\hfill\blacksquare}}
\newcommand{\relu}{\mathrm{ReLU}}

\newcommand{\reluf}[1]{\mathrm{ReLU}\left(#1\right)} 
\DeclareMathAlphabet{\mathbfsf}{\encodingdefault}{\sfdefault}{bx}{n}
\newcommand{\twoform}[1]{\mathbfsf{#1}}
\definecolor{myblue1}{rgb}{0.4, 0.7, 1}
\definecolor{myblue2}{rgb}{0, 0.4, 0.8}
\definecolor{myred1}{rgb}{1, 0.4, 0.4}
\newcommand{\wubs}[2]{w_{\textcolor{myblue1}{#1}}^{\textcolor{myred1}{#2}}}
\newcommand{\bubs}[2]{b_{\textcolor{myblue1}{#1}}^{\textcolor{myred1}{#2}}}
\newsavebox\MBox
\newcommand\Cline[2][black]{{\sbox\MBox{$#2$}%
  \rlap{\usebox\MBox}\color{#1}\rule[-1.5\dp\MBox]{\wd\MBox}{0.5pt}}}
\newtheorem{axiom}{Axiom}
\renewcommand\qedsymbol{$\blacksquare$}

\usepackage{titlesec}
\makeatletter
\@addtoreset{section}{part}
\makeatother
\titleformat{\part}[display]
{\normalfont\LARGE\bfseries\centering}{}{0pt}{}
\newcounter{magicrownumbers}
\newcommand\rownumber{\stepcounter{magicrownumbers}\arabic{magicrownumbers}.}
\newcolumntype{L}{>{\centering\arraybackslash}m{3cm}}
\usepackage{float} 
\makeatletter
\newcommand\newtag[2]{\def\@currentlabel{#1}\label{#2}}
\makeatother

\usepackage{verbatim}
\begin{document}
\normalem

\author{Brian D. Wood}
\title{Introduction to Engineering Mathematics and Analysis}
\subtitle{Modeling Physical Systems Using the Language of Mathematics}
\maketitle

\clearpage

\vspace*{\fill}
\noindent{ISBN-13 978-1-955101-33-2}\\

\noindent\includegraphics[height=1.5cm]{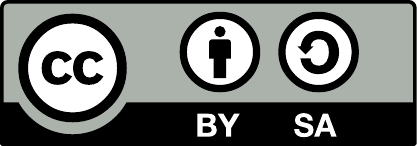}\\

\noindent This work is licensed under a Creative Commons Attribution 4.0 International License.

\noindent\url{https://creativecommons.org/licenses/by-sa/4.0/}\\

\noindent{Text Revision 1.3} \vspace{2mm}\\

\noindent{How to cite:}\vspace{2mm}\\
Wood, B. D. (2023). Introduction to Engineering Mathematics and Analysis: Modeling Physical Systems Using the Language of Mathematics. First Edition, Rev. 1.3. Oregon State University, Corvallis OR. DOI: \url{https://doi.org/10.5399/osu/1152}.
\clearpage

\frontmatter

%
%

\preface

\def\CHAP{chapter000_Introduction}
\theoremstyle{definition}


\addcontentsline{toc}{chapter}{Preface}  
\begin{minipage}[t]{.48\textwidth}
\begin{quote}
It has become almost a cliche to remark that nobody boasts of ignorance of literature, but it is socially acceptable to boast ignorance of science and proudly claim incompetence in mathematics...\\

--Richard Dawkins
\end{quote}
\end{minipage}

\subsection*{Why study mathematics?}

Scientists and engineers are, by definition, problem solvers.  But, it is curious to reflect on how they become problem solvers.  While it is possible that early in one's academic career, a course whose title included the words ``problem solving" might have been taken, but this is not likely to have been where students generally learn how to actually solve problems.  The topic of how engineers and scientists learn what they ultimately know is an area of research in education.  While that question is much to broad to address here, some comments about the use of mathematics in the education process is, however, within scope. 

\noindent The refrain ``But, when am I ever going to \emph{use} this?" is a familiar one to anyone who has ever studied (or taught) mathematics.  It underscores a particular problem with the way in which mathematics is presented.  In studying, for example, poetry, it would be unusual to hear those studying the subject to exclaim ``But when am I ever going to \emph{use} this poem?".  The problem begins to be clear when examining thinking about why mathematics is viewed as different than, say, poetry.  There is often an \emph{expectation} by learners that mathematics should be useful in a very applied sense.  This expectation is understandable because most of our early mathematical training (e.g., learning to add, subtract, and multiply; solving basic algebraic problems) has very obvious utility.  It is usually at the point where higher levels of mathematical thinking are introduced (e.g., formal linear algebra, calculus, differential equations) that the question of usefulness arises.  

Most of us study mathematics because it is \emph{useful}.  Few would argue that knowing how to add, subtract, multiply, and divide are skills that are not useful. Similarly, few would argue, for example, that spreadsheets are not useful; and spreadsheets are primarily a tool based in mathematical operations.  Most engineers and scientists have had the opportunity and need to apply the various rules of algebra to solve \emph{actual real-world problems}.  So, there are  subject areas for which the tools of mathematics are generally agreed to be useful.

One of the problems that occurs in the study of \emph{higher mathematics} is that it is not made clear \emph{in what sense} it is useful.  In part, this is because the utility of higher mathematics is not as directly obvious as it is for more basic mathematical reasoning.  It is worth making this as concrete as possible, however, to help better understand why the study of higher mathematics can be worth the effort that one has to put in. 

Learning mathematics is, in part, about learning algorithmic thinking.  Beyond this, however, the study of mathematics helps us learn and practice with algorithmic thinking.  Even if you never apply any of the mathematical tools learned, the \emph{process of algorithmic thinking} is translatable to problem solving in a very general context.  The study of mathematics, then, can help you become a better problem solver, even if mathematics is not used in solving the problem!  This presents a good rational for studying higher mathematics generally.  Not only does it teach one how to actually be conversant with new mathematical constructs (which, can arise in applications, depending on one's career path), but it also helps inculcate new skills for problem solving that are useful well beyond the intrinsic application of mathematical methods.  

Therefore, even if you never compute another derivative, or solve another partial differential equation after taking taking a course that covers these topics, the very act of learning the material will pay benefits.  The study of mathematics at all level helps establish and reinforce your ability to think about problems, and to generate algorithmic methods for solving the problems you encounter.  While this may or may not involve the formal mathematics that you have learned, the process of thinking about problems will be indelibly colored (to the positive!) by your experience in learning mathematics.   

\tableofcontents

\mainmatter
%
\abstract*{This is the abstract for chapter 000}

\begin{savequote}[0.55\linewidth]
\vspace{-10mm}
``In this sense, therefore, mathematics would appear to
be both more and less than a language for while being limited in its
linguistic capabilities it also seems to involve a form of thinking that has something in common with art and music."
\qauthor{{A. Ford and F.D. Peat}\newline{~from ``The Role of Language in Science"}}

``All models are wrong, but some are useful."

\qauthor{George Box, statistician.}
``The best material model of a cat is another, or preferably the same, cat.” 
\qauthor{{\quad\quad\quad Arturo Rosenblueth and} \newline{Norbert Wiener,}\newline{from ``The Role of Models in Science"}}\indexme{cat!Rosenblueth and Wiener}
\end{savequote}

\def\CHAP{chapter_00_models}
\theoremstyle{definition}


\chapter{Languages, Mathematics, and Models}

Before jumping into the topic of the various mathematical methods that are used for modeling in science and engineering, it is worth spending a little time discussing some of the more philosophical (or ``big picture") notions of what comprises mathematical modeling.  There is much to cover under the topic of \emph{mathematical modeling}, and some of these are best experienced through examples and applications.  However, there are also some good high-level questions to ponder regarding modeling generally.  Some of these questions are posed in the material that follows, and the discussion gravitates more towards the philosophy of science than is typical for the remainder of the text.  Despite the fact that the material is partly philosophical, that does not imply in any way that it is somehow unimportant in an applied sense.  Understanding what is implicitly embedded in the process of mathematical modeling is actually a very useful thing.  For example, explicitly noting (frequently unrecognized) assumptions is a key part of being a good problem solver!  In problem solving, to the extent possible, it is always a good idea to explain one's assumptions; not, perhaps, to the level of philosophical detail that is examined in the material following, but to the extent that it will help others (or yourself) understand how the problem was solved.

Because language (and in particular the language that is \emph{mathematics}) is part of \emph{modeling} it is worthwhile spending just a little effort attempting to understand exactly what languages are.  The next section begins to address this question.  While the question broadly enters some deep philosophical terrain, the summary discussion following helps at least expose the questions without becoming overly mired in philosophical structures.  This chapter as a whole is a collection of concepts regarding modeling and, in particular, modeling with mathematics.  The information that follows is material that helps better understand the process of modeling, and how mathematics helps the modeler toward that goal.

\section{What is a Language?}\indexme{language}

The purpose of discussing language in this introductory chapter is to explore the link between languages, mathematics, and models.  We cannot cover the topic of \emph{languages} in any depth, but even a cursory overview is useful for the purposes of this introduction.  

If you have ever studied languages, you may be familiar with the concepts of \emph{symbols}, \emph{syntax}, and \emph{semantics}.  Symbols are the way that we \emph{express} language (written or verbal).  Syntax explains how a set of words (symbols) are used together to \emph{according to rules} to form any of the following: (i) a statement/assertion, (ii) a question, (iii) a command, or (iv) an exclamation (these are the four options in most languages).  Semantics, on the other hand, seeks to assign \emph{meaning to syntax}.  

Languages can be divided into \emph{natural} \indexme{language!natural} languages and \emph{formal} languages.  A natural language is one that that has been developed instinctively and heuristically by humans to communicate with one another.  Such languages can be spoken, written, or both.  Natural languages tend to be highly \emph{expressive}; in other words, many kinds of concepts ranging from the subjective (``do you like this painting?") to the objective ("the toaster is broken").  One of the difficulties with natural languages is that it is difficult to codify the rules of the language.  So, while they are very expressive, it is difficult to know, explicitly, all of the rules of the language.  Most people who attempt to become fluent in a foreign language eventually encounter this problem; while it is not difficult to become proficient enough to communicate well, it is exceptionally difficult to inculcate the subtleties of every day language (e.g., the use of slang, idioms, and inflection are challenging) so that one is truly fluent.

A \emph{formal language} \indexme{language!formal} is one that is developed specifically to have a set of well-defined rules prescribing it.  These tend to be primarily written languages (for hopefully obvious reasons).  In these instances, one exchanges broad expressiveness with complex rules for restricted expressiveness, but with exceptionally well-defined rules.  The most familiar example of a formal language is probably a computer language such as Python or C++.  However, formal languages were developed well before computers existed.  As an example, the formal language known as first-order (predicate) logic was discovered in the late 1800s by a mathematician named C.S. Pierce, and brought (roughly) its present form by the famous mathematician David Hilbert in about 1915 \citep{ewald_first_order}.  First-order logic is widely considered to be adequate to allow the axiomatization of all ordinary mathematics.  In fact, this suggests that because mathematics arises from the language of first-order logic, mathematics itself has all of the necessary properties of a formal language.  The relationship between language and mathematics is a subject of study in its own right, and it is a fascinating topic because it co-mingles such disparate disciplines.  Interested students can find out more in the texts by \citep[][Chp.~13]{rosen2019discrete} listed in the bibliography.

While we will discuss modeling in the material that follows, it is worth pointing out that one of the most famous philosophers of the 20th century, Ludwig Wittgenstein, suggested that language itself reflected reality.  According to some scholars, this  implies that language is itself a model.  For example, \citet{mendie2019language} state of Wittgenstein's philosophy that  ``... a proposition (language), is a picture of reality, and a proposition (language), is also a model of reality as we imagine it."  We will not delve much more into the metaphysical aspects of language as a model that is used to create other models.  It is a worthwhile concept to ponder, however, when we discuss the concept of \emph{models} more generally below.

\section{What is Mathematics?}\indexme{mathematics}

It may seem odd to start out a textbook on mathematics and modeling with a (very rough) definition of what \emph{language} is.  But, because mathematics is itself a language, it helps us establish how to approach mathematics and mathematical thinking.  Like many broad concepts, actually \emph{defining} what mathematics can be challenging.  It is somewhat exceptional that something that we are so familiar with, and use so frequently is, at the same time, elusive to define.  None of us are alone in experiencing this difficulty, however.  Mathematicians, philosophers, logicians, linguists, and scientists of various disciplines have been attempting to define \emph{what} mathematics actually is for a very long time (for exammple, see the book \emph{What is Mathematics?} by the famous mathematician R. Courant, (Courant and Robbins, \citeyear{courant1996mathematics}) for one expert's opinion). 

From our perspective, the question of ``what is mathematics?" will be approached pragmatically.  We care about the question primarily in the hopes that understanding something about this question will also give us additional insight into the use and limits of mathematics.  It can also help better understand how to engage with mathematics (and mathematical thinking) to become more proficient problem solvers.

In the material above, we introduced the idea of mathematics as a \emph{language}.  Without getting to hung up on the intricate details (e.g., mathematics as a mode of thinking versus the particular way that it is expressed, see for example \citet{ford1988role} for more on this discussion), it is reasonably well accepted from a philosophical standpoint that mathematics is a formal language \citep[][\S I.2]{gowers2008princeton}.\indexme{mathematics!as a language}

While the status of mathematics as a language might sound like a primarily academic issue, it actually is one that has plenty of real-world relevance to learners of mathematics.  Generally, learning a new language is viewed as being a \emph{significant undertaking}, often requiring years of practice to reach any level of competence.  And, although mathematics is a formal language (so its rules are well-defined), becoming proficient in a strict rule-based language is also a challenge.  So, it is perhaps not that surprising that many people find mathematics a difficult topic-- the study of mathematics is, in a very real sense, the study of a foreign language.  It happens to be a very formal language with limited (and also exacting) structure, but it is a language nonetheless.  As a language, one can study it from an \emph{academic perspective} extensively, but still not be particularly good at using the language.  The reality is, like many languages, mathematics requires \emph{practice} to master.  

To many students, recognizing mathematics as a language can be somewhat reassuring.  First, it means that one can dismiss the idea that mathematics should somehow be easy, even if one has been studying it for many years.  Becoming proficient in mathematics is much like learning \emph{new ways of using and understanding a language}, even one that you already have some competence with.  For example, someone who is fluent in modern written and spoken English might still struggle when first encountering the language of Shakespeare.  Or one might be challenged by having to learn language that has specific disciplinary meaning, such as the jargon used in philosophy, history, or psychology.  

The importance here is that recognizing mathematics as a language can (and should) change the way that one learns about new mathematical ideas.  As a language, it takes practice (e.g., explicitly solving problems) and study to understand it.  However, it is not generally true that there are people who ``get" math and people who ``don't get" math, any more than the the same categories would not be made for, as an example, the language of Spanish or French.  Thinking about mathematics as a language allows one the flexibility to retrain their way of thinking about mathematics, and also to realize that practice (with attendant mistakes!) is an essential component of learning the language.

When one accepts that learning mathematics is much like learning a foreign language, it can put a fresh perspective on a topic that many otherwise approach unenthusiastically.  Sometimes the process of learning can involved \emph{un}learning patterns of thinking that were not productive.  This is often the case in learning mathematics, and hopefully the perspective of ``mathematics as a language" can be a organizational tool to help engage with the topic in new ways.

\section{What are Models?}\label{what?}\indexme{model}\indexme{modeling}

Defining what comprises a \emph{model} is a daunting task.  If one looks to the philosophy of science, there is no agreement whatsoever on the topic.  Fortunately, we have a more practical than philosophical need for defining models.  Thus, for the purposes of this text, we can define a model as follows.

\begin{definition}[model]\indexme{model!definition}
A model is a conceptual or physical \emph{abstraction} used to predict and explain a generally more complex feature (a thing), phenomenon (some natural entity that exists, such as gravity),  or process (such as, say, mass transport by diffusion) \emph{supposed} to occur in the physical world.  Usually, one thinks of a model as mapping (or explaining) how a set of inputs (independent variables) relates to a set of outputs (dependent variables), but this should be interpreted very liberally.  Importantly, a model can describe processes that may be only hypothetical.
\end{definition}
The focus on ``feature or process in the physical world" is purposeful.  Certainly one can imagine models for \emph{metaphysical} phenomena.    Here, we use the word \emph{metaphysical} in its most basic meaning, that is, outside of a possible description of the laws of physics as we understand them currently. As an example, one might develop a model for explaining the existence of ghosts; however, such a model would almost certainly not be within the bounds of physics as we currently understand it.  This still leaves many grey areas depending on context: for example, would a model of social-cultural phenomena, which is within the purview of the sciences,  be a valid one for us to consider?  Here, we can say that as long as there were a model that did not violate the laws of physics to explain the phenomenon of interest, then the model would be a valid one for us to discuss.  This will make more sense, perhaps, after the discussion of \emph{empirical models} provided below.

To be more compact, in the future we can refer to the \emph{feature}, \emph{phenomenon}, or \emph{process} as a \emph{physical system}\indexme{model!system}.  The definition of a system can be given as follows.

\begin{definition}[system]
A \emph{physical system}, or simply \emph{system},\indexme{system} is a group of interacting or interrelated elements (physical features, phenomena, or processes) that act according to a set of rules (which are may or may not be fully known) to form a unified whole.  A system may be described as being \emph{discrete} or \emph{continuous}. For discrete systems, each element is distinct and the total number of elements is, in principle, a unique integer.  For continuous systems, each element is defined (non-uniquely) as being part of the whole, but the division is only conceptual; the system cannot be though of as a unique sum of individual parts. \indexme{modeling!system}\indexme{system}
\end{definition}

\subsection{What are the Purposes of Models?}

The creation of models comes from a particular need or purpose.  Most readers of this text will already have had experience with models, and could probably arrive at some very good reasons that they are useful.   The following are reasons that model building is enacted are as follows, but the list is not necessarily exhaustive.

\begin{enumerate}
\item \textbf{For problem summary and simplification}.  One very common reason that models are generated are to \emph{simplify} an otherwise overwhelmingly complex system.  The case of an ideal gas (which will also be used as an example of complexity below) provides a great example.  At near standard temperature and pressure, a momentum balance on all molecules in a volume of gas can provide a very accurate value for the pressure.  This, of course, is quite a complicated computation because one mole of gas has $6.02 \times 10^{23}$ molecules.  Alternatively, the ideal gas law, $P=1/V (nRT)$ will also express the pressure provided one knows $V$, $n$, $R$, and $T$.  Computation of the pressure via the ideal gas law is clearly more efficient than the momentum balance computation.

\item \textbf{For understanding a system}.  Sometimes models are constructed to help better understand a system.  Again, we can use the ideal gas as an example.  While the ideal gas law was originally developed empirically (we will discuss \emph{empirical models} below), as science progressed there was a desire to understand more about how gases behaved.  Statistical mechanics is a branch of science that computes the statistical behavior of large numbers of bodies that obey the laws of mechanics (Newton's laws or the laws of quantum mechanics).  It turns out that one can show that applying Newton's laws to a large number of molecules in a fixed volume and known temperature, then the ideal gas law can be \emph{derived} as a result.  This is a case of proving an empirical \emph{macroscopic} law by computing averages over a well-defined \emph{microscopic} model.  Here the words ``macroscopic" and ``microscopic" are used only in a relative sense to establish the difference in length scales investigated in the two models.  The result is that now the ideal gas law can be shown to be consistent with both the original experiments that lead to the empirical law, and with a \emph{conceptual} model based on Newton's laws.  The fact that the ideal gas law can be developed by these to independent approaches simultaneously increases confidence in it, and also provides additional explanatory ability for our models.  

\item \textbf{For allowing prediction of system behavior}.  One of the most common use of models is to predict system behavior.  While there are many examples, one might consider something as simple as determining the forces on a truss that will be used as a bridge for pedestrian traffic.  While actually \emph{engineering} a bridge is a multistep process, it would at least start with applying Newton's laws to the proposed bridge structure to determine the distribution of forces in the members of the truss under typical loads.   This is something that is often done in undergraduate physics, or a course in statics.  With information about the distribution of forces, one could then begin to determine what materials (type of material, material shapes and sizes) would be required to function properly.  The advantage of this kind of approach is that predictions can be made for various kinds of options for the truss.  Different structural designs, different materials, and different loading conditions can all be done as various \emph{what if?} scenarios.  This is useful and efficient.  The prediction of system behavior from a model means that one does not need to physically build many different trusses and then test each of them.  The advantages of this kind of modeling are hopefully obvious.  
\end{enumerate}

\subsection{What Kinds of Models Can be Constructed?}

Restricting ourselves now to models that have some relevance to a physical system, we can identify several kinds of useful models and put them in categories.  Almost every attempt to categorize broad concepts leads to some lack of distinction for particular cases, but the general organization of ideas is still a useful one.  Hence, we categorize models as follows.

\begin{enumerate}[wide]
\item \textbf{Physical models}.  The quote at the start of this chapter ``The best model of a cat is another cat, preferably the same cat" is both meant to be somewhat humorous, but also to relate an essential feature of \emph{physical} models.  Physical models are, as their name suggests, models that are created out of matter (and, in some senses, energy) for representing a feature, phenomenon, or process.  The physical model can be an analogue; that is to say, it can be a phenomenon that simulates some physically completely different phenomenon.  For example, steady groundwater flow and electrical current can be described by the same differential equation.  Thus, in the past, researchers have used electrical analogues to model groundwater flow.  Physical models can also be based on the concept of \emph{similitude}.\indexme{model!similitude}  Here, the idea is to make a physical model that represents the actual physical system, but at a different scale.  The similarity between the physical model and reality is enforced by assuming that one or more dimensionless numbers are identical between the two.  The idea of similitude seems to have been developed by the fluid mechanist, Osborne Reynolds in his study of turbulence in channels \citep{reynolds1883xxix}.  At an extreme, a physical model can be the system of study itself.  At Oregon State University we have the H.G. Andrews Experimental Forest, which contains a number of experimental watersheds (a watershed is a geographically-defined region where rainfall and snowfall are channeled to an outflow).  In this case, ``the best model of a watershed is another watershed". 

One key aspect of physical models is that the model \emph{input} is the physical model itself, and a specific set of initial, boundary, and parametric conditions.  The \emph{output} of a physical model is almost always a sequence of \emph{measurements} in time, in space, or in both. 
In principle, a physical model is transparent model in regard to its assumptions, etc., and is sometimes called a \emph{white box} model (this is in contrast to a \emph{black box} model, described below).  

\item \textbf{Conceptual models}.  Conceptual models are models that are represented primarily by the use of some \emph{language} as defined above.  The kinds of conceptual models that we will be most interested in are ones that can be defined using the language of mathematics.  However, these are not the only kinds of conceptual models that are useful!   Explanations of systems using regular English (or any other natural language), using flow charts, or even drawings might all be classified as conceptual models.  The primary distinguishing feature of a conceptual model is that it is (i) \emph{not} a direct physical analogue to the system being modeled, and (ii) it is expressed using some kind of symbolic language (and here, we will extend the use of the words ``symbolic language" to include graphical representations such as flow charts or descriptive illustrations).  One advantage of conceptual models that are expressed using mathematics (or formal logic) as the language is that the communication of the model and the model inputs, structured, and outputs (or results) are about as clear and incontrovertible as is possible for a conceptual model.  Some conceptual models deal explicitly with uncertain or fuzzy data;  such models are sometimes referred to as \emph{grey box models}.

\item \textbf{Empirical models}. Empirical models are models that are based \emph{primarily} upon relating \emph{model inputs} to \emph{model outputs}, without necessarily understanding the intervening mechanisms that transform inputs to outputs.  Hence, such models are sometimes called \emph{black box} models (a term of uncertain origin, but one made popular in the context of modeling by  R. Ashby and N. Wiener around 1960).  Empirical models are probably the oldest kind of modeling that has been done in science, and has been successful in many different areas of the natural sciences.  There might be a tendency to view empirical models as being somehow \emph{less than} their non-empirical counterparts.  This is not generally true!  As one example, one might consider the original ``laws of friction" which, in its simplest form, stated that the friction experienced between two surfaces was linearly proportional to the normal force between them.  This, of course, has been a useful model that is still widely applied today.  It is also a model that has been based largely on empiricism.  While it is true that we now know that the observed friction arises from microscopic variations in the surfaces (and, depending on the scale of investigation, on other microscopic forces).  Yet, rarely do we attempt to explicitly model such phenomena.  Rather, the empirical evidence that has been built up over time provides a strong argument that the empirical model is both valid and accurate. 

Empiricism in science is even one branch of the philosophy of science.  Famously, Ernst Mach (of Mach number fame) was a staunch empiricist, and felt that science should be based as much as possible on what was strictly observable, and that interpretations should not be made via unobserved quantities.  While \emph{empirical models} still have a role in the natural sciences, rigidly subscribing to empiricist philosophies has fallen out of favor; in part, this has been because of some failures of the approach.  It was well known, for example, that Mach opposed Ludwig Boltzmann and others who proposed an atomic theory of physics.  His objections at the time were that atoms were not directly observable, and thus positing the existence of them was non-scientific.  This, of course, seems somewhat of a backward stance now that we routinely and directly measure atoms (e.g., in atomic force microscopy).  Thus, there must be some balance between purely conceptual modeling of the universe around us and purely empirical modeling.  A purely conceptual model that is incapable of \emph{ever} being measured (e.g., string theory) is somewhat vacuous.  A purely empirical model that does not dare to explain the phenomena that is responsible for the observations fails to move science forward at the pace it otherwise might.  

\item \textbf{Digital models}.  Digital models have been included last because of their unique status as models.  \emph{Digital models} are, for our purposes, models that rely on a digital computer to process input data and generate output data.  Such models, curiously, are in some ways a mixture of the three models specified above.  First of all, computers are, themselves, physical apparatus; much like analogue electrical models for groundwater flow, digital computers rely on physical hardware to construct a model.  Second, digital models are often \emph{approximations} of conceptual models.  For example, when programs uses a root-finding method on a computer, the underlying theory is entirely conceptual; the actual execution, however, is limited by the finite arithmetic of a computer, and is thus an approximation.  Finally, note that computers are complex machines that involve processors, memory, storage, etc.  Few users of computers actually understand, in detail, how a computer actually works.  Even someone who does understand how a computer works in detail usually cannot observe all processes that occur for a computational algorithm to transform input data into output data.  In this sense, we must view digital models as being, in some sense, empirical, if for no other reason than it functionally a ``black box".  Regardless of these limitations, digital models have quite literally revolutionized modeling in the natural sciences.   Most recently, the advent of \emph{machine learning} returns to sciences empirical roots in some ways by allowing users to make empirical sense of large sets of data by allowing a computational algorithm to fit or categorize data.  While there are many efforts to make such models understandable and interpretable by humans, most of them must be viewed currently as largely empirical.  This hardly means that they are not useful; rather, it means only that we have explanation that sometimes lacks deeper understanding.  

\end{enumerate}

\section{The Modeling Process}\indexme{model}\indexme{modeling}

One of the beneficial uses of mathematics for scientists and engineers is that it is a natural language for problem solving; although it is hardly true that all problem solving requires mathematics.   The process of problem solving involves the use of a language -- which may or may not be formal language for expressing mathematics or logic-- so that the problem solving process can be represented.  Assuming that a problem has been identified, solving the problem using a language (including the language of mathematics involves) at least the following steps.

\begin{enumerate}[wide]

    \item {\bf Abstraction.} \indexme{modeling!abstraction}  Once identified, the problem needs to be described in some terms that are \emph{simpler} than the actual system.  If the solution is to be a mathematical one, then the simplified system needs to be ``translated" into the appropriate mathematical expression.  Abstraction is often done in a universal sense (e.g., using \emph{variables} instead of numbers for constants) so that the abstracted representation can be generalized beyond the specific case of interest.  This is one of the powerful features of mathematics; a single mathematical statement representing the abstraction of a problem can actually \emph{represent} a whole class of similar problems.  As an example, think about the problem where there is a small deflection ($d$) of a end-loaded beam pinned at one end and free at the other (with the free end loaded).  With appropriate assumptions, the solution to this problem depends on the magnitude of the load (say, $W$), the length of the beam ($L$), the Young's modulus ($E$), and the moment of inertia of the cross section of the beam ($I$).  Under these circumstances, a whole \emph{class} of problems can be solved using mathematics.  To be concrete, the solution is of the form
    
    \begin{equation}
        d = \frac{W L^3}{3 E I}
    \end{equation}
    For all possible (and reasonable) values for $W$, $L$, $E$, and $I$.  This is very powerful statement that covers an enormous range of different physical situations.  This is part of why abstraction is an important concept; the results for an appropriately abstracted problem can be much more general than would be the case for specific instances of the problem (e.g., where $W$, $L$, $E$, and $I$ were all fixed values).   
    
    When possible, the abstraction step should include explicit statements about the assumptions that are imposed.  This step is an important one, because it provides a recognition of which variables are considered to be important to the model, and which are assumed not to be important.  If model revision is required (see step 4, below), a recognition of the variables that were not included can make it easier to create a revised model that may have more predictive power.

    Stating a problem in a language is an essential component to problem solving.  Stating a problem via mathematics is an especially powerful process.  For example, it allows one to communicate, reasonably exactly, the statement of the problem (with attendant simplifications) to someone else.  Because the language of mathematics is formalized and simplified, it has high fidelity (one is quite clear on the meaning of a particular mathematical statement, once formulated), but not as expressive as most natural languages (i.e., one does not generally use the language of mathematics to write poetry).
    
    \item {\bf Model building via algorithmic processing.}   There is also a deductive process that must occur so that the problem stated can eventually be used to determine a \emph{solution}, assuming that one exists.  The solution process is example of an algorithm, and the mental process of actually working through the steps to the solution is, as mentioned in the introduction, \emph{algorithmic thinking}.\indexme{modeling!algorithmic thinking}  In reality, there is no need use mathematics specifically for problem solving.  One could, for example, use modern English as a language and accomplish the same thing.  However, the compactness of the symbols used in mathematics make it a much more convenient tool for problem solving.  Another important feature is that, because mathematics is a highly constrained language used for a specific purpose, it is much less ambiguous than modern English would be. 
    
    \item {\bf Computation and concrete realization}.\indexme{modeling!computation}  Once a problem has been stated, and a solution process identified, a particular solution can be computed.  In the example above for a beam deflection, the concrete computational step would be the generation of the formula for deflection, and then the subsequent substitution for values of $W$, $L$, $E$, and $I$ for the particular case of interest.  The result would be a concrete \emph{number} expressing the predicted deflection in some specified \emph{units}.  While the \emph{results} of a model computation are concrete (in the sense that one obtains an answer), the results may not correspond well to the even more important behavior of the physical problem being solved (i.e., the ``real world" application of your model).  Thus, there are several steps to concrete realization: first, one must actually use their model to compute (or measure) a result produced by the model.   Second, assuming that one's model corresponds to a real, physical system, then the result of the model should compare \emph{well enough} to reality that the model is acceptable.  Here \emph{well enough} indicates a level of fidelity that is required for a particular application, and cannot be specified independently from knowledge of the fidelity needs for the application.  As an example, one might want to be quite certain that a beam design will not fail, because failure generally would represents a serious consequence; thus, high fidelity in the modeling may be required.  In another model, say a model for how many treats you predict your dog will eat on your next walk, may have few negative consequences if you are in error (however, your dog may feel differently about this); thus, lower fidelity is acceptable.
    
    \item {\bf Revision.} \indexme{modeling!revision}  A model is complete if it meets some pre-determined measure of fidelity for the application intended as discussed above.  Often, this means matching data, although other qualifications are possible.  We might find that our first efforts at modeling do not meet our goals for fidelity.  If this is the case, then the assumptions / neglected variables that were identified during the abstraction step can be very useful in helping to generate a revised model that has more predictive power when assessed by the the specific measures of fidelity that are used. There may be a number of iterations through steps 1-4 above before a successful model is attained.

    The revision step may also involve a process known as model \emph{validation}.  Model validation is the process of determining whether or not a particular model is appropriate for the intended purposes; thus, it has an obvious connection with the \emph{revision} process.  Usually, model validation requires that the model meet some performance metric within a prescribed level of accuracy or tolerance.  As an example, a performance metric might be that a conceptual model describes some observed data with a sum of squared deviations between the two being less than some value, $V$, determined by an external constraint (e.g., cost/benefit of the design, safety needs, etc.)  Model validation is an important step in the modeling process for many applications, and it represents an entire sub-discipline within modeling.
    
\end{enumerate}

The steps described above represent the process of \emph{modeling}.
Scientists and engineers use models all of the time.  Any time that we make some \emph{simple} system that is intended to represent the essential features of a more \emph{complex} system, we are making a model.  The particular kind of model depends upon the tools available and the kind of system being represented.

\section{Models and Units}

Finally, this discussion remains incomplete without addressing the problem of \emph{units}.\indexme{modeling!units}  While pure mathematics exists quite happily without the imposition of units, when modeling physical systems, the use of units becomes essential.  While most of us already know what is meant by the term \emph{units}, it is still useful to define them.  The use of a system of \emph{units} \indexme{units} allows one to identify and communicate the kinds of physical quantities that are being modeled.  It also allows one to establish the \emph{magnitude} of the quantities in some standardized sense. Formally, the definition might be given as follows.

\begin{definition}[base units]\indexme{units}\indexme{units!base units}
 A set of \emph{base units} provides a standardized name and magnitude for various kinds of physical quantities.  These form the basis of a measuring system, where other physical quantities can be expressed as multiples of the set of base units. 
\end{definition}

An important concept here is that there is a set of \emph{base} units which all other \emph{units} are constructed from.  For example, in the International System of Units (SI) \indexme{International System of Units (SI)}, the fundamental units are:

\begin{enumerate}[wide]
    \item The meter (symbol: \si{m}), used to measure length.
    \item The kilogram (symbol: \si{kg}), used to measure mass.
    \item The second (symbol: \si{s}), used to measure time.
    \item The ampere (symbol: \si{A}), used to measure electric current.
    \item The kelvin (symbol: \si{K}), used to measure temperature.
    \item The mole (symbol: \si{mol}), used to measure amount of substance or particles in matter.
    \item The candela (symbol: \si{cd}), used to measure light intensity.
\end{enumerate}
There are some curiosities involved in the \emph{writing} of units.  While many units are named after people, the full names of such units are not capitalized.  However, usually the \emph{symbols} for the associated units \emph{are} capitalized.  Thus, we have $1~\si{ampere}=1~\si{A}$.  Does this make sense?  Well, the question is not relevant: a formalism has been established, and thus it is clear \emph{what} to do here (even if it is not clear \emph{why} it is this way!)  The inclined student can examine the National Institute of Standards and Technology (NIST) Office of Weights and Measures (\url{https://www.nist.gov/pml/owm/writing-si-metric-system-units}) to learn more about the vagaries of unit names and conventions.

The primary point to be made here is that this \emph{base set} of units can describe all other possible physical quantities that can exist.  As a concrete example, we can think about the concept of \emph{voltage} (which is analogous in many ways to \emph{pressure} in fluid systems).  While we are used to expressing voltage in Volts (e.g., in North America, our single-phase household wiring is approximately 120 volts operating at a frequency of  60 \si{s^{-1}}).  However, the unit of volts is \emph{not} fundamental.  Instead, note the following.  One volt is equal to one joule per coulomb.  A joule is equal to a newton-meter.  A coulomb is equal to one Ampere-second.  Thus
\begin{equation}
1~ \si{volt} = 1~\frac{\si{J}}{\si{C}} = 1~ \frac{\si{N \cdot m}}{\si{A\cdot s}} = 1~ \frac{\si{(kg\cdot m/s^2) \cdot m}}{\si{A\cdot s}} = 1~ \frac{\si{kg\cdot m^2}}{\si{A\cdot s^3}}
\label{units}
\end{equation}
On the farthest right expression in Eq.~\eqref{units}, the volt is expressed in its base SI units.  The list of seven units provided above are sufficient to express \emph{all possible physical quantities that are known today}.  While one might be tempted to believe that ``units" would not continue to be researched today, it turns out that it is still a subject of much discussion!  There are efforts currently to establish all units in terms of ``universal constants" (e.g., the gravitational constant, $G$, or Plank's constant $\hslash$), and thus the topic is still one of evolving research (e.g., see \citet{borde2005base}).

\section{Mathematics, Models, Determinism, and Complexity}

When developing models, one usually makes certain kinds of very basic assumptions about the ultimate utility or purpose of the model.  For example, quite frequently we \emph{assume} that the models that we generate will allow us to make some kind of \emph{prediction}; this is indeed often the purpose that we bother to generate models in the first place.

While models have been successfully used by humankind for thousands of years, in (relatively speaking) more recent times, we have discovered that sometimes the models that we formulate provide us with information that is qualitatively (and quantitatively) different from what we generally have come to expect from models.  To be more specific, we learned that not all models can be expected to be \emph{deterministic}\indexme{determinism}\indexme{modeling!determinism} (a word that will be defined in more detail below).  Instead, we have found that the results of some models are so sensitive to their model parameters (for example a physical property, such as density) or to the conditions at their initial state or on boundaries, that in practice we cannot make deterministic predictions.  This does not mean that such models are not useful.  However, it does mean that the kinds of information that one expects from such models is different than what one expects from models that behave deterministically.  The concepts of determinism, and the related concept of complexity, are discussed in the material following. 

\subsection{Determinism in Models}
Determinism was a framework developed by Greek philosophers during the 7th and 6th centuries BCE.   In short, the philosophy was that cause and effect are bounded together in a way that could (in principle) be understood. It assumes that if an observer has \emph{sufficient information} about an object, that such an observer might be able to predict every consequent move of that object. In the mid- to late 1600's, famous philosophers and scientists (such as Ren/'e Descartes and Sir Isaac Newton) codified this idea as being a natural part of the physical world.  However, it was the mathematician-scientist Pierre-Simon Laplace who advanced the idea of determinism into the modern scientific discourse.  Laplace stated emphatically \citep{laplace1825essai}

\begin{quote}
We ought then to regard the present state of the universe as
the effect of its anterior state and as the cause of the one
which is to follow. Given for one instant an intelligence
which could comprehend all the forces by which nature is
animated and the respective situation of the beings who
compose it—an intelligence sufficiently vast to submit these
data to analysis—it would embrace in the same formula the
motions of the greatest bodies of the universe and those of
the lightest atom; for it, nothing would be uncertain and the
future, as the past, would be present to its eyes.
\end{quote}

This perspective was viewed with nearly the status of a physical \emph{law}.  However, in the late 1800s, Henri Poincar\'e suggested to the world that perhaps systems were not always as deterministic as had been understood.  Poincar\'e found that there were systems of differential equations (which, in turn, were models of physical systems) where the final solution was quite sensitive to the initial conditions.  Poincar\'e's observations led him to the following statement

\begin{quote}
A very small cause, which eludes us, determines a considerable
effect that we cannot fail to see, and so we say that
this effect is due to chance. If we knew exactly the laws of
nature and the state of the universe at the initial moment,
we could accurately predict the state of the same universe
at a subsequent moment\ldots But this is not always so,
and small differences in the initial conditions may generate
very large differences in the final phenomena. A small error
in the former will lead to an enormous error in the latter.
Prediction then becomes impossible, and we have a random
phenomenon.
\end{quote}

This observations has been examined and refined over the years; today we call the study of such sensitive systems (among other names) \emph{chaos theory}.  In short, what this means is that  \emph{both mathematical and physical systems} exist where, despite our best efforts and most advanced technology, we can never know enough about the system such that its future behavior is knowable with certainty.  While today we are somewhat familiar with this idea because of our understanding (even superficially) of quantum mechanics, the idea that even non-quantum-mechanical systems (or \emph{macroscopic systems}) can behave in ways we are not able to predict still comes as a surprise to many people.

The non-deterministic behavior of physical and mathematical systems led to somewhat of a revolution in the way that physicists, scientists, and mathematicians thought about and represented models of processes that occur in our world.  What we once thought of as \emph{knowable} if given the proper kinds and amount of information was now understood to be knowable only within certain bounds.  Even without the strangeness embodied in quantum mechanical systems, we found that some large-scale physical systems (such as the weather) could be so sensitive to the initial conditions or other physical properties, that their long-term behavior could not be modeled with certainty.  Primarily, this kind of behavior was observed for \emph{nonlinear} models.  The concept of nonlinearity is explained in more detail in the next chapter.  In short, however, one can think of nonlinearity in this application as meaning that small changes to the system can yield arbitrarily large responses.  This does not happen in so-called \emph{linear} systems, where small changes to the system always yield a response that is small in some sense. 

The gradual unraveling of the idea that the physical world and associated models of it should be \emph{deterministic} represented a significant change in the way that scientists and mathematicians thought about the models they used.  Because part of the purpose of a model is to describe the complicated universe by collecting and utilizing understandable ideas, it seemed that the lack of determinism created a problem.  This problem was specifically that \emph{some} reasonable models of physical processes yielded behavior that was inherently non-predictive in some ways.  The consideration of this problem led to new research that eventually became known as \emph{complexity science}\citep{phelan2001complexity}, or simply \emph{complexity}.

\subsection{Complexity in Models}
When we say something is \emph{complex} in English, we usually mean that it in the sense of ``not easily understood".  This could be a technical comment, or one that represents other factors involved in the real world (e.g., it would be safe to say that the geopolitical situation in the Middle East is complex, without intending any technical notion of complexity).  Similarly, in mathematics complex can also mean involving numbers that have an imaginary component.  In modeling, we usually mean neither of these things by the word \emph{complex} (although there are certainly cases where more than one of definitions of the word \emph{complex} might apply!).  While we will study only a few models that can be considered to be complex in this text, it is important to understand what complexity is as it relates to modeling.  \indexme{complexity}\indexme{modeling!complexity}

The \emph{complexity} of a model means different things in different disciplines.  However, nearly all definitions of the word, as it relates to modeling, involve answering the following question: \emph{how difficult is it to explain our model as it relates to predicting the behavior we are attempting to represent?}  One may be more concrete by replacing ``how difficult is it" with ``how much information is needed" (however, this also requires a subsequent definition of information, which we are not prepared to detail here!).     A key concept in answering this question involves an attempt to define \emph{regularities} \citep{phelan2001complexity} in our model and its solution.  Complexity science introduces  new ways to identify and study the regularities of even non-deterministic systems; importantly, these methods were represent new ways of looking at problems \citep{phelan2001complexity}.  While one generally thinks of finding regularities \indexme{modeling!regularity} in systems that are non-deterministic, this idea can also be useful for systems that are deterministic.  As an example, suppose that we wanted to measure the weight of 10,000 ball bearings.  While one deterministic option would be to \emph{measure} every bearing, one could introduce other strategies depending on what kind of information was actually needed.  If the \emph{actual} weight of every bearing was needed, we would indeed have to weigh each one (and keep track of them with labels).  However, suppose we needed only to know the weight with a particular confidence (say, with 95\% confidence)?  Then, we could measure the statistics of a randomly selected subset of the bearings, and make a statistical inference about the average and standard deviation (with enough measurements to assure 95\% accuracy).  Of course, when we do this, we are making a strong assumption about the \emph{regularity} of the system.  For example, in the case of the bearings, we are assuming that the population of all 10,000 bearings was regular enough such that standard statistics could be used to summarize them.  This is an example of a modeling assumption of system regularity that allows us to propose a deterministic model (even though the resulting model is statistical, it is deterministic in the sense we have defined for \emph{models}) for measuring the bearing weights without measuring each bearing individually.  Thus, the assumption and search for regularity \indexme{model!regularity} in models can be a very powerful tool!

So, how do we define the word complexity.  Again, there is no one universally agreed upon definition.  For the purposes of this text, however, we will define complexity in modeling as follows.
When applied to modeling, the word \emph{complex} (or \emph{complexity}) implies at least one of the following is true.\indexme{complexity!properties of}

\begin{enumerate}[wide]
    \item {\bf Large number of variables.} \indexme{model!large number of variables} Some models have so many \emph{degrees of freedom} (i.e., variables) that it is not practical or possible to solve the system.  An example here is the model of a mole of an ideal gas as a collection particles obeying classical mechanics (Newton's Laws) at a fixed temperature.  While such a model is indeed both deterministic and even accurate for certain gases and states, not many would find the tracking of 3 momentum variables and 3 position variables for each molecule (making a total of $(3+3)\times6.02\times10^{23}\approx 3.6\times10^{24}$ variables!) to be a very practical computation to make.  However, using statistical tools, and assuming certain kinds of regularity are manifest by the velocities of gas molecules, it is possible to actually derive the ideal gas law from the consideration of essentially Newton's laws applied to roughly $10^{23}$ particles.  The assumption of \emph{regularity} has clear power here: instead of needing on the order of $10^{24}$ variables to describe an ideal gas, we can use the classical ideal gas law $pV=nRT$.  For our case of a single mole of gas at a fixed temperature, the regularity assumption allows us to describe the system with variables ($p$, $V$) and one constant ($R$).  This represents a significant reduction in the difficulty we would have in explaining (or the ``information embodies by") our  model of the system!  \\
    
    \item {\bf Non-deterministic behavior.}  A second meaning for complexity in a model is the condition where the model behavior is so sensitive to initial, boundary, or other conditions (e.g., the exact size of a particular parameter) that it is nearly impossible to predict the behavior of the model in a classical sense.  Such models are called \emph{chaotic}\indexme{chaotic systems}, as discussed above.\indexme{model!chaotic} Colloquially, we can think of such models as ones that will predict a very different outcome with even a very tiny change in the conditions (initial condition, parameters, etc) the describe the model.  Because we seldom know initial, boundary, or parameters for real systems with high accuracy, the resulting model has very little predictive value because small errors in the conditions can yield wildly different results.  In short, many systems in science and engineering subscribed to the rule that ``small perturbations create small effects" (and this is always true in a sense for linear problems).  For chaotic problems, this rule is no longer true.\\
    
    \item {\bf Emergent behavior.}  This concept is not easy to describe concretely, but the term \emph{emergent}\indexme{emergent behavior} in describing the behavior of some models (usually nonlinear ones) is now so commonplace that it is important to understand what the essential features of \emph{emergent behavior} is.  In plain English, behavior describes the process where a model starts from a condition that is not very ``interesting" to one that has high structure, information content, or behavior that would not necessarily be expected.  In other words, the model transitions from one archetype to another in a surprising way.  While this is hardly a concrete definition, no concrete definitions seem to exist.  One of the first detailed studies for this kind of phenomena was conducted by the famous British mathematician Alan Turing (who is primarily known for his groundbreaking work on computers in the 1940s and 1950s).  In Turing's applications, the emergent behavior was given by a reaction-diffusion equation, and was proposed as  a possible explanation for how certain patterns arise in biology (e.g., the spots on a cheetah, or the pattern on some shells).  These patterns are sometimes called Turing patterns in honor of the discovery.   For the interested reader, a review of emergent behavior in general is given by \citet{krause_2018}.\indexme{model!emergent behavior}
\end{enumerate}
In each of these kinds of complexity, the search for regularity in the models and their solutions can yield understanding in areas where it would be otherwise lacking.  

\begin{figure}[t]
\sidecaption[t]
\centering
\includegraphics[scale=.5]{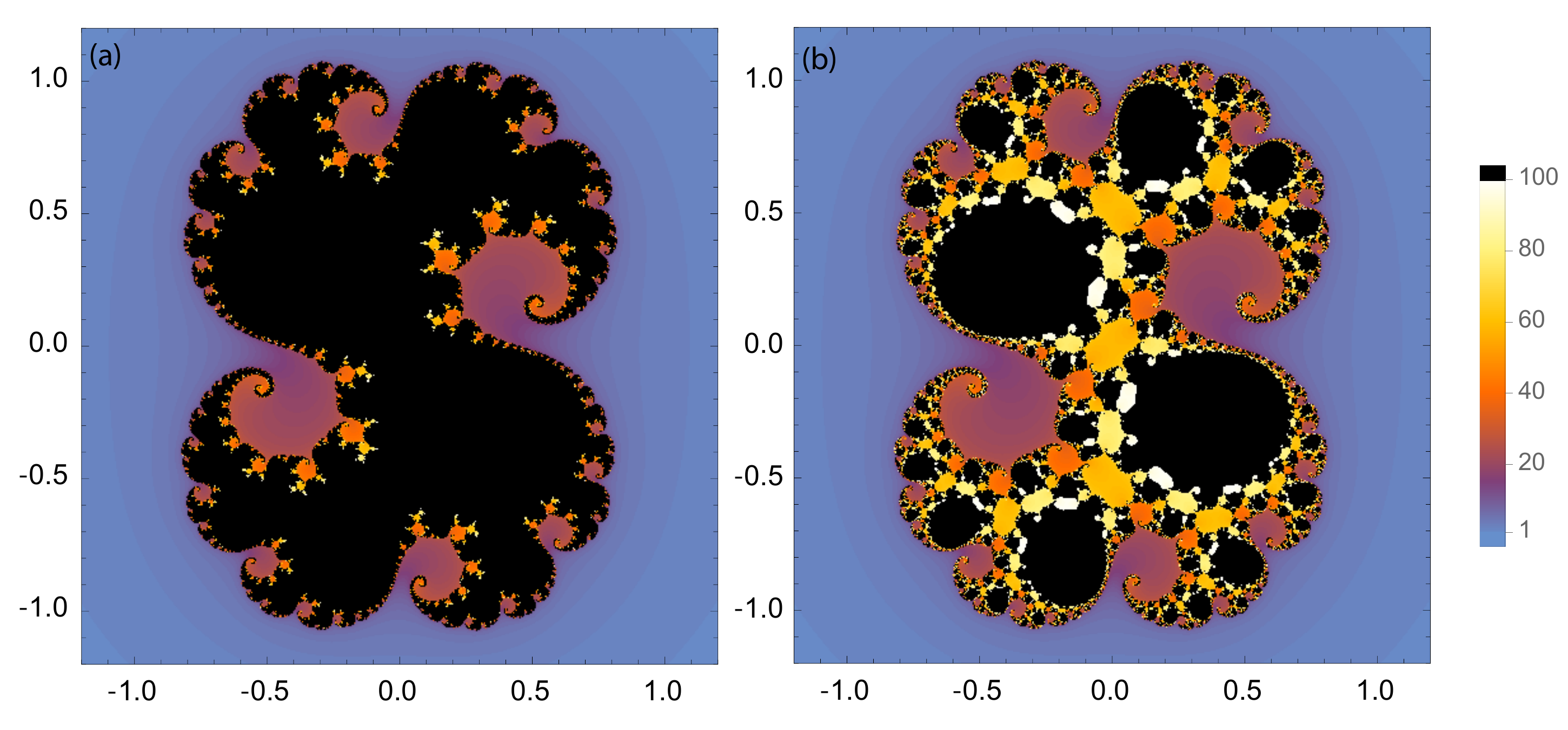}
\caption{Two representations of the Julia set. This set is determined iteratively for $z_j=x_j+i y_j$ by selecting a grid of $(x,y)$ of starting values, and then computing $z_{j+1}=z_j^2+c$ until convergence.  The value of $c$ is a constant number of the form $c=c_x + i c_y$.  Selecting different values of $c$ selects different Julia sets.  The images above are generated by computing the magnitude of each converged value of $z_\infty=x_\infty + i y_\infty$ on the grid, and computing its magnitude, $m$, where $m=(x_\infty^2+y_\infty^2)^{\tfrac{1}{2}}$. (a) $c=0.274-i 0.008$. (b) $c=0.275-i 0.008$.  For these images, a 0.4\% change in the real part of $c$ yielded dramatically different magnitudes for the resulting function.}\indexme{Julia set}
\label{julia}       
\end{figure}

The fairly well known example of complexity in mathematics can be found in the example of certain fractals.  As a specific example, one can point to the \emph{Julia set}, which is a chaotic function of the initial point selected for iteration.  Most of us have seen images of the Julia set, and these images have almost reached the level of pop culture; a cursory examination of Fig.~\ref{julia} might suggest that it is a familiar looking plot.  The Julia set is usually computed via iteration.  Even a small change in the initial condition for these iterations will yield a dramatically different set of results for the ``shape" of the resulting Julia set.  An example of high sensitivity for the Julia set is given in Fig.~\ref{julia}.  Here, a change of only 0.4\% of one of the parameters in the model led to dramatically magnitudes (represented by the color) for the set.\indexme{chaotic systems!Julia set}

As a second example of complexity, we can consider the third definition provided regarding the concept of emergent behavior.  Nonlinear reaction-diffusion systems are one of the oldest examples of systems which can show emergent behavior.  The diffusion reaction system given by Eqs.~\eqref{turingpat1}-\eqref{turingpat2} are known as a Turing model with Fitzhugh-Nagumo reactions; these equations represent the diffusive mixing and reaction of two chemical species (where the species concentrations are represented by the variables $u_1$ and $u_2$) that have the capacity to create self-organized non-homogeneous patterns in space in the steady state.  \indexme{emergent behavior!Turing models}

\begin{align}
    \frac{\partial u_1}{\partial t} &= D_1 \left(\frac{\partial^2 u_1}{\partial x^2}+\frac{\partial^2 u_1}{\partial y^2} \right)+(u_1 -u_2) -u_1^3 +\kappa \label{turingpat1} \\ 
    \frac{\partial u_2}{\partial t} &= D_2 \left(\frac{\partial^2 u_2}{\partial x^2}+\frac{\partial^2 u_2}{\partial y^2} \right)+(u_1 -u_2)
    \label{turingpat2}
\end{align}

\begin{figure}[t]
\sidecaption[t]
\centering
\includegraphics[scale=.5]{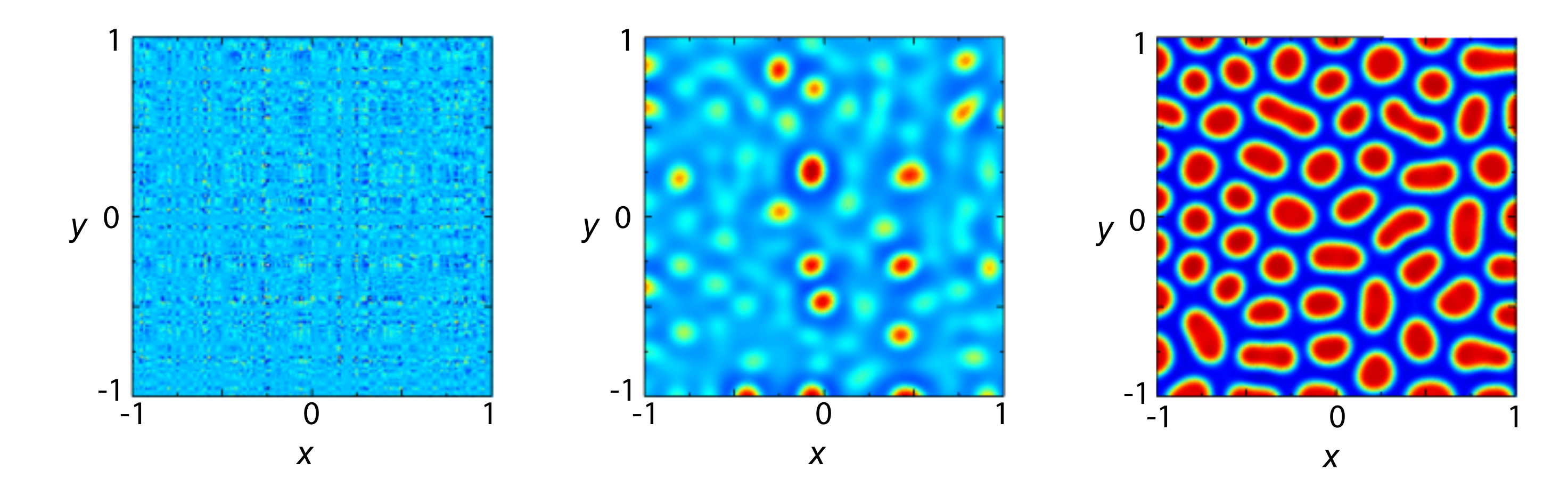}
\caption{Three spatial distributions of species $u_1$ reaction-diffusion system given by Eqs.~\eqref{turingpat1}-\eqref{turingpat2}.  These images represent the concentration of species $u_1$ for $t=0$, $100$, and $1000$ s, and illustrate emergent pattern formation.  (Left) At $t=0$, the initial configuration of the two species is a random function in space. (Middle) Pattern at $t=100$ s.  (Right) Pattern at $t=1000$ s, which represents near steady-state conditions.  The model was run with the following parameters: $D_1 =0.00028$ cm$^2$/s, $D_2 =0.005$ cm$^2$/s, $\tau = 0.1$, and $\kappa = -0.5$ mmol/(cm$^3\cdot$s).  For the color scale, red indicates a high concentration, whereas dark blue represents a low concentration.  ``Turing Bifurcation" images by Hubodeker are licensed under CC-BY-SA-3.0.}
\label{turing}       
\end{figure}

In Fig.~\ref{turing}, the spatial distribution of $u_1$ is plotted for three times.  At time zero (left panel),  the concentration of $u_1$ (indicated by the color) is random in space.  At the intermediate time (middle panel) one begins to observe spatial structure forming due the combination of diffusion and nonlinear reaction. Finally, at the near-steady-state condition (right panel), a spatial pattern has clearly formed showing high-concentration islands of $u_1$.  Although the conditions for $u_2$ are not plotted, they have similar initial conditions and evolution.

\subsection{Simplicity in Models}\indexme{modeling!simplicity}\indexme{simplicity}

Simplicity when applied to modeling is not the converse of ``complexity".  Model simplicity pertains to the idea that the most effective models impose only what is necessary to capture the phenomena modeled, and no more.  It turns out that this idea has been around in the philosophy of science for a very long time.  The idea that models should be as simple as possible has been around since at least the early 1200s \citep{franklin2001science}.  The principle is often called Occam's razor (where razor is a fanciful terminology indicating that the idea shaves away the unnecessary) after an English Franciscan friar named William of Occam.  While Occam did not originate the idea, he did make it a popular one in the philosophy of science.  One version of his statement (actually appearing in his works) is given by
\begin{quote}
    Plurality must never be posited without necessity
\end{quote}
In short, the idea is this: Given two competing models for a system, each with similar predictive power, the simpler model should be preferred.  Sometimes this idea is also called the \emph{principle of parsimony}, for obvious reasons.\indexme{parsimony}\indexme{modeling!Occam's razor}\indexme{Occam's razor}\indexme{modeling!parsimony}\indexme{parsimony}

The statement of Occam's razor has been repeated many times over the decades.  Some famous examples are as follows.

\begin{itemize}
    \item Isaac Newton stated in his famous text the \emph{Principia} ``We are to admit no more causes of natural things than such as are both true and sufficient to explain their appearances."
    \item Einstein can be quoted as stating (in a lecture at the University of Oxford, 1933) ``It can scarcely be denied that the supreme goal of all theory is to make the irreducible basic elements as simple and as few as possible without having to surrender the adequate representation  of a single  datum of  experience."
    \item A composer named Roger Sessions famously misquoted Einstein by stating ``I also remember a remark of Albert Einstein, which certainly applies to music. [Einstein] said, in effect, that everything should be as simple as it can be but not simpler!". 
    \item A memorable version of the statement was given by Nobel-prize-winning medical researcher Theodore Woodward.  His interpretation of the concept was stated by the aphorism ``When you hear hoofbeats, think horses, not zebras."  While the principle was intended to be applied to diagnosis of medical conditions, it is a colorful way of summarizing the core idea of Occam's razor.
\end{itemize}

While all of this sounds somewhat qualitative, it turns out that it can be adapted for use in a \emph{quantitative} environment. In the mid 1950s a statistician named Edwin Jaynes surprised the physics community by showing that Occam's razor, stated as the \emph{principle of maximum entropy} could allow one to derive statistical mechanical results with accuracy.  In conventional language, Jaynes was able to show that applying the statistical notion of maximum entropy in certain statistical mechanical computations was equivalent to inserting the minimum amount of assumptions.  Jaynes was able to show that the conventional assumptions required in classical statistical mechanics could be replaced by maximizing the entropy under known constraints.  Doing so would generate a solution that was \emph{maximally noncommital}\indexme{parsimony!maximally noncommital} to unjustified assumptions.  This represented somewhat of a revolution in the way that we think about modeling systems, and turned out to be very much a practical realization of Occam's razor.  Thus, what started out as a philosophical principle (Occam's razor) turned out to have utility in a quantifiable way.

\subsection{Overfitting Models}\indexme{modeling!overfitting}\indexme{overfitting}

An occurrence in modeling that is related to the concept of parsimony is the problem of overfitting a model.  The term \emph{overfitting} comes from the related disciplines of statistics and data science.   Conceptually, the idea is simple.  Overfitting means to model a particular set of phenomenological data very accurately, but in a way that does not generalize well for representing other data sets of the same phenomenon.  Frequently, this arises because of a lack of obeying Occam's razor.  It is easiest to explain this problem through an example. 

Drag on bodies moving in fluids has been studied for hundreds of years, and frequently one finds that there is a relationship between the drag force and the square of the velocity.  Suppose we run two independent sets of experiments where we measure the drag on a body in a wind tunnel for various velocities; we call these two sets of experiments A and B.  We use experiment A to calibrate a model that allows us to predict the drag for a velocity.  Then, we use our calibrated model from experiment A to predict the behavior of the data measured in experiment B.

\begin{figure}[t]
\sidecaption[t]
\centering
\includegraphics[scale=.3]{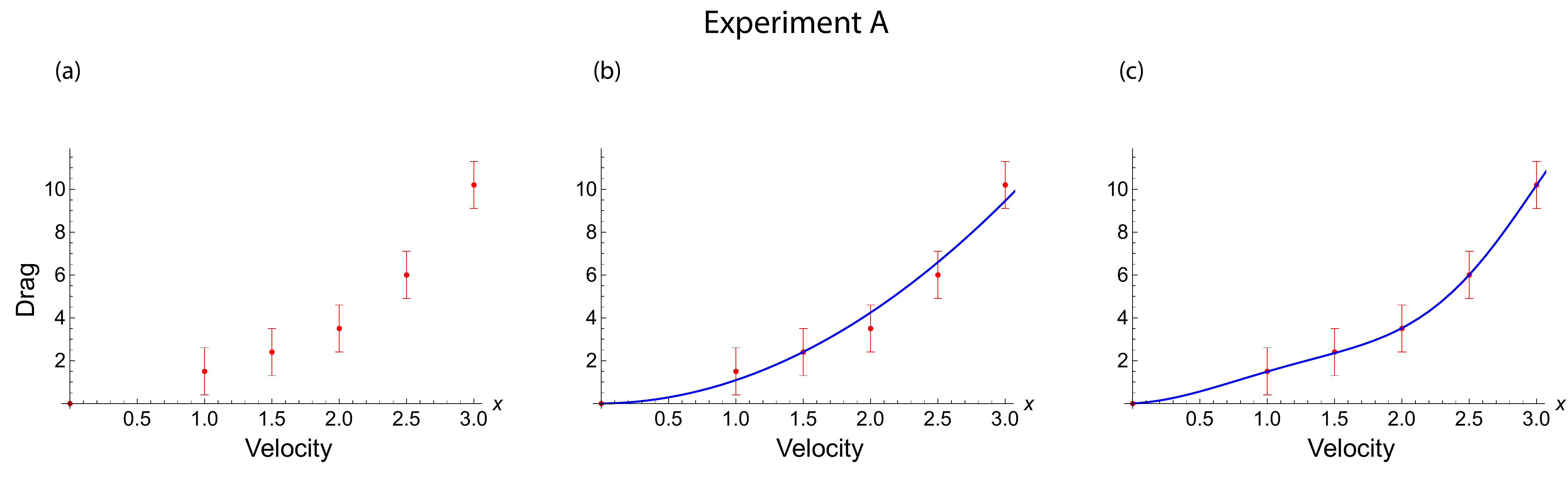}
\caption{Results of drag experiments, experiment A.  (a) The drag data for experiment A.  (b) A quadratic least-squares fit (constrained to be positive) of the form $F_D=1.05 U^2$.  (c) A fifth-order polynomial that fits the 6 data points perfectly; $F_D= -0.137553 x^5+1.17801 x^4-2.89406 x^3+2.66579 x^2+0.785579 x$}
\label{overfit1}       
\end{figure}

In Fig.~\ref{overfit1}, we have plotted the original data (complete with a 95\% confidence interval), and we have fit two different polynomials to the data.  First, our experience with drag laws would suggest that a quadratic should provide a reasonable fit to the data.  Our fitted quadratic (constrained to be everywhere positive) seems to match the data reasonably well, and is acceptable in that it falls within the 95\% confidence interval for each data point.  In Fig.~\ref{overfit1}(c), we show the results of a \emph{fifth-order} polynomial (constrained to be everywhere positive) fit to the data.  Including the origin, there are six data points total.  It is always possible to find a polynomial of order $n$ that will go through $n+1$ points.  Thus, the fit in this case is perfect.  But is it an \emph{optimal} fit?  This would depend, in some ways, upon our goals.  First, because we know that drag laws are generally quadratic in the velocity, a fifth-order polynomial goes against our ``prior data" indicating that we might expect a quadratic form.  And, while the fifth-order polynomial certainly fits the data perfectly, one has to wonder if the additional variation in the curvature represents anything \emph{physical}, or is simply a very accurate representation of the experimental error that we know to be part of the experiment.

\begin{figure}[t]
\sidecaption[t]
\centering
\includegraphics[scale=.3]{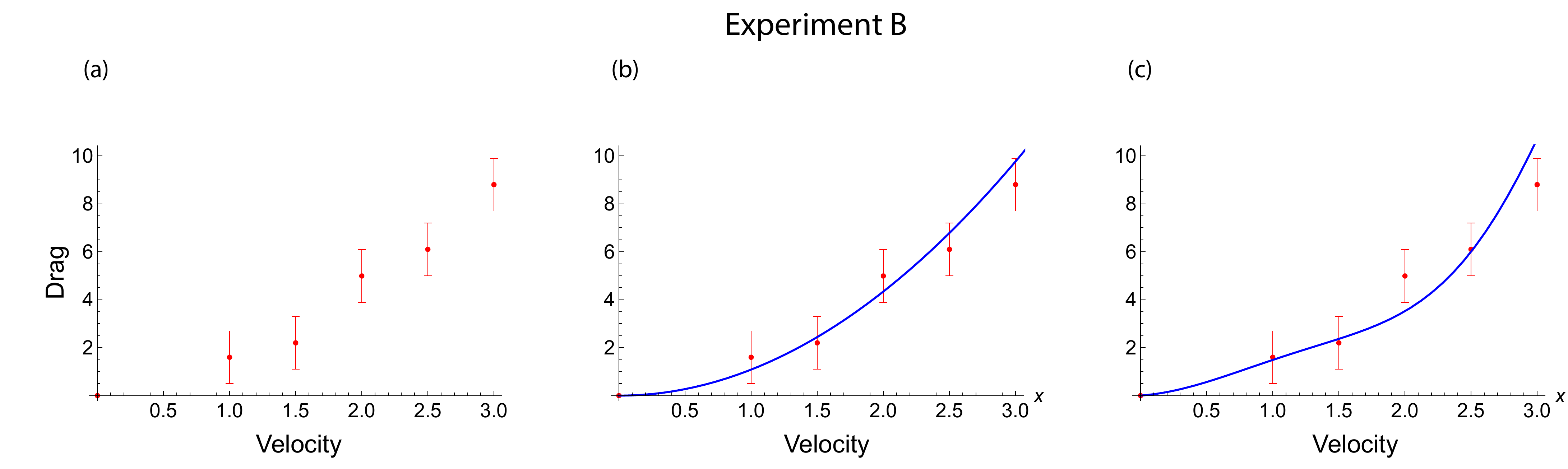}
\caption{Results of drag experiments, experiment B.  (a) The drag data for experiment A.  (b) The quadratic least-squares fit (constrained to be positive) of the form fit to model A, and tested on model B.  (c) The fifth-order polynomial fit from model A, and tested on model B.  One can see that for the data of experiment B, the quadratic model actually fits the data better than the fifth-order polynomial model. This is a result of over-fitting.}
\label{overfit2}       
\end{figure}

One way that we can begin to answer these questions is to look to see how well this model \emph{generalizes}.  In other words, we have two models derived from experiment A.  If these models generalize well, that means that if we take another set of experimental data (replete with its own measurement errors), the model will still provide a reasonable fit to the data (i.e., it ``explains" the data).  In Fig.~\ref{overfit2} we show the results of using the models from experiment A to fit the data from experiment B.  In Fig.~\ref{overfit2}(a) we show the data from experiment B, and in parts (b) and (c) we show the fits of the quadratic and fifth-order polynomial models to the data.  Here, it is no longer evident that the fifth-order polynomial does a better job of fitting the data.  In fact, because of the additional components of the polynomial, it actually explains the data less accurately than the quadratic model (i.e., there are two data points where the fifth-order model fails to be within the 95\% confidence interval of the data points).  
%

This is a good example of the problem of \emph{overfitting}.  What has happened here is that we have added unnecessary degrees of freedom (i.e., additional constants and their attendant polynomial functions $x^3$, $x^4$, $x^5$).  While these additional degrees of freedom allowed us to generate a more accurate model for a single instance of our drag data, the model failed to generalize well.  This is because the model was \emph{overly specific} to the data that we had.  We needlessly added additional degrees of freedom to our model, in contrast with Occam's razor.  This occurred because we ignored prior data (i.e., we knew that drag models are generally quadratic).  However, using the formalism of \emph{maximum entropy}, we could also have looked at many such experiments, and found that the errors induced by the fifth order model were not distributed by a normal distribution (which would be expected for the quadratic fit), but rather by some skewed distribution (because it systematically under-predicts the drag in the middle off the velocity range).  One can show that skewed distributions actually have \emph{less entropy} than do normal distributions.  In short, the more entropy a distribution has, the less parameters it takes to describe the distribution.  Thus, Occam's razor would suggest that we choose the distribution with the greatest entropy (and least assumptions), which, in this case, is the quadratic model.  It is interesting that Occam's razor has something quantitative to say about fitting models.  In fact, such entropy techniques have become important tools in the burgeoning field of machine learning.

\section{An Example of Model Building: Attneave's Cat}\indexme{cat!Attneave's}

The topic of how geometric shapes are processed by the human visual cortex has been a topic of both qualitative and quantitative research for many decades \citep{kohler1949cortical}, and it continues to be an active area of research for a variety of reasons ranging from the psychological (e.g., understanding the structure-function of the brain) to the practical (predicting how best to present information visually to an active user, such as in a ``heads up" display).  

In 1954, a psychological researcher named Fred Attneave proposed a model where he postulated that in 2-dimensional contour figures, the regions of high curvature would contain the greatest amount of information.  While this theory is still being discussed in the literature \citep{torfs2010,dewinter2008}, it appears to be one component of a possibly more complex understanding of human perception and vision.

It is interesting to examine the steps in problem solving listed in Section \ref{what?} as applied to Attneave's cat.  In the following list, thoughts about each of the steps is discussed. 

\begin{figure}[t]
\sidecaption[t]
\centering
\includegraphics[scale=.28]{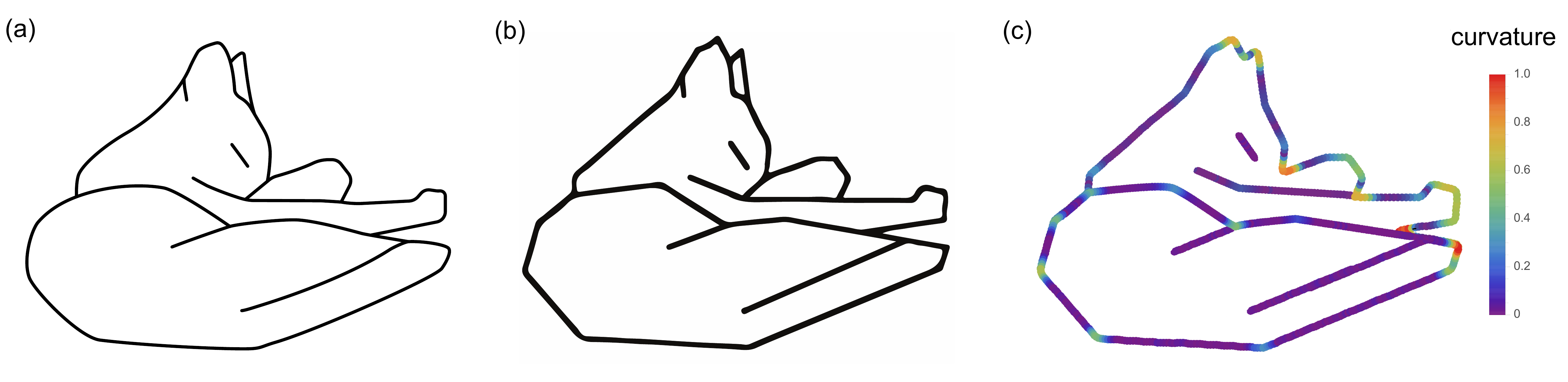}
\caption{Attneave's cat. (a) A smooth contour outline of a cat.  (b) The post-processed contour outline of Attneave's cat, expressed as straight lines connected by segments of high curvature.  This figure is generated by applying the algorithm described in the text.  (c) For reference, the rainbow color scale indicates the (normalized, absolute value) curvature  associated with the outline of the figure in (b).  Areas of high curvature are yellow to red; areas of low curvature are blue to purple.   }
\label{attneave}       
\end{figure}

\begin{enumerate}[wide]
    \item {\bf Abstraction.}  The abstraction step involved in Attneave's work is significant.  Attneave is not attempting to model a cat, but, rather, the image of a cat.  The abstraction step is the process by which Attneave asked (and attempted to answer) the question ``what is it that makes an image identifiable."  In particular, Attneave was interested in stripping away everything unnecessary so that only one or two elements could be focused on.  Ultimately, he chose to represent his images as a contour, with straight line segments joined by high-curvature segments.  The primary output from the abstraction step in this case was the formation of the hypothesis as a relatively simple  statement: regions of high curvature are the most informative in contour images.  In the abstraction step, it is important to identify one's assumptions and potentially additional or alternative variables that might be important in model revision (should that be necessary).  In Attneave's paper, this is addressed in the discussion and conclusions, where he suggests that other information (such as texture, color, contrast, and similarity) \emph{may potentially also be important}, but his results were focused specifically on curvature.  Identifying the simplifications imposed is an essential step in abstraction.  Models are rarely completed with a single effort; more frequently, a number of iterations are needed.  In order to facilitate these iterations, it is good to have a clear picture of what was included and what was left out of the model during the abstraction step.  One might then revisit these variables during revision, and decided that one or more of them is needed in order to produce a models with the fidelity desired.  
    
    \item {\bf Model building via an algorithmic process.}  The algorithm used to build the models is slightly more complex than that used simply to draw Attneave's cat.  First, Attneave wanted to test the hypothesis that curvature was indeed important.  This was done by experiments with a large number of subjects who looked at simple closed curves (think potato-shapes here), where subjects were asked to place points on the curve such that they could re-draw the curve using only the points selected.  From this, he was able to determine that his abstraction (segments with high curvature contain the most information) appeared to have validity.  This provided good evidence from which he could generate a ``model building" algorithm.  Without belaboring it too much, the algorithm would look something like the following:
    \begin{enumerate}
        \item For a given contour image, determine the regions of highest curvature.  Attneave did not investigate the process of generating contour images themselves in detail.  However, he did roughly sketch out how one might use edge detection to eliminate \emph{redundant information} to determine where contours for figures might be drawn.
        
        \item Identify the highest-curvature regions of the contour figure.  You may have learned in calculus that the curvature is the rate of change of the vector tangent to a curve as one translates along the curve.  For plane curves, this can be expressed by

        \begin{equation}
             \kappa ={\frac {|x'y''-y'x''|}{\left({x'}^{2}+{y'}^{2}\right)^{\frac {3}{2}}}}
        \end{equation}
        While Attneave did not compute the local curvature via this method, it is not difficult to compute.  In Fig.~\ref{attneave}(a), a version of Atteneave's cat is plotted.  In Fig.~\ref{attneave}(b) the absolute value of the curvature is plotted as the color scale along the outline of the figure known as Attneave's cat.  While this image represents the post-processed version of the contour image, it is nonetheless clear where the high-curvature regions of the figure reside.
        
        \item Keeping only the highest curvature regions in tact, replace the remainder of the components of the outline with straight-line segments connecting the endpoints of each high-curvature segment.  Again, the specifics of the method were not detailed by Attneave, but one could imagine using thresholding on the local curvature as a method to decide what to keep and what to eliminate.  Again, referencing Fig.~\ref{attneave}(b), it is clear that the purple regions represent nearly zero-curvature components of the outline.  
    \end{enumerate}
    
    While indeed a very simple algorithm, it is nonetheless an algorithm that can be applied with little uncertainty to nearly any contour image.  Some images (e.g., a circle, which has constant curvature) would be changed very little, whereas others (such as the image of the cat) would be significantly altered.
    
    \item {\bf Concrete realization.}  The realization process is, literally, just ``realizing" the model that has been constructed.  This means after identifying the appropriate abstractions and algorithm, and then using the algorithm to take the input data (in this case, a contour outline image) to generate an output (in this case, the simplified contour outline).  
    
    \item {\bf Revision.}  While no specific revision step is conducted in this particular modeling exercise, note that a number of important assumptions were identified by Attneave in his original (1954) paper.  Attneave's paper has been revisited many times by other researchers, and some of these have found that the neglected variables listed by Attneave were found to increase the descriptive ability of the models that they form.  So, although model revision was not a component of Attneave's model (at least, no discussion of alternative models was discussed in the paper; whether or not revisions were made in his initial work is not known), because he did describe a number of assumptions and neglected variables, other researchers have been able to extend his model to improve (to some extent) Attneave's original model.
    
\end{enumerate}

Attneave's cat provides a useful, and novel example of the problem solving process, with a focus on model development.  It is good to think about mathematics as being only one component of problems solving and modeling.  For many practical scientific and engineering problems, however, mathematics is the natural language for expressing models.  A solid background in mathematics will both provide useful modeling tools, as well as experience with the problem solving process in general.  Even if you never again use mathematics in problem solving (although this is hard to imagine for most scientists and engineers), the process of learning mathematics is, itself, a compelling method to practice the problem of problem solving.  Problem solving skills, once learned, have considerable transferability.  Learning and practicing mathematics will inevitably make you a better problem solver, regardless of the kind of problems that you are presented with.

\newpage
\section{Problems}

\begin{enumerate}[wide]
\item At the head of this chapter is the quote ``the best model of a cat is another cat, preferably the same cat".  While the quote is somewhat humorous (and has the benefit of mentioning cats), it also has substantial meaning.  In the context of thinking about systems, and abstractions of systems, write a few sentences explaining what this quote is meant to point out.

\item Another quote at the start of this chapter is ``all models are wrong; some are useful".  In the context of thinking about systems, and abstractions of systems, write a few sentences explaining what this quote is meant to point out. In particular, explain how it can be true that ``all models are wrong"?  Is there an example of any ``model" that is not wrong?  (Hint: refer to the quote repeated in problem 1.)
\setcounter{saveenumi}{\value{enumi}}
\end{enumerate}

\begin{minipage}[b]{0.4\textwidth}
\begin{enumerate}[wide,leftmargin=-0.45cm]
\setcounter{enumi}{\value{saveenumi}}
\item Circles have fascinated humankind for literally thousands of years.  Both the Babylonians and the Egyptians used approximations to pi, the ratio of the circumference of a circle to its diameter. But it was Archimedes of Syracuse (c.287 – c.212~ BC) who really got the ball rolling (so to speak) on the topic.  His idea was to inscribe regular polygons inside the circle.  With increasing numbers of sides, the area of these polygons would, he posited, get closer to the area of the circumscribing circle.  It is not difficult, using modern trigonometry, to show
\begin{equation}
A(r,n) = n r^2 \cos\left( \frac{\pi}{n}\right)  \sin\left( \frac{\pi}{n}\right)
\end{equation}

\vspace{0mm} 
\setcounter{saveenumi}{\value{enumi}}
\end{enumerate}
\end{minipage} 
\begin{minipage}[r]{0.6\textwidth}
\centering
\vspace{-65mm}
\includegraphics[scale=.6]{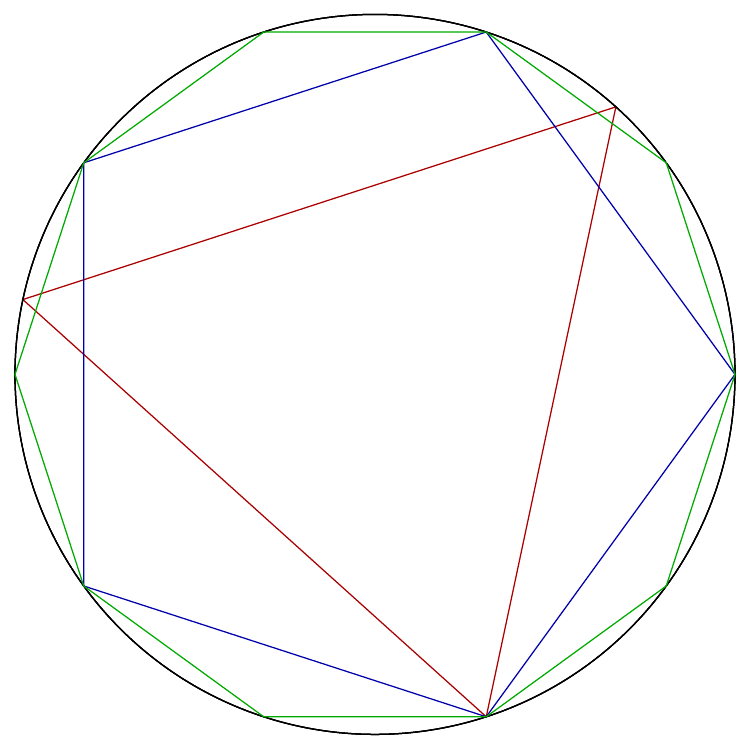}
\captionof{figure}{Polygons inscribed in a circle.}
\end{minipage}
\vspace{-5mm}

\noindent where $n\ge 3$ is the number of sides of the inscribed polygon.  We can think of $A(r,n)$ as being a \emph{model} for the area of a circle with radius $r$.  Compute the area of a circle of $r=1 \si{cm}$ for $n$ equal to 3, 4, 5, 10, 30 and 50, and compute the relative error $(A(1,n)-\pi)/\pi$ for those values.  Then answer the following questions.
\begin{enumerate}
\item[]

\begin{enumerate}
\item Does this model appear to be a sound and useful model?
\item Can you show that this model converges exactly to the area of a circle as $n$ gets arbitrarily large?  (Hint: for small $\theta$, $\sin(\theta)=\theta$).
\item Suppose you had a particular application in mind where you wanted to use $A(r,n)$ to estimate the area of a circle.  What kinds of considerations would you need to make in order to be able to choose the value of $n$ you would use?
\end{enumerate}
\end{enumerate}

\vspace{0mm}
\begin{enumerate}[wide]
\setcounter{enumi}{\value{saveenumi}}
\item \textbf{Chaotic Behavior.}  For certain kinds of phenomena, ranging from insect populations as a function of time to the wear rate of well-drilling bits \citep[cf.][]{li1975period}, the following model can relate the state at discrete intervals
\begin{equation}
x(n+1) = r x(n)\left[1-x(n) \right]
\label{logistics0}
\end{equation}
Where $r$ is a growth rate parameter, and $n$ indicates the number of the time period (each period with the same interval, say, $\Delta t$), and $x\in[0,1]$.  Equations of the form $x(n+1)=F(n)$ (like the equation above) relate the \emph{current state} of a system to the previous state.  They are called \emph{difference} equations; Eq.~\eqref{logistics0} in particular is sometimes called the \emph{logistics} equation.  

\hspace{-4mm}\begin{minipage}[l]{0.65\textwidth}
Starting in the 1960s, some researchers began to realize that such simple equations sometimes yielded surprisingly complicated (and \emph{complex} in the sense described earlier in the chapter) behavior.  In particular, there are values for $r$ where Eq.~\eqref{logistics0} becomes very sensitive to its initial value.  This expression is so simple that it can easily be coded up in a computer language, or even computed on a spreadsheet program.  Using a method of your choice, compute and plot the solutions up to $(n+1)=50$ given $x(1)$, In other words, you will need to compute a sequence of 49 values of Eq.~\eqref{logistics0} when given the first value.  For the parameters $x(1)$ and $r$ that should be used, refer to the table below. 
\end{minipage}
\hspace{2mm}\begin{minipage}[r]{0.35\textwidth}
    \centering
    \begin{tabular}{c|c|c}
      Case   & x(1) & r \\
      \hline
      1   &  0.5 & 1 \\
      2   &  0.5 & 2 \\
      3   &  0.5 & 3 \\
      4   &  0.5 & 3.99 \\
      5   &  0.5001 & 3.99 \\
    \end{tabular}
    \captionof{table}{Values for models using the Logistics difference equation.}
    \label{tab:my_label}
\end{minipage}

Once the solution is determined, compare cases 1 and 2.  What is different between them.  Case 3 has yet again different behavior-- how would you describe the behavior in time (i.e., steps) of Case 3?  Finally, note that Cases 4 and 5 are similar, except that the initial value is different by $0.001$ (or about 0.2\% of the initial value).  Do the two cases give similar or different behaviors in time?

\item \textbf{Models of probability.}  The concepts for understanding \emph{probabilistic} outcomes (e.g., the result of flipping a coin or tossing a six-sided die) was one of the first physical-conceptual systems for which intensive study of via models was applied.  The history of the topic indicates that improving ones odds at gambling was one of the motivations for serious study of the topic.  The statistician A. \citet{hald2003history} has stated
\begin{quote}
It was not until
the beginning of the 16th century that Italian mathematicians began to
discuss the odds of various outcomes of games of chance based on the
fundamental idea that the possible outcomes of a single game are equally
likely.
\end{quote}
While problems in probability theory can become quite challenging, here we propose a simple system for the purposes of thinking about \emph{modeling} the probabilities involved.  In particular, suppose one has an urn (the urn is, for some reason, the classical reservoir supposed in such problems!) that contains three balls: one yellow, one red, and one blue.  The question is: What is the \emph{probability} that the yellow ball is the first one picked if one selects three balls (one at a time) from the urn.

Now, this problem is really not about computing the probabilities involved (most people's experience would allow them to guess the probability to be $1/3$ without much additional thought), but to develop a model that \emph{explains the probability}.  While there are many different ``philosophies" regarding probability, we will not dwell on that here.  Instead, the following is proposed.  There are three balls, and they are selected one-at-a-time from the urn.  We assume that the order that the balls appear is important.  There are only a finite number of possible combinations that can be selected (e.g., blue, yellow, green is one; green, yellow, blue is another that is distinct from the first example).  \emph{If we assume that all outcomes of three balls are each equally likely}, then enumerating all possible outcomes will allow us to answer questions about the probabilities of each outcome.  For this problem, please do or answer the following.

\begin{enumerate}
\item Enumerate each possible outcome (hint: there are 6 total possible outcomes).  
\item Of the total number of outcomes, how many of them represent the cases where the yellow ball is drawn first?  
\item What is the ratio of the number of cases where the yellow ball is drawn first to the total number of cases?  
\item In this model, what is the primary modeling \emph{assumption} that allows us to actually compute the probabilities?  
\item Are there any additional assumptions regarding the process described that should be true so that the analysis is valid? 
\end{enumerate}
\label{probprob1}

\item Consider the situation described in problem \ref{probprob1}.  Now suppose that you want to know the probability of the selecting the balls in the following order: red, blue, yellow.  Can you describe the modeling process (i.e., write a short narrative) that explains how you arrived at your result?

\item Consider the situation described in problem \ref{probprob1}.  Now suppose that you want to know the probability of the following result: either the red \emph{or} the yellow ball is selected first.  What is the probability of this outcome?  Can you describe the modeling process (i.e., write a short narrative) that explains how you arrived at your result?  

\item Consider the situation described in problem \ref{probprob1}.  Now suppose that you want to know the probability of the following result: the first two balls selected are red and yellow, but in either order (yellow then red or red then yellow).  What is the probability associated with this outcome?  Can you describe the modeling process (i.e., write a short narrative) that explains how you arrived at your result?  

\item Consider the situation described in problem \ref{probprob1}.  Now suppose that you actually have a total of four balls: yellow, red, blue, and green.  Can you use the method described for \ref{probprob1} to determine the probability of the yellow ball being selected first?  While this method is convenient because it allows one to easily compute the probability of any possible outcome, can you see any disadvantages of this method (where each possible event is explicitly identified)?

\begin{minipage}[r]{0.6\textwidth}
\item Describe in words as concisely and clearly as you can the problem of \emph{overfitting} a data set.
\item Suppose you run two experiments that show the increase in temperature $(T-T_0)$ at ten evenly spaced times using infrared thermometry.  Your data sets are given in the Table \ref{thermdata}. 
\begin{enumerate}
\item Start with data set number 1, and fit a line ($f_1(x)=ax)$ and a quadratic (of the form $f_2(x)=a x^2 + b x$ to the data set, and compute the associated \emph{coefficient of determination} (also called the \emph{correlation coefficient}) $r^2$ for both fitted curves. You can do this in most spreadsheet programs, or in a computer language.  Which curve fits the data better?
\item Now use the functions fitted to the previous data ($f_1$ and $f_2$) to fit the data in experiment 2.  Again, compute the value of $r^2$ for both fitted curves.  Which curve fits the data better?  If you were planning to generalize a model from experiment 1 to use as a fit to subsequent experiments, on the basis of what you know, which one would you choose -- $f_1$ or $f_2$?
\end{enumerate}
\end{minipage}
\begin{minipage}[r]{0.4\textwidth}
\centering
\hspace{5mm}\begin{tabular}{c|c|c}
 ~Time~    & $~(T-T_0)~$ & $~(T-T_0)~$  \\
~      & Exp.~1     & Exp.~2 \\
\hline 
0     & -0.110     & -0.110 \\
1     & 1.100      &0.832  \\
2     & 1.835     & 5.035  \\
3     &  2.347    & 6.547  \\
4     &  2.009    & 6.809  \\
5     &  5.771    & 10.771  \\
6     &  6.208    & 11.008  \\
7     &  11.198    &13.765  \\
8     & 12.010     & 15.211  \\
9     & 16.986     & 17.215  \\
10    & 20.420     & 19.200  \\
\end{tabular}
\label{thermdata}
\captionof{table}{Two experiments measuring temperatures.}
\end{minipage}

\item \textbf{Information.} The concept of \emph{information} was brought up in the text.  One way of thinking about the information content of model is to consider how difficult it is to describe.  To make this concrete, think about a base-2 system of digits.  Consider 5-digit integer numbers represented in base-2 format.  Each such number has the form "XXXXX" where each ``X" is either a ``1" or a ``0".  If it is a ``1", then one sums up a unit of $2^{(n-1)}$, where $n$ represents the numerical value of the position of the digit starting from the far right.  Thus, the binary digits ``10010" represent the decimal number $2^4+2^1=17$.  The information in such a set of digits is just the number of digits (which is close to the logarithm of the maximum number expressible in the system-- e.g. $\textrm{log}_2 11111 = 4.95$, or approximately 5 binary digits or \emph{bits}).  With this in mind, please answer the following questions.

    \begin{enumerate}
    \item How much information is there in the binary number 10010?  How about 11111? 
    \item Which contains more information-- the binary number 11111 or the decimal number 31?
    \item Suppose I send a message containing a single binary digit that can be either a 1 or a 0.  How much information in bits is gained by the person who receives my message?
    \item Consider the following two examples of a 5000 digit binary number.  Case (i) each digit is generated randomly with a coin toss and written down as a 1 (heads) or 0 (tails)  until 5000 digits are created.  Case (ii) a number consisting of 5000 repetitions of the numeral ``1".  (Hint: Can the second one be compressed some way that uses less than 5000 digits to communicate with no loss of digits?  Can the first one?)
    \item In the answer to your previous question, does the concept of \emph{regularity} of the two numbers involved make a difference in your treatment of them?  
    \end{enumerate}

\item As briefly mentioned in this chapter, some philosophers view language itself as a model of reality.  In a few sentences, describe what this might mean-- that language is itself a model.  Reference the discussion on the properties of models and modeling, and attempt to describe the process of \emph{abstraction}, \emph{model building}, \emph{realization}, and \emph{revision} as it might apply to a natural language by providing examples. 

\item Read the paper titled \emph{What is complexity science, really?} by S.E. Phelan (\bibentry{phelan2001complexity}).  Phelan identifies three philosophical frameworks in which \emph{science} might be defined and interpreted.  What are these?  Does the paper suggest that any particular framework is best?

\item Read the paper titled \emph{This is not a universe: Metaphor, language, and representation} by Liliane Papin (\bibentry{papin1992not}).  Papin discusses some of the challenges that language poses in science.  In your own words, what do you think is the main point of the paper?  In other words, what is the message that Papin is attempting to communicate regarding language and science?

\end{enumerate}
\abstract*{This is the abstract for chapter 00}

\begin{savequote}[0.55\linewidth]
``A bad review is like baking a cake with all the best ingredients and having someone sit on it."

\qauthor{Danielle Fernandes Dominique Schuelein-Steel, fourth bestselling fiction author of all time}
\end{savequote}
\def\CHAP{chapter01_review}
\theoremstyle{definition}

\chapter{Mathematical Definitions, Concepts, and Review}\label{chaprev}

Almost everyone who uses this text will have some background in mathematics; it is assumed that this background includes advanced algebra, introductory calculus, linear algebra, and an introductory course in ordinary differential equations.  The purpose of this chapter is to review a host of definitions, concepts about mathematics, and some of the basic results learned in previous coursework.  The material presented here is not meant to be exhaustive, but is focused primarily on elements that will be useful in the remainder of the text.  

Because the presentation of this chapter is specific to ideas that arise in the remainder of the text, the presentation is somewhat an agglomeration of important topics rather than an exposition focused on one or two main ideas.  In a very few instances, there is the introduction of material that may not have been covered in  undergraduate mathematics; for this material, however, the necessary background is only what is described in this text.  Many of the topics in this chapter may be safely skipped by those who do not need a reminder about mathematical definitions, concepts, or a review of basic of linear algebra and calculus.

\section{Sets and Set Builder Notation}

While will not make extensive use of the concept of \emph{sets}, we will use \emph{set builder notation} as a convenient and compact way about discussing intervals on the real line (or higher dimensional Euclidean [sometimes called Cartesian] 2-D and 3-D spaces).   The version that is adopted for this text will be simpler than is possible for more general settings.  For example, there are symbolic representations for the conjunctions ``and" and ``or" in set builder notation; we will opt for simply using the words themselves to keep the new symbols to a minimum. The following are the elements of set builder notation that we will use.  While some of these may be defined further below (e.g., the word set), here the goal is simply to describe the notation.

\begin{enumerate}
    \item \textbf{Set}.  A set is denoted using curly braces, ``$\{ ~\}$".  Thus $\{a, b , f\}$ is a set containing three items, $a$, $b$, and $f$.
    \item \textbf{``Such that" symbol}.  In set builder notation, the colon ``:" is used to indicate the concept of ``such that" or ``with the property that" or sometimes ``as follows".  When one sees the colon after a variable (say, the variable $x$), it might be useful to think of this as indicating ``create the set of $x$ values such that..."
    \item \textbf{``Is an element of" symbol}.  We use ``$\in$" as a substitute for ``is an element of".  We will define element more completely below.
     \item \textbf{``Is not an element of" symbol}.  We use ``$\not\in$" as a substitute for ``is not an element of". 
    \item \textbf{Predicate}. A \emph{predicate}.  In English grammar, is a part of a sentence that states something about the subject.  In set theory, the word predicate has same notion.  Thus, a predicate can be interpreted as a \emph{rule} or \emph{formula} that must be applied.  Suppose the predicate is given by the equation (or inequality) $\Phi$, and we are asked to construct the set of all values of $x$ such that $\Phi$ is met.  To be concrete, let's assume the correspondence $\Phi \Leftrightarrow x>1$.  We could write such a set, call it $A$, as follows
    \begin{equation}
        A = \{x : x>0\}
    \end{equation}
    This is read as follows: ``create the set $A$ where $A$ is equal to the set of all numbers $x$ such that $x>0$".  There is an implication here that we know what \emph{kind} of number $x$ should be.  Let's assume that $x$ is meant to be an integer.  If we wanted to include this explicitly in the set statement, we could rewrite it to read the following.
        \begin{equation}
        A = \{x : x>0 \textrm{~and~} x\in\ \mathbb{Z}\}
    \end{equation}
    Here, $\mathbb{Z}$ is the special symbol reserved to indicate the set of integers.
\end{enumerate}
Although we will not use them much, one can define the statement ``for all" by $
\forall$ an the statement ``there exists" by $\exists$.  For example, the even integers greater than zero could be described in set builder notation as 
\begin{equation}
{\{n: (\exists k) k\in \mathbb {Z} \textrm{~and~} k>0 \textrm{~and~} n=2k\}}
\end{equation}
\noindent This would be read as ``create set of all values $n$ as follows: there exists $k$ values that are integers and are greater than zero.  The set of $n$ values is found by taking two times the set of these $k$ values."  This concludes the number of ideas and symbols that we will adopt via set builder notation.

Now that we have the notation established, we can proceed to define some of the basic features of sets (we will be ultimately be primarily interested in the applications of set notation to intervals on the real line).

\begin{definition}[Element]\indexme{set!element}
An \emph{element} of a set (also called a \emph{member} of a set) is the generic name associated any one of the objects contained by the set.  Often the symbol $x$ is used to denote the element of a set (although any symbol could be used).
\end{definition}

\begin{definition}[Set]\indexme{set}
A \emph{set} is a collection of objects; more specifically for our purposes, a collection of mathematical objects. Sets can contain a finite number of objects, or an infinite number of objects.
\end{definition}
\begin{definition}[subset]\indexme{set!subset}
A \emph{set} where the elements of of the set are also elements of some given set.  To be a \emph{proper subset}, the subset must not contain all of the elements of the given set.
\end{definition}
\begin{definition}[universal set]\indexme{universal set}
Many problems are created when attempting to define a completely general \emph{universal set}; in the broadest sense, they can lead to paradoxes in formal set theories.  For our purposes (which is primarily to discuss intervals on the real line), we adopt a universal set that is a well-defined entity that avoids such problems.  Specifically, we take $U$ to be the real line in one dimension (or the 2-D plane or 3-D space in higher dimensions).  One may also \emph{define} $U$ to be some subset $I$ of the real line, as long as all discussions relate to intervals that are proper subsets of $I$.
\end{definition}
\noindent Even though our use of set theory will be minimal, it is still useful to define the basic operations of the \emph{union} and \emph{intersection} of sets.  These are as follows.

\begin{definition}[union of sets]\indexme{set!union}
The union of two sets $A$ and $B$ is the set of elements $C$ belonging to either of the two sets (without repetition).  This is frequently denoted by $=A\cup B$.  If the elements of the sets are represented by $x$, this can be written as $A\cup B=\{x:x \in A \textrm{~or~} x\in B\}$.  Extensions to three or more sets can be done by operating on two sets at a time iteratively. 
\end{definition}
\begin{definition}[intersection of sets]\indexme{set!union}
The intersection of two sets $A$ and $B$ is the set of elements $C$ belonging to both of the two sets.  This is frequently denoted by $=A\cap B$.  This can also be written as $A\cap B=\{x:x \in A \textrm{~and~} x\in B\}$
\end{definition}
\noindent Note that the concepts of \emph{union} of sets and \emph{intersection} of sets correspond to the English conjunctions \emph{and} and \emph{or}, respectively.  Being explicit, the statements $D=A \cup B \cup C$ with $x \in D$ implies that $x$ is in set $A$ or in set $B$ or in set $C$ (which could include any intersections among those sets).  In contrast, the statement $D=A \cap B \cap C$ with $x \in D$ implies that $x$ is in each of $A$, $B$, and $C$.

\begin{figure}[t]
\sidecaption[t]
\centering
\includegraphics[scale=.7]{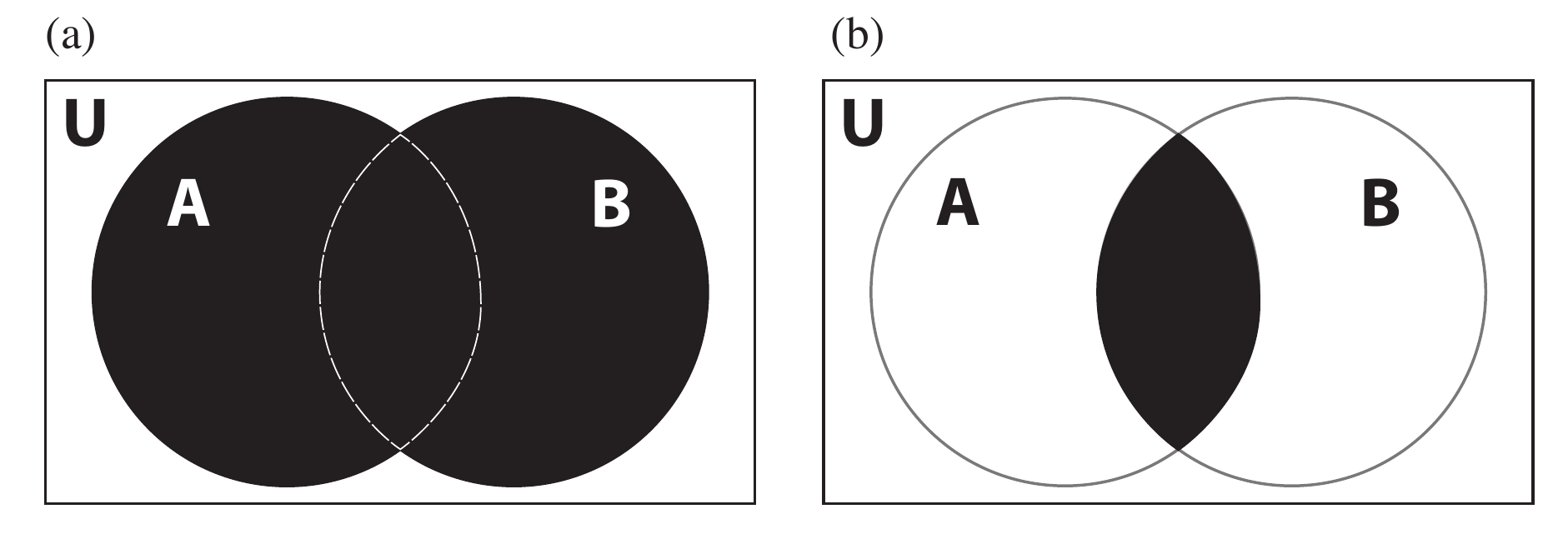}
\caption{(a) The union of two sets.  (b) The intersection of two sets. In these diagrams, ${A}$ and ${B}$ are sets, and these two sets are contained in the universal set $U$.  Thus, $A$ and $B$ are subsets of $U$}.
\label{sets}       
\end{figure}

\begin{definition}[compliment of a set]\indexme{set!compliment}
Suppose $C=A\cup B$. The compliment of $C$ in universal set $U$ consists of everything in $U$ that is not in set $C$.  This is often written as $C'=U \backslash C$ (where the backslash is read as subtraction or removal).  This can be alternatively written as $C'=\{x:x \in U \textrm{~and~} x \not\in A \textrm{~and~} x \not\in B\}$.  If $U$ defined is some subset of the real line (or higher dimensional spaces), then sometimes the compliment is called the \emph{compliment relative to} $U$.  
\end{definition}
\noindent These concepts are presented graphically in Fig.~\ref{sets}.  One additional definition is made here, in part because it helps to clarify a notation that is often observed in the mathematical literature.  

\begin{definition}[Cartesian product of sets]\indexme{set!Cartesian product}
The Cartesian product of two sets $A$ and $B$, is the set of all ordered pairs $(a, b)$ where $a \in A$ and $b \in B$.
\begin{equation}
A\times B=\{(a,b): a\in A\ {\mbox{ and }}\ b:\in B\}.
\end{equation}
\end{definition}
This latter notation will arise when we discuss the Euclidian plane and Euclidian space.

\subsection{Numbers}
The concept of numbers is a intuitive and natural thing that we have all been using since we were little and learning to count out objects on our fingers.  However, numbers themselves are an interesting topic; there is even an branch of mathematics called \emph{number theory} that studies the integers.  While the study of integers may not sound particularly exciting, it is actually an area of intense study!  For instance, number theory has been used to find very large prime numbers (numbers that are not the multiplication of two smaller numbers).  Such numbers are of great value in use in public key encryption schemes an other kinds of computer security.  In fact, there are several prize currently offered (of about $\$150$ K and as of 2021) to find a prime number with 10 or 100 million digits.  The largest prime number currently known (as of January 2022) is $2^{82,589,933}-1$, a number that has $24,862,048$ digits in base 10.

It turns out that the discussion of the various of \emph{sets of numbers} one might encounter in applied mathematics is a good way to discuss sets and intervals in a way that has familiarity.

The most frequently used sets of numbers are as follows.

\begin{figure}[t]
\sidecaption[t]
\centering
\includegraphics[scale=.2]{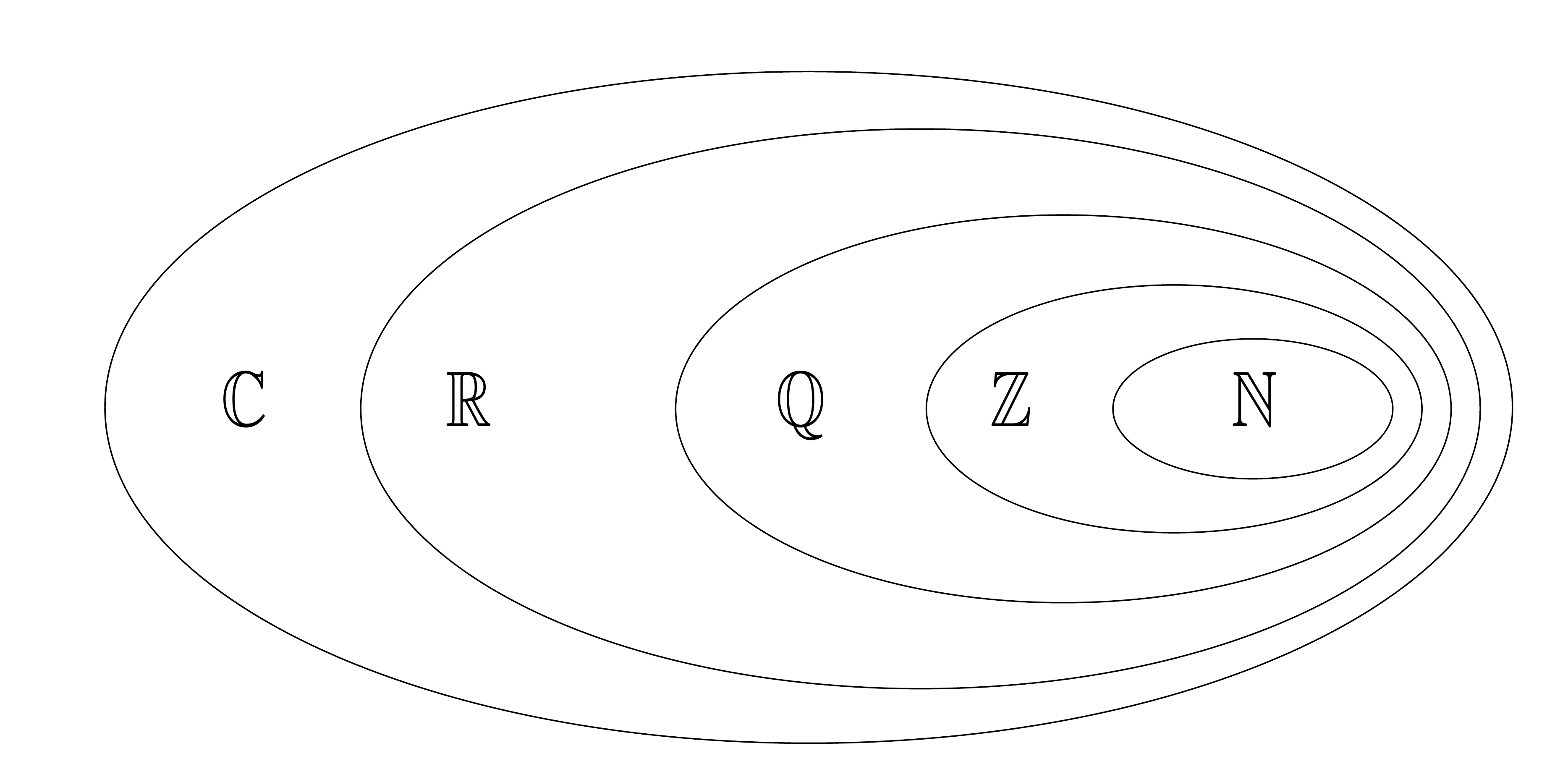}
\caption{The number systems, from the perspective of sets.  In this figure, each oval represents a set; an oval that is contains in another (larger) oval is a \emph{subset} of the containing set.}
\label{numbers}       
\end{figure}

\begin{itemize}
    \item $\mathbb{N}$ -- \textbf{The ``natural" numbers}.  These are the set of positive integers $\mathbb{N}=\{1,2,3,4, \ldots \}$.  Some definitions include zero; to denote that, the convention is to use the symbol $\mathbb{N}_0=\{0,1,2,3,\ldots\}$.\indexme{numbers!natural}\\
    \item $\mathbb{Z}$ -- \textbf{The integers}.  These are the set of all \emph{positive or negative} integers, including $0$ $\mathbb{Z}=\{\ldots -2, -1, 0, 1, 2, \ldots \}$.\indexme{numbers!integers}\\
    \item $\mathbb{Q}$ -- \textbf{The rational numbers}.  The ``$\mathbb{Q}$" here stands for \emph{quotient}, or, in other words, the ratio of two numbers.  The \emph{rational} numbers are formed by the set of all possible ratios of all integers $p$ and $q$ in $\mathbb{Z}$ so that all pairs for a ratio $p/q$ which is a member of $\mathbb{Q}$.  Note that there is a restriction; $q\ne0$.  Also, because we can have $q=1$, the integers are a subset of the rational numbers. \indexme{numbers!rational} \\
    \item $\mathbb{P}$ -- \textbf{The irrational numbers}.  In short, the set of all numbers that cannot be expressed as a rational number, i.e., they are not elements of $\mathbb{Q}$. \indexme{numbers!irrational} \\
    \item $\mathbb{R}$ -- \textbf{The real numbers}.  The real numbers are the set of all things that we might think of as conventional numbers; in other words, it contains as subsets all of the sets of numbers defined previously.  For our purposes, the real numbers can be though of most simply as the \emph{union} (or the combining of) the rational and the irrational numbers (in mathematics this might be written $\mathbb{R}=\mathbb{Q} \cup \mathbb{P}$). \indexme{numbers!real} The real numbers are what we use \emph{in principle} when doing computations relating to physical systems. The real line can also be though of as the 1-dimensional Euclidian space, where the Euclidian spaces are defined in the next section.\\
    
    \item $\mathbb{C}$ -- \textbf{The complex numbers}.  The complex numbers are a generalization of the concept of number that contain a \emph{real} and and \emph{imaginary} part.  Most of us have been introduced to the idea of complex numbers early in our mathematical career; however, they frequently instill much unease and consternation. As with many mathematical constructions, however, the ``reality" of complex numbers is not tremendously important.  What is important is that they meet certain algebraic necessities assuring mathematical consistency (they form a \emph{commutative} algebra), and, from a practical standpoint, they are useful.
\end{itemize}
A graphic indicating how these number systems (as sets) relate to one another is given in Fig.~\ref{numbers}.  Note that, for each of these number sets, the values $\pm \infty$ \emph{are not part of the set of numbers}!  There are number sets (known as the \emph{extended} reals) that contain $\infty$ as a value, but these sets will not be used in this text.  One extension to these sets are the sets of higher-dimension Euclidian spaces; this will be discussed below.  

These sets of numbers have interesting properties and histories of their own.  For example, the irrational numbers were first discovered by the Greeks.  In particular it is thought that a Pythagorean (i.e., a member of a Greek math cult following the tenants of Pythagoras), named Hippasus of Metapontum first proved this.  Strangely, although it is the Pythagorean theorem that allowed Hippasus to show that irrational numbers exist, this went in stark contradiction to the Pythagorean belief that all mathematics could be expressed through ratios of integers.  It is said that Hippasus discoverd this startling fact while at sea, and his fellow Pythagoreans were so upset by the revelation that they threw him overboard!  While the validity of this story is certainly suspect, what can be said is that revolutions in the understanding of numbers has frequently been met with substantial resistance.

While many think of the complex numbers (which will be reviewed below) as presenting challenges to concrete ways of thinking, even the real numbers are a bit stranger and philosophically more challenging than one might expect.  Some of the most contentious issues in defining modern mathematics has come from attempts to understand the set of real numbers, and many of these challenges arise via the subset of irrational numbers.  In particular the critical notions of \emph{limits} and \emph{completeness} of sets, topics that we will touch on (lightly) later, that took until the 19th and 20th centuries (see \citet{snow2003views} for more details) before the concepts were well understood.  Before moving onto other topics, the following example illustrates how the concept of real numbers is more complicated than it may appear on the surface.

\begin{svgraybox}
\begin{example}[$0.999\ldots  = 1$]
We have all be exposed the idea of repeating decimal numbers, but the characteristics of such numbers can be elusive.  For example, consider the sequence of numbers (we will define sequences more formally in \S \ref{seqseries}) that approach the value 1 as follows $s=\{0.9, 0.99, 0.999, 0.9999, \ldots\}$, or more generally $s=\{0.(9)_n, n\in \mathbb{N} \}$.  Clearly this sequence gets closer and closer to 1 as we add more repeats of the numeral 9 (i.e., increase $n$).  Also, if we fix any small number, $\epsilon$, then, no matter how small epsilon is, we can always take $n$ large enough so that we get closer than within $\epsilon$ to the value 1

\begin{equation}
    1-0.(9)_n < \epsilon
\end{equation}
In some sense, then, is $0.(9)_n$ equal to 1?  This seems like it would be a curious thing, but we can show that it is true. The following is an informal illustration.  There are proper, concrete proofs for the somewhat \emph{ad hoc} illustration below.  Nonetheless, the point is made, and the result is correct, as odd as the result may seem! \\

\begin{align}
    x&=0.9999\ldots \nonumber \\
    10x&=9.999\ldots \nonumber \\
     10x&=9+0.9999\ldots \nonumber \\
     10x&=9+x \nonumber \\
     9x &= 9  \nonumber \\
     x&=1
\end{align}

\noindent The explanation for this involves the idea of convergence of a sequence and completeness of a set to fully understand.  For now, we will put this notion aside and simply recognize that real numbers (particular those that are irrational or are repeating decimals) are trickier to understand than rational numbers are!  There is even some more philosophical discussion about in what sense real numbers are (objectively) real.  Because of the various concepts of infinity associated with the real numbers, they sometimes illustrate behavior that is counter intuitive (the fact that $0\times\infty$ is not defined is one such example; it is easy to construct examples where $0\times\infty$ is equal to any number that you choose).
\end{example}
\end{svgraybox}

\subsection{Euclidian Space}

Euclidan space might best be described, somewhat colloquially, as the 3-dimensional space that we are all familiar with.  Before continuing, in order to understand Euclid's definitions, we must first understand what an axiom is.

\begin{definition}[Axiom]\indexme{axiom}\indexme{axiom}
An \emph{axiom} is a statement about a physical or mathematical system that is assumed to be true, but cannot be \emph{proven} to be true within the system itself.  Sometimes axioms are called \emph{laws} or \emph{principles} or \emph{postulate}; these alternatives names are used historically or because of disciplinary differences in terminology (e.g., between physics and mathematics).
\end{definition}

Now we are in a position to define Euclidian space.   Euclidian space is a space that we are used to thinking about in, say, ordinary geometry.\indexme{geometry!Euclidian}  Euclidian spaces subscribe to the five axioms of Euclid's \emph{Elements}.  These are as follows.

\begin{enumerate}
\item  A straight line may be drawn between any two points.
\item  Any terminated straight line may be extended indefinitely.
\item  A circle may be drawn with any given point as center and any given radius.
\item All right angles are equal.
\item For any given point not on a given line, there is exactly one line through the point that does not meet the given line.
\end{enumerate}
where the language here is not identical to that of Euclid, but the concepts are (in particular, this version of axiom 5 is known as Playfair's axiom \citep[][p.~7]{playfair1795elements}).  So, why invoke Euclid at this juncture?   It turns out that most of the mathematics you have done have been done in Euclidan space.  In one dimension, Euclidian space is the real line; in two dimensions, the plane; and, in three dimensions what we often use as our model of space.  Using set notation, in 1-, 2-, and 3-dimensions, we have

\begin{align}
\textrm{1-dimension} && \mathbb{R}& = \mathbb{R}^1=\{x: x\in \mathbb{R}\} \\
\textrm{2-dimensions} && \mathbb{R}\times \mathbb{R}& = \mathbb{R}^2=\{(x,y): x\in \mathbb{R} \textrm{~and~}y\in \mathbb{R}\} \\
\textrm{3-dimensions} && \mathbb{R}\times \mathbb{R}\times \mathbb{R}& = \mathbb{R}^3=\{(x,y,z): x\in \mathbb{R} \textrm{~and~}y\in \mathbb{R} \textrm{~and~}x\in \mathbb{R}\} 
\end{align}
\noindent Here, the Cartesian product defined above has been used in defining $\mathbb{R}^2$ and $\mathbb{R}^3$.

Sometimes Euclidian geometry is also called \emph{flat} space; this is to say that it has no curvature to it.  Physicists (in particular A. Einstein) have shown that space is not, in fact, flat, but rather curved because of the effects of gravity.   Regardless, for most systems that are terrestrial (i.e., small enough) and not subject to need for extremely accurate measurements, we may think of space as being Euclidian.  What these spaces have in common is that our conventional notions of geometry (parallel lines do not intersect; all right angles are equal in angular measure) remain true in this model.  

\begin{figure}[t]
\centering
\includegraphics[scale=.2]{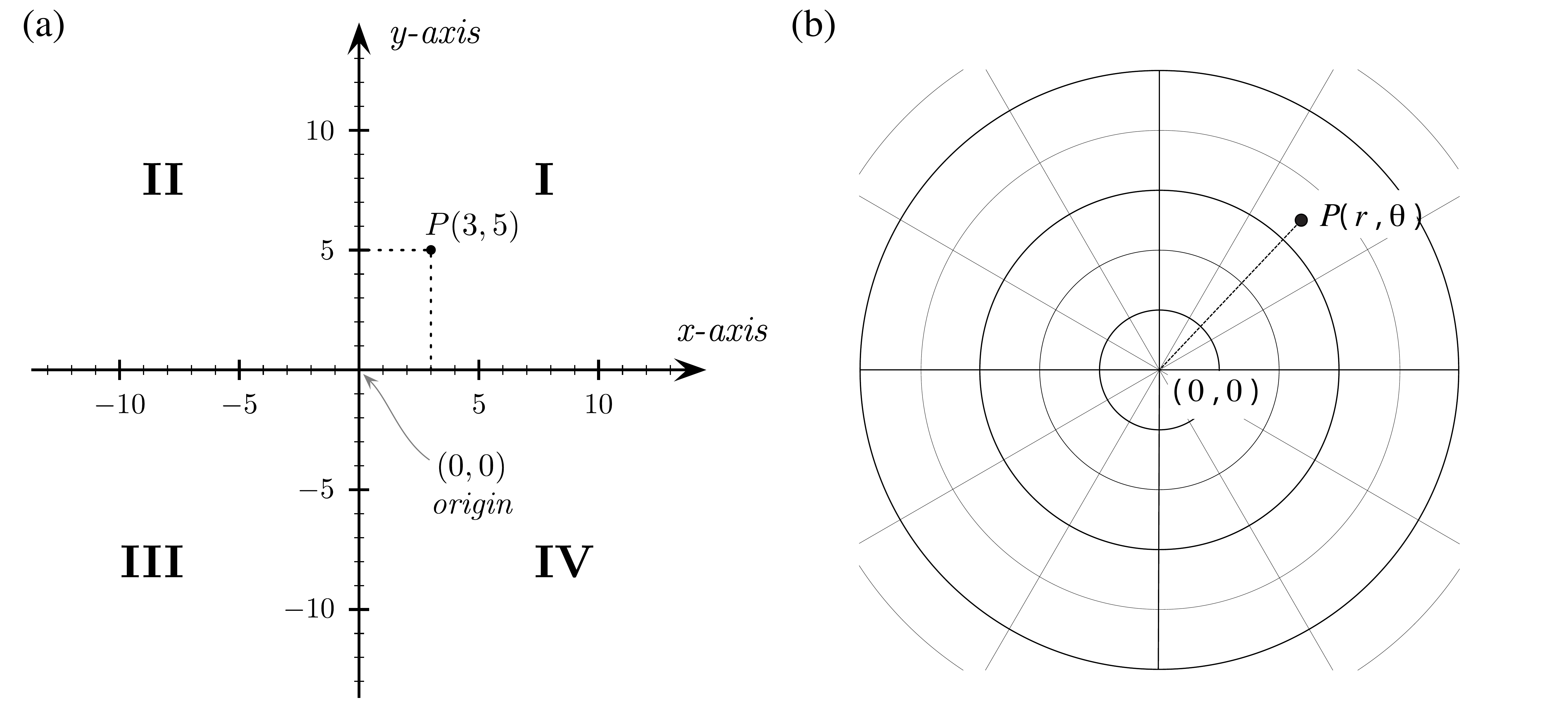}
\caption{Two coordinate systems for a 2-dimensional (2-D) Euclidian space. (a) The Cartesian coordinate system (CC BY-SA 3.0, \protect\url{https://commons.wikimedia.org/w/index.php?curid=603366}.  (b) The cylindrical coordinate system.  Regardless of the \emph{coordinate system} used to \emph{label} the points in Euclidian space, the space itself remains unchanged.  Thus, the five axioms of Euclidan geometry apply to both coordinate systems.  }
\label{fig:cartesian}       
\end{figure}

It should be noted that coordinate systems are \emph{independent} of the underlying space that they describe.  The conventional coordinate system that we use to describe Euclidian space is the Cartesian one, named after the French philosopher and mathematician Ren\'e Descartes, who first described it.  In fact, sometimes Euclidian space is called \emph{Cartesian space}.  However, some care is needed.  A Cartesian \emph{coordinate system} is the rectangular coordinate system that we are all familiar with (Fig~\ref{fig:cartesian}(a)).  The coordinate system may be a very convenient one, but it should not be conflated with Euclidian space itself.  Coordinate system are a convenience in which each point may be given a unique label; but the geometry of the coordinate system does not necessarily reflect the geometry of the underlying space.  For example, one may also adopt the familiar cylindrical coordinate system (Fig.~\ref{fig:cartesian}(b)) to describe Euclidian space.  This can be convenient, when, for example, one is modeling an object in Euclidian space that is itself cylindrical.  While the two coordinate systems are different, the underlying space (and the geometric principles that are assumed to be true on it) remain unchanged.  

While not an important concept for the material in this text, it is also possible that the underlying space itself is not \emph{flat} (i.e., it is not Euclidian).  For example, a spherical space is one that you may have heard of by analogy with the globe; in such a space, the five axioms of Euclidian geometry must be amended to create a consistent system.   The discovery of non-Euclidian spaces (e.g., such as the spherical space just mentioned, or the hyperbolic spaces) was made in the early 1800's by mathematicians in Europe and Russia.  The discovery of such unique geometries eventually diffused into popular culture, and represented a true revolution in the way that mathematicians, scientists, and the public at large viewed science and scientific discoveries.  As an example of the influence in popular culture, the ideas of non-Euclidian geometries were so widely known about that the author H. P. Lovecraft adapted the unique geometries to describe otherworldly settings in his fictional writing.  As with the development of set theory (discussed very briefly at the end of the next section), the discovery of non-Euclidian geometry helped spark the revolution in mathematics that started in the early 19th century and lasting  through the middle of the 20th century.

%

\subsection{Intervals on the Real Line}\label{real_line_subsec}

The notion of an \emph{interval} can be described in the language of \emph{sets}, and the concept of intervals is one that we will used frequently. Below the concepts of open, closed, finite, and unbounded intervals are discussed.  Analogous concepts apply for higher dimensions, but no discussion of those extensions is attempted at this juncture.  Similarly, set notation is discussed, but set theory (other than the basic notion that sets exist and have elements; and we cover the basic operations of unions and intersections). It is not an overstatement to say that set theory proper forms the basic underpinning of modern mathematics.  It is has also been one of the most difficult and contentious parts of mathematics, and contains controversies (or differences in approach and opinion) that are continue to generate discussion and research.  While we will not go into \emph{set theory} proper, or the details of why it has been such a challenging component of mathematics, we will briefly describe what gave rise to some of the difficulties at the end of this section.

\begin{definition}[Interval]\indexme{sets!interval}
An \emph{interval} is a set containing all of the real numbers between two specified real numbers $a$ and $b$; assume that $a<b$.  Intervals may be \emph{open}, \emph{closed}, or a mixture of open and closed.  A closed interval includes the end points, e.g.,
\begin{equation*}
     I=\{x: a \le x \le b \} \textrm{ usually denoted } I=[a,b]
\end{equation*}
whereas an open interval does not

\begin{equation*}
     I=\{x: a < x < b \} \textrm{ usually denoted } I=(a,b)
\end{equation*}

\indexme{interval}\indexme{interval!open}\indexme{interval!closed}
\end{definition}
One curiosity about open intervals such as $(a,b)$ is that they contain a maximum or a minimum real number.  By definition, the open interval does not include the end points $a$ and $b$.  Thus, for example, the open interval (1,2) contains all of the numbers greater than one, and all of the numbers less than 2, but not 1 or 2.  Another way to describe an open interval, then, is as follows.  For a moment, let's consider only the minimum part of the interval specified by $a$.  Now, suppose that we select a number in the interval as close to $a$ as we like; let's call this number $c$.  By definition, $c-a>0$.  We could then, say, select an even smaller number half way between $a$ and $c$, and call this $d=(c+a)/2$.  This number is also in the interval because $d-a=\tfrac{1}{2}(c-a)>0$.  This argument can be repeated indefinitely, and thus there is no smallest number in the open interval.  A similar argument can be made for the maximum number in the open interval.  

Intervals can also be half-open (or half-closed, which means the same thing) in the obvious way, e.g., $I=\{ x: a < x \le b\}$ is half open because the lower bound is not included.  \indexme{interval!half-open}. 

Finally, note that intervals can be \emph{unbounded}.  An unbounded interval is defined as follows.

\begin{definition}[Unbounded interval]\indexme{sets!unbounded interval}\indexme{interval!unbounded}
An \emph{unbounded interval} is a set containing all of the real numbers greater than or less than some specified real number, $a$.  Such intervals can be open or closed; thus $I=[a,\infty)$ is considered unbounded and closed, whereas $I=(a,\infty)$ is unbounded and open.  
\end{definition}

Accounting for various possibilities defined above for intervals (and defining the \emph{empty} and \emph{degenerate} intervals), then
intervals of the real numbers line can be classified into eleven different types (cf. \citet[][Chp.~3]{craig1969modern}) listed below. 
\begin{align*}
\intertext{\bf Empty:}
  &&  [b,a]&=(b,a)=[b,a)=(b,a]=(a,a)=[a,a)=(a,a]=\{\}=\varnothing \\
\intertext{\bf Degenerate:}
 && [a,a]&=\{a\}
\intertext{\bf Proper and bounded:}
        \textrm{Open:}&& (a,b)&=\{x: a<x<b\} \\
        \textrm{Closed:}&& [a,b]&=\{x: a\leq x\leq b\} \\
        \textrm{Left-closed, right-open:}&& [a,b)&=\{x: a\leq x<b\} \\
        \textrm{Left-open, right-closed:} && (a,b]&=\{x: a<x\leq b\} \\
\intertext{\bf Left-bounded and right-unbounded:}
        \textrm{Left-open:}&& (a,+\infty )&=\{x: x>a\} \\
        \textrm{Left-closed:}&& [a,+\infty )&=\{x: x\geq a\} \\
        \displaybreak[1]
\intertext{\bf Left-unbounded and right-bounded:}
        \textrm{Right-open:}&&  (-\infty ,b)&=\{x: x<b\} \\
        \textrm{Right-closed:} && (-\infty ,b]&=\{x: x\leq b\} \\
         \displaybreak[2]
\intertext{\bf Unbounded at both ends (both open and closed):} 
&& (-\infty ,+\infty )&=\mathbb{R}
\end{align*}
This last interval (the real line) is listed as both open and closed.  This is a technical detail that creates no end of discussion on mathematics forums on the internet.  The explanation is not easy without invoking additional mathematical structure (such as \emph{topology}).  In short, however, it might be described something like this:  The real line contains all possible real numbers; thus is it closed with respect to the real numbers.  However, there is no maximum nor minimum real number on the real line, thus (as for the discussion of the open interval) the real line is also open.  While this sounds like a paradox, it is also a statement that binary options (``true" or ``false") in mathematics are not always necessary; sometimes there is a \emph{third} option which is ``neither true nor false".  There is a saying in mathematics that ``A door must be either open or closed, and cannot be both, while a set can be open, or closed, or both, or neither!” \citep[][p.~91]{munkrestopology}.  While this may sound a bit noncommittal, the reality is that these concepts are both well defined and useful.  

\textbf{A short note about set theory}.  We have used set notation above, but not invoked set theory per se.  You may have learned a little about set theory even in grade school -- the ideas at that level are usually to discuss how one might group objects together (creating a set), and then compare the properties of various sets.  On the surface, this seems deceptively simple.  However, the goal of developing a robust set theory led to one of the most interesting periods of mathematical research in the history of mathematics.  By \emph{robust} here, something specific is meant.  This means that the theory must propose sufficient axioms such that it is both \emph{complete} and \emph{consistent}.  The word complete in this context means essentially that ``every true statement in the system can be proven from the axioms".  The word consistent means essentially ``every proper statement (or question if you prefer) within the system can be shown to be \emph{either} true \emph{or} false."  The hope was that because set theory was so fundamental, then, in principle at least, a powerful enough version of set theory could be the basis for \emph{deriving} all of mathematics (at least, in principle).  Quite unexpectedly, such a set theory was ultimately proven to be impossible.

Early work on set theory was done by many mathematicians; but in particular the work by two German mathematicians -- Richard Dedekind and George Cantor -- from about 1870 to 1900 paved the way for what is sometimes referred to as \emph{Naive set theory}.  A good and short introduction to this history can be found in \citep{ferreiros_encyclopedia}.  The modifier ``naive" here is not a slight; it means only that the theory was not cast in the language of  formal logic. While this set theory was incredibly useful, Cantor was alarmed to find that it was not consistent -- in other words, he found that paradoxes could be constructed in the theory.  Without recounting the entire history of set theory, the inability to generate a paradox-free set theory led to a crisis in mathematics in the early 1900's.  This crisis, while long in the making, was made most apparent by Bertrand Russel who described the paradox in 1901.  The Russel paradox is a mathematical statement much like the so-called \emph{liar's paradox}, which make the conflicting and self-referential statement ``This sentence is not true".

This kind of paradox would have deep and lasting ramifications for mathematics.  While set theory was ultimately repaired to some degree by the mathematicians Ernst Zermelo and Abraham Fraenkel in the 1920s (known in mathematics as the ZFC set theory).  The statement ``to some degree" must be qualified here.  The ZFC avoids the kinds of paradoxes that plagued earlier versions of set theory.  However, in the 1930s a mathematician named Kurt God\"el proved a rather unexpected result.  G\"odel was able to show, roughly, the following: Any axiomatic system which is complex enough to describe ordinary integer arithmetic is either \emph{incomplete} or \emph{inconsistent}.  In other words, if our axioms are consistent (i.e., do not lead to paradoxical or undecidable statements), then in every model of the axioms there is a statement which is \emph{true but not provable}.  While this result was somewhat of a blow to the idea of generating an overarching theory for mathematics, the actual impact is still debated.  First of all, the proof used by G\"odel was of a particular flavor known as \emph{first order logic}.  Thus, it is unclear if more capable logic languages suffer the same fate.  Secondly, while the premier set theory, ZFC, must by G\"odel's theorems either be inconsistent or have true but unprovable statements, it is not clear that this has any \emph{practical} import.  For example, if I state ``all non-pterodactyls are not non-dinosaurs", I have certainly uttered something true.  Regardless, if I am not interested in anything about dinosaurs (flying or not), then this true sentence is of no value to me.  So, while it might be impossible to prove all \emph{true} statements in ZFC, so far the such unreachable truths have not created any practical problems (possibly because, like the example, they convey truths that are not relevant).  And, to date, nobody has found an inconsistent statement within the context of ZFC.  Thus, ZFC has been a useful tool in mathematics.  The search for more powerful set theories continues as a topic of interest in mathematical research.

\subsection{Complex Numbers}\indexme{numbers!complex}

Unlike, say, the natural numbers (which can be illustrated by collecting actual objects), the complex numbers are a purely mathematical construct.  That does not mean that they are not useful or interesting, however.  There are many examples of concepts that exist only mathematically, but are nonetheless useful for many practical applications. 

The complex numbers are an extension of the real numbers.  Although almost everyone reading this text has probably encountered them previously, it is useful to recap their basic properties.

\begin{definition}
A \emph{complex number} assumes that there exists a mathematical object, called the unit imaginary number, $i$, such that $i=\sqrt{-1}$ so that $i^2=-1$.  Every complex number consists of two parts, a \emph{real} part, and an \emph{imaginary} part that is proportional to $i$.  The conventional form for a complex number is $a+b i$, where $a$ and $b$ are real numbers.  The set of all complex numbers is usually denoted by $\mathbb{C}$.
\end{definition}

The rules for addition, subtraction, and multiplication are slightly modified from those of real numbers as follows.  First, addition and subtraction are done by adding and/or subtracting the real and imaginary parts of a complex number independently.  Therefore

\[ (a+ i b) + (c+ id)= (a+c)+i (b+d)\]
\[ (a+ i b) - (c+ i d)= (a-c)+i (b-d)\]
Multiplication of two complex numbers is defined as follows
\[ (a+i b)\times (c+i d )= ac +i b c  +i a d i - b d  \]
or, equivalently, grouping terms
\[ (a+i b )\times (c+i d )= (ac-b d) +i (b c  +a d)   \]

The complex numbers are often represented as vectors on a plane, where the real part is plotted on the horizontal axis, and the complex part on the vertical axis.  In this representation, a complex number would be represented by a pair of points, i.e., $x+i y   \Leftrightarrow  (x,i y)$.  There are some advantages and disadvantages of this, but it generally improves interpretation for problems with physical significance.  When this formalism is adopted, a complex number is frequently represented typographically as a vector, as in ${\bf z}=(x, i y)$.  The representation of the complex numbers in this manner requires mildly re-defining the vector dot product.  We define this in the following.

\begin{definition}[The Complex Conjugate]
For a complex number ${\bf z}=a+i b =(a,i b)$, the \emph{complex conjugate} is defined by $\Bar{\bf z}=a- i b= (a,-i b)$.
\end{definition}

\begin{definition}[The Complex Dot Product]
For two vectors ${\bf z}_1=(a_1, i b_1)$ and ${\bf z}_2=(a_2, i b_2)$, the complex dot product between ${\bf z}_1$ and ${\bf z}_2$ is given by 
\begin{equation*}
    {\bf z}_1 \cdot \Bar{\bf z}_2 = (a_1 a_2+b_1 b_2)
\end{equation*}
\end{definition}
\noindent which has the advantage of being a single real number.  In particular, this means that a vector dotted with itself is defined by 

\[  {\bf z}\cdot \Bar{\bf z}= (a,i b)\cdot (a,-i b) = a^2+b^2 \] 
and the magnitude of a complex vector is defined by a rule that looks much like the rule for computing the length of a vector on a plane

\[ \left\| {\bf z} \right\|=\sqrt{{\bf z}\cdot\Bar{\bf z}}=\sqrt{a^2+b^2} \]

\begin{figure}[t]
\sidecaption[t]
\centering
\includegraphics[scale=.4]{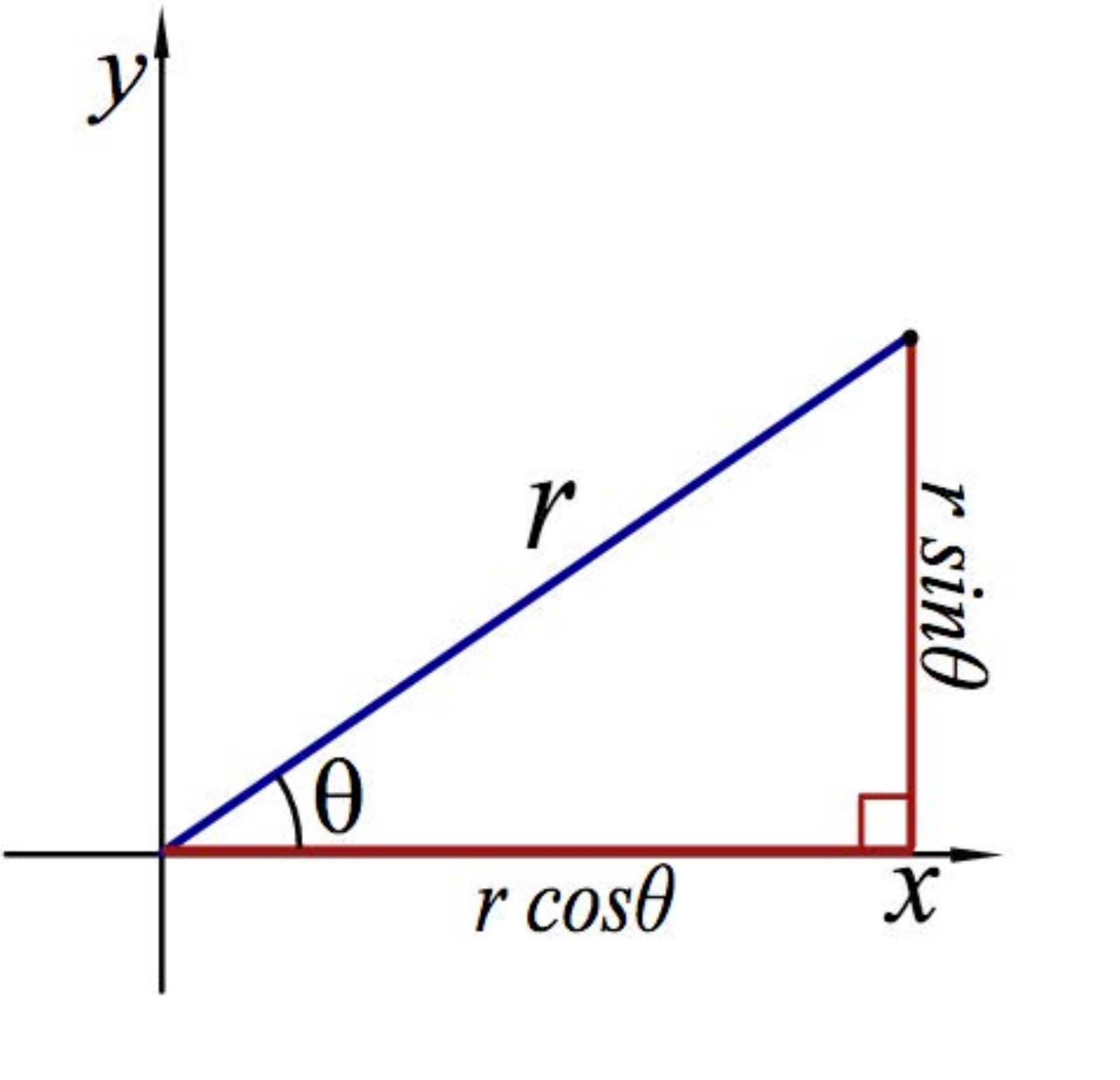}
\caption{The complex plane.  Complex numbers can be given a meaningful interpretation as vectors on a plane.}
\label{polar}       
\end{figure}

The recognition that complex numbers can be treated as vectors on a plane implies that they can be conveniently represented in polar coordinates.  In fact, there are some significant reasons for doing so.  To start, define the length  of the vector by the number $r= \left\| {\bf z} \right\|=\sqrt{a^2+b^2}$.  Then, we have the following relationships between the $(a,b)$ and $(r,\theta)$ coordinate systems (Fig.~\ref{polar})

\begin{align}
r&=\sqrt{a^2+b^2}&    a&=r \cos(\theta)& i b&= i r \sin(\theta) \\
\theta& = \tan^{-1}\left( \frac{b}{a}\right)
\end{align}
or,    
\begin{equation}
        {\bf z}= r(\cos(\theta),i \sin(\theta))=r(\cos(\theta)+ i \sin(\theta))
        \label{polareuler}
\end{equation}
In this representation, $r$ is sometimes called the \emph{magnitude} or \emph{modulus}\indexme{complex number!modulus}\indexme{modulus} or \emph{absolute value} \indexme{complex number! magnitude} of the complex number, and $\theta$ is called the \emph{argument} \indexme{argument}\indexme{complex number!argument} of the complex number.

You may have seen \emph{Euler's formula} before; it is the formula that leads to the famous relationship $e^{i \pi} = -1$ (where, recall, $e$ is the base for the natural logarithm).  More generally, however, Euler's formula takes the form

\begin{equation}
    e^{i \theta }= \cos \theta +i \sin \theta
\end{equation}
\noindent Although this formula looks truly remarkable, once we accept that the imaginary numbers are an acceptable extension of the reals, the proof of this result becomes quite simple (we will examine that further when we tackle infinite series).  For now, we adopt the formula without proof.  However, note that it allows a particularly simple representation of a complex number ${\bf z}$ as given by Eq.~(\eqref{polareuler}) in the form

\begin{equation}
        {\bf z}= r e^{i\theta}
        \label{polarform}
\end{equation}
As an interesting side note, this last formula allows one to make sense of the logarithm of a complex number.  Taking the natural logarithm of both sides of Eq.~\ref{polarform}, we have

\begin{align}
    ln({\bf z}) &=ln(r e^{i\theta})\\ 
    &=ln(r)+i\theta
\end{align}
Noting that $e^{i \theta} =e^{i(\theta +2k\pi)}, ~k=0,1,2,\ldots$, then it is clear that there are generally an infinite number of representations for the logarithm of a complex number

\begin{equation}
    ln({\bf z}) =ln(r)+i(\theta+2k\pi),~k=0,1,2,\ldots
\end{equation}

Although complex numbers have a helpful representation as vector quantities on the complex plan, they are technically just an extension of the real number system.  Therefore, it is not common to adopt a bold-face type to represent them (as we have done above).  In general, complex numbers are set in regular, italicized script, (e.g., $z$).  Generally, the context prevents there from being any confusion.  In future uses of complex numbers, we will not use bold-faced script to represent them.  Thus, the equation above would be more properly written

\begin{equation}
         z = r e^{i\theta}
\end{equation}

\section{Functions}
\label{sec:functions}

\begin{figure}[t]
\sidecaption[t]
\centering
\includegraphics[scale=.5]{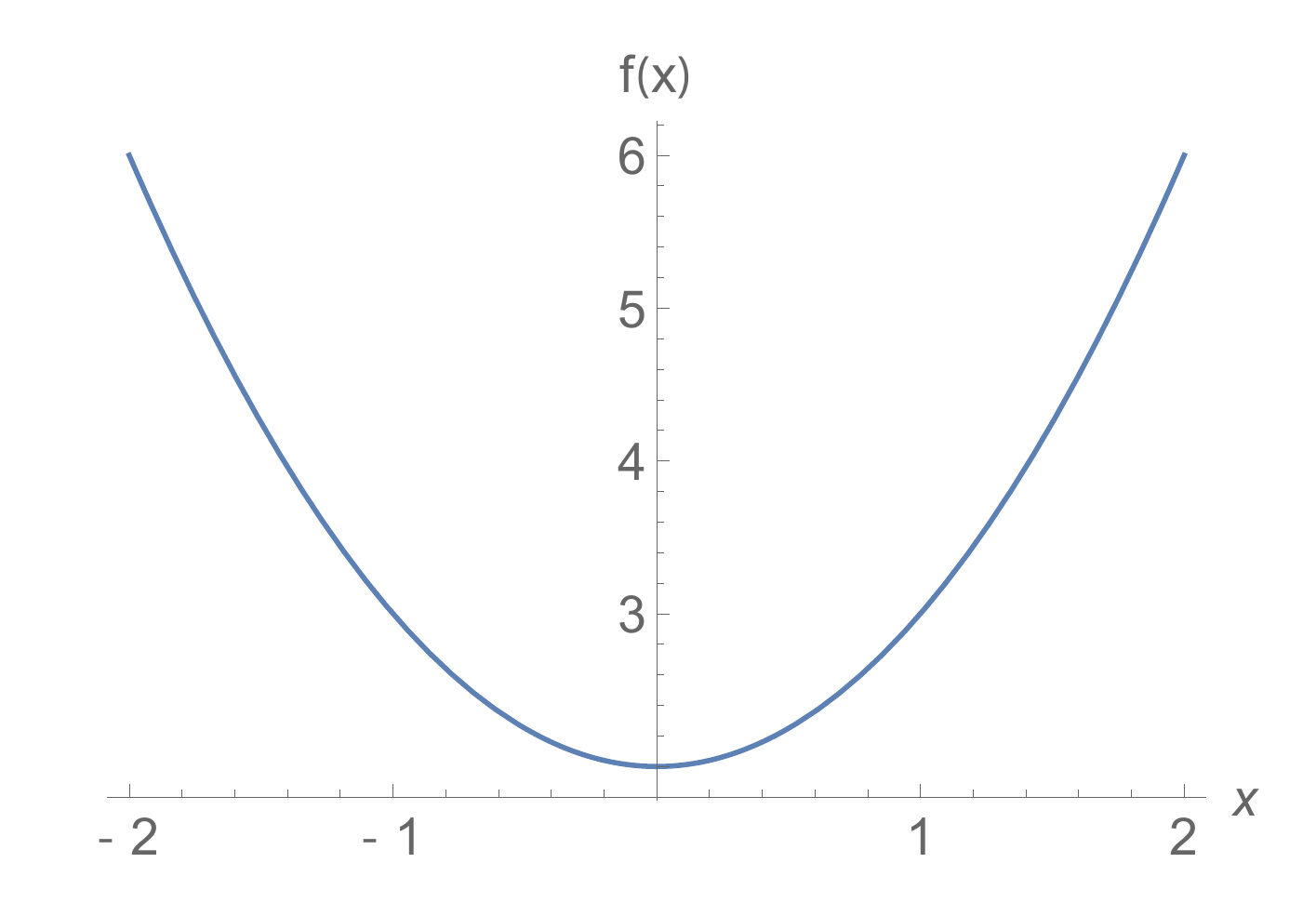}
\caption{A function.}
\label{fig:1}       
\end{figure}

We all have been introduced to the concept of functions.  When most people in college mathematics think of a function, the first thing that comes to mind is a relationship that looks something like 

\begin{equation}
    f(x) = x^2 +2,\qquad x\in[-2,2]
    \label{function}
\end{equation}
In this description, we are given a \emph{domain} (the numbers such that $-2 \le x \le 2$- recall, square brackets indicate \emph{inclusive} interval notation), and a \emph{range} for the function (the values $f(x)$ for the defined domain).
A graph of this function is plotted in Fig.~\ref{fig:1}.  Although this is how we typically think of a function, we can define it more generally.

\begin{definition}
A \emph{function} (or sometimes \emph{mapping}) is a relationship between two sets, $A$ and $B$, such that each element of $A$ is \emph{uniquely} associated with an element of $B$.
\end{definition}
\noindent The word \emph{uniquely} appearing in our definition here is really important: a function is always single-valued. 

\begin{figure}[b]
\sidecaption[t]
\centering
\includegraphics[scale=.7]{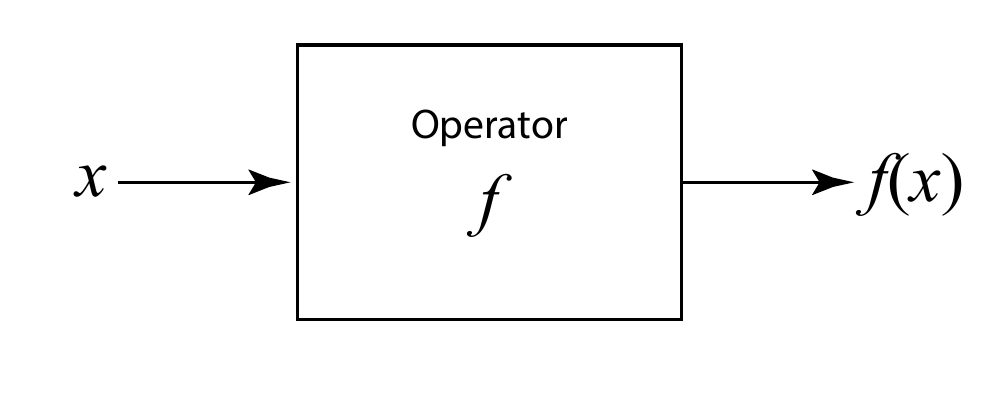}
\caption{An function, in conceptual form.  A function can also be considered a type of \emph{operator}.}
\label{fig:op}       
\end{figure}

This definition requires a bit of additional explanation for our purposes.  First of all, for us, the sets involved for a function are almost always intervals of the real number line. We even give these functions additional names and descriptions, as follows.  We usually refer to functions as having \emph{independent} and \emph{dependent} variables.  The independent variables are the \emph{inputs} to the function (or more properly, the set that is input to the function; this is the set $A$ in the definition).   Conversely, the \emph{dependent} variables are the set of values that are produced by the set of independent variables upon application of the function (this is the set $B$ in the definition). Each member of the set of independent variables is mapped to \emph{exactly one} member of the set of independent variables (this is the \emph{uniquely} part of the definition).  Fig.~\ref{fig:op} gives a pictorial representation of the process.

The concept of functions that strictly increase or decrease is one that is used routinely in mathematics; so much so that they have a special name: \emph{monotonic}.  These functions are defined as follows.

\begin{definition}
A function $f(x)$ is said to be \emph{monotonically increasing} on and interval, $I$, if the following are true:
\begin{enumerate}
    \item $f(x_1) \le f(x_2)$ for all $x_1 < x_2$ in $I$.
    \item $f(x)$ is not the constant function ($f(x)=C$, $C$ a constant). 
\end{enumerate}
Similarly, a function $f(x)$ is said to be \emph{monotonically decreasing} if the following are true:
\begin{enumerate}
    \item $f(x_1) \ge f(x_2)$ for all $x_1 < x_2$ in $I$.
    \item $f(x)$ is not the constant function. 
\end{enumerate}
\end{definition}
Obviously, functions in general are neither monotonically increasing nor decreasing, but a combination of these two concepts over subdomains of their total domain.  The constant function is neither monotonically increasing nor monotonically decreasing on any subdomain of its domain.

There is one final class of functions that requires definition because the term arises so frequently in applied mathematics.  These are the \emph{algebraically homogeneous}, or simply \emph{homogeneous}, functions.  Note that later on, the concept of homogeneous differential equations will be discussed; while the concepts are related, one should not confuse the two.  Here, we define algebraically homogeneous functions as follows:

\begin{definition}[(Algebraically) homogeneous function]
An algebraically homogeneous function is one with multiplicative scaling behaviour: if all its arguments are multiplied by some factor, $\alpha$, then the effect on the function is that its value is multiplied by some power of the factor $\alpha$. In other words, homogeneous functions display the following behavior
\begin{equation*}
    f(\alpha x) = \alpha^n f(x)
\end{equation*}
Here, $n$ is a real number, and it called the degree of the function.
\end{definition}

\begin{svgraybox}
\begin{example}[Homogenous functions.]
Below are a few examples of homogeneous functions.
\begin{enumerate}
    \item All polynomials of the form $f(x)=x^n$ are homogeneous.  To see this, just note the following:  $f(\alpha x)=(\alpha x)^n=\alpha^n x^2$.  In this case, the degree of the function is equal to the degree of the polynomial.
    
    \item No functions with additive constants are homogeneous.  This is also easy to demonstrate. Suppose $f(\alpha x)=\alpha^n f(x)$ is a homogeneous function.  Now consider the function $g(x)=f(x)+c$.  The function $g(x)$ cannot be homogeneous, because a multiplicative scalar $\alpha$ will not scale the constant $c$.  Functions like $g(x)$ (i.e., ones that contain an additive constant) are called \emph{affine}.  This just means that the two functions are connected by a linear scaling of the coordinate (the independent variable) followed by a translation (the additive scalar).
    
    \item  The function $f(x)=1/x^2$ is homogeneous.  To see this, note
    \begin{align*}
        f(\alpha x) &= \frac{1}{(\alpha x)^2}\\
        &= \frac{1}{\alpha^2}\frac{1}{x^2}\\
        &=\alpha^{-2} x^2
    \end{align*}
    So, this function is homogeneous, with degree equal to $-2$.
\end{enumerate}
\end{example}
\end{svgraybox}

\subsection{Boundedness}\label{boundedness}

In science and engineering, we usually think of the sets involved as being intervals of the real number line.  To put this in context, think of our example given by Eq.~\eqref{function} above.  Here, we can define of the \emph{input} set $A$ in set builder notation as $A=\{x: x\in[-2,2]\}$ (as a refresher, this is read as follows: ``The set A is constructed by collecting the all of the numbers $x$ such that $x$ is between -2 and 2, inclusive of the endpoints").  The set $B$ we can think of as all of the numbers defined by the function $f$.  Recall, we have $f(x)=x^2+2$.  The output set, $B$, is defined by $B=\{f: f(x)=x^2+2\textrm{~for~all~} x \in A\}$.  Note also that these intervals are either $closed$ or $open$ intervals.  The concepts of $closed$ and $open$ intervals (or, generally, sets) is actually much deeper and complex than one would assume (as described in \S\ref{real_line_subsec}.  Now, we combine the ideas of \emph{intervals} of independent variables, and \emph{functions} on those intervals.  The following is an example.

\begin{svgraybox}
\begin{example}[Domains of a Function as an Interval]
The function
\begin{equation}
f(x) =\frac{ \sin(x)}{x}, ~~x\in(0,\pi]
\end{equation}
has the domain $D=(0,\pi]$.  The function is monotonically decreasing on this domain.  While the function may not seem particularly unusual, it actually does have some strange behavior at $x=0$.  A plot of the function is given in the figure below.

{
\centering\includegraphics[scale=.45]{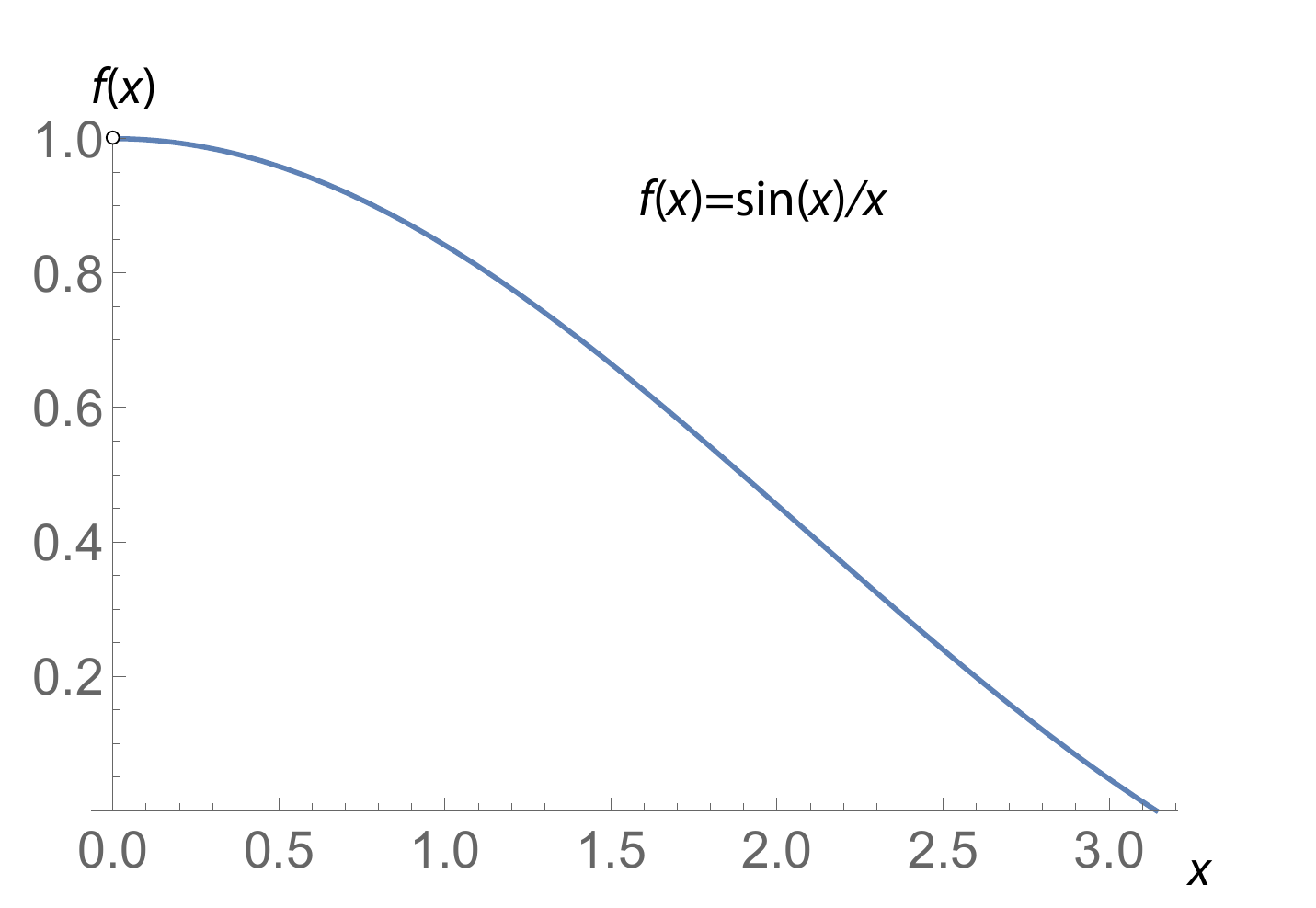}
\vspace{-3mm}
\captionof{figure}{The function $f(x)=(\sin x)/x$.\vspace{3mm}}
}

\noindent This function is perfectly well-defined all all points between $0$ and $\pi$ (that is, $D=(0,\pi]$), but it is \emph{undefined} at $x=0$.
\noindent One difficult concept arises in this kind of description, and this is encountered when you ask the following question:  What is the \emph{minimum} value of the set forming the domain?  Technically, the domain of this function has no minimum value!  The number $0$ is not in the domain (hence, the ``(" in the domain description).  However, for any small number $\epsilon$, say $\epsilon=1\times10^{-100000000}$, the value of the function $(\sin x)/x$ is very near 1.  It is only undefined at \emph{exactly} $x=0$.  In cases like this, where we can't technically use \emph{maximum} or \emph{minimum}, we say instead that the number 0 is the \emph{greatest lower bound} or \emph{infimum} of the domain.  Here, one needs to interpret ``greatest" as follows: of all of the lower bounds for for the interval where $f(x)$ is defined, $x=0$ is the largest of them.  Any larger value would be within the domain, and is thus not a bound.
\end{example}
\end{svgraybox}

The example above suggests that identifying the smallest number in an open set is not generally possible.  However, one can often identify the largest number that forms a lower bound for the set (in the example above, this is the number 0).  This leads to the following theorem (which we will not prove!), which we will use primarily as a definition.

\begin{theorem}[Least Upper Bound, Greatest Lower Bound]
If an interval has any lower bound, then it also has a \emph{greatest lower bound} (g.l.b., or \emph{infimum}).  That is to say, there is a unique number, $z_-$, which is a lower bound for the interval such no larger numbers are lower bounds.  Any numbers larger than $z_-$ are either (i) in the interval, or (ii) greater than or equal to an upper bound of the interval.

A similar statement can be made about the upper bound.  If an interval has any upper bound, then it also has a \emph{lease upper bound} (l.u.b, or \emph{supremum}).  The least upper bound is a unique number, $z_+$, which is an upper bound for the interval such no smaller numbers are upper bounds. 
\end{theorem}
\noindent To help make this more concrete, consider the interval we examined above: $(0,1]$.  For this interval, the greatest lower bound is the number 0, \emph{even though the number 0 is not in the interval itself.}  The least upper bound is the number 1.  No numbers, $z_0$, smaller than 1 are, obviously, upper bounds for the interval, because $z_0 <1$ and 1 is in the interval by definition.

Many of the functions that we encounter in applications are mappings that never tend toward infinite values.  In part, this is because physical systems never really have phenomena with infinite magnitude.  However, in many instances it can be useful \emph{model} a phenomenon as if it tended towards infinite value.  Additionally, some purely mathematical concepts require the definition of functions that tend toward infinity (e.g., the tangent function, $\tan(x)=\sin(x)/\cos(x)$, necessarily tends toward infinity as $\cos(x)$ tends toward zero).  It is useful to establish the concept of functions being bounded (and, hence, also being unbounded) as a matter of the vocabulary used in mathematics.

\begin{figure}[t]
\sidecaption[t]
\centering
\includegraphics[scale=.4]{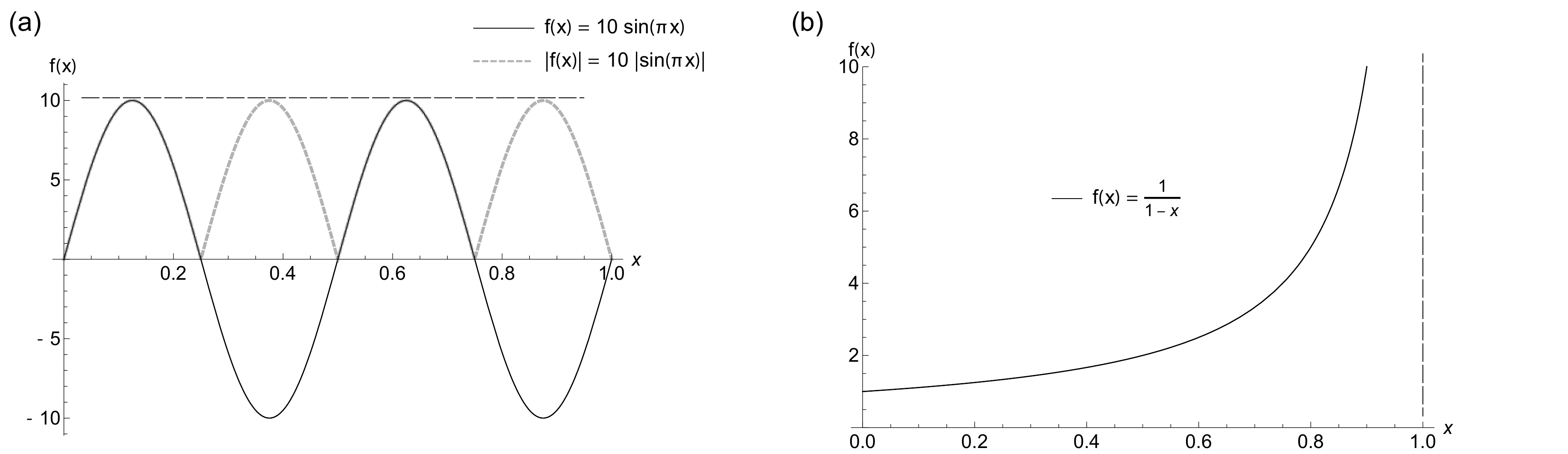}
\caption{(a) The function, $f(x) = \sin(4 \pi x)$, which is bounded on all possible intervals of the real line.  (b) The function, $f(x) = 1/(1-x)$, which is unbounded on the interval $x\in [0,1)$ as $x\rightarrow 1$.}
\label{boundedun}       
\end{figure}

\begin{definition}
A function, $f$, is \emph{bounded} on some interval $I=[a,b]$ if there is some number, $M$ such that $|f(x)| \le M$ for all $x \in I$.  
\end{definition}
\noindent Any function that is not bounded is unbounded.  An example of a bounded and unbounded function appears in Fig.~\ref{boundedun}.  In this figure, part (a) contains a bounded function.  The gray dashed line in the figure represents the absolute value of the function.  The horizontal line at $x=10$ indicates the minimum value of $M$ that bounds the function.  All values of the absolute value of the function are less than or equal to $M$, thus the function is bounded.  Note, in general there is no need to establish the \emph{minimum} value of $M$, as done here, to establish the boundedness of the function.  Had the choice been $M=12$, it would still be possible to show that all valued of the function are less than $M=12$, and the function is therefore bounded.  In Fig.~\ref{boundedun}(b), an the function $f(x) = 1/(1-x), ~0\le x < 1$ function is plotted.  This function is unbounded on this interval as $x\rightarrow 1$.  For every choice of $M$, we can always find a value of $x$ sufficiently close to 1 such that the value of $f(x)$ is larger than $M$.  Thus, even though the interval of definition does not contain the value of 1 (in which case one would find an infinity defined by $1/0$; this is sometimes called a \emph{singularity}), the function is still unbounded!  Note that while we often think of an unbounded function like $f(x)=1/(1-x)$ having a value of infinity at the singularity, it is more proper to think of the function as being \emph{undefined} at the point $x=1$.  The reason for this is the following: up to the point $x=1$, we can imagine a limiting behavior that is well defined.  In other words, it is correct to say 

\begin{equation}
    \mathop{\lim }\limits_{x\rightarrow 1} \frac{1}{1-x} =\infty
\end{equation}
because it indicates the limit of numbers that take the form $1/(1-x)$, where $x$ is any real number less than 1.  The ratio $1/(1-x)$ is always defined, because $x$ is never \emph{equal} to 1 in the interval defined.  This is in contrast to setting $x=1$; in this case, the function evaluates to $f(x=1) = 1/0$.  The form $1/0$ is not ``infinity" as it is sometimes (incorrectly) represented.  Technically, the expression $1/0$ has \emph{no definition} at all; it is an undefined operation.  While this last discussion is perhaps a bit overly-technical, there are a number of examples in applied mathematics where this seemingly technical point is enormously important in generating the proper understanding of a problem.

\subsection{Continuity}\label{continuity}

Functions are not necessarily always continuous.  For example, consider the following function.  Assuming $x\in[-2,2]$, define the function

\begin{equation}
f(x)=\left\{ 
{
\begin{array}{lll}
  5&~x=1 \\ 
   x^2 +2,&\mbox{otherwise} 
\end{array}
}
 \right.
    \label{function2a}
\end{equation}
\begin{figure}[b]
\sidecaption[t]
\centering
\includegraphics[scale=.5]{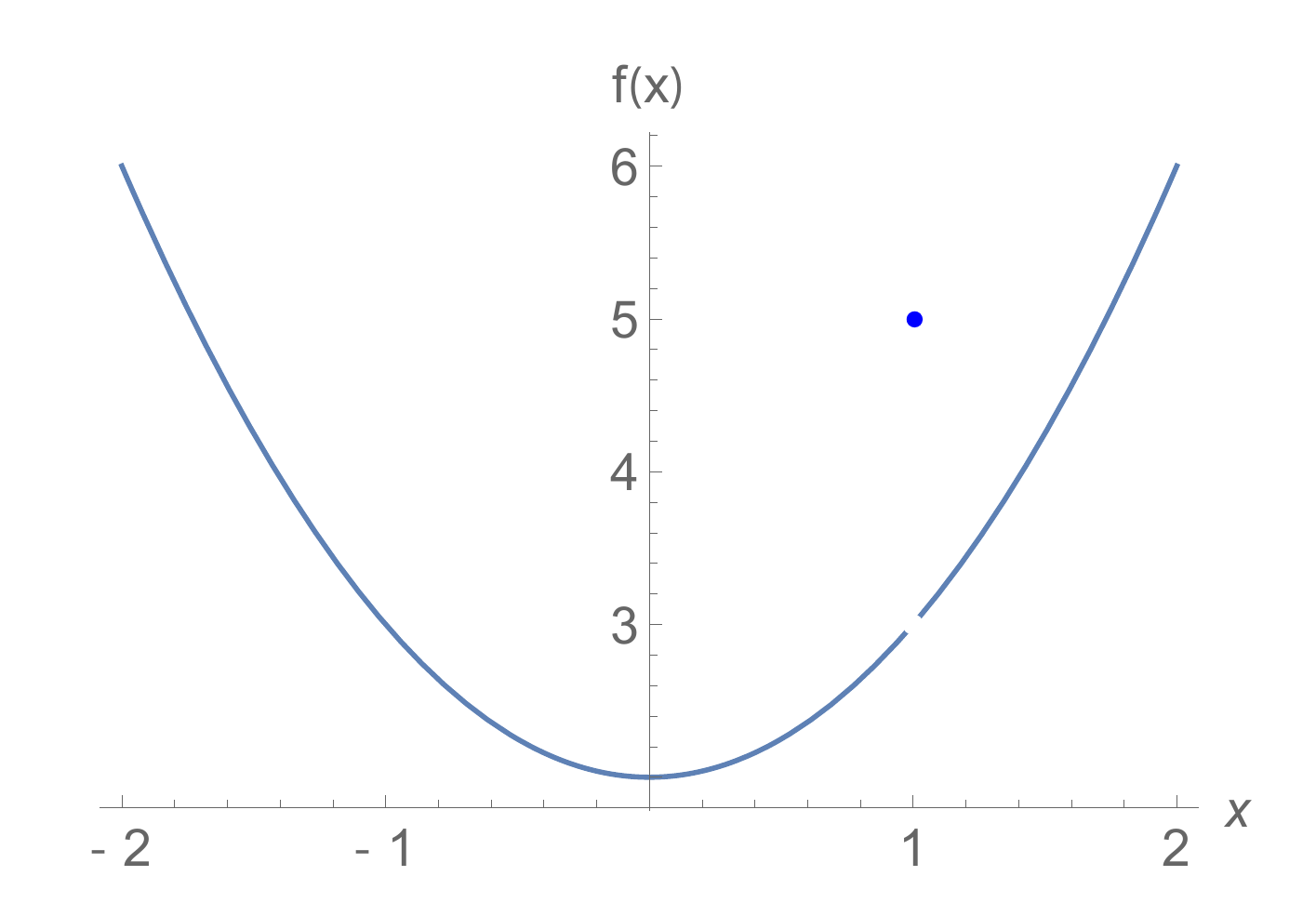}
\caption{A function.  This function has a single point of discontinuity.}
\label{fig:2}       
\end{figure}
This function is plotted in Fig.~\ref{fig:2}.  Although it is true that this function requires some extra handling (for the point $x=1$), the end result is not too dissimilar to what we are used to seeing.  However, it does indicate that it might be useful to further characterize functions on the basis of how continuous and smooth they are.  A few comments about the various labels that are applied to functions to describe how smooth they are are provided in the following.

A function $f(x)$ is \emph{continuous} if it changes gradually as the independent variable $x$ changes.  For a one-dimensional function, one can think of this as being able to draw the function with pen and paper without taking one's pen off the paper.  More formally, using the dreaded $\delta$-$\epsilon$ arguments, one make a more concrete statement.  Unlike most $\delta$-$\epsilon$ arguments, we are going to do this in slow steps.  First, suppose we are interested in the continuity of a function at a point $x=a$.  Now, we look at all of the points within a small distance, $\pm \delta$, around $a$.   We can denote those points by all of the values for $x$ such that $a-\delta < x < a+\delta$, or, equivalently $|x-a|<\delta$.  Now, for all of such points that we find on the $x-$axis (the domain), we can compute the absolute value $|f(x)-f(a)|$.  Suppose we do this, and we find that there is a number $\epsilon >0$ such that
\begin{equation}
    \left|f(x)-f(a)\right| < \epsilon~\mbox{for all $x$ such that } |x-a|<\delta
    \label{deltaepsilon}
\end{equation}
Now, for a continuous function, we expect that as we squeeze the interval $a-\delta < x < a+\delta$ by making $\delta$ smaller, we should also be able to correspondingly find a new value for $\epsilon$ that meets the criterion given in \eqref{deltaepsilon}.  In other words, as $|x-a|$ gets smaller, so does $|f(x)-f(a)|$.  This is essentially the way that continuity of a function is defined.  Formally, the argument is turned inside-out (which is common in the presentation of definitions and proofs in mathematics...), to read as follows

\begin{definition}
A function, $f$, is \emph{continuous} at a point $x=a$ if and only if (i) $a$ is in an open interval of the domain of $f$, and (ii) we can always pick an $\epsilon$  (as small as we like), and still always find some interval around $a$ (that is, $a-\delta < x < a+\delta$) such that $|f(x)-f(a)|<\epsilon$, regardless of how small we have chosen $\epsilon$ to be.
\end{definition}

\begin{definition}[Continuous Functions]
A function is \emph{continuous} on and interval, $I$, if it is continuous for \emph{each open sub-interval} of $I$.\indexme{function!continuous}
\end{definition}
In essence, the ``for each open sub-interval" part of the definition is to avoid any possible problems with defining continuity at the least upper or greatest lower bound of the interval.  In short, we don't need to worry about continuity at the end points of our domain.  This is designed specifically to avoid potential ambiguities as were discussed for the function $f(x)=1/(x-1),~0\le x < 1$ in Section \ref{boundedness}.

\begin{figure}[t]
\sidecaption[t]
\centering
\includegraphics[scale=.5]{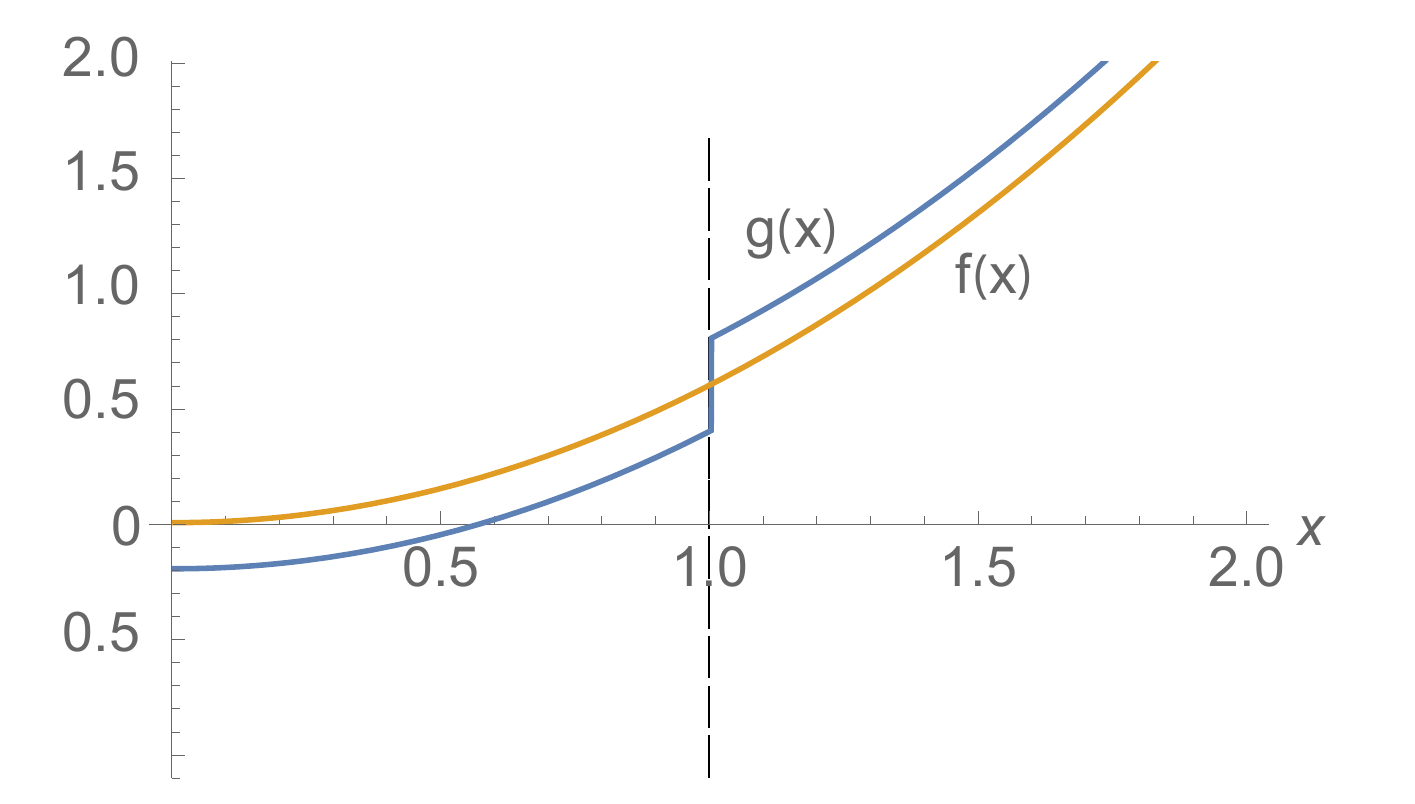}
\caption{Two functions, $f$ and $g$.  The function $g$ has a discontinuous jump in the value of $\tfrac{1}{2}$ at $x=1$.}
\label{fig:3}       
\end{figure}

These definitions can be a bit tricky to apply, and  they are not the only way to define continuity.  Alternatively, one can insist that the following limits, taken from the left and right sides of $a$ must be valid for each point in a domain for a function to be considered continuous.

\begin{align}
    \mathop {\lim }\limits_{x^+\rightarrow a}& |f(x)-f(a)|\implies 0 \\
    \mathop {\lim }\limits_{x^-\rightarrow a}& |f(x)-f(a)|\implies 0
\end{align}
And sometimes this alternative is easier to use in applications. Note.  The long double arrow ($\implies$) in mathematics should be interpreted as meaning ``this implies". 
Regardless of which method is used, the results are the same.  This is easiest to see with a simple example, as illustrated in Fig.~\ref{fig:3}.  \\

\begin{svgraybox}
\begin{example}[Continuity of Functions]

In Fig.~\ref{fig:3}, it is apparent from inspection that the function $g(x)=\tfrac{1}{2}x^2$ (the orange line on the Fig.~\ref{fig:3}) is continuous.  Suppose we consider the point $x=1$.  Clearly, we can pick a small number, say $\epsilon=1/1000$, such that there is a small enough interval around $x=1$ ($1-\delta<x<1+\delta$) where $|f(x)-f(1)|<1/1000$.  In fact, we can compute what this interval is as follows

\begin{align*}
    |f(1+\delta)-f(1)|&<1/1000 &&\\
    \intertext{Accounting for the absolute value, there are two options} \\
     \tfrac{1}{2}(1+\delta)^2-\tfrac{1}{2} & <1/1000 \mbox{ ~~~~or }& \tfrac{1}{2}-\tfrac{1}{2}(1-\delta)^2 & <1/1000\\
     \intertext{Solving these inequalities for $\delta$, we find}
     |\delta| &< \sqrt{\frac{1002}{1000}}-1& |\delta| &< 1-\sqrt{\frac{998}{1000}}
     \intertext{The first of these is the smallest, so it is safe to take $\delta$ equal to that value.}
\end{align*}
You can check this directly by noting that 
\begin{align*}
    \left|f\left(1+\sqrt{\frac{1002}{1000}}-1  \right) - f\left( 1 \right) \right|& = \frac{1}{1000} \\
     \left|f\left(1-\sqrt{\frac{1002}{1000}}-1  \right) - f\left( 1 \right) \right|& = 0.000999<\frac{1}{1000}
\end{align*}

No matter how small we make $\epsilon$, we can always find a value of $\delta$ so that the inequality is valid.  \\

Now consider the function $f$ (the blue line on Fig.~\ref{fig:3}).  Here, our scheme works reasonably well until we set $\epsilon$ to be less than about $1/2$.  Because there is a jump of $1/2$ right at $x=1$, we can easily find an $\epsilon$ that breaks our definition.  For example, set $\epsilon=1/10$.  Now we are looking for the values of $x$ in an interval $a-\delta \le x \le a+\delta$, such that $|f(x)-f(a)|<1/10$.  Of course, we can see just by looking at the graph that there is no such interval.  In any small interval around $x=a$, the value of $|f(x)-f(a)|$ is \emph{at least} $\tfrac{1}{2}$, but definitely never smaller than that value.  

The definition of continuous, then, aligns with our intuitive notion, even if the definition itself takes a bit of thinking to fully understand.
\end{example}
\end{svgraybox}

Functions that have a jump in them, such as the one in Fig.~\ref{fig:3} are called \emph{discontinuous}.  When there are a \emph{finite} number of discontinuities, sometimes the functions are called \emph{piecewise continuous}.  In other words, a piecewise continuous function is a function that has a finite number of jumps in it and doesn't blow up to $\pm \infty$ anywhere.   This is an important enough concept that it deserves a specific definition.

\begin{definition}[Piecewise Continuous Function]
A \emph{piecewise continuous function} is a function that is at continuous everywhere, except at a \emph{finite} number of points. \indexme{function!piecewise continuous}
\end{definition}

\begin{figure}[t]
\sidecaption[t]
\centering
\includegraphics[scale=.5]{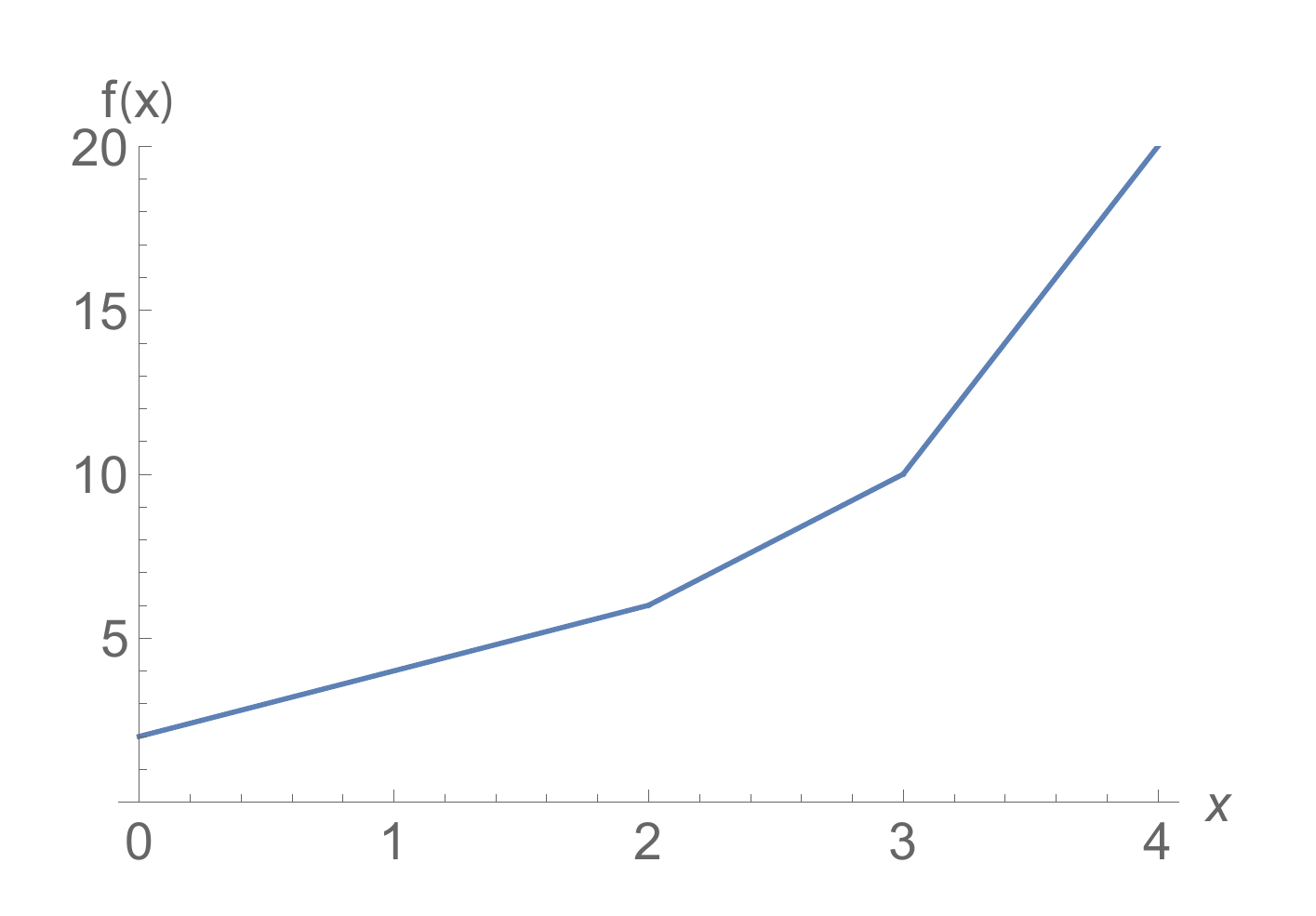}
\caption{A continuous function with discontinuous derivatives at a two points.  This kind of function is sometimes denoted $C^0$ to indicate that it is continuous, but it does not derivatives at all points in the domain.}
\label{fig:5}       
\end{figure}

Within continuous functions, the functions are sometimes by their derivatives.  As an example, examine Fig.~\ref{fig:5}.  This function is defined by

\begin{equation}
   f(x)= \begin{cases}
 2 x+2 & 0\le x<2 \\
 4 \left(x-\frac{1}{2}\right) & 2\le x<3 \\
 10 (x-2) & 3\le x \le 4
\end{cases}
\end{equation}

In general, the word \emph{smooth} is used to indicate a function that has a derivative at each point in its domain.  Thus, the function given by Eq.~\eqref{function} is smooth; the function given by Eq.~\eqref{function2a} is technically \emph{non-smooth}, but this arises because of the discontinuity imposed by a single point.  Extending the idea of piecewise functions, we can call the function given by Eq.~\eqref{function2} \emph{piecewise smooth}.  Most of the functions that we study in this text will be of this kind.

\begin{definition}[Piecewise Smooth Function]
A \emph{piecewise smooth function} is a function where the first derivative of the function is well-defined (i.e., it exists, and it is not infinite) everywhere, except at a \emph{finite} number of points.
\end{definition}

\begin{figure}[t]
\sidecaption[t]
\centering
\includegraphics[scale=.5]{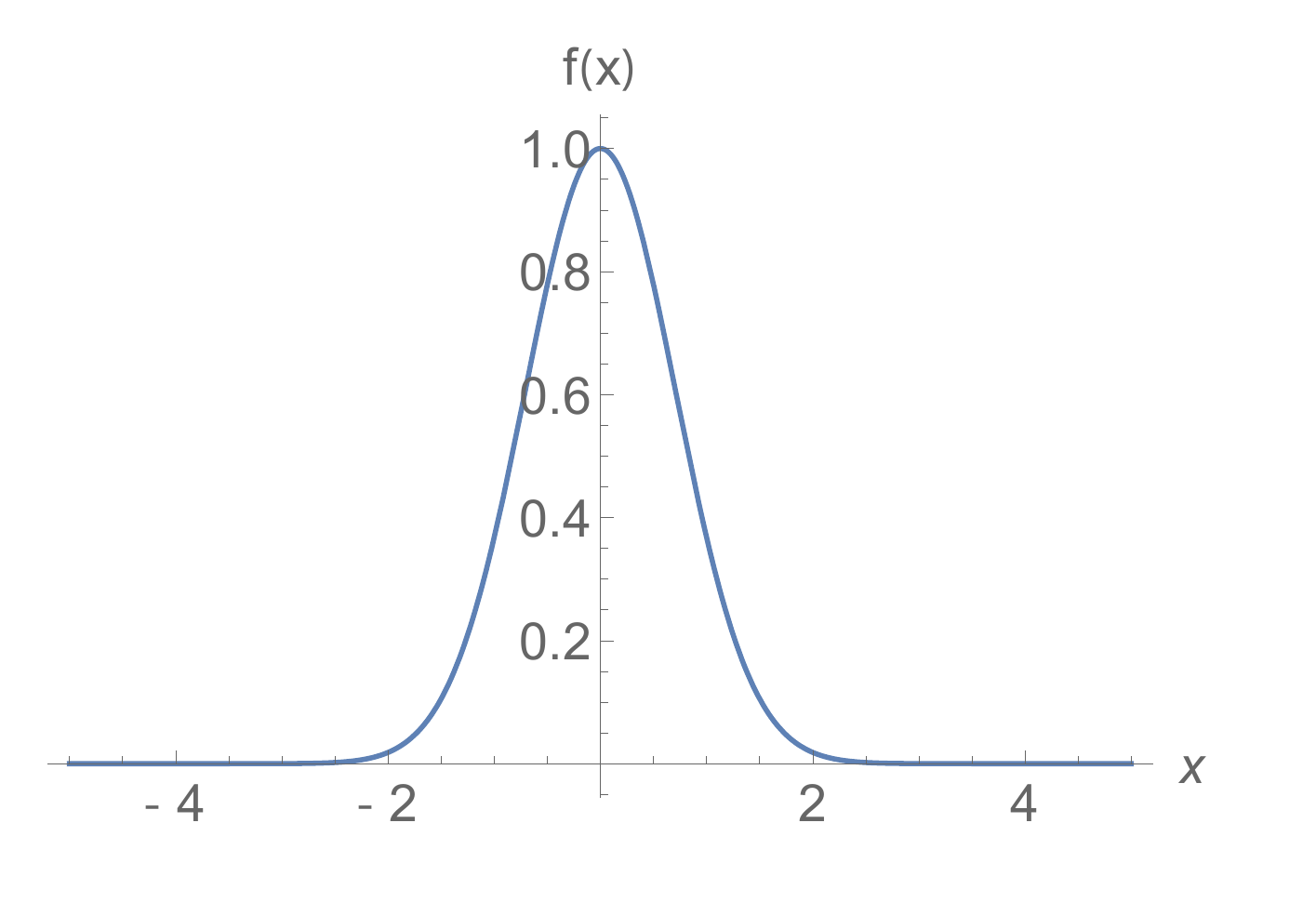}
\caption{A highly continuous function, $f(x) =1/\sqrt{\pi}\exp(-x^2)$.  This function  has derivatives of all order; thus it is called a $C^\infty$ function.  It also happens to be an analytic function (note, however, that there exist $C^\infty$ functions that are not analytic!).}
\label{fig:x2}       
\end{figure}

We can further characterize smoothness by how many derivatives a function has.  If a function has an infinite number of derivatives that exist, then it is called $C^\infty$.   An example here is the function $f(x)=\exp(-x^2)$.  This function is smooth (it is plotted in Fig.~\ref{fig:x2}), and is often called a \emph{Gaussian} function.  We can define particular $C^\infty$ functions that are called \emph{analytic} functions; however this will have to wait until after the review of derivatives and integration is presented.

\begin{figure}[t]
\sidecaption[t]
\centering
\includegraphics[scale=.5]{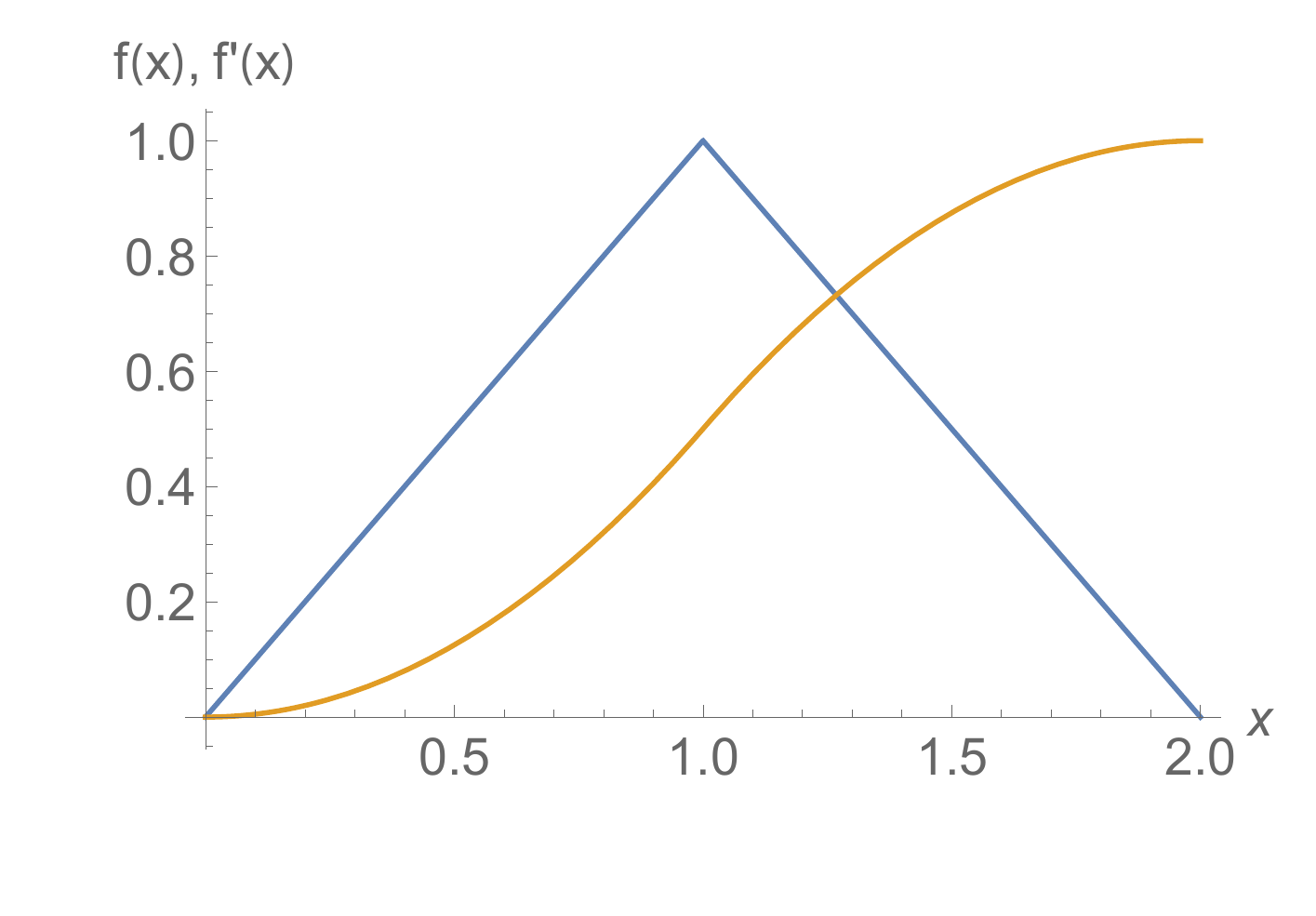}
\caption{A continuous function (orange) with continuous first derivatives (blue).  However, the second derivatives are not continuous at $x=1$.  This function is called a $C^1$ function to indicate that it has continuous first derivatives everywhere in the domain.}
\label{fig:06}       
\end{figure}

Some functions have only a limited number of derivatives that exist.  As an example here, consider the function

\begin{equation}
   f(x)= \begin{cases}
 \frac{x^2}{2} & 0\le x<1 \\
 - \frac{x^2}{2}+2 x -1 & 1\le x\le 2 
\end{cases}
\end{equation}
This function \emph{looks} innocuous enough- it is plotted as the orange curve in Fig.\ref{fig:06}.  However, its derivative (plotted in blue) has a cusp in it at the location $x=1$. Therefore, although its first derivative exists and is continuous, its second derivative is not defined at $x=1$.  To reflect this, functions with one continuous derivative are called $C^1$.  Although we have not yet reviewed derivatives, we note the following generalization of this idea.

\begin{definition}[$n^{th}$-order continuous function]
Suppose a function $f$ is such that the first $n$ derivatives (i.e., $f'$, $~f''$,$~f''' \ldots$ $~f^{(n)}$) are both \emph{bounded} and exist (i.e., there are no points such that the derivative of order $(n-1)$ generates a cusp or other discontinuity).  Then the function is continuous in the derivative up to order $n$, or, more simiply, $C^n$.
\end{definition}

\begin{figure}[t]
\sidecaption[t]
\centering
\includegraphics[scale=.5]{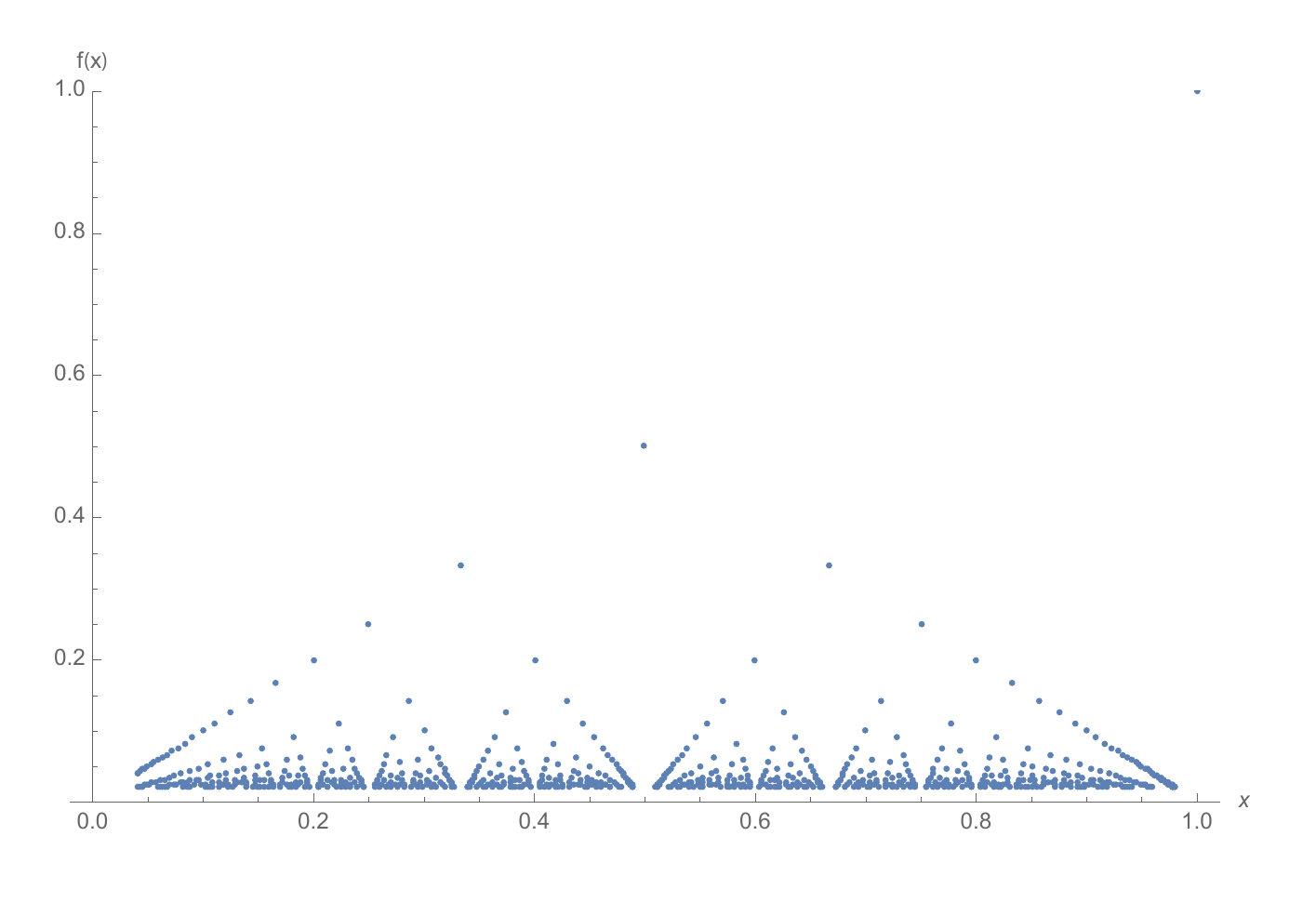}
\caption{A graph illustrating a collection of a values for the Thomae function using $2\le p\le 50$, $2\le q \le 50$.  This may look super weird, but it meets the criteria for a function.}
\label{thomae}       
\end{figure}

There is a special class of functions called \emph{analytic} that have a number of interesting properties.  Many of the familiar functions that we know about are analytic. Examples include
\begin{enumerate}
    \item All polynomials of finite degree.
    \item The exponential and logarithmic functions.
    \item The trigonometric functions.
\end{enumerate}
In the early days of mathematical evolution (say, through the early 1900s), analytic functions were synonymous with \emph{functions}.  Since that time, the notion of what constitutes a function has grown considerably; some of the unusual examples of functions that are not analytic in some part of their domain (or everywhere, in the case of the Thomae function given in Fig.~\ref{thomae}) are presented in the material later in the text.

\begin{definition}[Analytic Function-- Definition 1]
Suppose a function $f$ is defined on a domain, $D$. The function $f$ is \emph{analytic function} everywhere in $D$ if for every closed interval $K$ that is a subset of $D$ ($K\subset D$) there exists a constant $C$ such that for every point $x\in K$ and every non-negative integer, $k$, the following bound holds

\begin{equation}
\left | \frac{d^k f}{dx^k}(x) \right | \leq C^{k+1} k!
\end{equation}
\indexme{function!analytic} \indexme{analytic function}
\end{definition}
While this definition does make a mathematical mouthful, it is relatively easy to understand intuitively.  What the definition is trying to tell us is that, for any closed finite interval that is part of the domain of the function, that the function, nor any of its derivatives, go to infinity in that domain.  There are some technical issues regarding the domain of the function (and the closed subsets that are selected from it) that will not be covered here.  In general, however, we can think of a function as being analytic at a point in its domain if its value and the value of all of its derivatives do not tend toward infinity.

As a final note about the classification of functions, there is one additional classification that is useful to know about.  There are, in one sense, two different kinds of functions that we use in common practice.  There are polynomial functions, and the roots of polynomial functions to start with. A polynomial function is a polynomial (of any finite order) whose coefficients are also polynomials.  For example, a polynomial function of order 4 is defined by

\[ a_4(x)x^4+a_3(x)x^3+a_2(x)x^2+a_1(x) x + a_0(x) = 0 \]

\noindent For such a function, the  roots generate new functions, $f(x)$, involving (rational) fractional powers (including negative powers). Such functions are called \emph{algebraic} because they can be defined using only the rules of algebra applied in a \emph{finite} algorithm.  The following are all examples of algebraic functions.

\begin{align*}
    f(x)&=x\\
    f(x)&=x^{-1} = \tfrac{1}{x}\\
    f(x)&=\sqrt{x}\\
    f(x)&= \frac{\sqrt{1+x^3}}{x^{3/7}-\sqrt{7}x^{1/3}}
\end{align*}
The nice thing about these functions is that they can be entirely described mathematically simply by describing the polynomial that generates them. Algebraic functions are expressions involving only a finite number of terms, and using only the algebraic operations addition, subtraction, multiplication, division, and raising to a rational fractional power.   In a sense, we fully ``understand" these functions as long as we understand the operations that define them.  This is literally constructive- for many examples of algebraic functions, we can give someone an algorithm (with a finite number of terms) that explains what the function is and how to compute its values.  Importantly, all algebraic functions are given by the roots of some polynomial; however, for rational polynomials of degree 5 or higher, it is not true that all such polynomials have roots that are algebraic functions (the proof of this is known as the \emph{Abel–Ruffini theorem}).

It turns out, however, that there are many interesting and useful functions that are not the roots of any polynomial equation.  Familiar examples include the $\sin x$, $\cos x$, $\ln x$, and $e^x$ functions.  Such functions are called (and, yes, this is really the name) \emph{transcendental functions}.  The idea behind the name is that these functions transcend description by the discipline that we normally call algebra.  Transcendental functions will show up later in our studies of differential equations, and some will prove to be essential for describing solutions to such equations.  All transcendental functions are \emph{analytic} (a term we will define in by the Taylor series in detail later), so they all have convergent Taylor series representations.

There are functions that are even more \ldots \emph{interesting} to attempt to defined because their structure begins to challenge the concept of function altogether.  Consider the following function (sometimes called the modified Dirichlet function or the Thomae function)

\begin{equation}
f(x)=\left\{ 
{
\begin{array}{lll}
  0&~\mbox{$x=0$} \\ 
  \tfrac{1}{q}&~\mbox{$x\in\mathbb{Q},~x=p/q$} \\
  0& ~\mbox{$x\in\mathbb{P}$}
\end{array}
} \right.
    \label{function2}
\end{equation}
It is understood here that $x=p/q$ is expressed in fully reduced fraction form.  This is a really unusual function, and clearly it has no simple, closed algebraic formula. A plot of this function appears as Fig.~\ref{fig:3}.  In addition to its strange definition and look, it also has some other unusual characteristics.  It turns out that the real numbers are ``denser" than the rational numbers.  In a sense, there are  more (many, many more) irrational numbers than rational ones.  Thus, this function is continuous (and has a derivative) at each point where $x$ is an irrational number.  It is discontinuous at every point where $x$ is rational.  It is not terribly important to understand from this example all of the details of why this kind of behavior exists (although our discussion of intervals above is suggestive), but the essential idea is to build an intuition that the real (irrational) numbers constitute a much larger (denser) set than the rational numbers do.

\subsection{Equivalence Classes of Functions}\label{equivalence_class}
\begin{figure}[t]
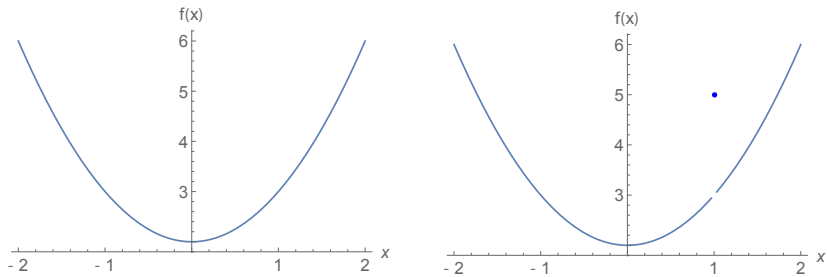

\sidecaption[t]
\centering
\includegraphics[scale=.4]{\CHAP/review_01.pdf}
\includegraphics[scale=.4]{\CHAP/review_02.pdf}
\caption{Two functions that belong to the same equivalence class.  The second differs from the first at only a single point.}
\label{fig:equiv}       
\end{figure}

Later on, when we discuss the concept of Fourier series, it will be helpful to understand the notion of \emph{equivalence classes} of functions.  In short, two functions $f$ and $g$ are considered to belong to the same equivalence class if they differ from each other by, at most, a finitely countable number of points.  An example of two functions that are of the same equivalence class is given in Fig.~\ref{fig:equiv}.  The mathematical representation of these two functions is given by \indexme{function!equivalence class}\indexme{equivalence class}

\begin{align*}
    f(x) &= x^2 \\
    g(x) &=\left\{ {\begin{array}{*{20}{c}}
{{x^2}~for~x \ne 1}\\
{5~for~x = 1}
\end{array}} \right.
\end{align*}

Note that we have the following limits for $g$ (the function on the right-hand side) near $x=1$

\begin{align*}
    lim_{x\rightarrow 1^+} ~g(x) &= 1 \\
        lim_{x \rightarrow 1^-} ~g(x) &= 1 
\end{align*}
where $1^+$ indicates the limit coming from the right-hand side, and $1^-$ represents the limit coming from the left-hand side.  Be sure to recall here that the limit exists in the sense of \emph{approaching} the value $x=1$, but not actually reaching that value.  If you need a refresher on one-sided limits, your undergraduate calculus text will definitely cover this material.

It turns out that, from the perspective of integration (including Riemann integration-- see \S\ref{integration}), the presence of a single discontinuous point does not affect the result of the integration. Speaking colloquially, a single point has no measure, so it does not add or subtract from the integral.  The following definitions makes this more formal.

\begin{definition}[Removable discontinuity]
A function is said to have a set of removable discontinuities if (a) the number of discontinuous points is finite, and (b) the left and right-hand limits at each point of discontinuity are equal.  The value of the function at the discontinuous points can be taken to be the value of the limit at that point. \indexme{discontinuity!removable discontinuity}
\end{definition}

\begin{definition}[Equivalence class of functions]
Two functions $f$ and $g$ are said to belong to the same equivalence class if (a) the two functions $f$ and $g$ have only removable discontinuities, and (b) $f-g$=0 for all points that are not located at a discontinuity in either $f$ or in $g$.  Alternatively, (b) could be stated by ``$f-g=0$ after all discontinuities are removed". 
\end{definition}

\section{Linear Spaces, Linear Operators, and Linear Algebra}

We will have only a few occasions to use linear algebra in this material, but it is worth briefly reviewing, and it also allows us to discuss the concept of an \emph{linear operator}.  Many of the problems that we encounter in the material to come involve linear operators, so spending some time understanding what they are is worthwhile.

Operators are just what they sounds like: They are mathematical constructs that operate on something (say, a function)  to generate something else (for example, a different function) (Fig.~\ref{fig:op}).  This is a very general concept, and it is difficult to give it a precise meaning.  As an example, all functions are operators.  To make this concrete, take a look at the following example.

\begin{svgraybox}
\begin{example}[Functions as operators.]

We are all familiar with the idea of a \emph{function}.  In the notation that we have are probably most familiar with, a function is defined as a one-to-one (i.e., unique) mapping between two sets.  In a more familiar setting, we think of a function as having a \emph{range} (which is usually an interval on the real line, plotted on the horizontal axis by convention) and a domain (which is usually some portion of the real line, plotted on the vertical axis by convention).  For example, a properly defined function would be any of the following:
\begin{align}
    f(x)&= 5x, &&0 < x < 10 \label{fx2}\\
    g(x)&=\sin[2 \sin(2 \sin\{2 \sin(x)\})]  &&\infty < x < \infty \\
    h(x)&=\sin[1/x] &&-0.1 < x < 0.1 
\end{align}
These functions are plotted in Fig.~\ref{f:functions}.  

\begin{centering}
\includegraphics[scale=1]{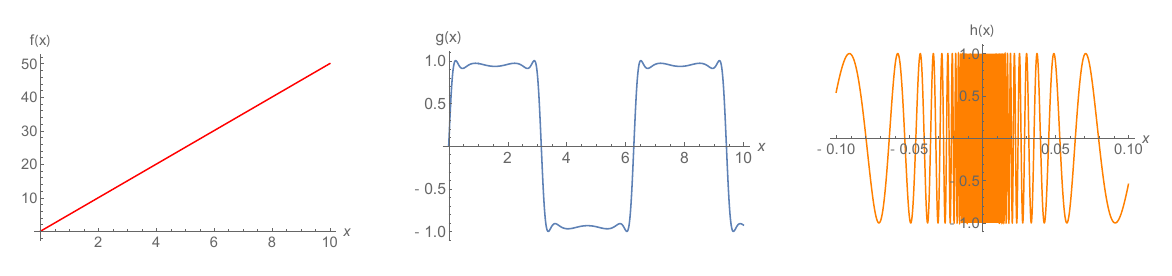}
\captionof{figure}{Figure for Example 1.2.  Three functions.}
\label{f:functions}  
\end{centering}
%
\noindent What is important to realize about these functions is that they are operations on their independent variables.  In some sense, a function exists \emph{separately} from their independent variables.  Thus, once I have defined, $f(x)$, as, for example, in Eq.~\eqref{fx2}, then what I have really done is to define an \emph{operation} to be performed on one set (the independent variables) to generate another set (the functional value). Often when we think of functions as being operators, we consider (1) the function to define the operation, and (2) the operation is done on whatever is put into the function.  The thing put into the function, when used in this sense, is often called the \emph{argument}.  The argument can be a simple independent variable defined over an interval (our conventional use of functions) or some other set.  

For example, suppose we define the following set: $X=\{x: x=1,3,4,5,42\}$.  We would read this in words as "The set $X$ is the set that contains values of $x$ such that $x$ is equal any of the values 1, 3, 4, 5, and 42".  The set $X$ has a finite number of values in it, but we can still use any of the functions above as operations on this set.  Let's interpret $f(X)$ the following way:\vspace{-5mm}

\begin{align*}
    f(X)&=\{f(1),f(3),f(4),f(5),f(42)\}\nonumber\\
    \intertext{or, computing the values using the expression $f(x)=5x$ \vspace{-5mm}}\nonumber\\
    f(X)&=\{5,~15,~20,~25,~210\} 
\end{align*}
We can even think about putting functions in other functions (creating a  as arguments when we use the operator idea.  As an example, consider the function $p(z)=z^2 + 1$.  Now, can we make sense of the operations $f(p)$ and $h(p)$?  Sure, we need only apply the operations to the  argument of the function, regardless of what the argument is. 

\begin{align*}
    f(p)&=5 p(z) & h(p)&=\sin(1/p(z))\\
    &= 5(z^2+1) &&=\sin\left(\tfrac{1}{z^2+1}\right) \\
    &= 5z^2 + 5 &&
\end{align*}
When we interpret a function $f(x)$ this way, we sometimes write the function without the independent variable (with the idea that the independent variable can be anything), as in simply $f$.  This is just a notational device used to be more compact, there is usually no deeper meaning associated with it. 
\end{example}
\end{svgraybox}

%
\subsection{Linear Operators}\label{lineop}
%
Now that we understand the basic idea of an operator, we can consider what a \emph{linear operator} is.  This is actually pretty simple at this point.

\begin{definition}[Linear Operator]
A \emph{linear operator}, $\mathscr{L}$, is an operator that subscribes to the properties of \emph{additivity} and \emph{homogeneity}.  Specifically, this means
\begin{align}
   additivity:& \mathscr{L}(f+g) = \mathscr{L}(f)+\mathscr{L}(g) \label{add}\\
   homogeneity:& \mathscr{L}(\alpha f)= \alpha \mathscr{L}(f)\label{homog}
\end{align}
\end{definition}

Linear operators are a subset of all possible operators.  In fact, most interesting phenomena in engineering and physics are nonlinear in general.  Linearity exists only as an approximation to the more general nonlinear behavior.  In this text, we will be generally concerned with (although not exclusively!) linear operations and linear operators. 

As discussed above, functions can themselves be though of as operators.  Thus, a linear function would subscribe to the properties of additivity and homogeneity above.  As and example, the following function is linear

\begin{equation}
    f(x) = 5 x
\end{equation}

To show that it is linear, note that we can prove both additivity and homogeneity in one step. If we make the substitution $x\rightarrow \alpha y + \beta z$, then we find

\begin{align}
    f(\alpha y+\beta z) &= 5(\alpha y+\beta z) \nonumber\\
    & = \alpha 5y + \beta  \nonumber \\
    &= \alpha f(y) + \beta f(z)~~~~(\therefore\textrm{~linear})
\end{align}
illustrating that the function is indeed linear.

\begin{figure}[t]
\sidecaption[t]
\centering
\includegraphics[scale=.5]{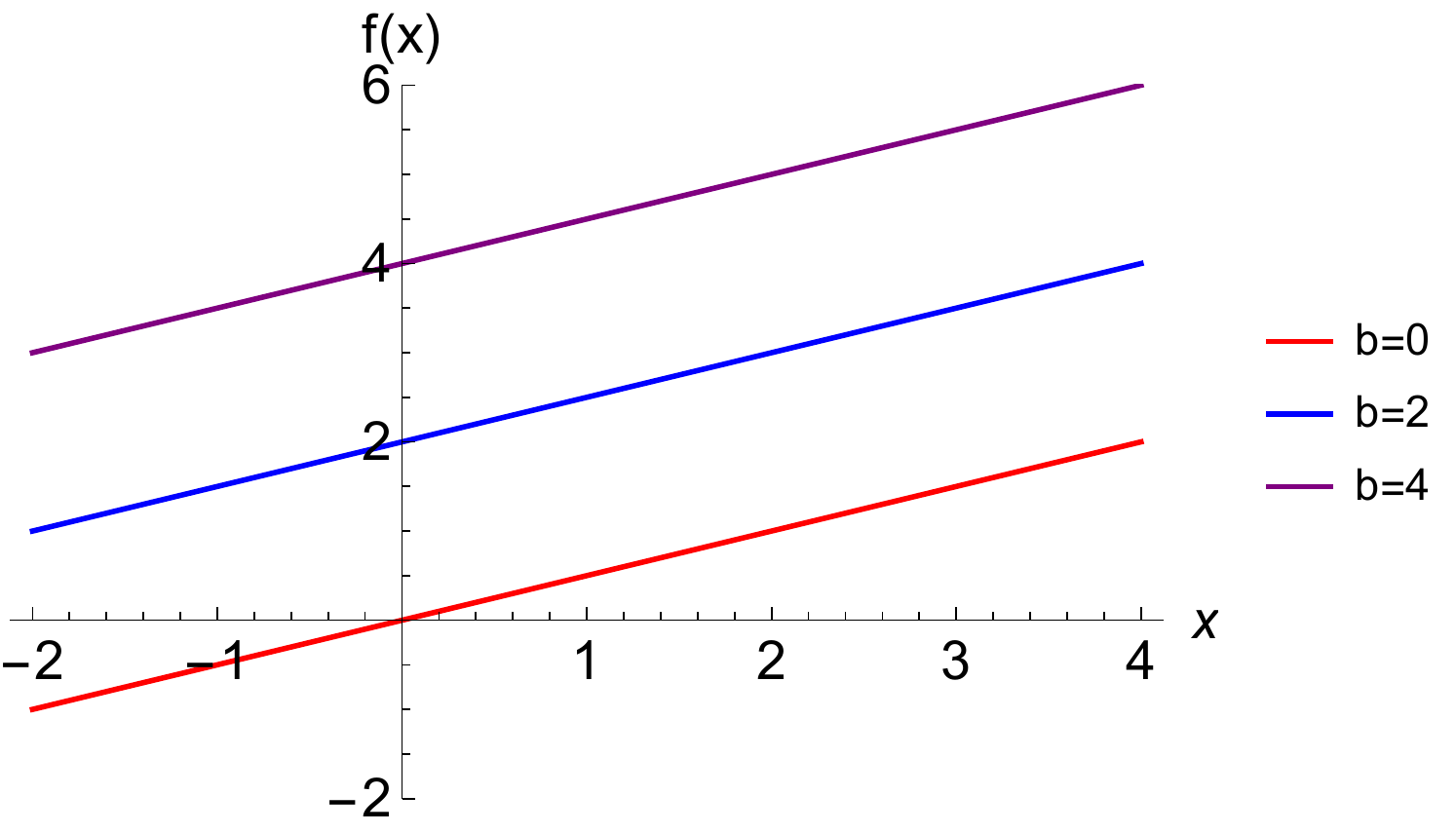}
\caption{Translations of the function $y=\tfrac{1}{2}x$.}
\label{fig:linearity}       
\end{figure}

There is a particular feature of linearity that can create some confusion.  We are used to calling functions such as $y=m x + b$ \emph{linear}.  However, it is not difficult to show that this function does not meet the properties of additivity and homogeneity (try it!).  There is a subtle reason for this failure.  In a very real sense, $y=m x + b$ is actually the \emph{translation} of a more fundamental function $y= mx$.   Take a look at Fig.~\ref{fig:linearity}.  Each of the functions illustrated is a translation of the function $y=\tfrac{1}{2} x$.  Such translations are called \emph{affine transformations}.  It turns out that when we call a function \emph{linear}, we actually mean that the fundamental function itself (where $b=0$) is linear.  The constant term only represents an affine transformation of this more fundamental function.  In fact, in each of the lines defined in Fig.~\ref{fig:linearity}, we could eliminate the constant term $b$ by simply making a transformation of the coordinate system.  For example, if we moved the coordinate system vertically by a distance of 2, then the blue line in the figure would then have $b=0$.  Thus, in a sense, all affine transformations of a function are the same as far as \emph{linearity} is concerned.  In order to assess linearity, the first step would be to first make an affine transformation of the coordinate system so that $b=0$.  In a more practical sense, we can essentially ignore constants when checking operators (including functions) for linearity.

\begin{svgraybox}
\begin{example}[A familiar example: The derivative as an operation.]
The process of differentiation is one example of a familiar operation.  Consider, for example, the following.
\begin{equation*}
\mathscr{L}(f(x)) = \frac{d}{dx}\left( f(x) \right) 
\end{equation*}
In this example, the \emph{operator} is defined by derivative notation, $\tfrac{d}{dx}$.  When we say $\mathscr{L} \equiv \tfrac{d}{dx}$, we are just defining the fact that the ``abstract" operator in this case is the derivative operation, denoted by $\tfrac{d}{dx}$.  
\end{example}
\end{svgraybox}
%
\begin{svgraybox}
\begin{example}[The derivative as a \emph{linear} operator.]
Is the derivative operator a \emph{linear} one?  We can check. Recall the definition.
\begin{equation}
    \mathscr{L}(\alpha f(x) + \beta g(x))=\alpha\mathscr{L}(f(x))+\beta \mathscr{L}(g(x)) \label{linop}
\end{equation}
Since our operator in this case is just $\tfrac{d}{dx}$, we need only to check that the linearity identity is met.  In other words, we need to evaluate
\begin{align*}
\mathscr{L}(\alpha f(x) + \beta g(x))&=\tfrac{d}{dx}(\alpha f(x) + \beta g(x))\\
    &=\frac{d}{dx}(\alpha f(x))+\frac{d}{dx}(\beta g(x)) \\
    &=\alpha\frac{d}{dx}(f(x))+\beta\frac{d}{dx}(g(x))
\end{align*}
This is exactly the form that a linear operation must take as defined by Eq.~\eqref{linop}, so the derivative is a linear operator.
\end{example}
\end{svgraybox}

A compound linear operator is just the sum or one or more linear operators. The easiest way to understand these is through some simple examples.  Although we will not use operator notation extensively in this text, it is a useful notation to understand.

\begin{svgraybox}
\begin{example}[Compound linear operators]
Compound linear operators act as follows.
\begin{equation}
    \left(\mathscr{L}_1+\mathscr{L}_2\right)( f(x))\equiv \mathscr{L}_1( f(x))+\mathscr{L}_2( f(x))
     \label{linopcomp}
\end{equation}
It is important to recognize that, although the notation looks as though it indicates multiplication, it does not!  Here is a more complex looking example that helps to understand the idea better.
\begin{align*}
\left( \frac{d}{dx}+\frac{d}{dy}\right)(x^2 y) &= \frac{d}{dx}(x^2 y)+\frac{d}{dy}(x^2 y) \\
&= 2xy+x^2
\end{align*}
Note that the operators \emph{do not} act like multiplication!
\begin{align*}
\left( \frac{d}{dx}+\frac{d}{dy}\right)(x^2 y+y^3) &= \frac{d}{dx}(x^2 y+y^3)+\frac{d}{dy}(x^2 y+y^3) \\
&= 2xy+x^2+3y^2
\end{align*}
\end{example}
\end{svgraybox}

Functions can be treated as operators also; but only \emph{homogeneous} functions are linear operators.

\begin{svgraybox}
\begin{example}[Functions as operators.]
Consider the functions defined in the previous example.  
\begin{align}
    f(x)&= 5x, &&0 < x < 10 \label{fx}\\
    g(x)&=\sin[2 \sin(2 \sin\{2 \sin(x)\})]  &&\infty < x < \infty \\
    h(x)&=\sin[1/x] &&-0.1 < x < 0.1 
\end{align}
Are they linear operators?  To check this, we need only try the operations on the quantity $\alpha y + \beta z$ to see if the conditions given by Eqs.~\eqref{add}-\eqref{homog} are met.  For the generic operator $L$, we will substitute our particular functional operators.  Thus, for the example of $f$, we have the result

\begin{align*}
    f(\alpha y + \beta z)&= 5(\alpha y + \beta z) \\
    &= \alpha 5y + \beta 5z \\
    &= \alpha f(y) +\beta f(z)
\end{align*}
Which indicates that $f$ is linear.  For the function $g$, we find

\begin{align*}
    g(\alpha y + \beta z)&= \sin[2 \sin(2 \sin\{2 \sin(\alpha y + \beta z)\})]\\
    &\ne \alpha g(y) +\beta g(z)
\end{align*}
so $g$ is not a linear operator (and also not a linear function!).  The final function $h$ is left to the reader to check.
\end{example}
\end{svgraybox}
%
There are two more concepts that are useful to introduce when discussing the (somewhat abstract) concept of linear operators.  These are the \emph{identity operator} and the \emph{inverse operator}.

\begin{definition}[Identity operator]
The identity operator, $\mathscr{I}$, is any operator such that $\mathscr{I}(f)= f$ for all admissible objects in the domain of $\mathscr{I}$ (e.g., numbers, functions, vectors, etc.) $f$.   
\end{definition}

\begin{svgraybox}
\begin{example}[Identity operators.]
Here are three common examples of identity operators that you have seen before.  
\begin{enumerate}
    \item The number ``1" is the identity operator in arithmetic.  For any number $a$, we have $1\cdot a = a$.
    
    \item We can define an identity function $g(f)$ in a domain $D=\{0<x<1\}$ as follows: $g(f(x))= f(x)$ for all $x$ in $D$.  It is not too difficult to recognize that this operator is equivalent to multiplying a function by the number 1.
    
    \item consider the following matrix multiplication (Note: linear algebra is covered in the next section).
    \begin{equation*}
       \twoform{I} = \left[ {\begin{array}{*{20}{c}}
1&0\\
0&1
\end{array}} \right]
    \end{equation*}
Now, for every possible vector ${\bf a}$ of dimension two, we have
    \begin{equation*}
       \twoform{I}\cdot {\bf a} = 
       \left[ {\begin{array}{*{20}{c}}
{{1_{}}}&0\\
0&{{1_{}}}
\end{array}} \right]\left[ {\begin{array}{*{20}{c}}
{{a_1}}\\
{{a_2}}
\end{array}} \right] = \left[ {\begin{array}{*{20}{c}}
{{a_1}}\\
{{a_2}}
\end{array}} \right]
    \end{equation*}
    Therefore, $\twoform{I}$ is an identity operator for the domain of vectors of dimension two.
\end{enumerate}
\end{example}
\end{svgraybox}

We won't be using the \emph{theorem-proof} format very often in this material, but occasionally it is helpful, especially when the proofs are short and clever.  The definition of the inverse operator can be done this way.

\begin{theorem}[Inverse operator]
 A linear operator can have an inverse, $\mathscr{L}^{-1}$, only if $\mathscr{L}(x) = 0$ implies that $x=0$.
\end{theorem}
\begin{proof} If $\mathscr{L}(x) = y$ then the inverse of $\mathscr{L}$ is the mapping which takes $y$ back to $x$. (AN ASIDE: Here, it might be helpful to think of a conventional function, such that $\mathscr{L}(x)=y(x)$, and $x$ is the domain of the horizontal axis.  For example, $\mathscr{L}(x)=x^2~for~0<x<1$ associates each independent value in the domain $x$ with a unique value for the result (which we call $y(x)$) in the range.)

Suppose now
that $\mathscr{L}(x_1) = y_0$ and $\mathscr{L}(x_2) = y_0$. Then by linearity $\mathscr{L}(x_1-x_2) = 0$. One of the following results must be true, either (a) $x_1-x_2 = 0$ (i.e., $x_1$ and $x_2$ are the same value), or (b) $x_1-x_2 \ne 0$, so that two values of $x$ in the domain of $\mathscr{L}(x)$ that would be mapped to the same value $y_0$ in the range.  
For the inverse function, the roles of the domain and range are interchanged.  Thus, for the inverse function, $y$ forms the domain, and  $\mathscr{L}^{-1}(y)=x$ is the range.  However, by definition, a function can have only one value in the range associated with a value in the domain.  Thus, the option where $x_1$ and $x_2$ are not equal is not possible (or the inverse would not be a function).  Thus we must have that $x_1=x_2$, and this means $\mathscr{L}(x)=0$ is only true for $x=0$.
\end{proof}

Without getting overly-technical, for our purposes an inverse linear operator $\mathscr{L}^{-1}$ will exist when we can show
\begin{equation}
    \mathscr{L}^{-1} (\mathscr{L}) = \mathscr{L}(\mathscr{L}^{-1})=\mathscr{I}
\end{equation}
and that $ \mathscr{L}^{-1}$ is never multiple-valued.

As usual, examples can really help make these ideas more clear.
\begin{svgraybox}
\begin{example}[Inverse operators.]
\begin{enumerate}
    \item Suppose we have $f(x) = x^2$ for $\{x: 0<x<1\}$.  Does this function have an inverse?  
    
    Try these functions
    \begin{align}
        y&=f(x)=x^2\label{bob1}\\
        x&=g(y)=\sqrt{y}\label{bob2}
    \end{align}
    Now, note that $f(g(y))=f(\sqrt{y})=y$.  But, from Eq.~\eqref{bob1}, we have $y=x^2$; thus, we  have $f(g(y))=x^2$.  Now, note that the identity operator is such that $\mathscr{I}(f(x))=\mathscr{I}(x^2)\equiv x^2$.  So, by definition, we must have that $f(g(x))=x^2 \equiv \mathscr{I}(f(x))$.  Therefore, $f(g(x))$ is the identity operator, and, by definition, $g$ is the inverse of $f$.
    
    Functions that are inverses of one another have an interesting graphical feature.  Consider the two functions above; if we plot them both as functions of $x$ (i.e., we plot the functions $y=x^2$ and $y=\sqrt{x}$, we obtain the plot below.  Functions that are inverses of one another have reflective symmetry about the line $y=x$.  
    
\begin{centering}
\includegraphics[scale=.5]{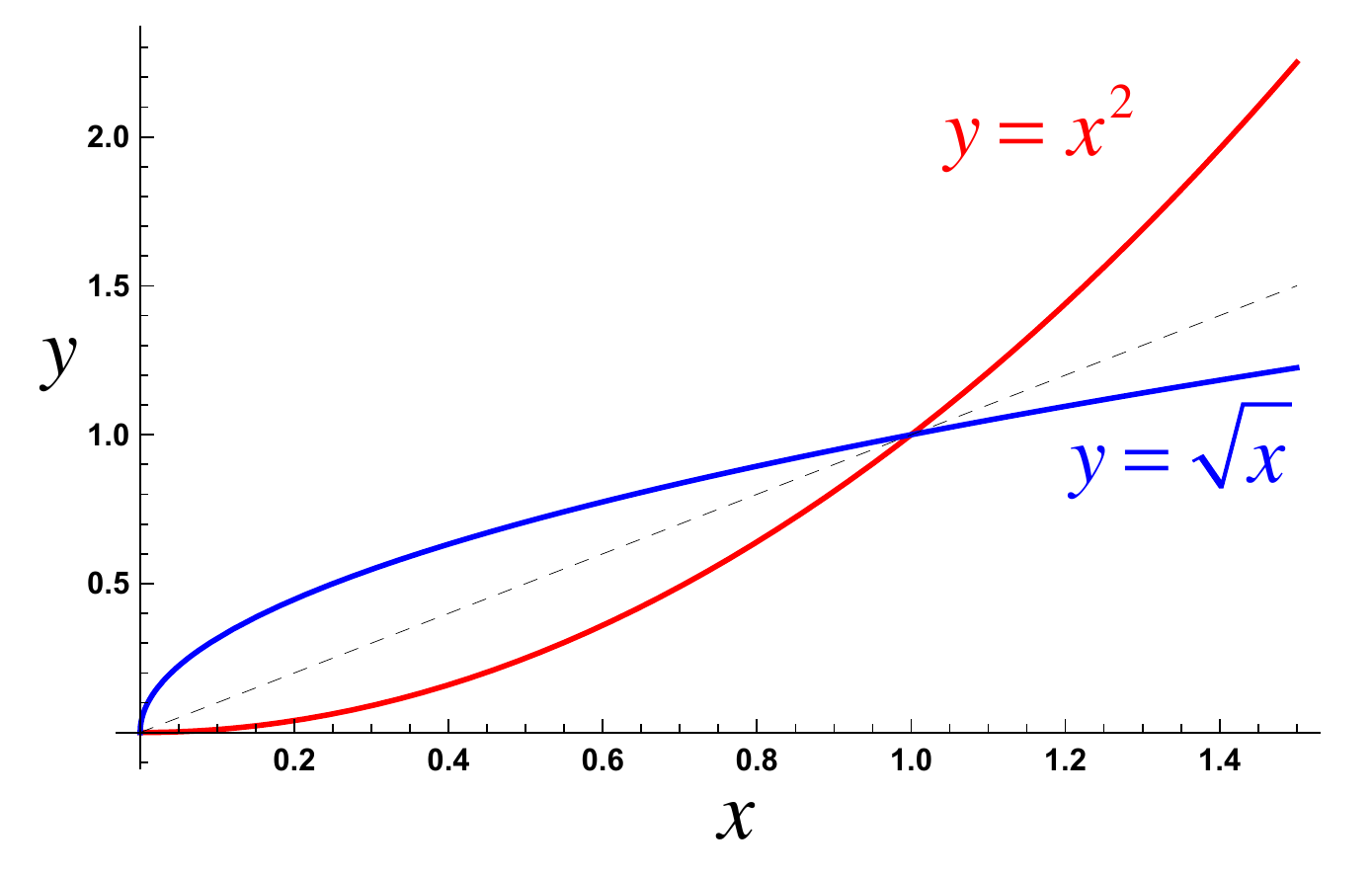}
\captionof{figure}{Inverse functions.}
\label{functions}  
\end{centering}

\end{enumerate}

\end{example}
\end{svgraybox}

There is one caveat that needs to be attended to when discussing the idea of linear operators and linear equations.  The idea of \emph{linearity} is usually associated with only the \emph{operator} part of the equation.  In particular, what this means is that constants involved in a function are excluded.  This is by convention, not necessity.  Consider the following linear equation.

\begin{align}
f(x) = 3x + 12
\end{align}
Obviously the operator $f$ must be linear, because $f(x)$ is a line!  However, if we attempt our normal process for checking linearity, we encounter a problem

\begin{align}
f(\alpha y + \beta z) &= 3(\alpha y + \beta z) + 12
&=\alpha 3 y + \beta 3 z + 12 \\
&\ne \alpha f(y) + \beta f(z)
\end{align}
So, what went wrong?  Well, the answer here is a bit tricky.  We should, technically, think of constants as not being part of the operator.  In other words, suppose we define $L(x)=3 x$.  Then, the equation above can be written

\begin{align}
f(x) = \mathscr{L}(x)+ 12
\end{align}
And the actual question can be posed as ``is the operator $\mathscr{L}$ linear"? And the answer to that question is obviously yes.  The important thing to remember here is that the \emph{linearity of the equation} has only to do with the operators involved, not with any particular constants that may be tacked on to the operators.  Thus, a linear equation is any equation whose operator (which, by definition does not include additive constants) is linear. 
%
\subsection{Linear Algebra}
%
Now that we know what a linear operator is, we can discuss linear algebra.  Most people's experience with linear algebra leaves them without an intuitive notion of what is really going on.  Because linear algebra is such a mathematically rich topic-- it is the first place where more formal mathematical analysis can sensibly be done in terms of vector spaces -- it tends to favor mathematical detail over process.  At any rate, for the material covered here, a simple review of linear algebra is sufficient.  To start, we will consider a set of linear functions; for the purposes of generating a concrete example, consider the following

\begin{align}
    f_1(x,y,z)&=3 x + 2 y + z -5 = 0 \\
    f_2(x,y,z)&=x + y -z -1  = 0\\
    f_3(x,y,z)&=2x + 2y +3z -10 =0
\end{align}
Obviously this is a set of equations, and obviously they are linear.  Do these equations have a solution?  Well, we know that there must be if the three equations are \emph{independent}.  That is to say, no equation is a linear combination of the other two; or, there are no coefficients $a,b\ne0$ such that $f_1=a f_2 + b f_3$.  Note that this one condition is sufficient (why?).  Assuming that the equations are linearly independent, and that there are as many equations as there are variables, then there is a solution to the set of equations (this is known from the \emph{fundamental theorem of linear algebra}- which we will not prove!)

One of the reasons for matrix algebra in the first place is that writing all of these equations down is somewhat repetitive.  Thus, we can compact the notation by defining matrix multiplication.  To start, note that the set of equations above can be written out more compactly (eliminating the function labels, which are unnecessary) 

\begin{align}
    3 x + 2 y + z  &= 5 \\
    x + y -z   &= 1\\
   2x + 2y +3z &=10
\end{align}
Now, suppose we define four vectors, as follows
\begin{align*}
    {\vec r}_1&= (3, 2, 1) \\
    {\vec r}_2&= (1, 1, -1) \\
    {\vec r}_3&= (2, 2, 3) \\
    {\vec x}&= (x, y, z) 
\end{align*}
Recall that the dot product between two vectors is given by ${\bf a}\cdot{\bf x}=(a_1,a_2,a_3)\cdot(x,y,z)= a_1 x+a_2 y + a_3 z$.  Noting this, we can write our equation more compactly as
\begin{align*}
    {\vec r}_1 \cdot  {\vec x} &= 5 \\
    {\vec r}_2 \cdot  {\vec x} &= 1 \\
    {\vec r}_1 \cdot  {\vec x} &= 10 
\end{align*}
Or, noting ${\vec x}=(x,y,z)$, using conventional notation for vectors
\begin{equation}
    \left[ {\begin{array}{*{20}{c}}
  {\bf r}_1 \\ 
   {\bf r}_2 \\ 
   {\bf r}_3
\end{array}} \right]\cdot
\left[ {\begin{array}{*{20}{c}}
  {x}\\ 
   {y}\\ 
   {z}
\end{array}} \right]
=
\left[ {\begin{array}{*{20}{c}}
  {5} \\ 
   {1} \\ 
   {10}
\end{array}} \right]
\end{equation}
This essentially \emph{defines} matrix multiplication by the notation

\begin{equation}
\left[ {\begin{array}{rrr}
  3&~2&1 \\ 
  1&~1&{ - 1} \\ 
  2&~2&3 
\end{array}} \right]\left[ {\begin{array}{*{20}{c}}
  x \\ 
  y \\ 
  z 
\end{array}} \right] = \left[ {\begin{array}{*{20}{c}}
  5 \\ 
  1 \\ 
  {10} 
\end{array}} \right]
\end{equation}
where to complete the multiplication, each row (taken as a vector) is dotted with the vector $\vec{x}=(x,y,z)$.  This is the easiest way to remember how to do matrix multiplication!

To solve this problem, there are a number of ways that we can proceed.  The easiest one is just to do what we would have done if we still had everything written out as three equations in three unknowns: eliminate variables from the equations simultaneously.  In short, we can do any of the following operations:
\begin{enumerate}
    \item Multiply any row by a constant.
    \item Add \emph{any} two rows, and replace either of those two by the result.
    \item Interchange any two rows.
\end{enumerate}
There is a caveat here-- whatever we do to the rows of the matrix, we need to also do to the vector on the right-hand side (they are equations after all!)
The simplest way to see this is to just do an example.  It is possible to do this in a super-orderly, algorithmic fashion, but with small matrices, it is often better to look for easy opportunities for simplifications.  To start, multiply row 2 by $-2$, and add that to row 3, the divide row 3 by $5$.   This gives

\begin{equation}
\left[ {\begin{array}{rrr}
  3&~2&1 \\ 
  1&~1&{ - 1} \\ 
  0&~0 &1 
\end{array}} \right]\left[ {\begin{array}{*{20}{c}}
  x \\ 
  y \\ 
  z 
\end{array}} \right] = \left[ {\begin{array}{*{20}{c}}
  5 \\ 
  1 \\ 
  {8/5} 
\end{array}} \right]
\end{equation}
Now, multiply row three by $-1$ and add to row 1; then add row 3 to row 2.  Now we have

\begin{equation}
\left[ {\begin{array}{rrr}
  3&~2&0 \\ 
  1&~1&{ 0} \\ 
  0&~0 &1 
\end{array}} \right]\left[ {\begin{array}{*{20}{c}}
  x \\ 
  y \\ 
  z 
\end{array}} \right] = \left[ {\begin{array}{*{20}{c}}
  17/5 \\ 
  13/5 \\ 
  {8/5} 
\end{array}} \right]
\end{equation}
Clearly, we are close.  Multipy row 2 by $-2$ and add to row 1.

\begin{equation}
\left[ {\begin{array}{rrr}
  1&~0&0 \\ 
  1&~1&{ 0} \\ 
  0&~0 &1 
\end{array}} \right]\left[ {\begin{array}{*{20}{c}}
  x \\ 
  y \\ 
  z 
\end{array}} \right] = \left[ {\begin{array}{*{20}{c}}
  -9/5 \\ 
  13/5 \\ 
  {8/5} 
\end{array}} \right]
\end{equation}
As the last step, multiply row 1 by $-1$, and add it to row 2

\begin{equation}
\left[ {\begin{array}{rrr}
  1&~0&0 \\ 
  0&~1&{ 0} \\ 
  0&~0 &1 
\end{array}} \right]\left[ {\begin{array}{*{20}{c}}
  x \\ 
  y \\ 
  z 
\end{array}} \right] = \left[ {\begin{array}{*{20}{c}}
  -9/5 \\ 
  22/5 \\ 
  {8/5} 
\end{array}} \right]
\end{equation}
The matrix on the left-hand side of this expression is called the \emph{identity} matrix, denoted $\twoform{I}$.  It is the matrix version of ``1".  Carrying out the matrix multiplication on the left-hand side leads to the solution

\begin{equation}
\left[ {\begin{array}{*{20}{c}}
  x \\ 
  y \\ 
  z 
\end{array}} \right] = \left[ {\begin{array}{*{20}{c}}
  -9/5 \\ 
  22/5 \\ 
  {8/5} 
\end{array}} \right]
\end{equation}

As a final note regarding linear algebra, there are a few more words to say about determining whether or not a set of equations is solvable or not.  Above, we mentioned that the equations needed to be linearly independent.  That is actually not an easy thing to check (even though it is easy to define).  It turns out that there is a characteristic number for a set of equations that indicates whether or not there is a solution to them.  This number is called the determinant, and it turns out to be somewhat difficult to define.  The following serves as a decent definition for the general case of an $n\times n$ matrix.

\begin{definition}
The determinant of a $n \times n$ matrix can be found by 
\begin{enumerate}
    \item Putting the matrix in upper or lower triangular form (i.e., conducting row reduction).
    \item Once in triangular form, multiplying the values on the diagonal gives the determinant.
\end{enumerate}
Here, \emph{triangular form} just means eliminating all of the entries in the matrix either above or below the diagonal.
There are a few caveats when using this approach to compute the determinant.  While you can always add a multiple of one row to another row, for other row-reduction operations more care is needed.  The following list provides guidance for the reduction process.

\begin{table}[h!]

\begin{center}
\caption{A list of operations that can be done for computing the determinant, and the associated effect of the operation.}
\label{t:determinant}
 \begin{tabular}{|c| c |c |} 
 \hline
 ~ & {\bf Type of operation} & {\bf Effect}  \\ [0.5ex] 
 \hline
 1 & ~~Add a multiple of one row to another row~~ & No effect \\ 
 \hline
 2 & Multiply a row by a constant, $c$ & ~~Determinant is multiplied by $c$~~ \\
 \hline
 3 & Interchange two rows & Determinant changes sign  \\
 \hline
\end{tabular}
\end{center}
\end{table}

\end{definition}

\begin{svgraybox}
\begin{example}[Determinants.]
Consider the matrix we just examined in the material above
\begin{equation}
    \twoform{A}=\left[ {\begin{array}{rrr}
  3&~2&1 \\ 
  1&~1&{ - 1} \\ 
  2&~2&3 
\end{array}} \right]
\end{equation}
To find its determinant, we first perform row operations (but not interchanging any two rows, and not multiplying any row by a constant) to eliminate the entries above or below the diagonal to make a triangular matrix.  A quick glance at the matrix $\twoform{A}$ indicates that either way will be reasonably easy.  One set of steps is as follows: (a) Multiply row 2 by $-2$ and adding it to row 3 (replacing row 3 with the sum), (b)  multiply row 1 by $-1/3$ and adding it to row 2 (replacing row 2 with the sum) should do it.  The result is

\begin{equation}
    \twoform{A}=\left[ {\begin{array}{rrr}
  3&~2&1 \\ 
  0&~1/3&{ - 4/3} \\ 
  0&~0&5 
\end{array}} \right]
\end{equation}
That wasn't too bad.  According to our definition, the determinant is just the product of the entries on the diagonal.
\begin{equation}
    \textrm{det}(\twoform{A})=\left|\twoform{A}\right| = 3\times 1/3 \times 5 = 5
\end{equation}
Here, two forms of notation for the determinant (``det" and the vertical bars) have been shown primarily for reference.  

\end{example}
\end{svgraybox}
%
There are several reasons that the determinant is handy.  For our purposes, it is useful because it can tell you whether or not a set of equations has any redundancies (i.e., whether or not it is linearly independent).

\begin{definition}
A square matrix $\twoform{A}$ is linearly independent if and only if it has a nonzero determinant, $\textrm{det}(A)\ne 0$.
\end{definition}
By the way, the \emph{if and only if} (sometimes \emph{iff}) statement means that the results apply both ways.  Above, for instance, it means ``if a square matrix has a nonzero determinant then it is linearly independent" \emph{and} ``If a square matrix is linearly independent, then it has a nonzero determinant."  

There is a method of solution of linear systems called \emph{Cramer's rule} that involves only computing determinants.  It is useful for small (2 by 2 or 3 by 3) matrices where the determinants are not difficult to compute. This method is discussed further in the problems.

\section{Calculus}

This review of calculus will be, like the sections before, short and focused more specifically on topics that are useful review for the material to come rather than an exhaustive summary of the subject.  The logical starting place for a review is the definition of the derivative.

\subsection{The Derivative}

\begin{definition}[The Derivative of a Function]
A derivative of a continuous (at least $C^1$, as defined above) function, $f$, at a point $x$ in the domain of the function is defined by the limit
\begin{align}
   f'(x)=\mathop{\lim }\limits_{t\rightarrow x}& \frac{f(t)-f(x)}{t-x} \\
    \intertext{or, equivalently, we can define}
     \vspace{-5mm}
     f'(x)=\mathop{\lim }\limits_{\Delta x\rightarrow 0}& \frac{f(x+\Delta x)-f(x)}{\Delta x} 
     \label{derivdef}
\end{align}
\end{definition}
There are at least two common notations for the derivative: $f'$ and $\tfrac{df}{dx}$.  For functions of a single variable, there is little chance for confusion.  Sometimes, for purposes of clarity or presentation, one of the two is preferable to the other, especially when multiple functions of different dependent variables are considered.

Note that the definition of the derivative automatically provides the definition of higher-order derivatives.  For example, consider the function $g(x)=f'(x)$.  Then we have

\begin{align*}
     g'(x)&=(f'(x))'=\mathop{\lim }\limits_{\Delta x\rightarrow 0} \frac{g(x+\Delta x)-g(x)}{\Delta x}  \\
     \intertext{Or, establishing the notation for the second derivative} \\
     f''(x)=\frac{d^2f}{dx^2}& =\mathop{\lim }\limits_{\Delta x\rightarrow 0} \frac{f'(x+\Delta x)-f'(x)}{\Delta x}\\
     \intertext{Or, applying the original definition of the derivative given by Eq.~\ref{derivdef}, we find}\\
     f''(x)=\frac{d^2f}{dx^2}&=\mathop{\lim }\limits_{\Delta x\rightarrow 0}\frac{f(x+2\Delta x)-2f(x+\Delta x)+f(x)}{(\Delta x)^2}
\end{align*}
This provides some explanation for why the second derivative is denoted by $\tfrac{d^2 f}{d x^2}$.

Generally, we do not derive the derivatives of functions from first principles, except perhaps in our introductory course on calculus.  However, it is useful to recall how this is done.

\begin{svgraybox}
\begin{example}[Derivatives.]
Computing derivatives directly from the definition of the derivative is not tremendously difficult, but sometimes it does require a little creativeness in determining the limit.  As an example, let's look at how to find the derivative of the function $f(x)=x^2$.  

\begin{align*}
    f'(x)& = \mathop{\lim }\limits_{t\rightarrow x}\frac{f(t)-f(x)}{t-x} \\
    &= \mathop{\lim }\limits_{t\rightarrow x}\frac{t^2-x^2}{t-x} \\
    &= \mathop{\lim }\limits_{t\rightarrow x}\frac{(t-x)(t+x)}{t-x}\\
    &= \mathop{\lim }\limits_{t\rightarrow x}(t+x)\\
    &= 2x
\end{align*}

\end{example}
\end{svgraybox}

Using the definition of the derivative, it is also possible to derive what the derivative of a \emph{product of two functions} is.  This is called the \emph{product rule for differentiation} or the \emph{Leibniz rule for differentiation}. 

\begin{theorem}[Product Rule for Differentiation]
\begin{equation}
 \frac{d}{dx}\left[f(x)g(x)\right]= f'(x)g(x)+f(x)g'(x)
 \end{equation}
\end{theorem}

\begin{proof}  
The proof for this is just a straightforward application of the definition of the derivative
\begin{align*}
   \frac{d}{dx}\left[ f(x) g(x)\right]&=\mathop{\lim }\limits_{t\rightarrow x} \frac{f(t)g(t)-f(x)g(x)}{t-x} \\
   &= \mathop{\lim }\limits_{t\rightarrow x} \frac{f(t)g(t)+[f(x)g(t)-f(x)g(t)]-f(x)g(x)}{t-x} \\
   &=\mathop{\lim }\limits_{t\rightarrow x} \frac{[f(t)-f(x)]g(t)}{t-x}+
  \mathop{\lim }\limits_{t\rightarrow x} \frac{[g(t)-g(x)]f(x)}{t-x} \\
  &= \frac{df(x)}{dx}g(x)+f(x)\frac{dg(x)}{dx} \tag*{\qedhere}
\end{align*}
\end{proof}

There is another rule for differentiation that is essential; this is the \emph{composition rule for derivatives}, more frequently called the \emph{chain rule for derivatives}.  This rule is handy when one has functions embedded in other functions; in other words, a \emph{composite} function.  

\begin{definition}[Composite Function]
A \emph{composite function} is a function whose independent variable is also a function.
\end{definition}
As an example, consider the relationship between position and velocity.  Suppose that you have a position function (in one dimension) such that $x=x(t)$. Then the velocity, which depends on position and time, $v(x(t),t)$, is a composite function.

The chain rule provides a method of computing the derivative of composite functions.  The proof of the chain rule is pretty complex, so it will not be presented here, but it is available in nearly every introductory text on calculus. The result, however, is important and will be used frequently in the material that follows.

\begin{theorem}[Chain Rule for Differentiation]
Suppose we have a function, $f$ whose argument is another function, $g$.  Assume that $g(t)$ has a derivative in the set of points $T=\{t: a<x<b\}$, and that  $f(y)$ has a derivative in the set of points $Y=\{g(t)\}$.  The composite function $f(g(t))$ is differentiable, and its derivative is given by

\begin{equation*}
    \frac{d}{dt}f(g(t)) = \frac{d f(g)}{dg}\frac{d g(t)}{dt}
\end{equation*}
\end{theorem}
Keeping the functions $g$ and $f$ straight sometimes causes confusion.  Some examples are helpful.
\begin{svgraybox}
\begin{example}[Product Rule, Velocity-Position Example.]

Above, we mentioned the relationship between position and velocity as an example that could be considered a composite function problem.  This is one example where one can compute derivatives both with and without the use of the chain rule.  In essence, it provides a way to validate that the chain rule leads to correct results.

Suppose you go out to run along the street, but in a very strange way (why you decide to do this is open to discussion). As you run along, you run in a way such that your position is the cube of the time you have been running (obviously, you can't keep this up forever, but for a short time it is possible.)  In particular, suppose your position, $x(t)$, is given by $x(t) = \tfrac{1}{3}t^3$.  Now, we know that velocity is defined as the time rate of change of position, and acceleration is the time rate of change of velocity, so 

\begin{align*}
    v(t)&=\frac{dx(t)}{dt} &a(t)&=\frac{dv(t)}{dt}&&(A) \\
    &= t^2&&= 2 t 
\end{align*}
From here, we are going to do something that seems a bit unusual, but it is necessary so that we can validate the chain rule.  First, note that we can express $x$ as a function of $v$, as follows (noting $v=t^2$ implied  $t=v^{\tfrac{1}{2}}$)

\begin{align*}
    x(v(t))&=\frac{1}{3}t^3= \frac{1}{3}v^{\tfrac{3}{2}}&&(B)
\end{align*}
Now, we have expressed $x$ as a composite function, $x(v(t))$.  We can compute the time derivative of this function by

\begin{align*}
    \frac{d x(v(t))}{dt}&=\frac{d x}{d v} \frac{d v}{dt}&&(C)
\end{align*}
We can compute $dx/dv$ by
\begin{equation*}
    \frac{d x}{dv} = \frac{d }{dv}\left(\frac{1}{3}v^{\tfrac{3}{2}}\right) = \frac{1}{2}v^{\tfrac{1}{2}}
\end{equation*}
and $dv/dt$ has already been computed above in Eq.~(A).  Combining these, we find
\begin{align*}
    \frac{d x(v(t))}{dt}&=\frac{d x}{d v} \frac{d v}{dt} \\
    &= \left(\frac{1}{2}v^{\tfrac{1}{2}}\right)\left( 2 t\right)\\
    \intertext{and substituting $v=t^2$ gives the result}
    \frac{d x(v(t))}{dt}&= t^2
    \intertext{which is identical to what is given in Eq.~(A) above.}
\end{align*}
\end{example}
\end{svgraybox}

\begin{svgraybox}
\begin{example}[Product Rule.]
Some kinds of catalyst can be deactivated by chemicals produced during the catalytic reaction, or by external factors such as UV light.  This process is sometimes called \emph{catalytic poisoning}.  Similar kinds of deactivation can happen to the enzymes in cell systems in biological reactors.  

Suppose that an experiment is run, and it is determined that the amount of product begin produced, under conditions of deactivation (whether cells or catalyst), is given by

\[ c(t) = c_0 \exp\left[-(k_0-k_1 t) t \right] \]
In other words, an first-order-like rate process has an effective rate constant that is a function of time, (that is, $k_{effective}=k_0=k_1 t$).  Determine the rate of reaction, $c'(t)$.\\

\noindent{\bf Solution.} 

Using the product rule requires that we first identify the composite functions.  Often, in practice, this is not done explicitly; rather, people just keep mental note of which function is which.  However, it is instructive to explicitly identify the functions when there is any potential for confusion.  For this problem, take

\begin{align*}
    g(t) &= -(k_0-k_1 t) t  \\
    f(g)&= c_0 \exp\left( g \right)
\end{align*}
We need to recall the definition for the derivative of the exponential 

\begin{equation*}
    \frac{d}{dg}exp(g)= exp(g)
\end{equation*}
The exponential is the only function (except the function 0) whose derivative is the same as the starting function!  With this, we have all we need.

\begin{align*}
    \frac{d g(t)}{dt} &= -(k_0-\frac{k_1 t}{2} ) \\
    \frac{d f(g)}{dg}&= \exp\left(g \right)\\
    \intertext{So, the result is \ldots} 
    \frac{d}{dt}f(g(t))&=\frac{d f(g)}{dg}\frac{d g(t)}{dt}\\
    &= c_0 \exp\left(g \right)\left[-(k_0-\frac{k_1 t}{2} )\right]
    \intertext{Substituting the function $g(x)$ (from above) into this result and rearranging gives the final result}
     \frac{d}{dt}f(g(t))&= -c_0(k_0-\frac{k_1 t}{2} )\exp\left[-(k_0-k_1 t) t \right] 
\end{align*}

\end{example}
\end{svgraybox}

\subsection{Partial Derivatives}
There is not much more to say regarding derivatives, except to explain the notion of a derivative when a function has multiple independent variables.  So, with out delay, we will define the partial derivative of a function with two independent variables.  The case of additional independent variables is identical, so no more than two is required for the definition.

\begin{definition}[Partial Derivatives]
For a function with two independent variables, a partial derivative is the derivative of the function with respect to only one of the two variables (the other variable being held constant.)  Assume that $f(x,t)$ a continuous (at least $C^1$) function of of $x$ and $t$ over a domain $\Omega$ ($\Omega$ could be an irregularly-shaped domain, so we will skip an effort to provide a more detailed set description of it).  The partial derivative of $f$ with respect to each variable is given by

\begin{align*}
   \frac{\partial f(x,t)}{\partial x}=\mathop{\lim }\limits_{\Delta x\rightarrow 0}& \frac{f(x+\Delta x,t)-f(x,t)}{\Delta x} \\
    \frac{\partial f(x,t)}{\partial t}=\mathop{\lim }\limits_{\Delta t\rightarrow 0}& \frac{f(x,t+\Delta t)-f(x,t)}{\Delta t} 
\end{align*}
\end{definition}
Similar results hold for functions of three or more variables, and the extension should be reasonably transparent based on the examples above.

\begin{definition}[The Chain Rule for Functions of Two or Three Variables]  Recall that a \emph{composite function} is a function whose independent variable is also a function.  When a function is dependent upon two or more variables that are themselves functions, the Chain Rule allows us to determine the derivative in the following form

\begin{align*}
    \frac{d }{d t} f(x(t), y(t)) = \frac{\partial f}{\partial x}\frac{d x}{d t} + \frac{\partial f}{\partial y}\frac{d y}{d t}
\end{align*}
Note that in this definition, the derivatives of $x(t)$ and $y(t)$ are not \emph{partial} derivatives, but conventional derivatives.  This is because $x(t)$ and $y(t)$ are \emph{functions of a single variable, $t$}.  Therefore, the conventional derivative is the correct form of the derivative for those quantities.
\end{definition}

\subsection{Integration}\label{integration}\indexme{Integration}

Integration is a simple concept, but it turns out to be quite deep in actual applications.  In introductory calculus, the concept of the Riemann integral is introduced.  This integral applies to most functions that are encountered in science and engineering; in particular, it is useful for evaluating functions that are \emph{piecewise smooth} as defined previously.  This is not the only kind of integral that can be defined.  For example, consider the following function (the Dirac function) 

\begin{equation*}
 f(x) =
\left\{ 
   {\begin{array}{{lll}}
  {1~\mbox{if x is irrational on} x\in[0,1] } \\ 
  {0~\mbox{if x is rational on } x\in[0,1] } 
\end{array}}\right.
\end{equation*}
This function is not piecewise smooth because it has an infinite number of holes in it (i.e., it has a hole at every possible fraction between 0 and 1!) The conventional Riemann integral cannot be used for such a function.  However, there are more general forms of the integral (such as the Lebesgue integral) that can be used to measure such functions.  We will not explore the Lebesgue integral in this text, but we will introduce (broadly) the ideas behind it at the end of this section on integration.  

For our purposes, the development of the integral will not be reviewed in detail.  Instead, a few important properties of the integral are presented.

\subsubsection{Riemann Sums, the Integral, and the Differential}\label{riemann}\indexme{integral!Riemann}

While a thorough treatment of integration theory is not needed here, some reminders about the definition of the conventional integral of continuous functions is useful.  The primary purpose here is to review the basic idea of the definition of the integral rather than to generate the most general notion possible. Therefore, it suffices for now to consider continuous functions (although functions with any finite number of discontinuities are also covered by this definition).  Suppose we want to compute the integral of a function, $f$ of a single independent variable, $x$, over some interval $x\in[a,b]$.  Recall, this corresponds to the area under the curve between $a$ and $b$.  If we were given a such a curve and asked to compute the area graphically, we might be tempted to estimate the area by constructing a sequence of rectangles approximating the curve (see Fig.) .  Suppose we do so, and we make each such rectangle have the same width on the $x-$axis.  Let $X=\{x_1, x_2, x_3 \ldots x_N\}$ such that $x_{i+1}-x_i = \Delta x$ (a constant) for all $i=1$ to $N-1$; hence, $\Delta x = (b-a)/(N-1)$.  The following sum is an approximation to the integral of $f$

\begin{equation}
    f_{\Delta x} = \sum_{i=1}^{i=(N-1)} f(x^*_i) \Delta x
\end{equation}
where $x^*_i$ is any value of $x$ such that $x_i \le x^*_i \le x_{i+1}$.  In this computation, then, $f(x^*_i)$ is an estimate of the height of the rectangle between $x_i$ and $x_{i+1}$.  Hence, the sum defining $f_{\Delta x}$ is an estimate of the area under the curve.  As $\Delta x$ is taken increasingly smaller, the sum naturally becomes increasingly accurate because a constant height becomes a better representation of the function in each rectangle.  The Riemann sum is then \emph{defined} by  

\begin{equation}
    \int_a^b f(x) \, dx = \mathop {\lim }\limits_{N \to \infty }\left( \sum_{i=1}^{i=(N-1)} f(x^*_i) \Delta x \right)
\end{equation}
%
\begin{figure}[t]
\sidecaption[t]
\centering
\includegraphics[scale=.33]{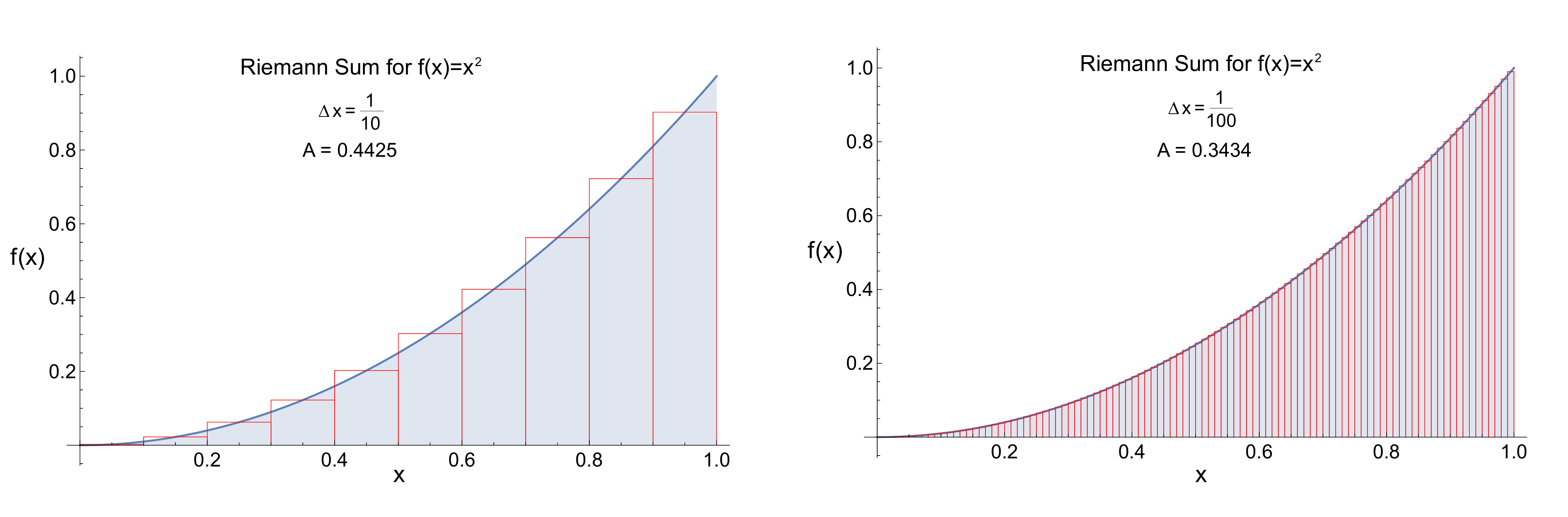}
\caption{Two examples of Riemann sums.  For these example, $x^*_i$ was taken at the midpoint of each rectangle.  On the left, a divisions of the domain with $N=10$, gives a total of 0.4424; this represents a 32\% over-estimate of the actual value. On the right  divisions of the domain with $N=100$, gives a total of 0.4424; this represents only a 3\% over-estimate of the actual value. }
\label{riemann0}       
\end{figure}
%
\noindent An example of the process of refining a Riemann sum is given in Fig.~\ref{riemann0}

A few notes are worth pointing out here.  First, note that the integral sign actually comes from an elongated ``S" as a reminder that the integral arises from a sum.  Second, the term $dx$ in the integral has a specific interpretable meaning.  It is commonly called a \emph{differential}, but the concept of a differential is not as straightforward as it might appear on the surface.  In fact, for many years even after the invention of calculus, the concept of the differential was not all that clear.  While there are many ways to make the concept of a differential formally concrete (it can be defined as a linear map from the real numbers to the real numbers \citet[][Chp.~7]{protter2012first}), for our purposes it is sufficient to consider it as follows.

A differential is a non-zero approximation to the change in a function with a change in the dependent variable being as small as needed to make the approximation attain the accuracy desired.  To be more formal, for a function $g(x)$, the differential is given by 

\begin{equation}
    dg = g'(x) dx
\end{equation}
where $dg$ stands for a ``small" change in the function $g$, and $dx$ stands for a small change in the independent variable.  In this case, the meaning of ``small" is a qualified one; it means small enough such that the error involved is less than some specified error.  Importantly, the concept is that the error can be driven to zero as $dx$ approaches zero.  While we have not yet reviewed the concept of the Taylor series (this appears in the material following), one can use a Taylor series to define the concept.  Recall

\begin{equation}
    g(x+\Delta x) = g(x) + \Delta x g'(x) +O[(\Delta x)^2] +\ldots
\end{equation}
Supposing that $\Delta x \ll 1$, then $(\Delta x)^2 \ll (Delta x)$, and the terms of $(\Delta x)^2$ and higher can be dropped relative to those involving $\Delta x$.  Rearranging we have

\begin{equation}
    g(x+\Delta x)-g(x) + =  \Delta x g'(x)
\end{equation}
Now, we note that the symbol $dx$ is used to indicate $\Delta x$ under the conditions that $\Delta x$ can be made arbitrarily small (this is, in essence, its definition); similarly, $dg = g(x+\Delta x)-g(x)$.  This gives us

\begin{equation}
    dg = g'(x) \Delta x
\end{equation}
Note that if we take the function $g(x) = x$ we end up with the relation $dx = Delta x$.  This explains the use of $dx$ in the representation of the integral.  By convention, the equality $dx = Delta x$ means that definition of the differential is given by 

\begin{equation}
    dg = g'(x) dx
\end{equation}

This is more than simply a formal manipulation of symbols.  Once we have defined differentials to be quantities whose error can be made as small as we like, then the resulting structure is essentially a linear one.  Some powerful methods in the analysis of, for example, non-Euclidian geometry.  They also arise in the study of differential equations.

\subsubsection{The Fundamental Theorem of Calculus}

The fundamental theorem of calculus says some really important things.  Primarily, it tells us the following.  

\begin{theorem}
Suppose a smooth function, $F$, is defined on the interval $[a,b]$.  Because the function is smooth, it has a derivative, $F'=f$.  Then
\begin{align*}
    \int\limits_{t=a}^{t=b} f(t) dt &= F(b)-F(a) \\
    \intertext{or, equivalently,}
        \int\limits_{t=a}^{t=b} \frac{dF(t)}{dt} dt &= F(b)-F(a) 
\end{align*}
\end{theorem}
\noindent The extension of this theorem to \emph{piecewise continuous} functions is straightforward.  It involves simply computing the integral over each of the (finite number) continuous intervals.\\

This is a very powerful theorem, and it essentially maps the problem of finding integrals on to the problem of finding derivatives.  That is to say, if we are given a function, $f$, and we happen to know a function $F$ whose derivative is equal to $f$, then we can compute the integral of $f$ with that knowledge.  The function $F$ is sometimes called the \emph{antiderivative} of $f$ for that reason.  Although this sounds somewhat circular, it is not.  Most of the ``known" integrals that exist do so because we have identified the antiderivative for the function.  

For many, a first course in calculus involves learning many ``techniques" to find the antiderivative of functions.  Most of this we will leave in the past, with the idea that we will all remember (or be able to look up) most of the fundamental integrals and derivatives that we encounter.  However, there is one ``technique" that is very useful in a number of applications, and it is one that we will have the opportunity to employ several times.  This is the rule for integration by parts.

\begin{theorem}
Let $f$ and $g$ be smooth functions on an interval $x\in[a,b]$.  Then
\begin{align*}
    \int\limits_{x=a}^{x=b} f(x) \frac{d g(x)}{d x} dx =
    f(x)g(x)\Bigg|_{x=a}^{x=b}\,\,-\,\, \int\limits_{x=a}^{x=b} \frac{d f(x)}{d x} g(x) dx\\
\end{align*}
\end{theorem}
\noindent Frequently, this rule is written in the easy-to-remember form
\begin{align*}
    \int\limits_{x=a}^{x=b} u dv =
     u v\Bigg|_{x=a}^{x=b}\,\,-\,\, \int\limits_{x=a}^{x=b} v du
\end{align*}

Integration by parts is particularly useful under the following circumstances: (i) there is an integrand of the form of a n{th}-order polynomial times some function that we know how to ingegrate $n$ times, or (ii) there is an integrand that contains an {th}-order derivative times some function that we know how to differentiate $n$ times, and we would like to eliminate the derivative.  This will be made clearer in the following example.

\begin{svgraybox}
\begin{example}[Integration by parts.]
Integrate the following functions over the interval $x\in[0,1]$
\begin{align*}
    (a)&& f(x)& = x \sin (\pi x) \\
    (b)&& g(x)&= \frac{d f}{dx} \sin (\pi x) \mbox{~~~~given $f(0)=0$}
\end{align*}
\noindent{\bf Solution.}
For (a), let $u = x$, $dv = \sin (\pi x) dx$.  Then, $du = dx$ and $v=-\pi \cos (\pi x)$ (noting $d/dx (-\cos x) = \sin x$).  Then

\begin{align*}
    \int\limits_0^1 x \sin (\pi x) \, dx &= \frac{x}{\pi} \cos (\pi x)\Big|_0^1 - \int\limits_0^1 -\cos (\pi x) dx \\
    &=\frac{1}{\pi}+\sin(\pi x)\Big|_0^1 \\
    &=\frac{1}{\pi}
\end{align*}

For (b), let $u = \sin (\pi x)$, $dv = \frac{d f}{dx} dx$.  Then, $du = 1/\pi cos(\pi x)$ and $v=f(x)$ (fundamental theorem of calculus).  Then
\begin{align*}
    \int\limits_0^1 \frac{df}{dx} \sin (\pi x) \, dx 
    &= \frac{x}{\pi} \cos (\pi x)f(x) \Big|_0^1 - \int\limits_0^1 f(x)\frac{1}{\pi} \cos(\pi x) dx \\
    &=- \int\limits_0^1 f(x)\frac{1}{\pi} \cos(\pi x) dx
\end{align*}
\noindent  In this latter example, we can go no further than this without knowing more about $f$; however, we have \emph{eliminated} the derivative, which is frequently a useful operation.

\end{example}
\end{svgraybox}

\subsubsection{Lebesque Integration }\label{lebesgue}\indexme{integral!Lebesgue}\indexme{Lebesgue integral}

In the material above, we presented an example of the Dirac function that was not integrable using the conventional (Riemann) integral.  Here, consider the following variation of that function

\begin{equation*}
 f(x) =
\left\{ 
   {\begin{array}{{lll}}
  {x~\mbox{if x is irrational on} x\in[0,1] } \\ 
  {0~\mbox{if x is rational on } x\in[0,1] } 
\end{array}}\right.
\end{equation*}

This is a function that maps all possible values of the independent variable $x$ on the interval $X=\{x: x\in [0,1]\}$ to some real number (in this case, this real number is also on the interval $f(x)\in [0,1]$).  In the mathematical literature, this function (mapping) is sometimes written as ``$f:X\rightarrow \mathbb{R}$", which just means (in written English) ``each element of $X$ is assigned exactly one value in the real numbers".  

\begin{figure}[t]
\sidecaption[t]
\centering
\includegraphics[scale=.5]{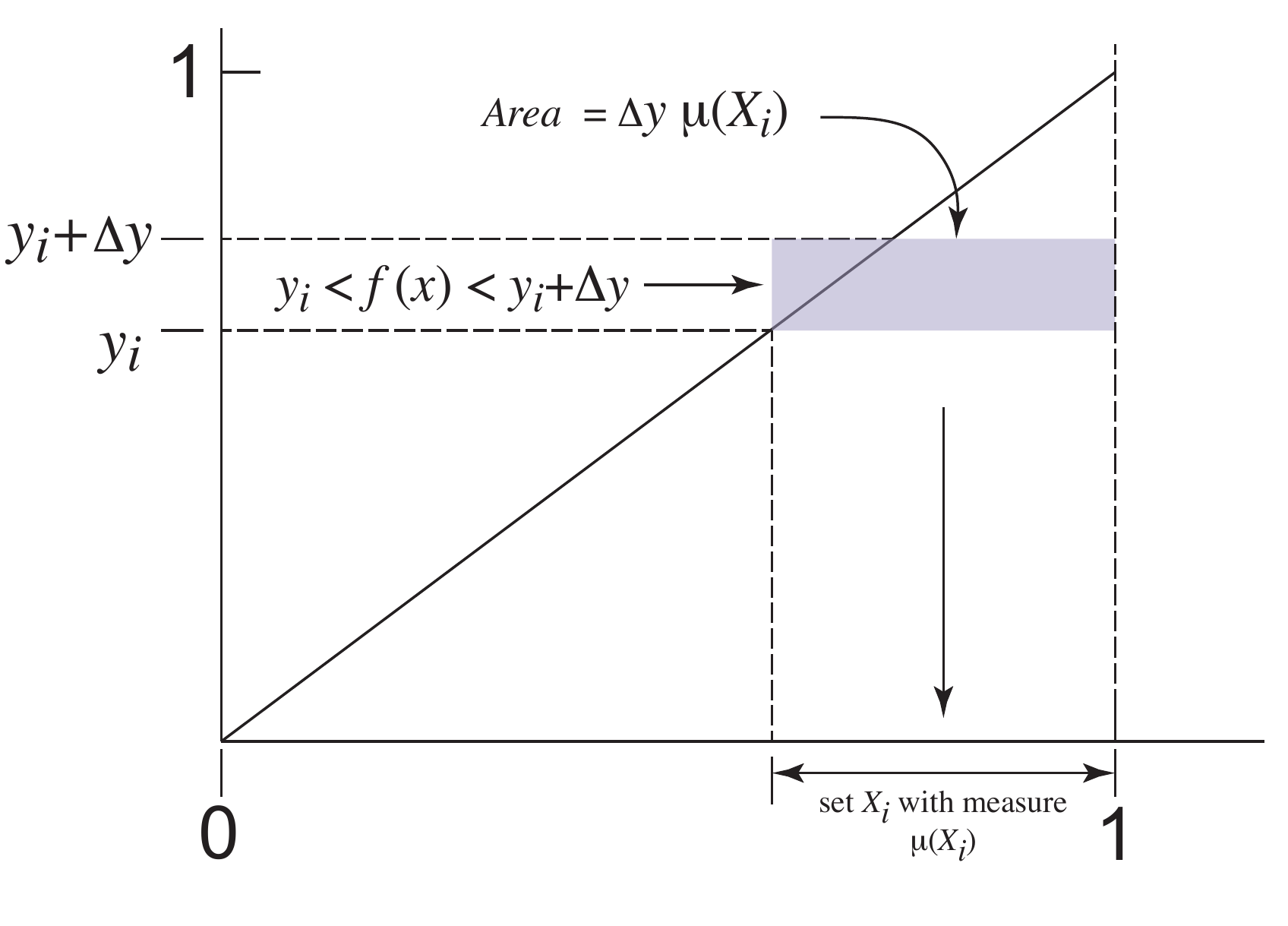}
\caption{The division of the range for a Lebesgue integral.}
\label{fig:lebesgue}       
\end{figure}

In the early 1900s, a French mathematician named Henri Lebesgue generalized the integral by thinking about it in a slightly more abstract way.  Instead of considering a ``nice" function in which we partition up the domain, and then consider various kinds of limiting operations that allows us to assign a value to the integration (which corresponds to ``area" under a curve for strictly positive functions), he considered partitioning the range instead.  In Fig.~\ref{fig:lebesgue}, and example of this process is provided for a function $f(x)$ where the range and domain are both on the intervals $[0,1]$.  Now, instead of vertical bricks being summed up over an interval, we have ``horizontal" bricks being summed up over the appropriate values of the range that correspond to $x\in[0,1]$.  

If we can assign a sensible \emph{measure} to the interval $x_i$, then we can compute the area of each such brick; the sum of these is, after the appropriate limiting process, the integral.  The problem now is to interpret what is meant by assigning a measure to the interval $x_i$.  For Riemann integrable functions, this measure corresponds to the conventional one that we think about for Riemann sums. 

\emph{Measure Theory} is the discipline within mathematics that deals with the question of ``how do we assign metrics to mathematical quantities?"  We are all familiar with the conventional Euclidian metric in $\mathbb{R}^3$-- it is just the \emph{length} of what we think of as vectors.  So, for a point $(x,y,z)$  with origin $(0,0,0)$, the Euclidian metric is

\begin{equation}
    d(x,y,z) = \sqrt{x^2+y^2+z^2}
\end{equation}
%

\begin{figure}[t]
\sidecaption[t]
\centering
\includegraphics[scale=.45]{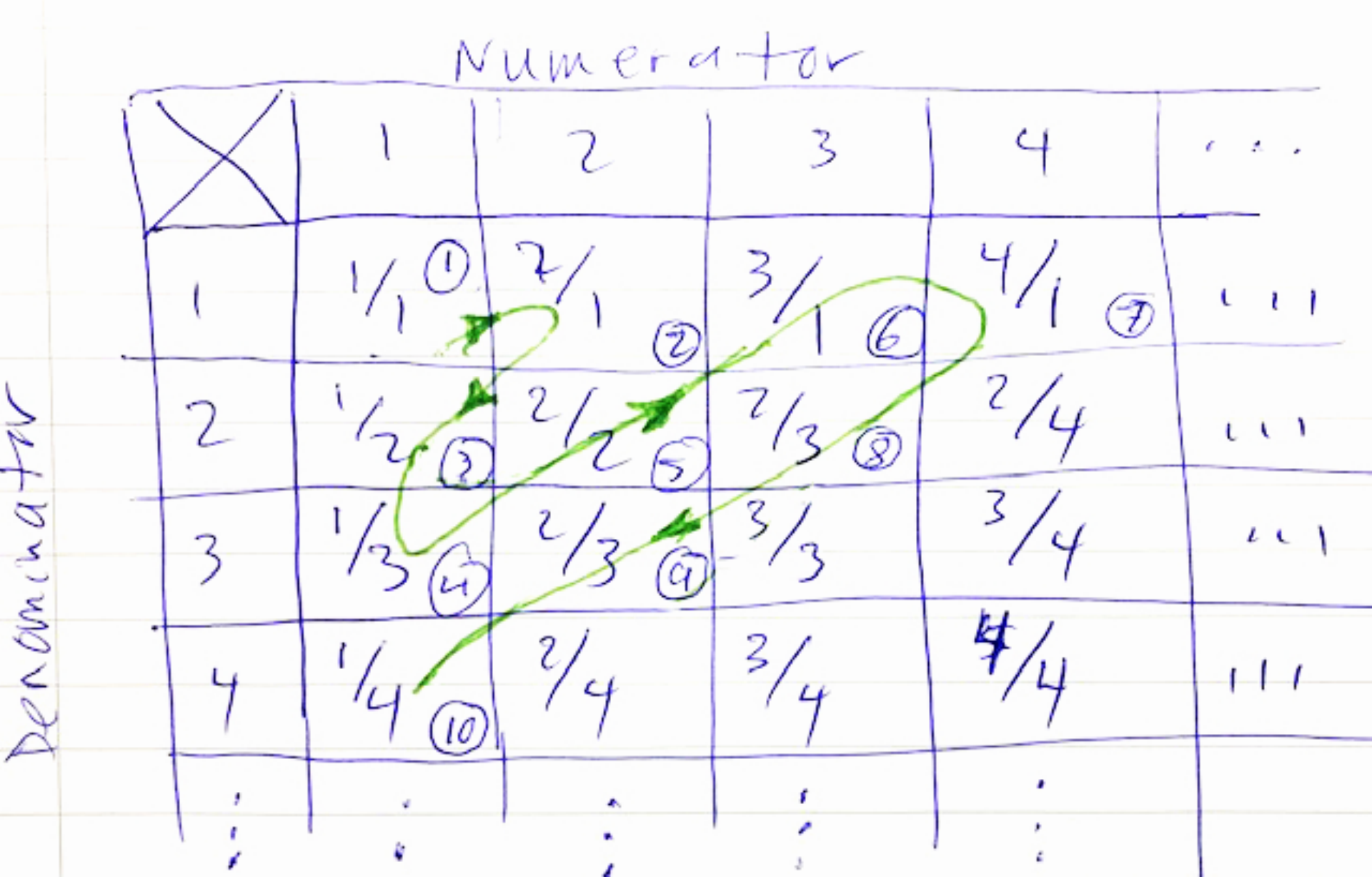}
\caption{Enumerating the rational numbers using the integers.}
\label{fig:rational}       
\end{figure}

For more exotic sets (such as the example of the rational numbers on the real line), a number of clever methods have been developed over time to assign a rational notion of ``measure" to such sets.  In particular, the measure that is used by the Lebesgue integral assigns a zero measure to isolated points (i.e., a single point on the number line is assigned zero distance, somewhat in accordance with our intuition).  It turns out that the rational numbers are much less \emph{dense} than the irrational numbers-- that is, in a sense there are many more irrational numbers than there are rational ones.  Thus, by the Lebesgue measure, the rational numbers look like isolated points that have a zero measure.  While I have summarily stated that such a notion exists, actually illustrating how to do this constructively is well beyond what we can accomplish here.  However, such considerations are not devoid of intuition. The the rational numbers can be put into a 1-to-1 correspondence with the real numbers; if you have never seen this construction, Fig.~\ref{fig:rational} contains the essence of the argument.  So, while the \emph{rational} numbers are denumerable by the integers, the irrational numbers are not.  And, these are the only two choices-- either a number is rational, or it is irrational.  So, the irrational numbers are much more dense than the rational ones.  Hence, we can legitimately think of the rational numbers as permeating the set of all possible real numbers with isolated holes (points).  Because the points have no measure, then the measure of any real interval, say $[a,b]$ is just $\mu=(a-b)$.  Again, it may appear as if we have accomplished nothing here, but this is not true!  The conventional Riemann sum for integrals does not exist for a set like the irrational numbers on $[0,1]$-- there is legitimately no way to even consider making intervals in such a case (What if one of the the interval end points lands on a rational number?  What happens as  you let $\Delta x$ tend toward zero-- it must pass though rational numbers!  There are many such problems...)  However, with the Lebesgue method of integration, this is no longer a problem.

%
For the purposes of this text, this is all we need to discuss about the Lebesgue integral (although there is much more that can be said about it!)  It is worth noting, however, that you may encounter this kind of integration in the mathematical literature.  If you do, you have most of the idea behind what such integrals actually represent, and you may still be able to read through the material.  As a final note, in the mathematics literature, one often encounters the concept of the integral of the square of a function.  This is usually called the $L_2$ metric, and it is understood that it is measured relative to the Lebesgue measure (hence the ``L" in the symbol for the metric).

\section{Sequences and Series}\label{seqseries}

Infinite sequences and series are incredibly useful, and they are the true workhorses of applied mathematics.  Many important functions (especially many of the transcendental ones that we use routinely, such as $\sin(x)$ and $\exp(x)$) are frequently \emph{defined} and/or \emph{computed} by use of an infinite series.  The coverage of sequences and series here is necessarily limited.  However, in the chapters following, we will have ample opportunity to revisit these concepts.

\subsection{Sequences}

Sequences are related to series.  One way to think of a sequence, is that it is just a indexed list of numbers or functional expressions; thus, sequences map the integers (the domain) to a another set (which may be constants, functions of the integer, or functions of the integer plus additional variables).  Thus, summing each term in an infinite sequence together would be one way to form an infinite series.  Although in introductory calculus courses, infinite sequences are often treated as functions of only the integers, they can also be functions of other independent variables unto themselves.  We will make use of the sums of such sequences in both the investigation of Power and Taylor series, and in the study of Fourier series appearing in later chapters.

To make the notion of sequences concrete, we have the following definition.


\begin{definition}[Finite Sequence]\indexme{sequence!finite}
A \emph{finite sequence} is a list of objects (elements) indexed by a consecutive subset of the natural numbers ($\mathbb{N}=\{1,2,3,4\ldots N\}$ or $\mathbb{N}_0=\{0,1,2,3,4\ldots N\})$.  In this text, a sequence of $N$ (or $N+1$ if the first index is 0) terms is denoted $A=(a_1,a_2,a_3,\ldots,a_N)$ (or $A=(a_0,a_1,a_2,\ldots,a_N)$ if the first index is 0).    The element in the position $i$ (where $i\in \mathbb{N}~or~\mathbb{N}_0$) is called the $i^{th}$ \emph{term} of the sequence.  Regardless of how the elements are defined, sequences are viewed as being \emph{functions over the natural numbers}.
\end{definition}
Note-- the primary feature of sequences is that they are a \emph{list} indexed by the real numbers.  One can think of this as a column of values on a spreadsheet program, where each entry has a unique number.  Also like a spreadsheet program, the objects in the list may be functions of the index (its numbered location), or of both the index and additional independent variables (representing, for example, spatial coordinate or time). 

\begin{definition}[Infinite Sequence]\indexme{sequence!infinite}
A \emph{sequence} is a list of objects (elements) indexed by the natural numbers, ($\mathbb{N}=\{1,2,3,4\ldots\}$ or $\mathbb{N}_0=\{0,1,2,3,4\ldots\}$).  This list may be specified by a functional rule, or by listing the elements explicitly.  A common notation for an infinite sequence is $A=(a_n)_{n=1}^\infty$ (or $A=(a_n)_{n=0}^\infty$ if the first index is zero).\\
\end{definition}
All of the notes regarding finite sequences apply to infinite sequences as well.  In particular, infinite sequences of functions that converge in a particular way to a specified function are relatively well-used concepts in applied analysis.  We will discuss this concept additionally in the material following.

\begin{figure}[t]
\sidecaption[t]
\centering
\includegraphics[scale=.33]{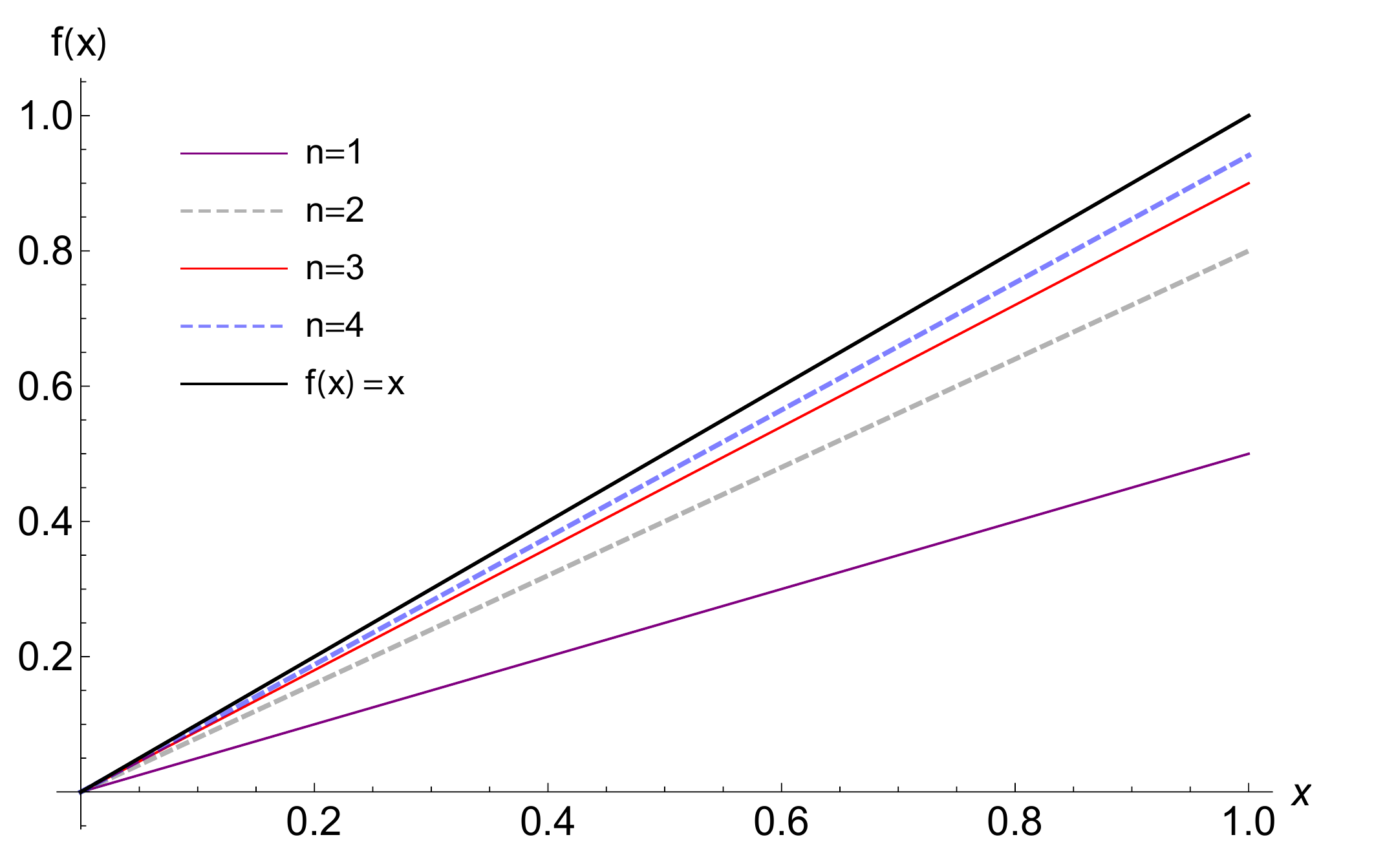}
\caption{A sequence of whole \emph{functions} of $x$.  In this case, the sequence itself is still indexed by the real numbers, $n$.  The sequence is given by $a_n(x) = $ }
\label{funcseq}       
\end{figure}

\begin{svgraybox}
\begin{example}[Examples of sequences]
Sequences are just lists of objects.  Thus, they represent a fairly general mathematical concept.  Here are a few examples.

\begin{enumerate}
    \item The following is a somewhat boring finite sequence: $A=(1,1,1,1,3,1)$.  It has a finite number of terms (six), and each element is the same, except for the fifth term.
    \item The following is an infinite sequence containing all of the even numbers greater than zero: $A=(2,4,6,8,\ldots)$.
    \item Here are two other ways to denote the very same sequence of the previous example: $A=(2n)_{n=1}^{\infty}$ and $S=(2n)_{n\in\mathbb{N}}$.
    \item Here is an infinite sequence that denotes a familiar irrational number: $A=(3,1,4,1,5,9,2,6,5,3,5,\ldots)$
    \item Here is a sequence whose terms converge to 1 as $n\rightarrow\infty$:
    $A=(1+1/n^2)_{n=1}^\infty =$
    $\left( 2,\tfrac{5}{4},\tfrac{10}{9},\right.$$\left.\tfrac{17}{16},\tfrac{26}{25},\ldots\right)$.  We can compute the following limit for this sequence:\\
    $\mathop {\lim }\limits_{n \to \infty } ({a_n})=1+\mathop {\lim }\limits_{n \to \infty }(\tfrac{1}{n^2}) = 1$.
    \item The following example is one which is a function of both the sets of integers (which form the sequence), and another independent variable (which, for concreteness, we can consider space)
    \begin{equation*}
        A= (1/(1+1/n^2) x)_{n=1}^{\infty}
    \end{equation*}
    Thus, this sequence, when computed term-by term, is given by $A=(x/2,x/(1+1/4)x, x/(1+1/9), x/(1+1/16),\ldots) $. A plot of the first four sequences of the function ($n=1,2,3$ and $4$) appear in Fig.~\ref{funcseq}.  A little though will indicate that this sequence eventually converges to the function $f(x)=x$.\indexme{sequence!sequence of functions}
\end{enumerate}
\end{example}
\vspace{-2mm}
\end{svgraybox}

In the examples above, we introduced the idea of \emph{convergence} for infinite sequences.  To be concrete, we need to establish what it means for a sequence to converge.  The issue that we need to consider is how to show that an infinite sequence actually approaches some limit that we define.  This is a somewhat tricky prospect.  Suppose we have a sequence denoted by $A=(a_n)_{n\in\mathbb{N}}$.  If you think about it,  there is no actual number $a_\infty$.  All we can really ask about an infinite sequence is the question ``what happens as $n$ becomes an arbitrarily large integer?"  This is exactly the notion behind a limit.  While we have discussed limits approaching a finite value in the material previous to this section, we have not yet discussed infinite limits.  

A \emph{mathematical} discussion of limits inevitably an $\epsilon$-$\delta$-like argument.  In this case, in stead of having a $\delta$ interval in which some error measure must be small, we instead have an integer for which all values of the sequence indexed by this integer or higher are sufficiently close to a limit.  To make this clear, here is the definition of the limit of a sequence.  Here, for clarity, we focus on sequences that are a function of only $n$ (unlike the final case described in the previous Example); however, extension of these concepts can be done on a pointwise basis for sequences of functions.  

\begin{definition}[Limits of a sequence: Sequence convergence]
Suppose $(a_n)_{n=1}^{\infty}$ is a sequence.  Then, the limit, $L$, exists if the sequence gets arbitrarily close to $L$ as $n$ increases.  In other words, for every $\epsilon >0$ that can be chosen, no matter how small the value of $\epsilon$, then it is also true that there is an integer, $N$, such that
\begin{equation}
    if~n\ge N, ~~then~~ |a_n-L|<\epsilon 
\end{equation}
In that case, we say that $(a_n)_{n=1}^{\infty}$ has a limit, $L$, and write

\begin{equation}
 \mathop {\lim }\limits_{n \to \infty } {a_n} = L
\end{equation}
\end{definition}
\indexme{sequence!convergence} \indexme{convergence!sequence}

If a sequence is not converging, does that imply that it must be diverging?  The answer here is no as evidenced by the counterexample $(a_n)_{n=1}^{\infty}=(-1)^n, n=1,2,3\ldots$.  This sequence neither converges nor diverges; it simply oscillates between $-1$ and $1$ periodically.   A definition for a diverging sequence can be made, and it is stated as follows.

\begin{definition}[Diverging sequences]
Suppose $(a_n)_{n=1}^{\infty}$ is a sequence.  Then, if $a_n$ tends to infinity as $n$ becomes arbitrarily large, the sequence is said to diverge.  In other words, if for every number $M<\infty$, no matter how large $M$, then there is still always an integer $N$ such that

\begin{equation}
    if~n\ge N, ~~then~~ |a_n| > M  
\end{equation}
In that case, we say that $(a_n)_{n=1}^{\infty}$ diverges (i.e., the limit tends to $\pm\infty$).

\begin{equation}
 \mathop {\lim }\limits_{n \to \infty } {a_n} =  \infty~~or~~ \mathop {\lim }\limits_{n \to \infty } {a_n} =  -\infty
\end{equation}
\end{definition}
\indexme{sequence!diverging}

One of the problems with deciding if a sequence is convergent using this definition is that one needs to have a limit, $L$ before the definition can be tested.  A way around this problem was described in a work by the French mathematician Augustin-Louis Cauchy in the early 1800s (although it was technically described by the mathematician Bernard Bolzano (from what is now part of the Czech Republic) first!)  The interesting insight that Cauchy (and Bolzano) had was to define convergence by the behavior of terms relative to one another rather than in an absolute sense.  Before presenting the definition, note that the logical terminology ``if and only if" (iff) is used in this definition.  This means only that the definition works in both directions.  This will be explained additionally after the statement.

\begin{definition}[Cauchy Sequences]
An infinite sequence $(a_n)_{m=1}^{\infty}$ is called a Cauchy sequence iff for every possible choice of $\epsilon>0$, there exists a positive integer $N$ such that\indexme{Cauchy sequence} \indexme{sequence!Cauchy}.

\begin{equation*}
    \| a_n -a_m \| < \epsilon~\textrm{for all } m>N ~\textrm{for all } n>N
\end{equation*}

\end{definition}
\noindent Just a note about the terminology ``iff"\indexme{iff and only iff (iff)}.  In the defintion above, the ``iff" indicates that if you specify a value for $\epsilon$, then you can always find values for $n,m>N$ such that the proof is true.  Also, if you specify a value of $N$, then for all $n,m>N$ you can always find a value of $\epsilon$ such that the proof is true.  So, the proof works both directions: an $N$ is implied to exist if the sequence is Cauchy.  In the other direction, if an $N$ is specified, then a value of $\epsilon$ is supposed to exist if the sequence is Cauchy.  

The usefulness of Cauchy sequences can be found in the following theorem, versions of which were proven by both Cauchy and Bolzano.

\begin{theorem}[Cauchy criterion for convergence of a sequence]
A necessary and sufficient condition for convergence of a sequence $(a_n)_{n=1}^{\infty}$ is that it is a Cauchy sequence. \indexme{sequence!Cauchy criterion for convergence}
\end{theorem}
Again, some additional explanation is needed here regarding the term ``necessary and sufficient"\indexme{necessary and sufficient}.  This terminology is related to ``if and only if" in a sense.  The terminology ``necessary and sufficient" is used to indicate whether or not the conclusions of the theorem always imply the initial statement (i.e., whether or not the proof is valid in reverse).  In short, ``necessary and sufficient" means that the proof is always true in the forward direction (all Cauchy sequences converge), but only sometimes true in the reverse direction (not all  convergent sequences are Cauchy sequences).  In this case, ``sufficient" is used to indicate that the observation of convergence is evidence to suggest that a sequence could be Cauchy (a non-converging sequence, therefore, cannot be Cauchy), but it is not sufficient evidence to prove that it is true.  Something else must be added (e.g., the definition given above) to determine if a convergent sequence is Cauchy. 

In summary, the concept of a Cauchy sequence is a powerful one because it allows one to investigate convergence properties of a sequence without first knowing the limit.  This closes an important logical gap (i.e., if one knows the limit of a sequence, then it is already obvious if the sequence converges or not!), and was an important landmark in abstract mathematical analysis.

\subsection{Series}

With the concept of sequences defined, it is relatively straightforward to define an \emph{infinite series}.  In short, an infinite series is the sum of some infinite sequence.  However, it is frequently useful to think about them in the following sense.

\begin{definition}[Series] \label{seriesdef}
Suppose $(a_n)_{n=0}^{\infty}$ is a sequence.  Now, define the \emph{partial sums} of the sequence by \indexme{series!definition} \indexme{series!partial sum}
\begin{align*}
    S_N &= \sum_{n=0}^{n=N} a_n \\
    \intertext{i.e.,}
    S_0&= a_0\\
    S_1&=a_0+a_1 \\
    S_2&=a_0+a_1+a_2\\
    S_3&=a_0+a_1+a_2+a_3\\
    &\ldots
\end{align*}
Then, the $S_n$ are known as the partial sums of the series.\indexme{series!partial sums}  As we allow $N\rightarrow\infty$, the partial sums define the following infinite series
\indexme{series}

\begin{equation}
    S = \sum_{n=0}^\infty a_n
\end{equation}
\indexme{series!infinite}
\end{definition}

There are two important things to note here.
\begin{enumerate}
    \item As noted above, it is not necessary that the lower bound of a series start at $0$.  However, most series used routinely in applied mathematics start at either $0$ or $1$.  If starting at some integer other than $0$, the definitions above would be modified in the obvious way.
    \item Not all series necessarily \emph{converge}.  In fact, it turns out that even some \emph{non-convergent} series are useful in applied mathematics.  In fact, if you have ever applied Stirling's approximation for the factorial
    \begin{equation}
        n! \sim \sqrt{2\pi n}\left(\frac{n}{e}\right)^n \left(1 +\frac{1}{12n}+\frac{1}{288n^2} - \frac{139}{51840n^3} -\frac{571}{2488320n^4}+ \cdots \right)
    \end{equation}
    then you have used a non-convergent series.  The concepts of convergence are one of the most important in applied mathematics, so some review of important results will be presented in the material following.
    
    \item Infinite series are often used to represent whole functions. As described under the material on sequences, we can think of the terms in the series to be the sum of a sequence of functions of both $n$ and some other independent variables. 
\end{enumerate}

\subsection{Series Convergence}

To start, it is important to define what it means for a series to converge.  As we have seen above, the list of partial sums form a sequence, so the study of convergence of sequences and series are substantially intertwined. 

\begin{definition}[Convergence of an Infinite Series]
Suppose we define the partial sums of an infinite series,  $S_n$, as we have above.  Then, the series is said to \emph{converge} if the partial sums tend toward a fixed limit $L$ as $N\rightarrow\infty$
\indexme{series!convergence} \indexme{convergence!series}

\begin{equation}
   S= \mathop {\lim }\limits_{N \to \infty } s_N = \sum_{n=0}^{N} a_n =  L
\end{equation}
\end{definition}
\begin{definition}[Absoslute Convergence of an Infinite Series]
Suppose we define the partial sums of an infinite series,  $S_n$, as we have above.  Then, the series is said to \emph{converge absolutely} if the partial sums tend toward a fixed limit $L$ as $N\rightarrow\infty$
\indexme{series!convergence} \indexme{convergence!series}

\begin{equation}
   S= \mathop {\lim }\limits_{N \to \infty } s_N = \sum_{n=0}^{N} |a_n| =  L
\end{equation}
\end{definition}

\noindent There are a number of \emph{convergence tests} that can be applied to determine if a sequence (or series) converges.  The tests below apply for any kind of series (except as noted), and include (1) the integral test, (2) the comparison test, (3) the limit comparison test, (4) the ratio test, (5) the root test, and (6) the Leibniz convergence test.   These are presented, without proof, as follows.

The following definitions and convergence theorems are provided below, without proof.  The theorems are all well-known, and you may have encountered in your studies of calculus.  Proofs can be found in any introductory text on calculus that covers infinite series.

\begin{theorem}[The integral test]
Let $(a_n)_{n=1}^\infty$ be a nonnegative (i.e., each term is positive or zero) sequence, and let $f$ be a continuous, monotonically decreasing function on $[0,\infty)$ defined such that
\begin{equation}
    f(n) = a_n~~\text{for }n\ge 1
\end{equation}
\indexme{series!integral test}
Then the series 
\begin{equation}
    \sum_{n=1}^{n=\infty} a_n
\end{equation}
converges if and only if the integral
\begin{equation}
  \int_{x=1}^{\infty} f(x) dx
\end{equation}
is finite.
\end{theorem}

\begin{theorem}[Comparison Test]
Suppose it is known that the series $\sum_{n=0}^{\infty} |b_n|$ converges.  If, for a second series $\sum_{n=0}^{\infty} |a_n|$ the condition $|a_n|\le |b_n|$, then the second series in $a_n$ also converges. 
\indexme{series!comparison test}
\end{theorem}

\begin{theorem}[Limit Comparison Test]
Suppose the series $\sum_{n=0}^{\infty} |b_n|$ converges.  For a second series, $\sum_{n=0}^{\infty} |a_n|$, we wish to determine the convergence properties.  If the condition 

\begin{equation}
    \mathop {\lim }\limits_{N \to \infty }  \left| \frac{a_n}{b_n}\right| =  L
\end{equation}
where $L$ is some finite number, then the series associated with $a_n$ also converges.
\indexme{series!limit test}
\end{theorem}

\begin{theorem}[Ratio Test]
Suppose $a_n$ is not equal to zero for all values of $n$.  We can say that a series converges (absolutely) if \indexme{series!ratio test}
\begin{align}
   \underset{n\rightarrow\infty}{\lim} \frac{|a_{n+1}|}{|a_{n}|} < 1
\end{align}
This latter expression is called the (general) ratio test (which you may have learned in a course on calculus).  Note that if $r=1$, it is not obvious whether the series converges or diverges.
\indexme{series!ratio test}

\noindent The converse is also true; that is, a series is said to \emph{diverge} if 
\begin{align}
  \underset{n\rightarrow\infty}{\lim}  \frac{|a_{n+1}|}{|a_{n}|} > 1
\end{align}
If the limit of this quantity is exactly equal to 1, then nothing can be said about convergence of the series. 
\end{theorem}

\begin{theorem}[Root Test]
Suppose that one computes the limit
\begin{align}
   \underset{n\rightarrow\infty}{\lim} ({|a_{n}|})^{\tfrac{1}{n}} = L
\end{align}
and that $L$ is a finite number such that $L<1$.  Then the series converges (absolutely).  If $L>1$ then the series diverges.  If $L=1$, the convergence of the series is undecidable by the root test.
\indexme{series!root test}
\end{theorem}

\begin{theorem}[Leibniz convergence test]  

Suppose the sequence $(|a_n|)_{n=1}^{n=\infty}$ is
\begin{enumerate}
    \item A monotonically decreasing sequence ($|a_{n+1}| < |a_{n}|$ for all values of $n$), and
    \item $\underset{n\rightarrow\infty}{\lim} a_{n}=0$
\end{enumerate}
Then the alternating series
\begin{align}
    \sum_{n=0}^\infty (-1)^n a_n \\
     \sum_{n=0}^\infty (-1)^{n+1} a_n
\end{align}
converge.
\end{theorem}
\indexme{series!Leibniz convergence test}

\subsection{Power Series}

A \emph{power series} is just a label applied to a particular kind of series.  A power series is defined by any series of the form

\begin{equation}
 f(x) =   \sum_{n=0}^{\infty} c_n x^n
 \label{power1}
\end{equation}
\indexme{series!power}
One very significant distinction that separates power series from series more generally is that power series are defined by the sum of a \emph{sequence of functions}.  Note that comparing a power series with the general series described above, we have the important correspondence

\begin{equation}
    a_n = c_nx^n
\end{equation}

Note that for a power series (as with series in general) the lower index can start at any number, although this number will generally be 0 or 1.  Power series are one of the most frequently used devices in applied mathematics.  Note that the power series can also be shifted to be defined around any point, $a$, in the domain as follows

\begin{equation}
  f(x) = \sum_{n=0}^{\infty} d_n (x-a)^n
  \label{power2}
\end{equation}
\indexme{series!power} \indexme{power series}
where now we have the correspondence

\begin{equation}
    a_n = d_n(x-a)^n
\end{equation}
This shifted form of the power series can be a convenient notation to use.  Note, however, one could in principle expand the series, and regroup terms to recover Eq.~\eqref{power1}.  Therefore, Eqs.~\eqref{power1} and \eqref{power2} represent two different ways of \emph{expressing the same series} (i.e., they can be equivalent).  

In analogy with the discussion above, one can define the \emph{partial sums} of a power series as a sum truncated at some finite term numbered $N$.  Being explicit, note that the first few partial sums of the power series defined above are given by

\begin{align*}
    S_0 &= c_0\\
    S_1 &= c_0 + c_1 x \\
    S_2 &= c_0 + c_1 x +c_2 x^2\\
    S_3 &= c_0+ c_1 x + c_2 x^2 + c_3 x^3\\
    & \ldots 
\end{align*}
It is useful to compare these partial sums with those given in Definition \ref{seriesdef}.

In the next example, we provide the power series representation of a few familiar functions.

\begin{svgraybox}
\begin{example}[Power Series.]
Trigonometric functions are examples of transcendental functions (as introduced above). By definition, this means that there is no \emph{finite} set of algebraic steps that one can use to generate such functions (e.g., there is no polynomial that will reproduce them exactly).   If you think about trigonometric functions, you will realize that you have probably never been asked to compute the value of, say, the sine function yourself.  This is because the sine function must be defined by an algorithm that has, in principle, an infinite number of steps!

Although functions like $\sin(x)$ and $\cos(x)$ cannot be represented by a finite number of algebraic steps, they \emph{can be represented by an infinite number of algebraic steps}.  The following are the power series for these functions.

\begin{align}
    \cos(x) &= \sum_{n=0}^\infty \frac{{(-1)}^n}{(2n)!} x^{2n} = 1- \frac{x^2}{2!}+\frac{x^4}{4!}-\frac{x^6}{6!}+\ldots\\
     \sin(x) &=  \sum_{n=0}^{\infty} \frac{(-1)^nx^{2n+1}}{(2n+1)!}= x - {\frac{x^{3}}{3!}} + {\frac{x^{5}}{5!}} - \cdots 
\end{align}
While it is true that, in principle, the value of these functions is defined by an infinite number of terms, it is often the case with series that the terms defined by increasingly large $n$ become smaller and smaller in magnitude.  In a practical sense, this means that it is often possible to get a good representation for a series at a point of interest by using a finite number of terms.  As an example, in the figure below we plot the approximation to the sine function using $N=2,5,10,$ and $20$ terms.  It is clear that even for $N=2$, the power series actually gives a good approximation to the function in the interval of approximately $0<x<1.5$.

\begin{centering}
\includegraphics[scale=.5]{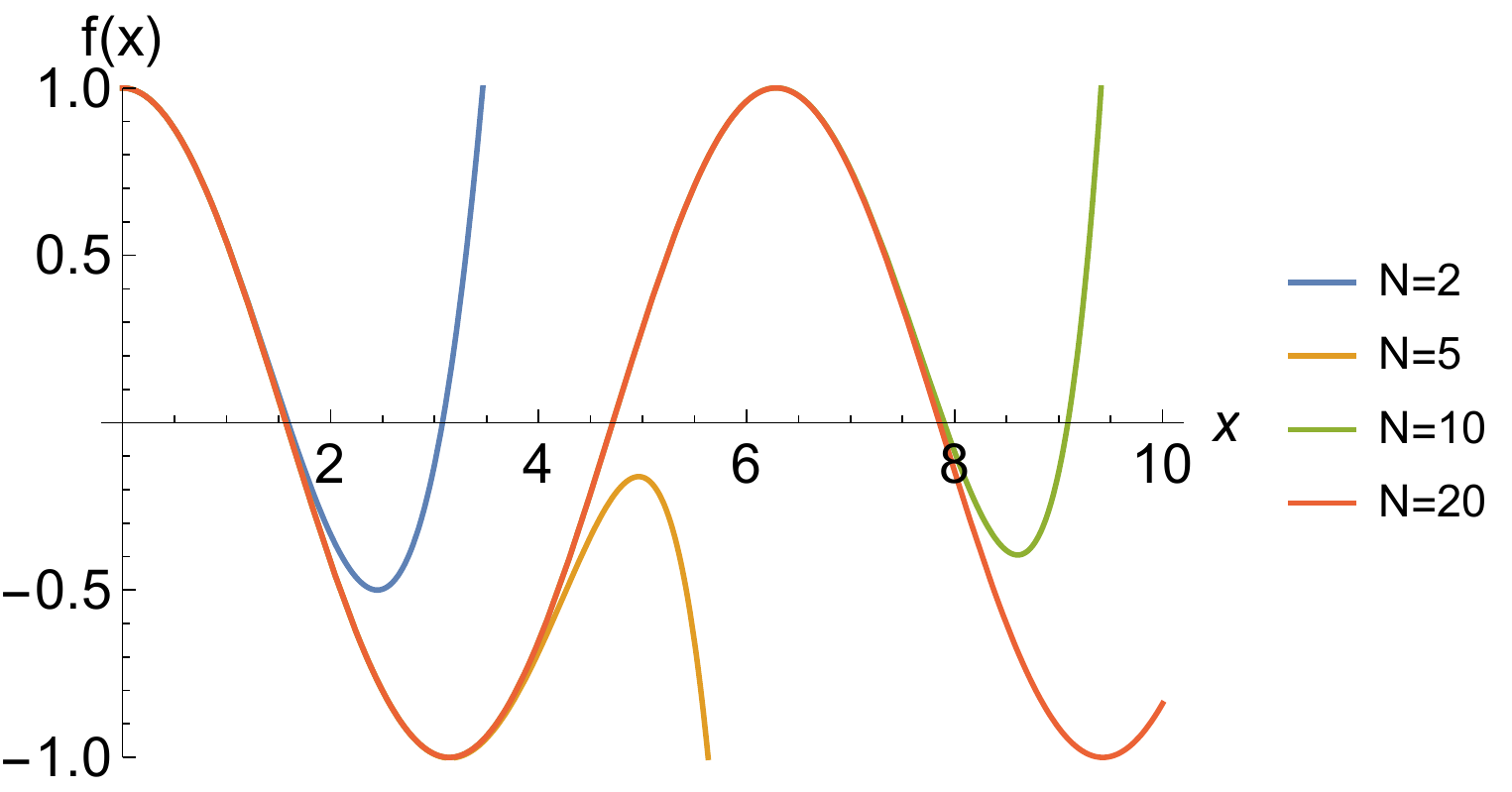}
\captionof{figure}{The sine function approximated by its power series, using an increasing number of terms.}
\end{centering}

\end{example}
\end{svgraybox}

\section{Taylor's Theorem and Taylor Series}

A Taylor series is a particular power series named after the British mathematician, Brook Taylor, who made it famous.  The genius of Taylor was to work out how one might compute the terms in a power series for any function for which one knows the derivatives.  Interestingly, Newton himself also developed what we now call the Taylor series \citep{ferraro2008rise}, although his original manuscripts on the topic remained unpublished in his lifetime.  Even more curiously, a Scottish mathematician named Colin Maclaurin promoted Newton's unpublished methods on infinite series (with Newton's approval) in a mathematics textbook; in this form the series were expanded around zero and are called \emph{Maclaurin} series.  Newton, Maclaurin, and Taylor were all members of the Royal Society of London contemporaneously.   Both Maclaurin and Taylor knew Newton personally \citep{oconnor}, and it seems plausible that Maclaurin and Taylor very likely knew one another too.  While these three scientist-mathematicians knew of each other, much of the work by them on series contains unique contributions from each of them, and, in many cases, there was even ``re-invention" of the same ideas independently.  This perhaps underscores the idea that science (especially in the days where communication was not immediate as it is today) is often a pattern of discovering and rediscovering.  There is an adages that says something to the effect of ``the person credited with a discovery is usually just the one that was able to describe it well enough that it need not be rediscovered again."  Taylor series turned out to be one of the most important developments in analytical mathematics; therefore, it should not be too surprising that the ideas were built up over a period of time, and not exclusively Taylor's.  

\subsection{Taylor Series}

Taylor series are a particular form of power series.  What makes them \emph{distinct} from generic power series is that the coefficients $c_n$ are defined by combinations of the derivatives of the function that is to be represented in series form.   To start out, lets constructively define the \emph{Taylor series} as follows. \indexme{Taylor series}

\begin{definition}[Taylor Series]
For a function that has an infinite number of derivatives at a point, $x=a$, that function has a Taylor series expansion that represents the value of the function at a location $x=(a+\Delta x)$ as follows (where $\Delta x = [x-a]$ or, equivalently, $x=a+\Delta x$)

\begin{align*}
    \underbrace{f(a+\Delta x)}_{f(x)}& = f(a) + \Delta x \left.\frac{d f}{dx}\right|_{x=a}+\frac{1}{2!} (\Delta x)^2 \left.\frac{d^2 f}{dx^2}\right|_{x=a}+\frac{1}{3!} (\Delta x)^3\left. \frac{d^3 f}{dx^3}\right|_{x=a}+ \ldots \\
    \intertext{or, equivalently, the following notation is frequently used}
    f(x) &= f(a) + (x-a) f'(a)+\frac{1}{2!} (x-a)^2 f''(a)+\frac{1}{3!} (x-a)^3 f'''(a)+ \ldots 
    \end{align*}
Finally, note this last form can be written in the familiar summation notation for the Taylor series
    \begin{align*}
    f(x)& = \sum_{n=0}^\infty \frac{f^{(n)}(a)}{n!} (x-a)^n
\end{align*}
\end{definition}
\indexme{series!Taylor}\indexme{Taylor series}
\noindent where $f^{(n)}$ is the n$^th$ derivative of the function $f$, and the factorial of an integer, $n$, is defined by $n!= n\cdot(n-1)\cdot(n-2)\ldots 2\cdot 1$, and $0!=1$.
Note that with this definition, we have that a Taylor series is a (shifted) power series of the form

\begin{equation*}
    d_n = \frac{f^{(n)}(a)}{n!}
\end{equation*}
where, recall, the coefficients $d_n$ are defined for the shifted power series given by Eq.~\eqref{power2}, and these coefficients are distinct from those for the non-shifted power series.  

This definition says nothing about whether or not this series \emph{converges} or not; the topic of convergence is addressed material that follows.  In the introductory material on Series above, the notion of partial sums of sequences was defined as one type of series.  This notion is now extended to partial sums of the Taylor series; each of these partial sums represents an entire \emph{function} rather than simply a number.

\begin{definition}[Partial Sums]  The partial sum of a Taylor series is the sum truncated at some positive integer $N$ such that 

\begin{equation*}
    f_N(x) =  \sum_{n=0}^N \frac{f^{(n)}(a)}{n!} (x-a)^n
\end{equation*}
\end{definition}

Note, that every infinite series can also be described by the \emph{sequence of its partial sums}.  This is an important point, because it provides and illustration of how a \emph{function} might be described by a sequence of functions, $f_n$, that converges to a desired function, $f$.  This will be shown in the next example.  However, first we need to discuss what \emph{kinds} of functions have convergent Taylor series.  These functions are the \emph{analytic functions} that were described earlier.  Here, a second (and more fundamental) definition of the term analytic is provided.

The Taylor series provides us with an intuitive tool to define \emph{analytic} functions (which were defined previously in \S \ref{continuity}).  This second definition is as follows.  Note that $x=a+\Delta x$ implies that $\Delta x = (x-a)$.  Thus for $a=0$, we have $\Delta x = x$ (a Maclaurin series).

\begin{definition}[Analytic Functions-- Definition 2]
An \emph{analytic function} is a function which (i) is defined (i.e., it has a computable value) on an open interval ($I$) of the real line ($I$ is the set of values, $x$, such that $x\in(a,b); \mbox{ or, in set notation } I=\{x: a < x < b\}$), (ii) its Taylor series converges to a definite value on the interval $I$, and (iii) the definite value it converges to is equal to $f(x)$ for all values of $x$ in the interval $I$.
\indexme{function!analytic} \indexme{analytic function}
\end{definition}
\noindent Admittedly, that is quite a mouthful of requirements; however, many functions that of are of interest in science and engineering are analytic.  For example, all polynomials on a finite domain are analytic.  The functions $\sin x$, $\cos x$, and $\ln x$ are analytic on appropriate domains.  So, although the notion of being analytic seems stringent on the surface, most of the functions we commonly encounter are analytic on some portion of their domains.  

One of the nice properties of all analytic functions is that \emph{all} of the derivatives of an analytic function are bounded.  In other words, no derivative of an analytic function can grow to infinity anywhere in the domain where the function is analytic.  This is strongly suggested by the fact that all analytic functions have convergent Taylor series (i.e., any derivative that approached infinity would prevent the Taylor series from converging to a finite value).  Nonetheless, the following theorem is stated without proof.

\begin{theorem}[Derivatives of analytic functions]
Let $x_0$ be a point, and $r$ a positive value.  Let a function, $f$ be analytic in an interval around the point $x_0$; that is, for some positive $r$, the function $f$ is analytic in $\{I: x_0-r < x < x_0+r\}$ 
Then there exists a positive real number $M$, such that for every $n \in \mathbb{N}$

\[  |f^{(n)}(x)|  \le \frac {M r \, n!} {{(r - |x - x_0| )}^{(n + 1)} }\]

In other words, the derivatives of $f$, $f^{(n)}$, are bounded in the interval, and the derivative cannot grow arbitrarily large within the interval.
\end{theorem}

The next example is one that shows that all polynomials are analytic functions.\\

\begin{svgraybox}
\begin{example}[Polynomials are Analytic.]
Suppose we have a polynomial on an open interval other than $I=(-\infty, \infty)$.  The polynomial is analytic everywhere in its domain.

This is not a proof, but it could easily be turned into on.  Instead, consider a concrete polynomial, for example 
\begin{equation*}
    f(x) = x^3 + 2x^2 +2 x +1
\end{equation*}
defined on the \emph{open} interval $0<x<1$. Suppose we want to determine the Taylor series for this polynomial around the point $a=0$.  We can expand this polynomial in a Taylor series around the point $x=a=0$ to determine the value of the function at the point $x=(a+\Delta x)$.  To do so, we first need to compute all of the derivatives of $f$.  This seems as though it might be nearly impossible, until we note
\begin{align*}
    f(x) &= x^3 + 2x^2 +2 x +1\\
    f'(x)&=3 x^2 + 4 x + 2 \\
    f''(x) &= 6 x + 4 \\
    f'''(x) &= 6 \\
    f^{(4)}(x) &= 0\\
    f^{(5)}(x) &= 0
    \ldots&
\end{align*}
Now, according to the formula, the Taylor series is given by
\begin{align*}
f(a+\Delta x) &= f(a) + \Delta x \frac{d f}{dx}+\frac{1}{2!} (\Delta x)^2 \frac{d^2 f}{dx^2}+\frac{1}{3!} (\Delta x)^3\frac{d^3 f}{dx^3}+\ldots 
\end{align*}
Noting that, for this example, $\Delta x=x-a=x-0=x$, then we have
\begin{align*}
    f(0) &= 1\\
    f'(0)&= 2 \\
    f''(0) &= 4\\
    f'''(0) &= 6 \\
    f^{(4)}(0) &= 0\\
    f^{(5)}(0) &= 0
    \ldots&
\end{align*}

\begin{align*}
f(x)&= f(0) + x \left.\frac{d f}{dx}\right|_{x=0}+\frac{1}{2!} x^2 \left.\frac{d^2 f}{dx^2}\right|_{x=0}+\frac{1}{3!} x^3\left.\frac{d^3 f}{dx^3}\right|_{x=0}+\ldots \\
\intertext{Substituting the derivatives above yields}
f(x)&= 1 + x (2) +\frac{1}{2!} x^2 (4)+\frac{1}{3!} x^3 6+0+0+0+\ldots
\intertext{Simplifying terms, we find}
f(x)&=x^3 + 2x^2 +2 x +1
\end{align*}
which was the original polynomial.  Thus, this polynomial is its own Taylor series!  This is actually true for all polynomials, and it is not hard to prove using the principle of induction.  The ancillary conclusion that can be reached is that all polynomials defined on a (open, non-infinite) domain are analytic functions.
\end{example}
\end{svgraybox}

\begin{svgraybox}
\begin{example}[Partial Sums Example.]

In the previous example, the polynomial $f(x) = x^3 + 2x^2 +145 x -1$, when expanded around $x=0$, had the following partial sums

\begin{align*}
    f_0(x)&=1\\
    f_1(x) &= 2 x +1 \\
    f_2(x) &= 2 x^2+ 2 x +1 \\
    f_3(x) &= x^3 + 2x^2 +2 x +1\\
    f_4(x) &= x^3 + 2x^2 +2 x +1\\
    &\ldots
\end{align*}
It is interesting to see how the sequence $(f_1, f_2, f_3,f_4, f_5, \ldots)$ converges to the function $f(x)$ exactly in this instance.  To see this, we plot each of the partial sums on $x\in[0,1]$ giving

\begin{centering}
\includegraphics[scale=.5]{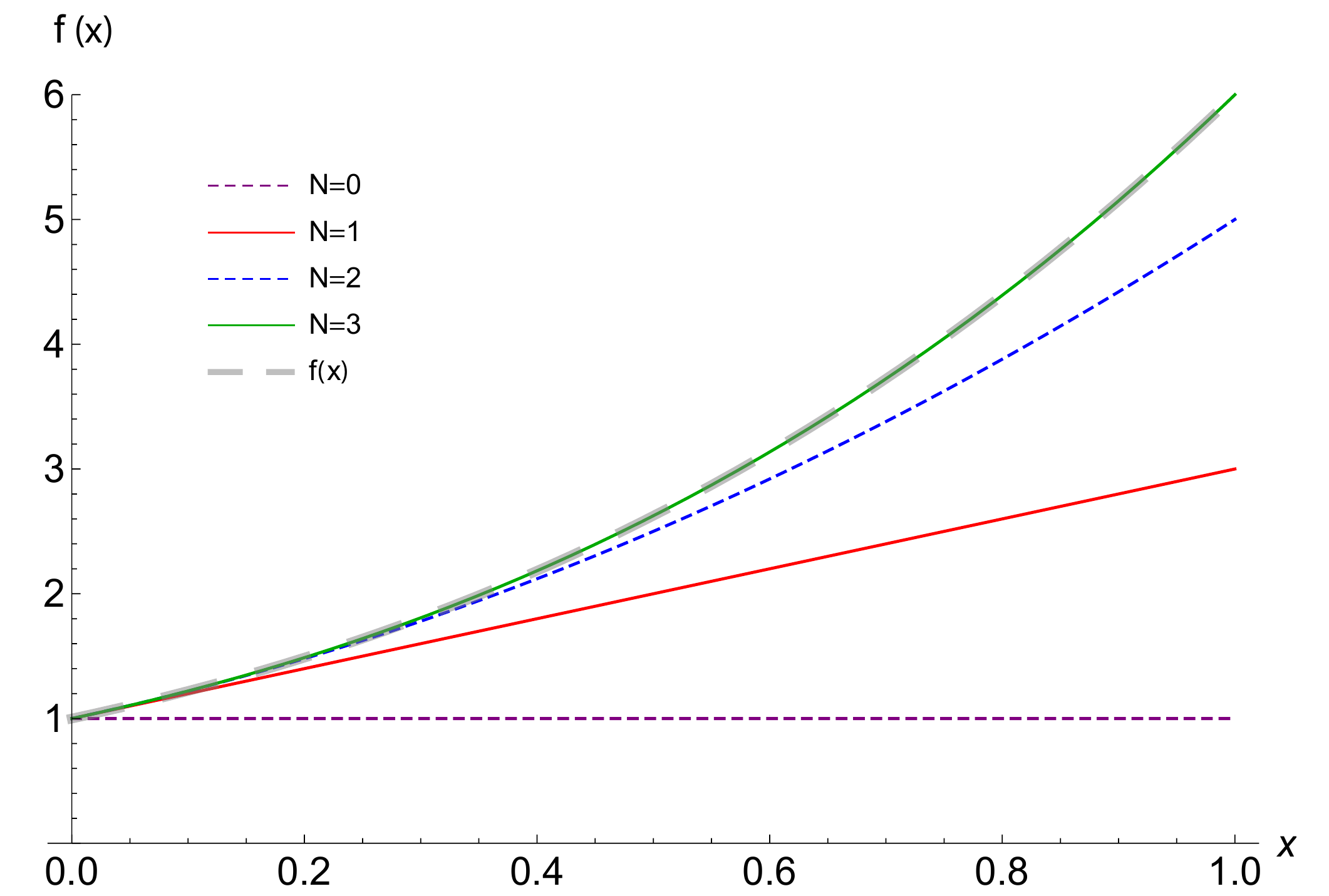}
\captionof{figure}{The partial sums (Taylor polynomials) of the Taylor series of $f(x) = x^3 + 2x^2 +2 x +1$.}
\label{taylor_poly}  
\end{centering}

\end{example}
\end{svgraybox}

\subsection{Taylor Series Construction}

So, we have established that the Taylor series can represent all analytic functions, and we examined a particular example of a polynomial that illustrated that the series actually did what we hoped it would.  One might, at this juncture, wonder: \emph{why does the Taylor series work?}   It turns out, that this is reasonably easy to prove.  Again, we will not do this as a formal proof, but rather as a construction that illustrates the method.  A formal proof can be constructed from this outline.

To start, recall the fundamental theorem of calculus for the function $f$ in the form

\begin{equation}
\int\limits_{x=a}^{x=a+\Delta x} \frac{df}{d\xi} d\xi = f(a+\Delta x) - f(a)
\end{equation}
\noindent where here we have used $\xi$ as the variable of integration.  This is already reasonably suggestive.  Rearranging, note that we have
\begin{equation}
f(a+\Delta x)= f(a) + \int\limits_{x=a}^{x=a+\Delta x} \frac{df}{d\xi} d\xi 
\end{equation}
Now, we use integration by parts once (setting $u=df/d\xi$, $dv= d\xi$) to give

\begin{equation}
f(a+\Delta x)= f(a) + \left. (a+\Delta x) \frac{df}{dx}\right|_{a+\Delta x}-a\left.\frac{df}{dx}\right|_{a}- \int\limits_{x=a}^{x=a+\Delta x} \xi \frac{d^2f}{d\xi^2} d\xi
\label{taylor1}
\end{equation}
Now note, also by the fundamental theorem of calculus
\begin{equation}
    \int\limits_{x=a}^{x=a+\Delta x} (a+\Delta x) \frac{d^2f}{d\xi ^2} d\xi
    =(a+\Delta x) \left.\frac{df}{dx}\right|_{a+\Delta x}-(a+\Delta x) \left.\frac{df}{dx}\right|_{a}
\end{equation}
\noindent Rewriting this as

\begin{equation}
    (a+\Delta x) \left.\frac{df}{dx}\right|_{a+\Delta x}=\int\limits_{x=a}^{x=a+\Delta x} a \frac{d^2f}{d\xi ^2} d\xi
    +(a+\Delta x) \left.\frac{df}{dx}\right|_{a}
\end{equation}
\noindent Finally, substituting this into  Eq.~(\ref{taylor1}) gives
\begin{equation}
f(a+\Delta x)= f(a) + \left. \Delta x\frac{df}{dx}\right|_{a}- \int\limits_{x=a}^{x=a+\Delta x} (\xi-a)\frac{d^2f}{d\xi^2} d\xi 
\end{equation}
Repeating this process $n$ times yields Taylor's formula to term $n$.

\subsection{Uniform Versus Pointwise Convergence for Power and Taylor Series}\label{convergence_types}

Both power and Taylor series are unique to the study of series because they involve an infinite sum of {functions indexed by the integers} $\ge 0$. Assuming, for example, that we have a single independent variable in space, these series can then be functions of both $x$ in some interval $I$ and $n\in \mathbb{N}_0$.  More explicitly, we have a sequence of functions whose partial sums are
\begin{equation}
    f_n(x) = \sum_{n=0}^N d_n (x-a)^n
\end{equation}
or, assuming a Taylor series where $d_n = \frac{f^{(n)}(a)}{n!}$

\begin{equation}
    f_n(x) = \sum_{n=0}^N \frac{f^{(n)}(a)}{n!}(x-a)^n
\end{equation}
Thus, as far as convergence is concerned, we must consider the convergence of these series for each point $x$ as $n\rightarrow \infty$.  There are many senses (or modes) of  convergence that have been adopted in mathematics.  Here, the discussion will be focused on two that are particularly relvant to series representations of functions; they apply to the Taylor series described in this section as well as the Fourier series that will be introduced in future chapters.

Before stating the definitions for these two kinds of convergence, it is possible to generate some intuition about them.  The convergence known as \emph{uniform} is defined by the idea that one can make the series approximation to the function exhibit a maximum error, $\epsilon$, for \emph{every} point in the domain.  In contrast, pointwise convergence is not as strict.  It suggests that, while every point in the domain must converge, one cannot necessarily ascribe a singe error, $\epsilon$, which is the maximum error.  Instead, the error associated with each point must be determined independently.  Before continuing, recall that a sequence of functions $f_n(x)$ (e.g., the partial sums of a Taylor series) is said to converge if

\begin{equation}
    \underset{n\rightarrow\infty}{\lim} f_n(x) = f(x)
\end{equation}
Now, we are prepared to present the definitions of pointwise and uniform convergence.

\begin{definition}[Pointwise Convergence of a Sequence of Functions]
A sequence of functions, $f_n(x)$ is said to converge \emph{pointwise} on some interval $I$ if for all $x\in I$ and for \emph{every} specified error $\epsilon(x)$, we can always find a value of $N$ such that \indexme{convergence!pointwise} \indexme{series!pointwise convergence}

\begin{equation}
    |f_n(x) -f(x)| < \epsilon(x), \textrm{ for all } n\ge N
\end{equation}
\end{definition}

\begin{figure}[t]
\sidecaption[t]
\centering
\includegraphics[scale=.45]{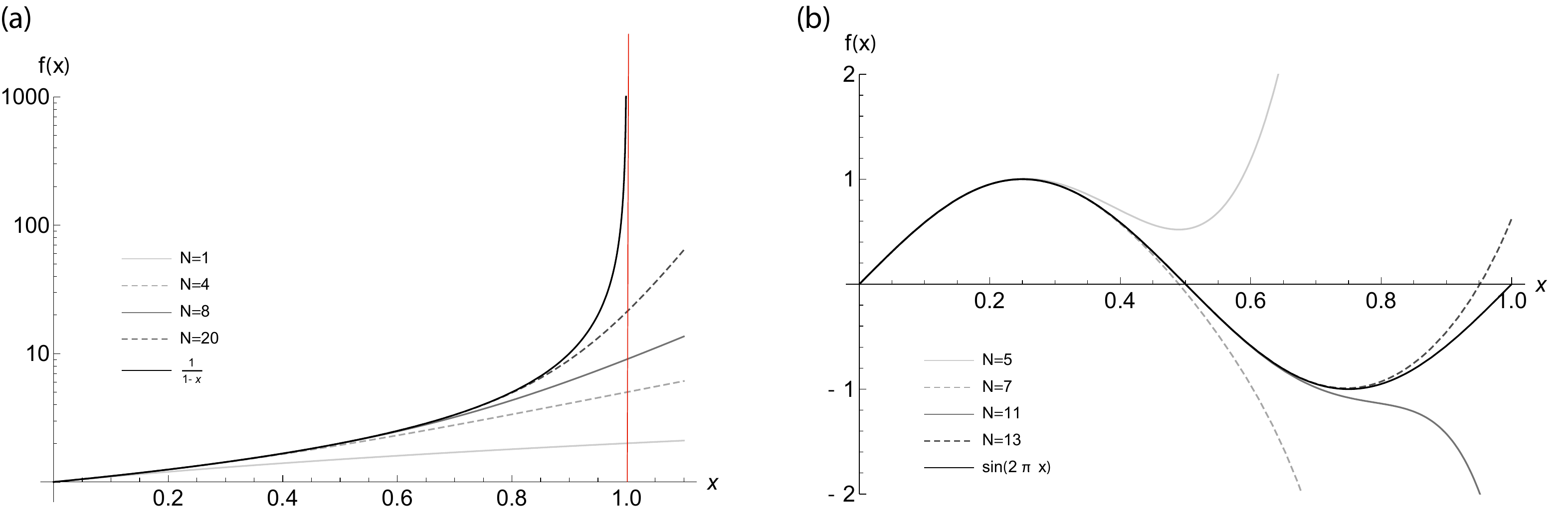}
\caption{Series convergence types.  (a) A uniformly converging series. (b) A series that converges pointwise.  }
\label{converge}       
\end{figure}

The critical thing to note here is that the value of $N$ needed depends upon both the location, $x$, and the value of $\epsilon(x)$ selected at that point.  For this case, there is no guarantee that there is a single combination of both $\epsilon(x)$ and $N$ that will be valid for all points.  This can be a bit confusing to ponder, but an example can help.  In Fig.~\ref{converge}(a), a Taylor series approximation to the function $f(x) = 1/(x-1)$ on the interval $x\in[0,1)$ is illustrated.  It is clear from this figure that increasing the number of terms in the partial sums leads to improved approximations for the series.  However, note the following.  At the point $x=1$, the function tends to infinity.  There are no polynomials that tend toward infinity at a some finite value of $x$; thus, the error in the vicinity of the point $x=1$ is complicated.  While it is true that for any $x<1$ (no matter how close we come to the point $x=1$), it is always possible to find a value of $N$ large enough that the error is smaller than any value of $\epsilon(x)$ that we pick, we can only determine this value of $N$ once both $x$ and $\epsilon(x)$ are specified.  However, there is no finite value of $N$ that works for all points!  For example, suppose we find the minimum value of $N$ that provides $\epsilon < \tfrac{1}{10^6}$ for $x=999/1000$.  This value of $N$ would be undoubtedly large; however, it would not also give $\epsilon < \tfrac{1}{10^6}$ if we chose $x=9999/10000$.  Thus, the necessary value of $N$ is always achievable, but it also always depends upon the values of $x$ and $\epsilon(x)$ chosen.

In contrast to pointwise convergence, uniform convergence can be defined as follows.

\begin{definition}[Uniform Convergence of a Sequence of Functions]
A sequence of functions, $f_n(x)$ is said to converge \emph{uniformly} on some interval $I$ if for all $x\in I$ and for \emph{every} specified error $\epsilon$, we can always find a value of $N$ such that \indexme{convergence!uniform} \indexme{series!uniform convergence}

\begin{equation}
    |f_n(x) -f(x)| < \epsilon, \textrm{ for all } x\in I \textrm{ and}\textrm{ for all } n\ge N
\end{equation}
\end{definition}
Note the distinction here in the domain of the convergence.  Here, a single value of $\epsilon$ can be selected, and it is guaranteed that there is some value of $N$ large enough so that the difference between the series approximation and the function is smaller than $\epsilon$ for every point in the domain simultaneously.  In Fig.~\ref{converge}(b), a function that converges uniformly (the sine function) is plotted. While there is always some error in the Taylor polynomial approximation to the function, it is also possible to make this error as small as we like  everywhere in the domain $I$ simply by taking $N$ to be large enough.  In other words, the error has a concrete upper bound, and that bound depends only upon the size $N$, but not on the location $x$ that is chosen.

The topic of uniform versus pointwise convergence will arise again in the study of the Fourier series representation of functions.  While the distinction between the two modes of convergence appears to be subtle, the actual ramifications are large!  In a practical sense, uniform convergence means that one can often find a good approximation to a function that works everywhere in a domain of interest with some finite number of terms.  For a function that converges only in the pointwise sense, the approximation of the function in an applied sense can be much more complicated. 

\subsection{$^\star$The Taylor Series in Approximation Theory}\label{approxtheory}

Taylor series, and power series in general, do not always converge; even if they converge, they may converge only on limited domains.  For the series

\begin{equation}
 f(x) = \sum_{n=0}^\infty \frac{f^{(n)}(a)}{n!} (x-a)^n
\end{equation}
there are three possibilities for convergence, as follows.
\begin{enumerate}
    \item The series converges only for $x=0$ (in which case, it converges to zero there).
    \item The series converges for all $x$. 
    \item There exists a number $R$ such that $\sum_{n=0}^\infty \frac{f^{(n)}(a)}{n!} (x-a)^n$ converges for $|x|<R$ and diverges for $x>R$.  For $x=R$, the series might converge or diverge.
\end{enumerate}

Most of the useful Taylor series, as one might imagine, are convergent.  This means that they either (1) converge everywhere on the real line, or (2) converge in some interval $|x|<R$.  Convergent Taylor series are interesting in that they look like polynomials of infinite-order.  In fact, there is no such thing as an infinite-order polynomial (the proper term there would be a power series or Taylor series), but every finite approximation of a Taylor series is some finite order polynomial.  To see that, suppose that we set $c_n = {f^{(n)}(a)}/{n!} $, and that we conduct the expansion about $x=0$ (we can do the case for $a\ne0$, and all the following results are the same, it is just more complicated to explain).  Then, the Taylor series takes the form

\begin{equation}
 f(x) = \sum_{n=0}^\infty c_n x^n
\end{equation}
We can make finite approximations to every Taylor series, simply by truncating it at some value $N$.  Then, we have the partial sums discussed earlier, which we represent by 

\begin{equation}
 f_N(x) = \sum_{n=0}^N c_n x^n
 \label{Taylor2}
\end{equation}
These partial sums represent a \emph{sequence of functions}, and for convergent Taylor series this sequence of functions gets to be a better and better approximation for $f$ as $N$ increases.  

Suppose we consider the first five terms in a Taylor series of the form of Eq.~\eqref{Taylor2}.  A little thought will indicate that this series must take the form

\begin{equation}
    f_5(x) =  c_0 +c_1 x +c_2 x^2 + c_3 x^3 + c_4 x^4 + c_5 x^5
    \label{taylorpoly}
\end{equation}
But, this is, as noted above, just a polynomial!  Apparently, the existence of a Taylor series convergent over some interval indicates that the function being approximated can be represented by a polynomial.  It also indicates that we can make the error associated with this polynomial representation as accurate as we like simply by taking a sufficient number of terms in the expansion (this comes from the very definition of \emph{convergent}).  The type of convergence that one observes will be dependent upon the function investigated, but nonetheless one can find a (finite) polynomial that can approximate the function with as small an error as one likes.

Expressions like Eq.~\eqref{taylorpoly} are called Taylor polynomials.  One reason that this is interesting is because it shows, in effect, that any function with a convergent Taylor series has an expansion in \emph{basis functions} that are polynomials.  Specifically, the set of basis functions is given by 

\begin{equation}
    1,~ x, ~x^2, ~x^3, ~x^4, ~x^5 \ldots
\end{equation}
Suppose, for a moment, we call these functions by a symbolic name.  Let $p_0 = 1,~ p_1(x)=1,~p_2(x)=x^2,~ etc.$.  Then, the set of functions $\{p_i\}$ is sufficiently ``rich", that linear combinations of these basis functions can provide an approximation to any function with a convergent Taylor series.  From the example above with an expansion of five terms, we would have the linear combination given by

\begin{equation}
   f_5(x)= c_0 p_0(x) + c_1 p_1(x) + c_2 p_2(x) + c_3 p_3(x) +c_4 p_4(x) + c_5 p_5(x)
\end{equation}
This is quite an amazing result.  In essence it says that every function with a Taylor series that converges on some interval can be approximated to any accuracy that we like by some polynomial!  An example will be helpful here. \\

\begin{svgraybox}
\begin{example}[Polynomial approximations via Taylor series.]

For this example, we consider the familiar function $f(x)=e^x$.  The Taylor series for this function (computed for a=0) is easily calculated, and fairly well known.

\begin{equation}
    e^x = 1 + x + \frac{1}{2} x^2 + \frac{1}{6} x^3 + \frac{1}{24}x^4+ \frac{1}{120}x^5 +\ldots
\end{equation}
Once nice property of the Taylor series for $e^x$ is that it converges for every $x$ on the real line.  
It is interesting to see how well different polynomials provide estimates for $e^x$.  Suppose we consider the interval $0\le x \le 3$.  In the plot below, the first 6 Taylor polynomials are plotted for the exponential series.  To be clear, the polynomials are explicitly listed below.  Note that in each case, the function is a polynomial of the form $\sum_0^N c_i p_i(x)$, where the $p_i(x)$ are the polynomial basis functions, $p_0(x)=1$, $p_1(x)=x$, $p_2(x)=x^2$, etc.

\begin{align*}
    f_0(x) &= 1 \\
    f_1(x) &= 1 + x \\
    f_2(x) &= 1 + x + \frac{1}{2} x^2\\
    f_3(x) &= 1 + x + \frac{1}{2} x^2 + \frac{1}{6} x^3\\
    f_4(x) &= 1 + x + \frac{1}{2} x^2 + \frac{1}{6} x^3+ \frac{1}{24}x^4\\
    f_5(x) &= 1 + x + \frac{1}{2} x^2 + \frac{1}{6} x^3 + \frac{1}{24}x^4+ \frac{1}{120}x^5
\end{align*}
From this plot, it is clear that the sequence of Taylor polynomials are functions, and this sequence of functions appears to converge to the function $f(x)=e^x$ (as we had expected).  Additionally, we see that as the order of the Taylor polynomial increases, the approximation improves.  For the Taylor polynomial $f_5(x)$, the fit is quite good over the entire range plotted, with the maximum error being only about 8\%.  

\begin{centering}
\includegraphics[scale=.4]{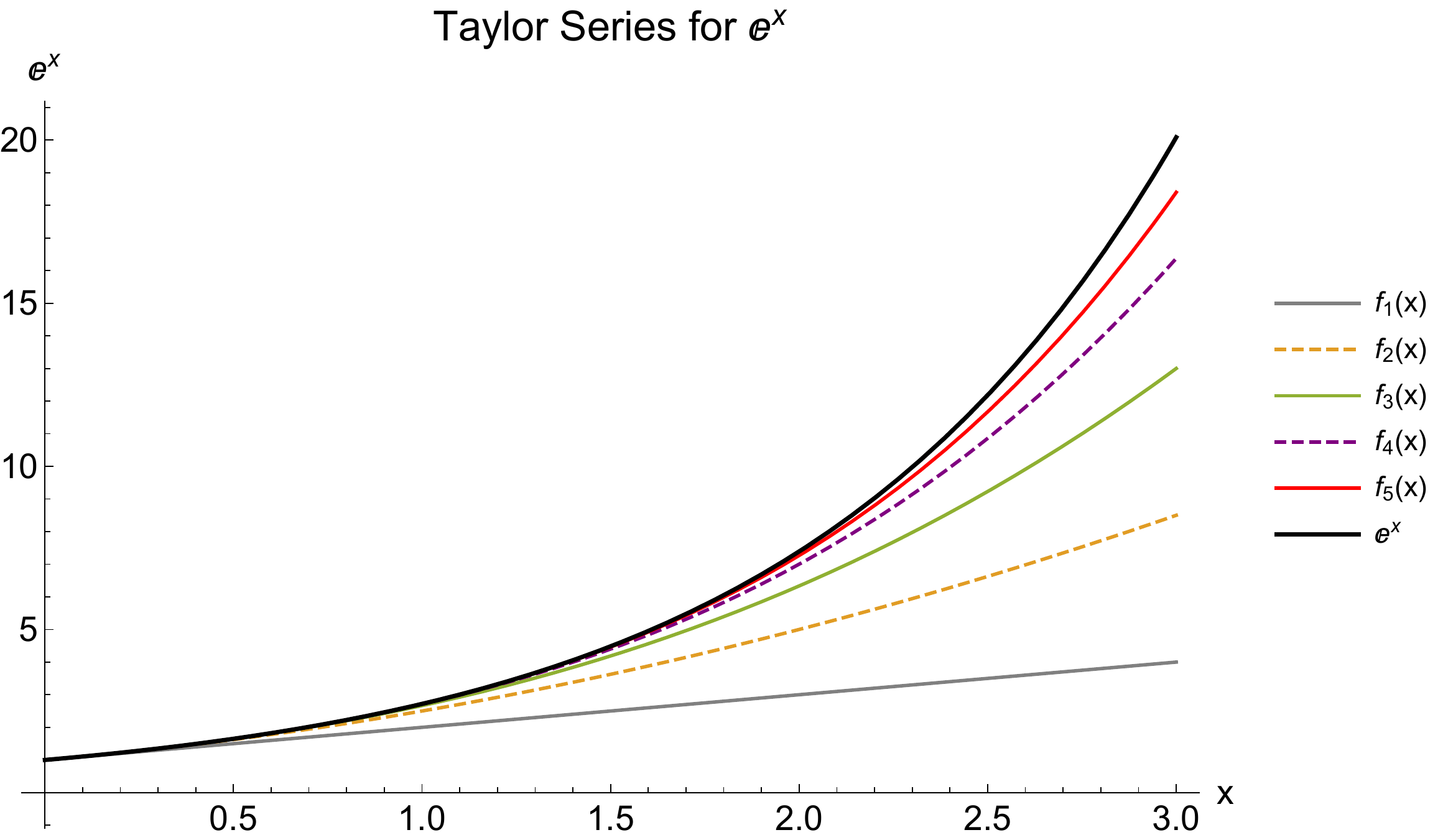}
\captionof{figure}{The first six Taylor polynomials for the function $f(x)=e^x$ plotted for the interval $0\le x\le 3$.} 
\end{centering}

\end{example}
\end{svgraybox}

It is really quite remarkable to realize that any analytic function can be approximated by a polynomial to any accuracy that we might desire.  The example of Taylor polynomials provides two interesting concepts worth reflecting on.

\begin{enumerate}
    \item Functions can be created by linear combinations of an entire set of other functions called \emph{basis functions}.  Here, the basis functions were the polynomials $p_0(x)=1$, $p_1(x)=x$, $p_2(x)=x^2$.  The linear combination was formed by taking the appropriate weighting of these functions, and summing them up.  This is done automatically during the process of computing the Taylor series.  However, it is useful at this juncture to realize that we can take any analytic function and, using a Taylor series, \emph{decompose} this function into the weighted sum of polynomial basis functions.  We will return to this idea of decomposing functions in terms of basis functions when we discuss Fourier series.
    
    \item Second, we have touched on what is known in mathematics as \emph{approximation theory}.  Approximation theory is concerned with the process of representing a function by a linear combination of basis functions, as we have done with the Taylor polynomials.   There is a much stronger statement regarding the approximation of functions with polynomials known as the Weierstrass approximation theorem.  While it also proposes polynomials as approximations to continuous functions, it is able to make even stronger statements regarding the behavior of the convergence of those approximations to the function than can be said in general for Taylor series.
\end{enumerate}
This discussion has touched on approximation theory.  In some senses, approximation theory more generally forms the basis for all of modern real (and complex!) analysis (i.e., almost every topic one can imagine in applied mathematics).  Later on in the text, the study of Fourier series is introduced. Fourier series is one of the most important tools in all of applied mathematics; and, ultimately, it is also the study of approximation theory for functions on bounded intervals.  It is important to look for patterns when undertaking the study of mathematics.  A well-understood approach of one type (e.g., approximation theory using Taylor series) can be an amazingly helpful analogue when studying new problems with similar structure (e.g., approximation theory on finite intervals using Fourier series).

\section{$^\star$Functionals and Integral Transforms}

In this section, a concept that has already been ``seen" by most students is re-evaluated for the purposes of later use.  In particular, this section is concerned with the notion of \emph{functionals}\indexme{functional}\indexme{function!functional}, which, in short, are just functions of functions.  Recall that a \emph{function} uniquely maps every element of a set $A$ (the domain) to a single element in set $B$ (the range).  We usually think of the set $A$ as being some portion of the real line (or of $\mathbb{R}^2$, or the complex plane, etc.) defined by a number.  Similarly, we usually think of the output of the function as being a number; thus a function maps numbers (or vectors of numbers in multiple dimensions or for complex numbers) to numbers. 

For a \emph{functional}, the same definition is true, except both the input set, $A$, and the output set, $B$, are not numbers, but entire \emph{functions}.  The following is a definition

\begin{definition}[Functional]
A function whose \emph{domain} is a set of functions, and whose \emph{range} is another set of functions (including the constant functions) is known as a \emph{functional} (cf., \citep[][Chp. IV]{courant1953methods})
\end{definition}

Some examples will make the notion of a functional much clearer.  

\begin{svgraybox}
\begin{example}[A simple functional.]

Suppose we consider the set of all possible functions, $f(x;a) = a x^2$, where $a$ is some constant parameter, and $x \in [-1,1]$.  In this notation, we have indicated the parameters associated with a function by listing them (preceded by a semi-colon) with the independent variables.  This is a common notation that is used when the parameters of a function might also be considered to vary.   

We can define the following functional, $F[f(x;a)]$ that relates each function $f$ to a single real number as follows.

\begin{align*}
    F[f(x;a)] &= \int_{x=-1}^{x=1} a x^2\, dx
\end{align*}
In this case, we can even determine the particular real number that each function, $f$, is mapped to by the functional $F$.

\begin{align*}
    F[f(x;a)] &= \int_{x=-1}^{x=1} a x^2\, dx \\
    &=a \left. \frac{x^3}{3}\right|_{x=-1}^{x=1}\\
    &=\frac{2}{3} a 
\end{align*}
The action of the functional, $F$ is now clear: it maps every function of the form $f(x)=ax^2$ \emph{uniquely} to the number F[f(x;a)] = 2/3 a.  The mapping is unique because each value of $a$ defines a different member of the set of parabolas we have defined, and each such parabola is linked to a real number that depends only on $a$.  
\end{example}
\end{svgraybox}

According to our definition, functionals can also map \emph{complete functions} to \emph{new functions}.  An example here can be seen by defining the conventional moving-box average.  

\begin{svgraybox}
\begin{example}[A functional resulting in a new function: Moving average]

Suppose we define the following box function

\begin{equation*}
    B(y) =
    \begin{cases}
    2 & 0 \le y \le \tfrac{1}{2}\\
    0 &otherwise
    \end{cases}
\end{equation*}
Now, for any function $f$, we can define the following averaging operation given by a functional, $F$

\begin{equation*}
    F[f(x),x] = \int_{y=-\tfrac{1}{4}}^{y=\tfrac{1}{4}} f(x+y)B(y)\, dy
\end{equation*}
This functional has been constructed to (1) for any point $x$, take the (uniformly weighted) average of the function $f(x)$ between $x-\tfrac{1}{4}$ and $x+\tfrac{1}{4}$, and (2) assigns this new value to the point $x$.  (Note: this means that the function needs to be defined on the interval $x-\tfrac{1}{4}<x<x+\tfrac{1}{4}$).  Thus, the function assigns a uniformly weighted average of the function near $x$ to to the point $x$.  To make this a bit more clear, the following figures can help.  First, note that the averaging function is plotted in Fig.~\ref{avefunc}.  

\begin{centering}
\includegraphics[scale=.35]{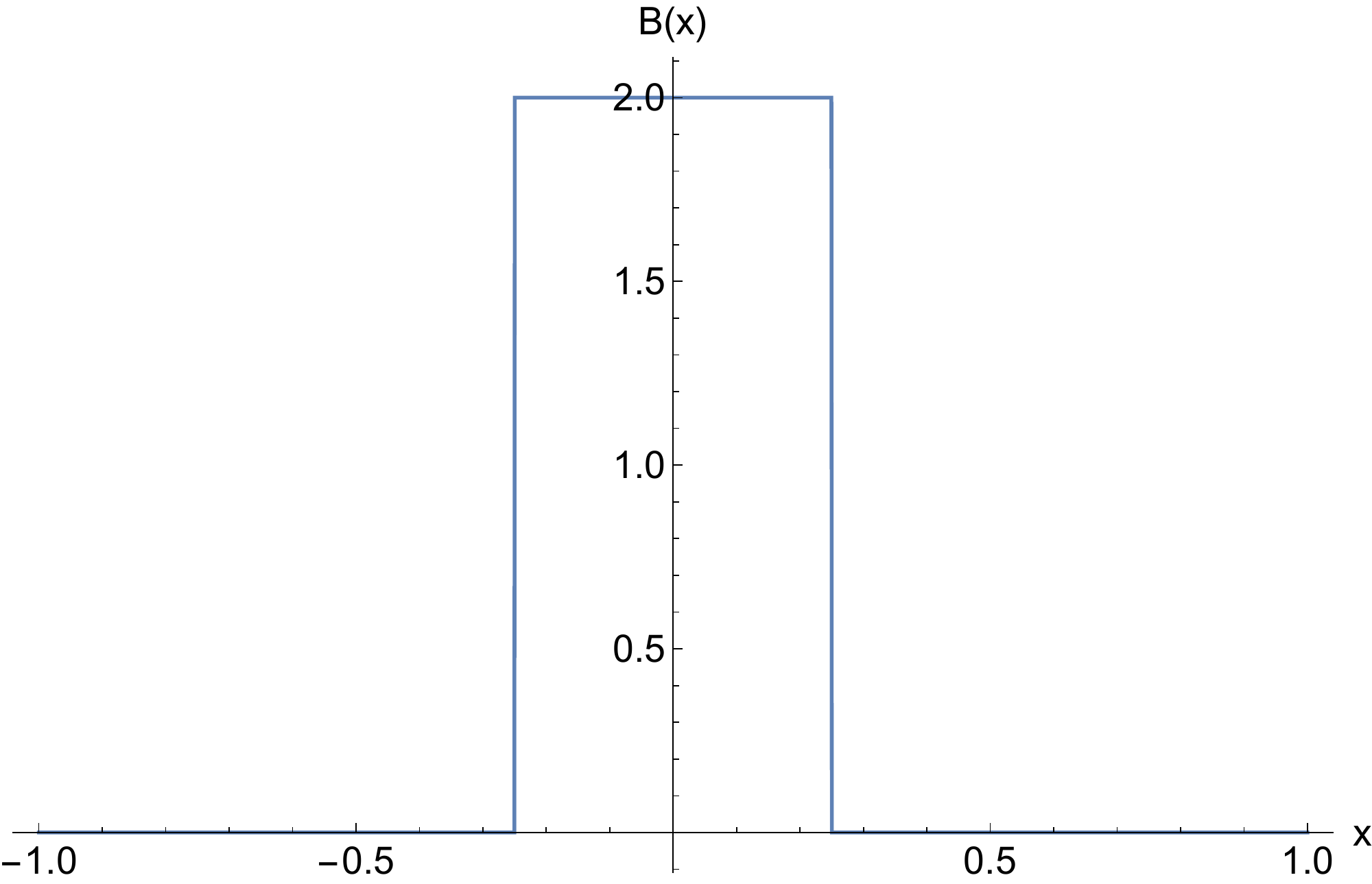}
\captionof{figure}{The averaging function; a simple uniform box with area equal to 1.}
\label{avefunc}  
\end{centering}

In Fig.~\ref{noisefunc}, we have plotted a roughly linear function that has been subjected to random noise.  Often, when one wants to remove nose from a function, this can be done by passing an averaging function over the noisy function to smooth out the random fluctuations.  In this case, we will be using the box function defined above.  

\begin{centering}
\includegraphics[scale=.35]{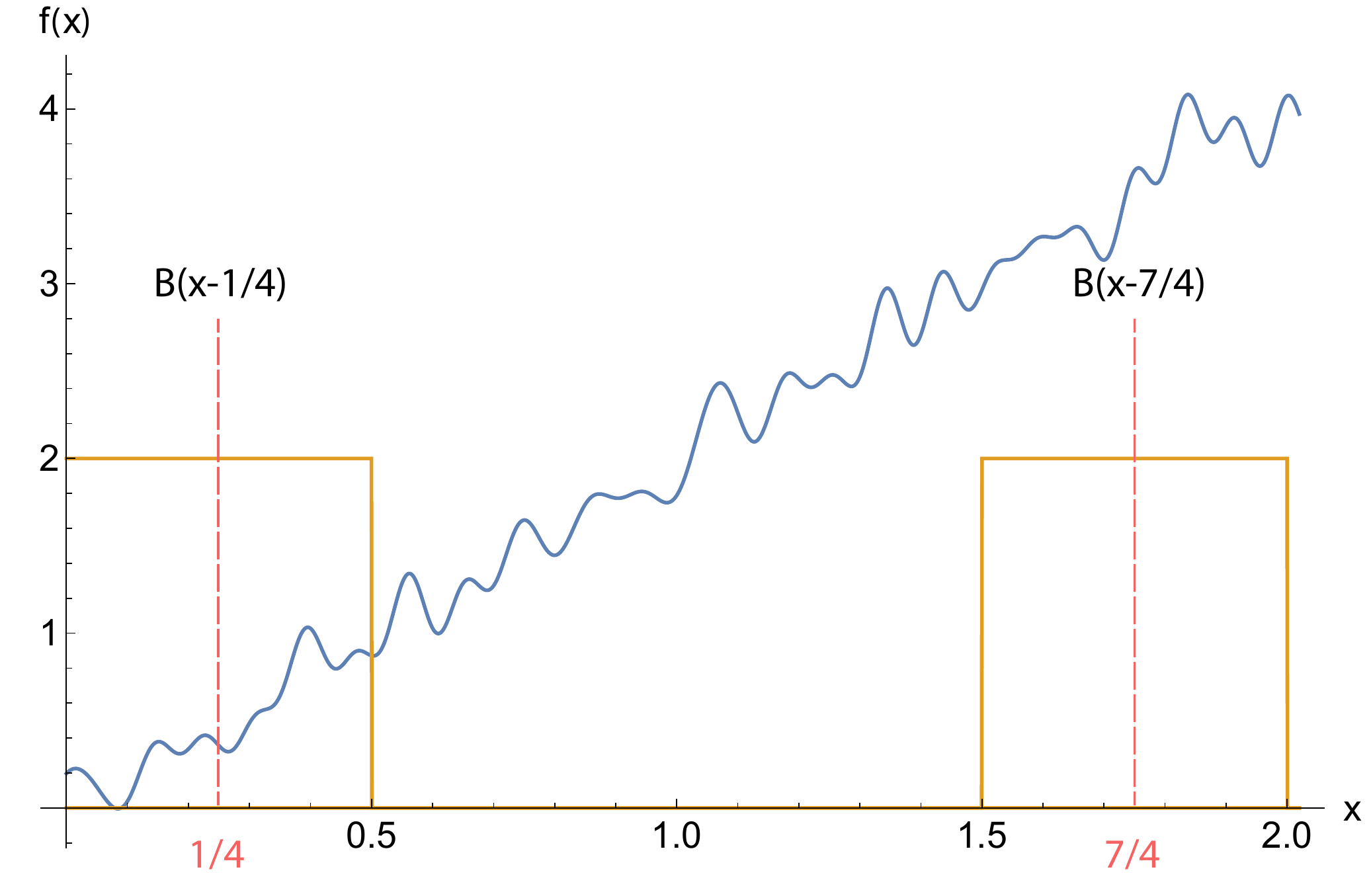}
\captionof{figure}{A noisy function (in blue) on the interval $x\in[0,2]$ to be smoothed. The box functions $B(x-1/4)$ and $B(x-7/4)$ are plotted for reference.  Note that the functional, $F(x)$ is only defined for $\tfrac{1}{4} \le x \le \tfrac{7}{4}$ for this case, because the box function would extend beyond the bounds of definition for the function otherwise.}
\label{noisefunc}  
\end{centering}

Upon applying the moving box average via the functional $F$ defined above, we have as a result a smoothed version of our function.  This smoothed function can be observed in Fig.~

\begin{centering}
\includegraphics[scale=.35]{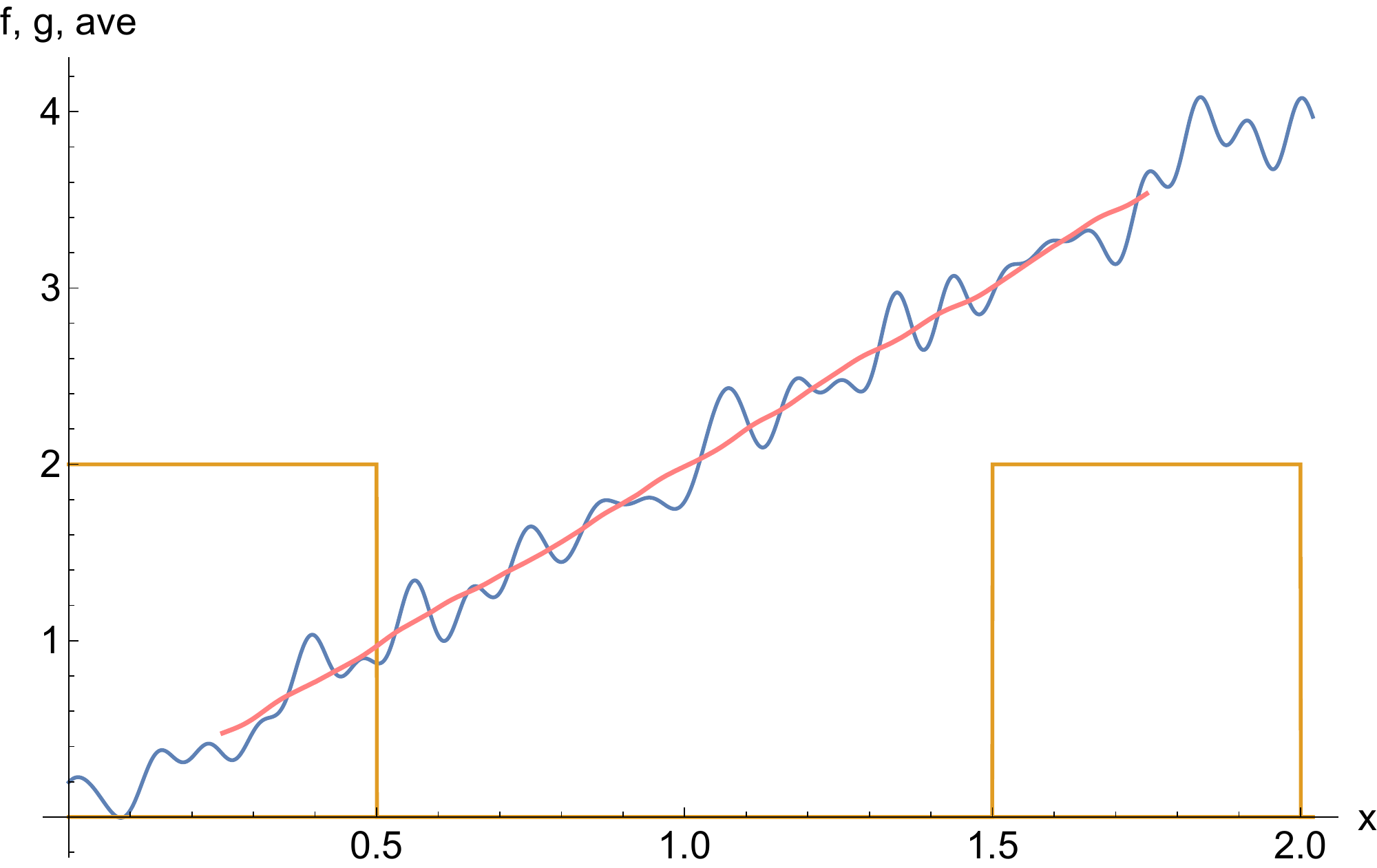}
\captionof{figure}{A noisy function (in blue) on the interval $x\in[0,2]$ to be smoothed. The line in red indicates the smoothed function.  Note that the smoothing is only defined between $\tfrac{1}{4} \le x \le \tfrac{7}{4}$, because it is only within these bounds that the averaging function is within the domain of the noisy function $f(x)$.}
\label{noisefunc2}  
\end{centering}

\end{example}
\end{svgraybox}

This second example is one where the functional specifically takes the form

\begin{equation}
    g(x) = \left.F[f]\right|_x = \int_{y=a}^{y=b} K(x,y) f(y) \, dy
\end{equation}
where here an alternate notation has been used ($\left.F[f]\right|_x$)to indicate independent variable of the \emph{resulting} functional. This should not be confused with the independent variable for the function $f$ (which, in this case, is represented by the variable $y$ under the integral sign).  Such integral transforms are used widely in applied mathematics.  In this use, the function $K(x,y)$ is usually called the \emph{kernel} \indexme{integral!kernel}\indexme{functional!kernel}\indexme{kernel} of the integral.  One of the most well-known uses of the functional defined by an integral with a kernel operation is that the the Fourier and Laplace transforms.  These two integral transform methods are presented in detail in later chapters.

\newpage
\section*{Problems}
\begin{enumerate}

\item  Use set builder notation to make the following statement: ``create the set of numbers $n$ in the integers that are less than the value three".  

\item Use set builder notation to make the following statement: ``create the set of numbers $n$ in the integers that are less than the value three and greater than the value of negative twenty-seven". 

\item The perfect squares are the squares of the natural numbers; let's assume that $0$ is a perfect square.  Using set builder notation, build an express for the set of all perfect squares (including zero).    

\item Suppose we would like to define the interval $I=(a,b]$.  Use set builder notation to develop a set expression for this interval in the form $I=\{x:$ ...

\item Suppose $I_1 = (0,1]$ and $I_2=[1,2]$.  Now define $I_3=I_1\cup I_2$  Use set builder notation to describe $I_3$ in the form $I_3=\{x:$ ...

\item Suppose $I_1 = (0,1]$ and $I_2=[1,2]$.  Now define $I_3=I_1\cap I_2$  Use set builder notation to describe $I_3$ in the form $I_3=\{x:$ ...

\item Suppose $I_1 = (0,1]$ and $I_2=(1,2]$.  Now define $I_3=I_1\cap I_2$  Use set builder notation to describe $I_3$ in the form $I_3=\{x:$ ...

\item The determinant for a $2\times2$ matrix is relatively straightforward to find.  Consider the matrix
    \[\twoform{A} = \left[ {\begin{array}{*{20}{c}}
  a_{11}&a_{12} \\ 
  a_{21}&a_{22} 
\end{array}} \right]\]
where $a_{ij}$ indicates (by convention) the value of the $i^{th}$ row and the $j^{th}$ column.

Now, do the following.  (i) Use row operations to make the matrix in upper triangular form, and (ii) compute the determinant from the product along the diagonal of the triangular matrix.  The result should be $\textrm{det}(\twoform{A})=a_{11}a_{22}-a_{12}a_{21}$.\\

\item As long as we have computed the determinant for a $2\times2$ matrix, we may as well learn about the solution to a set of two equations in two unknowns.  There is a method, known as \emph{Cramer's rule} that allows you to find the solution to \emph{any} set of linear equations as long as you can find the appropriate determinants.  It turns out that for an $n\times n$ system of equations, you need to find $n$ determinants.  In general, this is not so easy to do; for small matrices, however, it can be very convenient. 

For two equations in two unknowns, the coefficient matrix, $\twoform{A}$ is $2 \times 2$.  For a $2\times 2$ system of the form
\[ \twoform{A}\cdot \vec{x}=\vec{b} \]
or
\[ 
\left[ {\begin{array}{*{20}{c}}
  a_{11}&a_{12} \\ 
  a_{21}&a_{22} 
\end{array}}\right]
\cdot
\left[ {\begin{array}{*{20}{c}}
  x \\ 
  y 
\end{array}}\right]
=\left[ {\begin{array}{*{20}{c}}
  {b_1} \\ 
 {b_2} 
\end{array}}\right]
\]
Cramer's rule says that the solution can be found from the following determinants where the right-hand side has been substituted into the columns of $\twoform{A}$

\begin{align*}
    x&=\frac{1}{\textrm{det}(\twoform{A})}\left| {\begin{array}{*{20}{c}}
  {{b_1}}&{{a_{12}}} \\ 
  {{b_2}}&{{a_{22}}} 
\end{array}} \right|  &
y&=\frac{1}{\textrm{det}(\twoform{A})}\left| {\begin{array}{*{20}{c}}
  {{a_{11}}}&{{b_1}} \\ 
  {{a_{21}}}&{{b_2}} 
\end{array}} \right| 
\end{align*}
where both $\textrm{det}(\cdot)$ and $| \cdot |$ indicate the determinant.
Suppose we have the set of equations
\begin{align*}
    2 x + 4 y = 2 \\
    2 x + 2 y = \frac{3}{2}
\end{align*}
or, in matrix form
\[ 
\left[ {\begin{array}{*{20}{c}}
  2&4 \\ 
  2&2 
\end{array}}\right]
\cdot
\left[ {\begin{array}{*{20}{c}}
  x \\ 
  y 
\end{array}}\right]
=\left[ {\begin{array}{*{20}{c}}
  2 \\ 
  \frac{3}{2} 
\end{array}}\right]
\]
Use Cramer's rule to solve for $x$ and $y$.\\

\item Linearity is an important concept to be familiar with.  For the following operators, $\mathscr{L}(\cdot)$, determine if they are linear or nonlinear
    \begin{enumerate}
        \item $\mathscr{L}(x)=a x^2 + bx$
        \item $\mathscr{L}(x) = \sqrt{a^2 x} $
        \item $ \mathscr{L} = \frac{d^2}{dx^2}+ a \frac{d}{dx}, \mbox{ so that }~\mathscr{L}(f) = \frac{d^2f}{dx^2}+ a \frac{df}{dx} $
        \item $ \mathscr{L}(f) = \int_x f(x) dx $
    \end{enumerate}

\item Integration is done one of two ways: either using an indefinite integral, or by using a definite integral.  The fundamental theorem of calculus connects these for us.  One way to help remember the fundamental theorem of calculus is to compute the area under a curve using the two methods.  So, please try this.  For the function $f(x)=x$ on the interval $x\in[1,2]$ compute the area two ways.
\begin{enumerate}
    \item Compute the indefinite integral.  This will contain an unknown constant; call it $C$.  The value of $C$ is fixed when we specify which interval we want to compute the area for.  So, suppose that we want to compute the area for the following two intervals: (i) $0\le x\le 1$, and (ii) $0 \le x \le 2$.  Compute the area under the curve $f(x) = x$ for intervals (i) and (ii) by drawing a sketch of the function, and then using simple geometry to compute the area of the associated regions.  
    
    \item Now, compute the numerical value of the indefinite integral for this problem as follows.  You know the proper areas that each interval represents.  You have also computed the indefinite integral with unknown constant, $C$ in part (a).   For the intervals (i) and (ii), compute the value of the unspecified constant, $C$.  It should be the same for the two intervals.
    
    \item Note that we can get the area of the interval $1\le x\le 2 \Leftrightarrow x\in[1,2]$ by subtraction using the two areas above.  Draw this on a sketch, and also compute what the area is numerically (from your results in part (b)).
    
    \item Finally, compute the definite integral
    \begin{equation}
        \int_{x=1}^{x=2} x \, dx
    \end{equation}
    You should (kind of obviously perhaps at this point) get exactly the same value.  The fundamental theorem of calculus states that the area for some fixed interval is just the different between two indefinite integrals, as you have shown.
\end{enumerate}

\item Use integration by parts to compute the following integrals (i.e., find an answer in terms of functions no longer involving integrals) for $x\in[0,1]$.  Note: You may have to apply integration by parts more than once to get a final result!.
    \begin{enumerate}
        \item $f(x) = x \sin(\pi x)$
        \item $f(x) = x \cos(\pi x)$
        \item $f(x) = x^2 \sin (\pi x)$
        \item $f(x) = x^2 \cos(\pi x)$
    \end{enumerate}
    
\item For the functions $f$ and $g$ in Fig.~\ref{fig:3}, consider the domain ${\bf x}\in[0,1.5]$.  What are the \emph{ranges} of the two functions (approximate- estimate as best you can)? Are either of the functions piecewise continuous?  Are either of the functions piecewise smooth? \\

\item  The exponential function has the unusual property that its Taylor series converges everywhere on the real line.  Given that information, is the following function \emph{analytic} on the domain $D=[0,1]$?  

\begin{equation}
    f(x) = \exp(\pi (x-1)), 
\end{equation}

Why or why not?\\

\item {\bf Euler's rule}.  Euler's rule, given by 
\begin{equation*}
    e^{i x}=\cos x + i\sin x
\end{equation*}
is actually pretty easy to prove.  What you need to know to do that is the following (and that $x^0=1$ and $0!=1$ by definition)

\begin{align*}
     \exp(x) & = \sum_{n=0}^\infty \frac{x^n}{n!}=\frac{x^0}{0!} +\frac{x^1}{1!}+ \frac{x^2}{2!} + \frac{x^3}{3!} +\frac{x^4}{4!}+\ldots \\
    \cos x &= 1-\frac{1}{2!}x^2+\frac{1}{4!}x^4-\ldots \\
    \sin x &= x-\frac{1}{3!}x^3+\frac{1}{5!}x^5-\ldots
\end{align*}
Use these series, in addition to the properties of the imaginary number $i=\sqrt{-1}$ to prove Euler's rule.\\

\item Every polynomial is its own Taylor series.  To see this, determine the Taylor series for the polynomial

\begin{equation*}
    y(x) = 3 x^3 + 4 x^2 + 2 x + 5
\end{equation*}
Around the point $x=0$.  Start by first computing $y(0)$, $y'(0)$, $y''(0)$, $y'''(0)$, and $y^{\textrm{\scriptsize(IV)}}(0)$.  Then use the formula for the Taylor series to arrange these terms in the proper form, and compute the coefficients for each power of $x$.

\item Consider the following function, defined piecewise over the interval $0\le x \le 2$.
\begin{equation*}
f(x)=
    \begin{cases}
    0 \textrm{~for } 0\le x \le 1 \\
    x \textrm{~for } 1 < x \le 2 
    \end{cases}
\end{equation*}
 \begin{enumerate}
        \item Is this function continuous on the interval $x\in[0,2]$?
        \item Is this function $C^1$ continuous on the interval $x\in[0,2]$?
        \item Can this function be described as piecewise $C^1$ continuous?  Why or why not?
        \item Does this function have a conventional derivative defined for all $x\in(0,2)$?
    \end{enumerate}

\item Consider the following function, defined piecewise over the interval $0\le x \le 1$.
\begin{equation*}
f(x)=
    \begin{cases}
    2 \textrm{~for } 0\le x \le \tfrac{1}{2} \\
    -2 \textrm{~for } \tfrac{1}{2} < x \le 1 
    \end{cases}
\end{equation*}
 \begin{enumerate}
        \item Is this function continuous on the interval $x\in[0,1]$?
        \item Is this function $C^1$ continuous on the interval $x\in[0,1]$?
        \item Can this function be described as piecewise $C^1$ continuous?  Why or why not?
        \item This function is integrable.  To integrate it, simply break the domain of integration into two parts, one covering $x\in[0,\tfrac{1}{2}]$, and the other covering $x\in(\tfrac{1}{2},1]$.  
        
        Plot the integral of this function as a function of $x$, where $x$ represents the upper bound of the integral.  You will have to break the problem into two cases: one for the case that $x\le \tfrac{1}{2}$, and the other for $x>\tfrac{1}{2}$.
    \end{enumerate}

\item Using your own words (and graphics if you like) describe your understanding of the \emph{differential} $dy$.  How does it relate to the derivative $y'(x)$, where $x$ is the independent variable?

\item  Which of the following functions are \emph{homogeneous}?  Show your work.
    \begin{enumerate}
        \item $f(x) = x $
        \item $f(x) = x^2 \sin (\pi x)$
        \item $f(x) = x^\pi$
        \item $f(x)=\sqrt{5 x}$
        \item $f(x) = x^2 \sqrt{x}$
    \end{enumerate}
    
\item Does the following series converge? Which test for convergence can be used to establish the convergence or divergence of this series?

\begin{equation}
    P(n) = \sum_{n=1}^{\infty} \frac{1}{n^2}
\end{equation}

\item Is the following function analytic at $x=0$?  Is it continuous at $x=0$?  Assume that the domain is\\ $-\infty < x < \infty$.
\begin{equation}
f(x) = \exp\left(-\frac{1}{x^2} \right)
\end{equation}

\item Taylor series are central to what are called \emph{asymptotic} expansions.  For such expansions when the expansion is for values of $x$ that are much smaller than the radius of convergence, $R$.  These series often generate a very accurate representation of the function (in a limited domain) with just a few terms.

One frequently used example of an asymptotic series is for the function

\begin{equation}
    f(x) = \frac{1}{1+x}, ~ 0\le x < 1
\end{equation}
The radius of convergence for this problem is $|x|<1$. 

For this problem, do the following.
\begin{enumerate}
    \item Determine the general Taylor series for this function.  Note: This function has a radius of convergence $R<1$, and it does not converge at the boundary $R=1$.  In other words, this series converges for $|x|<1$.
    
    \item Compute the first five terms of this series explicitly (i.e., work out the Taylor polynomial for it).  Using Mathematica (or whatever software is convenient for you) plot the series approximation $f_5(x)$ and the original function $f(x)=1/(1-x)$ over the interval $x\in[0,1]$.  
    
    \item  This series converges for $|x|<1$.  In the figure below, the Taylor polynomials for $N=3,~4,~9,~10$, and $50$ are shown.  Regardless of the size of $N$, each polynomial gives a value of 1 or 0, depending on whether the final term is even or odd.  Thus, the series does not converge at $x=1$, no matter how many terms are computed.  \\
    
    It can be shown that the series \emph{does} converge for all $x$ less than 1.  Is the convergence of this series \emph{uniform convergence} or \emph{pointwise convergence}?  in the domain $x\in[0,1)$?  Provide a sentence explaining your choice.\\
    
        \begin{figure}[t]
    \sidecaption[t]
    \centering
    \includegraphics[scale=.4]{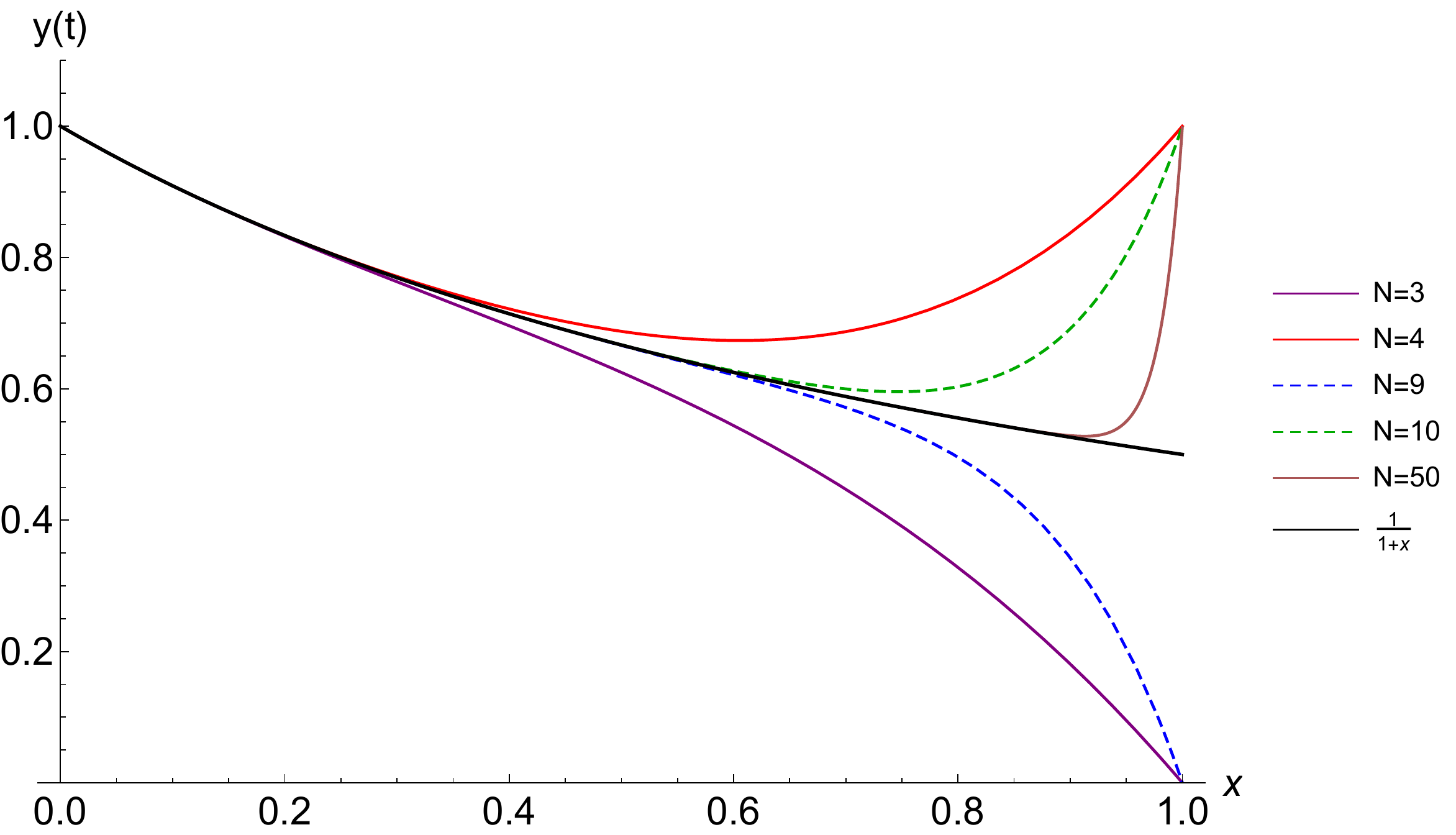}
    \caption{Several Taylor polynomials for $f(x)=1/(1+x)$ plotted on the interval $x\in[0,1]$.  }
    \label{one_over}       
    \end{figure}
\end{enumerate}

\item {\bf The Divide and Conquer method for finding the square root.}
The ancient Babylonians had a method to compute the square root of any positive real number, $a$.  It was an iterative algorithm that started with an initial \emph{guess} of the root, and then provided the method to find the exact root.  Mathematically, what they did was to develop the following \emph{implicit} relationship for the square root.
\begin{align}
    x^2&=a&&(by~definition) \nonumber\\
    2x^2&=x^2+a&&(add~x^2~both~sides) \nonumber\\
    2x&=x+\frac{a}{x}&&(divide~by~x) \nonumber\\
    x&=\frac{x+a/x}{2}&&(solve~for~x,~implicit~equation) \nonumber
\end{align}
Here, $a$ is the starting number, and $x$ is the square root being sought.

Given an initial guess, $x_{1}$, then a better estimate could be found by 
\begin{equation}
    x_2=\frac{x_1+a/x_1}{2}
\end{equation}
or, in general
\begin{equation}
    x_{i+1}=\frac{x_i+a/x_i}{2}, \qquad i=1,2,3,\ldots (or~i\in \mathbb{N})
\end{equation}

Use this method to develop a sequence of estimates for $\sqrt{2}$ (i.e., $a=2$), starting with the initial guess of $x_1=1$ .  You can use a calculator (or spreadsheet) to do the work, but please show the result at each step. Also compute the relative error at each estimate step $i$, $\epsilon_i = (estimate_i - actual)/actual$, where the value of $estimate_i$ is the value computed at the $i^{th}$ step of the algorithm above, and the value of $actual$ is the value computed by a calculator or spreadsheet.  Stop the algorithm when the relative error is less than $1/100$.\\ 
    
\end{enumerate}

%
\abstract*{This is the abstract for chapter 1.}

\begin{savequote}[0.55\linewidth]
``In the broad light of day mathematicians check their equations and their proofs, leaving no stone unturned in their search for rigour. But, at night, under the full moon, they dream, they float among the stars and wonder at the miracle of the heavens. They are inspired. Without dreams there is no art, no mathematics, no life."
\qauthor{---Michael Atiyah, English Mathematician and co-discoverer of the Atiyah–Singer index theorem (widely used in counting the number of independent solutions to differential equations).}
\end{savequote}
\begin{savequote}[0.55\linewidth]
``Classification of mathematical problems as linear and nonlinear is like classification of the Universe as bananas and non-bananas." 
\qauthor{---Author Unknown.}
 \end{savequote}


\chapter{First and Second Order Ordinary Differential Equations}\label{ODEChap}
%
\def\CHAP {chapter02_odes}
%
\section{Terminology}
In the study of ordinary differential equations (ODEs), there are a few definitions and some vocabulary that are helpful to establish up front.  Some of these are repeated from the review chapter.   The repetition of ideas when learning is not actually redundant!  Most of us do not learn by seeing something a single time.  The more fundamental and important an idea or a definition is, the more it may be worth repeating, ideally using slightly different perspective to appeal to different ways of thinking about the topic.

\begin{itemize}

\item \textbf{Domain.}\indexme{domain} A \emph{domain} is the set of possible values that can be selected for a particular application.  Domains can be represented as either a set of discrete values, or by continuous intervals.  Values not in the domain are not generally valid for the application in question.  As an example, suppose one conducted an experiment which collected four temperature measurements (whose values were represented by, for example, $T_1, T_2, T_3, T_4$ at for different times, $t_1, t_2, t_3, t_4$.  Then, the values $t_1, t_2, t_3, t_4$ represent the domain of the experimental results collected. \\
\item \textbf{Range.}\indexme{range} The \emph{range} is the set of possible values that is associated with the particular domain in an application.  As an example, the range of values from the temperature experiment described above would be the temperatures that were measured: $T_1, T_2, T_3, T_4$.  Each of these temperatures is uniquely associated with a single point from the domain.\\
\item \textbf{Independent variable.} \indexme{independent variable} \indexme{variables!independent} An independent variable is one in which represent \emph{causes} or \emph{inputs} to the mathematical description of a process.  These variables do not, in principle, depend on any other quantity-- in other words, one is free to select them from any valid domain.  In the temperature experiment example above, the independent variables are the times that were selected to make temperature measurements: $t_1, t_2, t_3, t_4$.  The particular times selected to make measurements are constrained only by factors external to the process itself.  These might be governed by convenience, or possibly constrained by physical considerations (e.g., one might not want to make 10,000 measurements over a 4 day period if system reaches a steady condition after 5 minutes).  Sometimes the independent variable is ``tagged'' with the list of dependent variables to keep the relationship clear.  So, in the example above, one might indicate that the temperature is considered to be a function of time by writing $T(t)$.  Note the relationship between the independent variable and the domain of a function.\\ 
\item \textbf{Dependent variable.} \indexme{dependent variable} \indexme{variables!dependent} The dependent variable is the one that is computed once the independent variable is specified.  In problems of practical interest, the dependent variable is the variable of interest, representing the physical quantity that one wants to predict. Note the relationship between the dependent variable and the range of a function.\\

\item \textbf{Derivative.}\indexme{derivative!definition}  A derivative is the differential rate of change of a dependent variable variable as one of its independent variables changes.  If there is a single dependent variable with a single independent variable, then the derivative corresponds to what we conventionally think of as the slope of a simple function when plotted on two axes (with the vertical axis representing the independent variable, and the horizontal axis representing the dependent variable).  In calculus, we all learned that the derivative was given by
\begin{align}
  \frac{{dy}}{{dt}} & \equiv \mathop {\lim }\limits_{\Delta t \to 0} \frac{{\Delta y}}{{\Delta t}} \hfill \\
   &\equiv \mathop {\lim }\limits_{\Delta t \to 0} \frac{{y(t + \Delta t) - y(t)}}{{\Delta t}} \hfill
\end{align}
Note that although the notation \smash{$\tfrac{dy}{dt}$} is in common usage, it is also very common to see the notation for the derivative expressed using the \emph{prime} notation (see \ref{note01} below), where the derivative is expressed as $y'(t)$.  \\

The \emph{order of a derivative}\indexme{derivative!order} is just a count of how many times the derivative of a function has been taken.  Thus, when a function has its derivative taken twice, this is denoted $y''(t)$ or \smash{$\tfrac{dy^2}{dt^2}$}, and this is also called a \emph{second-order derivative}.\\

\item \textbf{Ordinary differential equation.}\indexme{ordinary differential equation}   An \emph{ordinary} differential equation (ODE) is one in which there is only one independent variable.  All derivatives appearing in the equation are taken with respect to the single independent variable.  The term ordinary is used in contrast with the term partial differential equation which may be with respect to more than one independent variable.  Note that one can have more than one independent variable in the case of coupled systems of equations.  For example, consider the conventional multiplicative reaction rate for fully mixed chemical reactors. Suppose there are two chemical species, $y_1$ and $y_2$ who must come in contact for the reaction \smash{$y_1 ~{\mathop\to \limits^{k_1}} ~y_2$} to happen.  The reactions can be represented by
\begin{align}
\frac{d y_1}{dt} &=-k_1 y_1 y_2 \\
\frac{d y_2}{dt} &= k_1 y_1 y_2 
\end{align}
Although there are \emph{two coupled} differential equations here, there is still only \emph{one} independent variable.  Thus, this set of equations are still ODEs. \\
\item \textbf{Order of an ODE.}  The \emph{order} of a differential equation is the order of the highest derivative in the equation.  Examples.
\begin{align}
y'' + 2 y' &= 0 &(second-order ODE) \\
y''' + 3y &= 42 &(third-order~ODE) \\
y^{(iv)} &= f(t)  &(fourth-order~ODE) \\
\frac{dy^2}{dt^2} +\left(\frac{dy}{dt} \right)^3 &= 0 &(second-order~ODE)
\end{align}
Note that in many texts, the prime notation is altered using Roman or Arabic numerals in parentheses to signify derivatives of order four or more.  Thus $y''''(t) = y^{(iv)}(t)=y^{(4)}(t)$.  You can see why this might be the case- writing a $7^{th}$-order derivative by $y'''''''(t)$ starts to get to be a bit ridiculous.\\

\item \textbf{Homogeneous ODE.}  A \emph{homogeneous} \indexme{ordinary differential equation! homogeneous} ODE is one in which there are no terms that do not involve the independent variable. For example, examine the following
\begin{align}
y''+3y' +7y &= 0& &homogeneous \\
y''+3y' +7y &= 5t& &nonhomogeneous \\
y''(t)+3t y'(t) +7 t^2 y(t) &= 0& &homogeneous \\
y'(t) +7 y(t) &= \pi& &nonhomogeneous \\
\end{align}
More generally, a linear homogeneous equation of order $n$ takes the form
\begin{equation}
a_n(t) y^{(n)}(t) + a_{n-1}(t) y^{(n-1)}(t) + \ldots + a_1 y(t) = 0
\end{equation}
where the $a_n(t)$ are coefficients of the ODE (which may, themselves, be functions of the independent variable $t$), and $y$ is the independent variable.\\

\item \textbf{Linear ODE.}  A \emph{linear ODE}\indexme{ordinary differential equation!linear} is one which satisfies two properties for the {\it homogeneous} part of the ODE:
\begin{enumerate}
	\item If the function $y(t)$ is a solution to the ODE, then multiplying the solution by any constant, $\alpha$, generates a function $y_2$ that is also a solution; in other words $y_2(t) =\alpha y(t)$ is a solution.  Somewhat confusingly, this aspect by itself is sometimes known as algebraic \emph{homogeneity}; this was defined and reviewed in section \ref{sec:functions}. 
	
    \item If one has two linearly independent solutions to the ODE, call them $y_1$ and $y_2$, then the linear combination of the two is given by $y(t) = y_1(t) + y_2(t)$; this linear combination is also a solution to the ODE.  This property is sometimes called \emph{additivity}.
\end{enumerate}

\item \textbf{Nonlinear ODE.} Any ODE that is not linear is nonlinear. \indexme{ordinary differential equation!nonlinear} 

\item \textbf{Analytic function.}  An analytic function is one that has a local representation as a convergent power series.  In other words, for any point given by the independent variable $x$, the function $f(x)$ is given by a power series that converges over some radius of convergence.

\item \textbf{Transcendental function.} \indexme{function!transcendental} A transcendental function is an analytic function that can not be specified by a (finite) polynomial.  In other words, it cannot be expressed exactly by a finite sequence of algebraic operations.  Such functions are often approximated by the appropriate truncation of an infinite series representation.

\item \textbf{Ancillary condition}. An ancillary condition is one where the \emph{dependent variable} (or one of its derivatives) given at a specified value of its \emph{independent variables}.  

\item \textbf{Initial condition}. \indexme{ordinary differential equation!initial condition}  An initial condition is an ancillary condition that is given at a specific time; this is assumed to be the initial or starting time for the process.  As an example, to find the concentration during a reaction, one would need to know the initial concentration at some time.  An ODE could then be used to determine its subsequent locations as a function of the time elapsed since the initial condition.

\item \textbf{Boundary condition}\indexme{ordinary differential equation!boundary condition} A boundary condition is the spatial equivalent to an initial condition.   As an example here, the steady heat flow in an insulated rod will lead to a temperature distribution that is a function of the ancillary conditions specified for the two ends; these are the boundary conditions.

\end{itemize}

\begin{svgraybox}
\begin{note}[Notation for derivatives.]\label{note01}

\noindent There are a number of notations for derivatives in use, and this comes from the long history of the evolution of calculus at nearly the same time but by different people.  The two most famous (and rightfully so) names in the development of calculus are Newton and Leibniz (also spelled \emph{Leibnitz}); both developed the fundamental ideas about calculus in the 1660s.  Newton preferred to indicate derivatives by a dot over the variable; thus, a derivative of a function $y(t)$ with respect to $t$ would be indicated by $\dot y(t)$, or, when the independent variable was obvious in the formulation, just $\dot y$.  Leibniz adopted a notation that contained more intuitive content, but was a bit more unwieldy; he used the notation \smash{$\tfrac{dy}{dt}$} to indicate the derivative. One of the advantages to this notation was that it was \emph{always} clear what the independent variable was!  Despite the fact that this was only notation, there was actually a significant rivalry between mathematicians in England and those in Europe over this topic.  To further complicate things, other notations cropped up, including the notation of the French mathematician Lagrange (who introduced the prime to indicate derivatives, as in $y'(T)$) and the Swiss mathematician Euler (who introduced the notation $Dy(t)$ to indicate the derivative). Much of this problem came to some (minor) resolution in the early 1830s, when at the University of Cambridge (England) a group calling themselves the \emph{Analytical Society} was formed specifically to promote the Leibniz notation for the derivative.  Although this may seem somewhat absurd out of context, there was a real point to this transition.  Because the Leibniz notation was being used throughout Europe, students in England were not able to easily understand the newest work being done on European mathematicians.  This was seen as putting England behind, and thus the Analytical Society was born.  The switch to the Leibniz notation helped foster better communication between England and Europe.  

In the present day, one finds both the Lagrange prime notation, and the Leibniz notation used frequently.  In addition, in some disciplines such as continuum mechanics and physics one still finds the dot notation used as a reference to the influence that Newton had on the development of these areas.
\end{note}
\end{svgraybox}

\section{Introduction}

As soon as one defines the derivative, it turns out that one also poses a \emph{differential equation}.  A differential equation is exactly what it sounds like:  it is a mathematical expression that involves the derivative of some unknown function rather than, say, a purely algebraic expression of some unknown function.  The word \emph{ordinary} is used as a prefix to describe differential equations that contain exactly one independent variable.  Differential equations that involve two or more independent variables are called \emph{partial differential equations}; those will be the subject of later chapters.

It was the German mathematician Gottfried Wilhelm Leibniz of calculus fame who is given credit for formally writing down the first differential equation with solution in 1675 \citep{sasser1992history}; however, Isaac Newton also had examples of solutions as early as 1671, although these were not published until 1736 \citep{krishnachandran2020}.  Leibniz's realization was essentially the one explained by the first sentence of this introduction.  Suppose we have a function $y(x) = 1/2 x^2$.  The derivative of this function is given by the equation

\begin{equation}
    \frac{dy}{dx} = x
    \label{firstode}
\end{equation}
While thinking about this in the ``forward" direction represents taking the \emph{derivative} of the function $y(x) = 1/2 x^2$, there is a second interpretation.  That interpretation involves viewing Eq.~\eqref{firstode} as an equation with an unknown function $y$ and solving it.  Or, put as a question, the problem is ``what function $y(x)$ exists whose derivative is specified by the quantity $dy/dx = x$ everywhere in its domain?" Of course, we know the answer to \emph{this} question, because we developed it from knowledge of $y(x)$ in the first place.

While viewing the derivative problem in reverse may seem trivial, the change in perspective is actually  dramatic.  Once one admits the idea that equations expressed in terms of derivatives might be somehow solved to determine whole functions, then a world of possibilities for new kinds of equations (and questions) is suddenly opened up.  Instead of simply knowing pairs of functions and their derivatives, one could now formulate problems strictly in terms of the rate-of-change of something physical, and then go ``seeking" the solution without necessarily knowing it beforehand.  This is a very powerful skill when applied to problems in science and engineering, and it has particular applications in problems where the rates of change of various quantities are given via some kind of conservation law.

Our review of ODEs will cover only the following: (1) Linear and (some) nonlinear first-order ODEs, and (2) linear second-order ODEs with constant coefficients.  The more general case of linear second-order ODEs with nonconstant coefficients is an interesting topic (and we will encounter those when studying Sturm-Liouville problems), but the topic involves series solutions; these solutions can get a bit tedious without adding substantially to developing new concepts for better understanding ODEs, so the choice has been made to omit them in this review.  With a basic understanding of the solution to linear second-order problems with constant coefficients, motivated students can learn to handle problems with nonconstant coefficients fairly readily with a little self-study of the topic.

\section{First-Order ODEs}

First order ODEs are those which contain derivatives of order one as the highest order of derivative in the problem.  The most general form for a first-order ODE is (in the independent variable $x$) \indexme{ordinary differential equation!first order}

\begin{align}
    \frac{d y}{dx} &= f(x,y) \\
    \intertext{or, equivalently}
    y'(x) &= f(x,y)
    \label{first_order}
\end{align}
Here, the unknown $y(x)$ is a function.  Thus, if one were to translate Eq.~\eqref{first_order} in plain English, it would be the statement ``$y$ is the function whose first derivative with respect to the independent variable $x$ is the function f(x,y) everywhere in the domain of $x$."  Because $f(x,y)$ can be any (reasonable) function of $y$, the equation as a whole can be linear or nonlinear.  For example, if we chose $f(x,y) = \alpha y(x) +\beta y^2(x) +x$, we can determine (using the conventional tests for linearity) that the problem is nonlinear for this choice of $f(x,y)$.

\subsection{Checking a Solution to a Differential Equations}

First-order differential equations appear in many applications that are familiar to scientists and engineers.  For example, the classical first-order reaction problem for a  fully-mixed batch reactor is an example.  If $c$ is the concentration, then we think of the time-history $c(t)$ for the reactor as being specified by an initial configuration (initial condition) $c_0$, and a rate expression.  Together these form the description of the system taking the form
\begin{align}
\frac{d c(t)}{dt}&=-k_1 c(t) \label{first1.1}\\
c(0)&=c_0 \label{first1.2}
\end{align}
In this expression, $k_1$ is the \emph{first-order reaction rate constant}.  The system of equations given by Eqs.~(\ref{first1.1})-(\ref{first1.2}) represents a first-order differential equation for the mass balance of the chemical species $c$ as a function of time $t$.  Equation (\ref{first1.2}) is often called an \emph{initial} condition when the problem has time as an independent variable, or a \emph{boundary} condition when space is the independent variable.  For our purposes, we will sometimes simply call them \emph{ancillary} conditions, because these expressions provide additional data about the problem that must be determined from considerations that come from outside of mathematics.  We will discuss that in additional detail later.

Note that, by convention, the variable $c$ is called the \emph{dependent} variable because its value is viewed as being a function of time. Similarly, $t$ is called the \emph{independent} variable because its values can be selected freely. Sometimes the dependent variable is ``tagged'' with the independent variables to keep this dependence clear; hence, one might write $c(t)$ to assure that the relationship between $c$ and $t$ is kept clear. 

Although we have not yet discussed how to find solutions, one might recall that the solution to this problem is an exponential function in time
\begin{equation}
c(t)=c_0 exp\left(-k_1 t \right) \label{solution1.1}
\end{equation}
We can check to see that this is a solution by noting the following.
\begin{enumerate}
\item Taking the derivative of both sides of Eq.~(\ref{solution1.1}) and then substituting the result for both the derivative and the solution back into Eq.~(\ref{first1.1}), we can determine that the equation is \emph{valid}.
\item We also note that at $t=0$, we have that $exp(0) = 1$, so that $c(0)=c_0$ and thus meets the initial condition.  This verifies that the solution given by Eq.~(\ref{solution1.1}) is both a solution, and meets the ancillary conditions required for the solution.
\end{enumerate}
Below, the steps in this process are detailed.
\begin{svgraybox}
\begin{example}[Determining if a given solution is a valid one]
Often you will want to check to see if a proposed solution is actually a valid one-- in other words, that it meets both the ancillary conditions (inital or boundary) and the ODE itself.  The solution may be one that is given to you, or one that you derive yourself.  \emph{When working problems, it is always a good idea to check your solutions following the steps below to determine that your solution is correct!}.  The steps in checking the solution are the same regardless of the ODE and ancillary conditions.  For the problem of the first-order reaction with a single specified initial condition, the steps are as follows.

\begin{itemize}
\item \textbf{Determine the necessary derivatives for the solution}.  In our example, we have a first order equation, so we need derivatives of only order 1.  
\begin{align}
c(t) &= c_0 exp(-k_1 t) \\
c'(t) &= - k_1 c_0 exp(-k_1 t)
\end{align}

\item \textbf{Substitute the solution and its derivatives into the ODE}.  Substituting the solution and its first derivative (computed above) into the ODE, we find

\begin{equation}
- k_1 c_0 exp(-k_1 t)=-k_1 c_0 exp(-k_1 t) 
\end{equation}
Clearly, the two sides of this equation are equal, and this verifies that the proposed solution is indeed correct.\\

\item \textbf{Determine if the ancillary conditions are met}.  In this case, we have only one ancillary condition (an initial condition).  Substituting the value of the independent variable at the time specified in the initial condition ($t=0$), we find
\begin{align}
c(t=0) &= c_0 exp(-k_1 (0)) \\
c(t=0) &= c_0
\end{align}
Clearly, the ancillary condition is met for the time specified.  Since ODE is validated and the specified condition is met, then Eq.~(\ref{solution1.1}) is validated as being the particular solution we are looking for.
\end{itemize}

\end{example}
\end{svgraybox}

\subsection{Directly Integrable First-order ODEs}\indexme{ordinary differential equation!integrable}

Integrable linear first-order ODEs represent an interesting case because of their direct link to calculus. You may recall from your first course in calculus the much lauded (but rarely remembered) \emph{fundamental theorem of calculus}.  There are a number of presentations for the theorem; two of them are given below.

\begin{theorem}[Fundamental theorem of calculus]
Suppose $f$ is a continuous function on any open interval $I$ of the real line.  Then, there exists a function $F$ defined by 
\begin{equation}
F(t) = \int_x f(x) \,dt +C
\end{equation}
called the \emph{antiderivative} or \emph{indefinite integral} of $f$ on the interval $I$.  Here, $C$ is usually called a \emph{constant of integration}.  The function $F$ has the property that

\begin{equation}
F'(t) = f(t)
\label{FOC}
\end{equation}
for all values of $t$ in the interval $I$.
\end{theorem}
The interesting thing to note here is that Eq.~\eqref{FOC} is actually a first-order ODE!  Therefore, one interpretation of the fundamental theorem of calculus is that it defines the relationship between the integral and the derivative by solving the first-order ODE given by Eq.~\eqref{FOC}.

The primary feature of the fundamental theorem of calculus is that it establishes the idea that differentiation is the inverse operation to integration (hence, the term \emph{antiderivative} is often used for indefinite integrals).  In a more applied context, there is a corollary that is perhaps more familiar; this corollary relates the indefinite integral to the definite integral (i.e., the area under a curve in a closed interval).  

\begin{corollary}[The definite integral]
If $f$ is a continuous function on the closed interval $[a,b]$, and $F$ is its antiderivative, then the area under the function $f$ in the interval is given by a quantity called the \emph{definite integral}, and it is defined by
\begin{equation}
\int_{t=a}^{t=b} f(t)\, dt = F(b)-F(a)
\end{equation}
\end{corollary}

With this tool in hand, we are now ready to consider the solution to the first case of ODEs.  The simplest first-order ODEs are linear, directly integrable problems.  For these problems, the function $f$ depends on \emph{only} the independent variable.  We use the terminology $f(y,t)=Q(t)$

\begin{equation}
\frac{d y(t)}{dt}=Q(t)
\label{1storderlinear}
\end{equation}
Note that, although we are using the independent variable $t$ here, this variable could stand for anything (not necessarily time!).  In applications, we will often use the variable $t$ when dealing with problems where the independent variable is meant to represent time, and the variables $x$ in problems where it represents space.  However, in the study of ODEs it is understood that there is only \emph{one} independent variable, and the particular symbol used to represent that variable has no inherent significance on its own.

Solving Eq.~(\ref{1storderlinear}) is a straightforward application of the fundamental theorem of calculus.  Dividing through by $P(t)$ and integrating both sides of this equation, we find
\begin{equation}
\int\frac{d y(t)}{dt}\,dt=\int Q(t)\,dt
\end{equation}
The fundamental theorem of calculus then immediately gives the result
\begin{equation}
y(t)=\int Q(t)\,dt +C_1
\label{linearsolution}
\end{equation}
Where here, we have chosen to put the constant of integration, $C_1$, on the right-hand side of the result (N.B., the sign of the arbitrary constant at this juncture is somewhat irrelevant; when it is evaluated on the basis of an ancillary condition, the correct sign will be imposed).

\begin{svgraybox}
\begin{example}[First-order separable equations]\label{skier}
The conventional problems of computing the position, velocity, and acceleration of an object are usually expressed by systems of first-order, separable ODEs.  Recall from your physics classes that, if position is expressed as a function of time, $x(t)$, we have:
\begin{equation*}
   v(t)=\frac{dx(t)}{dt}\qquad a(t) = \frac{dv(t)}{dt} = \frac{d^2 x(t)}{dt^2}
\end{equation*}

We can use these relationships to solve interesting dynamical problems.  Consider the following.  Suppose we want to determine the velocity as a function of time for a skier starting from a standstill, and accelerating down slope as shown in the figure.  We also want to know how long in seconds it takes to get to the bottom of the slope, and the velocity of the skier when they get there.

{\bf Assumptions.}  It is always a good idea to \emph{list your assumptions} when solving problems.  That way, it is clear under what circumstances the result you obtain is valid.  It also provides you with something to look back on for problem improvement/refinement should your result come out less reasonable than you would have hoped.  For this problem, we will make some \emph{very rough} approximations that may not be altogether reasonable.  These are (1) the resistance between the skis and the snow will be neglected, and (2) air resistance will be neglected.  Note: if we wanted to include such effects, we might hope to find an ``overall" drag term that takes the form $drag~force=\alpha v^2(t)$.

To start the problem, we do a simple force balance on the skier, noting that $v(t=0)=x(t=0)=0$.  The acceleration in the down slope direction is found to be

\begin{align*}
    F&=m a(t) \Rightarrow a(t)=F/m\\
    \intertext{But, we can compute the force to be}
    F&=m g \sin \theta\\
    \intertext{Combining these, we have an expression for the acceleration}
    a(t) = F/m = g \sin \theta
\end{align*}

Note that, in this case, the acceleration is actually a constant.  To find the velocity, we need only integrate.  This integration could be done as a definite integral (where the bounds of integration would incorporate the initial condition), or as an indefinite integral (where the integration constant would be determined by using the initial condition as a constraint).

\begin{centering}
\includegraphics[scale=.55]{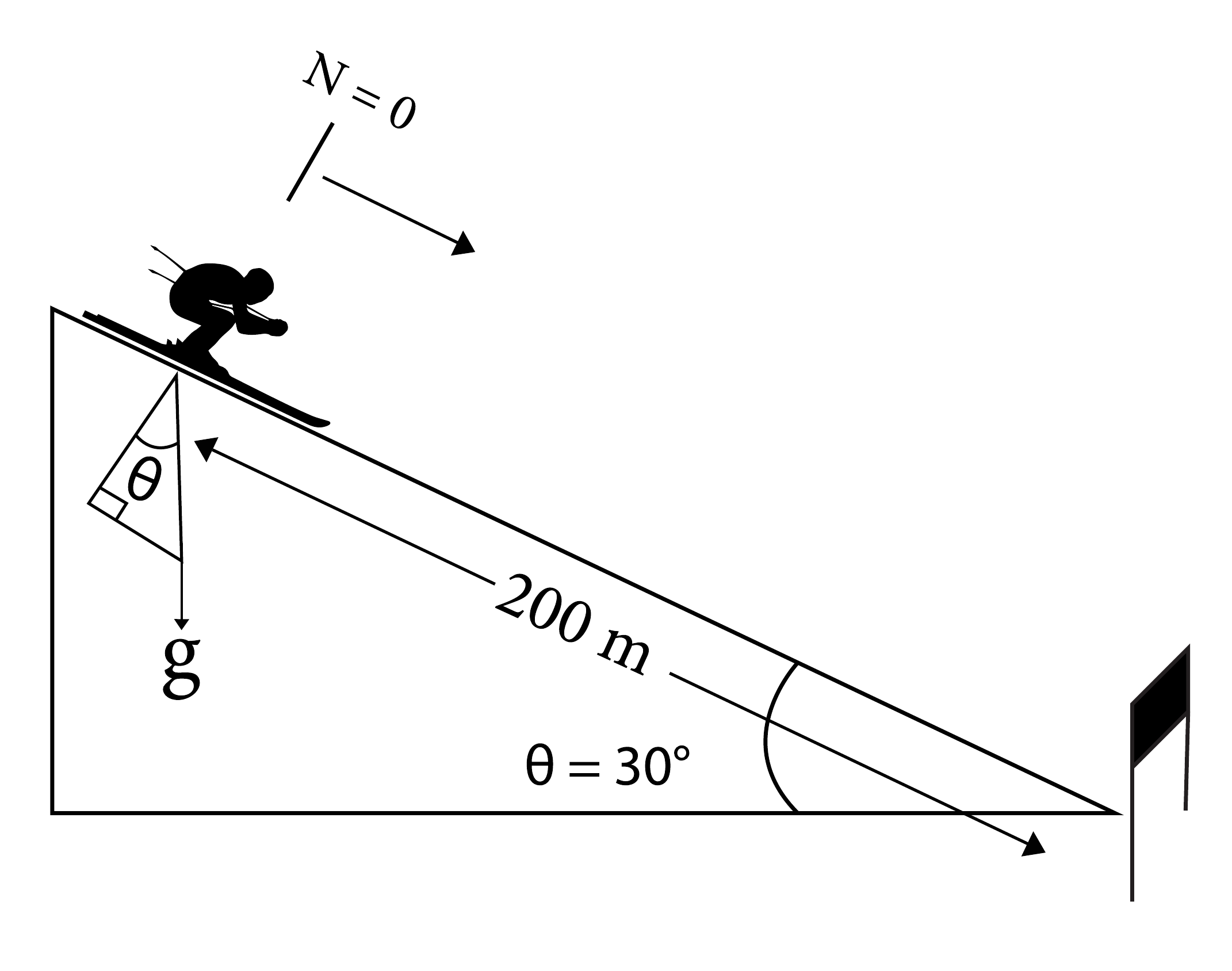}
\captionof{figure}{A skier.}
\label{fig:ski}  
\end{centering}
\vspace{5mm}
\begin{align*}
    \textrm{Recalling} ~~~a = \frac{d v(t)}{dt}\\
    \intertext{Upon substituting for $a$}
    \frac{d v(t)}{dt} &= g \sin \theta\\
    \int_t \frac{d v(t)}{dt} dt &= \int_t g \sin \theta dt\\
    v(t) &= g t \sin \theta +C_1\\
    \intertext{Using the fact that $v(0)=0$, then $C_1=0$.}
    v(t) &= g t \sin \theta
\end{align*}
The position as a function of time is found by a similar integration.
\begin{align*}
    \textrm{Recalling} ~~~v(t) = \frac{d x(t)}{dt}\\
    \intertext{Upon substititing for $v(T)$}
    \frac{d x(t)}{dt} &= g t \sin \theta\\
    \int_t \frac{d x(t)}{dt} dt &= \int_t g t \sin \theta dt\\
    x(t) &=  \tfrac{1}{2} g t^2 \sin \theta +C_2\\
    \intertext{Using the fact that $x(0)=0$, then $C_2=0$.}
    x(t) &= \tfrac{1}{2}g t^2 \sin \theta\\
    \intertext{or, solving for $t$}
    t &= \left(\frac{x(t)}{\tfrac{1}{2}g  \sin \theta}\right)^{\tfrac{1}{2}}
\end{align*}
With this last expression, we can solve for the time given the final position $x(t_f) = 200$m.  Using $g=9.81$ m/s, the result is $t=9.0$ s.  From the expression for $v(t)$, we can find the velocity at the end of the run to be $44$~ m/s (equivalent to 99 miles per hour!).

Certainly, this would be the absolute upper limit that one could reach on 200 meters of ski slope pitched at $30^{\circ}$.  Note that a $30^{\circ}$ slope is equivalent to 67\% slope; this actually quite steep!
\end{example}
\end{svgraybox}

\subsubsection{First-Order Separable Equations}\indexme{ordinary differential equation!separable}
%
As a reminder, the most general form for a first-order ODE is

\begin{equation}
    \frac{dy}{dt} = f(y(t), t)
\end{equation}
where $f(y(t), t)$ can be \emph{any} linear or nonlinear function of $y(t)$.  In general, there is no solution to all problems of this form.  However, for a reasonably common subset of this equation, solutions do exist; for example, if $f$ is a linear function of $y$ then solutions are certain to exist.  For example

\begin{equation}
    \frac{dy}{dt} = y(t)
\end{equation}
is a linear problem.  It represents the statement that we are looking function whose derivative is the same as the function itself.  A little thought will indicate that the exponential is the only such function (i.e., $y(t) = e^t$).  

This form is called a \emph{separable} first-order ODE, and it occurs when the function $f$ can be expressed as the product of two functions, one explicitly involving only $t$, and the other explicitly involving only $y$.  In other words, the problem takes the form

\begin{equation}
    \frac{dy}{dt} = r(y)Q(t)
    \label{separableeq}
\end{equation}
In general separable equations like this can be linear or nonlinear, and solutions for both cases are possible assuming that we can complete the necessary integration.  First-order ODEs are one of the few areas of mathematical modeling where one can obtain solutions for relatively general conditions.  Most physical systems actually contain some nonlinearity.  It is only our models that approximate such systems as being linear.  Nonlinear problems are the rule rather than the exception. 

Nearly all discussions of problems of this form write it in the following way, where $R(y)=1/r(y)$.

\begin{equation}
R(y) \frac{d y}{dt}=Q(t)
\label{1stordernonlinear}
\end{equation}

There are few nonlinear problems where a general solution exists.  Problems of the form of Eq.~\eqref{1stordernonlinear} do have solutions in general, and therefore represent a really unusual situation in mathematics!  A few notes are in order here.  First, equations of the form of Eq.~\eqref{1stordernonlinear} can, in principle, be linear (take $R(y) = 1$) or nonlinear; regardless, the solution method discussed in the following is valid.  Second, in the material that follows, it is important to be aware of potential confusion between independent and dependent variables.  One can compute an integral that is written in terms of dependent variables.  This can happen when a transformation of variables allows an integration of the dependent variable to be expressed entirely in terms of the dependent variable.  

This discussion can be put in context with an example solution for a particular nonlinear first-order ODE.  The general (abstracted) case will be developed following the example.

\begin{svgraybox}
\begin{example}[First-order nonlinear separable equations]\label{sepnon}
Consider the following nonlinear ODE (offered without a physical motivation).

\begin{equation}
    y^2 \frac{dy}{dt} = t^2
    \label{examplenon}
\end{equation}
This expression is clearly nonlinear in $y$.  Just as a reminder about testing for nonlinearity: the two required properties are (1) additivity and (2) homogeneity as discussed in \S \ref{lineop}; failing either one of these tests indicates nonlinearity.  If we substitute $y=\alpha y_1(t)$ into our problem we find

\begin{align}
\intertext{\it Checking for nonlinearity}
    (\alpha y_1)^2 \frac{ d(\alpha y_1)}{dt} &= t^2 \\
   \alpha^2 y^2_1 \alpha \frac{dy_1}{dt} & = t^2 \\
   \alpha^3 y^2_1 \frac{dy_1}{dt}  = t^2
\end{align}
thus, this problem fails homogeneity and is nonlinear.  Regardless, the problem is a separable one because it is in the form of \eqref{separableeq}.  To solve it, we start by integrating both sides with respect to $t$.  This gives

\begin{align}
    \int y^2 \frac{dy}{dt} dt = \int t^2 \,dt
    \label{integralsep}
\end{align}
In Section \ref{riemann}, the differential was defined.  There, we showed that it was appropriate to make the following correspondence

\begin{align}
    dy = y'(\tau) d\tau \\
    \intertext{or, equivalently}
    dy = \left(\frac{dy}{d\tau}\right) d\tau
    \label{differential2}
\end{align}
while it is tempting to think of this as if the differential quantities $dt$ can be ``cancelled out", this is not the most useful way of thinking about this expression.  In reality, we \emph{must} mean something different from that interpretation because the derivative $dy/dt$ represents the \emph{limit of a ratio as $dt$ tends to zero}.  However, we cannot think of $dy$ or $dt$ individually as being defined by such a limit; otherwise, these two quantities would be identically zero, and Eq.~\ref{differential2} would not be particularly useful.  Instead, we think of this relationship between differentials as being described as follows.
\begin{quote}
    When the increment of the independent variable, $dt$ is made sufficiently small, then the function $y(t)$ changes by the amount $dy$, with an error $\epsilon$.  The error $\epsilon$ in this approximation that can be made as small as we like by decreasing $dt$.
\end{quote}

With this recognition, we are in a position to solve Eq.~\eqref{integralsep}.  First we make the substitution $dy = (dy/dt) dt$ to give  

\begin{align}
    \int y^2 dy = \int t^2 \, dt
\end{align}

Integrating this, we find 

\begin{equation}
    \frac{1}{3} y^3 = \frac{1}{3} t^3 +C_1
\end{equation}
Noting that $3C_1$ is just a new constant, we can write this in the form

\begin{equation}
  y^3 = t^3 +C_2
\end{equation}
There are no ancillary conditions given for this problem, so this solution is complete.  Usually, if it is possible, we try to express a solution in its \emph{explicit} form.  Doing so here yields

\begin{equation}
 y(t) = (t^3 +C_2)^{\tfrac{1}{3}}
\end{equation}

\end{example}
\end{svgraybox}

Now that we have seen a concrete example, it is pretty easy to see how these problems can be solved generally.  Here is the solution method, but done with arbitrary functions.  The problem is

\begin{align}
\frac{dy}{dt}&= r(y)Q(t)\\
y(t_0) & = y_0 
\intertext{or, recalling $R(y) = 1/r(y)$}
R(y)\frac{dy}{dt}&= Q(t) \\
y(t_0) & = y_0 
\end{align}
Now, we need just integrate both sides with respect to $t$.  This gives

\begin{align}
\int R(y)\frac{dy}{dt}\, dt= \int Q(t)\, dt
\end{align}
As mentioned in the previous example problem, in Chapter 1 (Section \ref{riemann}) we reviewed the concept of the differential as part of the review of integration theory.  There, the case was made that while the derivative was defined by a limiting process of the ratio of two quantities that approach zero, the differential was defined as a \emph{nonzero} quantity.  In essence, to think about differentials, we assume that we can treat ``dy" and ``dt" not by a limiting process, but by thinking about them as being a ratio of two quantities that approximate the derivative with an error, $\epsilon$ that can be made as small as we like.  With this model in mind, we can properly treat first-order derivatives almost as if they represent algebraic quantities.  Thus we can write in general

\begin{equation}
    dy = \frac{dy}{dt} dt
\end{equation}
giving

\begin{align}
\int R(y) dy= \int Q(t)\, dt+C
\label{separablesoln}
\end{align}
where the final constant represents the constant of integration to be determined by the ancillary conditions (if any) that are provided.

\begin{svgraybox}
\begin{example}[A first-order, linear equation]
Suppose we have the problem
\begin{equation}
y \frac{dy}{dt} = t^2 
\end{equation}
with ancillary condition
\begin{equation}
y(0)=1
\end{equation}
What is the explicit solution for $y(t)$?\\
~\\
{\bf Solution.}
We can use the expression given by Eq.~(\eqref{separablesoln}) to solve this problem.  Clearly, the solution is given by
\begin{equation}
\int y \,dy   =\int t^2\,dt+C
\end{equation}
or, upon computing the integrals
\begin{equation}
y^2(t) -y^2(0) = \tfrac{2}{3}t^3+C
\end{equation}
Noting the initial condition, this can be expressed
\begin{equation}
y^2(t) =( 1+\tfrac{2}{3}t^3)
\end{equation}
where we must have $C=0$ (it is easy to verify that this expression matches the ancillary condition at $t=0$).  This is an implicit equation for $y$, since the expression does not provide us with the value of $y$, but, rather, the value of $y^2$.  In this case we are lucky-- the expression is a quadratic, and we can use the quadratic formula (or, in this case, simply taking the square root of both sides of the expression) to get an explicit expression for $y$.
\begin{equation}
y(t) = \pm \,\sqrt[]{1+\tfrac{2}{3}t^3}
\end{equation}
The derivative of the solution is
\begin{equation}
\frac{dy}{dt}=\pm\, \tfrac{1}{2}(1+\tfrac{2}{3}t^3)^{-\tfrac{1}{2}}\, 2t^2
\end{equation}
And with this information, it is easy to verify that the original nonlinear ODE is also met by the solution.

Note that we can also solve this problem using \emph{definite} integration.  When we do so, we no longer have a constant of integration because this constant is automatically accounted for by the information in the bounds.  To solve the problem using definite integrals, we integrate in time from where information is known (at $t_0$) to any unknown time (symbolized by simply $t$) as follows

\begin{equation}
\int_{y(0)}^{y(t)} y \,dy   =\int_{\tau= 0}^{\tau = t} \tau^2\,d\tau
\end{equation}

Note that here we are using the variable of integration, $\tau$, to represent the time variable within the integral.  Recall, this is a conventional notation for integration that helps prevent confusion of variables of integration and variables within the bounds of integration.  Integrating this, we find the solution

\begin{equation}
\left. y(t) \right|_{y(0)}^{y(t)} = \left. \tfrac{2}{3} \tau^3\right|_{\tau = 0}^{\tau=t}
\end{equation}
Noting the initial condition $y(0)=1$, we have the solution

\begin{equation}
y(t) = \pm \,\sqrt[]{1+\tfrac{2}{3}t^3}
\end{equation}
which is identical to the solution via indefinite integration.
\end{example}
\end{svgraybox}

The previous problem was one that happened to be linear in the function $y$ (Question: can you show this?).  For nonlinear separable equations, the solution to nonlinear problems proceeds essentially the same way.  The only real issue one faces with separable problems is whether or not the two sides, once separated, can be integrated.  

In the next example, a nonlinear problem in the dependent variable is solved using separation.  For this problem, it is not possible to express the result in an explicit form (i.e., in the form $y(t) = function(t)$.  Instead, we leave the result in implicit form.

\begin{svgraybox}
\begin{example}[Michaelis-Menten kinetics]
Michaelis-Menten kinetics describe kinetic processes that are rate-limited by the total number (or amount) of compound that causes the reaction to occur.  As an example, reactions due to a catalysts on the surface of a carrier particle (of the form $A\xrightarrow{{Catalyst}}B$ are sometimes rate-limited because once all of the catalysts sites are occupied by the reactant ($A$).  The reaction rate  is given by
\begin{equation*}
    \frac{dc_A}{dt}=- \frac{k c_A}{c_A+K_A}
\end{equation*}
where $k$ is the rate constant (mol/(L$\cdot$s) and $K_A$ is the half-saturation constant (mol/L).  Michaelis-Menten kinetics are sometimes also called \emph{saturation} kinetics, because the kinetic rate achieves a maximum value as the concentration of the reactant species increases.  This is illustrated

\begin{centering}
\includegraphics[scale=.15]{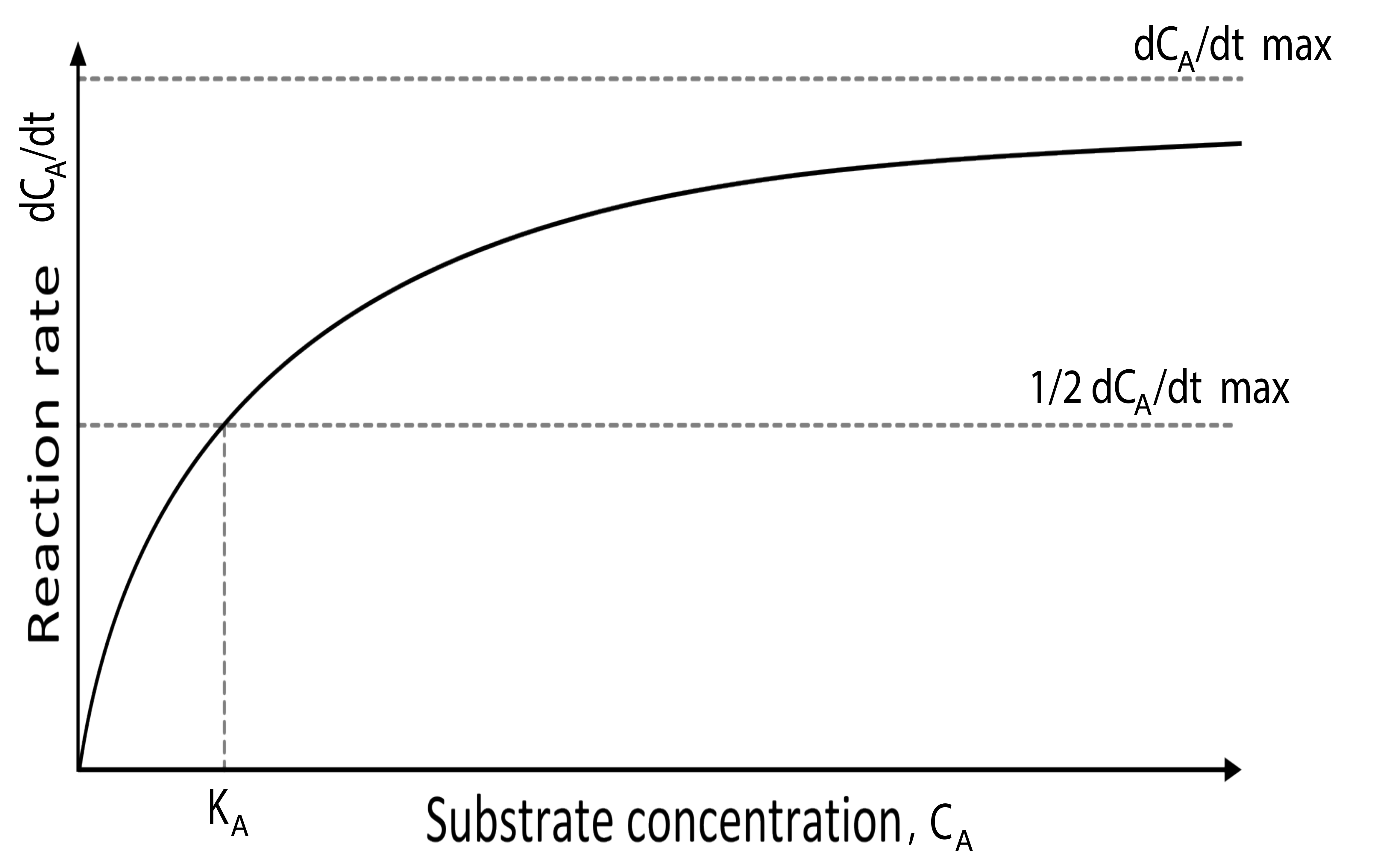}
\captionof{figure}{Michaelis-Menten kinetics.  Two interesting features about these kinetics can be seen through the geometry of the curve of rate versus concentration.  At concentrations that are much larger than the half-saturation constant, $K_m$, the rate of reaction is approximately equal to the kinetic rate parameter, $k$.  The half-saturation constant, $K_A$ can be found by first finding the rate equal to half of the maximum rate; then, one find the corresponding value on the concentration axis.  This value is equal to $K_A$. }
\label{fig:1.1}  
\end{centering}

Determine a functional relationship for $c_A(t)$ assuming this rate law, and the initial condition $c_A(0)=C_0$.  How long does it take for $c_A$ to get to one-half of its initial value?

\noindent \emph{Solution}.  Start by doing separation of variables, leading to the form

\begin{align*}
    -\frac{c_A+K_A}{c_A} dc_A= k dt \\
     -\int_{C=c_A(0)}^{C=c_A(t)}\frac{C+K_A}{C} dC= \int_{\tau = 0}^{\tau=t} k d\tau \\
\end{align*}
(Note that here, we have been careful to use variables of integration when writing out the \emph{definite integral}.  Sometimes this is helpful to do; it reinforces what the correct bounds for the integration are.  However, for compactness and expedience, once one is familiar with indefinite integration, the explicit switching to a separate dummy variable for integration is not necessary.  The resulting expression is a bit ``sloppy", but as long as it is clear what is being done, it usually causes no problems.)

Performing the integrations, we find
\begin{align*}
    -\left. C\right|_{c_A(0)}^{c_A(t)} -K_A \left. \ln{C}\right|_{c_A(0)}^{c_A(t)}= kt
\end{align*}
Evaluating this result at the bounds, and using the properties of logarithms, the final result is

Performing the integrations, we find
\begin{align*}
    t=\frac{1}{k}\left[ C_0 - c_A(t) +K_A  \ln{\frac{c_A}{C_0}}\right]
\end{align*}
As a check on this result, note two properties that we would expect from this solution (1) The solution meets the initial condition at $t=0$ (i.e., for $t=0$, $c_A(0)=C_0$), and (b) the result is positive for all values of $c_A$; we expect this latter property because time can not be a negative quantity!  Also note that argument of the logarithm is \emph{dimensionless}.  It is a good practice to assure that in your \emph{final results} all logarithmic arguments should be dimensionless!  They arguments may not be dimensionless during intermediate steps, but if they cannot be made to be free of units at the end of the problem, this may indicate that there has been an error.  Generally speaking, it is improper to take the logarithm of a quantity with units.

To solve the second part of the problem, we substitute $c_A(t)=C_0/2$, and solve for time; the result is
\begin{align*}
    t=\frac{1}{k}\left[ \tfrac{1}{2}C_0 +K_A  \ln{2}\right]
\end{align*}
\end{example}
\end{svgraybox}

\begin{svgraybox}
\begin{example}[Solution to a linear first-order ODE by separation of variables]
Suppose we have the classical first-order reaction in a completely-mixed reactor problem discussed above, and specified by the set of equations
\begin{align}
\frac{dy}{dt} &= -k_1 \,y \\
y(t=0) = c_0
\end{align}
What is the solution?\\

\emph{Answer.} Separating variables we have
\begin{align}
\frac{1}{y} \,dy&= -k_1\, dt
\end{align}
Integrating both sides of this immediately yields
\begin{align}
\int \frac{1}{y} \,dy&= -k_1 \int \, dt
\end{align}
or
\begin{align}
\ln(y) &= -k_1 t + C \\
y(t) &=C\, \exp[-k_1 t]
\end{align}
Substituting the ancillary equation into this result, we find
\begin{align}
y(0) &= C \\
&or\\
C&= c_0
\end{align}
so that the final result is the familiar exponentially-decreasing function of time
\begin{equation}
y(t) =c_0\, \exp[-k_1 t]
\end{equation}
\end{example}
\end{svgraybox}

\begin{svgraybox}
\begin{example}[Solution to a nonlinear first order ODE by separation of variables]
Suppose instead of the classical first-order reaction in a completely-mixed reactor problem, we have a second-order problem specified by
\begin{align}
\frac{dy}{dt} &= -k_1 y^2 \\
y(t=0) & = c_0
\end{align}
What is the solution?\\

\emph{Answer.} Separating variables as above, we have
\begin{align}
\frac{1}{y^2}\,dy &= -k_1 \,dt \\
\end{align}
Integrating both sides of this immediately yields
\begin{align}
\int \frac{1}{y^2} \,dy&= -k_1 \int \, dt
\end{align}
or
\begin{equation}
-y^{-1} = -k_1 t + C 
\end{equation}
Upon solving explicitly for $y(t)$ we have
\begin{equation}
y(t) =-\frac{1}{-k_1 t + C}
\end{equation}
Finaly, substituting the ancillary equation into this result, we find
\begin{align}
y(0) &= -\frac{1}{C} \\
&or\\
C&=-\frac{1}{c_0}
\end{align}
so that the final result is the following function of time
\begin{equation}
y(t) =\frac{c_0}{k_1 c_0 t +1}
\end{equation}
Noting that at $t=0$, we have $y(0)=c_0$, we find that the ancillary condition is met.  Also noting that the derivative is given by
\begin{equation}
\frac{dy}{dt} = -\frac{k_1 c_0}{(k_1 c_0 t +1)^2}
\end{equation}
it is easy to validate that the solution not only meet the initial condition, but it meets the original ODE also.
\end{example}
\end{svgraybox}
%
\subsection{Non-Separable Linear First-Order ODEs}
%
In the previous sections, first-order ODEs were solved by direct use of the fundamental theorem of calculus.  This approach required that the equations be \emph{separable}.  However, there are plenty of examples of equations that are not immediately separable.  The canonical form for these equations is

\begin{equation}
\frac{dy}{dt}+ P(t) \,y = Q(t)
\end{equation}
There is no obvious way to integrate both sides of this equation to get a solution. Doing so would result in an integral with integrand $P(t)y(t)$; because we do not \emph{know} $y(t)$, we are stuck with an implicit integral solution for $y(t)$.  While technically not a disaster, it is not a solution that is easy to compute the values of the dependent variable.

To handle problems of this general form for linear first-order ODEs requires first transforming this problem.  This idea is a key element in analytical modeling and in mathematical analysis generally; it is important enough to highlight the idea specifically.
\begin{svgraybox}
\textbf{Key Idea 1.1}. When encountering a problem that one does not know how to solve, it is worthwhile considering if it is possible to transform the problem into a form that does have a known solution.
\end{svgraybox}
\noindent This idea behind this statement is simultaneously extremely useful, and potentially frustrating!  It suggests what must be attempted, but it provides no information about how this is to be accomplished.  This brings up a second key idea that arises in mathematical analysis.

\begin{svgraybox}
\textbf{Key Idea 1.2}.  Mathematical analysis often involves the creative use of intuition or experimentation to find solutions.
\end{svgraybox}
An example of these two ideas is exactly what is embodied by the following solution to the general first-order linear ODE.  It is often most useful to see such developments as an example before stating the more general and abstract result.  

\begin{svgraybox}
\begin{example} \label{cstrex}
Suppose we have a continuously-stirred tank reactor (CSTR) undergoing a first-order reaction process $A\rightarrow B$ with kinetic rate constant $k_1$.  If the flow rate in and out are both equal to $Q$, the concentration of species $A$ is $c_0$ at the inlet, and the initial concentration of species $A$ is $y_A(0)=0$ then the net balance for the concentration of species $A$ reactor is given by the ODE

\begin{align}
\frac{dy_A}{dt}&=\frac{Q}{V} c_0 - \frac{Q}{V} y_A-k_1 y_A \label{cstr}\\
y_A(0) &= 0
\end{align}\\

\noindent Here, we note that the quantity \smash{$\frac{Q}{V}$} is related to the \emph{residence time} for the reactor, $\tau$
\begin{equation}
\frac{Q}{V}=\frac{1}{\tau}
\end{equation}

{\centering
\includegraphics[scale=.1]{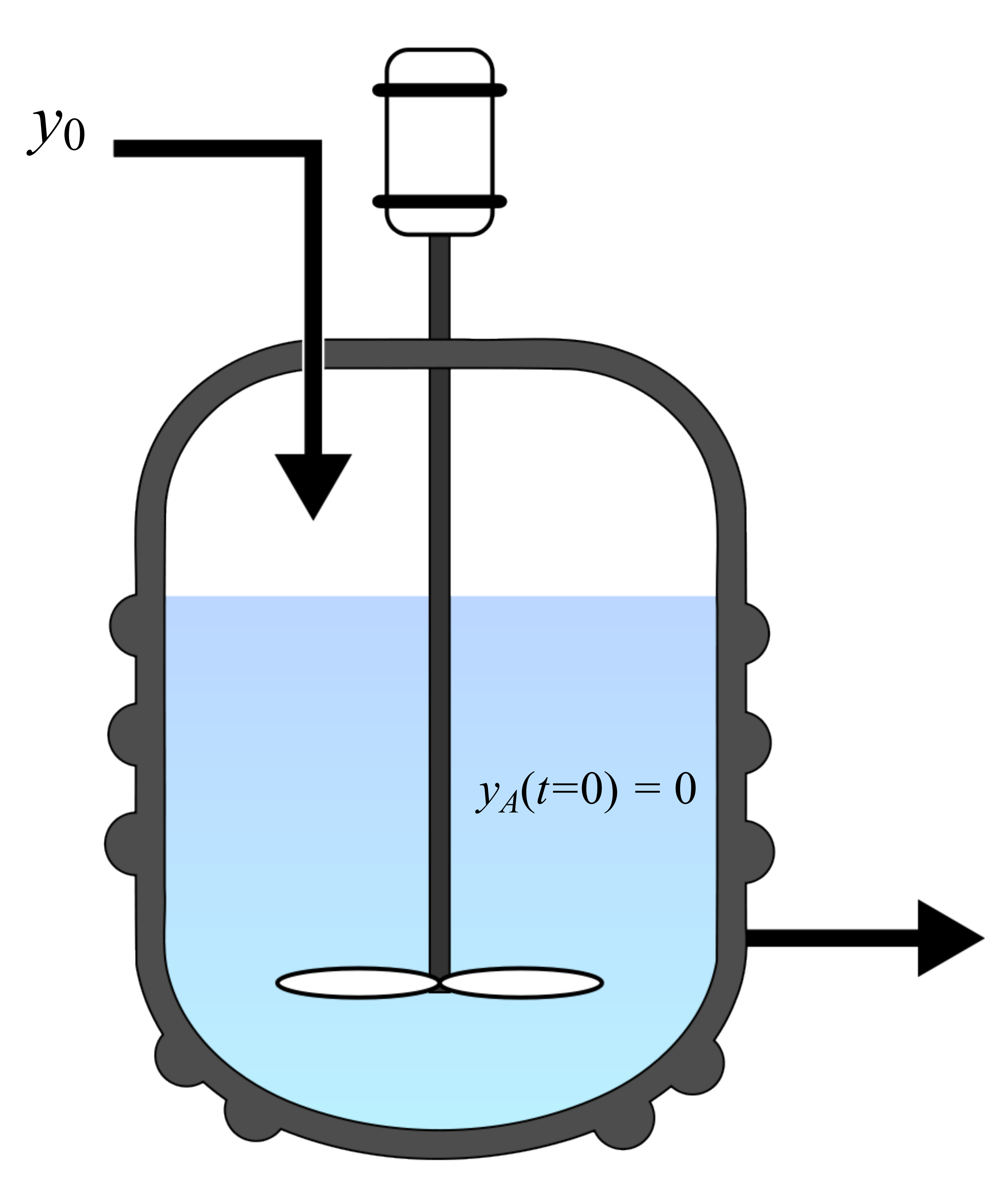}
\captionof{figure}{A continuously-stirred tank reactor. The tank initially has zero concentration of species $A$ or $B$, but does contain a suspended catalyst (as pellets that do not leave the system).  At time $t=0$, species $A$ begins to enter the reactor, and the catalytic reaction proceeds.}
}
\vspace{5mm}

To solve this problem, we start by rewriting the equation in the following form

\begin{equation}
\frac{dy_A}{dt}+(\tau^{-1}+k_1)y_A= \tau^{-1} y_0 
\end{equation}
Now, it is clear that this expression is not of the form of one that has been covered in the material above.  However, we note that if we set $\beta = (\tau^{-1}+k_1)$, and then multiply this expression through by $e^{\beta t}$, we have

\begin{equation}
\frac{dy_A}{dt}e^{\beta t}+ y_A \beta e^{\beta t}= e^{\beta t}\tau^{-1} y_0 
\end{equation}
Although it may not be obvious that this has improved things, it has once we make an observation: The left-hand side of the expression is now a total derivative.  In other words, we note

\begin{equation}
\frac{d}{dt}(y_A e^{\beta t}) = \frac{dy_A}{dt}e^{\beta t}+y_A \beta e^{\beta t} 
\end{equation}
But this is exactly the left-hand side of the equation above!  Thus, we can substitute as follows
\begin{equation}
\frac{d}{dt}(y_A e^{\beta t})= e^{\beta t}\tau^{-1} y_0 
\end{equation}
Now, we have transformed the original equation into one we know how to solve- this expression is separable in the conventional way that was discussed previously.  If this is unclear, we can make the following temporary transformation of variables
\begin{equation}
\psi(t) = y_A e^{\beta t}
\end{equation}
Then, we have 
\begin{equation}
\frac{d\psi}{dt}= e^{\beta t}\tau^{-1} y_0 
\end{equation}
Separating variables and integrating, we find
\begin{equation}
\int d\psi= \int e^{\beta t} \tau^{-1} y_0 \, dt +C
\end{equation}
This yields
\begin{equation}
\psi(t)= \frac{1}{\beta} e^{\beta t}\tau^{-1} y_0 \, dt +C
\end{equation}
Returning to the original variables, we have

\begin{equation}
y_A \exp\left({(\tau^{-1}+k_1) t}\right) =\frac{1}{(\tau^{-1}+k_1)} \exp\left({(\tau^{-1}+k_1) t}\right) \tau^{-1} y_0 + C
\end{equation}
Finally, solving for $y_A$ gives the expression
\begin{equation}
y_A(t)  =\frac{1}{(1+\tau k_1)}  y_0 +\exp\left(-(1+\tau k_1) \frac{t}{\tau}\right)C
\end{equation}
Setting $t=0$ allows us to use the initial condition to solve for $C$
\begin{equation}
C  =-\frac{1}{(1+\tau k_1)}y_0 
\end{equation}
Thus, the final solution is

\begin{equation}
y_A(t)  =\frac{y_0}{(1+\tau k_1)} \left[1-\exp\left(-(1+\tau k_1) \frac{t}{\tau}\right) \right]
\end{equation}

\end{example}
\end{svgraybox}

%
In Example \ref{cstrex}, the function $P(t)$ that was incorporated into the exponential was a simple constant, $\beta = \tau^{-1}+k_1$.  Recall, however, the general form is 

\begin{equation}
\frac{dy}{dt}+ P(t) \,y = Q(t)
\label{general1st}
\end{equation}
The solution here is actually no more difficult than when $P(t)$ is a constant.  For this more general case, take 
\begin{equation}
s(t) = \int P(t) \,dt
\end{equation}
This transformation can be done as an indefinite integral.  The constant for the indefinite integral can be taken as zero (and \textbf{always} will be when we use an integrating factor), because any such constant eventually can be eliminated from the resulting equations.  Note that 
\begin{equation}
\frac{ds}{dt} =  P(t)
\end{equation}
By the fundamental theorem of calculus.

Now, multiplying Eq.~(\ref{general1st}) through by $\exp[\beta(t)]$ we find
\begin{equation}
\frac{dy}{dt}\exp[s(t)]+ P(t)\exp[s(t)] \,y = \exp[s(t)] Q(t)
\end{equation}
This expression can be rewritten as
\begin{equation}
\frac{d}{dt}\left(\exp[s(t)] \,y\right) = \exp[s(t)] Q(t)
\label{imintegrable}
\end{equation}
Now, set \smash{$\psi(t) =\exp[s(t)] \,y $}, separate variables, and integrate both sides as \emph{definite} integrals.  This results in
\begin{equation}
\int_{\psi(t'=t_0)}^{\psi(t'=t)} \,d\psi= \int_{t'=t_0}^{t'=t}\exp(s(t')) Q(t') dt'
\end{equation}
Integrating this yields
\begin{equation}
\psi(t) -\psi(t_0) = \int_{t'=t_0}^{t'=t}\exp(s(t')) Q(t') dt'
\end{equation}
Upon returning to the original variables, we find
\begin{equation}
\exp[s(t)] \,y(t) -\exp[s(t_0)] \,y(t_0) = \int_{t'=t_0}^{t'=t}\exp[s(t')] Q(t') dt'
\end{equation}
And solving this for $y(t)$ gives us the general formula for the solution to the separable first-order linear ODE
\begin{equation}
\boxed{
~y(t)  =\exp[-s(t)]\left[\int_{t'=t_0}^{t'=t}\exp[s(t')] Q(t') dt'+ y(t_0)\exp[s(t_0)]\right]~
}
\label{parts}
\end{equation}
where recall
\begin{equation}
s(t) = \exp\left(\int P(t) dt\right)
\end{equation}
This expression gives the solution, then, to \emph{any} linear non-separable first-order ODE.  

The entire analysis above can also be done using indefinite integration.  The entire process is exactly the same up to Eq.~\ref{imintegrable}.  From there, we can integrate both sides as indefinite integrals, but remembering to add the appropriate constant of integration.  Thus, we would have the result (starting right after Eq.~\ref{imintegrable})

\begin{equation}
\int \,d\psi= \int\exp(s(t)) Q(t) dt
\end{equation}
\begin{equation}
\psi(t)= \int\exp(s(t)) Q(t) dt +C_1
\end{equation}
Where $C_1$ is a constant of integration from the integration of the left-hand side.  Note, that the integral on the right-hand will also generate a constant of integration. However, the two constants can be combined, so we can ignore the constant coming from the remaining integral. 
Recalling that $\psi(t) =\exp[s(t)] \,y $, then this can be simplified to give $y(t)$ as follows

\begin{equation}
\exp[s(t)] \,y(t)= \int\exp(s(t')) Q(t') dt' +C_1 
\end{equation}
\begin{equation}
\boxed{
~y(t) = \exp[-s(t)]\left[\int\exp(s(t')) Q(t') dt' +C_1 \right]~
}
\label{secondvers}
\end{equation}
Where, again, remember that \emph{only one constant}, $C_1$, is generated for this problem in total, which is consistent for the fact that it is a first-order problem.
By comparing the two boxed results above, it is clear that $C_1=y(t_0)\exp[s(t_0)]$.  However, some people prefer to use the simpler looking version given by Eq.~\eqref{secondvers}, and then simply evaluate the unknown constant using the ancillary data.  The results will be the same regardless of which approach is taken.

\begin{svgraybox}
\begin{example}[CSTR with catalyst deactivation.]
Solid catalysts are used in a huge variety of reactors, from high-tech reactors to make pharmaceuticals to reactors for ion-exchange in the treatment of drinking water.  One potential problem with solid catalysts is that they can become deactivated over time.  There are any number of deactivation processes (e.g., particle sintering, catalyst poisoning, and catalyst coking are three types of deactivation that are common); all of them reduce catalyst effectiveness, and have significant economic impact.

One model for catalyst inactivation involves an inverse-time function that decreases the net rate of reaction with increasing time.  Suppose we have a reaction that proceeds in the presence of a catalyst with the simple form \smash{$A\mathop  \to \limits^C B$}, where $A$ is the reactant, $B$ is the product, and $C$ is the catalyst.  
Returning to the case of a CSTR from the previous example, the molar mass balances for a catalytic reaction with deactivation for chemical species $A$  can be specified by

\begin{align}
\frac{dy_A}{dt} &= \tau^{-1}y_0 -\tau^{-1} y_A -\left(\frac{1}{1+k_d t} \right) k_1 y_A
\label{deactive}\\
y_A(0)&=0 \\
y_B(0)& = 0 \\
\end{align}
We have assumed that the catalyst concentration is constant in the CSTR (catalyst neither enters nor leaves the reactor).
Note that because initially there is no species $B$ in solution, the concentration of species $B$ at any time can be computed from how much of species $A$ has been converted.  While it is possible to also get the concentration of species $B$ at any time, for now we will focus only on the balance for species $A$.

Here, we can see that Eq.~(\ref{deactive}) is identical to Eq.~(\ref{cstr}), except the reaction term includes an additional function of time that multiplies the first-order reaction.  This function decreases the net rate of reaction as time increases.
This problem would be somewhat intimidating if we had not previously developed the integration factor formula given by Eq.~(\ref{parts}) (and even with it, it may still be challenging!)  Rewriting the mass balance equation, we find
\begin{equation}
\frac{dy_A}{dt} +\left[ \tau^{-1}+\left(\frac{k_1}{1+k_d t} \right) \right] y_A= \tau^{-1}y_0 
\end{equation}
In this form, it is clear that we have
\begin{equation}
s(t) = \int \left[ \tau^{-1}+\left(\frac{k_1}{1+k_d t} \right) \right] dt
\end{equation}
or, after completing the integral
\begin{equation}
s(t) =  \left[\frac{t}{\tau}+\frac{k_1}{k_d} \ln(1+k_d t)\right] 
\end{equation}
This makes rather short work of writing down the solution; recalling Eq.~(\ref{parts})

\begin{equation}
    y(t)  =y(t_0)\exp[s(t_0)]\exp[-s(t)]+ \exp[-s(t)]\,\int_{t'=t_0}^{t'=t}\exp[s(t')] Q(t') dt'
\end{equation}
Making the appropriate substitutions (and noting that $y(t_0)=0$), we have the solution
\begin{align}
y_A(t) &= \exp\left[-\left(\frac{t}{\tau}+\frac{k_1}{k_d} \ln(1+k_d t)\right)\right]
\int_{t'=0}^{t'=t} \exp\left[\frac{t}{\tau}+\frac{k_1}{k_d} \ln(1+k_d t)\right] \frac{y_0}{\tau} dt'
\end{align}
Computing this last integral is somewhat of a challenge; however, the integral is tabulated in tables, and can be computed by symbolic mathematics programs like Mathematica.  The result is
\begin{align}
y_A(t) &= \tau  e^{-\frac{1}{k_d \tau }} \left(\frac{1}{k_d \tau
   }\right)^{-\frac{k_1}{k_d}} \left[\Gamma
   \left(\frac{k_1+k_d}{k_d},-\frac{k_d t+1}{k_d \tau
   }\right)-\Gamma \left(\frac{k_1+k_d}{k_d},-\frac{1}{k_d \tau
   }\right)\right]
\end{align}
Plots of the solution for species $A$ and $B$ are provided in Fig.~\ref{fig:1.2}.\\
\indent The Gamma function, $\Gamma$ is a transcendental function that corresponds to a generalization of the factorial operation to the real (or complex) numbers.  The function is widely used, and is usually an intrinsic function in both symbolic mathematics programs like Mathematica, and in interpreted coding languages like MATLAB. \\

{
\centering\includegraphics[scale=.7]{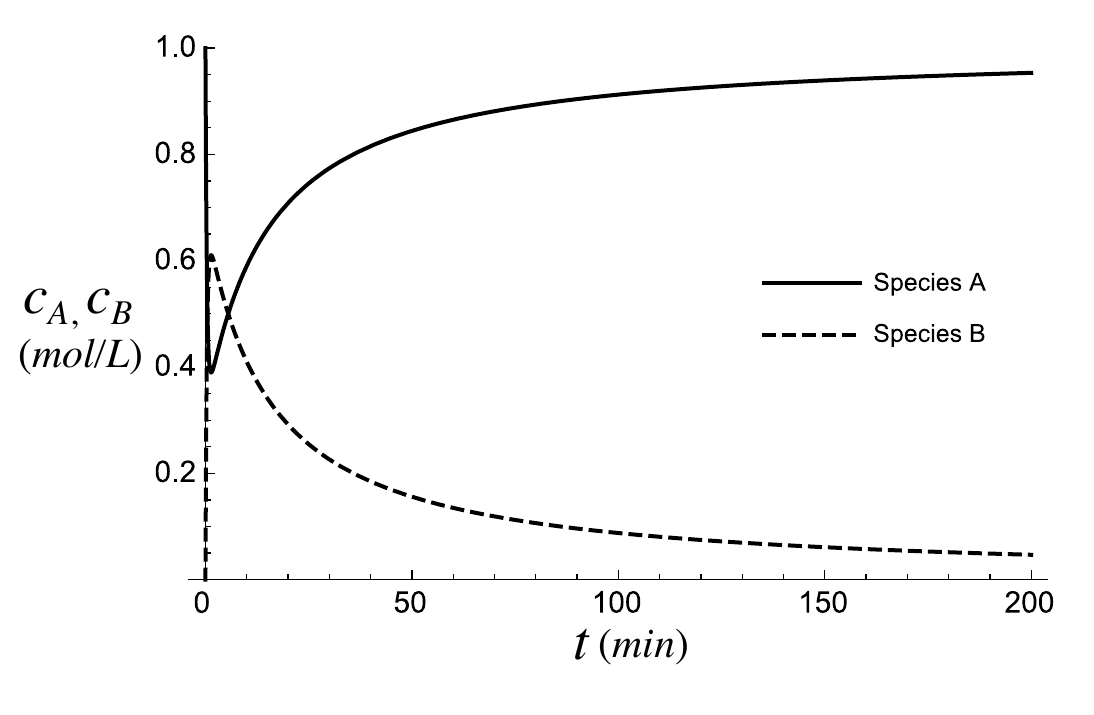}
\captionof{figure}{Solutions for species $A$ and species $B$ as a function of time for a catalytic reaction with deactivation.  For this solution, $k_1=2~min^{-1}$, $k_d=0.1~min^{-1}$, $\tau=1~min^{-1}$, $y_0=1~mol\cdot m^{-3}$.}
\label{fig:1.2}  
}
\vspace{5mm}
Just for completeness, note that the two Gamma functions are defined by the following integrals
\begin{align}
Gamma&&\Gamma_1(x) &= \int_{t=0}^{t\rightarrow \infty} t^{x-1} e^{-t} dt \\
Incomplete~Gamma&&\Gamma_2(a,x) &= \int_{t=x}^{t\rightarrow \infty} t^{a-1} e^{-t} dt 
\end{align}
More information on the Gamma function can be found in many texts.  The Gamma function can be evaluated using common computational platforms such as MATLAB or Mathematica, both of which have built-in routines for computing these functions.
\end{example}
\end{svgraybox}

\section{Second-Order ODEs with Constant Coefficients}\indexme{ordinary differential equation!second order}

We will cover exactly two kinds of second-order ODEs: (1) directly integrable second-order ODEs, and (2) general second-order ODEs with constant coefficients.

Linear second-order ODEs are functions of a single independent variable that take the form
\begin{equation}
    \frac{d^2 y}{dx^2}+ b(x)\frac{d y}{dx} + c(x) y  = g(x)
    \label{general_second}
\end{equation}

In general, solutions of these equations (where the coefficients are functions rather than constants) require series solution methods which we will not be reviewing.  However, there are a few other interesting problems that are solvable. These are the (1) directly integrable problems, and (2) the case of linear second-order ODEs with constant coefficients.  

\subsection{Directly Integrable Second-Order ODEs}

If a second-order can be put in the form

\begin{equation}
    \frac{d}{dx}\left( p(x) \frac{d y}{dx} \right) = h(x)
    \label{second}
\end{equation}
then, in principle, it can be solved directly by two integrations through the application of the fundamental theorem of calculus (FTC).

Suppose the ODE above is defined on $x\in[a,b]$, it has two appropriate ancillary conditions imposed, and $p(x) \ne 0$ anywhere in the domain.  Then, integrating both sides of Eq.~\eqref{second} gives (using the fundamental theorem of calculus), and then dividing by $p(x)$

\begin{equation}
      \frac{d y}{dx}  = \frac{1}{p(x)}\int_x h(\xi) d\xi
\end{equation}
where, here we have used a dummy variable of integration to avoid confusion in later developments.
A second integration gives a solution for $y(x)$

\begin{equation}
      y(x)  = \int_x \frac{1}{p(x)} \int_x h(\xi) d\xi dx
\end{equation}
Note that these two integrations will generate two constants of integration.

\begin{svgraybox}
\begin{example}[Directly Integrable Second-Order ODE.]
Suppose we have an integrable of the form
\[ 
 \frac{d}{dx}\left( (x^2+1) \frac{d y}{dx} \right) = 3 x^2, ~~x\in[0,1]\]
 with 
 \[ y(0)=0, y(1)=0 \]
 What is the solution?
 
 \noindent{\bf Solution.}  To start, integrate both sides as an indefinite integral and divide through by $x^2+1$
 
 \[ 
 \frac{d y}{dx}  = \int_x \frac{x^3+c_1}{x^2+1} \, dx \]
Now, this integral is not (at least apparently) easy to do.  One could try something like the method of partial fractions (usually introduced in introductory calculus) to integrate it.  Alternatively, there are tables of integrals, and various symbolic mathematics languages that can be used to find the value of this integral.  We are going to take the latter option here, noting that the integral can be computed as 
\vspace{-2mm}
\[ 
 \frac{d y}{dx}  = \frac{1}{2} [x^2 +2 c_1 \arctan(x) -\ln(1+x^2)]
 \]
A second integration yields
\begin{align*}
y(x) & = {c_1} \left(x \arctan(x)-\frac{1}{2} \log
   \left(x^2+1\right)\right)+\frac{x^3}{6} \\
   &-\frac{1}{2}x \log \left(x^2+1\right)+ x- \arctan(x)+c_2
\end{align*}
Applying the first ancillary condition, $y(0)=0$, we find $c_2=0$. Applying the second ancillary condition $y(1)=0$ gives
\[ 
0 = c_1\left(\frac{\pi}{4}-\frac{1}{2}  \log (2)\right)-\frac{\pi
   }{4}+\frac{7}{6}-\frac{\log (2)}{2}
\]
or
\[ c_1 = \frac{ \frac{\pi}{4}-\frac{7}{6} + \frac{1}{2}\log (2)}{\pi -2 \log (2)} \]
So, the final solution is, somewhat unbelievably, the function
\vspace{-2mm}
\begin{align*}
y(x)  &= \left(\frac{ \frac{\pi}{4}-\frac{7}{6} + \frac{1}{2}\log (2)}{\pi -2 \log (2)} \right)
\left(x \arctan(x)-\frac{1}{2} \log
   \left(x^2+1\right)\right)+\frac{x^3}{6} \\
   &-\frac{1}{2}x \log \left(x^2+1\right)+ x- \arctan(x)
 \end{align*}

{
\centering\includegraphics[scale=.45]{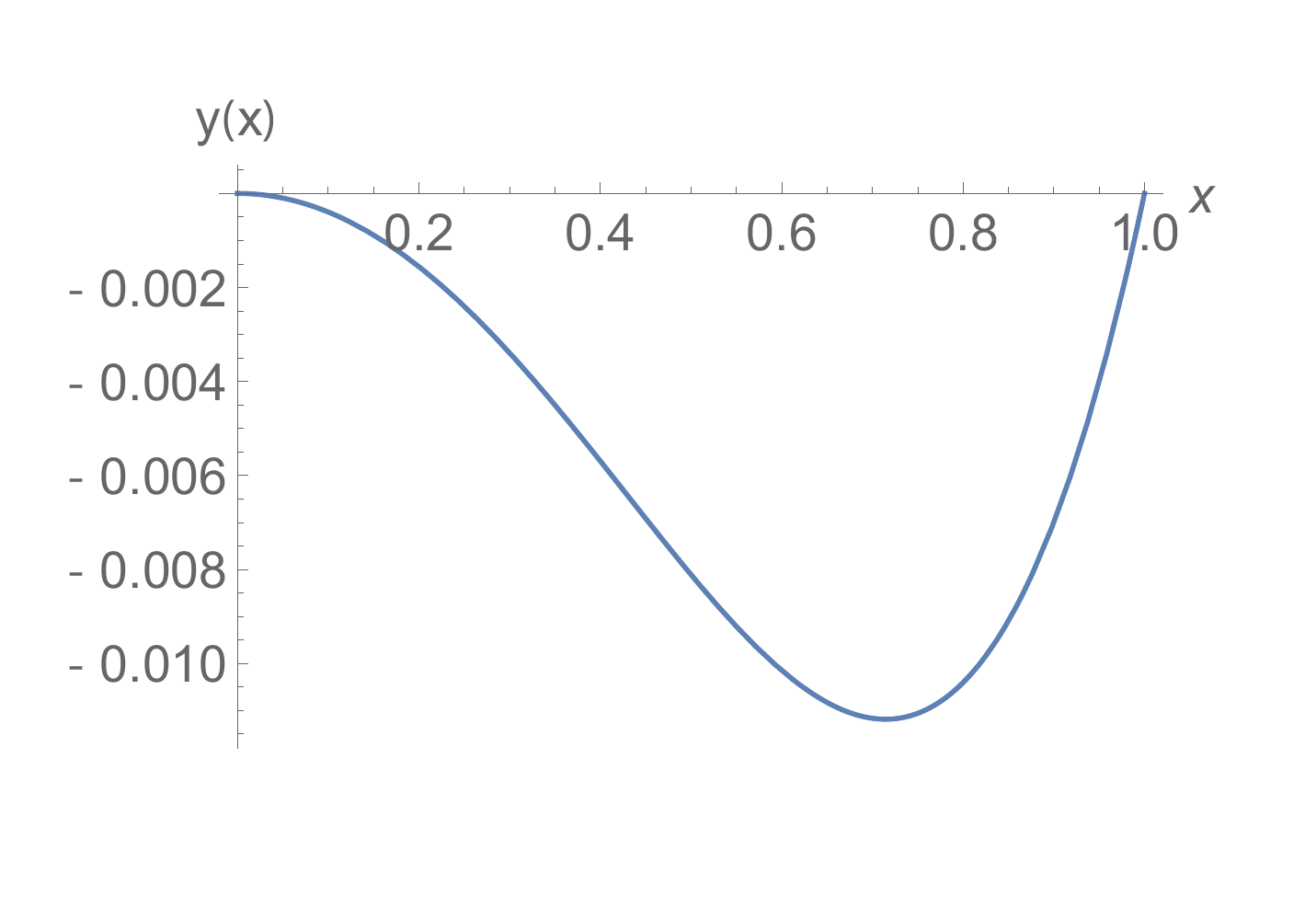}
\vspace{-5mm}
%
\captionof{figure}{Plot of solution, $y(x)$}
}

It is possible to verify that $y(0)=y(1)=0$ for this result, suggesting strongly that it is a correct one.
\end{example}
\end{svgraybox}

\subsection{Thinking about Solutions to Second-order ODEs}

\begin{figure}[t]
\sidecaption[t]
\centering
\includegraphics[scale=.5]{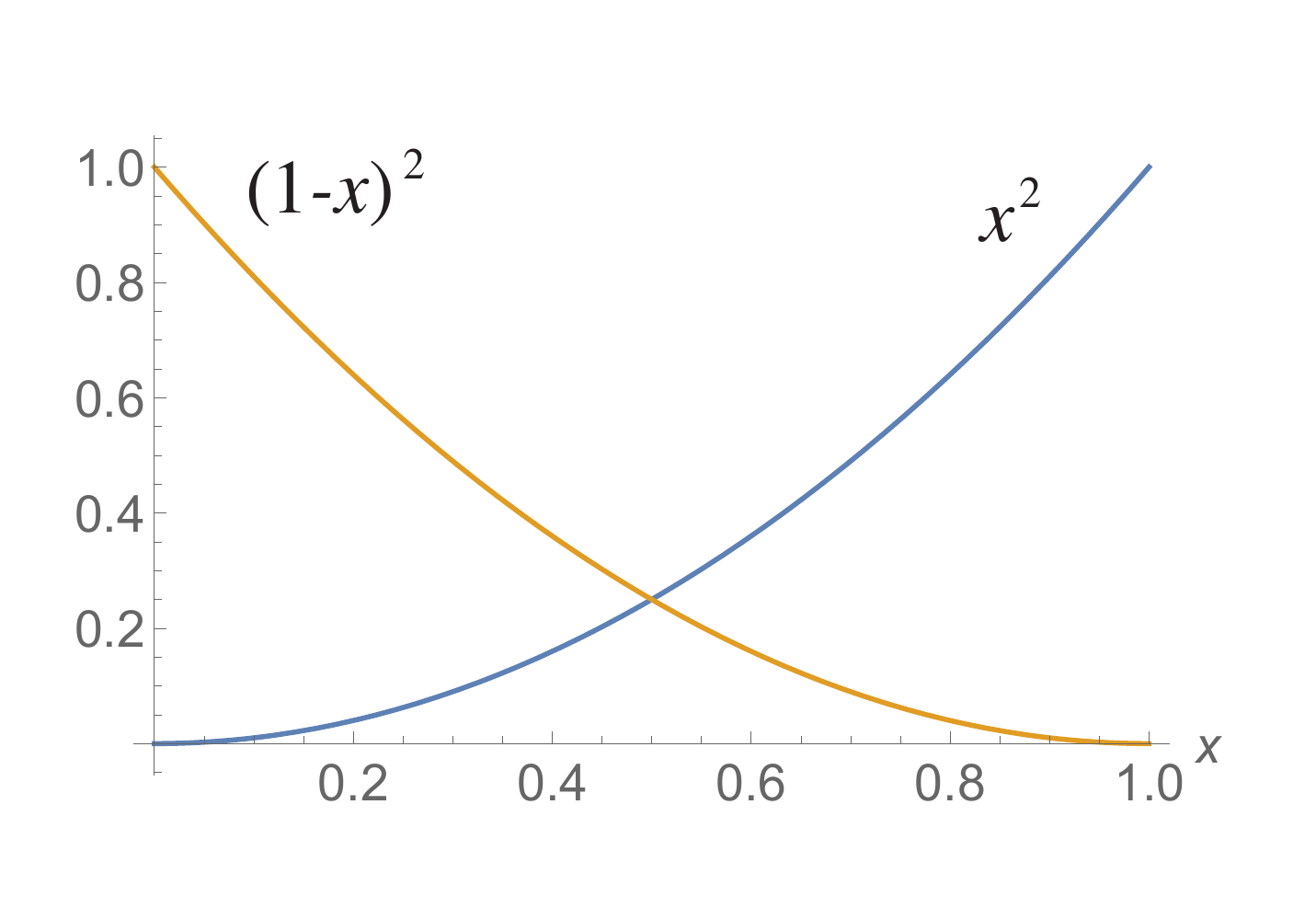}
\caption{The functions $y(x)=x^2$ and $y(x)=(1-x)^2$ }
\label{squares}       
\end{figure}

Before doing any additional analysis, it is worth doing a little thinking to build our intuitive understanding of problems of this sort.  So, let's think about this problem entirely backward for a few moments.  To start, consider the function

\begin{align}
    y(x) &= \alpha y_1(x)  +\beta y_2(x)\label{gen0}   ~~~ x\in [0,1] 
    \end{align}
    \begin{align*}
    \intertext{where}
    y_1(x) &=  x^2 \\
    y_2(x) & = (1-x^2)
\end{align*}
This function is clearly the sum of two independent functions; these independent functions are plotted in Fig.~\ref{squares}

These functions have some interesting features.  First, note that $y_1$ is zero at one end, and 1 at the other; $y_2$ has the same feature, but at opposite ends of the domain.  Also, the derivative $y_1$ is zero at $x=0$, and the derivative of $y_2$ is zero at $x=1$.  These features give lots of flexibility if one wants to develop new functions by creating linear combinations of these two.  For example, suppose I want a function on $[0,1]$ that is (1) a quadradic, (2) has the value 3 on the left-hand side, and has the value of 2 on the right-hand side.  A little though will indicate that we can \emph{create} this function by the following linear combination

\begin{align*}
    y(x) &= 2 y_1(x) + 3 y_2(x)  \\
    &= 2 x^2 + 3(1-x)^2 \\
    &= 5x^2 - 6x +3
\end{align*}

\begin{figure}[t]
\sidecaption[t]
\centering
\includegraphics[scale=.5]{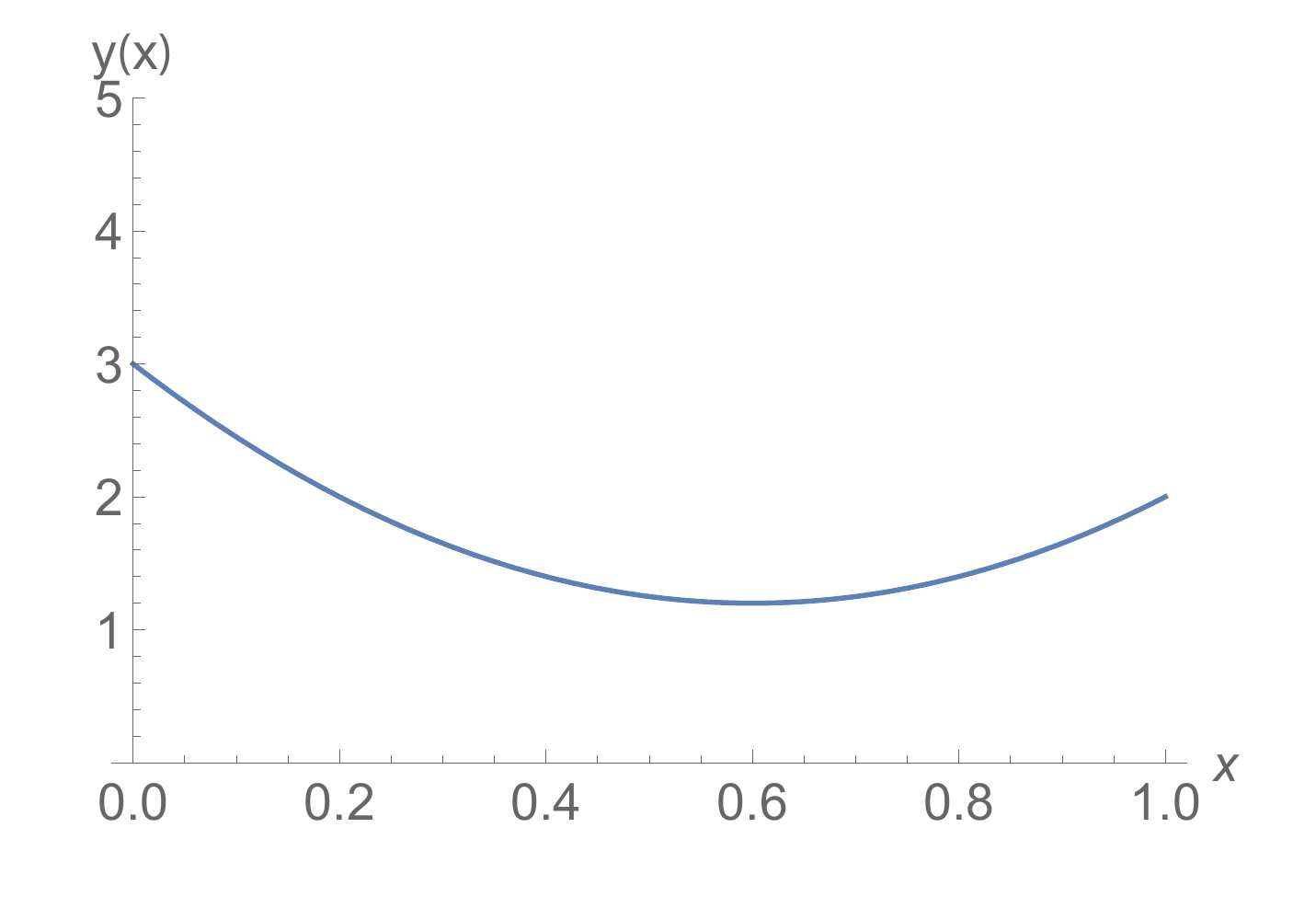}
\caption{The functions $y(x)$ composed specifically so that $y(0)=3$ and $y(1)=2$. }
\label{squares2}       
\end{figure}
If we want to be really clever, we can even indicate the functional \emph{value} at one end, and the \emph{slope} at the other end.  To do this, first note

\[
y'(x) = 2 \alpha x  -2\beta (1-x)   ~~~ x\in [0,1] \\
\]
So, now to meet two conditions (one for the value at one end, one for the slope at the other end) we can solve these two equations simultaneously.  For example, suppose we wanted $y(0)=1, ~y'(1)=20$. This implies
\begin{align*}
    1&=0+\beta(1-0) &&\rightarrow \beta=1&& \mbox{~(from the equation for $y(x)$)}\\
    20&=2\alpha+0 &&\rightarrow \alpha=10&& \mbox{~(from the equation for $y'(x)$)}
\end{align*}

\begin{figure}[t]
\sidecaption[t]
\centering
\includegraphics[scale=.48]{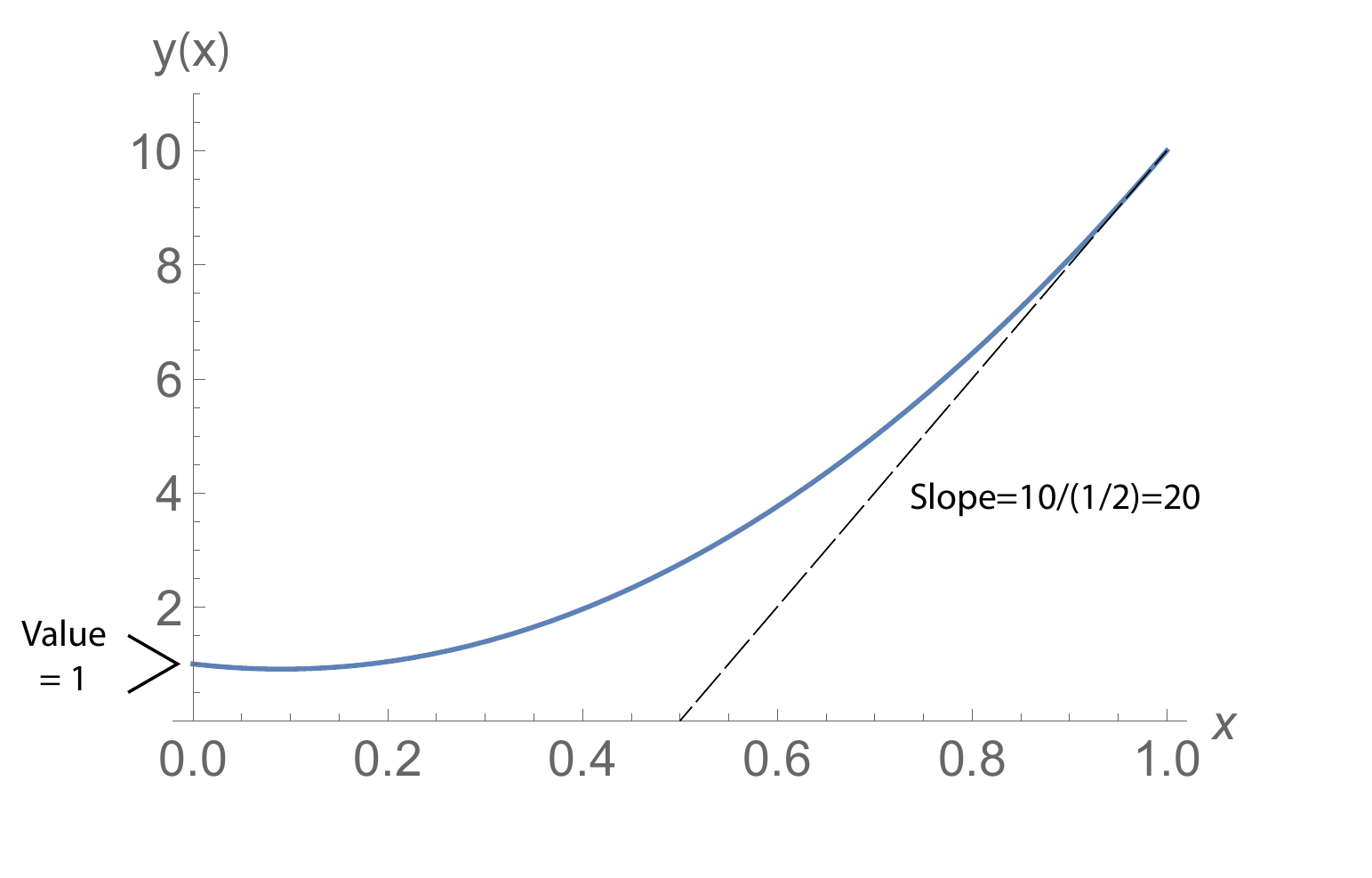}
\caption{The functions $y(x)$ composed specifically so that $y(0)=1$ and $y'(1)=20$. }
\label{squares3}       
\end{figure}

So far, we have not said a word about the solutions to ODEs.  However, the functions that we just examined have relevance here.  To see this, we just have to consider the solutions to the following nonhomogeneous problem.  This is a directly integrable problem (as described in the previous section), but with a particularly simple form.

\begin{align*}
    \frac{d^2 y}{dx^2} = k_0
\end{align*}
where $k_0$ is some specified constant.

Without any ancillary information, we can still solve this problem.  To do so, we integrate both sides to give

\begin{align*}
    \int_x\frac{d^2 y}{dx^2}dx &= \int_x k_0 dx \Rightarrow \\
    \frac{d y}{dx} &=  k_0 x + C_1
\end{align*}
and, integrating a second time, we find

\begin{align}
    \int_x  \frac{d y}{dx} dx &= \int_x (k_0 x+C_1) dx \Rightarrow \\
     y(x)& =  \frac{k_0}{2} x^2 + C_1 x + C_2
     \label{gen2}
\end{align}
It turns out that the function that we examined earlier, $y(x) = \alpha x^2 +\beta (1-x)^2$, is equivalent to this solution. To see this, note that if we expand and regroup Eq.~\eqref{gen0}, we have

\[ 
y(x)=(\alpha + \beta)x^2 -2 \beta x + \beta
\]
Comparing the two solutions at the $x^2$ term, we find that we must have $k_0=2(\alpha+\beta)$ for the two solutions to match.  Suppose we re-examine the case where we had $y(0)=3,~y(1)=2$.  Solving Eq.~\eqref{gen2} with these constraints gives

\begin{align*}
    y(0)&=3 \Rightarrow C_2 = 3 \\
    y(1) &=2 \Rightarrow C_1 = -(1+k_0/2)
\end{align*}
The constraint we have on compatibility between the two equations does not allow \emph{any} values for $\alpha$ and $\beta$; the value of these are related to $k_0$ by $2(\alpha+\beta)=k_0$.  This means that our problem with $y(0)=3,~y(1)=2$, which had $\alpha = 3$, $\beta = 2$, corresponds to $k_0 = 10$. Thus, we have $C_1 = -6$.  It is easy to show that this solution is the same solution as plotted in Fig.~\ref{squares2}.  Alternatively, we can say that the first solution that we plotted (where we forced the right-hand value to be 2 and the left-hand value to be 3)

\begin{align*}
    y(x)&=2x^2 +3(1-x)^2 \\
    \intertext{Corresponds to the second-order ODE}
    \frac{d^2 y}{dy^2}&= 10 \\
    y(0) &= 3 \\
    y(1)&= 2
\end{align*}

So, why did we bother with this long example?  Well, for one thing, it showed us, empirically at least, that we do indeed generate two constants of integration for second-order ODEs.  It also showed us that we can think of the solution to such problems as the linear combination of \emph{two linearly independent} solutions.  In the example above, we thought of these two solutions as $y_1(x)=x^2$ and $y_2(x)=(1-x)^2$.

\subsection{Linear Second-Order ODEs with constant Coefficients: Homogeneous Case}\indexme{ordinary differential equation!second-order homogeneous}

Now, let's examine the general case of linear second-order equations with constant coefficients.  Recall the general form for a \emph{linear} second-order ODE given by Eq.~\eqref{general_second}.  When the term $c(x)\equiv 0$, then the equation is called \emph{homogeneous}.  Sometimes, functions like $g(x)$ are called \emph{source} terms.  Because functions like this are not dependent on $y$, and because the entire term is independent from $y$, the action of these functions can be though of as ``sources" or ``sinks" that drive the problem.

We will study the particular \emph{homogeneous} case where the coefficients of the differential equation are constants rather than functions.  This means that the equations of interest to us take the form

\begin{equation}
    \frac{d^2 y}{dx^2}+ b \frac{d y}{dx} + c y  = 0
    \label{homog_second}
\end{equation}
First, we note that \emph{every} homogeneous second-order ODE with constant coefficients can be put in this form.   Second, we note that
the solution to an equation like this must involve two unknown constants.  We can think of there being an unknown constant generated for each ``integration" that needs to be done to invert the derivatives.  Here, we are using the term integration somewhat loosely, but it is a correct interpretation.  Because the highest order of derivative is two, this means two integrations need to be done.  Finally, because of our previous experience, we might expect that the solution to this equation is the linear combination of two independent solutions.  We have not yet proved this (we will offer a proof of sorts later), but this idea is consistent with our observations so far.  

We expect that there are two linearly independent solutions to the general homogeneous problem with constant coefficients, and we know that we can make new solutions by making linear combinations of them.  However- as to what the solutions actually might be, we have not much to go on.  We also know that for the directly integrable homogeneous problem

\begin{equation}
    \frac{d^2 y}{dx^2}=0
\end{equation}
the solution must be linear of the form $y(x)=C_1 x + C_2$.  So, whatever more general solution we find for Eq.~\eqref{homog_second}, it must have the linear case as one of the possibilities (we will see, ultimately, that it does).  

So, to proceed further, we are going to use a process that is not discussed very often in mathematics:  We are going to use our intuition.  If we \emph{look} at Eq.~\eqref{homog_second}, we note right away that we are looking for some function $y(x)$ such that it is at least possible to take its derivatives, multiply those by constants, and sum them to get zero.  In fact, for the case $b=0$, we even know that we must have a function whose second derivative is some constant multiple of the function that we started with.  Clearly all polynomials are out as possibilities.  Only a few potential functions that behave this way come immediately to mind.  First, we know that functions like $\sin x$ and $\cos x$ behave this way; so those seem like reasonable possibilities.  Another function that can behave this way is the the exponential function with complex exponent.  In fact, because of Euler's identity ($e^{i x}= \cos x + i \sin x$), the exponential function even subsumes the $\sin x$ and $\cos x$ functions.  So, as literally a \emph{guess} (albeit, and educated one), trying $y(x) = e^{s x}$ (where $s$ could potentially be a real or a complex number) seems like a good place to start.  Trying this, we find the following result

\begin{align*}
    &\frac{d^2}{dx^2}e^{s x}+ b \frac{d}{dx}e^{s x}+ c e^{s x} =0 \\
    &s^2 e^{s x} + b s e^{s x} + c  e^{s x} = 0
\end{align*}
So far, we are off to a good start.  Now, we just have to determine under what conditions this last expression might possibly ever be true.  Because $s^2 e^{s x}$ can never be zero at any point (even if $s$ is complex), then we can safely divide both sides by this value. This simplifies things dramatically 

\begin{align*}
    &s^2  + b s  + c   = 0
\end{align*}
We are trying to determine if, given constant (real number) values for $b$ and $c$, there is a solution to the above equation in terms of $s$.  But the expression is now just a quadratic equation in $s$, so we can use the quadratic formula to solve it.  Recall, the solutions are (noting that for the form $a x^2 + bx +c$, we have $a=1$)

\begin{align}
    s_1 = \frac{-b + \sqrt{b^2-4 c}}{2} \\
    s_2 = \frac{-b - \sqrt{b^2-4 c}}{2}
\end{align}
This result tells us that, except for the case of a repeated root (which we will cover later), we do indeed have two independent solutions for $y(x)$.  These are $y_1(x) = e^{s_1 x}$ and $y_2(x) = e^{s_2 x}$.  We can combine these two solutions linearly to generate the most general solution possible.  To do this, we need the \emph{principle of superposition} proved below.

\begin{theorem}[Principle of Superposition for Second-Order ODES]
Suppose $y_1(x)$ and $y_2(x)$ are solutions to the ODE 
\[
    \frac{d^2 y}{dx^2}+ b \frac{d y}{dx} + c y  = 0
\]
Then, the linear combination $y(x) = \alpha y_1(x)+\beta y_2(x)$ is also a solution.
\begin{proof}
Let $L=\frac{d^2 }{dx^2}+ b \frac{d }{dx} + c $.  Then $L[\alpha y_1+\beta y_2] = \alpha L[y_1]+\beta L[y_2]$ by linearity.  Because $L[y_1]$ and $L[y_2]$ are both solutions, they are both equal to zero.  This proves that $L[\alpha y_1+\beta y_2] =0$.
\end{proof}
\end{theorem}

Given all the information above, we can now write our general solution to be of the form

\begin{equation}
    y(x) = \alpha e^{s_1 x}+\beta e^{s_2 x}
\end{equation}
A little more thought about this solution indicates that we can actually say a bit more.  The possibility for roots in a quadratic equation are as follows.

\begin{enumerate}
    \item Case 1.  The roots $s_1$ and $s_2$ are distinct real numbers.
    \item Case 2.  The  roots $s_1$ and $s_2$ are real, but not distinct ($s_1=s_2$).
    \item Case 3.  The roots are complex conjugates of the form $s_1=u+i \lambda$, $~s_2=u-i \lambda$.
\end{enumerate}
These are the only three options.   Because the roots of the equation $s^2  + b s  + c   = 0$ basically define the kind of solution you have, this equation is frequently called the \emph{characteristic equation}.  The easiest way to proceed from here is to consider each option in sequence.

\subsubsection{Case 1: The Roots are Real and Distinct}

For this case, the solution is relatively uncomplicated; it is just a linear combination of two  exponential functions with real exponents.  Thus, the solution is

\begin{equation}
    y(x) = \alpha e^{s_1 x} + \beta e^{s_2 x}
\end{equation}
where $\alpha$ and $\beta$ are two constants that are determined by the two ancillary conditions.  

\begin{svgraybox}
\begin{example}[Contaminant Degradation in a River]
When organic carbon is dumped into a river, it causes a decrease in oxygen in the river as the carbon is degraded by microorganisms in the river water.  Suppose a vegetable processing plant is allowed to put its effluent into a river at a concentration of 20 mg/L of organic carbon (See Figure).  Assume that organic carbon is discharged constantly (24 hours a day), mixed immediately, and that the effective concentration in the river of this readily-degradable organic carbon is $u(x=0)=u_0=$2 mg/L at the point of discharge.  Also assume that the convection-dispersion-reaction equation is applicable, and that all processes have reached steady state.

How far downstream will the carbon affect oxygen levels?  We will interpret the ``affecting" oxygen levels, to mean that the organic carbon concentration is equal to or above $u=0.05$ mg/L. Assume that the velocity is $v_0=1.2$ km/h, the dispersion coefficient is $D_0 = 0.7$ km$^2$/h, and first-order degradation rate is equal to $k_1=0.3$ h$^{-1}$.  Determine a set of reasonable boundary conditions based on the physical organization of the problem in order to develop a solution.

{
\centering\includegraphics[scale=.65]{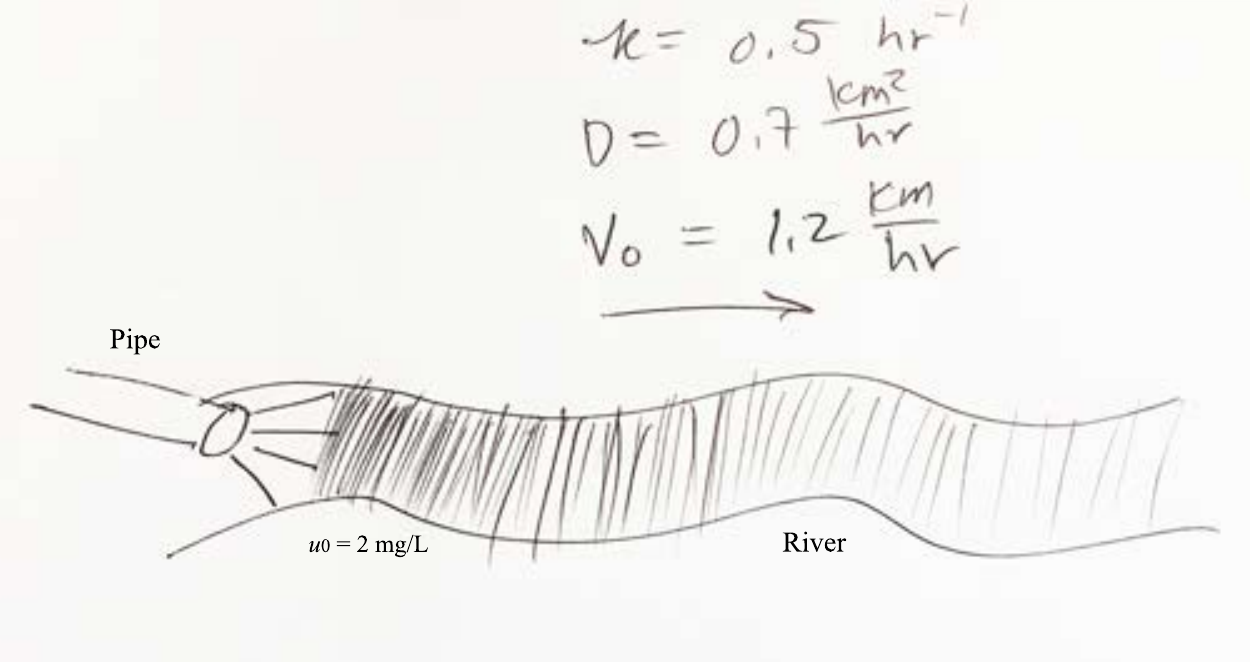}
\captionof{figure}{A stream.  Carefully drawn by hand.}
\label{stream1}  
}

\noindent{\bf Solution.}  The problem statement suggests that the appropriate balance equation takes the form

\begin{align*}
    D_0 \frac{d^2 u}{d x^2} - v_0 \frac{d u}{d x}- k_1 u(x) &= 0 \\
    \intertext{Or, upon rearranging a bit,}
    \frac{d^2 u}{d x^2} - \frac{v_0}{D_0} \frac{d u}{d x}- \frac{k_1}{D_0} u(x) &= 0 \\
\end{align*}
Taking $a=1$, $b=-v_0/D_0$, and $c=-k_1/D_0$.  Using the quadratic formula, we have the roots
\begin{align*}
    s_1 &= \frac{1}{2}\left(\frac{v_0}{D_0} + \sqrt{\left(\frac{v_0}{D_0}\right)^2+4\frac{k_1}{D_0}}\right) \\
    s_2 &= \frac{1}{2}\left(\frac{v_0}{D_0} - \sqrt{\left(\frac{v_0}{D_0}\right)^2+4\frac{k_1}{D_0}}\right)
\end{align*}
The general solution is, then
\begin{align*}
    u(x) &=\alpha \exp\left(\frac{v_0}{2 D_0} x \right) \exp\left(\frac{1}{2}x\sqrt{\left(\frac{v_0}{D_0}\right)^2+4\frac{k_1}{D_0}} \right) \\
    &+
    \beta  \exp\left(\frac{v_0}{2 D_0} x\right) \exp\left(-\frac{1}{2}x\sqrt{\left(\frac{v_0}{D_0}\right)^2+4\frac{k_1}{D_0}}\right)
\end{align*}
In order to determine the solution, we need to impose some ancillary conditions to more fully define the problem (in this case we might call them \emph{boundary} conditions since they apply to space).  The condition at $x=0$ is not too difficult to work out.  Clearly, we would like the concentration there to be $u_0$, so we can set $u(0)=u_0$.  For the second ancillary condition, things are not necessarily as clear.  However, a little thinking would indicate the following: The concentration downstream from the source can only \emph{decrease} not increase.  This is because the problem is essentially one of conservative transport at a constant velocity (with dispersive spreading) plus a decay reaction.  Or, more simply, we might impose the condition that our concentration can not become arbitrarily high downstream from the source, although it can go to zero.  Both of these conditions are consistent with the idea that the constant for the first exponential term must be zero.  Why?  Because if it is not zero, then the concentratio will grow exponentially large as one goes downstream from the source- this does not make any \emph{physical} sense based on our intuition about the problem.  A more ``mathematical", but otherwise equivalent statement might be this: we expect the concentration of carbon to become arbitrarily small as $x$ becomes arbitrarily large.  Either way, we now have the solution
\begin{align*}
    u(x) &=\beta \exp\left(\frac{v_0}{2 D_0} x \right) \exp\left(\frac{1}{2}x\sqrt{-\left(\frac{v_0}{D_0}\right)^2+4\frac{k_1}{D_0}} \right) \\
    u(0)&=u_0
\end{align*}
Evaluating this solution at $x=0$ and using the boundary condition at $x=0$, we rapidly find the final solution
\begin{align*}
    u(x) &=u_0 \exp\left(\frac{v_0}{2 D_0} x \right) \exp\left(-\frac{1}{2}x\sqrt{\left(\frac{v_0}{D_0}\right)^2+4\frac{k_1}{D_0}} \right) \\
\end{align*}
A plot of this function appears below.

{
\centering\includegraphics[scale=.5]{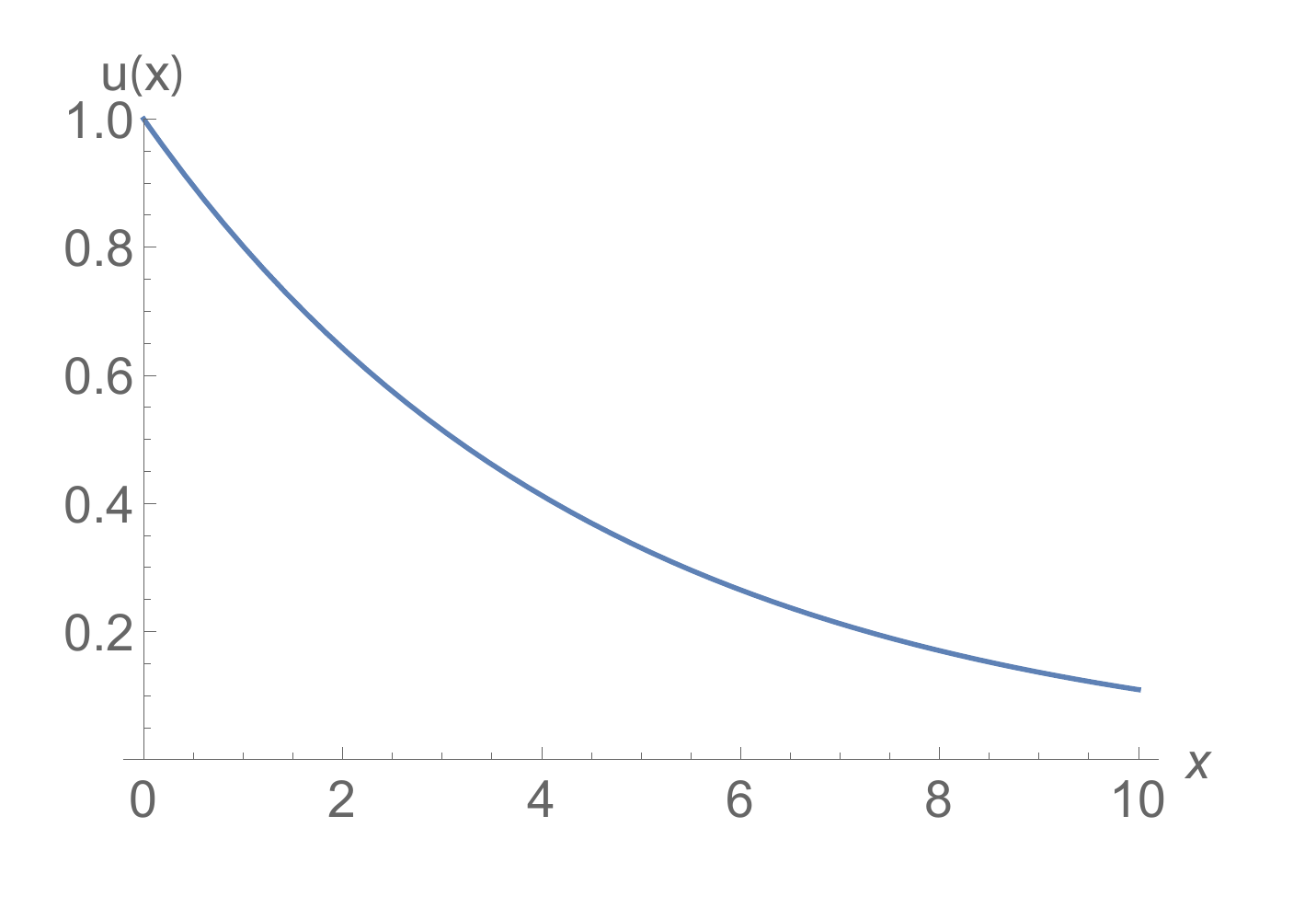}
\vspace{-5mm}
\captionof{figure}{Plot of solution, $u(x)$}
\label{stream2}  
}

The last part of the question asks for the distance at which the carbon content becomes acceptable- i.e., less than 0.05 mg/L.  We can determine this by setting

\begin{align*}
    \frac{u(x)}{u_0}&= \exp\left(\frac{v_0}{2 D_0} x  -\frac{1}{2}x\sqrt{\left(\frac{v_0}{D_0}\right)^2+4\frac{k_1}{D_0}} \right) \\
    \intertext{Or, taking the natural log of both sides}
    \ln\left(\frac{u(x)}{u_0}\right) &= \left(\frac{v_0}{2 D_0} x  -\frac{1}{2}x\sqrt{\left(\frac{v_0}{D_0}\right)^2+4\frac{k_1}{D_0}} \right)
    \intertext{Solving for $x$ gives}
     x&= \ln\left(\frac{u(x)}{u_0}\right)\left(\frac{v_0}{2 D_0}  -\frac{1}{2}\sqrt{\left(\frac{v_0}{D_0}\right)^2+4\frac{k_1}{D_0}} \right)^{-1}
     \intertext{And, substituting in the parameter values from above, we have the solution}
     x&= \ln\left(\frac{0.05}{2}\right)\left(\frac{1.2~km/h}{2\times 0.7~ km^2/h}  -\frac{1}{2}\sqrt{\left(\frac{1.2~km/h}{0.7~km^2/h}\right)^2+4\frac{0.5~h^{-1}}{0.7~km^2/h}} \right)^{-1} \\
     x&=16.62~km
\end{align*}
\end{example}
\end{svgraybox}

\subsubsection{Case 2: The Roots are Real, but Not Distinct}
This is the case where the terms in the radical of the characteristic equation $s^2  + b s  + c   = 0$ ,  is zero.   Or, equivalently, we have that $b^2=4 c$, so that the only root is $s=-b/2$.  This presents a little bit of a problem, because we know that we need \emph{two independent} equations to form the general solution to a second-order ODE.  

The resolution to this problem is to find a second solution.  Recall, the way that we found the exponential solutions to begin with was somewhat empirical; we thought about what kinds of functions could possibly be combined in a way that the archetype for a second-order homogeneous ODE might possibly be valid.  The exponential was a candidate that turned out to work.  We do know that the first solution generates exponentials of the form $y_1(x)=\alpha e^{-b/2 x}$.  Lacking any definite direction, we might hope that our second solution could be one that is proportional to this first solution.  For example, we could suggest

\[ y_2(x) = f(x) e^{-b/2 x} \]

And hope that this works.  Well, to be quite honest, this is roughly how the process unfolded.  If we are guessing a functions that might work for $f(x)$, we might try the simplest thing we can think of.  

\[ y_2(x) = x e^{-b/2 x} \]

Of course, this suggestion wouldn't have been made if it wasn't going to work.  To see that it does, we need only take some derivatives

\begin{align*}
    y_2(x) &=  x e^{-x b/2 } \\
     y_2'(x) &=  e^{- x b/2 } -x\frac{b}{2}e^{-x b/2 }\\
      y_2''(x) &=  -\frac{b}{2} e^{-x b/2 } -\frac{b}{2}e^{-x b/2 } +x\left(\frac{b}{2}\right)^2e^{-x b/2 }\\
      &=-be^{- x b/2 } +x\frac{b^2}{4}e^{-x b/2 }
\end{align*}
Now, substitute this into
\begin{equation}
    \frac{d^2 y}{d x^2}+ b \frac{d y}{dx}+ cy = 0
\end{equation}
To give
\begin{align*}
   & -be^{-b/2 x} +x\frac{b^2}{4} e^{-b/2 x}
    +b e^{-b/2 x} -x\frac{b^2}{2}e^{-b/2 x}
    +c x e^{-b/2 x} = 0
    \intertext{And, after simplifying, this is}
    & 
    \left[-\frac{b^2}{4}+c 
    \right] x e^{-b/2 x}=0
\end{align*}
Finally, recalling that for this case we have $b^2=4c \Rightarrow -b^2/4+c=0$, we have proved that $y_2$ is a solution.  

\begin{svgraybox}
\begin{example}[Case 2 Example]
Consider the following steady diffusion problem 
\begin{align*}
    D\frac{d^2 u}{d x^2} &= 0 \\
    u(0) &= u_0 ~~~~\text{(This is a specified concentration boundary condition)}\\
    \left. -D\frac{du}{dx}\right|_{x=L} &= 0~~~~~~\text{(This is a \emph{no-flux}  boundary condition.)}
\end{align*}
Show that this solution is a special case of the general solution for Case 2. \\

\noindent{\bf Solution.} For this example, we have $a=1$ and $b=c=0$.  Therefore, the only roots of the characteristic equation are $s_1 = s_2 =0$.  This is a repeated root (even though it is zero), meaning that the general solution is
\begin{align*}
    u(x) &= \alpha e^0+\beta x e^0 \\
    &= \alpha + \beta x
\end{align*}
Note that this solution is actually a line, which is indeed a special case of the general exponential solution.  Using the two boundary conditions, we have
\begin{align*}
    u_0& = \alpha \\
    0  &= \beta
\end{align*}
Thus, our final solution is 
\[ u(x) = u_0 ~\text{(a constant)} \]
\end{example}
\end{svgraybox}

\subsubsection{Case 3: The Roots Are Complex Conjugates}

For this case, the characteristic equation ($s^2 + b x + c$) is such that $b^2-4 c <0$; therefore, one is faced with complex roots coming from the quadratic equation.  i.e., 
\begin{align*}
    s_1 &= -\frac{b}{2}+\frac{1}{2}\sqrt{b^2 -4 c} \\
    s_2 &= -\frac{b}{2}-\frac{1}{2}\sqrt{b^2 -4 c}
\end{align*}
Because $b^2-4 c <0$, the radical generates an imaginary number.  For convenience, set $u=-b/2$ and $i \lambda = 1/2\sqrt{b^2-4c}$.  Then, the two roots can be put in the form

\begin{align*}
    s_1 &= \mu+i \lambda \\
    s_2 &= \mu-i \lambda
\end{align*}
which is a complex conjugate pair ($s_1 = \overline{s_2}$).  Our general solution is, therefore 

\begin{align}
    y(x) & = \alpha e^{\,\mu x+i \lambda x}+  \beta e^{\,\mu x-i \lambda x} \nonumber \\
    &= e^{\,\mu x}\left( \alpha e^{i \lambda x}+ \beta e^{-i \lambda x} \right) 
\end{align}
Recalling Euler's identity, we can write the two independent solutions as 

\begin{align*}
    y_1(x)=e^{\,\mu x}e^{i x}= e^{\,\mu x}\left[\cos (\lambda x) + i \sin (\lambda x)\right] \\
    y_2(x)=e^{\,\mu x}e^{-i x}= e^{\,\mu x}\left[\cos (\lambda x) - i \sin (\lambda x)\right]
\end{align*}
Now, if we were satisfied with solutions that, in general, contained complex numbers (and, in some cases, such as quantum mechanics, we might very well be!), then we would be done.  However, for the vast majority of problems in science and engineering, the problems are such that we expect a real number as an answer.  This is actually no problem.  We can, if we like, construct two new solutions from the existing one by creating linear combinations of them.  Recall, that when we make linear combinations, we can multiply either of the two equations by any constant we like- including complex numbers.  Suppose, then, that we construct two new independent solutions as follows

\begin{align*}
    y_A(x) &= e^{\,\mu x}\frac{1}{2}\left[ y_1(x)+y_2(x)\right] \\
    &= e^{\,\mu x}\frac{1}{2}\left[(\cos (\lambda x)+i \sin (\lambda x))+(\cos (\lambda x)-i \sin (\lambda x))  \right] \\
    &=e^{\,\mu x}\frac{1}{2}(2 \cos (\lambda x)) \\
    &=e^{\,\mu x}\cos (\lambda x)
    \intertext{and}
    y_B(x) &=e^{\,\mu x}\frac{i}{2}\left[ y_1(x)-y_2(x)\right] \\
    &= -e^{\,\mu x}\frac{i}{2}\left[(\cos (\lambda x)+i \sin (\lambda x))-(\cos (\lambda x)-i \sin (\lambda x))  \right] \\
    &=-e^{\,\mu x}\frac{i}{2}(2 i \sin (\lambda x)) \\
    &= e^{\,\mu x}\sin (\lambda x)
\end{align*}
With a little thought, we might have realized this from the start.  Our two independent solutions can be considered to be either $e^{\,\mu x}\sin x$ and $e^{\,\mu x}\cos x$ or $e^{\,\mu x}e^{ix}$ and $e^{\,\mu x}e^{-ix}$.  It actually makes no difference mathematically- one set of solutions can always be constructed from the other.  However, for applications where we desire solutions that are real functions (rather than complex functions), we are better off using the trigonometric solutions.  Thus, we can take our general solution to be in this case

\begin{align}
    y(x) &= \alpha y_A(x)+\beta y_B(x) \nonumber \\
    &= \alpha e^{\,\mu x} \cos (\lambda x) +\beta e^{\,\mu x} \sin (\lambda x) \nonumber\\
    \intertext{Or, factoring out the exponential}
    y(x)&=e^{\,\mu x}\left( \alpha  \cos (\lambda x) +\beta  \sin (\lambda x) \right)
\end{align}

\begin{svgraybox}
\begin{example}[Harmonic Oscillator.]
\emph{Harmonic oscillator} is just a very complicated term for describing something that bounces periodically.  This could be the pendulum of a clock, a weight on a spring, or the CEO or Pfizer (slogan: ``Pfizer Quality!") inexplicably bouncing on a trampoline in the Mojave Desert on a moonless Wednesday night in late August.  Suppose we have the weight on a spring version (sorry).  Initially, the spring is pulled by the amount $10~cm$.  Suppose the weight has a mass of $1 kg$, and the spring constant for the spring is $K=3~N/cm$.  Assume that the system is ideal so that it does not ever lose energy (which, admittedly, is impossible; but, it is often very close to reality for short periods of time!).  What is the function defining the oscillations of the weight?  Assume that the motion takes place on a frictionless horizontal plane (say, on an air hockey table) so that you do not need to worry about gravity.\\

{
\centering\includegraphics[scale=.5]{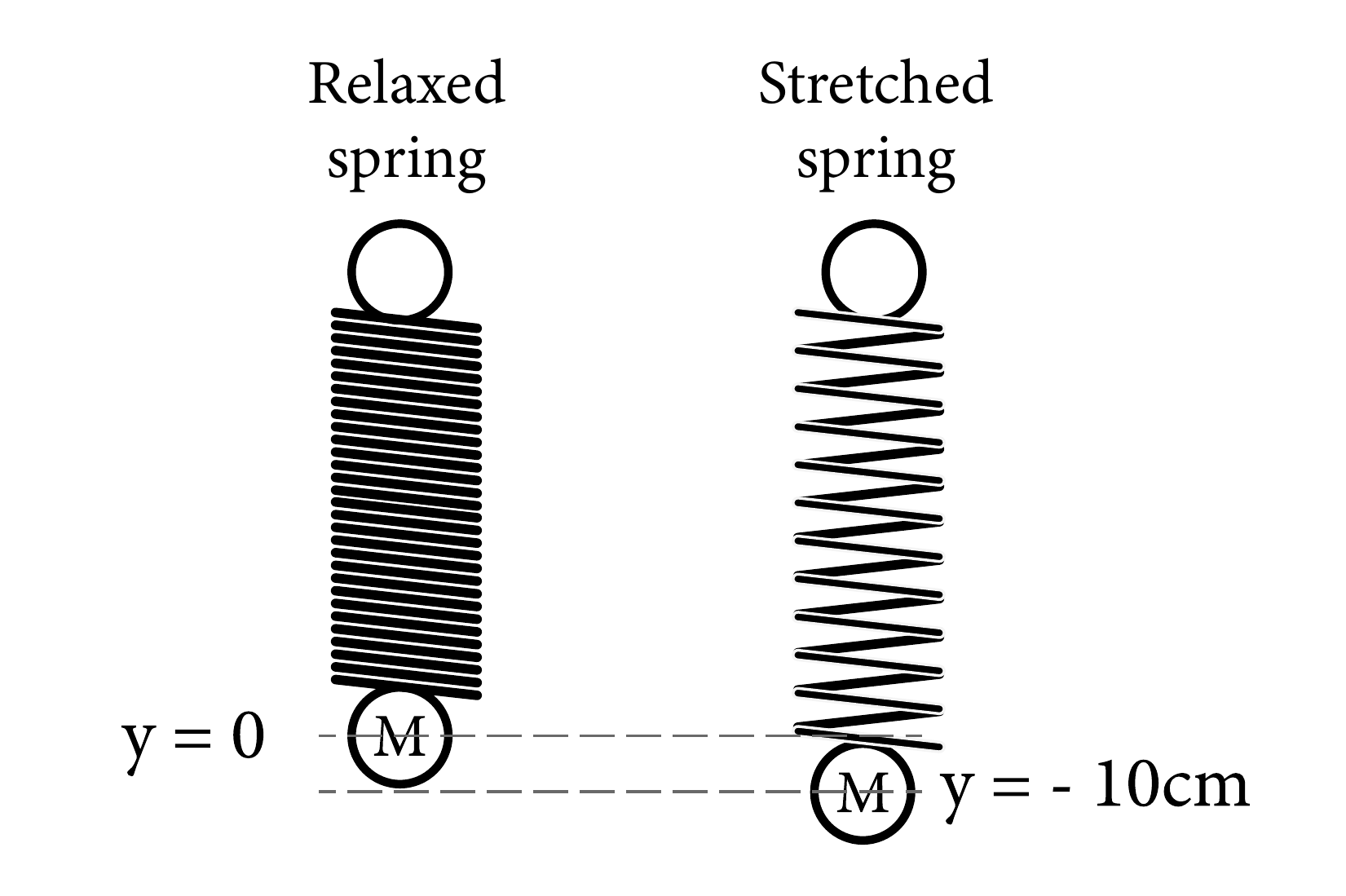}
\captionof{figure}{A weight on a spring- with motion in the horizontal plane.  On an air hockey table.}
\label{spring}  
}

\vspace{5mm}
\noindent {\bf Solution.} A simple force balance on the mass just as it is let go is as follows
\begin{align*}
    m\frac{d v}{dt}=  K \Delta y
\end{align*}
where $dv/dt$ is the acceleration (time rate of change of velocity), $g$ is the gravitational constant ($9.81~m/s^2$), and $\Delta y$ is the initial displacement $\Delta y = 0-(-10)~cm=10~cm$. Because the equilibrium position is $y(0)=0$, we have that $\Delta y=-y(t)$. Recall that $v(t)=dy/dt$.  With this in mind, our force balance can be written (with two ancillary conditions based on the physical system)
\begin{align*}
  &&  \frac{d^2 y}{dt^2}+\frac{K}{m} y(t) &= 0&& \\
 & condition~1&  y(0)&=-10 &&\text{initial position at y=-10 cm}\\
 & condition~2& y'(0)&=0 &&\text{zero initial velocity}
\end{align*}
The characteristic equation for this problem is  $s^2 + K/m = 0$ ($a=1$, $b=0$, $c=K/m$).  Thus the two roots are
\begin{align*}
    s_1=0+\frac{1}{2}\sqrt{-4\frac{K}{m}} = i\sqrt{\frac{K}{m}} \\
    s_2=0-\frac{1}{2}\sqrt{-4\frac{K}{m}} =-i\sqrt{\frac{K}{m}}
\end{align*}
Using the notation established above, we have $\mu=0$ and $\lambda = \sqrt{K/m}$.  Therefore, the solution that we want is
\begin{align*}
    y(t)& =e^0\left(\alpha \cos (t\sqrt{K/m} ) + \beta \sin(t\sqrt{K/m} ) \right) \\
    \intertext{Taking the derivative (for use with the second ancillary condition)}
    y'(t) &=\left(-\sqrt{K/m}\alpha \sin (t\sqrt{K/m} ) + \sqrt{K/m}\beta \cos(t\sqrt{K/m} ) \right) \\
    \intertext{Using the two ancillary conditions, we have}
    -10 &= \alpha \cos(0)+\beta \sin(0)\Rightarrow \alpha=-10 \\
    0&=-\sqrt{K/m}\alpha \sin (0)+ \sqrt{K/m}\beta \cos (0)\Rightarrow \beta = 0
\end{align*}
So, our solution is
\begin{equation}
    y(t) = -10 \cos(t\sqrt{K/m} )
\end{equation}
Does this solution make sense?  We can check.  First of all, we should check to see if the initial conditions (the two ancillary conditions) are met.  Note: $y'(x)=-10\sqrt{K/m} \sin(t \sqrt{K/m})$.
\begin{align*}
    y(0)&=-10 \cos(0) \\
    y'(0)& = -10 \sqrt{K/m} \sin(0) = 0
\end{align*}
So, yes, our solution does meet the initial conditions.  It is easy to verify that it also meets the ODE (take derivatives, and substitute back into the original ODE). A plot of the oscillations for $0 < t < 60~s$ is given below. 

{
\centering\includegraphics[scale=.75]{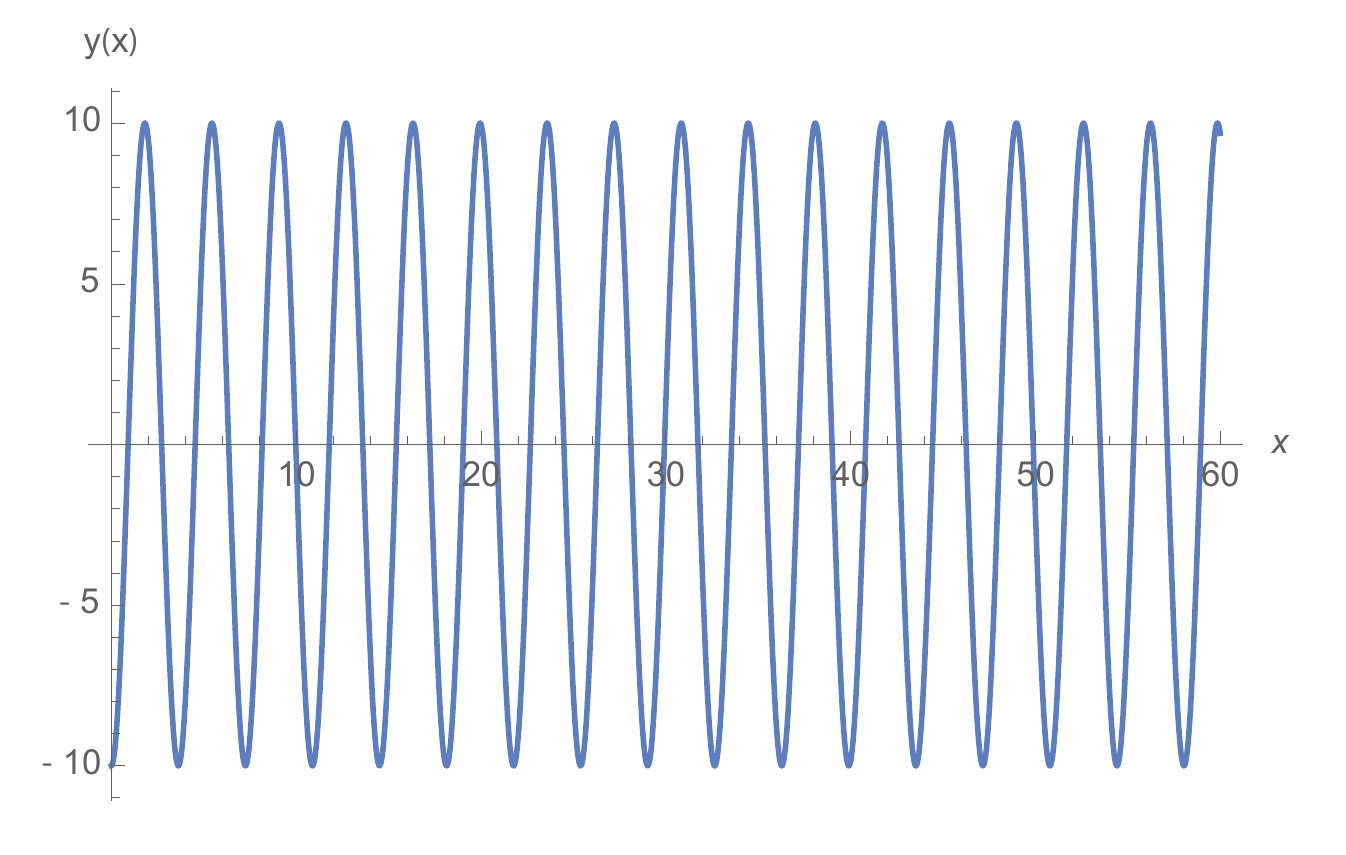}
\vspace{-5mm}
\captionof{figure}{Plot of solution for the weight on a spring.}
\label{bounce}  
}
\end{example}
\end{svgraybox}

\subsection{Solutions for Nonhomogeneous Second-Order ODEs with Constant Coefficients}\indexme{ordinary differential equation!second-order nonhomogeneous}\indexme{ordinary differential equation!variation of parameters}\label{variationofparams}

We will study only one approach for nonhomogeneous second-order ODEs.  In particular, we are getting around to solving the problem given by 

\begin{equation}
    \frac{d^2 y}{d x^2}+b\frac{d y}{dx}+ cy(x) = g(x)
    \label{nonhomog}
\end{equation}
As a matter of notation, we call any solution to this problem \emph{a particular} solution, $y_p(x)$.

The method that we are going to use is called \emph{variation of parameters}, and it was developed by two somewhat famous mathematicians: Leonard Euler (of Euler identity fame), and Joseph-Louis Lagrange.  The basic idea behind the method is as follows.  Suppose we have the two solutions to the homogeneous second-order equation

\begin{equation}
    \frac{d^2 y}{d x^2}+b\frac{d y}{dx}+ cy(x) = 0
    \label{homog2}
\end{equation}
We call the general solution to this problem the \emph{homogeneous} solution \index{ordinary differential equation!homogeneous solution}.  Because we will make use of $y(x)$ to indicate the solution to Eq.~\eqref{nonhomog}, we will use $y_h(x)$ to indicate the homogeneous solution
\begin{equation*}
    y_h(x) = \alpha y_1(x) +\beta y_2(x)
\end{equation*}
Recall that we previously defined the operator notation for general second-order ODEs to be 

\begin{equation}
  L =  \frac{d^2 }{d x^2}+b\frac{d }{dx}+ c = 0
\end{equation}
so that $L[y_h]=0$ gives us Eq.~\eqref{homog2}.  This is really nothing but a shorthand notation.  So, now we consider the problem $L[y_p]=g(x)$-- the nonhomogeneous problem.   There are an infinite number of particular solutions, $y_p$.  To see this, suppose we have \emph{any} particular solution to Eq.~\eqref{nonhomog}.  We need only examine the sum 
\begin{align*}
    y(x) &= y_p(x) + \gamma y_h(x) \\
    &= y_p(x) + \alpha y_1(x)+\beta y_2(x) 
\end{align*}
where $\gamma$ is a constant.  Then note

\begin{align*}
    L[y(x)] &= L[y_p(x)] + \gamma L[y_h(x)] \\
    \intertext{But, by definition $L[y_h]=L[\alpha y_1(x)+\beta y_2(x) ]=0,$ and $L[y_p]=g(x)$, so}
    L[y(x)]&= g(x)  
\end{align*}
So, if we have any particular solution, we can always add any multiple of the homogeneous solution, and we have a new solution.  The specific reason that we might need to add in the homogeneous solution is that the particular solution that we find might not meet the necessary ancillary conditions.  However, we know we can find specific values of $\alpha$ and $\beta$ in the homogeneous solution to meet the possible boundary conditions.  Thus, we can use the added homogeneous solution as a way of adjusting our particular solution so that it meets the necessary boundary conditions.

The technique of variation of parameters can be explained in terms of the theory of Green's functions, but we will not be taking that approach here (for one thing, this requires some understanding of distribution theory).  Alternatively, let's think about the properties of the solution that we seek, $y(x)$.  Clearly $y(x)$ cannot be a constant multiple of either $y_1(x)$ or $y_2(x)$, because then it would be a solution to the homogeneous equation.  In other words, we realize that 

\[ \frac{y(x)}{y_1(x)} \ne C_1 ~~~~~~~~~~ \frac{y(x)}{y_2(x)} \ne C_2 \]
where $C_1$ and $C_2$ are constants.  The argument at this juncture seems ridiculously trivial.  If the ratios above cannot be equal to constants, then they must be equal to functions.  Specifically, functions of $x$, since that is the only independent variable we are considering.  Thus, we argue

\[ \frac{y(x)}{y_1(x)} = u_1(x) ~~~~~~~~~~ \frac{y(x)}{y_2(x)} = u_2(x) \]
For some (currently unknown) functions $u_1$ and $u_2$.  This is equivalent to stating that a solution for $y_p(x)$ must be of the form

\[ y_p(x) = u_1(x) y_1(x) + u_2(x)y_2(x)  \label{particularform}\]
Now, if you really think about it, at this juncture we haven't really done much.  In fact, if we have any three bounded and continuous functions,$f,g$, and $h$, it is not difficult to show that we can always find two functions $u_1$ and $u_2$ such that $h=u_1 f+u_2 g$.  The real brilliance of Eq.~\eqref{particularform} is that we already know that $L[y_1]=L[y_2]=0$.  So, working with Eq.~\eqref{particularform} to find the functions $u_1$ and $u_2$ is going to be made much simpler by this fact.  This becomes clearer when we actually try to do this.  We are looking for the solution

\begin{align*}
    L[y_p] &= g(x) \\
    \intertext{or, substituting}
    L[u_1(x) y_1(x) + u_2(x)y_2(x)]= g(x) \\
\end{align*} 
Or, writing this all out more painfully
\begin{align}
    &\frac{d^2}{dy^2}[u_1 y_1+u_2 y_2]+ b \frac{d}{dy}[u_1 y_1+u_2 y_2] + c[u_1 y_1+u_2 y_2] = g(x)
    \label{painful}
\end{align}
Let's start by noting
\begin{align*}
     \frac{d}{dy}[u_1(x) y_1(x)]&=u_1'y_1 + u_1 y_1' \\
     \frac{d}{dy}[u_2(x) y_2(x)]&= u_2'y_2 + u_2 y_2'
\end{align*}
so that
\begin{align*}
    \frac{d}{dy}[u_1(x) y_1(x)+ u_2(x) y_2(x)]&=u_1'y_1 + u_1 y_1' + u_2'y_2 + u_2 y_2'\\
    &= (u_1'y_1 +u_2'y_2)+(u_1 y_1' +  u_2 y_2')
\end{align*}
and
\begin{align*}
     \frac{d^2}{dy^2}[u_1 y_1+u_2 y_2]&=\frac{d}{dy}[(u_1'y_1 +u_2'y_2)+(u_1 y_1' +  u_2 y_2')] \\
     &=(u_1'y_1 +u_2'y_2)' + u_1'y_1'+u_1 y_1''+u_2' y_2' +u_2 y_2''
\end{align*}

Substituting these results back into Eq.~\eqref{painful}, we have

\begin{align*}
    &(u_1'y_1 +u_2'y_2)' + u_1'y_1'+u_1 y_1''+u_2' y_2' +u_2 y_2'' \\
    +&b[(u_1'y_1 +u_2'y_2)+(u_1 y_1' +  u_2 y_2')] \\
    +&c[u_1 y_1+u_2 y_2] = g(x) 
\end{align*}
From here, we collect all of the terms that are multiplied by $u_1$ or $u_2$.  This result is

\begin{align*}
&u_1(y_1''+b y_1'+ c y_1)+ u_2(y_2''+b y_2'+ c y_2) \\
+& u_1'y_1'+ +u_2' y_2'\\
+&b(u_1'y_1 +u_2'y_2)\\
+&(u_1'y_1 +u_2'y_2)'= g(x) 
\end{align*}
This simplified immediately because $y_1$ and $y_2$ are solutions to the homogeneous problem $y''+by'+cy=0$, leaving us with
\begin{align*}
 u_1'y_1'+ +u_2' y_2'
+b(u_1'y_1 +u_2'y_2)
+(u_1'y_1 +u_2'y_2)'= g(x) 
\end{align*}
At this juncture, we have one equation, but two unknown functions $u_1$ and $u_2$.  This kind of a problem has an infinite number of solutions; we only need one solution, and any particular solution will do.  So, we are free to add any additional equation that we like (as long as it is independent) to make our system solvable.  A particularly convenient requirement is to set $u_1' y_1 +u_2'y_2=0$; you can see that doing so is reasonably clever because it eliminates \emph{two} groups of terms in the equation above!  In summary, then, we have the two equations

\begin{align*}
u_1'y_1 +u_2'y_2&=0\\
 u_1'y_1'+ +u_2' y_2'&= g(x) 
\end{align*}
which is two equations in the unknowns $u_1'$ and $u_2'$.  In fact, we can arrange this in a matrix 

\begin{equation*}
    \left[ {\begin{array}{*{20}{c}}
  {{y_1}}&{{y_2}} \\ 
  {{{y'}_1}}&{{{y'}_2}} 
\end{array}} \right]\left[ {\begin{array}{*{20}{c}}
  {{{u'}_1}} \\ 
  {{{u'}_2}} 
\end{array}} \right] = \left[ {\begin{array}{*{20}{c}}
  0 \\ 
  {g(x)} 
\end{array}} \right]
\end{equation*}
and then use Cramer's rule to solve for $u_1'$ and $u_2'$.  First note that 
\[\left| {\begin{array}{*{20}{c}}
  {{y_1}}&{{y_2}} \\ 
  {{{y'}_1}}&{{{y'}_2}} 
\end{array}} \right|=y_1 y_2'-y_2 y_1' \] 
So, Cramer's rule gives us

\begin{align*}
\frac{du_1}{dx}&=\frac{1}{y_1 y_2'-y_2 y_1'}\left| {\begin{array}{*{20}{c}}
  0&{{y_2}} \\ 
  g(x)&{{{y'}_2}} 
\end{array}} \right| 
& \frac{du_2}{dx}&=\frac{1}{y_1 y_2'-y_2 y_1'}\left| {\begin{array}{*{20}{c}}
  y_1&{0} \\ 
  y_1'&{g(x)} 
\end{array}} \right|\\  
\intertext{or}
\frac{du_1}{dx}&=\frac{-y_2 g(x)}{y_1 y_2'-y_2 y_1'} 
& \frac{du_2}{dx}&=\frac{y_1 g(x)}{y_1 y_2'-y_2 y_1'} 
\end{align*}
Although not essential, we can put this in integral form by integrating both sides of this expression with respect to $x$.  Note that, because we are looking for \emph{any} particular solution, we can discard any constants of integraiton generated at this point.  

\begin{align*}
    u_1(x)&=\int_x \frac{-y_2 g(x)}{y_1 y_2'-y_2 y_1'}\,dx
& u_2(x)&=\int_x \frac{y_1 g(x)}{y_1 y_2'-y_2 y_1'} \, dx
\end{align*}
This completes our solution.  However, to summarize, recall that the two solutions above allow us to generate the following particular solution

\begin{equation*}
    y_p(x)=u_1(x)y_1(x)+u_2(x)y_2(x)
\end{equation*}
To this particular solution, we must add the general form of the homogeneous solution to get our final solution for the nonhomogeneous equation

\begin{align*}
    y(x) = y_p(x) + y_h(x)
    \intertext{or}
    y(x) = y_p(x) + \alpha y_1(x)+\beta y_2(x)
\end{align*}
This last step is crucial!  Without the addition of the homogeneous solution, we do not have two free constants to specify via two ancillary conditions.  As a reminder, because $L[y_h]=0$, adding the homogeneous solution means that we still meet the nonhomogeneous second order ODE 
\[
    \frac{d^2 y}{d x^2}+b\frac{d y}{dx}+ cy(x) = g(x)
\]
because adding the homogeneous solution does nothing but add zero to this ODE! An example is helpful here.

\begin{svgraybox}
\begin{example}[Solution for a nonhomogeneous second-order ODE]
The previous example of a harmonic oscillator is an interesting one to use for this example, in part because we already have the homogeneous solution.  A \emph{forced} harmonic oscillator is one that is forced to oscillate at a frequency different from its \emph{natural} frequency.  Consider the system on the air hockey table described in the previous example.  Suppose our weight is made of magnetic material, and we put a giant electromagnetic array under the table that can create a magnetic force which fluctuates in time according to
\[ g(x) = -10 \cos(10 t) \]
This is different from the natural frequency that we found above that was equal to $-10 \cos(\sqrt{K/m} t)=-10 \cos(\sqrt{3} t)$.  What happens when we do this?  The system of equations that we want to solve is now

\begin{align*}
  &&  \frac{d^2 y}{dt^2}+\frac{K}{m} y(t) &= -10 \cos(10 t)&& \\
 & condition~1&  y(0)&=-10 &&\text{initial position at y=-10 cm}\\
 & condition~2& y'(0)&=0 &&\text{zero initial velocity}
\end{align*}
Recall, our homogeneous solution for this case was
\[    y(t) =\alpha \cos (t\sqrt{K/m} ) + \beta \sin(t\sqrt{K/m} ) \]
Note that here, we have $y_1(t) = \cos (t\sqrt{K/m})$, $y_2(t)=\sin(t\sqrt{K/m} )$.
Following the example above, we are now seeking a solution to the nonhomogenous problem of the form

\[ y_p(t) =u_1(t) \cos (t\sqrt{K/m} ) + u_2(t) \sin(t\sqrt{K/m} ) \]
To find this solution, we need the functions $y_1, ~y_2, ~y_1', ~y_2'$ and $g(t)$.  It is also handy to have the combination $y_1 y_2'-y_2 y_1'$.

\begin{align*}
    y_1&=\cos (t\sqrt{K/m}) \\
    y_2&= \sin(t\sqrt{K/m}) \\
    y_1'&= -\sqrt{K/m}\sin(t\sqrt{K/m}) \\
    y_2'&= \sqrt{K/m}\cos(t\sqrt{K/m}) \\
    y_1 y_2'-y_2 y_1' &= \sqrt{K/m}\cos^2 (t\sqrt{K/m})  +\sqrt{K/m}\sin^2(t\sqrt{K/m})=\sqrt{K/m}\\
    g(t)&=-10 \cos(10 t)
\end{align*}
At this juncture, the solution is a bunch of busy work
\begin{align*}
    \frac{d u_1}{dt }&=\frac{10\sin (t\sqrt{K/m}) \cos(10 t)}{\sqrt{K/m}} \\
    \frac{d u_2}{dt }&=\frac{-10\cos (t\sqrt{K/m}) \cos(10 t)}{\sqrt{K/m}}
\end{align*}
These can be integrated using standard trigonometric techniques, or by looking the integrals up on a table.  The results are
\begin{align*}
    u_1(t) &= 10 \sqrt{\frac{K}{m}} \left(-\frac{\cos \left(t
   \left(\sqrt{\frac{K}{m}}-10\right)\right)}{2
   \left(\sqrt{\frac{K}{m}}-10\right)}-\frac{\cos \left(t
   \left(\sqrt{\frac{K}{m}}+10\right)\right)}{2
   \left(\sqrt{\frac{K}{m}}+10\right)}\right) \\
   u_2(t) &= -10 \sqrt{\frac{K}{m}} \left(\frac{\sin \left(t
   \left(\sqrt{\frac{K}{m}}-10\right)\right)}{2
   \left(\sqrt{\frac{K}{m}}-10\right)}+\frac{\sin \left(t
   \left(\sqrt{\frac{K}{m}}+10\right)\right)}{2
   \left(\sqrt{\frac{K}{m}}+10\right)}\right)
\end{align*}
Well, worse things have happened... Recalling that $K/m=3$ and simplifying, our final solution is
\begin{align*}
    y(t) &= \alpha  \cos \left(\sqrt{3} t\right)+\beta  \sin \left(\sqrt{3} t\right)+\frac{30}{97}
   \cos (10 t) \\
   \intertext{and its derivative is}
y'(t)&=-\sqrt{3} \alpha  \sin \left(\sqrt{3} t\right)+\sqrt{3} \beta  \cos \left(\sqrt{3}
   t\right)-\frac{300}{97} \sin (10 t)
   \end{align*}
From the first ancillary condition we have
\[-10 =\alpha \cos(0)+ \frac{30}{97}\cos(0) \Rightarrow \alpha = -\frac{1000}{97} \]
From the second ancillary condition we have
\[ 0 = \sqrt{3}\beta \cos(0) \Rightarrow \beta=0 \]
Apparently, the solution we want is
\begin{align*}
    y(t) &= -\frac{1000}{97}\cos \left(\sqrt{3} t\right)+\frac{30}{97}
   \cos (10 t)
\end{align*}
\end{example}
\end{svgraybox}

\newpage
\section*{Problems}
\subsection*{Practice Problems}
For the problems below, keep in mind that $y=y(x)$ is the \emph{dependent} variable, and $x$ is the \emph{independent} variable.  Note also that the choice of what we call the independent variable is not important.  Thus, some problems might be posed in terms of other independent and dependent variables (e.g., $u(t)$ and $t$). \\
  
\noindent For the following problems, determine the \emph{general} solution by separating variables.  Also list the dependent variable in each case.  Please solve explicitly for the dependent variable only if it is both possible, and does not require unusual effort (e.g., it would not be expected that you would solve for the roots of any polynomial solution beyond quadratic).

\setlength{\columnsep}{2cm}
\begin{multicols}{2}
\begin{enumerate}[topsep=8pt,itemsep=4pt,partopsep=4pt, parsep=4pt]
    \item $\frac{dy}{dx} = \frac{x}{y}$
    \item $\frac{2y}{y^2+1}\frac{dy}{dx} = \frac{1}{x^2}$
    \item $\frac{dy}{dx} = \frac{y}{x}$
    \item $\frac{dy}{dx} = \frac{\sin x -\cos x}{y^2+y}$
    \item $(y^2-3)\frac{dy}{dx} = 1$
    \item $(1+x^2)\frac{dy}{dx} = (1+y^2)$
    \item $\frac{1+e^t}{1-e^{-y}} \frac{dy}{dt} + e^{t+y} = 0$
    \item $\exp(t+y) \frac{dy}{dt}= \frac{t}{y}$, for $t>0$, $y(t)>0$
    \item $z^3 \frac{du}{dz} = \sqrt{z^2 -u^2 z^2}$
    \item $u'(t) = 4 t u^2(t) $
    \item $y'(t) = m c^2 t^{\pi} y(t)$ where $m$ and $c$ are constants
        \saveenumerate
\end{enumerate}
\end{multicols}

\noindent For the following problems, determine the \emph{general} solution by whatever method is required.  If ancillary conditions are specified, then also find the \emph{particular} solution for those conditions.

\setlength{\columnsep}{2cm}
\begin{multicols}{2}
\begin{enumerate}[topsep=8pt,itemsep=4pt,partopsep=4pt, parsep=4pt]
\restoreenumerate
\item $\frac{dy}{dx} + \frac{1}{x^2}y = 0$, $x\in(0,1]$
\item $\frac{dy}{dt} +2y = 4$, $~y(0)=0$, $~x\in[0,2]$
\item $\frac{dy}{dx} -a y = f(x)$, $~x\in[a,b]$\vspace{1mm}\\
where $f(x)$ is some integrable function.  You will not be able to solve this explicitly!  You will have to leave it in integral form.
\item $\frac{dy}{dz} = 2y + z$, $~z\in [0,\infty]$\vspace{2mm}\\ $y'(0)=1$
\item $y'+  x^5 y = x^5$
\item $y'(t)=y(t) + \frac{1}{1-e^{-t}}$
\item $\frac{dy}{dt} + y \cos t = \cos t$, $~t\in[0,\infty]$, $y(0)=5$
\item $\frac{dy(t)}{dt} - y(t) -t \sin(t) = 0$
 \saveenumerate
\end{enumerate}
\end{multicols}

\noindent For the following problems, determine the \emph{general} solution for the second order ODE.

\setlength{\columnsep}{2cm}
\begin{multicols}{2}
\begin{enumerate}[topsep=8pt,itemsep=4pt,partopsep=4pt, parsep=4pt]
\restoreenumerate
\item $y''(t) - 5 y'(t) -14 y(t) = 0$
\item $y''(x) -25 y(x) = 0$
\item $\frac{d^2 y}{d x^2} + 2 \frac{d y}{dx} -24 y = 0$
\item $\frac{d^2 y}{d x^2} - 6 \frac{d y}{dx} +9 y = 0$
\item $\frac{d^2 y}{d x^2} + 10 \frac{d y}{dx} +25 y = 0$
\item $\frac{d^2 y}{d x^2} + 2\sqrt{2} \frac{d y}{dx} +2 y = 0$
\item $\frac{d^2 y}{d x^2} + 3\frac{d y}{dx} +3 y = 0$
\item $\frac{d^2 y}{d x^2} + 2 \frac{d y}{dx} +5 y = 0$
\item $6\frac{d^2 y}{d t^2} -4 \frac{d y}{dt} = 0$
\item $6\frac{d^2 y(t)}{d t^2} -4 y(t) = 0$
 \saveenumerate
\end{enumerate}
\end{multicols}

\noindent For the following problems, determine the \emph{particular} solution for the second order ODE and the two given ancillary conditions.

\setlength{\columnsep}{2cm}
\begin{multicols}{2}
\begin{enumerate}[topsep=8pt,itemsep=4pt,partopsep=4pt, parsep=4pt]
\restoreenumerate
\item $\frac{d^2 y}{d x^2} - 2 \frac{d y}{dx} -15 y = 0$, $~y(0)=1$, $~y'(0)=-1$, $~x\in[0,\infty]$
\item $\frac{d^2 y}{d x^2} +16  y = 0$, $~y(0)=3$, $~y'(0)=12$, $~x\in[0,\infty]$
\item $\frac{d^2 y}{d x^2} - 10 \frac{d y}{dx} +25 y = 0$, $~y(0)=0$, $~y(1)=0$, $~x\in[0,1]$
\item $4y''(t) -4 y'(t) +y(t) =0$,$~y(0)=1$,$~y'(0)=0$,$~t\in[0,\infty]$
\item Show that if $y_1(x)$ and $y_2(x)$ are each solutions to the ODE
\begin{equation*}
    \frac{d^2 y}{d x^2}+ b \frac{d y}{d x}  +cy = 0
\end{equation*}
then, for any constants $C_1$ and $C_2$, the function $y(x) = C_1 y_1(x) + C_2 y_2(x)$ is also a solution.  To start, substitute this last expression for $y(x)$ directly into the ODE.  Group terms containing $y_1$ and terms containing $y_2$. 
 \saveenumerate
\end{enumerate}
\end{multicols}

\noindent For the following problems, determine the solution for the \emph{nonhomogeneous} second order ODE.  If ancillary conditions are provided, solve the problem for the unknown coefficients.

\setlength{\columnsep}{2cm}
\begin{multicols}{2}
\begin{enumerate}[topsep=8pt,itemsep=4pt,partopsep=4pt, parsep=4pt]
\restoreenumerate
\item $y''(t)-y(t)'-2y(t) =2e^{-t}$, $y(0)=1$, $y'(0)=0$
\item $\frac{d^2 y}{d x^2} +2\frac{d^2 y}{d x^2}+y = 3e^{-t}$, $y(0)=0$, $y'(0)=1$
\item $\frac{d^2 y}{d x^2} -\frac{d^2 y}{d x^2}+\tfrac{1}{4}y = 4e^{t/2}$, $y(0)=1$, $y'(0)=1$
\saveenumerate
\end{enumerate}
\end{multicols}

\subsection*{Applied and More Challenging Problems}

\begin{enumerate}
    \item Repeat the skier problem (Example problem \ref{skier}), but with the following modification.  The net drag due to friction from both the interface with the snow and from air resistance will be accounted for by a term proportional to the velocity:  $F_{drag}=-\alpha m v(t)$.  Add this force into the force balance, and then re-solve the problem using $\alpha = 1$ s$^{-1}$. 
    

   \item For this problem involving a CSTR (of volume $\forall$) with flow rate $Q$, and with first-order reaction, determine how a sudden change in the influent concentration, $y_{in}$ would be manifest in the CSTR given that it has an initial concentration equal to $y_0$.  
    \begin{enumerate}
        \item The mass balance and initial condition are given by 
        \begin{align}
            \frac{dy}{dt} +(\frac{1}{\tau}+k)y &= \frac{1}{\tau} y_{in}\\
            y(t=0)=y_0
        \end{align}
        Recall here that $\tau=\forall/Q$ (the hydraulic residence time), $\forall$ is the volume of the reactor, $Q$ is the flow rate at both the inlet and the outlet, $y_{in}$ is the influent concentration, $y_0$ is the initial concentration in the tank, and $k$ is the first-order reaction rate constant.
        
        Solve this problem for $y(t)$.  You may find it helpful to make the substitution $\beta = (\frac{1}{\tau}+k)$ to simplify the analysis (and prevent mistakes from dealing with complex algebraic computations), but do not forget to convert back to the original variables if you do so.
        
        \item For the conditions $\forall=$10 m$^3$, $Q=$1 m$^3$/hr, $k=$0.4 hr$^{-1}$, $y_0=$2 moles/m$^3$ and $y_{in}=$1 moles/m$^3$, plot the solution to the problem for $0 \le t \le 100$ hr.  On the same plot, plot the solution with identical parameters, except with $k=0$ hr$^{-1}$.  Comment on the primary difference in the two solutions.  What is the steady-state solution ($t\longrightarrow \infty$) for each of the problems?
    \end{enumerate}

    \item Catalysts are sometimes subject to degradation.  This can be for many reasons. Some examples include thermal degradation, poisoning by other reactive chemicals (that changing the oxidation state of catalytic site or bind to the site irreversibly), and fouling by components in the treatment stream (physically altering the catalyst).  
    
    Suppose we have the catalytic reaction $A\mathop  \to \limits^C B$), where species $A$ is converted to species $B$ in the presence of the active catalyst sites, represented by species $C$.  Assume that the catalysts is being degrades due to the reaction process occurring at a low pH.  
    
    We can model a catalytic reaction by creating balances for both the catalyst ($y_C$) and the chemical being catalysed (chemical species $A$).  A possible set of balances for catalytic degradation during the production of the three chemical species (reactant $A$, product $B$, and catalyst $C$) is as follows.
    \begin{align}
        \frac{d y_C(t)}{dt} &= -k_C y_C(t) \\
        y_C(0) &= C_0 \\
        \frac{d y_A(t)}{dt} &= -k_A y_A(t) y_C(t) \\
        y_A(0) &= A_0 \\
        \frac{d y_B(t)}{dt} &= k_B y_A(t) y_C(t) \\
        y_B(0) &= B_0
    \end{align}
    We can solve this problem in a step-wise fashion.  First, note that the balance for the catalyst, $C$, does not depend at all on the concentrations of the other two species.  Thus, an explicit expression for the concentration $y_C(t)$ is straightforward to evaluate.  For the problem, do the following.
    \begin{enumerate}
        \item Solve the balance for the catalyst to determine an expression for $y_C(t)$.
        \item Substitute your expression for $y_C(t)$ into the right-hand side of the balance for species $A$.  Solve this expression for $y_A(t)$.
        \item With the solution for $y_A(t)$ and $y_C(t)$ determined, the final step is to substitute these into the right-hand side of the balance for species $B$.  Solve the resulting expression for $y_B(t)$.  You will probably want to use a symbolic integration software like Mathematica to compute the associated integral.
        \item Plot the resulting function on the interval $0\le t \le 20$ min for the following list of constants.  To help make the graph easier to view, plot the normalized values $y_A/A_0$, $y_B/A_0$, and $y_
        C/C_0$.
        \begin{align}
            &k_A = 0.2 ~\textrm{min}^{-1}  && A_0= 10 ~\textrm{mmol/L} \\
            &k_B = 0.3 ~\textrm{min}^{-1}  && B_0= 0 ~\textrm{mmol/L} \\
            &k_C = 0.25 ~\textrm{min}^{-1}  && C_0= 1 ~\textrm{mmol/L}
        \end{align}
    \end{enumerate}

    \item Sometimes problems are separable, but the resulting solution is \emph{implicit} in the dependent variable.  As an example, consider the following problem for degradation of a chemical species by enzyme kinetics in a batch reactor.  If Michaelis-Menten type kinetics apply, an initial value problem for a well-stirred batch system can be stated by
    \begin{align}
        \frac{d y}{dt} &= - \mu \frac{y}{y+K}\\
        y(0) &= y_0
    \end{align}
    where $y$ is the concentration of the chemical species (mass/volume),  $k$ is the reaction rate parameter (mass/time) and $K$ is the half saturation constant (mass/volume).  For this problem, do the following.
    \begin{enumerate}
        \item Solve the problem by simple separation of variables.  Include the initial condition by either solving using definite integration, or by evaluating the constant of integration from the initial condition if indefinite integration is used.
        \item The result is explicit in the variable $t$, but implicit in the variable $y$.  Suppose we have $y_0=20$ mmol/L, $\mu=0.1$ mmol$\cdot$s$^{-1}$, and $K=8$ mmol/L.  
        Suppose you wanted to solve the equation for values of the concentration between $0\le y\le 20$, for the time period of approximately $0 \le t\le 50$ s.  How could you do this?  \emph{Hint.  Try reversing the roles of the dependent and independent variables, treating time as if it were the dependent variable)}.  
        \item Make a plot of your results using whatever plotting software is convenient for you.  If you use Mathematica, you may find the plotting command \texttt{ParametricPlot} useful.  You can read about how to use that in the help pages for Mathematica.
    \end{enumerate}

    
    \item Suppose we have a completely stirred reactor that, because of a problem with the pump, the flow rate decreases in time.  Specifically, assume \label{transientQ}
    \begin{equation}
        Q(t) = \frac{Q_0}{\alpha + \gamma t}
    \end{equation}
    Assume that the initial condition is $y(0)=0$.
    Where $\alpha$ and $\gamma$ are parameters that control the decrease in flow rate over time.  The revised mass balance equation can be taken for this case as
    \begin{equation}
        \frac{dy}{dt} = \frac{1}{\tau}\left(\frac{1}{\alpha + \beta t}\right) y_{in} -\frac{1}{\tau} \left(\frac{1}{\alpha + \beta t}\right) y(t) -k y(t)
    \end{equation}
    \begin{enumerate}
        \item Begin by rewriting this expression so that you have a single term on the left-hand side that can be associated with $P(t)$, and a single term on the right hand side associated with $Q(t)$.  Separate from the ODE, write down the explicit functions representing $P(t)$ and Q(t).  
        
        \item Determine the integrating factor by first computing $s(t)$ from
        \begin{equation}
            s(t) = \int P(t) dt
        \end{equation}
        Then, compute explicit expressions for $e^{s(t)}$ and $e^{-s(t)}$.  You should need to use the logarithmic identity $a \ln (x) = \ln (x^a)$.
        
        \item Finally, using Mathematica (or similar symbolic integration software), determine the result
        \begin{equation}
            y(t) = e^{-s(t)}\left[ \int Q(t) e^{s(t)} dt + C_1\right]
        \end{equation}
        Computing this integral is challenging, and will lead to a solution in terms of the Gamma function that was introduced earlier in the text.  Note that in Mathematica, you will want to compute this integral by putting conditions on the constants and variables, as follows.\\
        
        \texttt{Integrate[Q(t), \{t, 0, t\}, 
 Assumptions $->$\{\{t, k, t\} $\in$ Reals \&\& \{t, k, tau\} $\ge$ 0\}]}\\
 
 \noindent(where, here, you will substitute the appropriate expression for $Q(t)$).
 
        \item  Noting the initial condition $y(0)=0$, find the appropriate value for $C_1$.
        
        \item Finally, using the values $k=0.1$ min$^{-1}$ and $\tau$ = 0.5 min, $a=1$ min, and $b=1$, plot the solution over the range $0\le t \le 100$ min using whatever software is conveneint.  Your solution should look like the plot in Fig.~\ref{f:transientQ}.    Does the resulting function match your intuition about what should happen to the effluent concentration over time?  Why or why not?
    \end{enumerate}
\begin{figure}[t]
\sidecaption[t]
\centering
\includegraphics[scale=.4]{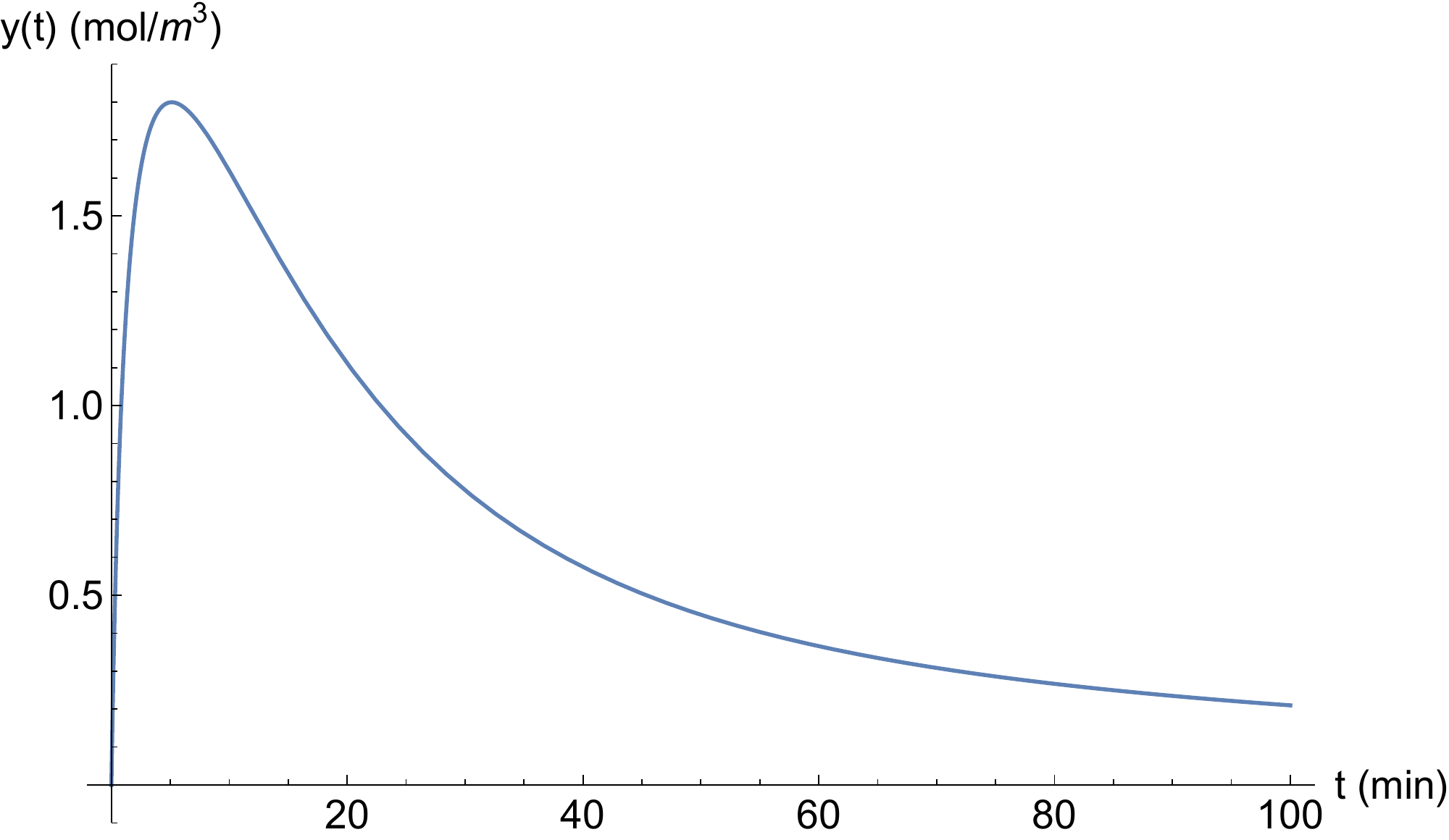}
\caption{Solution to Applied Problem \ref{transientQ}. }
\label{f:transientQ}       
\end{figure}

\item Dispersion is the spreading that occurs during chemical species transport resulting from the variations in the velocity field plus the effects of molecular diffusion.  Sometimes for a plug flow reactor, one models the reactor \emph{without} dispersion; the idea is that for plug flow reactors, this is one of the goals (to reduce spreading).\\

For this problem, assume we have a 1 m long plug flow reactor that catalyses a reactant $A$ to a product $B$, and model the system at steady state both with and without the process of hydrodynamic dispersion.  Suppose for now, we are interested only in the concentration evolution of species $A$ (the reactant).  The problem is outlined as follows.

 Without hydrodynamic dispersion, the steady state differential mass balance equation is given by the first order expression
    \begin{align}
       - u_0 \frac{d y}{dx} &= k y(x) \\
       y(0) &= y_0
    \end{align}
    
 Including hydrodynamic dispersion, the steady state differential mass balance equation is given by the first order expression
    \begin{align*}
      D \frac{d^2 y}{dx^2} - u_0 \frac{d y}{dx} &= k y(x) \\
       y(0) &= y_0 \\
       y'(L) &= 0
    \end{align*}
(note, we need one additional ancillary (boundary) condition to determine the constants, because now the problem is second-order). To begin, for each problem divide through by whatever constant makes the coefficient of the \emph{highest derivative} equal to 1.  This step often makes problems a bit easier to think about.   Please remember not to substitute any numeric values into the solution until the solution is complete.
    
\begin{enumerate}
\item Solve the first problem representing convection without dispersion.  Note that to solve this problem, it is useful to consider the residence time of the system.  The system is $L$ units long with a fluid velocity of $u_0$. The time it takes a parcel of fluid entering the system at $t=0$ to exit the system is 

\begin{equation}
    L = u_0 \tau
\end{equation}
where $\tau$ is defined as the residence time for a plug flow reactor.  Clearly, we have that the residence time is just the length divided by the velocity, $\tau = L/u_0$.  Note also that there is a relationship between the variables $x$ and $t$ in this problem.  Any location $x$ can be converted into an equivalent residence time (up to that point) by the relationship $t =x/u_0$.  Thus, for a purely convective system, it is possible to convert between times and spatial locations. 

\item Using $D=1\times 10^{-9}$~m$^2$ s$^{-1}$, $u_0=1\times 10^{-8}$ m s$^{-1}$, $k=1\times 10^{-8}$ s$^{-1}$, and $y_0=100$ mmol/L, solve this problem over the interval $0 < x < L$, with $L=1$ m.  Plot the two solutions using whatever software is convenient.  

\item How do the two solutions compare?  Why do you think they are different from one another?
\end{enumerate}


\item Steady state heat transport in a radiator fin, assuming that the fin is thin compared to its length, can be described by the following second-order ODE.  Here, we are imposing boundary conditions at the left of a specified temperature ($y(0) = y_0 = 400$ K), and the right boundary as a zero gradient condition.

\begin{align*}
    K\frac{d^2 y}{dx^2}-\frac{h}{\ell} y &= -\frac{h}{\ell} y_\infty  \\
    y(0) &= y_0 \\
    \left.\frac{d y}{dx}\right|_{x=L}&= 0
\end{align*}
Physically, this equation represents the temperature, $y$, in the fin assuming that there is both heat conduction, and convective cooling by an external air source blowing over the radiator fin.  

Assume the fin is $L=0.1$ m long, with thermal conductivity $K=385$ W m$^{-1}$K$^{-1}$, and with heat transfer coefficient $h=150$W m$^{-2}$K$^{-1}$.  The thickness of the fin is $\ell=0.002$ m.   Assume the far-field temperature, $y_\infty$ is 293 K.  Solve this problem.  Plot the solution over the domain and range $0<x<L$ and $293 < T < 400$.  Please remember not to substitute any numeric values into the solution until the solution is complete.
\end{enumerate}
\abstract*{This is the abstract for chapter 3.}

\begin{savequote}[0.55\linewidth]
``Fourier's Theorem \ldots is not
only one of the most beautiful results of modern analysis, but
may be said to furnish an indispensable instrument in the treatment of nearly every recondite question in modern physics."
\qauthor{Lord Kelvin, Treatise on Natural Philosophy (1912)}
\end{savequote}


\chapter{Fourier Series Part I: Introductory Concepts}\label{FS_1}\indexme{Fourier series}
%
\def\CHAP {chapter03_fourier_series1}
Fourier series are one of the most widely used, and powerful tools in all of applied mathematics, with applications that range from signal analysis, to compression algorithms for making smaller computer files. Before investigating Fourier series in detail, it will be useful to revisit some basic concepts about infinite series in general.  Again, it is helpful to start with some terminology.

\section{Terminology}
In the study of Fourier series, a number of new concepts arise.  For convenience, some of the more important definitions and vocabulary are summarized here.

\begin{itemize}

\item \textbf{Fourier series.} A series expressed in terms of a weighted sum of sine or cosine functions (or both).  The sine and cosine functions are of increasing frequency, and the weight functions represent the corresponding amplitude for that function.  As an example, a general Fourier sine series on $x\in[0,1]$ is given in terms of the amplitudes $B_n$ and frequencies (modulated by the integer $n$) as follows
\[ f(x) = \sum_{n=1}^\infty B_n \sin(n Pi x) \]
\\
%
\begin{figure}[t]
\sidecaption[t]
\centering
\includegraphics[scale=.5]{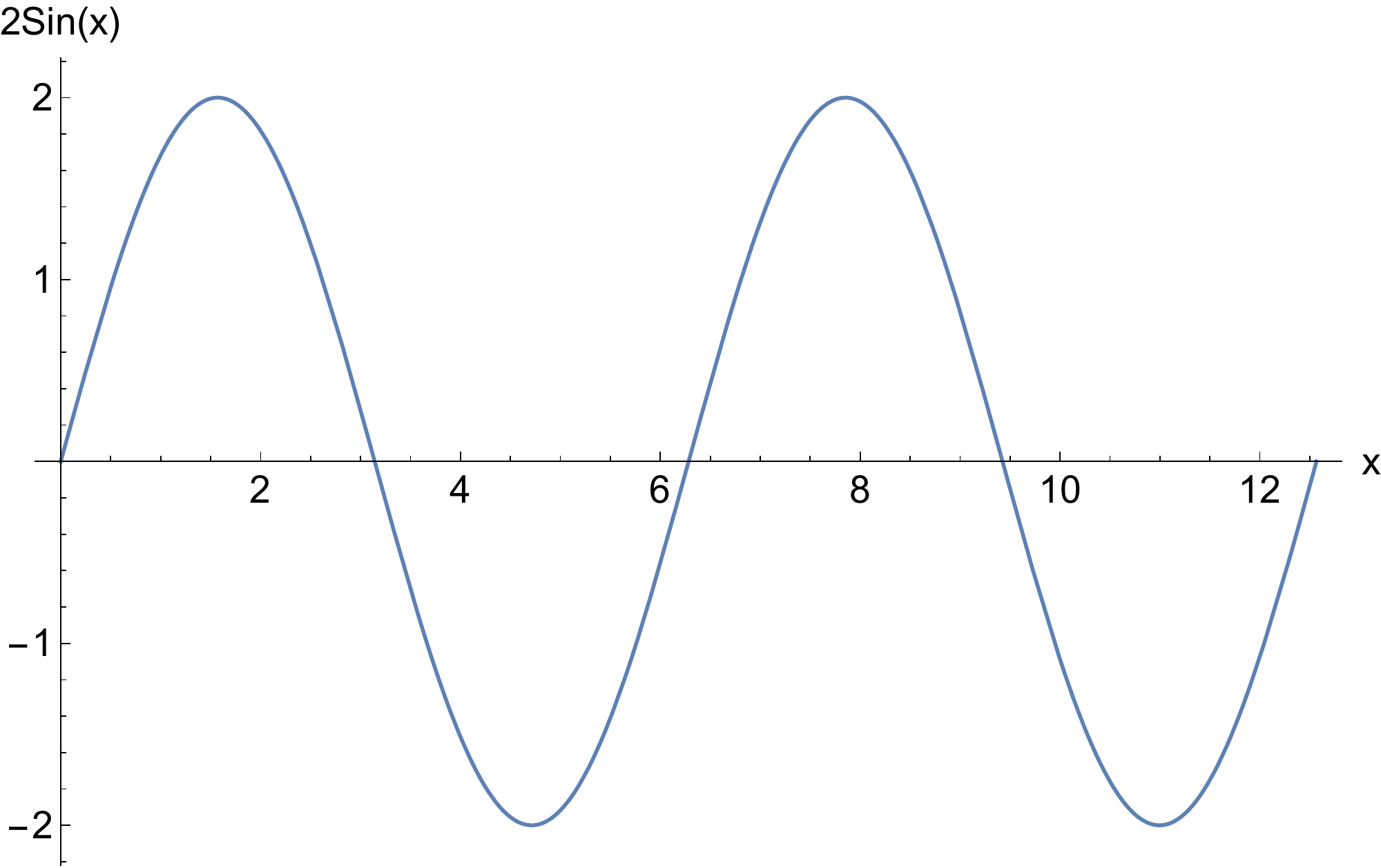}
\caption{ The periodic function $f(x)=\sin(x)$ plotted over the interval $0\le x \le 4 \pi$.  The amplitude of this function is $B=2$. The period of this function is $P=2\pi$; the function repeats itself every multiple of $2\pi$ so that $f(x)=f(n 2\pi x)$ for all positive integers, $n$. }
\label{sine}       
\end{figure}
\item \textbf{Amplitude.}  The amplitude of a periodic function is the maximum value of its vertical magnitude.  For a sine function in the form $B \sin x $, the coefficient $B$ represents its amplitude.  For example, if $B=2$, then the function oscillates between $\pm 2$ with changing $x$.\\

\item \textbf{Period.}  A periodic function, $f(x)$, is said to have a period $P_\lambda$ if $f(x)=f(x+nP)$ for all positive integers $n$.  In other words, the period is the distance (space or time) along the $x-$axis that defines the repeating unit of the periodic function.  When the independent variable is space, sometimes $P_\lambda$ is called the \emph{wavelength} and given the symbol $\lambda$.  The subscript $\lambda$ is added to the the symbol for the period ($P_\lambda$) and the frequency ($F_\lambda$, defined below) so that these symbols are defined uniquely (and as a reminder of their connection with the concept of wavelength).  For sine and cosine functions in the form $\sin(n \pi x)$ or $\cos(n \pi x)$, the period is given by 
\begin{equation*}
    P_\lambda = \frac{2 \pi}{n}
\end{equation*}
\\

\item \textbf{Frequency.} If the independent variable is time, the period is sometimes expressed as the \emph{frequency}.  Because of its frequent use, it has become common to speak of the frequency even when the independent variable is not time.  Frequency is best thought of as the number of complete cycles of the periodic function per unit length or unit time.  The relationship between period, $P_\lambda$, and frequency, $F$ is

\[F = \frac{1}{P_\lambda}\]

Therefore, if the Period is $P_\lambda=2 \pi$, the frequency is $F_\lambda=1/(2 \pi)$.  This indicates that there is one full cycle (one period) occurring every $2\pi$ units of the independent variable.  For sine and cosine functions in the form $\sin(n \pi x)$ or $\cos(n \pi x)$, the period is given by 

\[F = \frac{1}{P_\lambda}=\frac{n}{2 \pi}\]

The SI units for frequency are the Hertz (Hz) which measures the number of periods completed per second.  However- take caution!  The terminology \emph{period} is sometimes used even when the independent variable is something other than time.\\

\item \textbf{Convergence.} Convergence of series is a complicated topic.  We discussed the concept of convergence briefly for Taylor series.  For Fourier series, the functions that we examine will all converge in some useful sense.  For continuous functions, the series will converge pointwise.  This means that for any $x$ chosen in the domain, the partial sums evaluated at $x_0$, ($S_n(x_0)$),  get as close to the value of the function at $x$ ($f(x_0)$) as we like by increasing $n$ sufficiently.  For discontinuous functions, the Fourier series converge in a sense that will be discussed later in the chapter.

\item \textbf{Spectrum}.  We will find that Fourier series for a function $f$ are are (in general) infinite series containing weighted sum of sine functions, cosine functions, or both.  For each such function, the series is indexed by an integer, $n$, as described above where $n$ specifies both the frequency and the amplitude.  The \emph{spectrum} of the function $f$ is then given by the discrete (but infinite...) plot of amplitude versus the value of $n$.  

\item \textbf{Convergence}.  The property of a series (or sequence) to attain a finite limit as the number of terms tends towards infinity.  In Chapter 2, two modes of convergence for series were discussed-- \emph{pointwise} convergence and \emph{uniform} convergence, with uniform convergence being the stricter mode.  
\end{itemize}

\section{Review: Power series}

Most of us are familiar with infinite series because they are introduced in introductory calculus courses.  There is a very intuitive way to think about Taylor series.  To set the stage, consider first the problem of fitting a polynomial to a set of data.  For example, we all know that a unique line can be fit through two points; similarly,  a quadratic through three points, and cubic through four points, etc.   This result is supported in a more general context of fitting $n+1$ points exactly by the following theorem (which we state, but do not prove)

\begin{theorem}
For any set of $k+1$ points $(x_i,y_i)$ where
\begin{itemize}[leftmargin=\dimexpr 26pt+.1in]
    \item $i$ is an index such that $i=1,\ldots k+1$
    \item $x_i\ne x_j$ for all $i \ne j$
\end{itemize}
then there is a \emph{unique} polynomial of degree \emph{at most} $k$ such that $f(x_i)=y_i$ for all $i=1,\ldots, k+1$.
\end{theorem}
\noindent  This theorem raises an interesting question.  Suppose we have a continuous function $f(x)$ (let's assume that all orders of derivatives also exist and are bounded, so that $f(x)$ is a $C^\infty$ function) defined on some interval $I=(a,b)$.  Now, assume we sample $n+1$ points in the interval $I$.  Apparently, we can get some kind of approximation to the function $f(x)$ by fitting an $n^{th}$-order polynomial to the $n+1$ sampled points.  This can be done, in practice, by a method such as Lagrange polynomial interpolation (which we will not discuss).

Now, by definition, all polynomials are of finite order. The order of the highest power can be as \emph{large} as you like, but it cannot be infinity; thus there are a finite number of terms.  The interesting question that might occur after thinking about this problem is whether or not one can consider a kind of ``infinite" polynomial that could fit continuous curves exactly.  Being somewhat loose, we might wonder of something like an ``infinite" polynomial would correspond to fitting an ``infinite" number of points on the interval $I$.  We would expect such a polynomial to provide an exact (in some sense) representation of the continuous function $f(x)$.  

While this is not necessarily a rigorous mathematical way of posing the question, it does offer a kind of analogy that provides some basis for working with what are called power series.  Recall, a polynomial of order $n$ is given by

\begin{align*}
    f(x)&= \sum_{n=0}^{n=N} a_n (x-b)^n
\end{align*}
where, here a \emph{shift} of the amount $b$ has been included in the polynomial.  A power series, then, is the generalization of this idea to the limit $n\rightarrow\infty$.
\begin{align*}
    f(x)&= \sum_{n=0}^{n=\infty} a_n (x-b)^n
\end{align*}
Note that any finite number of terms from this series is a polynomial; however, a power series extends the sum to infinity.  Despite their similarity, power series and polynomials do exhibit some very different kinds of behavior in appropriate limits.  Specifically, we note the following.

\begin{enumerate}
    \item First, all polynomials converge to a \emph{finite value} in every finite interval, $I=(a,b)$, where $a$ and $b$ are real numbers.
    \item Conversely, it is not true that every power series converges.  While finite polynomials have finite values on any finite interval, ``infinite" polynomials might not converge on any interval.
    \item It is not difficult to show that any polynomial tends toward plus or minus infinity as $x$ tends toward plus or minus infinity.  In other words, for a polynomial $f(x)$ of $n^{th}$ degree, we have
    \begin{equation*}
        \left| {\mathop {\lim }\limits_{x \to \infty } f(x)} \right| = \infty 
    \end{equation*}
    \item  A power series, in contrast to a polynomial of order $n$, can converge to \emph{any value}, including finite values, as $x$ tends toward infinity.  This is actually quite a strange result in some ways.  For example, suppose we have the following behavior for the power series representing $f(x)$
    \begin{equation*}
            {\mathop {\lim }\limits_{x \to \infty } f(x)} ={\mathop {\lim }\limits_{x \to \infty } \sum_{n=0}^{n=\infty} a_n (x-b)^n} = 0
    \end{equation*}
    Certainly, power series with this kind of behavior exist (for example, the power series for $f(x)=\exp{(-x^2)}$ behaves this way).  However, for \emph{any finite approximation to the power series} using only $N$ terms, we have (regardless of how large $N$ is)
     \begin{equation*}
           {\mathop {\lim }\limits_{x \to \infty } \left|\sum_{n=0}^{n=N} a_n (x-b)^n\right| } = \infty
    \end{equation*}
    This is very strange behavior because the power series tends to zero as $x\rightarrow\infty$, however \emph{every} finite approximation of that series diverges to $\pm \infty$.  It is frequently the case though that extensions of concepts from the finite to the infinite lead to non-intuitive behavior!
\end{enumerate}

\section{Review: Taylor Series}

The Taylor series should be a familiar concept from introductory calculus.  It was named after the English mathematician Brook Taylor who worked on problems of calculus in the early 1700's.  Although primarily known for his work on calculus, his interest ranged widely (as was typical of the times), even authoring a treatise under the unusual title \emph{On the Lawfulness of Eating Blood} which was discovered, upon his demise, to be among his unpublished papers.

The Taylor series is actually derivable directly from the definition of the power series.  Note, if we take the $n^{th}$ derivative of the power series we find the following

\begin{align*}
    & first~derivative & f'(x) &= a_1 +2 a_2 (x-b) + 3 a_3 (x-b)^2 +\ldots \\
    & second~derivative & f''(x) &= 2 a_2 + 6 a_3 (x-b) + \ldots \\
    & n^{th} derivative & f^{(n)}(x) &= n! a_{n} + (n+1)! a_{n+1} (x-b) + \ldots
\end{align*}
Starting from this result, if we now set $x=b$, we find that all of the terms except the first are identically zero.  The coefficients for the power series above are then defined by

\begin{align*}
 f'(b) &= a_1 \\
     \frac{1}{2!}f''(b) &=  a_2  \\
   \frac{1}{n!}f^{(n)}(b) &=  a_{n}
\end{align*}
Substituting this result into the power series above, gives us the \emph{Taylor series} around the point $x=b$

\begin{equation*}
     f(x)= \sum_{n=0}^{n=\infty} (x-b)^n \frac{f^{(n)}(b)}{n!} 
\end{equation*}
Sometimes the series written around $b=0$ is called the \emph{Maclaurin} series. Setting $b=0$ gives this series

\begin{equation*}
     f(x)= \sum_{n=0}^{n=\infty} x^n \frac{f^{(n)}(0) }{n!}
\end{equation*}
Any function that has a convergent Taylor series a set of points, $x\in X$ is called \emph{analytic} on $X$.  Many series (such as the exponential function) converge for all possible values of $x$, and are thus analytic everywhere.  Any finite truncation of a Taylor series is called a \emph{Taylor polynomial}.  Note that every finite truncation of a Taylor series is actually a polynomial!  Also true is that analytic functions converge nicely: as one increase the number of terms in the sum approximating the function $f$, the result gets uniformly closer to the actual value of $f$. An immediate consequence of these properties is that, for every function that is \emph{analytic} in some domain, $X$, there is a polynomial that can represent that function as closely as we like.  Essentially, this follows from the definition of an analytic function, and the Taylor polynomial.

\begin{svgraybox}
\begin{example}[Taylor/Maclaurin series examples]
The function $f(x)=\exp(x)$ is an interesting example to consider for a Taylor series expansion, in part because it converges exactly everywhere in the domain $x\in(-\infty,\infty)$.  To start the analysis, note that we first need an infinite number of derivatives of $f(x)$.  While this may seem like a daunting task, often derivatives of functions exhibit a pattern of behavior that can be exploited (using the principle of induction).  Consider the following
\begin{align*}
    f(x) &= \exp(x) \\
    f'(x) &= \exp(x)\\
    f''(x)&= \exp(x) \\
    \ldots
\end{align*}
Here, we do not need to work too hard to see the pattern.  For \emph{any derivative of order $n$}, the derivative is given by the equation $f^{(n)}(x) = \exp(x)$.  Now, suppose that we expand the series around the point $x=0$ (i.e., we choose the shift parameter, $b$, to be zero, which is also the definition of a Maclaurin series).  Then, each derivative is equal to unity $f^{(n)}(0) = 1$.  The resulting Taylor series is 
\begin{align*}
    \exp(x)&= \sum_{n=0}^{n=\infty} x^n \frac{1}{n!}\\
      &= 1+x+\frac{x^2}{2!} + \frac{x^3}{3!}+\frac{x^4}{4!}+\frac{x^5}{5!}+\frac{x^6}{6!}+\ldots
\end{align*}

As a second case, suppose we would like to find the Taylor series for $\sin(x)$.  
\begin{align*}
    f(x) &= \sin(x) \\
    f'(x) &= (-1)^2 \cos(x)\\
    f''(x)&= (-1)^3 \sin(x) \\
    f'''(x)& =(-1)^4 \cos(x)\\
    \ldots
\end{align*}
In principle, we can define the Taylor series to be expanded around any real number. And there are good reasons to do that in some cases; for example, if you know you need a particular expansion to be accurate in a specific interval, $I=(a_1,a_2)$, it would be good to choose $a_1 < b < a_2$ for improved accuracy in the expansion. In this example, no particular point was stated, so we adopt $x=0$ for convenience.  Under those conditions, we find
\begin{align*}
    f(0) &= 0 \\
    f'(0) &= 1\\
    f''(0)&= 0 \\
    f'''(0)& =-1 \\
    f^{(4)}(0)& =0\\
    f^{(5)}(0)& =1\\
    \ldots
\end{align*}
In other words, we want only the \emph{odd} terms, and those terms alternate sign.  A little thought, will indicate that the function $\textrm{odd}(n)=2n+1,~n=0,1,2,\ldots$ counts by odd numbers only.  To switch signs we multiply by powers of $-1$; we can define a function that alternates sign by $\textrm{alt}(n)=(-1)^n$.  Putting this together, we find the following result valid for any value of $n$
\begin{equation*}
    f^{(2n+1)}(0) = \textrm{alt}(n)=(-1)^{n} 
\end{equation*}
Note the following examples
\begin{align*}
n&=0 & 2n+1&=1 & f^{(1)}(0) &= (-1)^0 = +1 \\
n&=1 & 2n+1&=3 & f^{(3)}(0) &= (-1)^1 = -1 \\
n&=2 &2n+1&=5& f^{(5)}(0) &= (-1)^2 =  +1
\end{align*}
Recalling the definition of the Maclaurin series
\begin{equation*}
     f(x)= \sum_{n=0}^{n=\infty} x^n \frac{f^{(n)}(0) }{n!}
\end{equation*}
and noting again that only the odd terms are non-zero, the series representation for $\sin(x)$ is
\vspace{-4mm}
\begin{align*}
     f(x)&= \sum_{n=0}^{n=\infty} x^{2n+1} \frac{f^{(2n+1)}(0) }{(2n+1)!} \\
     & = \sum_{n=0}^{n=\infty} x^{2n+1} \frac{(-1)^n }{(2n+1)!}\\
     &= x-\frac{x^3}{3!}+\frac{x^5}{5!}-\frac{x^7}{7!}+\ldots
\end{align*}
Plots of finite Taylor polynomial (Taylor series truncated at order $n$) are given in the figure below.

{
\centering\fbox{\includegraphics[scale=.6]{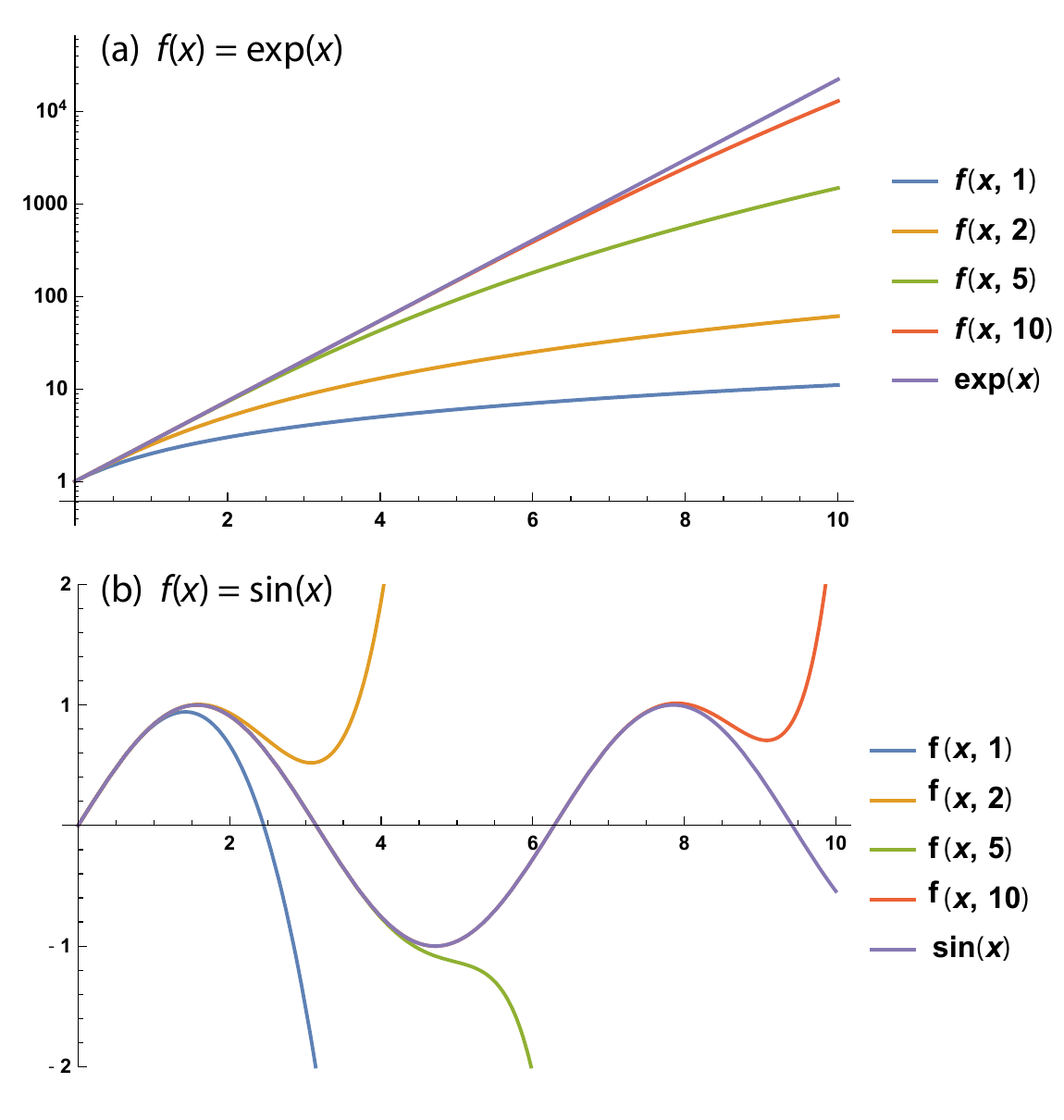}}
\vspace{-1mm}
\captionof{figure}{(a) Exponential, and (b) Sine functions plotted for $1,~2,~5,$ and $10$ terms in the approximating Taylor polynomial.}
\label{functions2}  
}

\end{example}
\end{svgraybox}

Taylor series for specific functions are widely available, and are also not unreasonably difficult to compute directly from the definition.   Note, however, that all Taylor series have a \emph{radius of converge}.  In other words, there is a domain that a Taylor series will converge (i.e., tend toward a definite number as the number of terms in the sum increase, and that definite number is equal to the value of $f$ at that point).  It is generally easy to determine if a series converges or not, but much more difficult to determine its radius of convergence.  We will not discuss convergence properties of series in this chapter, but we will be careful to list their radius of convergence if it is known.  A few well-known Taylor series are given by

\begin{align*}
 e^{x}& =\sum _{n=0}^{\infty }{\frac {x^{n}}{n!}}&&=1+x+{\frac {x^{2}}{2!}}+{\frac {x^{3}}{3!}}+\cdots  && (converges~for~all~x) \\
\log(1+x)&=\sum _{n=1}^{\infty }(-1)^{n+1}\frac {x^{n}}{n}&&=x-\frac {x^{2}}{2}+\frac {x^{3}}{3}-\cdots  &&(converges~for~| x | < 1) \\
\frac {1}{1-x}&=\sum _{n=0}^{\infty }x^{n} &&=x+x^2+x^3+\ldots && (converges~for~| x | < 1) \\
{\frac {1}{(1-x)^{2}}}&=\sum _{n=1}^{\infty }nx^{n-1} &&=1+2x+3x^2+4x^3+\ldots && (converges~for~| x | < 1)\\
 \sin x&=\sum _{n=0}^{\infty }{\frac {(-1)^{n}}{(2n+1)!}}x^{2n+1}&&=x-{\frac {x^{3}}{3!}}+{\frac {x^{5}}{5!}}-\cdots && (converges~for~all~x)\\
\cos x&=\sum _{n=0}^{\infty }{\frac {(-1)^{n}}{(2n)!}}x^{2n}&&=1-{\frac {x^{2}}{2!}}+{\frac {x^{4}}{4!}}-\cdots && (converges~for~all~x)\\
\end{align*}

One of the most useful properties of series like this is for computing the value of \emph{transcendental} functions (Recall from Chapter 1, a transcendental function is (in short) an analytic that cannot be expressed by a finite polynomial).   Thus, the Taylor series give one of the few methods that are available to compute approximations to transcendental functions.  If you think about it, how else would you compute, for example, the value of $\sin \pi/12$?  You could estimate it graphically (by, say, drawing a giant unit circle) like the ancient Greek's did, but a method to do it using real numbers requires a series solution (or some other algorithm).  We get accustomed to hitting the ``sin" button on our calculators, but computing the actual values of the $\sin$ function is actually very difficult!  In the distant past (before calculators- gasp!), series solutions were used to compute the values for useful transcendental functions (like $\sin$, $\cos$, $ln$, etc.), and they were published in very large books.  Engineers and scientist would have to \emph{look up these numbers in immense tables} whenever they wanted to know a particular value for $\sin$ or $\cos$.  Let's hear it for calculators!

In summary, this discussion about series is meant to highlight a few key concepts.
\begin{enumerate}
    \item Power series can be interpreted as an extension of the concept of polynomials to infinite degree.
    \item Taylor and Maclauren series are a special cases of power series.
    \item Power series may converge only on some finite interval; it is also possible that they do not converge at all.  
    \item Because the sum is an infinite one, the limiting behavior of a series (e.g., as $|x|\rightarrow \infty$) might not be well approximated by a finite truncation of the series.
    \item Infinite sums frequently thwart intuition.  Infinite sums can create behavior that would otherwise not be expected!  
\end{enumerate}

\section{Trigonometric Series}

\begin{figure*}[t]
\centering
\includegraphics[scale=0.5]{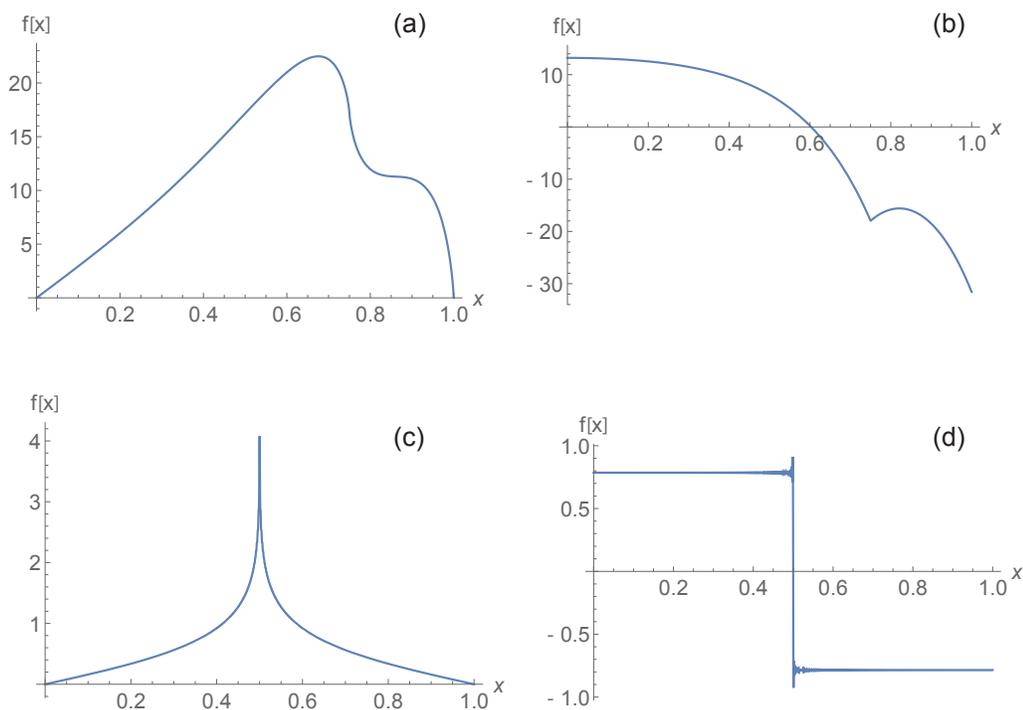}
\caption{The functions given by (a)-(d) above. }
\label{fseries}       
\end{figure*}

Power series are not the only kind of series that provides useful results.  In fact, expansions in \emph{trigonometric series} are probably more widely used than any other kind of series expansion.
A trigonometric series is one that contains a trigonometric function that varies with the index of the sum.  A few examples illustrate some examples of trigonometric series on the interval (0,1).  In each of these examples, $n$ is an integer greater than or equal to zero.

\begin{align*}
    &(a)& B_n &= -50\cos\left(\frac{n \pi}{r}+\frac{\pi}{4} \right)&& f(x) = \sum_{n=1}^\infty B_n \sin\left( n \pi x\right) 
    \end{align*}
    \begin{align*}
    &(b)& A_n &= -50\cos\left(\frac{n \pi}{r}+\frac{\pi}{4} \right))&& f(x) = \sum_{n=1}^\infty A_n \cos\left( n \pi x\right) \\
    &(c)&B_n&=\frac{\sin \left(\frac{\pi  n}{2}\right)}{n}&& f(x) = \sum_{n=1}^\infty B_n \sin\left( n \pi x\right) \\
    &(c)&A_n&=\frac{\sin \left(\frac{\pi  n}{2}\right)}{n}&& f(x) = \sum_{n=1}^\infty A_n \cos\left( n \pi x\right) \\
\end{align*}
The plots of these series appear in Fig.~\ref{fseries}
where here we have used $A_n$ to indicate the coefficients of cosine series, and $B_n$ to indicate the coefficient for sine series.

As you can see in these results, trigonometric series can allow the representation of some very peculiar functions.  Some of the interesting features about the particular functions plotted above are as follows.
\begin{itemize}
    \item Although the two functions used to create the new functions are $C^\infty$ (i.e., infinitely smooth analytic functions), the series seem to give functions appear to be, in some cases, non-differentiable at every point.
    
    \item Although the two functions used to create the new functions are $C^\infty$ (i.e., infinitely smooth analytic functions), the series seem to give functions appear to be, in some cases, discontinuous.
    
    \item The sine and cosine series give very different series solutions even when $A_n=B_n$.
\end{itemize}
Although all of these examples are interesting, it is not clear that any of them correspond to any classical function that we recognize- either polynomial or transcendental.  So, although we certainly can generate many interesting trigonometric series, the questions arises: ``if I am given a specific, known (polynomial or transcendental function) on a finite interval, can I determine a trigonometric series for that?"  This is a question that occurred to Joseph Fourier, a French engineer, in the early 1800's.  The answer to this question, remarkably, is a resounding yes for almost any function that one can imagine.  In fact, modern mathematics was dramatically shaped by the quest for the answer to this question.  Not only did it cause mathematicians to re-think the notion of what a function is, but it caused them to refine the mathematical methods that have led to modern mathematical analysis.

\section{Fourier Series} \indexme{Fourier series}

For power series, one can think of the expansions as being in an infinite set of polynomials (although, technically, there is no infinite-order polynomial).  For trigonometric series, the series expansions are, not surprisingly, in trigonometric functions.  

For the development of the Taylor series, we were able to develop a scheme in which we (1) proposed an expansion in an infinite series with an infinite number of (unknown) coefficients, $a_n$, and (2) determined the infinite number of coefficients by finding repeating patterns in the derivatives needed, so that any derivative of order $n$ could be explicitly computed if $n$ were specified.  

For Fourier series, we will do something similar.  However, instead of the unknown parameters being a function of an infinite number of derivatives, we will be able to express the unknown parameters as an infinite number of \emph{integrals}.  While this sounds on the surface to be dire, like the case for Taylor series, we will find that we can determine repeating patterns in the integrals so that we can derive closed-form expressions for the integrals in terms of the series index $n$.

A critical component for the development of the Fourier series is the extension of the concept of \emph{orthogonality} to continuous functions.  While on the surface this may not immediately make intuitive sense, it can be made intuitive by analogy with the familiar concept of orthogonality for finite vectors.

\subsection{Orthogonality Revisited}\label{orthog}

\begin{figure*}[t]
\centering
\includegraphics[scale=0.5]{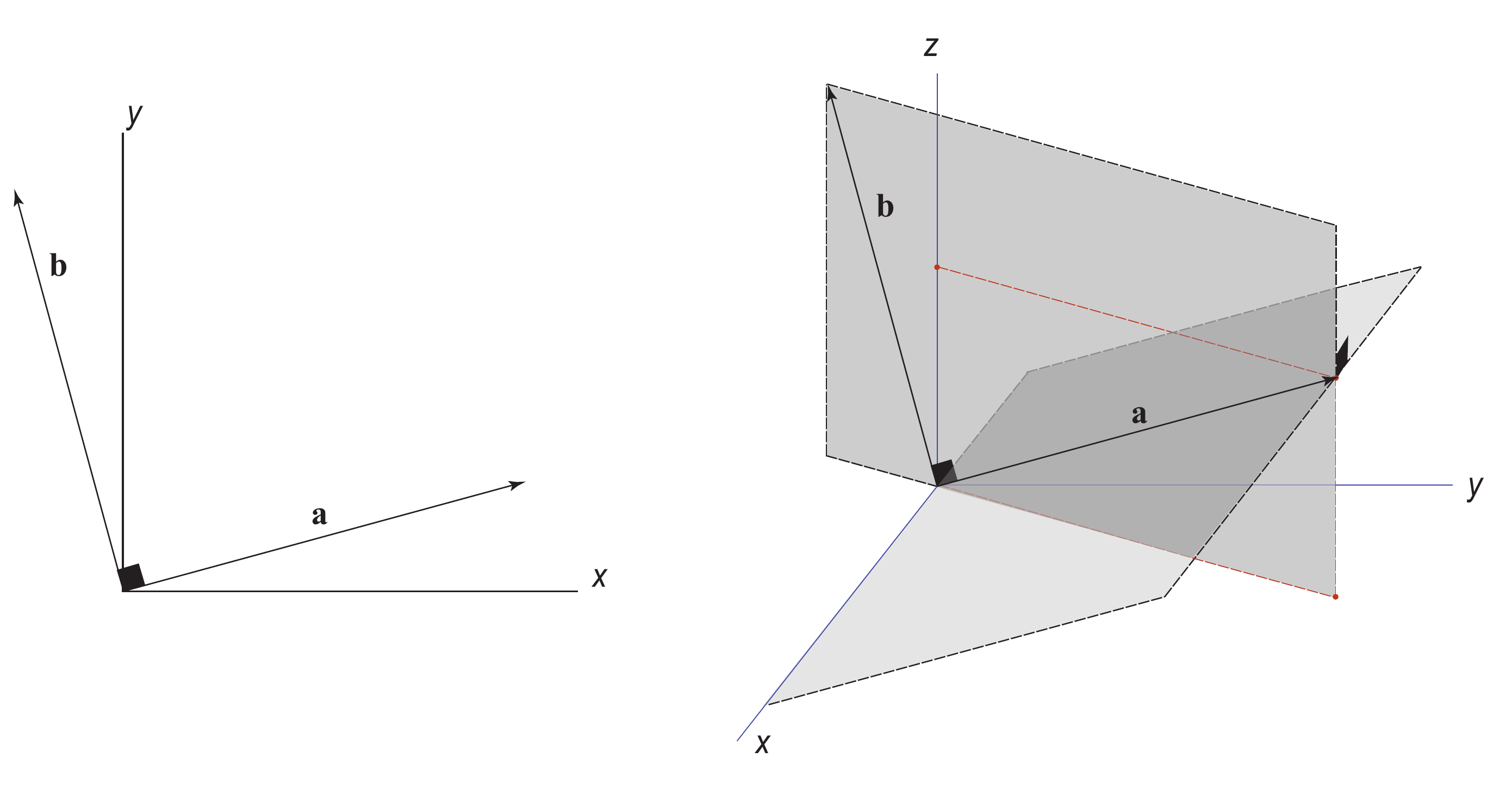}
\caption{Perpendicular vectors ${\bf a}$, and ${\bf b}$ in 2- and 3-dimensions. }
\label{perp}       
\end{figure*}

At some point in our mathematical education, most of us have encountered the concept of \emph{orthogonality} for two vectors.  Perpendicular vectors in both 2- and 3-dimensions are illustrated in Fig.~\ref{perp}.  There are two related concepts that help us define and describe perpendicular vectors in 2- and 3-dimensions.

\begin{enumerate}
    \item First, any two vectors ${\bf a}$ and ${\bf b}$ are said to be perpendicular if the dot product between the two vectors is zero.  That is, two vectors are perpendicular if ${\bf a}\cdot{\bf b} = 0$.  Recalling that the dot product of two vectors is equal to their magnitudes times the cosine of the angle between them, we have
    \begin{equation*}
        {\bf a}\cdot{\bf b}=\|{\bf a}\| \|{\bf b}\|\cos{\theta}
    \end{equation*}
    where $\theta$ is the angle between the two vectors in the plane that contains them both.  Note that this is zero, exactly when $\theta = \pi/2 = 90^\circ$.
    
    \item Second, there are a set of mutually perpendicular \emph{basis vectors} that can be used to define any arbitrary vector as a linear \emph{weighted sum} of the basis vectors.  These vectors are often given the symbols ${\bf i}$, ${\bf i}$, and ${\bf i}$ corresponding to unit vectors in the $x-$, $y-$, and $z-$ directions, respectively.  Thus, in the Cartesian coordinate system, we have 
    \begin{align*}
        {\bf i} = (1,0,0)\qquad {\bf j} = (0,1,0)\qquad {\bf k} = (0,0,1)
    \end{align*}
    Using these basis vectors, any vector ${\bf a}=(a_1,a_2,a_3)$ can be represented by the following weighted sum of the basis vectors
    ${\bf a}=a_1{\bf i}+ a_2{\bf j}+ a_3{\bf k}$.
\end{enumerate}
%

In summary then, every vector can be decomposed into a linear weighted sum of its basis vectors, and each of these basis vectors is, by definition, perpendicular to the others.  These ideas are relatively familiar and easy to grasp in 2- and 3-dimensions.  Of course, these concepts can be extended to any finite number of dimensions.  We cannot visualize such extensions, but mathematically the concepts remain valid.  We will discuss in what sense two \emph{functions2} can be perpendicular by making an analogy of functions being approximated by finitely-long vectors.

\begin{figure*}[t]
\centering
\includegraphics[scale=0.65]{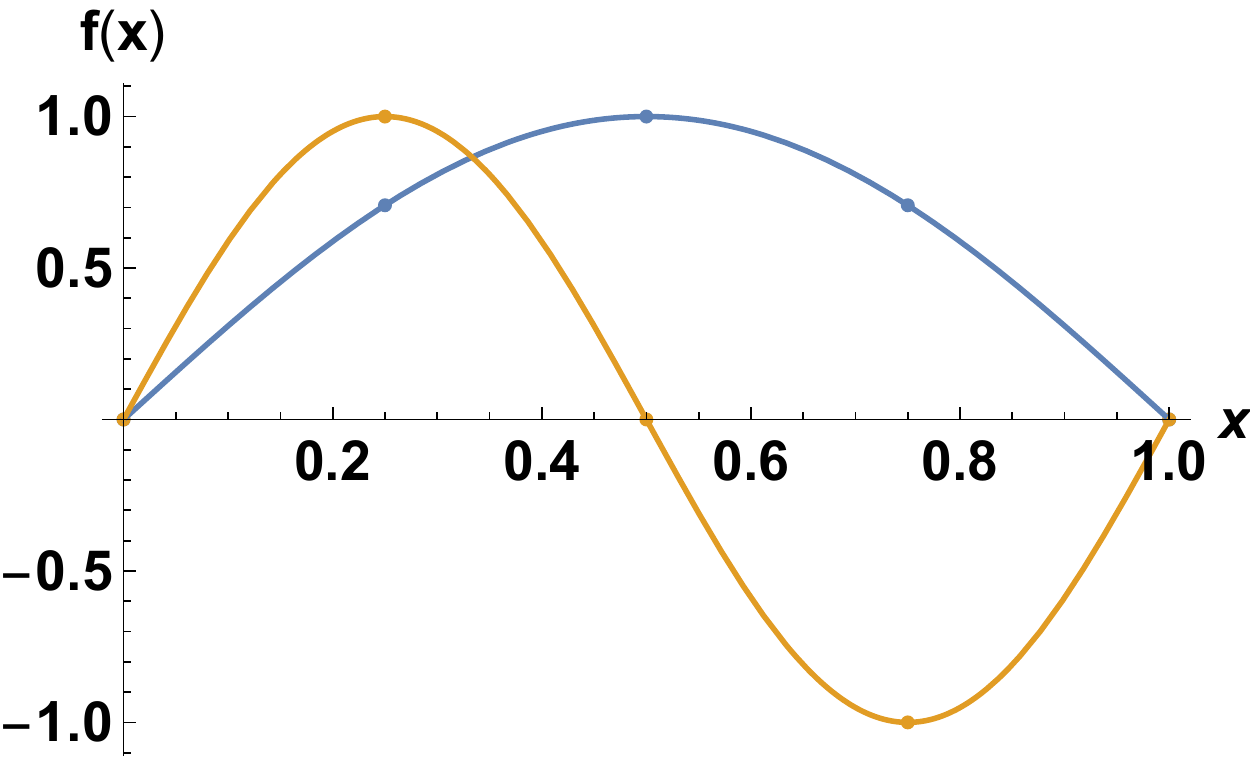}
\caption{The functions $f(x)=\sin(\pi x)$ and $g(x)=\sin(2\pi x)$, approximated by four discreet sampling intervals (defined by 5 points, ${\bf x}=(0,1/4,1/2,3/4,1)$).  This representation results in  5-dimensional vectors representing the values of the sine functions.  The continuous representation of these two functions is shown for comparison. }
\label{sin_4points}       
\end{figure*}

Now, consider two \emph{finite} approximations to the functions $f(x)=\sin(\pi x)$ and $g(x)=\sin(2 Pi x)$ on the interval $I=[0,1]$.  By finite approximation, in this case we mean 
\begin{enumerate}
    \item Segmenting the interval into $N$ pieces represented by the vector ${\bf x}$, where the spacing is given by $\Delta x_i=1/N$, and the components of the vector are given by the recursive relationship $x_{i+1}=x_{i}+\Delta x_i,~i=0,1,2,\ldots, N$. 
    \item Determining the values of $\sin(\pi x_i)$ and $\sin(2\pi x_i)$, for $i=0,1,2,\ldots,N$.
\end{enumerate}
As a simple (and crude) example, consider the case of $N=4$.  For this case, the values of ${\bf x}$ are given by ${\bf x}=(0, 1/4, 1/2, 3/4, 1)$.  For notation, define the vector of values for the $\sin(\pi x_i)$ function by $_{1\pi}{\bf S}=\sin(\pi x_i),~ i=0,1,2,3,\ldots,N$. So, for this example, the result is $_{1\cdot\pi}{\bf S}=(0,\sin(\pi/4),1,\sin(\pi/4),0)$.  We will use a similar definition for the $\sin(2\pi x_i)$, so that $_{2\pi}{\bf S}=(0,1,0,-1,0)$.  The two sine functions, and the points for the discrete approximation of them, are illustrated in Fig.~\ref{sin_4points}.

\begin{figure*}[t]
\centering
\includegraphics[scale=0.65]{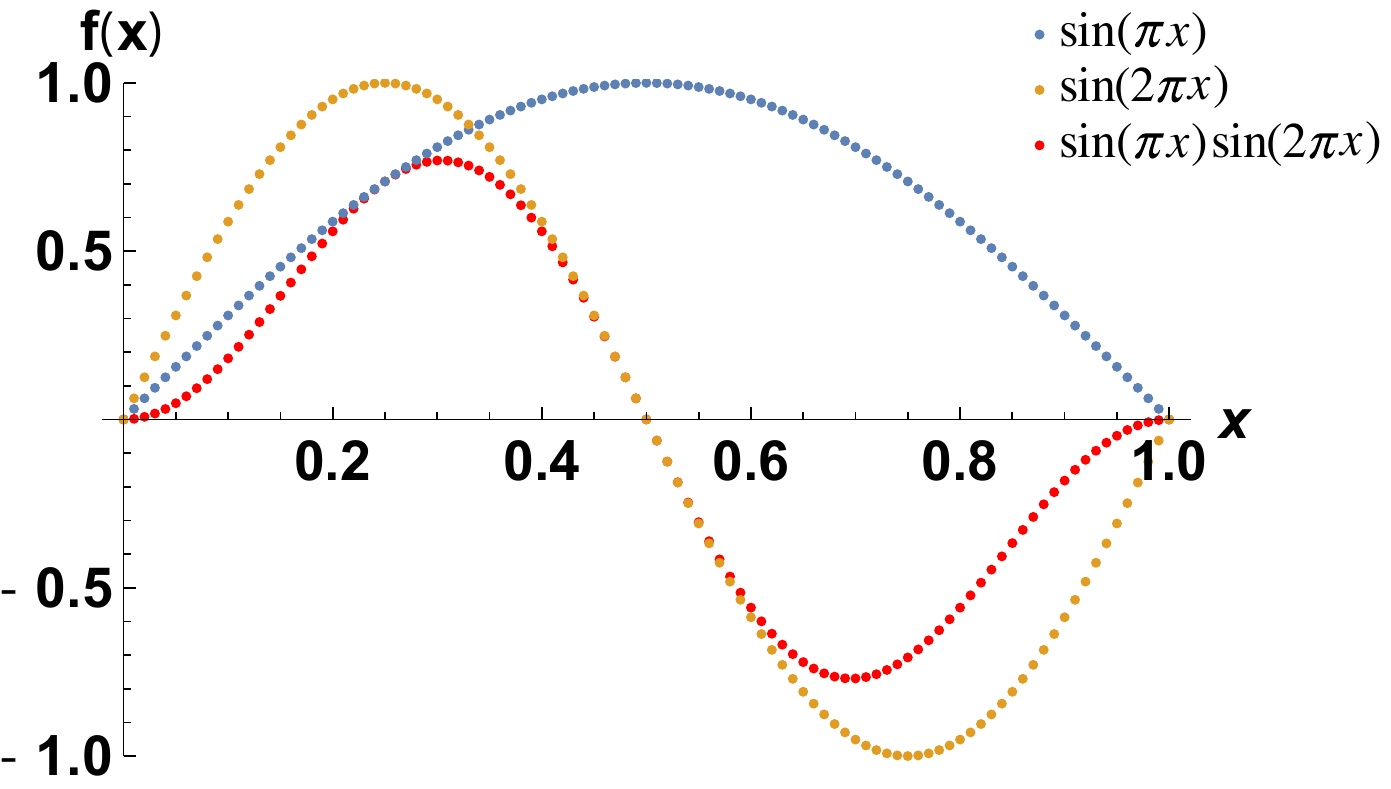}
\caption{The functions $f(x)=\sin(\pi x)$ and $g(x)=\sin(2\pi x)$, approximated by four discreet sampling intervals (defined by 101 points).  This representation results in  101-dimensional vectors representing the values of the sine functions.}
\label{sin_100points}       
\end{figure*}

%
Although this is a rather crude approximation to the two sine functions, note that for our carefully chosen sampling points, we have the following result

\begin{align*}
    _{1\pi}{\bf S}\cdot_{2\pi}{\bf S} &=(0,\sin(\pi/4),1,\sin(\pi/4),0)\cdot(0,1,0,-1,0)\nonumber\\
    &=0\cdot0+\sin(\pi/4)\cdot1+\sin(\pi/4)\cdot(-1)+0\cdot0\nonumber\\
    &=0
\end{align*}
So, even though this is an approximate representation of the two sine functions, we find that the finite vectors representing the functions themselves are orthogonal. This is a bit more compelling when we use a larger number of points, say $N=100$.  For that case, the approximate discrete functions (each with 101 points) are illustrated in Fig.~\ref{sin_100points}.

%
For this discrete representation, we can see that the number of points used is high enough that we really are capturing the features of the two sine functions.  We have a large list of numbers for ${\bf x}$, $_{1\pi}{\bf S}$, and $_{2\pi}{\bf S}$, as follows

\begin{align*}
    {\bf x}&=(0,1/100, 2/100, 3/100, \ldots, 1) \\
    _{1\pi}{\bf S}&=(0,\sin(\pi/100), \sin(2 \pi/100),\ldots, 0)\\
    _{2\pi}{\bf S}&=(0,\sin(2\pi/100), \sin(4 \pi/100),\ldots, 0)
\end{align*}
While we can verify that $_{1\pi}{\bf S}\cdot _{2\pi}{\bf S}$ by direct computation, a look at the graph in Fig.~\ref{sin_100points} makes this somewhat unnecessary.  The red points plotted in this graph represent the product of the two sine functions. A little though will indicate that the symmetry of the problem guarantees that the sum of the red points (which represents the dot product $_{1\pi}{\bf S}\cdot _{2\pi}{\bf S}$) will sum to zero.  

While this is an approximate method to understand how two continuous functions can be though of as being orthogonal, it also is to some extent rigorous.  For example, consider the quantity

\begin{align*}
    \sum_{i=0}^{i=N+1} {_{1\pi}{\bf S}} \Delta x_i&= \sum_{i=0}^{i=N+1} \sin(\pi x_i) \Delta x_i\\
    \intertext{In the appropriate limit, this gives}
    \mathop {\lim }\limits_{N \to \infty } \left(\sum_{i=0}^{i=N+1} \sin(\pi x_i) \Delta x_i\right)&= \int_{0}^{1} \sin(\pi x) \,dx
\end{align*}
And, in a similar fashion, we have

\begin{align*}
    \sum_{i=0}^{i=N+1} {_{1\pi}{\bf S}}\cdot {_{2\pi}{\bf S}} \Delta x_i&= \sum_{i=0}^{i=N+1} \sin(\pi x_i) \sin(2 \pi x) \Delta x_i\\
    \intertext{In the appropriate limit, this gives}
    \mathop {\lim }\limits_{N \to \infty } \left(\sum_{i=0}^{i=N+1} \sin(\pi x_i) \sin(2 \pi x) \Delta x_i\right)&= \int_{0}^{1} \sin(\pi x_i) \sin(2 \pi x) \,dx = 0
\end{align*}
In particular, it is easy to validate that this last integral is identically zero.  Thus, the concept of two continuous functions being \emph{orthogonal} to one another really does have a direct, demonstrable connection to the case of finite-sized vectors being orthogonal.

\subsection{Fourier Sine Series}\indexme{Fourier series! sine series}

Fourier series are a special kind of trigonometric series that have somewhat astounding properties.  Fourier series use only the sine and cosine functions to expand a function $f(x)$ as a series.  In starting the investigation of Fourier series, we will consider first only the interval $I=\{x: x\in [0,1]\}$; this interval is sometimes called the \emph{unit interval}.  We can expand the definition of the Fourier series to other intervals once we understand how they work on the unit interval.  

To start, consider the most basic question that we can ask.  If we have a simple, analytic function on an interval $x\in [0,1]$, can we determine the Fourier series for it?  To be concrete, let's suppose we have a specific function, say $f(x)=e^{-x}$, and we decide we would like to find a sine series for that function.  The question is, can we find a series of the form

\begin{align}
    f(x)&=\sum_{n=0}^{\infty} B_n \sin(n \pi x)\\
    \intertext{Or, being specific to the example given, $f(x)=e^{-x}$}
    e^{-x}&=\sum_{n=0}^{\infty} B_n \sin(n \pi x)
    \label{expo}
\end{align}
where for both expressions, $n$ is an integer, and $B_n$ is an infinite sequence of constants (i.e., the sequence $B_n=(B_1, B_2, B_3, B_4, \ldots ).$  There are two primary questions that we need to address about such a proposed series.  These are
\begin{enumerate}
    \item Is there a method (an algorithm or constructive proof) that allows us to determine the infinite sequence of constants $B_n$?
    
    \item If we can find the constants $B_n$, can we show that the resulting series converges to the function $f(x)$?
\end{enumerate}

\begin{figure*}[t]
\centering
\includegraphics[scale=0.85]{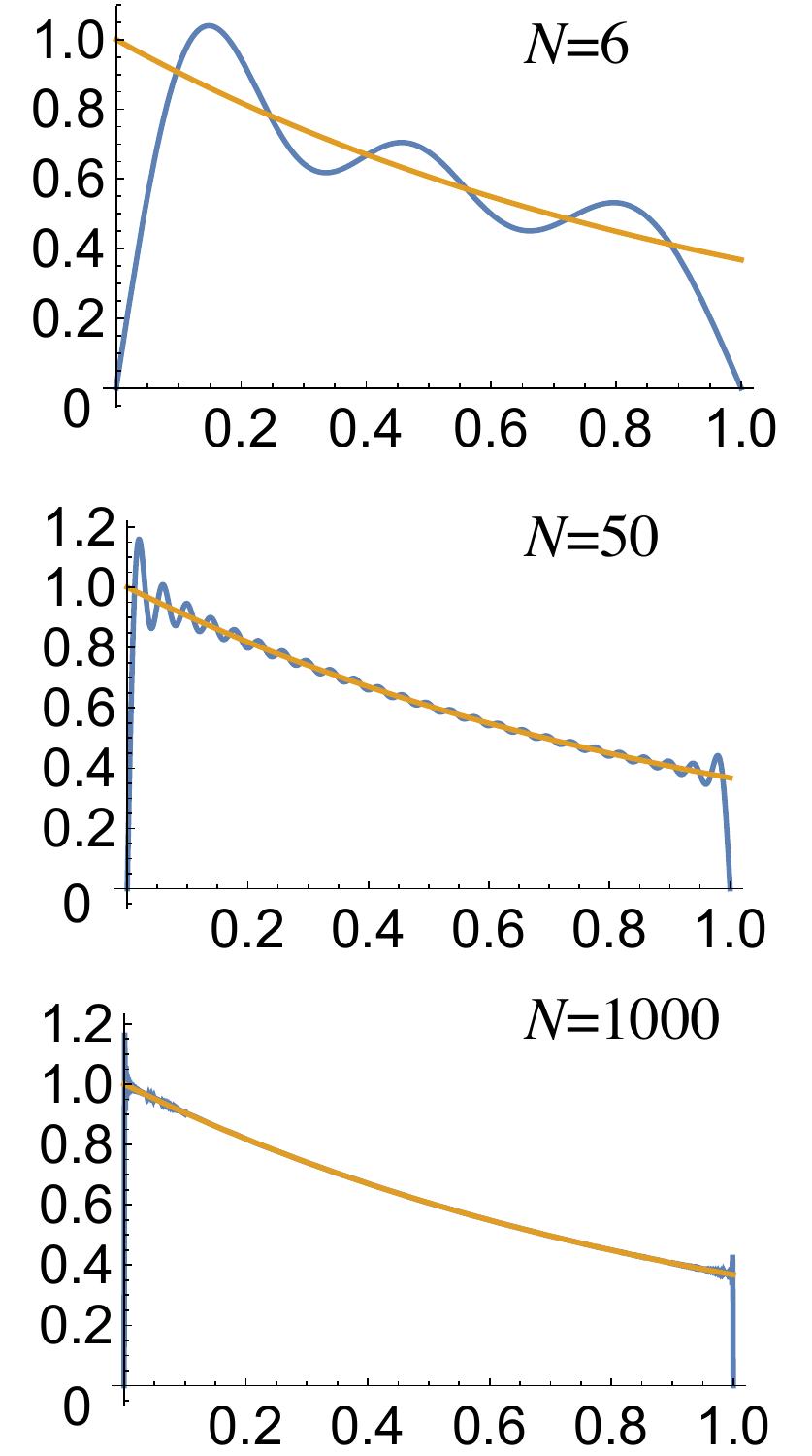}
\caption{The partial sum of the series, for the total number of terms, $N$, equal to $N=6,50, \text{ and } 1000$. }
\label{graph2}       
\end{figure*}

Before we proceed, we need to make a few notes.  First, we note that because $\sin(0 \pi x) =0$, technically we do not need to start the series at $x=0$; we could start at $x=1$.  By convention, sine series are usually written as

\begin{align*}
    e^{-x}=\sum_{n=1}^{\infty} B_n \sin(n \pi x)
\end{align*}
It is not wrong to start the series at $n=0$ though!  The function $\sin(0 \pi x)$ is one of the basis functions for the sine series; it just does not add anything if it is maintained.  Secondly, note the following identity (which we discussed in the context of orthogonality in the previous section)

\begin{align}
    \int\limits_{x=0}^{x=1} \sin(n \pi x) \sin(m \pi x) dx &=
    \left\{ {\begin{array}{*{20}{c}}
  {0~if~m\ne n} \\ 
  {\frac{1}{2}~if~m=n} 
\end{array}} 
\right.
\label{orthogonal}
\end{align}
In these expressions, both $n$ and $m$ are integers, not necessarily equal.  Because of the direct analogy with the dot product in the case of finite-length vectors, we call two functions \emph{orthogonal} when they meet a condition like Eq.~\eqref{orthogonal}.   With this information in hand, we now note the following reasonably amazing result.  If we multiply both sides of Eq.~\eqref{expo} by $\sin(m \pi x)$ and integrate, we find

\begin{align*}
    \int\limits_{x=0}^{x=1} e^{-x} \sin(m \pi x) dx =\sum_{n=1}^{\infty} B_n  \int\limits_{x=0}^{x=1} \sin(n \pi x) \sin(m \pi x) dx
\end{align*}
where we have taken the integral inside the sum (and this is always allowable for our purposes; the same is not necessarily true for differentiation!)
The reason that this is interesting, is because of what it does to the right-hand side of the equation.  When we started, we had an infinite number of values of $B_n$ to contend with.  But, because the integral on the right-hand side is only non-zero when $m=n$, then we must have the following

\begin{align*}
    \int\limits_{x=0}^{x=1} e^{-x} \sin(n \pi x) dx = B_n  \int\limits_{x=0}^{x=1} \sin^2(n \pi x)  dx
\end{align*}
In other words, all of the terms in the series, after integration, are zero, except for the \emph{one term} where $n$ and $m$ are equal.  For that one term, the integral becomes the integral of $\sin^2(n \pi x)$ over the interval $x\in[0,1]$.  By the identities above, this is exactly $\tfrac{1}{2}$.  Thus, we end up with the result (rearranging a little) that allows us to compute $B_n$:
\begin{align*}
   B_n =  2\int\limits_{x=0}^{x=1} e^{-x} \sin(n \pi x) dx 
\end{align*}
This integral can be done by parts (and then using the fact that it generates a repeating function that can be collected on one side), from a table of integrals, or from software like Mathematica.  The result is 

\begin{equation*}
    B_n = \frac{\pi  n-e \pi  n \cos (\pi  n)}{\pi ^2 n^2+1}
\end{equation*}
And, while this result looks a bit strange, we have nonetheless apparently computed the value of $B_n$ for all possible values of $n$.  In other words, we have exactly the coefficients we need to compute the result for our series.  Our series apparently takes the form

\begin{align*}
    e^{-x}=\sum_{n=1}^{\infty} \frac{\pi  n-e \pi  n \cos (\pi  n)}{\pi ^2 n^2+1} \sin(n \pi x)
\end{align*}
Can this possibly be true?  Is there some way that a bunch of sine curves, properly weighted, and added together can somehow form an exponential?  Well, let's find out.  Below, I have computed the first 10 values of $B_n$ 
\begin{align*}
    B_n=&\{0.790704, 0.196239, 0.287042, 0.099972, 0.173461, 0.0668818, \\
&0.124146, 0.050223, 0.0966368, 0.0402013\}
\end{align*}
These values represent the \emph{amplitude} or the height of each of the sine functions up to $n=6$.  To see the functions themselves, we need only compute the values of $B_n \sin(n \pi x),~n=1,2,3,4,5,6$.  These are plotted on Fig.~\ref{graph2}.  As a matter of terminology, note that the functions $\sin(n \pi x)$ are called \emph{basis} functions.\indexme{Fourier series!basis functions}\indexme{basis functions}  For each such basis function, the integer $n$ controls the \emph{frequency} of the basis function, and the associated coefficient $B_n$ controls the \emph{amplitude} of the basis function.  

\begin{figure*}[t]
\centering
\includegraphics[scale=0.7]{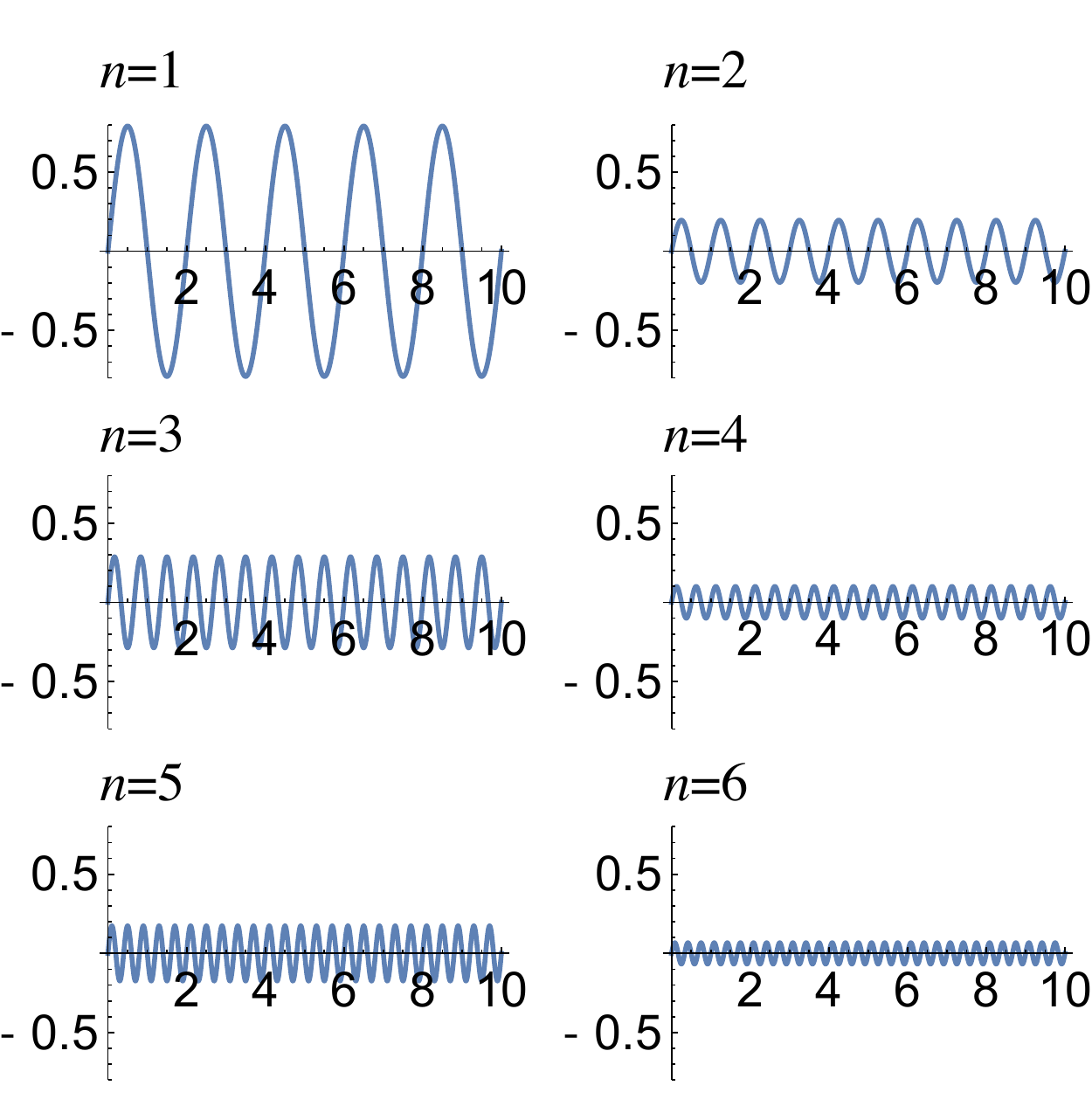}
\caption{The functions $B_n \sin (n \pi x)$ for $n=1,2,3,4,5,6$. }
\label{graph}       
\end{figure*}
%
It is clear from the plot in Fig.~\ref{graph} that each of the terms in the series is a sine function, but each with a different amplitude, and each with a frequency that increases with $n$.  The most remarkable thing, however, is what happens when we add these components together.  Adding the first 6 terms (as we have above) gives us an approximation to the function $f(x)=e^{-x}$ that is not necessarily good (although, it is not necessarily bad either!).  If we compute the first 50 terms, things begin to look much nicer.  Finally, with the first $100$ terms, the function and the series are almost indistinguishable at most points.  

The set of functions $~E=\{0,~\sin(\pi x),~\sin(2 \pi x), ~\sin(3 \pi x), \ldots\} $ can be thought of as a set of \emph{orthogonal basis functions} from which new functions can be built.  There is an analogy with \emph{orthogonal basis vectors} here that was introduced in \S\ref{orthog}.

For infinite-dimensional vectors (functions), the  equivalent of the dot-product is an integral that is called the \emph{inner product} or sometimes the $L^2$ (pronounced "ell-two") inner product.  Although we will not make extensive use of it, there is even a specific symbol used in mathematics to indicate the $L^2$ inner product.  On the interval $x \in [0,1]$ we would have

\begin{equation*}
    \langle \sin(n \pi x),\sin(m \pi x)\rangle = \int\limits_{x=0}^{x=1}  \sin(n \pi x)\sin(m \pi x)\, dx
\end{equation*}

To summarize the main results for the sine series, we have the following.  For any smooth function, $f(x)$, on $x \in [0,1]$, the Fourier sine series is given by

\begin{align}
    f(x)&=\sum_{n=0}^{\infty} B_n \sin(n \pi x) \\
   B_n& =  2\int\limits_{x=0}^{x=1} f(x) \sin(n \pi x) dx, ~~ n=0,1,2,3,\ldots
\end{align}
or, equivalently, because $\sin(0 \pi x)=0\, \Rightarrow B_0=0$
\begin{align}
    f(x)&=\sum_{n=1}^{\infty} B_n \sin(n \pi x) \\
   B_n& =  2\int\limits_{x=0}^{x=1} f(x) \sin(n \pi x) dx, ~~ n=1,2,3,\ldots
\end{align}
We can think of the set of functions $E_{sin}=\{0,\sin(\pi x), \sin(2 \pi x), \sin(3\pi x),\ldots\}$ as being the basis functions from which any function on the unit interval $I=[0,1]$ can be reconstructed. The reconstruction requires an infinite sum to be made, and the sum is a weighted one, where the weights (amplitudes) are given by the values of $B_n$.  Note that there is an almost exact correspondence here to the case of reconstructing an arbitrary finite vector ${\bf a}$ by computing the sum of a weighted set of basis functions: ${\bf a}=a_1{\bf i}+a_2{\bf j}+a_3{\bf k}$.\\


\begin{svgraybox}
\begin{example}[Fourier Sine Series Example: The Heaviside Function]\label{heavisideex}

It can be quite useful to define a function that is equal to zero to the left of some point (say, $x_0$), and equal to 1 to the right of $x_0$.  Such a function is known as a \emph{step} function or \emph{Heaviside} function.  This function will be discussed more in the chapter titled ``The Step and Delta Functions" (Chp.~\ref{deltachap}).  For now, we define the function as follows.
\begin{equation}
H(x) = \begin{cases}
0 & \textrm{~for~} x<0 \\
1 & \textrm{~for~} x\ge 0 \\
\end{cases}
\end{equation}
Note that the step function can be shifted in the standard way; that is, the step function that is zero to the left of $x_0$ and one to the right is given by $H(x-x_0)$.

Even though this function is not a continuous one, we can still find its Fourier series.  This raises an interesting fact regarding the convergence of Fourier series that was not initially understood: The infinite sum of continuous and smooth functions can be discontinuous!  We will put further discussion of this off for now, and simply go forward with the idea that even discontinuous functions can have well-defined Fourier series that converge to the function in some (yet to be specified) sense.  Suppose that we have the Heaviside function shifted to the right by one-half.  Recall, by definition we have

\begin{align}
    f(x)&=\sum_{n=0}^{\infty} B_n \sin(n \pi x) \\
   B_n& =  2\int\limits_{x=0}^{x=1} H(x-\tfrac{1}{2}) \sin(n \pi x) dx, ~~ n=0,1,2,3,\ldots
\end{align}
Thus, noting that this integral is identically zero from zero to $\tfrac{1}{2}$, and unity for $x\ge\tfrac{1}{2}$ we have
\begin{equation}
    B_n=2 \int\limits_{x=\tfrac{1}{2}}^{x=1} \sin(n \pi x) \,dx 
\end{equation}
This integral is not difficult to do; the result is
\begin{equation}
    B_n = \frac{2}{n \pi} \left[\cos\left( \frac{n \pi}{2}\right)-\cos(n \pi)\right]
\end{equation}
The function and its approximation are plotted in Fig.~\ref{heavisideplot0}.  

{
\centering\includegraphics[scale=.75]{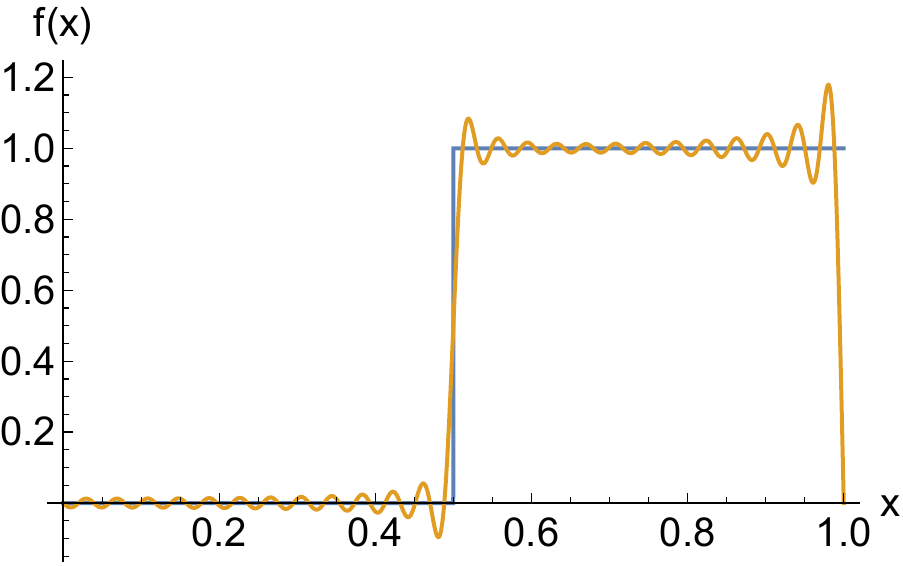}
\captionof{figure}{ The Heaviside function $H(x-\tfrac{1}{2})$ (blue) and the approximating cosine series  with $n=50$ terms (orange).} 
\label{heavisideplot0}
}
\vspace{4mm}
%
You may notice that at the discontinuities, this approximation seems to overshoot and undershoot the actual value by a factor of nearly 10\%; this is illustrated by the ``spikes" observable at $x=1/2$ and $x=1$.  The reason for this phenomenon is somewhat complicated, but it is a consequence of the use of a finite number of terms to approximate the series.  While the overshoot and undershoot are not eliminated by the use of more terms, the \emph{width} of the spikes diminishes to zero as $n\rightarrow \infty$.  This is known as \emph{pointwise} convergence (as discussed in the review given in Chapter 2-- see \S \ref{convergence_types}); recall that this is distinctly different from \emph{uniform} convergence.  While we will not study convergence extensively in this chapter, it is important to understand the concept of convergence, and how different modes of convergence might look on graphs of the function with increasing $n$.
\end{example}
\end{svgraybox}

\section{Fourier Cosine Series}\indexme{Fourier series!cosine series}

In the previous section, we showed that for ``nice" smooth functions, we can use nothing other than calculus to develop the representation for the Fourier sine series.  As you might imagine, the same thing is possible for the cosine, and the procedure is roughly the same.  

For the sine series, we found that the set of basis functions $~E_{\sin}=\{\sin{0 \pi x},~\sin(\pi x),$ $~\sin(2 \pi x), ~\sin(3 \pi x), \ldots\} $ allowed us to construct a sine series on $x\in[0,1]$ of the form
\begin{equation*}
    f(x)=\sum\limits_{n=1}^{n\rightarrow \infty} B_n \sin(n \pi x)
\end{equation*}
And, using the \emph{orthogonality} of the functions $\sin(n \pi x)$ and $\sin(m \pi x)$ ($m,n~\in~\,\mathbb{N}$), we were able to determine an integral equation for \emph{all} of the amplitude coefficients, $B_n$.  Note that, although we excluded it from the sum, the zero function is, technically, one of the functions in this list of basis functions (given that it is the function associated with $\sin(0 \pi x)$).  Making this analogy, we might guess that the set of basis functions for the cosine series might be $~E_{\cos}=\{\cos{0 \pi x},~\cos(\pi x),$ $~\cos(2 \pi x), ~\cos(3 \pi x), \ldots\}$.  This points out one important distinction between the sine and cosine series: the basis function associated with $n=0$ is nonzero for the cosine series! To be specific, note that $\cos(0 \pi x) = 1$  This is a non-symmetry that creates a few headaches, mostly because we are always having to keep track of this fact.  However, we can continue forward as before.  To start, we suggest the following form for the cosine series

\begin{equation}
        f(x)=\sum\limits_{n=0}^{n\rightarrow \infty} A_n \cos(n \pi x)
\end{equation}
where now the series \emph{must} start at zero.  Again, by convention, this series is often written in a form that allows the sum to start at $n=1$ (just as is done for the sine series).  Somewhat ridiculously, I am even going to keep the $\cos(0 \pi x)$ function as such in the series

\begin{equation}
        f(x)=A'_0\cos(0 \pi x)+\sum\limits_{n=1}^{n\rightarrow \infty} A_n \cos(n \pi x)
   \label{coser}
\end{equation}
Here, I have adopted the unusual notation of $A'_0$ for the first term in the Fourier cosine series.   The reason for this is that in this text, we will represent this term in a way that is not conventional, but it is useful and leads to fewer mistakes. This will be discussed further below. 

In exactly an analogous way as for the sine series, we also have an orthogonality condition for the inner product for the cosine functions.  Note

\begin{align}
    \int\limits_{x=0}^{x=1} \cos(n \pi x) \cos(m \pi x) dx &=
    \left\{ {\begin{array}{*{20}{c}}
  {0~~if~m\ne n~~~~~~~} \\ 
  {\frac{1}{2}~~if~m=n\ne 0} \\
   {1~~if~m= n=0} \\ 
  \end{array}} 
\right. \\
\end{align}
\label{cosident}
Note that the case $n=m=0$ leads to a result that is different from the rest of the cases ($n,m>0)$.  As we did for the sine series, we are now going to work out a method to determine the values of $A_n$ (and $A'_0$) for the cosine series.  To start, we can multiply both sides of Eq.~\eqref{coser} by $\cos(m \pi x)$ and integrate.  This gives

\begin{equation*}
    \int\limits_{x=0}^{x=1}  f(x)\cos(m \pi x) \, dx=A'_0 \int\limits_{x=0}^{x=1}\cos(0 \pi x)\cos(m \pi x) \, dx+\sum\limits_{n=1}^{n\rightarrow \infty}  \int\limits_{x=0}^{x=1}A_n \cos(n \pi x)\cos(m \pi x) \, dx
\end{equation*}
Using the results from Eq.~\eqref{cosident}, we have (and, specifically, recalling that every term in this integrated version of the series is \emph{zero} except for the one term where $n=m$) the following two cases.

\begin{itemize}
\item Case 1.  $n=m=0$
\begin{align*}
    &\int\limits_{x=0}^{x=1}  f(x)\cos(0 \pi x) \, dx=A'_0  \\
    \intertext{so that}
    &A'_0= \int\limits_{x=0}^{x=1}  f(x) \, dx
\end{align*}

\item Case 2.  $n=m, \text{ and } n,m>0$
\begin{align*}
    &\int\limits_{x=0}^{x=1}  f(x)\cos(n \pi x) \, dx=\frac{1}{2}A_n  \\
    \intertext{so that}
    &A_n= 2 \int\limits_{x=0}^{x=1}  f(x)\cos(n \pi x) \, dx
\end{align*}
\end{itemize}
Here, we have used $m=n$ to express the final result in terms of $n$ (again, this is the conventional notation).  Now, if we look at the results for $A'_0$ and $A_n$ we notice something that is a little bit annoying.  The general form for $A_n$ actually works for $A'_0$ also, except the result would be two times too large.  So, there are two equally reasonable ways to proceed.  

\begin{enumerate}
\item We just remember the formula for $A'_0$ and for $A_n$ ($n>0$) separately, and go along on our way.

\item We remember only the formula for $A_n$, with the idea that it also works for $n=0$; but then we also have to redefine $A'_0$ as follows

\[ A'_0 = \frac{A_0}{2}  \]
\end{enumerate}
The reality is that both are fine, and both require that you remember two different things.  The use of $A_0$ is one adopted by most texts by convention. However, in practice, this is not the best approach.  Instead, I prefer to think of the first term of the Fourier cosine series as being the \emph{average value} of the function over the interval.  For the unit interval, the interval length is $L=1$.  Thus we have the corresponding average

\begin{equation*}
    A'_0 = \frac{1}{1}\int\limits_{x=0}^{x=1}  f(x) \, dx
\end{equation*}
Later, when we allow the intervals to be of arbitrary length, $I=[0,L]$, the first term in the Fourier cosine series will \emph{still be} the average value, computed for this more general case by

\begin{equation*}
    A'_0 = \frac{1}{L}\int\limits_{x=0}^{x=L}  f(x) \, dx
\end{equation*}

In summary, we have

\begin{align}
        f(x)&=\frac{A_0}{2}+\sum\limits_{n=1}^{n\rightarrow \infty} A_n \cos(n \pi x) \\
        A_n&= 2 \int\limits_{x=0}^{x=1}  f(x)\cos(n \pi x) \, dx, ~~n=0,1,2,\ldots \\
        \intertext{or, preferably }
        f(x)&={A'_0}+\sum\limits_{n=1}^{n\rightarrow \infty} A_n \cos(n \pi x) \\
        A'_0& = \int\limits_{x=0}^{x=1}  f(x)\, dx \\
        A_n&= 2 \int\limits_{x=0}^{x=1}  f(x)\cos(n \pi x) \, dx, ~~n=1,2,\ldots
\end{align}
where we \emph{always} compute the first term in any cosine series as the average value over the appropriate interval.

These results are essentially identical to those for the sine series, except that we always will have to deal with the $n=0$ term of the cosine series (which is identically zero for the sine series) whenever we use the cosine series.   This all makes cosine series just a little less fun, but ultimately it will be worth it.  We will discuss later reasons that we might prefer one series over another.

\begin{svgraybox}
\begin{example}[Fourier Sine and Cosine Series Compared]\label{sincoscompare}
The function $f(x)=x$ is a frequent example for sine and cosine series, mostly because the expressions for $A_n$ and $B_n$ are integrable.  To start, let's compute the sine series.  This is given by
\begin{align*}
    B_n&=2 \int\limits_{x=0}^{x=1} x \sin(n \pi x) \,dx \\
    \intertext{Integrating by parts gives}
    B_n&= -2\frac{ \cos(n \pi)}{n \pi}=-\frac{2(-1)^n}{n \pi}
    \intertext{So, the result is}
    f(x) &= \sum\limits_{n=1}^{n\rightarrow \infty} B_n \sin(n \pi x)\\
    x &= \sum\limits_{n=1}^{n\rightarrow \infty} -\frac{2(-1)^n}{n \pi} \sin(n \pi x)
\end{align*}
Now, for the cosine series we have
\begin{align*}
    A_n&=2 \int\limits_{x=0}^{x=1} x \cos(n \pi x) \,dx \\
    \intertext{Integrating by parts gives}
    A_n&= \frac{2(-1+(-1)^n)}{n^2 \pi^2}
    \intertext{Note that for $n=0$ this gives us 0/0, which is undefined.  Rather than deal with that problem, we go back to the original definition for the first term in the series}
    A'_0& = \int\limits_{x=0}^{x=1} x  \,dx = \frac{1}{2} \\
    \intertext{So, in conclusion, we have}
    f(x) &= A'_0 + \sum\limits_{n=1}^{n\rightarrow \infty} A_n \cos(n \pi x) \\
    x &=\frac{1}{2}+ \sum\limits_{n=1}^{n\rightarrow \infty} \frac{2(-1+(-1)^n)}{n^2 \pi^2} \cos(n \pi x)
\end{align*}
Thus, we have two \emph{different} series for the same function.  Is there any practical difference between them?  We can check by plotting these functions up for $n=10$ and comparing.

{
\centering\includegraphics[scale=.5]{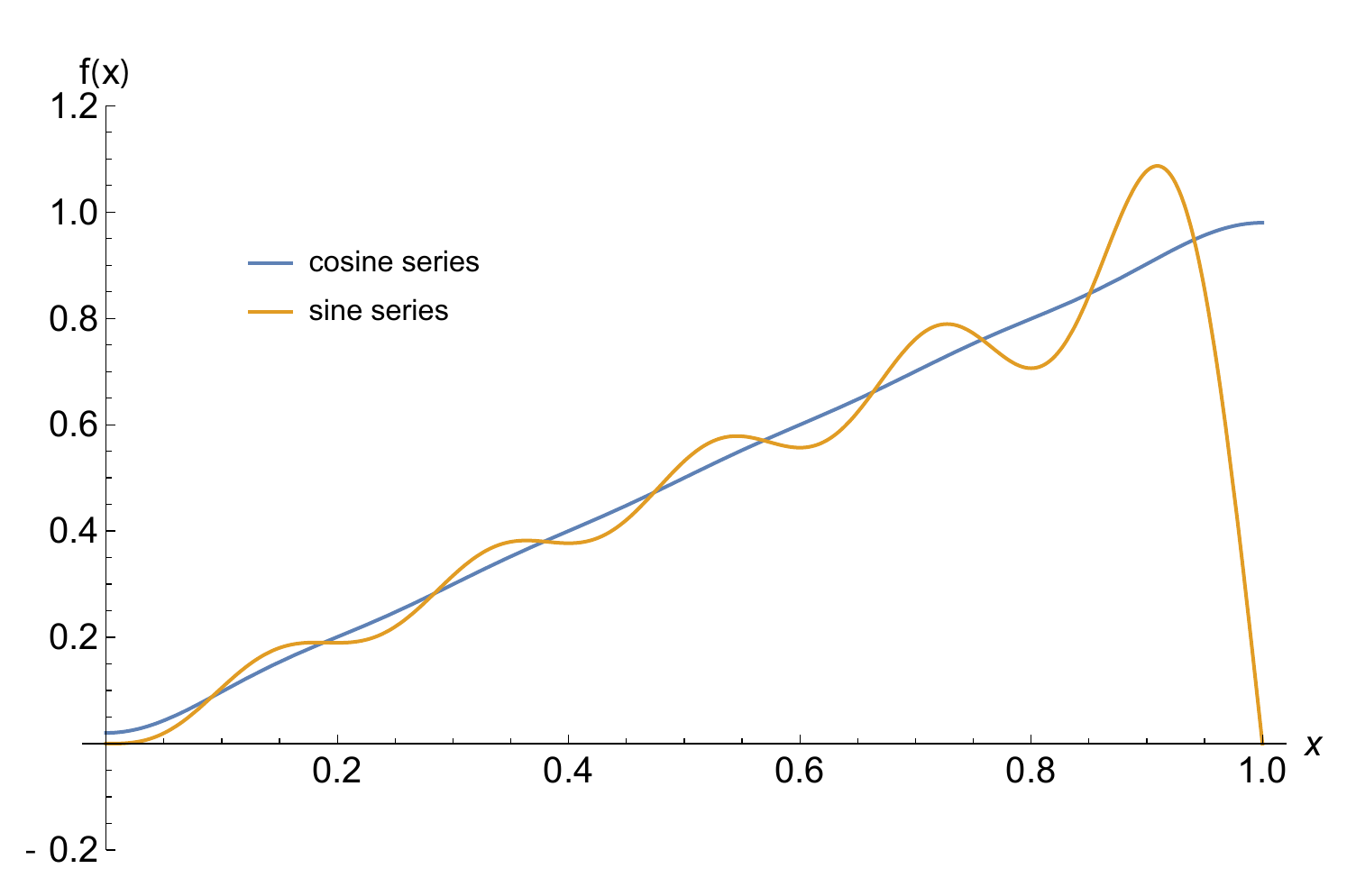}
\captionof{figure}{Cosine and sine series for the function $f(x)=x$, each with 10 terms.}
\label{trigfunctions}  
}
\vspace{4mm}
We have two different series for this function, and it is clear from observation that the behavior of the two series is quite different.  First of all, the series using the cosine basis functions seems to converge much better than the one using the sine series.  Second, the sine series is zero at the location $x=1$, which is not what we want.  As we increase the number of terms, we can see that this behavior never really goes away.  Here is a plot of the same two series with a total number of terms equal to 100 for each series

{
\centering\includegraphics[scale=.5]{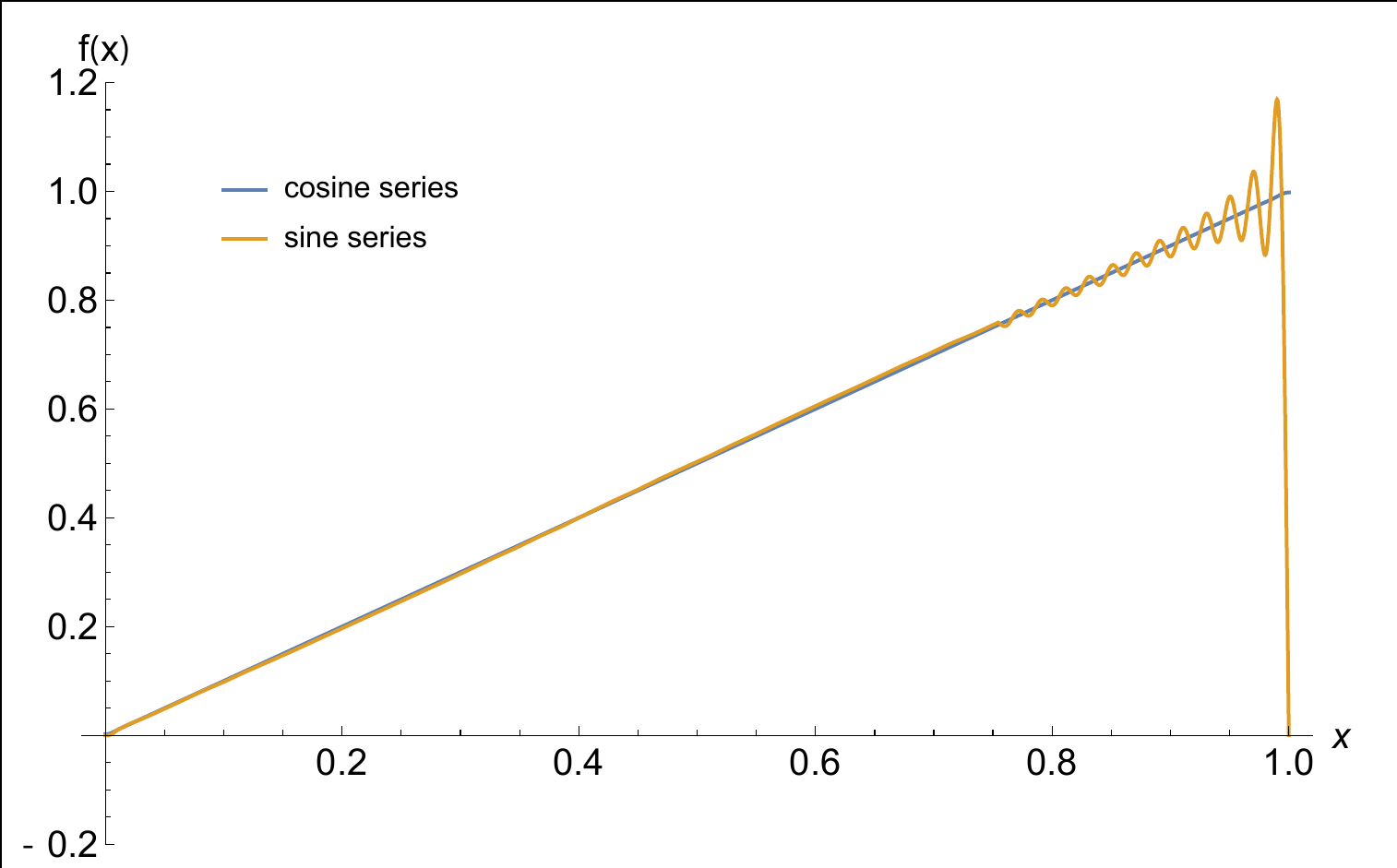}
\vspace{2mm}
\captionof{figure}{Cosine and sine series for the function $f(x)=x$, each with 100 terms.}
\label{trigfunctions2}  
}
\vspace{4mm}
As  you can see, something slightly odd happens near $x=1$ with the sine series.  The explanation is, at least in part, this: The sine series is forced to be zero at all multiples of $\pi$; thus each term in $\sin(n \pi x)$ is \emph{identically} zero when $x=1$.  The best that the sine series can do is to try to create a discontinuity that drops precipitously from 1 to zero as you approach $x=1$ from the left.  Which is what it, in fact, does do.  There are a few more problems that arise (for instance, one can observe some oscillations near the point $x=1$ with the sine series) that will be discussed later.
\end{example}
\end{svgraybox}

\section{Comparison of Basis Functions}\indexme{Fourier series!basis functions}

Both the Fourier sine and cosine series have a similar features.  For each, we can think of any function $f$ on a finite interval $D=[a,b]$ as being decomposed into its constituent amplitudes and frequencies, as represented by the weighted sine or cosine functions.  These functions form a \emph{basis}, as mentioned above, although it is a basis with an infinite number of basis ``vectors" (in this case, the basis vectors are the sine and cosine functions).  

It is helpful to see some of the parallels between the sine and cosine series.  In Table \ref{t:basis_table}, the basis functions, as a function of $n$, are listed.  There are relatively clear parallels between the two kinds of basis functions, with the exception of the $n=0$ case.  For that case, there is a lack of parallelism.  For the sine function, we have that $\sin(0)=0$, whereas for the cosine function,  $\cos(0)=1$.  So, technically, there is a basis function for the sine functions for $n=0$, however, because that function is just the zero constant function, it adds nothing to the corresponding series.  For the cosine series, we have a different constant function given by $\cos(0)=1$.  This function does contribute to the series by allowing the series to be translated by a constant.  This actually can make it easier for some problems; as an example, the cosine series for $f(x)=1$ is just the $n=0$ term (all other terms are zero).  For the sine series the expansion for $f(x)=1$ is substantially more complicated.  This distinction is made explicit by Eqs.~\eqref{constsine}-\eqref{constcos}.

\begin{align}
    &\textrm{sine series:}& &1 = \sum_{n=1}^\infty B_n \sin(n \pi x) =  \sum_{n=1}^\infty \frac{1-(-1)^n)}{n\pi} \sin(n \pi x)\label{constsine}\\
    &\textrm{cosine series:}& &1 = A'_0+\sum_{n=1}^\infty A_n \cos(n \pi x) = 1, \nonumber\\
    &&&\quad\textrm{ (i.e., $~A'_0=1$, $A_n=0$ for $n>0$)}\label{constcos}
\end{align}

\begin{table}[ht!]
\caption{Basis functions for sine and cosine Fourier series.}
\centering
\def\arraystretch{1.5}%
\setlength{\tabcolsep}{2mm} 
\begin{tabular}{|c|c|c|}
\hline
$~~n~~$   & $~~\sin(n \pi x)~~$  & $~~\cos(n \pi x)~~$\\ 
\hline
 0& $0$       & $1$  \\
 1&   $\sin(\pi x)$   & $\cos(\pi x)$ \\
 2&   $\sin(2\pi x)$  & $\cos(2\pi x)$\\
 3&   $\sin(3\pi x)$  & $\cos(3\pi x)$\\
 4&   $\sin(4\pi x)$  & $\cos(4\pi x)$\\
 \vdots & \vdots & \vdots \\
 $n$ & $\sin(n\pi x)$ & $\cos(n \pi x)$ \\
 \hline
\end{tabular}
\label{t:basis_table}
\end{table}

\section{Fourier Series Convergence: The Dirichlet Conditions}\indexme{Fourier series!Dirichlet conditions}

It would not be an overstatement to say that the understanding of what conditions are needed for Fourier series to converge were part of what motivated many components of modern mathematical analysis.  A full understanding of the convergence of Fourier series required new tools and definitions in mathematics that had not been required for analytic functions.  For example, the entire notion of what constitutes a \emph{function} was opened to question by the existence of Fourier series.  As an example, the Heaviside function is not an analytic function (hence, it has no Taylor series expression), but it could be constructed by Fourier series.  This required a re-examination as to what kinds of mappings could legitimately be called a function.  A similar existential crisis also developed for the notions of what it meant for an infinite series to \emph{converge}.  For analytic functions, the notions of convergence of Taylor series to the appropriate function were reasonably well understood.  However, for Fourier series, come conceptual curiosities were created.  Again, thinking of the Heaviside function,  it is clear from our example of the Fourier sine series expansion that it is possible to develop an \emph{approximation} to the Heaviside function by a Fourier series.  However, each term of the Fourier series expansion is formed by a weighted $C^\infty$-smooth sine function.  Thus there arises the following conceptual disconnect: how can a sum of \emph{smooth} functions generate a function that has a discontinuity (with no derivative at the point of discontinuity)?  Such realizations were part of what initiated a deeper understanding of what it meant for a series to converge.  It also led to the realization that when dealing with \emph{infinite} sums, the results may confound one's intuition.  

The conditions for which Fourier series can be understood to converge have been extended well beyond even the cases that originally caused consternation with mathematicians; interested readers can find very general treatments in the texts by \citet{zygmund1955}, \citet{lighthill1970}, and \citet{hardy1999}.  For our purposes, the most general possible extensions are not necessary.  One of the very first researchers on the topic of convergence of Fourier series was a German mathematician named Peter Gustav Lejeune-Dirichlet (for whom the ``Dirichlet" boundary condition is named), who in 1829 published a set of criterion for the convergence of Fourier series that covered a wide class of practical functions.  Dirichlet not only determined the conditions such that certain Fourier series converge, but he also developed a more modern concept of the function.  Notably, the Dirichlet conditions can be described concisely, and have a good physical interpretation.  Dirichlet's conditions can be stated as follows.

\begin{enumerate}
    \item The function $f(x)$ is defined on the range $x\in[-L,L]$, and bounded on this interval.  Note that functions defined on the half interval $x\in[0,L]$ can be periodically extended to cover the symmetric interval around zero.
    
    \item The function $f(x)$ has, at most, a finite number of points $x_{0,i}, i=1 \ldots N$ of discontinuity (i.e., a countable number of points where the left and right limits at $x_{0,i}$ approach $f(x_{-,i})$ and $f(x_{+,i})$ respectively).  At such points, $f(x_{0,i})$ is given the value
    \begin{equation}
        f(x_{0,i})=\frac{f(x_{-,i})+f(x_{+,i})}{2}
    \end{equation}
    that is, just the average value at the discontinuity. 
    
    \item The function $f(x)$ has a countable number of maxima or minima on the interval $x\in[-L,l]$.   
\end{enumerate}
The first condition is sometimes stated by requirement that the function be \emph{absolutely integrable}, i.e., 

\[ \int_{-L}^{L} |f(x)|\, dx < \infty  \]
and this definition will prove to be useful when the topic of Fourier transforms arises in Chapter \ref{Fouriertransforms}.
Requirement 3 above is sometimes extended to functions that have an infinite number of maxima and/or minima on the interval and have a property known as \emph{bounded variation}.  However, this extension was not part of Dirichlet's conditions, and this extension was due to the French researcher Camille Jordan \citep{jordan1881series}.  For the material following, we will assume that we investigate at least meet the Dirichlet conditions (unless otherwise specified).

\section{The Spectrum}\indexme{Fourier series!spectrum}

The spectrum of a Fourier series is simply the collection of amplitude coefficients $A_n$ or $B_n$ plotted versus the index, $n$.
Examination of the spectrum of a Fourier series has many practical applications (such as filtering in the Fourier domain), and it also helps us to better understand the physical meaning of the Fourier series.   

Suppose we look more at both the sine and cosine series for the function
\begin{equation}
    f(x)=x, \quad x\in[0,1]
\end{equation} 
that we examined above.  Recall, we had the following list of amplitudes as a function of $n$, where $n\in\mathbb{N}$ (i.e., $n=1,2,3,\ldots$) for the sine and cosine series for this function

\begin{align*}
&\text{cosine}&&
\left\{ {\begin{array}{*{20}{l}}
  {A'_0=\dfrac{1}{2}} \\ 
  {}\\
  {A_n=\dfrac{2(-1+(-1)^n)}{n^2 \pi^2}~~~~n=1,2,3,\ldots} 
 \\
 {}
\end{array}} \right.\\
&\text{sine}&&   ~~~~~B_n = \frac{-2(-1)^n}{n \pi}~~~~n=1,2,3,\ldots
\end{align*}
There is an interesting fact about these two series.  Although the functions that they represent are infinite dimensional (in the \emph{uncountable} infinity sense of the word infinite), the series representations are \emph{countably} infinite.  In other words, if you know all of the amplitudes, you can reconstruct the function.  In a sense, this represents a form of compression of the information embedded in a function.  In fact, Fourier series are sometimes used for exactly that purpose (e.g., in the compression of images in the JPEG format!)

\begin{figure*}[ht]
\sidecaption[t]
\centering
\includegraphics[scale=0.45]{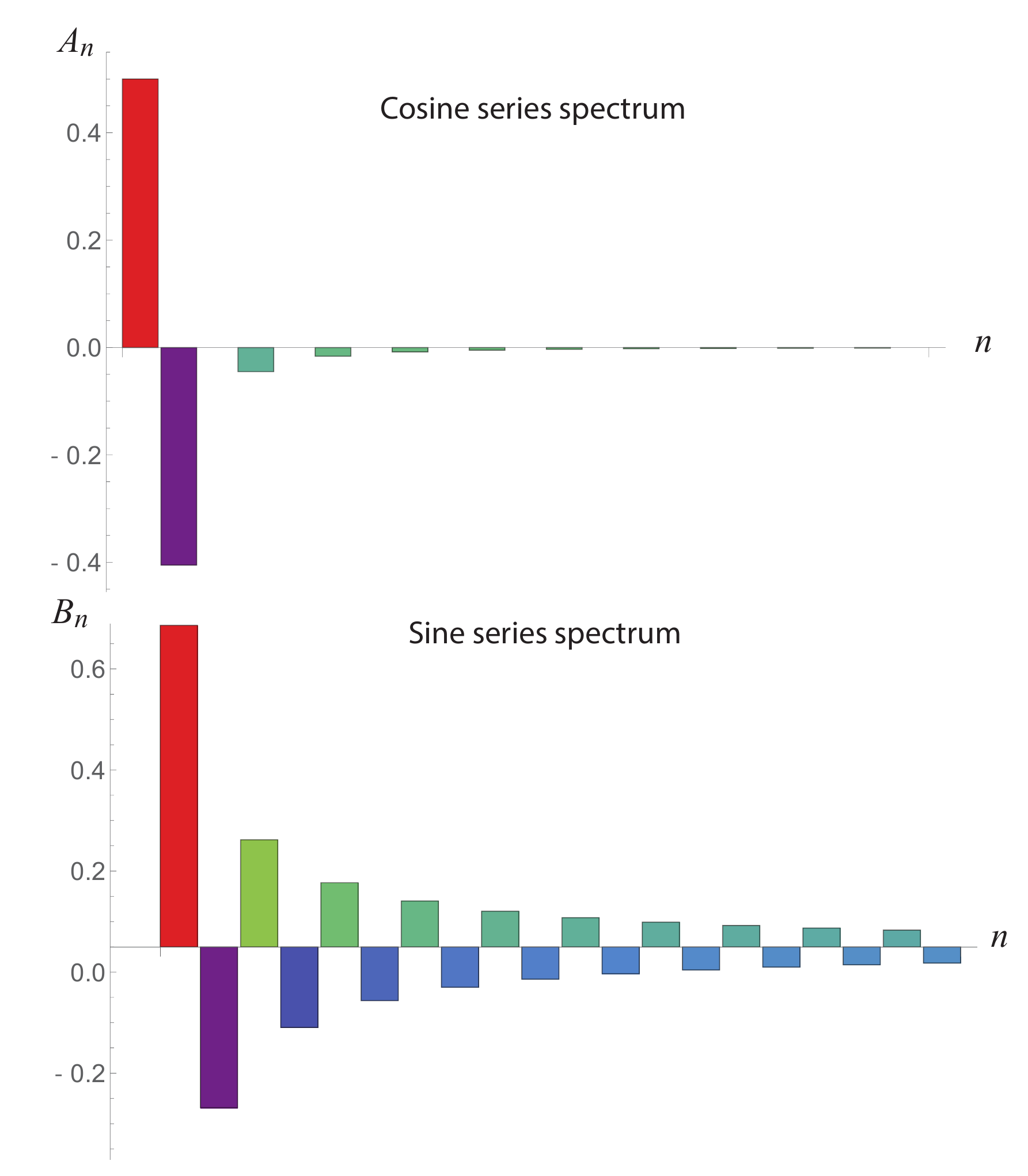}
\caption{The spectra for the cosine series (top) and sine series (bottom). The first 20 amplitudes are shown.}
\label{spectra}       
\end{figure*}

The amplitudes of the sine and cosine functions are useful for a number of reasons.  The main one is that they can give you some sense for how quickly a function converges.  The graph of $A_n$ or $B_n$ versus $n$ represents the graph of the amplitudes, and it is called the \emph{spectrum} of a sine or cosine series.  As an example, let's look at the two spectrum plots for the sine and cosine series for $f(x)=x~$ (Fig.~\ref{spectra}).

\begin{figure}[ht]
\sidecaption[t]
\centering
\includegraphics[scale=0.48]{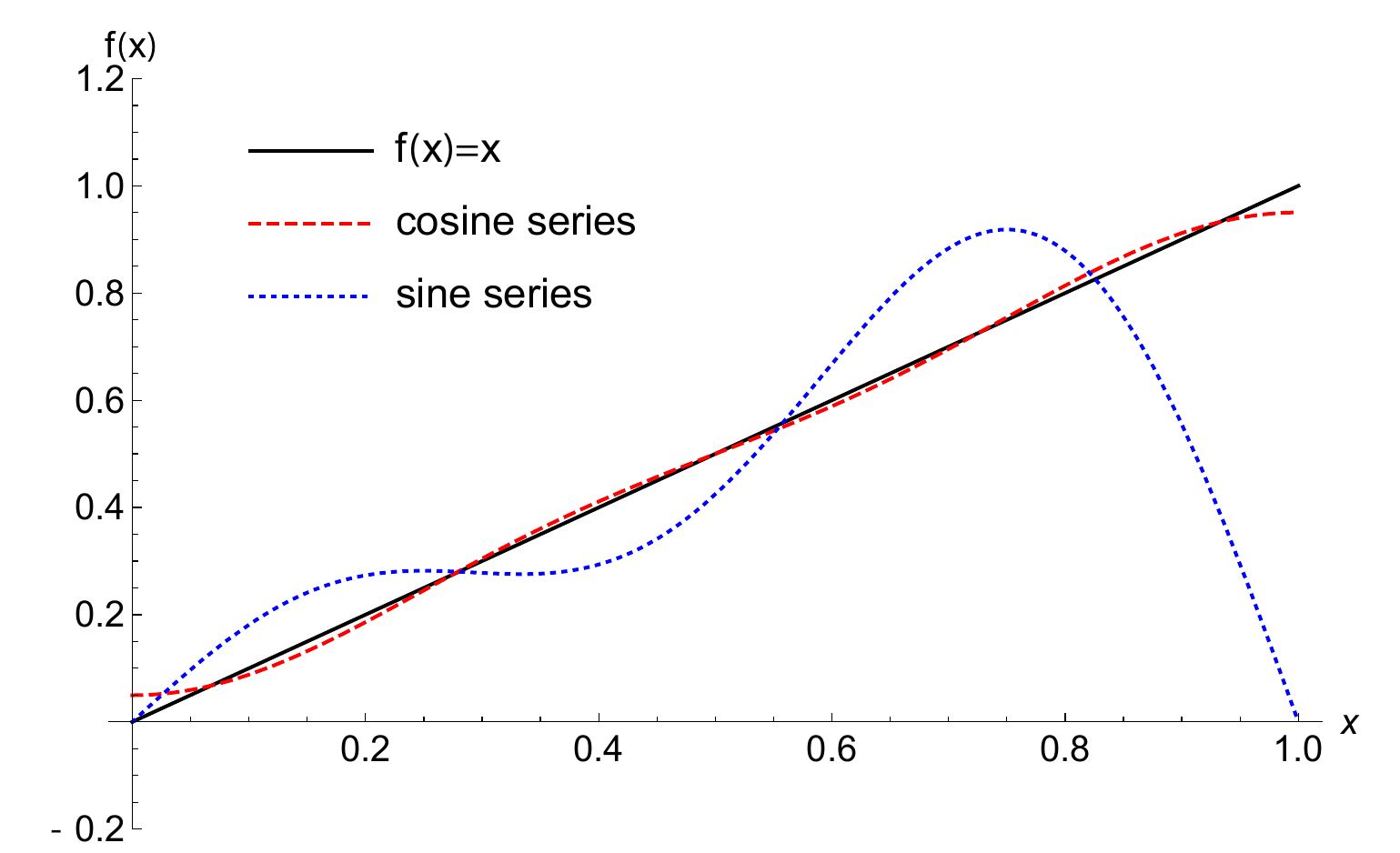}
\caption{Plots of the cosine and sine series, using only the first three nonzero amplitudes.}
\label{compare}       
\end{figure}

Looking at these two spectra, one can see several important differences.  The primary one is that the cosine series for this function is dominated by the first two or three (nonzero) amplitudes; the remaining amplitudes are quite small in comparison, and they decay (become smaller) very quickly.  For the sine series, although it is true that the first two amplitudes are the largest (in magnitude), they don't exactly dominate the others, and they do not seem to decay as quickly.  In fact, it is not clear that any of the first 20 amplitudes for the sine series might be neglected compared to the first two.  This gives us some guidance as to how well the two series converge.  By examining these plots, one might expect that the cosine series would give good results even if it included only the amplitudes $A_0,~A_1$ and $A_3$ ($A_2$ and all other even amplitudes being zero).  Whereas the same cannot be said for the sine series.  To check this, the plot of the cosine and sine series constructed from the first three (nonzero) amplitudes are given in the plot shown in Fig.~\ref{compare}.
Examining the results of Fig.~\ref{compare} suggests that our assertion about the relative sizes and the rate of decay of sizes of the amplitudes was largely correct.

\vspace{-2mm}

The spectrum also has significant practical applications.  Whenever one is tackling problems that contain time-series that have periodic components (e.g., tide heights over several weeks, atmospheric pressure over several days, hourly temperatures collected over a month), the spectrum can tell you much about what frequencies (frequency =$n/2$) are important.  Similarly, one can use a spectral analysis to help analyze data for \emph{noise} and even build filters to remover noise.  This is explored in the next example.

\begin{svgraybox}
\begin{example}[Noise filtering of a DC voltage]
Suppose you are working in a lab, and a piece of equipment is turned off and on using a direct current (DC) signal.  When the equipment gets senses a large voltage spike (say, $20<V<50$ volts) it turns on.  When it senses a small voltage spike ($2<V<10$) it turns off.  Voltage changes smaller than about $2 V$ are ignored. To make this piece of equipment work, you need a nice clean source of DC voltage that can send voltage spikes of various sizes.  Suppose that you have a DC source, and you send out low- and high-voltage spikes about every 1 second.  You measure its voltage output for (what you expect to be) a sequence of $25V$ and $5V$ voltage spikes occurring every 1 s or so.  You record this information (an oscilloscope is a device that can do this), and you find that your DC signal is actually really not very clean at all.  In fact, it is contaminated with all kinds of noise.  Suppose it looks like the plot in Fig.~\ref{voltage1}.  This is kind of a disaster, because there is so much noise in the signal that some of the small pulses sometimes generate voltage that is over $10V$ (thus, they would not turn off the device!)

You show your results to the local person-who-can-build-anything, and they say that you just need a DC filter to remove the high-frequency contamination coming from the rectifier (which is a device that converts AC current from an outlet into DC current). They can build one for you as long as you can suggest what frequencies to filter out.

{
\centering\fbox{\includegraphics[scale=.35]{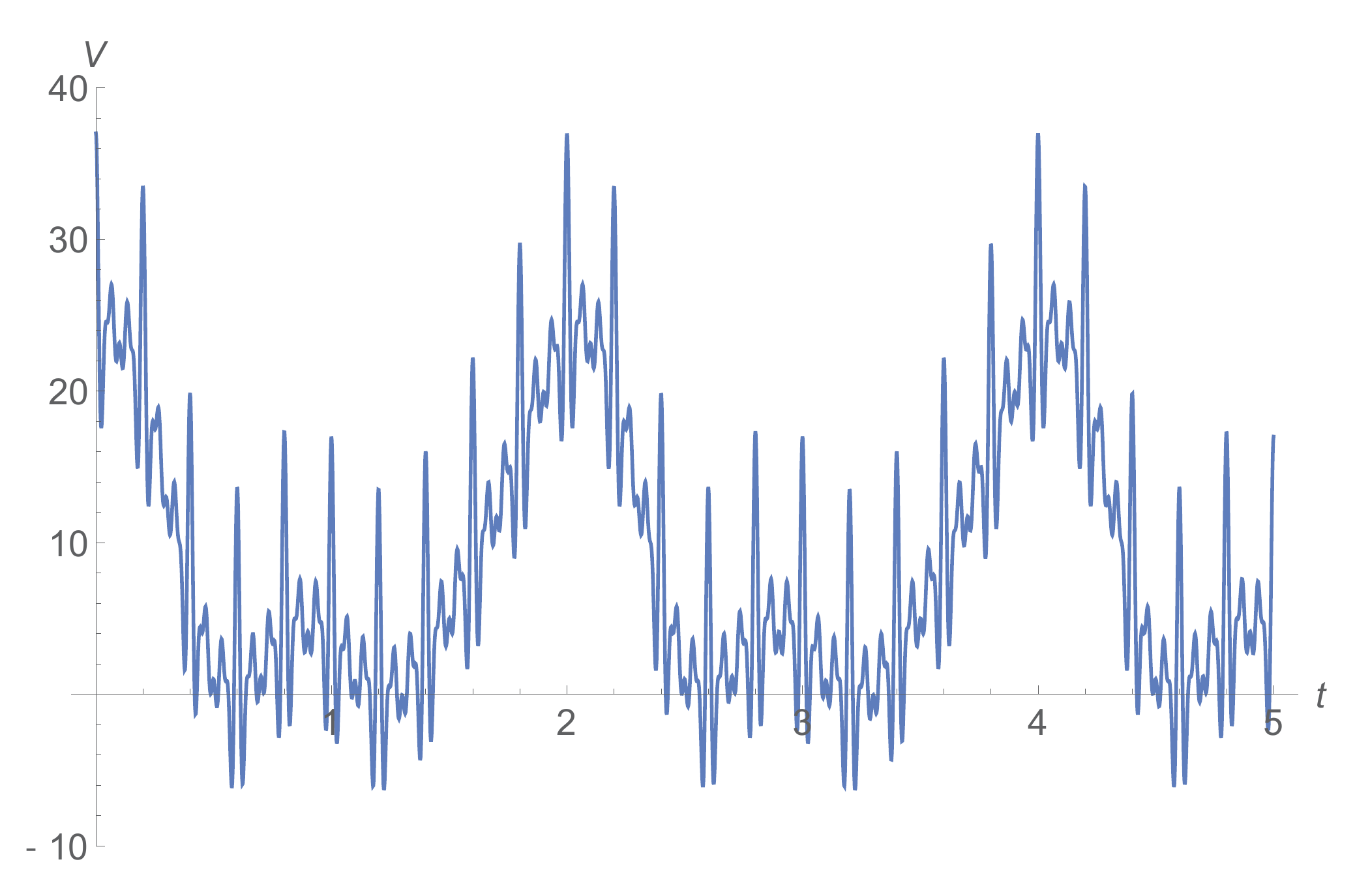}}
\vspace{-1mm}
\captionof{figure}{A noisy voltage-versus-time signal.}
\label{voltage1}  
}
{
\centering\fbox{\includegraphics[scale=.44]{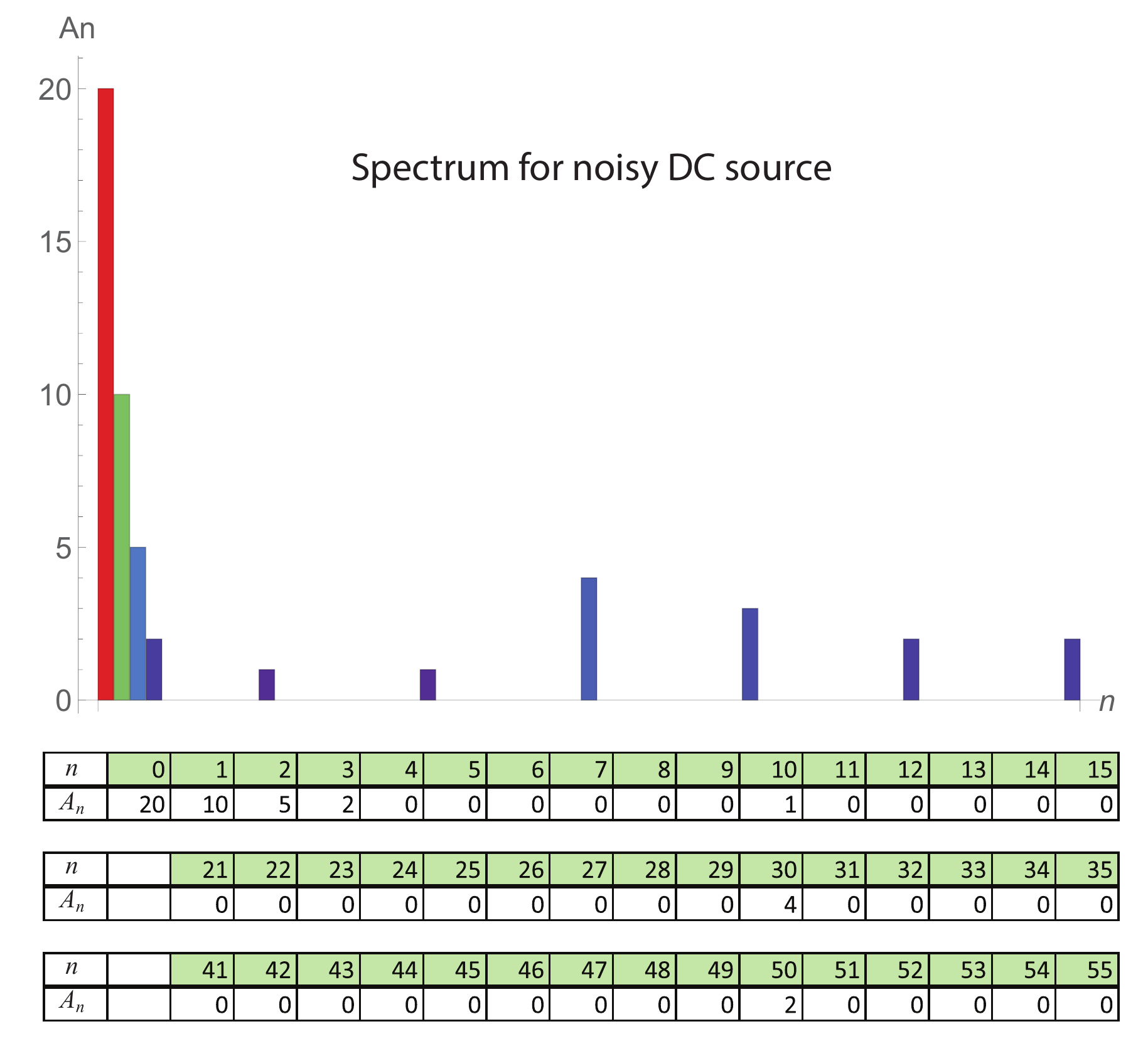}}
\vspace{-1mm}
\captionof{figure}{A DCT spectrum for your noisy current data.}
\label{spectrum2}  
}

Knowing what you know about Fourier series, you decide that what you really need is the spectrum for your DC current.  There are many ways of computing spectra directly from raw data; the most popular one (which will not be described in detail!) is call the \emph{Fast Fourier Transform} or FFT.  Actually, you will use a fast version of what is called the discrete cosine transform (DCT) so that you only get cosine amplitude-frequency information.  Being a Mathematica whiz, you compute the spectrum, and you find the spectrum shown in Fig.~\ref{spectrum2}.  Looking at this, you immediately see the problem!  You expect the signals to have a period around $1~s$, (i.e., $n\approx 2$, but you see that there are a number of components with an $n$ much greater than this value.  The solution seems clear: you just need to develop a filter to ``cut off" all of the high frequencies.  Looking at the plot, you decide it is safe to cut off all of the values where $n>10$ (which means that you decided to keep the first little bump in the spectrum plot).  To check to see if it this is sufficient, you even compute the new signal based on this filtering.  You can basically reconstruct the correct cosine signal by taking the first four non-zero values of $A_n$.  You find the result

\begin{equation*}
    V(t) = 9 + 10 \cos (\pi  x)+5 \cos (2 \pi  x) + 2 \sin (3 \pi  x)+\cos (10 \pi  x)
\end{equation*}
Plotting this will allow you to assess how well your filter will work.  The plot appears as Fig.~\ref{voltage3}.  Looking at these results, you decide that your cut off filter suggestion will work.  You get "on" and "off" peaks of the right magnitude, and the voltage fluctuations are small enough that they will be ignored.

{
\centering\fbox{\includegraphics[scale=.5]{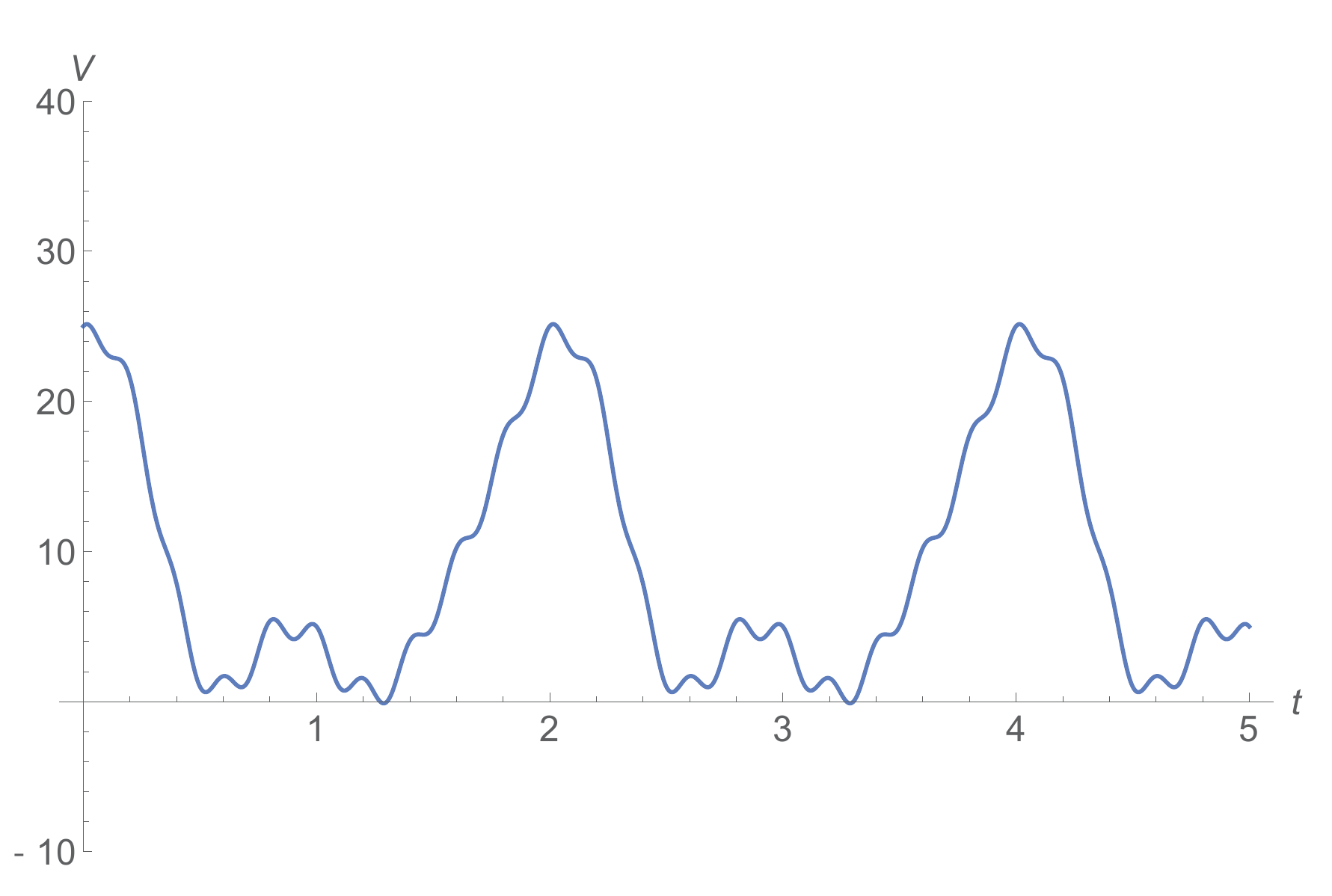}}
\vspace{2mm}
\captionof{figure}{A nicely filtered voltage signal!}
\label{voltage3}  
}

\end{example}
\end{svgraybox}

\section{Change of Interval}\indexme{Fourier series!change of interval}

\begin{figure}[t]
\sidecaption[t]
\centering
\includegraphics[scale=0.6]{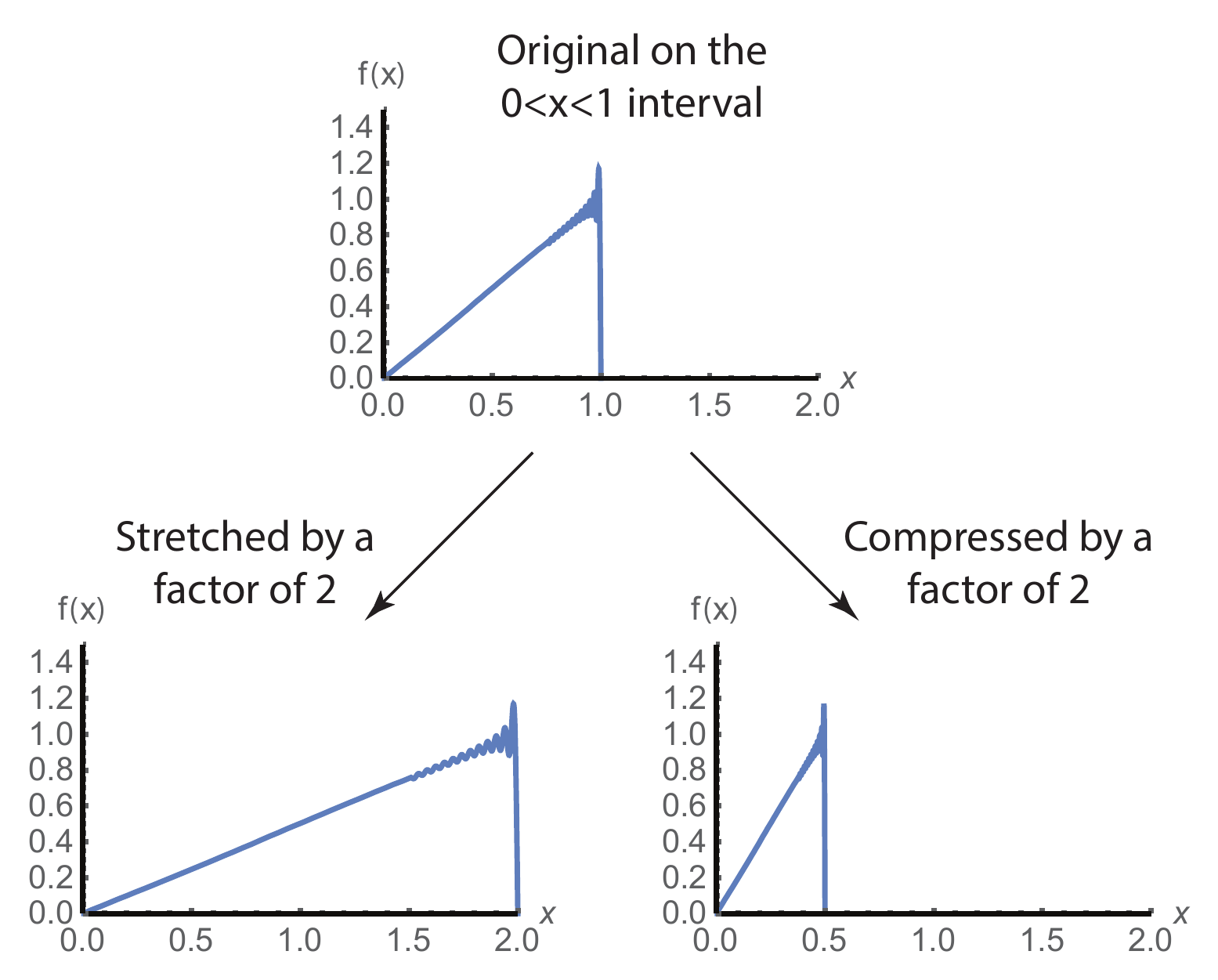}
\caption{Rescaling the domain of the Fourier series.}
\label{stretch}       
\end{figure}

In the examples above, we have examined the Fourier series exclusively on the interval $x\in[0,1]$.  Of course, we might want to define the Fourier series on some more general interval, say  $x\in[0,L]$.  This change is actually not all that difficult.  When we think about it, changing to the new interval basically means that the function that occurs on $x\in[0,1]$ is either stretched or compressed along the $x$-axis so that it now occurs on $x\in[0,L]$; the vertical behavior of the function remains unchanged (except that it is mapped to these stretched or compressed coordinates).  Fig.~\ref{stretch} gives a rough idea of the process. 

Making this change in the domain is not as difficult as it might seem.   To start, think about the following variable

\[  x=\frac{z}{L} \]
where $z\in[0,L]$.  While $z$ on the right-hand side goes from $0$ to $L$, the variable $x$ on the left-hand side only goes from $0$ to $1$.  We have essentially \emph{mapped} the $[0,L]$ onto the $[0,1]$ interval with this transformation.  Why is this useful? Well, let's consider the Fourier sine series for the function $f(x)=x$ that we derived above for $x\in[0,1]$.  The result was

\begin{align*}
    f(x) &= \sum\limits_{n=1}^{n\rightarrow \infty} -\frac{2(-1)^n}{n \pi} \sin\left(n \pi x\right),~~x\in[0,1] \\
    \intertext{Now, consider the series }
      f\left( \frac{z}{L}\right) &= \sum\limits_{n=1}^{n\rightarrow \infty} -\frac{2(-1)^n}{n \pi} \sin\left(n \pi \frac{z}{2}\right),~~x\in[0,2] \\
\end{align*}
This is basically the same series, except the interval is now twice as long.  In fact, recognizing that the variable $z$ is just a symbol for the independent variable, we could re-label the variable in this new function back to $x$ if we like.  Also, we technically do not need to put parameters like $L$ inside our notation for the function, just the list of independent variables.  To make this clear, we can write

\[ f\left( \frac{x}{L}\right) = g(x) \]
So, our longer-interval series can be equivalently written (in terms of the independent variable $x$)

\begin{align*}
      g(x) &= \sum\limits_{n=1}^{n\rightarrow \infty} -\frac{2(-1)^n}{n \pi} \sin\left(n \pi \frac{x}{2}\right),~~x\in[0,2] \\
\end{align*}
Note, this is no longer the Fourier series for $f(x) = x$; it is the Fourier series for we have defined the Fourier series for a \emph{new} function

\[ g(x) = \frac{1}{2} x \]
If you think about what happens here, the value of $g(1)$ is the same as the value of $f(1/2)$; similarly $g(2)$ is the same as the value of $f(1)$.  So, the net result is the original function $f$ is now stretched out over an interval that is twice as long.  This actually is the case that is plotted in Fig.~\ref{stretch} (upper plot and lower left plot).  So, that is really all that there is to the process of changing the interval.

If you think about what has happened here, it is really not all that complicated.  In short, we have replaced the set of orthogonal basis functions $E_{\sin}=\{0,~ \sin(1 \pi x),~ \sin(2 \pi x),\cdots,~x\in[0,1]\}$  with the new set of basis functions $E_{\sin}=\{0,~ \sin(1 \pi x/L),~ \sin(2 \pi x?),\cdots,~x\in[0,L]\}$.

\begin{figure}[t]
\sidecaption[t]
\centering
\includegraphics[scale=0.6]{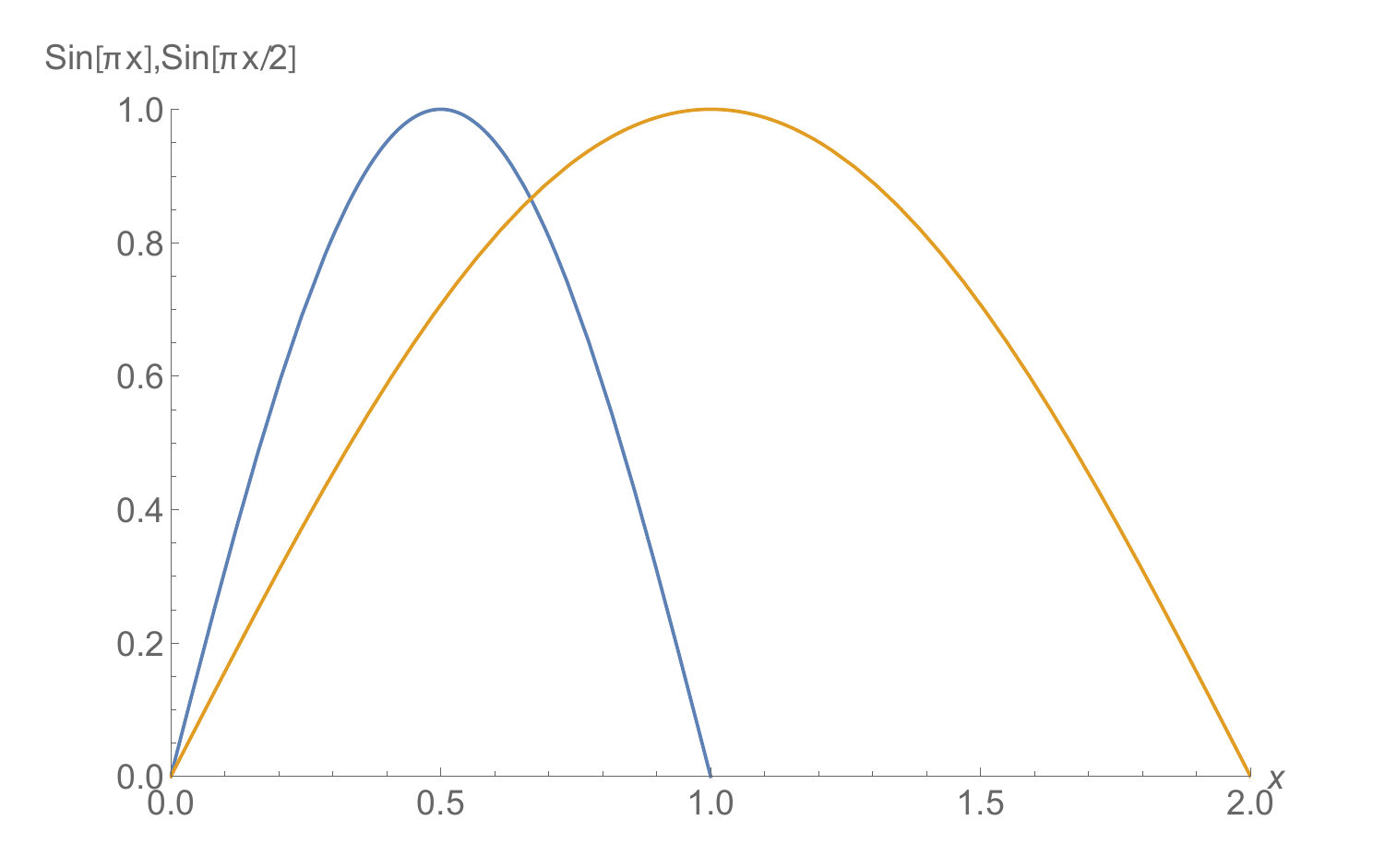}
\caption{Two sine functions with $n=1$.  One spans twice the distance of the other, but they both correspond to the first \emph{frequency} of the sine series.}
\label{twosine}       
\end{figure}

Now, we have on little detail to clean up.  Suppose you are given a series on an interval that is not $x\in[0,1]$ to begin with.  You could (a) convert the function to an interval that covered $x\in[0,1]$ by an appropriate transformation of variables, and and carry on as we did above, or (b) we could just re-derive the expressions for $B_n$ for the larger interval.  Generally, the second of these two options is going to be much more convenient.  To do this, we just re-create the process that we did with the unit interval.  We end up with the following integrals we need to evaluate

\begin{align}
    \int_{x=0}^{x=L} \sin(n \pi x/L) \sin(m \pi x/L) = 
    \left\{ {\begin{array}{*{20}{c}}
  &0 &&n\ne m\\ 
  &\frac{L}{2}&&n=m
\end{array}} \right.
\end{align}
And, we ultimately find that for $x\in[0,L]$
\begin{align*}
    B_n&=\frac{2}{L} \int\limits_{x=0}^{x=L} f(x) \sin\left(n \pi \frac{x}{L}\right) \,dx
\end{align*}
Note that this reduces to the same expression for the unit interval when $L=1$!

\section{Fourier Series on Symmetric Intervals around Zero}\indexme{Fourier series!symmetric intervals around zero}

Although we have been examining Fourier sine and cosine series on positive intervals such as $I=[0,L]$, in actuality the natural domain for Fourier sine and cosine series are symmetric intervals around zero, $I=[-L,L]$.  To understand this, we need primarily just two concepts

\begin{enumerate}
    \item The concept of  \emph{odd} versus \emph{even} functions.
    \item The concept that the functions $sin(n \pi x/L)$ and $\cos(n \pi x/L)$ each repeat with a period of $2L$ (N.B., for $n=1$, the functions repeat exactly once in an interval of $2L$; for $n>1$ the functions repeat multiple times in an interval of $2L$).
\end{enumerate}
We start the discussion with definitions of odd and even functions.

\subsection{Even and Odd Functions}\indexme{Fourier series!even and odd functions}\indexme{function!even and odd}

So far, we have discussed only Fourier series on the interval $[0,L]$, for some positive value for $L$.  In general, there is absolutely nothing preventing a Fourier series to be defined on \emph{any} finite interval $I=\{x:a < x < b\}$ where $a$ and $b$ are any real numbers.  However, including the origin in the domain has the potential to create some technical problems for finding Fourier sine and cosine series.  Before starting the discussion on Fourier series on more general intervals, we will discuss the concept of \emph{odd} and \emph{even} functions, and their relationship to Fourier series.

\begin{definition}
Suppose a function $f$ is defined on some \emph{symmetric interval} around origin ($x=0$), i.e.,  $I=[-L,L]$, where $a$ and is some positive real number.  Then the function may be \emph{even}, \emph{odd}, or \emph{neither even nor odd} as follows.
\begin{itemize}
    \item The function $f$ is said to be \emph{odd} on $I$ if $f(x)=-f(-x)$ for all $x\in I$.
    \item The function $f$ is said to be \emph{even} on $I$ if $f(x)=f(-x)$ for all $x\in I$.
    \item If neither of these two conditions is true, then $f$ is neither even nor odd on the interval $I$.
\end{itemize}  
\end{definition}
The concepts of even and odd functions are more than just a obscure mathematical property.  It turns out that, for every single function defined on a symmetric interval around zero can be decomposed into an even part and and odd part whose sum gives the original function.  The proof of this is intuitive and clever, so it will be explained in a few lines.

\begin{figure}[t]
\sidecaption[t]
\centering
\includegraphics[scale=0.55]{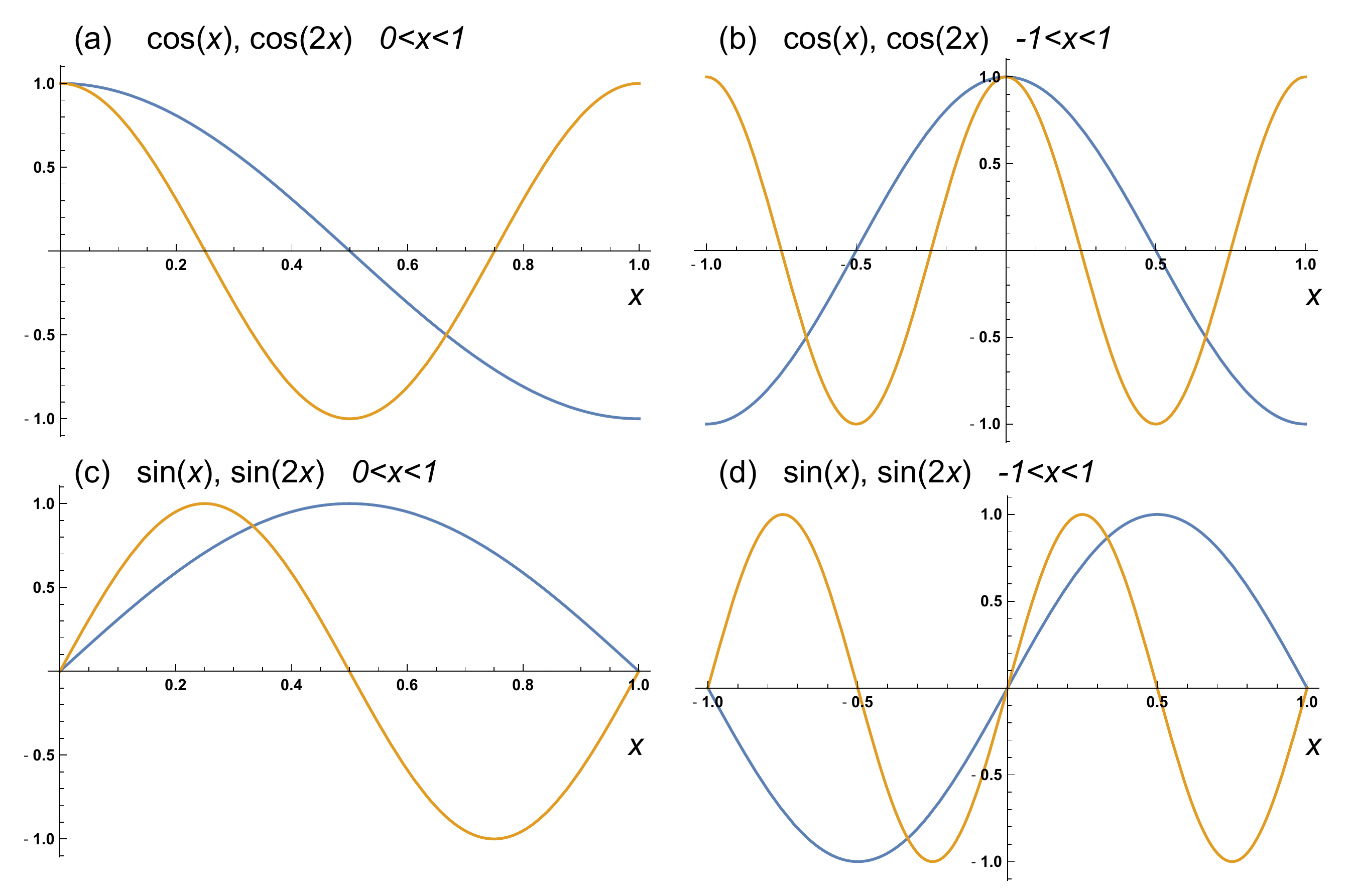}
\caption{Sine and cosine functions on the intervals $0<x<1$ and $-1<x<1$.  From the definitions in the text (and by examining parts (b) and (d) of this figure), it should be clear that cosine is an even function and sine is and odd function on the interval $-1<x<1$.}
\label{sinecosine}       
\end{figure}

\begin{theorem}
Every function, $f$, that is neither even nor odd on a symmetric interval $I$ around zero can be decomposed into a sum of two functions on that interval, one of which is even, and the other odd.
\end{theorem}
\begin{proof}
Define two new functions as follows
\begin{equation*}
    f_{even}(x)=(f(x)+f(-x))/2 \qquad  f_{odd}(x)=(f(x)-f(-x))/2
\end{equation*}
It is easy to verify that $f_{even}(-x)=f_{even}(x)$ and $f_{odd}(-x)=-f_{odd}(x)$, and thus these two functions are even and odd, respectively.  It is also easy to verify that $f_{even}(x)+f_{odd}(x)=f(x)$.  
\end{proof}
It turns out that this decomposition into even and odd components is also \emph{unique}; that is, there is no other such decomposition, only the one that we constructed above.  Examples of the sine function and the cosine function on the interval $-1<x<1$ is given in Fig.~\ref{sinecosine} ((b) and (d)).  This graphical presentation allows us to note that the functional definitions given above lead to geometrical properties of even and odd functions that can be observed in their graphs.  In particular, note the following

\begin{itemize}
    \item The cosine function is an \emph{even} function ($\cos(-\pi x) = \cos(\pi x)$).  Graphically, we can generate the cosine function on $[-1,0]$ by reflecting the cosine function on $[0,1]$ about the vertical axis.  
    
    \item The sine function is an \emph{odd} function ($\sin(-\pi x) =- \sin(\pi x)$).  Graphically, we can generate the sine function on $[-1,0]$ by reflecting the sine function on $[0,1]$ about the vertical axis, and then a second time about the horizontal axis.   
\end{itemize}

There are a few additional features about even and odd functions that should be mentioned.  The first is an assessment of what happens when one multiplies even and odd functions.  As an example of even and odd functions to keep in mind, think of the even function $f_{even}(x)=x^2$ and the odd function $f_{odd}(x)=x^3$ on $x=[-1,1]$.  Now we state the following theorems (without proof)

\begin{figure}[t]
\sidecaption[t]
\centering
\includegraphics[scale=0.4]{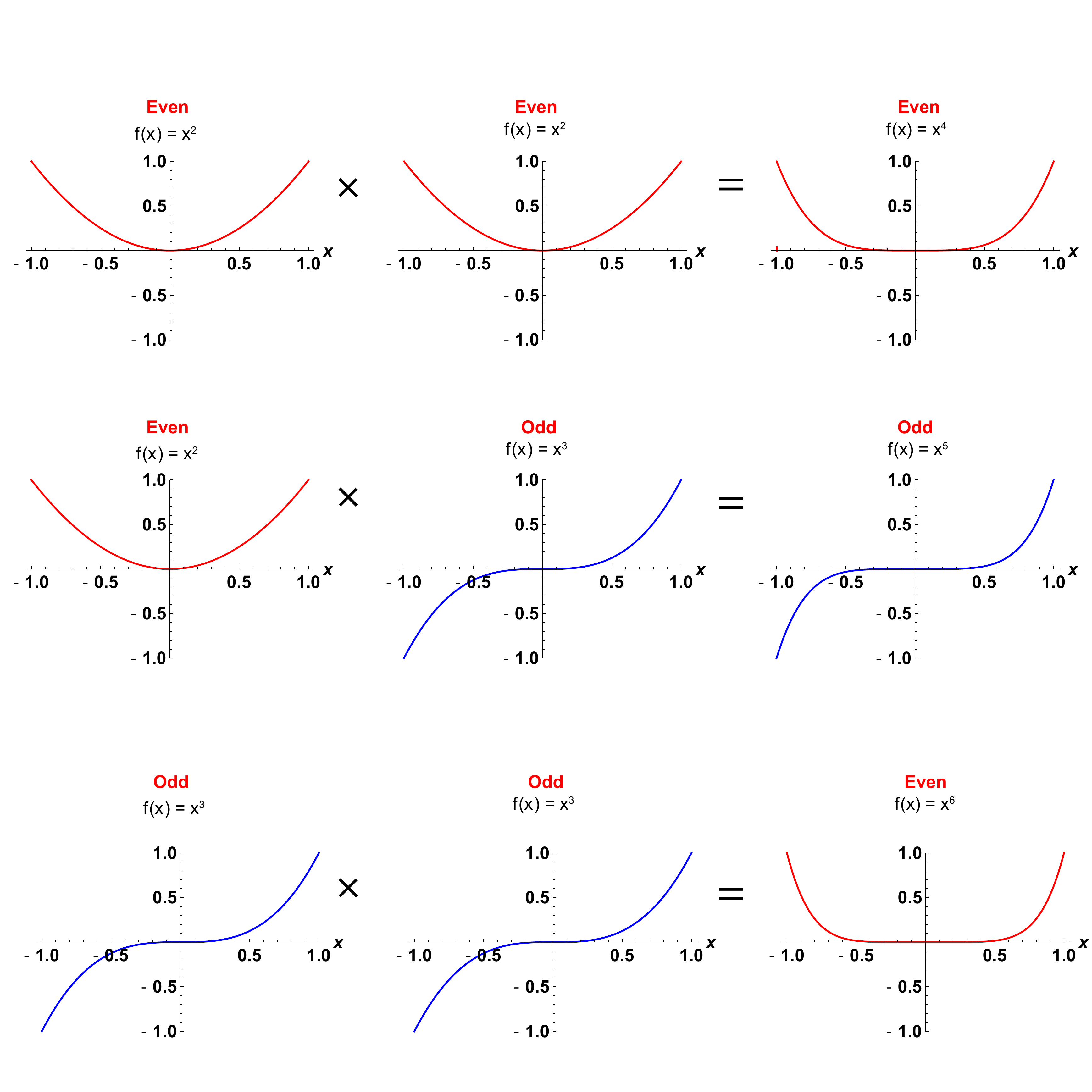}
\caption{The result of multiplying even and odd functions.  In this figure, the multiplication happens from left to right, across rows. }
\label{evenodd}       
\end{figure}

\begin{theorem}[Products of even and odd functions]
Suppose we have and even function $f_{even}(x)$ and an odd function $f_{odd}(x)$ defined on the interval $[-L,L]$.  Then, the products of these functions have the following characteristics
\begin{enumerate}
    \item $f_{even}\times f_{even}$ is an even function.
    \item $f_{even}\times f_{odd}$ is an odd function.
    \item $f_{odd}\times f_{odd}$ is an even function.
\end{enumerate}
\end{theorem}
For easy visualization, these results are summarized graphically in Fig.~\ref{evenodd}.

\begin{theorem}[The integral of an even and odd functions]
Suppose we have and even function $g_{even}(x)$ and an odd function $g_{odd}(x)$ defined on the interval $[-L,L]$.  Then the integral of $g_{even}$ over  the interval $[-L,L]$ is non-zero.  The integral of $g_{odd}$ over  the interval $[-L,L]$ is zero.
\end{theorem}
The reason that these two theorems are important, is because they explain what we can and cannot do with Fourier sine and cosine series.  In particular, recalling that the sine series is odd, and the cosine series is even, then we have the following results.  Suppose $f(x)$ function on the interval $I=[-L,L]$.  Then, the following things must be true
\begin{enumerate}
    \item If $f(x)$ is even, then it can be expanded only as a cosine series.  A sine series will not work because each integral would involve an even times an odd function, which is odd; the resulting integrals defining $B_n$ would be zero.
    \item If $f(x)$ is odd, then it can be expanded only as a sine series.  A cosine series will not work because each integral would involve an odd times an even function, which is odd; the resulting integrals defining $A_n$ would be zero.
    \item If $f(x)$ is neither even nor odd, then it can be thought of as being the sum of two functions $f(x)=f_{even}(x)+f_{odd}(x)$.  The Fourier series must have both sine and cosine components to represent the entire function.
\end{enumerate}

%
Now that we have an understanding of even and odd functions, we can consider finding the Fourier transform of general functions on symmetric intervals around zero.

\subsection{Fourier Series on the Interval $I=[-L,L]$}
Now consider the symmetric interval $I=[-1,1]$.  It turns out that the relevant basis functions for this are $\sin(n \pi x)$ and $\cos(n \pi x)$, just as they were for the unit interval $I=[0,1]$.  To see this more clearly, examine the behavior of both the sine and cosine functions on the interval $I=[-1,1]$ as shown in Fig.~\ref{sinecosine}.  It is easy to see that the interval $I=[-1,1]$ allows the sine and cosine functions to go through an integer number of periods; in fact, the functions on the interval $[-1,0]$ are just reflections of the functions on the interval $[0,1]$ as discussed above.  The basis functions on the interval $I=[-1,1]$ are exactly the same as those on $I=[0,1]$; specifically, the functions $\sin(n \pi x)$ and $\cos(n \pi x)$.  

There is one slightly tricky point about this change of interval.  Note that for a symmetric interval $I=[-L,L]$, we can map this interval to the interval $I=[-1,1]$ simply by dividing independent variable by $L$.  In other words, the new basis functions become $\sin(n \pi x/L)$ and $\cos(n \pi x/L)$.  This is worth noting, because we are still dividing by $L$, which corresponds to only one-half of the full interval width (the interval width is $2L$).  This is unlike the case of the change of interval on $I=[0,L]$, where we divided by the entire width of the interval, $L$.  

The best way to think about the relationship between Fourier series on $I=[0,1]$ versus those on $I=[-1,1]$ is to think of them as the \emph{same} series, where we are simply ignoring the negative part of the interval.   In fact, this is exactly what we are doing when we compute a Fourier series for a function $f(x)$ on an interval $I=[0,1]$.  If we expand $f(x)$ as a sine series, then we are ignoring an \emph{odd} extension of the solution on the interval $[-1,0]$.  If we expand $f(x)$ as a sine series, then we are ignoring an \emph{even} extension of the solution on the interval $[-1,0]$.  This can be made a bit more clear through an example.  Before proceeding, however, we note the following integrals on the general interval $I=[-L,L]$

\begin{align}
    \int_{x=-L}^{x=L} \sin(n \pi x/L) \sin(m \pi x/L) = 
    \left\{ {\begin{array}{*{20}{c}}
  &0 &&n\ne m\\ 
  &L&&n=m
\end{array}} \right.
\end{align}
\begin{align}
    \int\limits_{x=-L}^{x=L} \cos(n \pi x/L) \cos(m \pi x/L) dx &=
    \left\{ {\begin{array}{*{20}{c}}
  {0~~if~m\ne n~~~~~~~} \\ 
  {L~~if~m=n\ne 0} \\
   {2L~~if~m= n=0} \\ 
  \end{array}} 
\right.
\end{align}
Note that the interval is now $2L$, but the bounds of the integral go from $-L$ to $L$. This sometimes creates some confusion, so be aware of the details here.

These two integrals are still \emph{orthogonal}, but the case for $n=m$ now evaluates to twice the value we had for the unit interval.  This makes sense, since the total domain of integration is twice that of the unit interval.  Following this process through to the evaluation of $A_n$ and $B_n$ leads to the (hopefully unsurprising) results

\begin{align*}
    A'_0 &= \frac{1}{2L}\int_{-L}^{L} f(x) dx \\
    A_n &=\frac{1}{L}\int_{-L}^{L} f(x)\cos(n \pi x/L) dx,~~n=1,2,3,\ldots \\
    B_n &=\frac{1}{L}\int_{-L}^{L} f(x)\sin(n \pi x/L) dx~~n=1,2,3,\ldots\\
    \intertext{with}
    f(x)&={A'_0}+\sum_{n=1}^{\infty} A_n \cos(n \pi x/L) \\
    \intertext{or}
    f(x)&=\sum_{n=1}^{\infty} B_n \sin(n \pi x/L) \\
\end{align*}
Recall, it is easy to produce the formula for the first term of the cosine series ($A'_0$) if you remember that it is just the \emph{average} of the function.
\begin{svgraybox}
\begin{example}[Sine series for \mbox{$f(x)=x$ on $I=[0,1]$ versus on $I=[-1,1]$}]\label{intervalcompare}
Lets return to the familiar case of the sine series expanded fo the function $f(x)=x$.  To begin, lets compute the sine series for this function on the interval $I=[0,1]$.  Following the details above, we have
\begin{align*}
    B_n=2\int_{0}^{1} x \cos(n \pi x) dx = \frac{-2(-1)^n}{n \pi}
\end{align*}
Which gives us the series
\begin{align*}
    x=\sum_{n=0}^{n=\infty}\frac{-2(-1)^n}{n \pi}\sin(n \pi x),~ ~~\textrm{for $0<x<1$}
\end{align*}
This should be familiar; we have computed this series before.  Now, let's try the steps that we did previously for the unit interval for the interval $I=[-1,1]$. We start with the definition of the series expansion on $I=[-1,1]$
\begin{align*}
    x=\sum_{n=0}^{n=\infty}B_n \sin(n \pi x), ~~~\textrm{for $-1<x<1$}
\end{align*}
As done previously, we use orthogonality to determine the $B_n$.  Multiplying both sides of the last equation by $sin(m \pi x)$ and integrating gives us
\begin{align*}
\int_{-1}^{1}x\sin(m \pi x)\,dx&=\sum_{n=0}^{n=\infty}B_n \int_{-1}^{1}\sin(n \pi x)\sin(m \pi x)\, dx\\
\intertext{and, using the orthogonality of the sine functions, only one term in the sum on the right hand side is non-zero; this corresponds to the case where $n=m$}
B_n&=\int_{-1}^{1}x\sin(m \pi x)\,dx 
\intertext{Finally, computing this integral gives}
B_n&=\frac{-2(-1)^n}{n \pi}\\
\intertext{which gives us the series}
 x&=\sum_{n=0}^{n=\infty}\frac{-2(-1)^n}{n \pi}\sin(n \pi x),~ ~~\textrm{for $-1<x<1$}
\end{align*}
\vspace{-5mm}
{
\centering\fbox{\includegraphics[scale=.65]{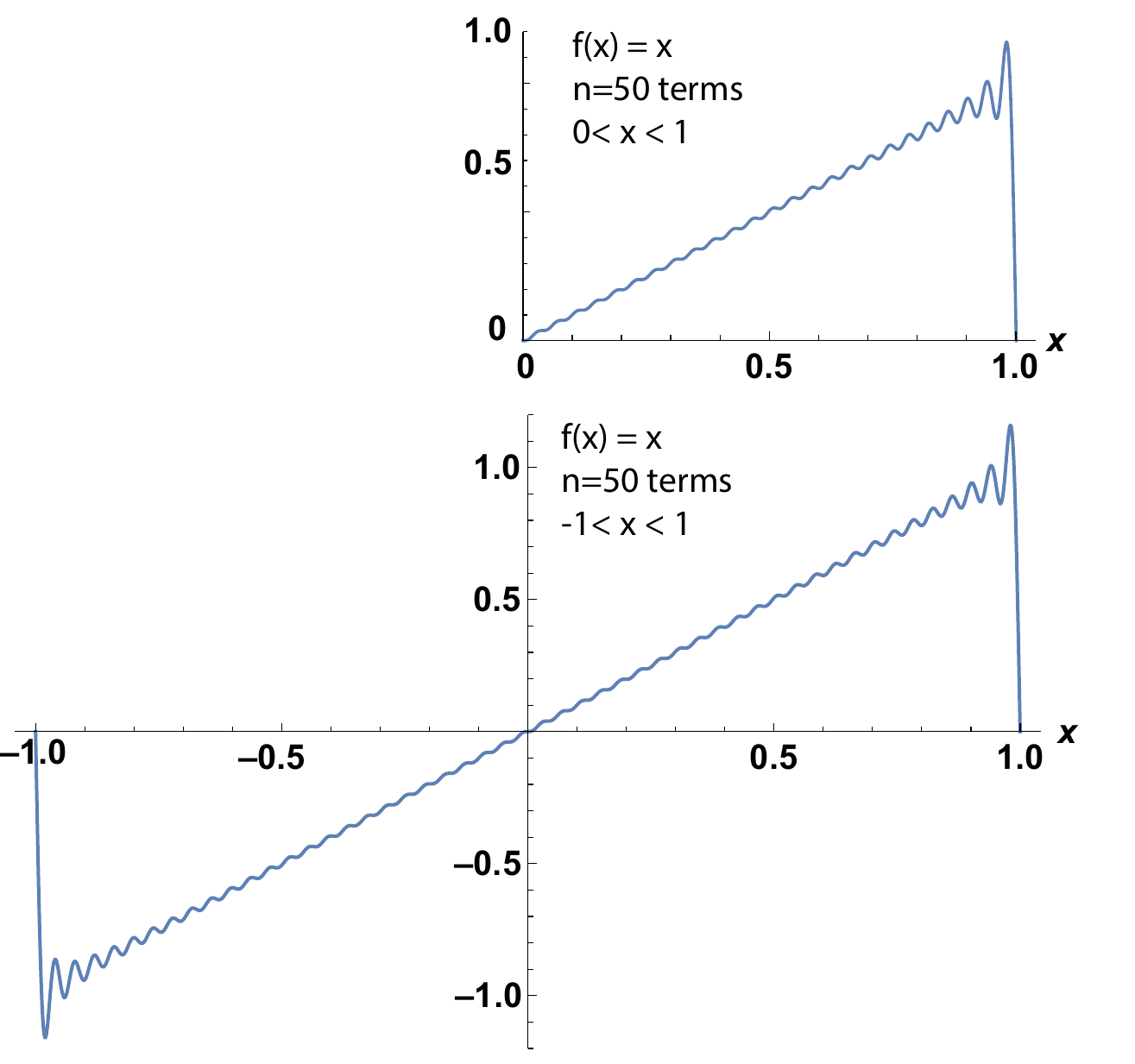}}
\vspace{2mm}
\captionof{figure}{Comparison of the Fourier sine series for the function $f(x)=x$ on the domains $I_1=[0,1]$ and $I_2=[-1,1]$.  Note that sine is and odd function, thus the function on the symmetric domain $I_2=[-1,1]$ is also odd.\vspace{5mm}}
\label{intervals}  
}
\noindent Comparing this with the result for the unit interval, we find that the two series are identical!  Thus, when we expand a function as a sine series on the interval $I=[0,1]$, we really are generating the same result as for the interval $I=[-1,1]$, it is just that we are ignoring the negative portion of the solution (because, by definition of the domain being $I=[0,1]$, it is not in the domain that is of interest to us).  A plot can help compare the solution over the two domains; such a plot given in fig.~\ref{intervals}. 
\end{example}
\end{svgraybox}

The case where a function is neither even nor odd is an interesting one.  To start, we will \emph{attempt} to expand the function using only a sine series and then only a cosine series.  We will find that \emph{neither} series is able to expand the function by itself.  However, the \emph{sum} of the two series will reproduce the function.  In short, the sine series reproduces the \emph{odd} part of the function, and the cosine series reproduces the \emph{even} part of the function.  An example is really useful here.

\begin{svgraybox}
\begin{example}[Cosine and Sine series for \mbox{$f(x)=e^{-x}$ on $I=[-1,1]$}]

The function $f(x)=e^{-x}$ is neither even nor odd on $I=[-1,1]$; this is easy to see by looking at a plot of the function on this interval.

{
\centering\fbox{\includegraphics[scale=.45]{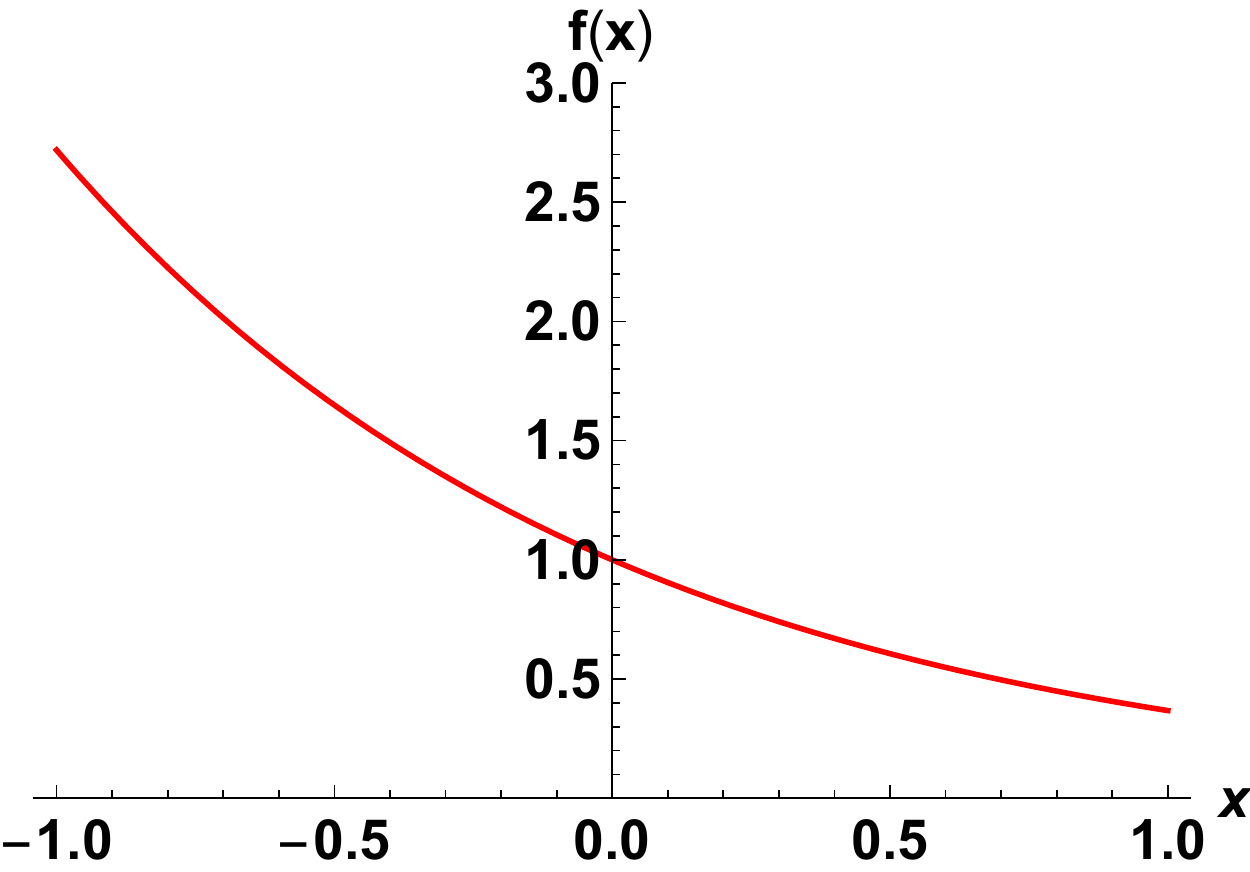}}
\vspace{2mm}
\captionof{figure}{The function $e^{-x}$ is neither even nor odd.\vspace{5mm}}
\label{exp}  
}
%
In principle, we could \emph{split} the function $f(x)=e^{-x}$ into its even and odd parts, and then compute the Fourier series for each of these parts.  However, the computation of the series coefficients $A_n$ and $B_n$ will actually take care of this for us; the multiplication by cosine or sine filters out the even and odd parts of the function (respectively) during the computation of these integrals. So, let's try computing the two parts as separate operations, starting with the cosine series.  For reinforcement of the ideas of orthogonality, this example will be shown with all steps developed (rather than simply relying on the formulas derived earlier).

\begin{align*}
    e^{-x} &= \sum_{n=0}^{n=\infty} A_n \cos(n \pi x), ~~\textrm{where for this problem $L=1$.}\\
    \intertext{Multiplying both sides by $\cos(m \pi x)$ and integrating gives}
    \int_{-1}^{1} e^{-x} \cos(m \pi x)dx &= \sum_{n=0}^{n=\infty} A_n \int_{-1}^{1}\cos(n \pi x)\cos(m \pi x)dx\\
\end{align*}
For the integral on the right, there are three cases that can occur: (1) $m\ne n$, in which case the integral is zero, (2) $n=m\ne 0$, in which case the integral is equal to 1, and (3) $n=m=0$, in which case the integral is equal to 2 (the length of the interval!)  Regardless of case, for each value of $m$ chosen, there is only \emph{one} nonzero value of the sum as $n$ goes from zero to infinity.  Thus, we can rewrite this equation as (and, remembering that because $n=m$, we are free to set the index to $n$ everywhere)
\begin{align*}
    \int_{-1}^{1} e^{-x} \cos(n \pi x)dx &=  A_m \int_{-1}^{1}\cos(n \pi x)\cos(n \pi x)dx\\
\end{align*}
We treat the $n=0$ case independently.  For that case, we find that the integral on the right gives the result of 2, and, solving for $A'_0$, we find

\begin{equation*}
    A'_0=\frac{1}{2} \int_{-1}^{1} e^{-x} \cos(n \pi x)dx =\sinh{(1)}
\end{equation*}
For the remaining values of $n$, we find

\begin{equation*}
    A_n= \int_{-1}^{1} e^{-x} \cos(n \pi x)dx =\frac{(-1)^n(-1+e^2)}{e(1+n^2 \pi^2)}
\end{equation*}
Together, these yield the series

\begin{equation*}
   f_{even}(x)=\sinh{(1)}+\sum_{n=1}^{n=\infty} \frac{(-1)^n(-1+e^2)}{e(1+n^2 \pi^2)} \cos(n \pi x)
\end{equation*}
Here, the function is labeled ``$f_{even}(x)$ because it is not in fact the function $f(x)=e^{-x}$; this can be seen in Fig.~\ref{fevenfodd}.  The associated Fourier cosine series \emph{must be} an even function, and we can verify that it is.  Conceptually, this funtion represents the ``even" component of the function $f(x)=e^{-x}$.

For the sine series, we will not reproduce the entire sequence of operations leading to the resulting series.  Instead, we perform operations analogous to those above, and find the result (presented previously)

\begin{equation}
    B_n = \int_{-1}^{1} e^{-x} \sin(n \pi x)dx =\frac{(-1)^n(-1+e^2) n \pi}{e(1+n^2 \pi^2)}
\end{equation}
%
{
\centering\fbox{\includegraphics[scale=.4]{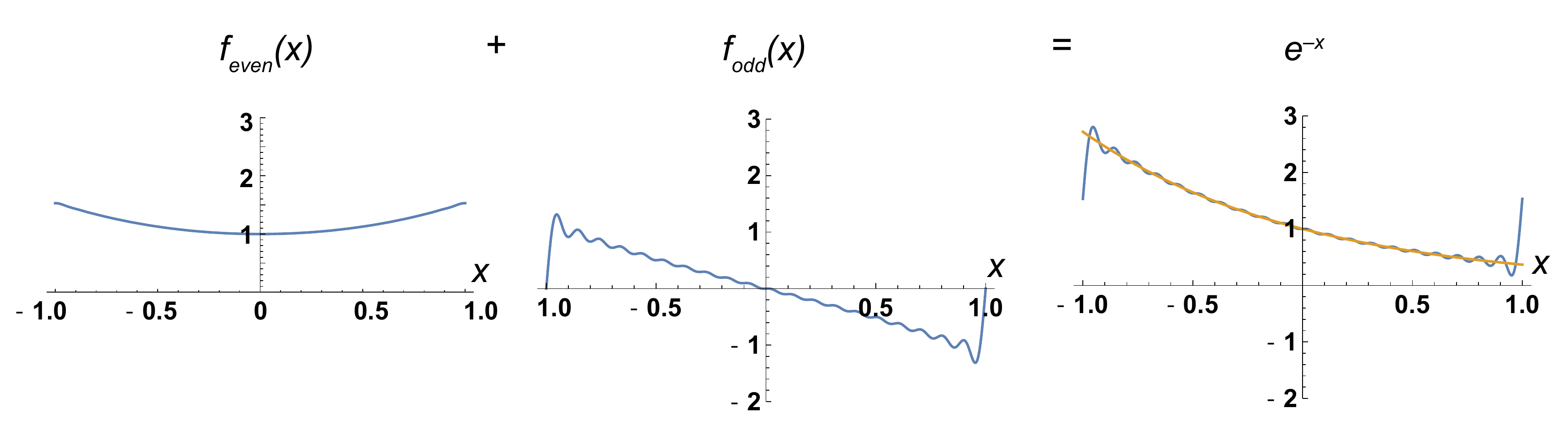}}
\vspace{2mm}
\captionof{figure}{The function $e^{-x}$ is the sum of an even (cosine) Fourier series plus and odd (sine) Fourier series \vspace{5mm}}
\label{fevenfodd}  
}
And, the corresponding series is

\begin{equation*}
     f_{odd}(x)=\sum_{n=1}^{n=\infty} \frac{(-1)^n(-1+e^2)n \pi}{e(1+n^2 \pi^2)} \sin(n \pi x)
\end{equation*}
This function is also plotted in Fig.~\ref{fevenfodd}. By inspection, it is clear that the sine series is in fact an odd function.

Finally, note that the sum of these two series actually reproduces the function that we seek

\begin{align*}
e^{-x} &=  f_{odd}(x)+ f_{even}(x)\\
&= \sinh{(1)}+\sum_{n=1}^{n=\infty} \frac{(-1)^n(-1+e^2)}{e(1+n^2 \pi^2)} \cos(n \pi x)
+\sum_{n=1}^{n=\infty} \frac{(-1)^n(-1+e^2)n \pi}{e(1+n^2 \pi^2)} \sin(n \pi x)
\end{align*}
Figure ~\ref{fevenfodd} illustrates the sum of these two series, showing that the sum of the sine and cosine series expansions for the function results in a total series that reproduces the function.
\end{example}
\end{svgraybox}
From the example above, we can draw some conclusions.
\begin{enumerate}
    \item Functions on symmetric intervals around zero are either even, odd, or neither even nor odd.  For the latter case, it is \emph{always} possible to decompose the function into unique even and odd parts, i.e., $f(x)=f_{even}(x)+f_{odd}(x)$.
    \item Even functions are represented exclusively by Fourier cosine series.
    \item Odd functions are represented exclusively by Fourier sine series.
    \item For functions that are neither even nor odd, we need to sum both the Fourier cosine and sine series so that both the even and odd components of the function are represented.  This  expression takes the general form
    
    \begin{align*}
        f(x)&=f(x)+f_{even}(x)+f_{odd}(x) \\
        f(x)&= A'_0 +\sum\limits_{n=1}^{n\rightarrow \infty} A_n \cos(n \pi x) 
        +\sum\limits_{n=1}^{n\rightarrow \infty} B_n \sin(n \pi x)
    \end{align*}
\end{enumerate}

\subsection{Periodic Continuation}\indexme{Fourier series!periodic continuation} \label{periodic_contin}

The examples above illustrated that even if we compute a Fourier series for $x\in[0,L]$, the Fourier series is actually defined for $x\in[-L,L]$.  The shape of the function in the negative part of the interval will depend on whether a sine or a cosine series is used to expand the function.  This is exactly what we saw in Example \ref{intervalcompare}.

It turns out that the Fourier series on an interval $[-L,L]$ is actually a \emph{periodic} function.  That is, outside of the interval $[-L,L]$, the function simply repeats by generating copies of the function that is defined for $[-L,L]$.  

\begin{figure}[t]
\sidecaption[t]
\centering
\includegraphics[scale=0.6]{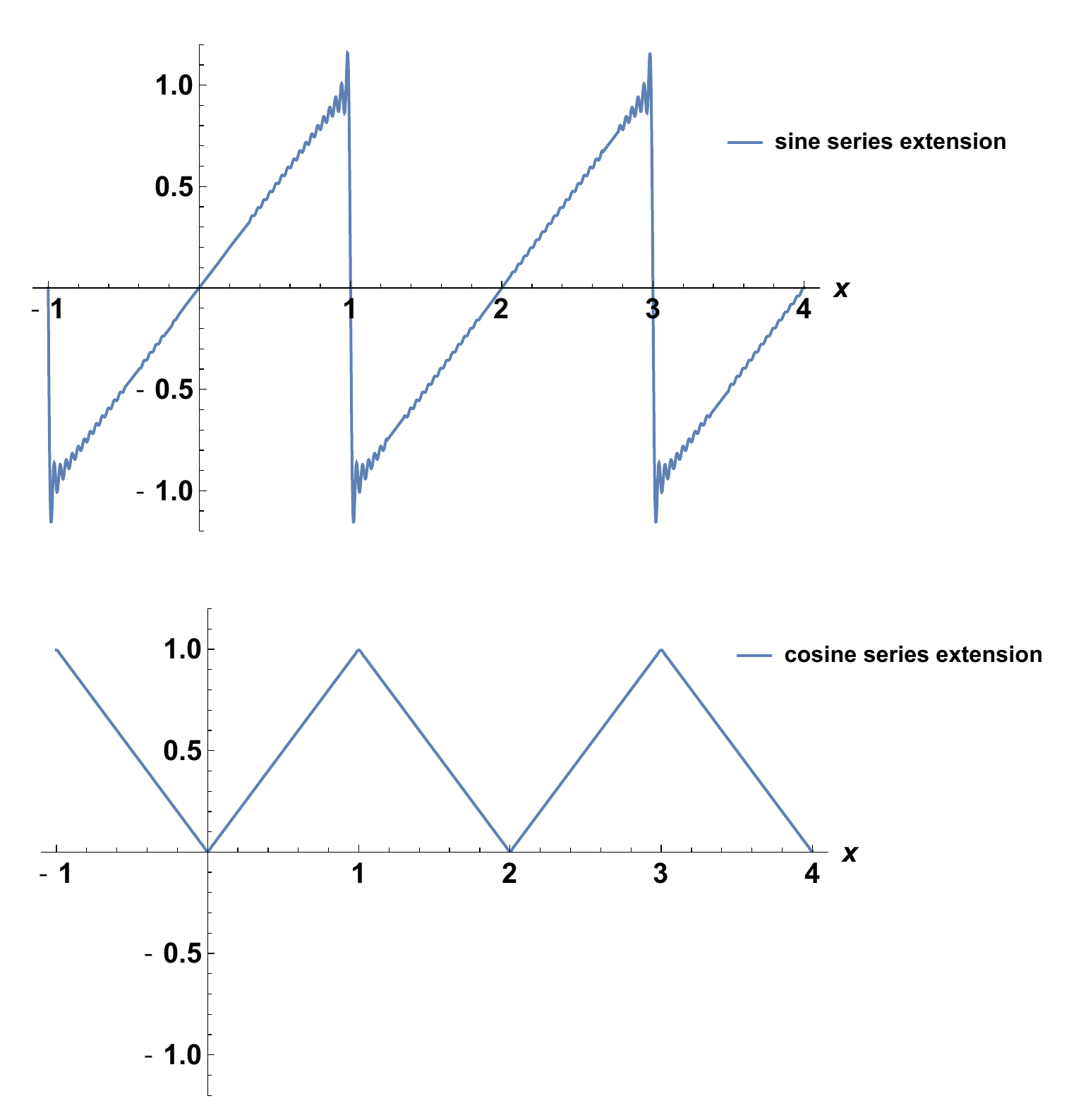}
\caption{The extensions of the Fourier sine and cosine series for the function $f(x)=x,~x$ to the interval $x\in[-1,4]$.}
\label{extension}       
\end{figure}

As an example, consider our cosine and sine series developed for $x\in[0,1]$ in Example \ref{sincoscompare}. Recall, the two functions that we found for the unit interval were

\begin{align}
&\textrm{cosine function}&&\nonumber\\
&&    x &=\frac{1}{2}+ \sum\limits_{n=1}^{n\rightarrow \infty} \frac{2(-1+(-1)^n)}{n^2 \pi^2} \cos(n \pi x)\\
 &\textrm{sine function}&&\nonumber\\ 
 &&x &=\sum\limits_{n=1}^{n\rightarrow \infty} \frac{-2(-1)^n)}{n \pi} \sin(n \pi x)
\end{align}
Now, when we plot them over $x\in[-1,4]$ (Fig.~\ref{extension}), we find that each of the functions repeats what is represented on $x\in[-1,1]$.  This makes sense-- each of the basis functions ($\sin(\pi x),~\sin(2\pi x), \ldots$) are defined for any value of $x$, and by construction, each of the trigonometric functions have a period of $2/n$. Thus, the lowest frequency functions involved are $\sin(\pi x)$ and $\cos(\pi x)$, and both of these repeat exactly with a period equal to 2.  All other functions repeat even more frequently, but each repeats $n$ times over an interval of 2. Therefore, the function constructed from the sum of these has no choice but to also repeat with a period of 2.  Comparing the sine and cosine series for $f(x)=x$ in Fig.~\ref{extension}, it is clear that this is in fact the case.  This example also underscores again the differences between the odd (sine) series expansions, and the even (cosine) series expansions. 

\section{Introduction to Convergence of Fourier Series}

In calculus, the issue of the convergence arises in the study of infinite series.  The determination of whether or not a particular series will converge in some sense is a central problem in applied mathematics.  

It would not be an understatement to say that the study of the convergence properties of Fourier series actually ushered in a new era of analysis in mathematics.  When Fourier series were first introduced, Fourier claimed without proof (starting around 1807, when his initial results were presented) that the method would work for \emph{any} function.  It should be noted, however, that even the notion of what properly constituted a \emph{function} and the definition of \emph{integration} were not well-defined at the time.  In that context, one can better understand Fourier's overly-enthusiastic proclamations.  

Over the ensuing time since Fourier proposed his theory, it has gradually been realized that the theory is much more robust than might have initially been envisioned.  While the details of convergence can involve very technical mathematical concepts, it is not incorrect to state that many of the functions that arise naturally in physics and engineering converge in some useful and intuitive sense.  

The earliest serious study of the convergence of Fourier series was done by P. Dirichlet, who in 1829 proposed the first theorem regarding the convergence of Fourier series.  While by modern standards, this theorem is substantially more restrictive than technically necessary, it is a useful touchpoint for understanding the convergence of Fourier series.  The conditions proposed by Dirichlet are \emph{sufficient} to guarantee \emph{at least point-wise} convergence, but are stricter than is technically necessary.  In other words, there are series for functions that do not meet these conditions, but the series still converge in some defined sense.  Regardless, the Dirichlet conditions cover many of the important functions that would arise in more applied problems.  While a more detailed discussion of convergence properties of Fourier series will be delayed, we state the Dirichlet conditions as follows.

\begin{theorem}[Convergence of Fourier Series: Dirichlet conditions.]  The Fourier series for a function $f$ converges pointwise on an interval $I$ if the following conditions are met.
\begin{enumerate}
    \item The function is \emph{bounded} over the domain $I$.  This is sometimes stated more rigorously by the idea that the function f is \emph{absolutely integrable}, which means
    \[ \int |f(x)| \,dx \, < \infty\]
    \item The function is piecewise $C^1$ continuous.  This means, by definition, that there exist a \emph{finite number} of points where the right and left derivatives are not equal (i.e., a finite number of discontinuities in the derivative).  
    \item The function $f$ has only a \emph{finite} number of maxima or minima in the interval.  In other words, the function can not oscillate infinitely fast anywhere in the interval.
\end{enumerate}
\end{theorem}
As a note, the finite number of maxima or minima prohibits functions such as  $f(x)=\tfrac{1}{\sin(1/x)}$ from meeting the conditions because such functions oscillate infinitely fast as $x\rightarrow 0$.  If the the three conditions above are met, then the function $f$ can be said to have a Fourier series that converges everywhere in the domain, $I$.  The series converges pointwise to its actual value $f(x)$ except at points of discontinuity.   At points of discontinuity, $x_d$ the series converges to the average value on the two sides of the discontinuity, i.e., 

\[ f(x_d) =\frac{1}{2}[f(x_d^+)+f(x_d^-)] \]
Here, $f(x_d^+)$ represents the value of the function at $x_d$ approaching from the right, and $f(x_d^-)$ represents the value of the function at $x_d$ approaching from the left.


\section{$^\star$Appropriate bases make a difference}\indexme{basis functions!comparison of}

This section is a bit of a tangent, but one that is worthwhile.  Now that we have a handle on Fourier series, we are in a position to discuss the idea that some bases are more appropriate than others for expressing particular functions.  The point of this short section is not intended to be particularly rigorous, but, rather, to give students an opportunity to \emph{think about} how different kinds of series representations might be better or worse for particular applications.

We have learned about two kinds of expansions.  First, we discussed the power series of the form

\begin{equation}
    f(x)=f(x_0)+(x-x_0)f'(x_0) +\tfrac{1}{2!} (x-x_0)^2 f''(x_0)+\ldots
\end{equation}
Now, even though we do not always think about the power series as being an expansion in polynomials, it actually is!  Suppose that we take the expansion around $x_0=0$.  It becomes pretty clear that the resulting series looks like an expansion in \emph{monomials} (the variable $x$ raised to some power is a monomial).  

\begin{equation}
    f(x)=f(0)+xf'(0) +\tfrac{1}{2!} x^2f''(0)+\ldots
\end{equation}
Here, we can think of the functions $1$, $x$, $x^2$, etc., as the basis functions for the expansion, and $f(0)$, $f'(0)$, $f''(0)$, etc., are the coefficients of the expansion.

\begin{figure}[t]
\sidecaption[t]
\centering
\includegraphics[scale=0.6]{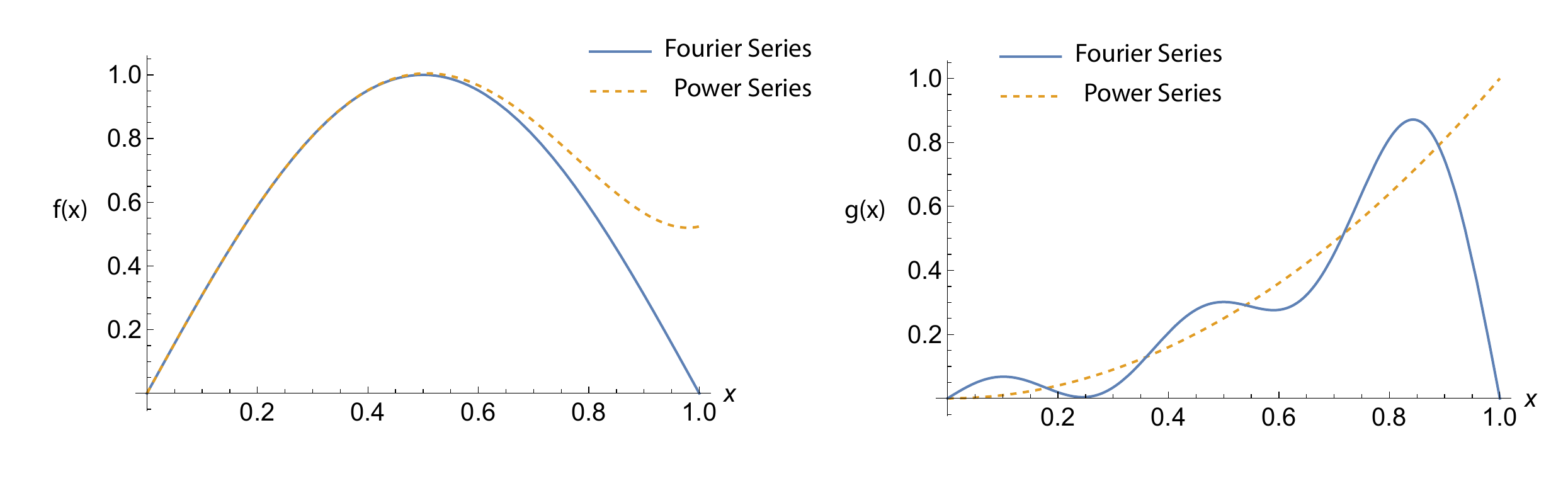}
\caption{Comparison of the series expansions for $f(x)=\sin(\pi x)$ and $g(x)=x^2$ by the Fourier series and power series expansions.}
\label{fig:compare1}       
\end{figure}

Now, lets think about the series expansion of two different functions on $x\in[0,1]$:  $f(x)=\sin(\pi x)$ and $g(x)=x^2$.  In the material in this chapter, we have learned how to compute the Fourier series for each of these functions.  And we already know how to expand functions as power series.  Here, we will just list the results of the expansions for comparison.  The first seven terms in both expansions for the functions $f(x)$ and $g(x)$ are given below.

\begin{align}
f(x)=\sin(\pi x):&&&\nonumber\\
Fourier& &f(x) &=1\cdot \sin(\pi x)+ 0 + 0 + 0 + 0+0+\ldots\\
Power&  & f(x) &= 0 + \pi x +0 + \frac{\pi^3 x^3}{6}+0 \frac{\pi^5 x^5}{120} + 0+\ldots\\
g(x)=x^2:~~~~~~~~&&&\nonumber\\
Fourier& & g(x) &=\frac{\left(2 \pi ^2-8\right) }{\pi ^3}\sin (\pi  x)-\frac{1}{\pi
   }\sin (2 \pi  x)+\frac{\left(18 \pi ^2-8\right) }{27 \pi ^3}\sin (3 \pi  x)-\frac{1}{2
   \pi }\sin (4 \pi  x)\nonumber \\
   &&&+\frac{\left(50 \pi ^2-8\right) }{125 \pi ^3}\sin (5 \pi  x)-\frac{1}{3 \pi }\sin (6 \pi 
   x)+\ldots\\
Power& &  g(x) &= 0+0+x^2+0+0+0 +\ldots
\end{align}
Plots of these series appear in Fig.~\ref{fig:compare1}.  For the function $f(x)$, it is clear that the Fourier series is very efficient at capturing the function (in this case, it is exact), and this requires only one non-zero term.  However, the power series expansion contains significant deviations near $x=1$ even with six terms.  The converse can be see for the function $g(x)$.  For that function, the power series gives an exact result, whereas the Fourier series gives an approximation that is not particularly good when only the first six terms are used. 

So, it becomes clear that different kinds of basis functions can give dramatically different results.  When one is attempting to\emph{fit} an expansion to a function, there are choices among expansions.  While we have examined only two here, there are many other kinds of basis functions that can be used to approximate a target function.  When one would like a good and efficient approximation, exploring the opportunities for choice of basis functions to use can make a difference!

\newpage
\section*{Problems}
\subsection*{Practice Problems}
\begin{enumerate}
\item We first encountered the idea of expanding a function in terms of a set of basis functions when we examined the Taylor series for analytic functions (Section \ref{approxtheory}).  In that material, it was illustrated that a Taylor series truncated at $n=N$ was called a Taylor polynomial of order N.  In such expansions, we have a sum of terms proportional to the following basis set of polynomials

\begin{equation*}
    \mathscr{p} =\{1, x, x^2, x^3, x^4, \ldots\}
\end{equation*}
where we use the notation $p_0(x) = 1$, $p_1(x) = x$,$p_2(x)=x^2$, etc.
Thus, a Taylor series can be represented by
\begin{align*}
f(x)&= \sum_{n=0}^\infty d_n p_n(x) \\
\intertext{or, equivalently}
    f(x) &= \sum_{n=0}^\infty d_n x^n
\end{align*}
where $d_n$ represents the appropriate coefficient of term $n$ of the Taylor series, and $p_n$ are the basis functions.

The finite Taylor polynomials of order $N$ are given by 

\begin{equation*}
    f_N(x) = d_0\cdot 1 + d_1 x + d_2 x^2 +d_3 x^4 +\ldots + d_N x^N
\end{equation*}
It is clear from this representation that a Taylor polynomial is just a linear combination of the polynomials in the set $\mathscr{P}$.  Although the polynomials themselves are not linear (except, of course, for the cases $N=0$ and $N=1$), their sum is in fact a linear combination of the basis functions.

These polynomials can be very useful approximations for functions that otherwise have no easy analytical expression.  For example, the function

\begin{equation}
    f(x) = \frac{1}{1-x}
\end{equation}
converges on the interval $-1 < x < 1$ (i.e., its radius of convergence, $R$, is $|x|<1$
When doing analysis, encountering functions like this can create serious difficulties because they are not linear. However, for $x$ small enough, one can use the first few terms of a Taylor series to very accurately represent the function.  This gives us a good way to introduce the idea of generating an approximation to an analytic function using a polynomial (extracted from its Taylor series).  

For this problem, do the following to illustrate how a linear combination of basis functions can generate an approximation to a function that is analytic in some domain.  Do the following.

\begin{enumerate}
    \item Determine the general Taylor series for this function.  Note: This function has a radius of convergence $R<1$, and it does not converge at the boundary $R=1$.  In other words, this series converges for $|x|<1$.
    
    \item Compute the first five terms of this series explicitly (i.e., work out the Taylor polynomial for it).  Using Mathematica (or whatever software is convenient for you) plot the series approximation $f_5(x)$ and the original function $f(x)=1/(1-x)$ over the interval $x\in[-1, 1]$; recall, $f_5(x)$ represents the sum of the series for $n=0$ to $n=5$ (i.e., the first six terms of the series).  
    
    Compute the fractional error by
    \begin{equation*}
        error(x) = [f_5(x)-f(x)]/f(x)
    \end{equation*}
    and plot this on the same plot with $f_5(x)$ and $f(x)$.
    \item For $x$ small enough, the first two terms in the Taylor series will be accurate, i.e., the Taylor series can be given by 
    \begin{equation*}
        f_1(x) = 1+x
    \end{equation*}
    Determine the maximum value for $x_{max}$ such that the $f_1(x)$ approximation for $f(x)$ is less than 1\% within the domain $0 \le x \le x_{max}$.  Note-- you do not need to find an exact value of $x_{max}$ such that the error is 1\%!  Just find a reasonable value for $x_{max}$ where the maximum error is in the range of $0.08$ to $1.2$\%.  Hint: Try plotting the difference $[f_1(x)-f(x)]$ to find the needed value of $x_{max}$.
\end{enumerate}
\item Much like the fact that a polynomial is its own Taylor series, any trigonometric function in terms of the Fourier basis functions is its own Fourier series.  To see this, consider the following function
    
    \begin{equation*}
        f(x) = 6 \sin(5 \pi x), ~~x\in[0,1]
    \end{equation*}
    Using the definitions we have learned for the Fourier series, show that this function is its own Fourier series expansion.
   \saveenumerate
\end{enumerate}
  
\noindent For the following problems, find the \emph{Fourier sine series} on the interval $x\in[0,1]$ for the function indicated.  You do not have to re-derive the formula for the coefficients $B_n$ (although you are certainly welcome to if that works better for you).   Plot the function using 50 terms.  Also plot the associated spectrum for the first 50 terms.

\setlength{\columnsep}{2cm}
\begin{multicols}{2}
\begin{enumerate}[topsep=8pt,itemsep=4pt,partopsep=4pt, parsep=4pt]
\restoreenumerate
    \item $f(x)=1$
    \item $f(x)=x$ \\(you will have to use integration by parts)
    \item $f(x)=x^2$ \\(you will have to use integration by parts twice here)
    \item  $f(x)=\sqrt{x}$\\
    (A symbolic mathematics program can be used to compute the integral defining the coefficients $B_n$)
    \item  $f(x)=\sin(x)$ \\
    (A symbolic mathematics program can be used to compute the integral defining the coefficients $B_n$)
    \item  $f(x)=\cos(x)$\\
    (A symbolic mathematics program can be used to compute the integral defining the coefficients $B_n$)
     \item $f(x)=\left\{\begin{array}{l}1 \textrm{ for } x\leq\frac{1}{2}\\0\ \textrm{ for }\frac{1}{2}< x\leq1\end{array}\right.$
    \item $f(x)=\left\{\begin{array}{l} 0 \textrm{ for } x\leq\frac{1}{2}\\1\ \textrm{ for }\frac{1}{2}< x \leq 1 \end{array}\right.$
    \item $f(x)=\left\{\begin{array}{l}x \textrm{ for } x\leq\frac{1}{2}\\0\ \textrm{ for }\frac{1}{2}< x\leq1\end{array}\right.$
     \item $f(x)=\left\{\begin{array}{l}1-2x \textrm{ for } x\leq\frac{1}{2}\\0\ \textrm{ for }\frac{1}{2}< x\leq1\end{array}\right.$
     \item $f(x)=\left\{\begin{array}{l}0 \textrm{ for } x<\frac{1}{3}\\1\ \textrm{ for }\frac{1}{3}\leq x\leq \frac{2}{3}\\
     0 \textrm{ for }\frac{2}{3}< x\leq1\end{array}\right.$
        \saveenumerate
\end{enumerate}
\end{multicols}

\noindent For the following problems, find the \emph{Fourier cosine series} on the interval $x\in[0,1]$ for the function indicated.  You do not have to re-derive the formula for the coefficients $B_n$ (although you are certainly welcome to if that works better for you).   Plot the function using 50 terms.  Also plot the associated spectrum for the first 50 terms.

\setlength{\columnsep}{2cm}
\begin{multicols}{2}
\begin{enumerate}[topsep=8pt,itemsep=4pt,partopsep=4pt, parsep=4pt]
\restoreenumerate
    \item $f(x)=1$
    \item $f(x)=x$ \\(you will have to use integration by parts)
    \item $f(x)=x^2$ \\(you will have to use integration by parts twice here)
    \item  $f(x)=\sqrt{x}$\\
    (A symbolic mathematics program can be used to compute the integral defining the coefficients $A_n$)
    \item  $f(x)=\sin(x)$ \\
    (A symbolic mathematics program can be used to compute the integral defining the coefficients $A_n$)
    \item  $f(x)=\cos(x)$\\
    (A symbolic mathematics program can be used to compute the integral defining the coefficients $A_n$)
     \item $f(x)=\left\{\begin{array}{l}1 \textrm{ for } x\leq\frac{1}{2}\\0\ \textrm{ for }\frac{1}{2}< x \leq 1\end{array}\right.$
    \item $f(x)=\left\{\begin{array}{l} 0 \textrm{ for } x\leq\frac{1}{2}\\1\ \textrm{ for }\frac{1}{2}< x \leq 1 \end{array}\right.$
    \item $f(x)=\left\{\begin{array}{l} x \textrm{ for } x\leq\frac{1}{2}\\0\ \textrm{ for }\frac{1}{2} < x \leq 1\end{array}\right.$
     \item $f(x)=\left\{\begin{array}{l}1-2x \textrm{ for } x\leq\frac{1}{2}\\0\ \textrm{ for }\frac{1}{2}< x \leq 1\end{array}\right.$
     \item $f(x)=\left\{\begin{array}{l}0 \textrm{ for } x <\frac{1}{3}\\1\ \textrm{ for }\frac{1}{3} \leq x\leq \frac{2}{3}\\
     0 \textrm{ for }\frac{2}{3}< x \leq 1\end{array}\right.$
        \saveenumerate
\end{enumerate}
\end{multicols}

\noindent For the following problems, find the \emph{Fourier series} for the function on the interval indicated.  Use the full series expansion (involving both sine and cosine components) unless otherwise noted.
\setlength{\columnsep}{2cm}
\begin{multicols}{2}
\begin{enumerate}[topsep=8pt,itemsep=4pt,partopsep=4pt, parsep=4pt]
\restoreenumerate
    \item $f(x)=1$, $x\in[-1,1]$
    \item $f(x) =x$, $x\in[-1,1]$
    \item $f(x)= |x|$, $x\in[-1,1]$
    \item $f(x)=\exp(2x)$, $x\in[-2,2]$.
    \item $f(x)=\exp(2x)$, $x\in[0,2]$; do this as a sine series.
    \item $f(x)=\exp(2x)$, $x\in[0,2]$; do this as a cosine series.
    \item $f(x) = \frac{1}{1+x}$
    \saveenumerate
\end{enumerate}
\end{multicols}

\bigskip
\subsection*{Applied and More Challenging Problems}
\begin{enumerate}
\restoreenumerate

\item {\bf Even and odd extensions.}  The interval  $x\in[0,L]$ is sometimes called the \emph{half interval}, because both the Fourier sine and cosine series on this interval are also defined on the interval $[-L,L]$.  In fact, both the Fourier sine and cosine series are defined on $[-L,L]$, and they are periodic with period $2L$.  For the sine series, the extension to $x\in [-L,L]$ should result in an odd function since sine is an odd function on this interval.  For the cosine series, the extension to $x\in [-L,L]$ should result in an even function since sine is an even function on this interval.\\

For the following problems, develop \emph{both the Fourier sine and cosine series} for the interval $x\in [0,L]$ indicated.  Then, show that both series are defined on $[-L,L]$ by using the the result you have obtained, but plotting it over $x\in[-L,L]$ using $N=50$ terms.  You should find that the sine series is defined by an \emph{odd} extension of the function that appeared for $x\in[0,L]$, and the cosine series is defined by an \emph{even} extension of the function that appeared for $x\in[0,L]$.
\begin{enumerate}
    \item $f(x) = x,~~x\in[0,1]$
    \item $f(x) = x^2~~x\in[0,2]$
    \item $f(x) = \left\{\begin{array}{l}1 \textrm{ for } x\leq\frac{1}{2}\\0\ \textrm{ for }\frac{1}{2}< x \leq 1\end{array}\right.$
    \item     \vspace{-0mm} $ f(x)=
    \begin{cases}
 0 & 0<x<\frac{1}{4} \\
 8 x-2 & \frac{1}{4}\leq x\leq \frac{1}{2} \\
 6-8 x & \frac{1}{2}<x\leq \frac{3}{4}\\
 0 & x>\frac{3}{4} \\
\end{cases}
$

(Note: This implicitly indicates that the function is defined on the interval $~~x\in[0,1]$.)
\end{enumerate}

\item {\bf Periodic extensions.}  
This extension of Fourier sine and cosine series on $x\in[-L,L]$ for subsequent periodic intervals is called \emph{periodic extension}. Because both the sine and cosine function are periodic, with period equal to $2L$, when the series developed for $x\in[-L,L]$ is extended in either the positive or negative direction, the solution repeats because of the periodicity of the underlying trigonometric functions.

For the following problems, develop \emph{both the Fourier sine and cosine series} for the interval $x\in [0,L]$ indicated.  Then, show that the that the series is actually a periodic one that repeats every $2L$ by plotting the series for $x\in [-5L,5L]$. Use $N=50$ terms.
\begin{enumerate}
    \item $f(x) = x,~~x\in[0,1]$
    \item $f(x) = x^2~~x\in[0,2]$
    \item $f(x) = \left\{\begin{array}{l}1 \textrm{ for } x\leq\frac{1}{2}\\0\ \textrm{ for }\frac{1}{2}< x \leq 1\end{array}\right.$
    \item $ f(x)=
    \begin{cases}
 0 & 0<x<\frac{1}{4} \\
8 x-2 & \frac{1}{4}\leq x\leq \frac{1}{2} \\
 6-8 x & \frac{1}{2}<x\leq \frac{3}{4}\\
 0 & x>\frac{3}{4} \\
\end{cases}
$

(Note: This implicitly indicates that the function is defined on the interval $~~x\in[0,1]$.)
\end{enumerate}

\item  Consider the following function, defined piecewise over the interval $0\le x \le 1$ (Fig.~\ref{piecewise}).\label{prob_3_23}
\begin{equation*}
f(x)=
    \begin{cases}
 0 & 0<x<\frac{1}{4} \\
 16 x-4 & \frac{1}{4}\leq x\leq \frac{1}{2} \\
 12-16 x & \frac{1}{2}<x\leq \frac{3}{4}\\
 0 & x>\frac{3}{4} \\
\end{cases}
\end{equation*}
\begin{figure}[b]
\sidecaption[t]
\centering
\includegraphics[scale=0.3]{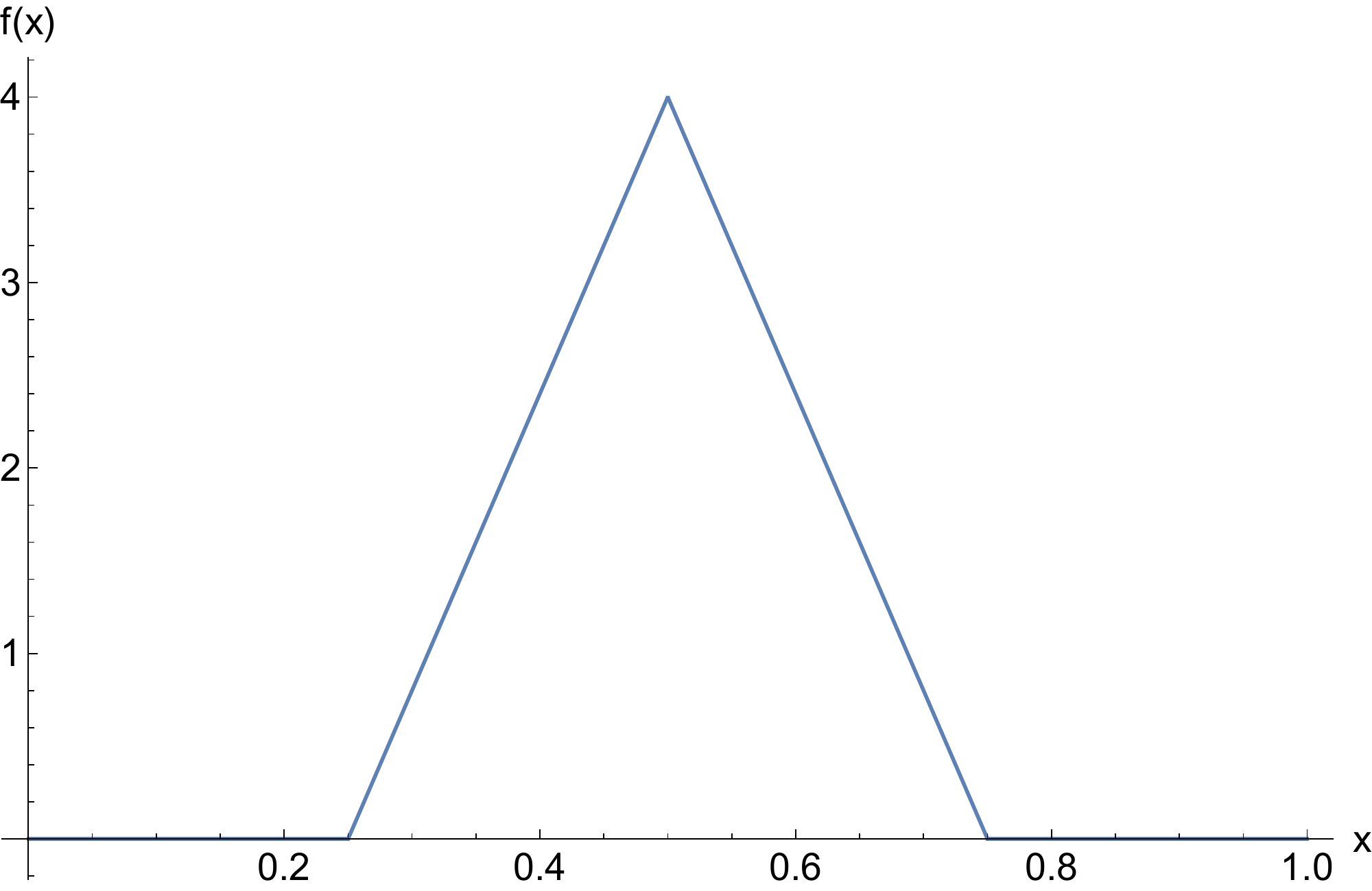}\hspace{2mm}
\includegraphics[scale=0.3]{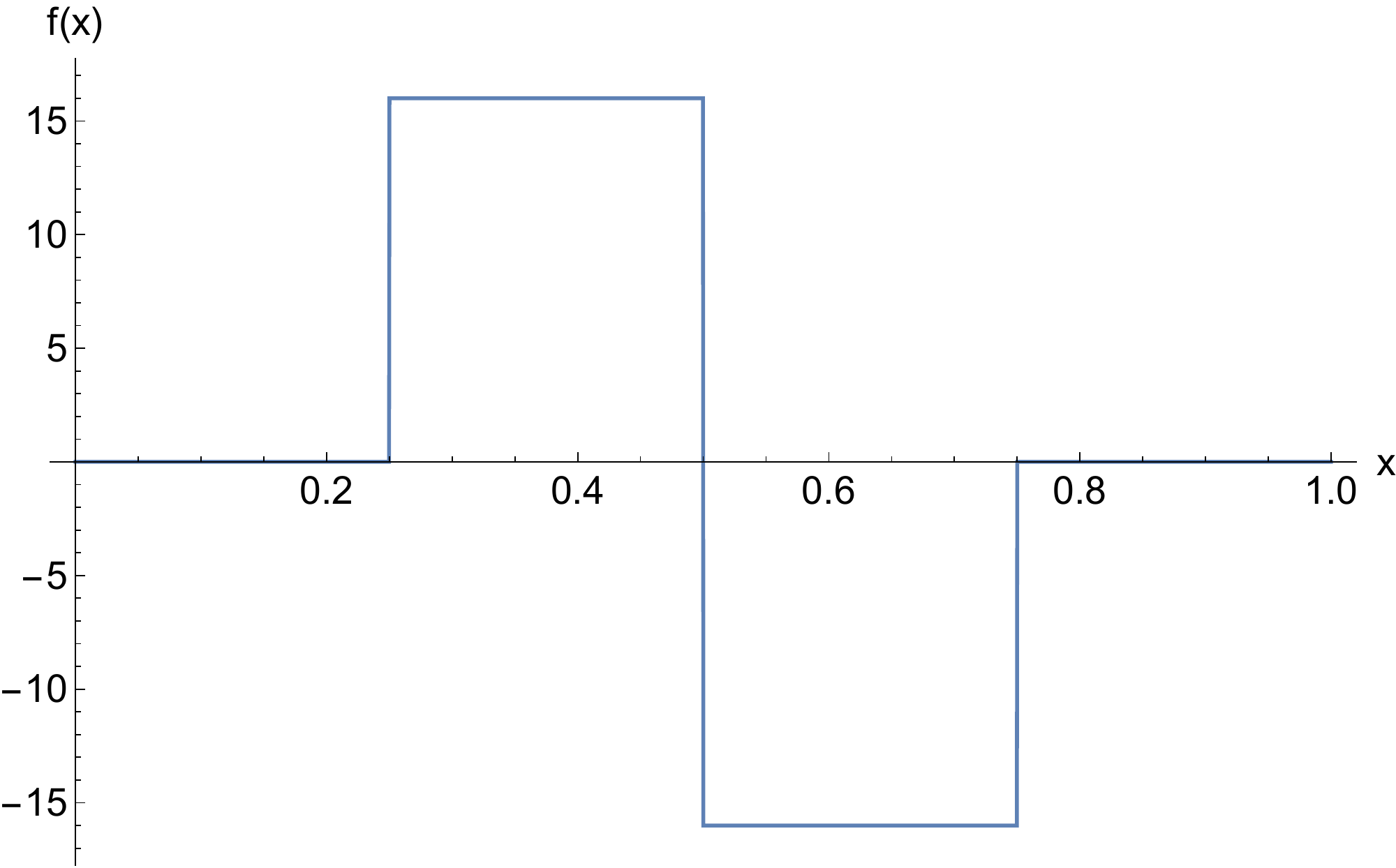}
\caption{The piecewise function defined for problems \ref{prob_3_23} (left) and \ref{prob_3_24} (right).}
\label{piecewise}       
\end{figure}
Find the Fourier sine series for this function.  Plot the function using the first 100 terms in the series.\\

\item Consider the following function, defined piecewise over the interval $0\le x \le 1$ (Fig.~\ref{piecewise}).  Note that this function is the derivative of the piecewise function given in problem \ref{prob_3_23}. \label{prob_3_24}
\begin{equation*}
f(x)=
    \begin{cases}
 0 & 0<x<\frac{1}{4} \\
 16 & \frac{1}{4}\leq x\leq \frac{1}{2} \\
 -16  & \frac{1}{2}<x\leq \frac{3}{4}\\
 0 & x>\frac{3}{4} \\
\end{cases}
\end{equation*}

Find the Fourier sine series for this function.   Plot the function using the first 100 terms in the series.

\item  There are many requirements as to when a Fourier series can be differentiated term-by-term.  While this question is a complex one in general, what we can say is that on the interval $x\in[0,1]$, if $f(0)=f(1)=f'(0)=f'(1)=0$ \emph{and} the function $f$ is piecewise $C^2$, then both the Fourier sine and cosine series for $f(x)$ can be differentiated to obtain the appropriate series for $f'(x)$.  Show that this is the case by doing the following.
\begin{enumerate}
    \item Compute the termwise derivative of the Fourier series for the function given in problem \ref{prob_3_23}.
    \item Plot the function using the first 100 terms in the resulting series.
    \item Compare this sum to the exact solution given by the piecewise function given in \ref{prob_3_24}.
\end{enumerate}

\item For some functions, the evaluation of the sine or cosine series requires extra care and thought because the function $f(x)$ itself is one of the series terms.  For instance, consider the following piecewise function \label{prob_3_25}
\[
f(x)=
\begin{cases}
 \sin (2 \pi  x) & 0<x\le \frac{1}{2} \\
 0 & \frac{1}{2}<x \le 1
\end{cases}
\]
For this function, plotted in Fig.~\ref{piecewise2}, the set of coefficients for the Fourier sine series, $B_n$, contain an apparent singularity (the denominator goes to zero).  Thus, the problem requires some additional handling.  To see one resolution to this problem, complete the following steps.
\begin{enumerate}
    \item Compute the Fourier sine series for this problem.  Verify that the result for the coefficients is
    \[B_n = \frac{4 \sin\left(\tfrac{1}{2} n \pi \right)}{(-4 +n^2)\pi} \] 

Note that for $n=2$, this result presents a problem.  The numerator is zero, and the denominator is zero, leading to a $0/0$ indefinite form.

\item There are a few ways in which this problem can be handled.  While the result for $B_n$ is correct, the indeterminate form arises because we have computed the integral for the general case.  The indeterminate form is avoided if we consider the integral for the $n=2$ case directly, because we can make obvious simplifications that are not obvious if we consider the more general case.  In particular, note that the integration for the $n=2$ case is
\[ \int_{0}^{\tfrac{1}{2}}  \sin(2 \pi x) \sin(2 \pi x) \,dx = \int_{0}^{\tfrac{1}{2}}  \sin^2(2 \pi x)\, dx   \]
Once one combines the two functions into $\sin^2(2 \pi x)$, the indefinite form no longer occurs.  Show that the integration above leads to the result $B_2 = \tfrac{1}{2}$.
\begin{figure}[t]
\sidecaption[t]
\centering
\includegraphics[scale=0.6]{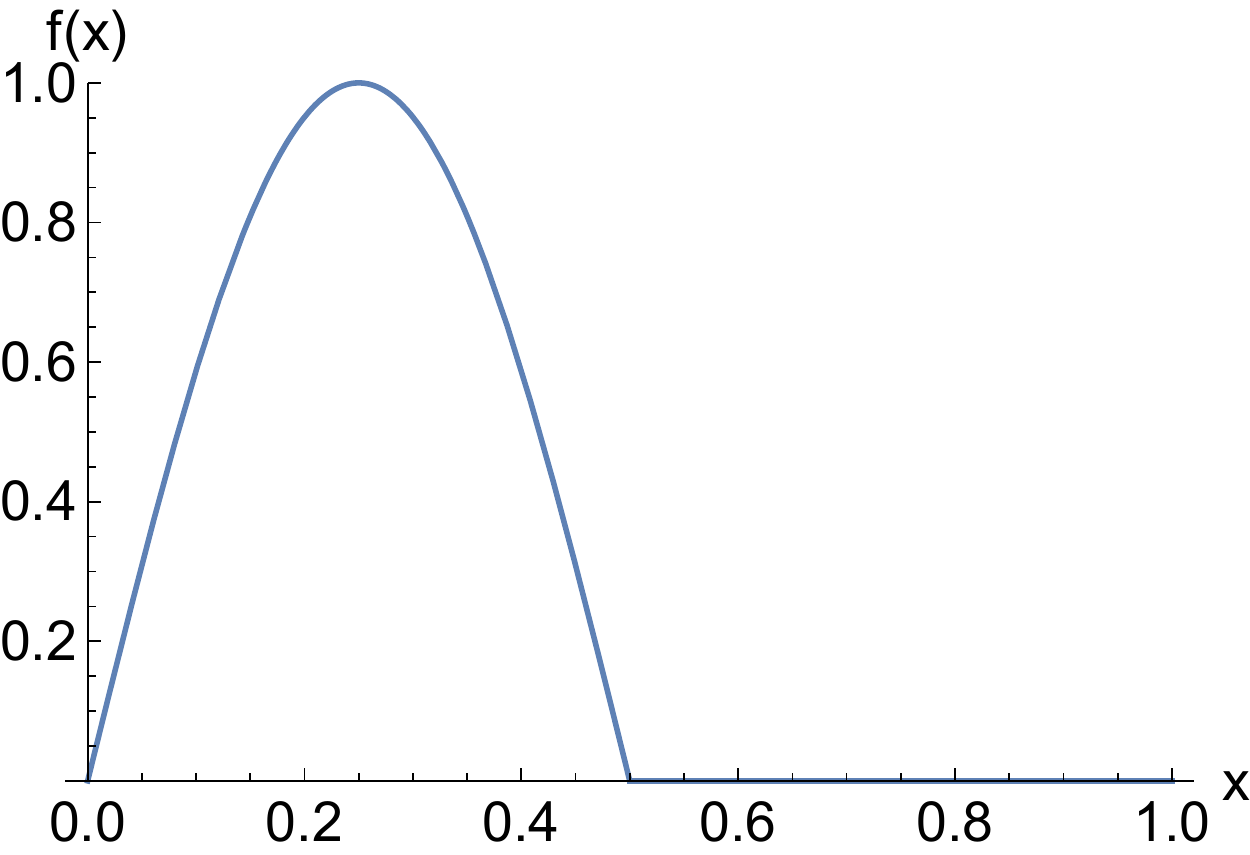}\hspace{2mm}
\caption{The piecewise function defined for problems \ref{prob_3_25}.}
\label{piecewise2}       
\end{figure}

\item There is a second, equivalent way to handle this problem.  Recall from your study of calculus L'H\^opital's rule, which states that for an $0/0$ indefinite form
\[  {\mathop {\lim }\limits_{n \to n_0 }} \frac{f(n)}{g(n)} = {\mathop {\lim }\limits_{n \to n_0 }} \frac{f'(n)}{g'(n)} \]
Even though the function $B_n$ is not a continuous one, we can certainly treat it as though it were.  Show that, treating $n$ as a continuous variable, you can use L'H\^opital's rule on $B_n$ in the form
\[ B_2 = {\mathop {\lim }\limits_{n \to 2 }} \frac{4 \sin\left(\tfrac{1}{2} n \pi \right)}{(-4 +n^2)\pi} \]
to show that $B_2 = \tfrac{1}{2}$.

\item Plot the Fourier sine series using $N=10$ terms, and the piecewise function to compare them.
\end{enumerate}

\item The function $f(x)=\exp(-(a x)^2)$ (this is the classical \emph{Gaussian} function used in statistics and physics) has a curious property.  Its spectrum has very close to the same shape as the function itself (assuming that $a$ is large enough such that $f(x)<<1$ as $x\rightarrow L$, $x)\in[-L,L]$.   For this problem, assume that the interval of interest is $x\in [-1,1]$.

\begin{enumerate}
\item To see that the spectrum looks like the original function, compute the Fourier cosine series for this function for the case. To do so, you will need to note the following integral

\[ \int_{0}^1 \exp(-a^2 x^2) \cos(n \pi x) \, dx = 
\begin{cases}
\tfrac{1}{2a}\sqrt{\pi}\, \textrm{erf}(a), &n=0\\
\tfrac{1}{a} \sqrt{\pi}\exp\left(-\tfrac{\pi^2 n^2}{4 a^2}  \right) \textrm{erf}(a),&n>0
\end{cases}
\]
\item Assume that $a=3$, so that $f(x) = \exp(-(3 x)^2)$. 
On two separate plots, plot (i) The original function $f(x)$ on $x\in[0,1]$, and (ii) $B_n=B(n)$.  For the latter function, we can using the following plotting trick to extend the function to an interpolated version that is smoother looking.  The trick is as follows: instead of plotting $B(n)$, plot the coefficient $B(n/100)$.  For this example, plot $B(n/100)$ over the interval $n\in[0,600]$.  This should give an interpolated plot that looks smooth (although the values of $B_n$ appear at every division of 100 units on the horizontal axis). \\ 

In Mathematica, one command to plot this latter function is given by 

\small{
\begin{verbatim}
    DiscretePlot[B[n/100], {n, 1, 600}, ExtentSize -> Full, 
    ColorFunction -> "Rainbow", PlotRange -> All, AxesLabel -> {"n", "B(n)"}, 
    ImageSize -> Large, PlotRange -> All]
\end{verbatim}
}
\end{enumerate}

\item Proving that that the sine and cosine functions are \emph{orthogonal} is something that was not shown in the text, but it is not difficult to illustrate.  To do so, one needs to use the following trigonometric identities

\begin{align*}
    \cos(A)\cos(B) &= \frac{1}{2}\left[\cos(A-B)+\cos(A+B)\right] \\
    \sin(A)\sin(B) &= \frac{1}{2}\left[\sin(A-B)+\sin(A+B)\right] \\
\end{align*}
Using these identities, show that the following are true.
\begin{align*}
    \int_{x=0}^{x=1} \cos\left( n \pi x\right) \cos\left( m \pi x\right) \, dx &= \begin{cases}
    \frac{1}{2} \textrm{ for } n=m \\
    0 \textrm{ for } n \ne m\end{cases}\\
    n=1,2,3,\ldots \\
    \int_{x=0}^{x=1} \sin\left( n \pi x\right) \sin\left( m \pi x\right) \, dx &= \begin{cases}
    \frac{1}{2} \textrm{ for } n=m \\
    0 \textrm{ for } n \ne m\end{cases}\\
    n=1,2,3,\ldots \\
\end{align*}
NOTE: Please work this out carefully by hand.

\end{enumerate}
%

%

\abstract*{This is the abstract for chapter 00}

\begin{savequote}[0.55\linewidth]
``It was a sort of act of faith with us that any equations
which describe fundamental laws of Nature must have
great mathematical beauty in them. It was a very
profitable religion to hold and can be considered as the
basis of much of our success."

\qauthor{Paul Dirac on the laws of quantum mechanics, where the delta function arises as an important feature}
\end{savequote}
\def\CHAP{chapter04_the_delta_function}
\theoremstyle{definition}

\chapter{The Step and Delta Functions}\label{deltachap}
%
\def\CHAP {chapter04_the_delta_function}

The concept of what constitutes a \emph{function} was briefly covered in Chapter \ref{chaprev}.  While the concept itself seems simple enough, this is in part because generally we have been generally been exposed to the concept throughout our education.  However, the word \emph{function} itself did not even exist until the late 1600's, when the mathematicians Gottfried Leibniz (of calculus fame) and Johann Bernoulli (a Swiss mathematician, who's son Daniel was famous for the Bernoulli principle of fluid mechanics, and for the gamma function, discussed below) began to develop the concept of function more formally.  Even into the early 1800's, the concept of \emph{function} was still thought by many to apply only to \emph{analytic} functions (see Chapter \ref{chaprev}).  

Many mathematicians contributed to the generalization of the concept of what defined a function.  However, it is not an overstatement to say that it was the theory of Joseph Fourier who prompted modern efforts for defining functions.  Through the use of Fourier series (see Chapter \ref{FS_1}), Joseph Fourier was able to \emph{construct} functions that behaved very much unlike the functions that mathematicians were used to contemplating.  For example, Fourier series existed which defined functions with discontinuities in the value (i.e., jumps) within the function's domain.  More alarmingly, these discontinuous functions were constructed entirely from infinite series of \emph{continuous} functions, which seemed to present a conceptual paradox.

While we will not study function theory as a separate topic, we will in this section cover two important functions that come up frequently in applications: (1) The step (or Heaviside) function, (2) and the delta function.  The first of these was briefly discussed in the chapter on Fourier series.  The second function mentioned, the delta function, is technically not a function; it is more correctly called a \emph{generalized function} or a \emph{distribution}. However, the terminology \emph{delta function} is so thoroughly ingrained in science and engineering that we will adopt this (slightly incorrect) terminology.   

You may have encountered the concept of the delta function previously.  The delta function is a good example of a mathematical concept whose justification was very much inspired by the fields of engineering and physics.  In fact, the function is often referred to as the Dirac delta function in honor of the physicist Paul Dirac (8 August 1902 – 20 October 1984), who used the delta function as an important component of his description of quantum mechanics.   Another common application of the delta function is to represent point charges in electrostatics.  In engineering mechanics, the delta function is used to represent point loads, as described above for the case of the force acting on a beam.  As mentioned above, the delta function was routinely in used applications long before it was understood mathematically.  In the 1950's, the more general theory of generalized functions (which includes the delta function) was finally established by a mathematician named Laurent Schwartz.  While the theory extended the notion of what constituted a function to new mathematical constructs, we will not pursue that course here.  Instead, we will focus specifically on the delta function and the related step function.  The ideas will be developed using primarily the tools of calculus.

\section{Terminology}

\begin{itemize}

\item {\bf Generalized function.}  A notion that expanded the definition of functions.  In particular, a generalized function is a mathematical object that may not meet the definition of a regular (classical) function, but can be described by its action on other functions via integration.  All regular functions are generalized functions, but some generalized functions are not regular functions.   The delta function is the most well known example of a generalized function.  \indexme{generalized function}  \\

\item {\bf Delta function.}  A generalized function that physically represents a concentrated source over a very small (relative to other scales of the problem) time or space interval.  Mathematically, the delta function posed many difficulties; ultimately, these unusual mathematical objects were given a sound mathematical framework generally known as the theory of \emph{generalized functions}.  The word generalized was chosen because these new objects had function-like utility, but were not functions in any mathematically conventional sense.    While the delta function does not behave like any known classical function (e.g., it is non-zero at only one point, but its integral is nonetheless equal to 1), it does arise from applied, physically-based considerations.  As an example, point charges in the theory of electrostatics are representable by delta functions.  Similarly, point forces in the mechanics of materials can be represented by delta functions. \indexme{delta function}\indexme{generalized function!delta function}\\

\item {\bf Delta sequence.}  A function, $f$ that is indexed by an integer, $n$, such that the sequence of functions $(f_1, f_2, f_3, \ldots)$ becomes increasingly peaked and narrow as $n\rightarrow\infty$.  An example is the conventional Gaussian function, written in the form
\[ G(x; a, n) = \frac{1}{\sqrt{\pi (a/n)^2}} \exp\left(- \frac{x^2}{(a/ n)^2}\right) \]
This function is illustrated for several values of $n$ in Fig.~\ref{gaussians}.  It is easy to see in this figure that as $n$ increases, the function becomes more peaked and narrower.  In the limit, this sequence of functions becomes a \emph{delta function} (hence its name).

\begin{figure}[!ht]
\sidecaption[t]
\centering
\includegraphics[scale=.48]{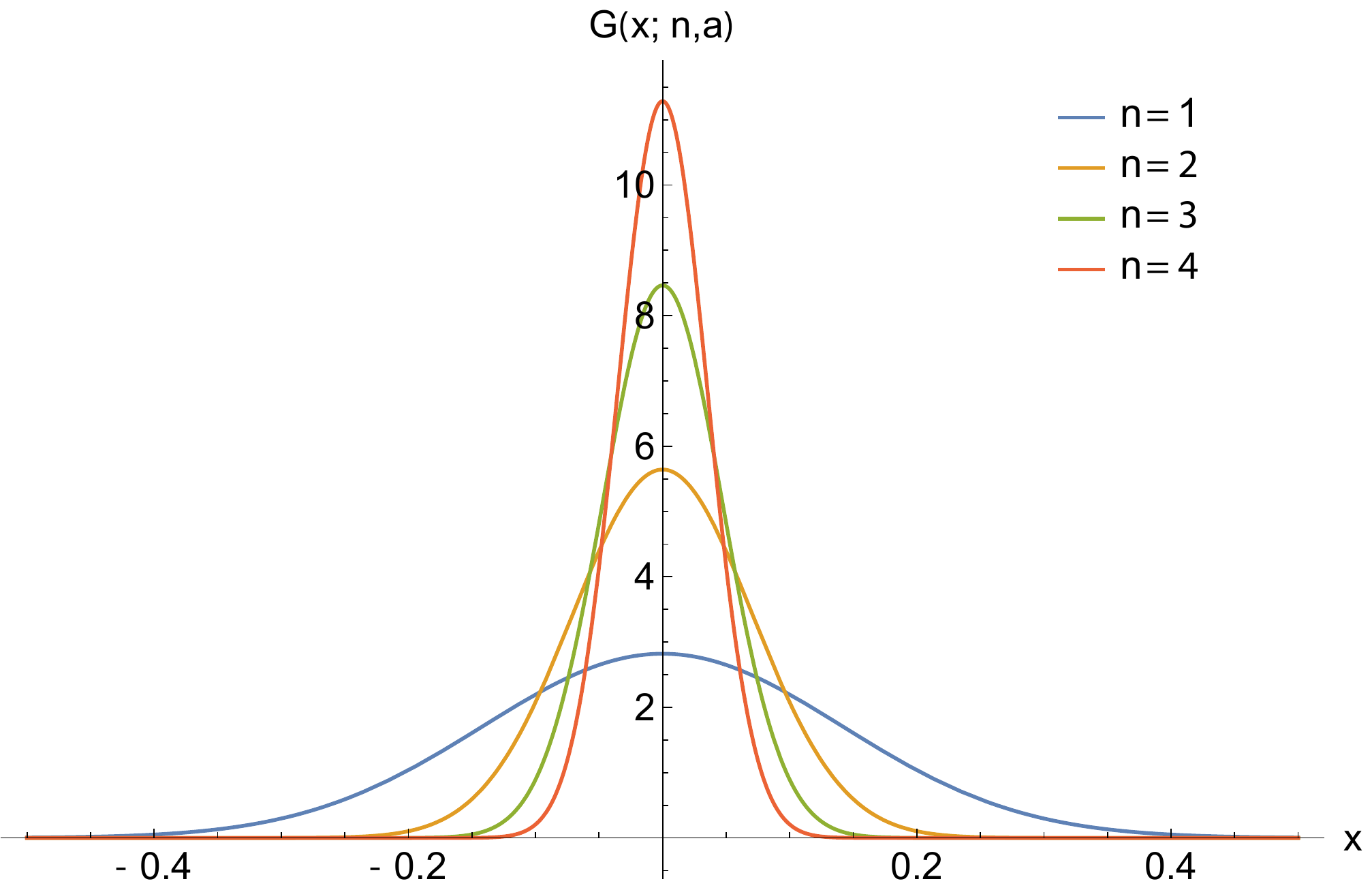}
\caption{{\bf A delta sequence.}  The Gaussian $G(x;\tfrac{1}{5},n)$ with standard deviation $\sigma = \tfrac{1}{5} n,~n\in\mathbb{N}$ is an example of a \emph{delta sequence}.  As $n$ increases, the function becomes increasingly peaked and narrow, even though the area stays constant at 1. }
\label{gaussians}       
\end{figure}

\item {\bf Step function or Heaviside function.}  Step functions are functions that are discontinuous at a single point, $x_0$.  The function is zero for $x<x_0$.  At the point of discontinuity, the step function jumps from 0 to the value 1, and remains at that value for $x>x_0$.  At $x=x_0$, the function is technically not defined; however, a better way of thinking about this is that the function can take on any value, $c$ between 0 and 1 at $x_0$.  As such, the step function is actually always an entire \emph{equivalence class} (see \S \ref{equivalence_class}) of functions, where the value chosen for the function at $x_0$ defines the particular function in the class.\\

\item {\bf Generalized derivative.}  The derivative of a generalized function. This concept will be defined by appealing to integration by parts.  The ideas is to give concrete \emph{meaning} to what it means to, for example, take the derivative of the step function.  Classically, the step function has no derivative defined at the point of discontinuity (i.e., the function jumps a finite distance over an interval of size zero, so its classical derivative is infinite.) The theory of generalized functions gives meaning to the ``derivative" of the step function.  While we will not explore the theory of generalized functions \emph{per se}, we will nonetheless adopt an approach that gives an intuitive (but mathematically rigorous) notion of the derivative of the step function. \indexme{generalized function! generalized derivative}\indexme{derivative!generalized}\\

\item {\bf Function compact support.} A \emph{function of compact support} is (for a single independent variable) a function that is nonzero only on some closed interval.  While it is hard to conjure up images of such function, there are many examples that are even $C^\infty$ (differentiable any number of times).  One routinely adopted example is
\[ f(x) = \exp(1/x)\exp[1/(1-x)] \]

This function is plotted in Fig.~\ref{compact_func}.  Functions of compact support that are also $C^\infty$ are sometimes called \emph{bump functions} or \emph{mollifying functions}.\\
\begin{figure}[t]
\sidecaption[t]
\centering
\includegraphics[scale=.48]{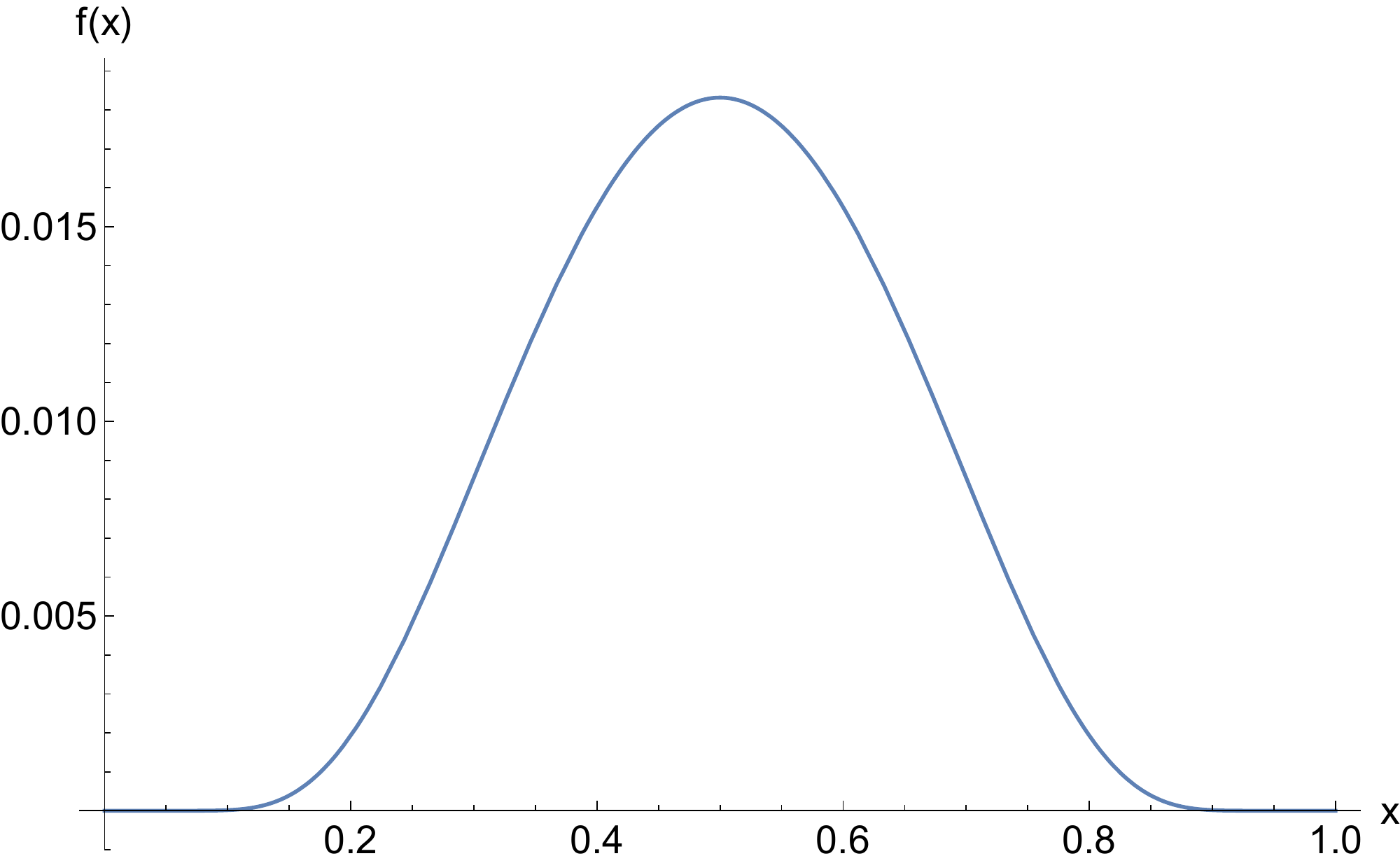}
\caption{{\bf A compact function.}  This function is non-zero only for $0\le x \le 1$, and is $C^\infty$.  Functions like this that are both compact and $C^\infty$ form an important class of functions called \emph{bump} functions.}
\label{compact_func}       
\end{figure}

\item {\bf Gamma function, $\boldsymbol{\Gamma}$.}  \indexme{function!gamma function} The gamma function is an analytic function on the real line. One interpretation of the gamma function is that it extends the notion of the factorial to the real numbers.  The history of the gamma function reads like a ``who's who" of classical mathematics, with important contributions regarding the applications and propertries of the function being established by Leonhard Euler (1707--1783), Adrien-Marie Legendre (1752--1833), Carl Friedrich Gauss (1777--1855), Joseph Liouville (1809--1882), Karl Weierstrass (1815--1897), and Charles Hermite (1822--1901), among others. 

\end{itemize}

\section{The Step and Delta Functions: The Basic Idea}

\begin{figure}[t]
\sidecaption[t]
\centering
\includegraphics[scale=.75]{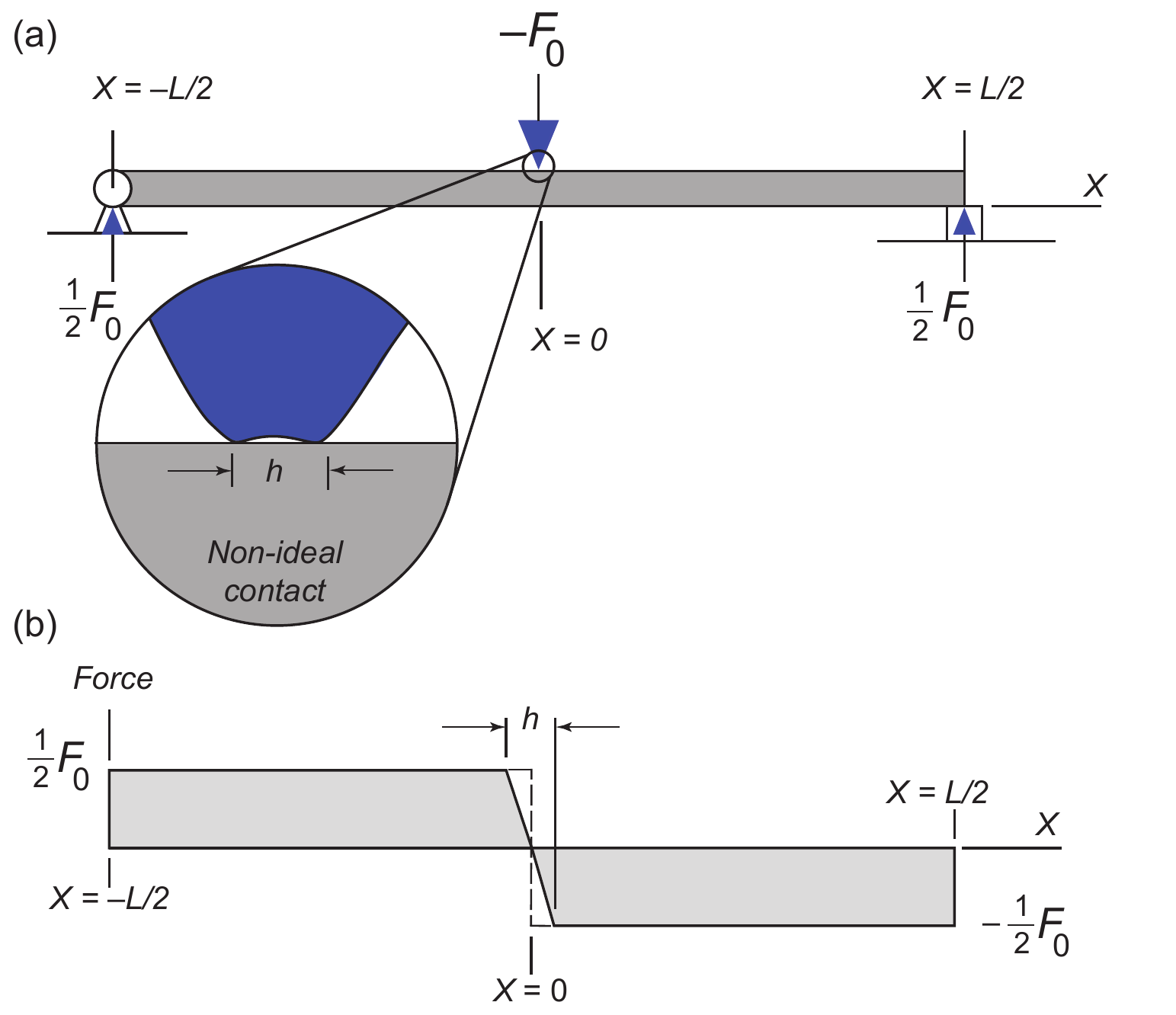}
\caption{{\bf Idealization of non-ideal system.}  Mathematical modeling always requires a level of abstraction from ``reality".  The trick in modeling is to make abstractions that are good approximations to the real problem, and also lead to mathematically tractable results.  The use of the delta function is an example of an abstraction (i.e., an representation that is not exact, but generates the correct behavior anyway) for many physical quantities of interest.  In the example in part (a) of this illustration, the actual force applied to the beam is not ideal.  This means that it does not occur at a single point, but is a complex function of $x$ that, in general, we would not know explicitly (see expanded view).  In part (b) of this figure, we show what the force diagram (also known as a \emph{shear diagram}) would look like for a force that is distributed over a small area of size $h$.  Here, we have approximated the force as being uniformly distributed over the distance $h$, even though the actual distribution may be more complex.}
\label{beam}       
\end{figure}

In this section, the purpose is to build the basic intuition about the delta function and the step function.   The intuition about the delta function is actually fairly simple to visualize.  Physically, the delta function represents some quantity (e.g., a force, a concentration, heat energy) that is \emph{idealized} as if it were concentrated at a single point.  There is a purpose for this model of a physical system. For example, consider the force of a sharp wedge bearing down on a rigid beam.  In Fig.~ \ref{beam}, a force acting on a beam is illustrated at two levels of resolution.  From afar, the idea of representing the force as if it were applied at a single point seems reasonable.  However, we have to consider what we actually mean by applying a force to a single point.  Physically, this is impossible; all physical contacts must occur through an interaction of finite areas.  Physically, a \emph{point} does not exist; a point is a strictly mathematical abstraction.  Thus, if we zoom in on the region where the force is applied, it is clearly not occurring at a single point; in fact, the application of the force happens over a small area, and, in general, such interactions will occur in a complex way that we will not be able to know in detail.

\subsection{Mathematical Modeling of a Point Source}

While it is true that \emph{physically} we cannot apply a force to a single point, there are instances where \emph{mathematically} we would like to model a quantity that is applied to a small area (relative to some other dimension; in this case, the length of the beam) as if it were applied to a point.  The reason that we might want to make such a model are as follows.

\begin{enumerate}
\item Generally, we would not know the \emph{actual} distribution of forces concentrated over a small region.  Thus, we would not know how to represent the distribution of forces exactly.  

\item Even if we could represent the distribution of forces exactly, there is the sense that, because of the difference in scales, $h \ll L$, it seems unlikely that the particular distribution of forces over $h$ will be relevant.

\item There is some hope that a mathematical model that represents the force as if it were applied to a single point would offer mathematical simplifications than would a more detailed model. 
\end{enumerate}
At this juncture, then, we need to address the question ``\emph{is it possible to generate a sensible model of the force as if it applied to a single point?}"  This is our goal for this chapter.\\

\begin{figure}[t]
\sidecaption[t]
\centering
\includegraphics[scale=.75]{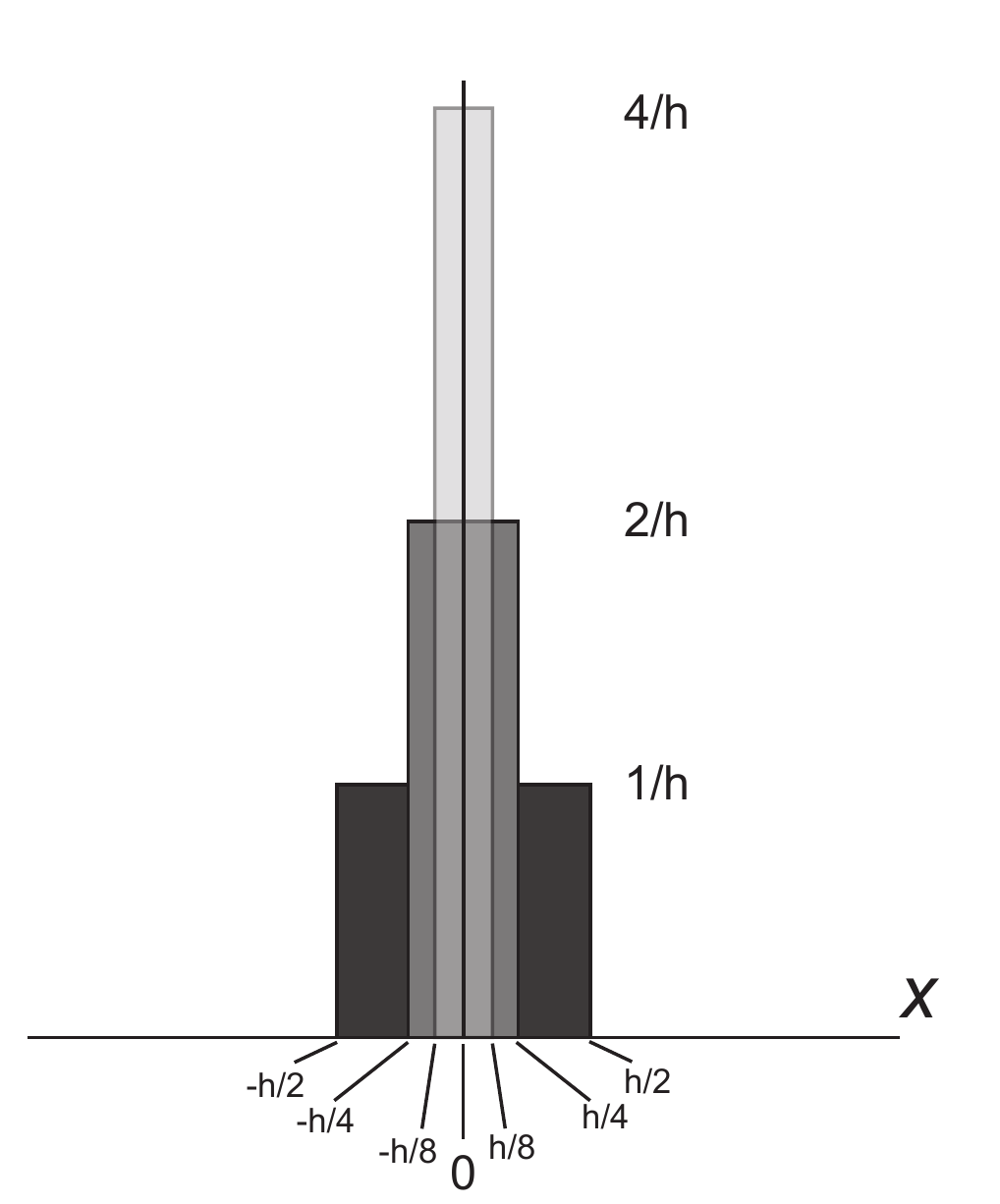}
\caption{{\bf Approximations to the delta function.}  A \emph{sequence} of rectangle functions, indexed by $n$.  The area under the rectangle is fixed at a value of one.  As $n$ increases, the rectangle gets thinner and taller, but the area always remains one. }
\label{rect_delt}       
\end{figure}

To start our thinking about the problem, consider the set of functions illustrated in Fig.~\ref{rect_delt}.  Each of the rectangles shown has the same area-- an area of unity (i.e., $area=1$).  Clearly, for each integer increase in $n$, the width of the rectangle decreases by a factor of 2, and the height increases by a factor of 2.  Thus, as $n$ increases, the rectangle becomes increasingly peaked and narrow.  This function is given explicitly by

\begin{equation}
F_0{\delta _n}(x) = 
\begin{cases}
F_0/(h/n) & \mbox{for } - \frac{h}{2 n} < x < \frac{h}{2 n} \\ 
  0 &  otherwise 
  \end{cases}
  \label{delta_force}
\end{equation}
Some additional interpretation is helpful here.  Consider the case $n=1$ and $n=2$.  These case represents the force of the beam, were it to be spread uniformly over the interval of size $h$ and $h/2$ respectively. 

\begin{align}
F_0{\delta _1}(x)& = 
\begin{cases}
F_0 \frac{1}{h} & \mbox{for } - \frac{h}{2} < x < \frac{h}{2} \\ 
  0 &  otherwise 
  \end{cases}\\
  \label{delta_force2}
F_0 {\delta _2}(x)& = 
\begin{cases}
F_0 \frac{2}{h} & \mbox{for } - \frac{h}{4} < x < \frac{h}{4} \\ 
  0 &  otherwise 
  \end{cases}
\end{align}
Increasing values of $n$ affect the force density in obvious ways.  Note additionally, the constant $F_0$ is distinct from the delta function, $\delta_n$; its role is to modify the \emph{magnitude} of the delta function.  

As a concrete example of these functions, suppose we were to integrate, $F_0 \delta_1$, over the interval $-h/2 \le x \le h/2$. We would find

\begin{equation}
  \int_{\xi=-h/2}^{\xi=h/2} {F_0}\frac{1}{h} \delta_1(\xi) \,d\xi =  F_0 \int_{\xi=-h/2}^{\xi=h/2} \frac{1}{h} \delta_1(\xi) \,d\xi = F_0 \frac{1}{h}h=F_0
\end{equation}

We now have a model in which the area of interaction can be adjusted to be as larger or small as we like, while the total force applied to the beam remains constant.  With this representation, we can compute the force on the beam at any location $x$ as follows

\begin{equation}
    F(x) = \frac{F_0}{2}-\int_{\xi=-L/2}^{\xi=x} F_0 \delta_n(\xi) \,d\xi
    \label{shear}
\end{equation}
This result generates a distribution of shear forces in the beam as illustrated in Fig.~\ref{beam}.  We can interpret the Eq.~\eqref{shear} describing this distribution of shear forces as follows.  On the left-hand side, the support provides an upward force of $F_0/2$.  As we approach the region near $x=0$, the force per unit width is spread uniformly (by our model) over the small distance $h$. The force expression given by Eq.~\ref{shear} represents this by an increase in the total downward force as $x$ increases in this region, such that integrating the force density across the width $-h/(2n)< x < h/(2n)$, one recovers the downward force $F_0$.  Finally, as $x$ becomes greater than $h/(2n)$, then the total downward force applied is the amount $F_0$.  A little thought on the physics of the system will hopefully match your intuition.  

\subsection{The Path Forward: Limiting Behavior}

Now, although our model for the force distribution is much simpler than what the actual distribution might be, it is still has not achieved the notion of being applied at a single point.  In fact, for those who have studied beams in physics or a course in statics, the force diagram illustrated in Fig.~\ref{beam} don't look like those that you have seen before.  Usually, the transition at $x=0$ is treated \emph{as if the force were in fact applied at a single point}.  In that case, the force diagram near $x=0$ would resemble the curve given by the dashed line rather than the solid line, which contains an observable finite interval over which the force increases.  Interestingly, we can make our representation for the force as close to the idealized dashed line as we like, simply by making the value of $h$ smaller (or, looking at Eq.~\eqref{delta_force}, making $n$ larger).  This begins to give us some insight as to the nature of the delta function.  It appears that we might think of this function as being a \emph{limit} of the function described above as the width of that function tends toward zero.

The problem now is making mathematical sense out of an idea that is motivated by a physical model.  Quite frequently, mathematics and science have interacted to reinforce one another.  This was certainly the case for the delta function, where physicists and engineers used the delta function long before it had a solid mathematical theory underpinning it.  Because the delta function had such obvious intuition as the limit of a process that had a sensible interpretation, it was widely adopted by scientists and engineers.  Correspondingly, mathematicians realized that the delta function was not a function in any conventional sense of the word.  There was no theory that described the mathematics of such a function, and its use lead to quite a number of mathematical conundrums.  However, as mentioned previously, a full understanding of the mathematics of this problem was finally developed in the 1950s.  Now, such functions are used with confidence in both the physical sciences and in mathematics.

\section{A Construction for the Step and the Delta Functions}\indexme{function!step} \indexme{step function}\indexme{function!Heaviside}\indexme{generalized function!delta function}\indexme{delta function}\indexme{Heaviside function}

The description given above was a relatively intuitive presentation of the delta function.  Now, we will firm things up a bit by illustrating how to \emph{construct} delta functions from appropriate sequences.  Once we see how this is done and under what conditions it is possible to do so, we can simply adopt the delta function as a ``shorthand" notation for a more complex process.  Fortunately, the complexities of the process do not need to be repeated for every problem; once the delta function is understood, it can be used by adopting an intuitively appealing set of rules.

\begin{figure}[t]
\sidecaption[t]
\centering
\includegraphics[scale=.8]{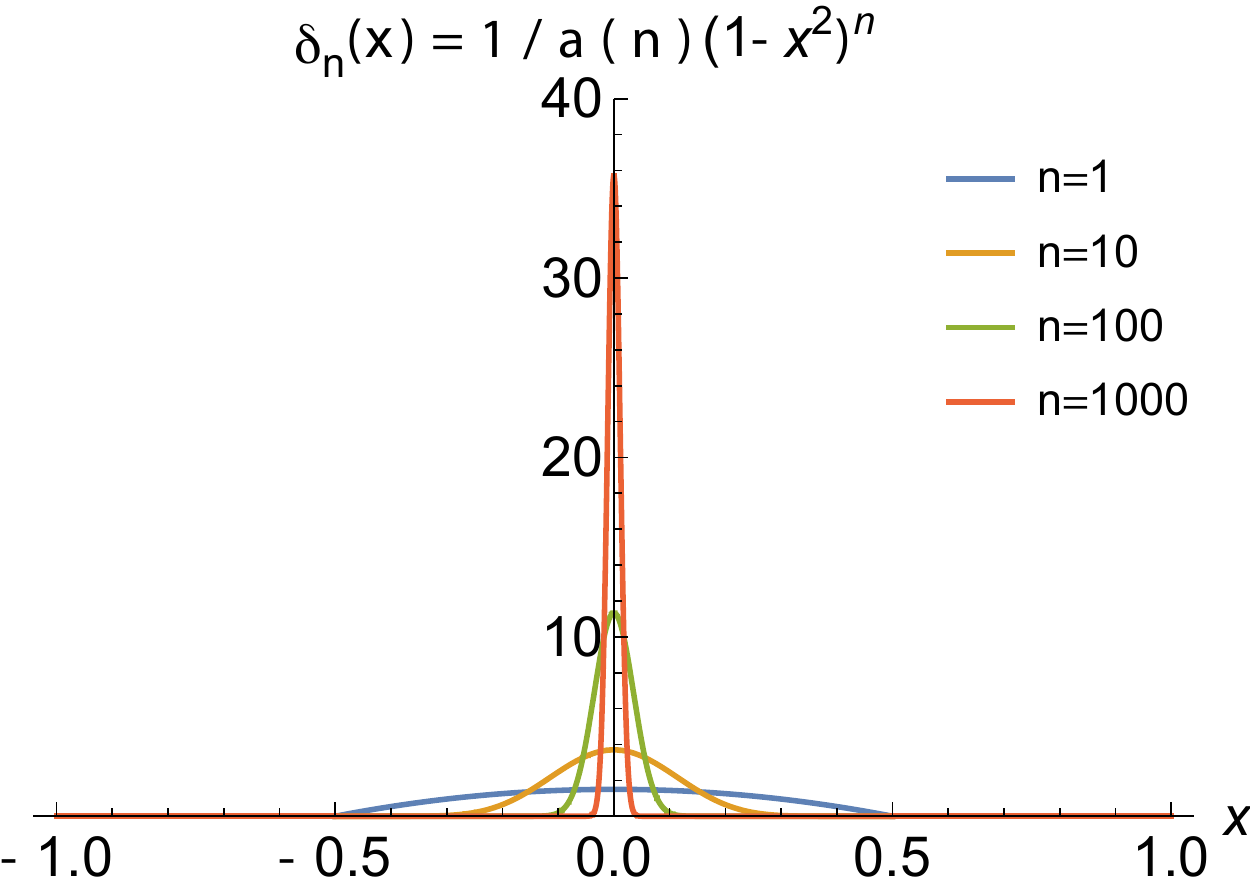}
\caption{{\bf Approximations to the delta function.} A \emph{sequence} of functions that approaches the delta function in the limit.  These sequences are constructed using the function $\delta_n(x)=1/a(n)(1-x^2)^n$, which is normalized by $a(n)$ to always have unit area.  As the value of $n$ increases, the function rapidly becomes narrower and more peaked.}
\label{deltaa}       
\end{figure}

\begin{figure}[t]
\sidecaption[t]
\centering
\includegraphics[scale=.8]{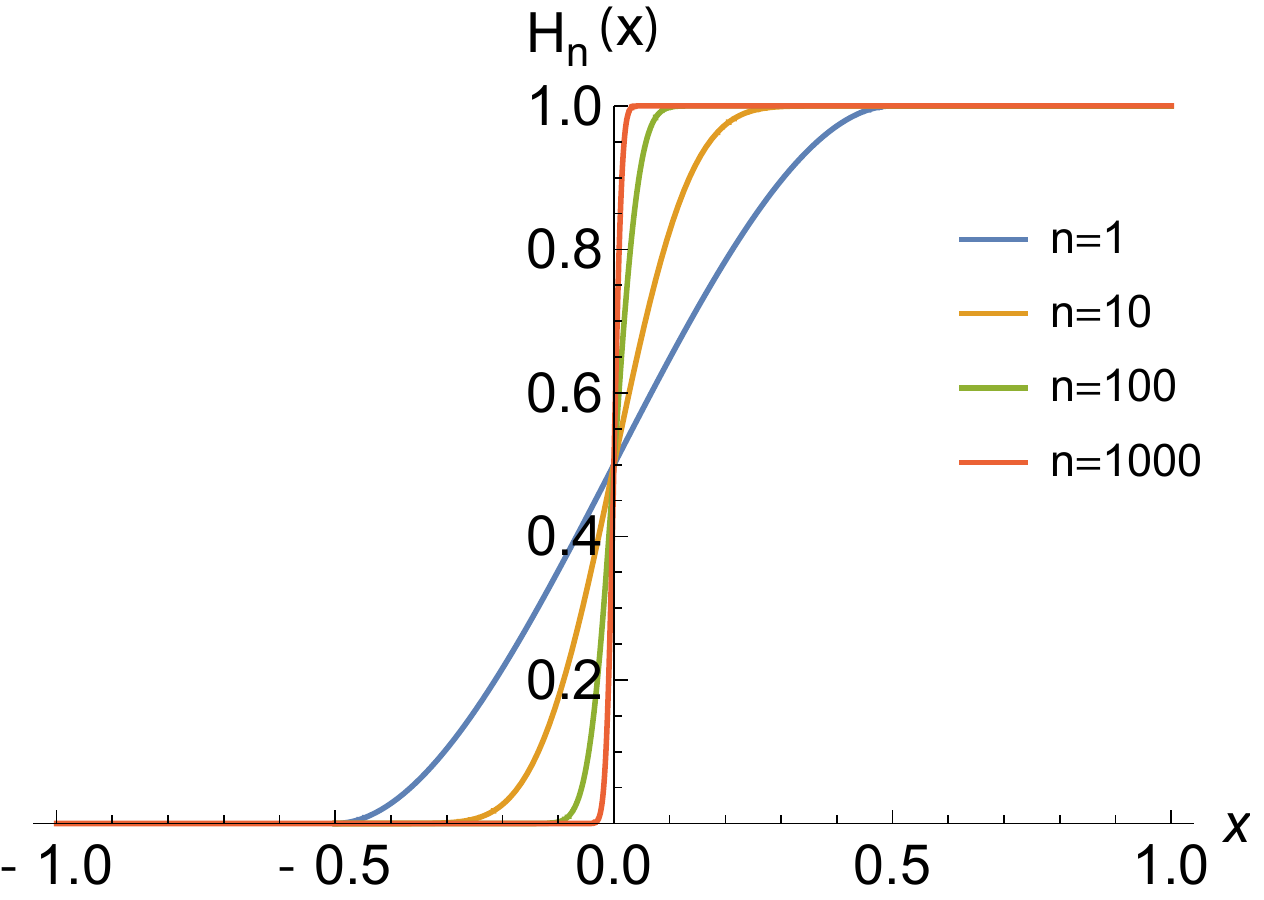}
\caption{{\bf Approximations to the step function.}  A \emph{sequence} of functions that approaches the \emph{integral} of the delta function in the limit.  These sequences are constructed using the function $H_n(x)=1/a(n)$$\int_{-1/2}^x (1-\xi^2)^n d\xi$.  As the value of $n$ increases, the function rapidly becomes abrupt and step-like.  The step function is also called the Heaviside function in honor of the electrical engineer Oliver Heaviside who popularized its use.}
\label{stepa}       
\end{figure}

To start the exploration of the delta and step functions, we are going to revise the previous definition to use functions that are compact and differentiable on an interval.   Consider the following function and its integral.  The function $\delta_n$ is technically compact on $x\in(-1/2,1/2)$.  The exponent $n$ is an integer greater than zero.  The function below is normalized by its area, $a(n)$, so that the integral of the function is always unity (i.e., area of the normalized function is $area=1$) regardless of the value of $n$.  The function is also differentiable $n+1$ times.

\begin{align}
    \delta_n(x) =
    \begin{cases}
    0, &x<-\tfrac{1}{2}\\
     \frac{1}{a(n)}
    (1-x^2)^n, & -\tfrac{1}{2}\le x \le \tfrac{1}{2}\\
     0, &x>\tfrac{1}{2}
    \end{cases}
    \label{deltafun}
\end{align}
\begin{align}
    H_n(x) &= \int_{\xi=-1/2}^{\xi=x} \,\,(1-\xi^2)^n\, d\xi
    \label{stepfun}
\end{align}
\begin{align}
    a(n) &=H_n(1/2)= \int_{\xi=-1/2}^{\xi=1/2} \,\,(1-\xi^2)^n\, d\xi
    \label{areafun}
\end{align}
The expression for $H_n(x)$ is the integral of $\delta_n(x)$ to the value $x$; note that this means that $H_n(1/2)$ is the area under the curve of the function $\delta_n(x)$. This means that $\delta_n(x)$ is \emph{normalized} so that the area under the curve is always one.   The mathematical expression for $H_n$ is quite complicated, and not technically necessary at this point for the discussion.  To help with the visualization of this function and its integral, these are plotted in Figs.~\ref{deltaa} and \ref{stepa}  for various values of $n$.

A few things are apparent from these plots.  First, the plot of the function $\delta_n(x)$ shows that $\delta_n(x)$ is a compact function (i.e., its domain is a finite one, $x\in[-1/2, 1/2]$) that becomes increasingly narrow and tall as $n$ becomes large enough.  It appears that, as we increase $n$, the maximum of $\delta_n(x)$ increases in value, and the values of the function not immediately near the origin become arbitrarily small, comparatively.  In fact, we can determine the \emph{second spatial moment about the vertical axis} (the variance) of the function $\delta_n(x)$ by integration by parts. This is a very messy computation, so the details will not be shown.  However, the result is the following. 

\begin{align}
\sigma(n)& = \frac{1}{F_\Phi(1/2,n)} \int_{x=-1/2}^{x=1/2} x^2 (1-x^2)^n\, dx \nonumber\\
&=\frac{1}{4 (2 n+3)}
\end{align}
Knowing the characteristic width of this function, we can see that as $n$ increases without bound, we have the somewhat strange result
\begin{align}
\mathop {\lim }\limits_{n \to \infty } \sigma(n) &= 0
\end{align}

With a little bit more work, we can find the area under the curve (also known as the \emph{zeroth spatial moment}) for $\delta_n(x)$

\begin{equation}
    a(n) =  \int_{\xi=-1/2}^{\xi=1/2} \,\,(1-\xi^2)^n\, d\xi= \frac{\sqrt{\pi } \Gamma (n+1)}{2 \Gamma \left(n+\frac{3}{2}\right)}
\end{equation}
The function defined by the symbol $\Gamma$ is known as the \emph{gamma function}; it is described in more detail in the gray text box for Example \ref{gamma} (The gamma function).

\begin{svgraybox}
\begin{example}[The gamma function]\label{gamma}

We are all familiar with the idea of the factorial, $n\!$ that applies to any integer.  This function is plotted for the first few integers in Fig.~\ref{gammafact}.  One might wonder if there is concept similar to the factorial, but for any real number.  It turns out that there is, and this function is called the gamma function.  The gamma function is one of the most interesting and important functions in applied mathematics, and it was derived originally by Daniel Bernoulli.  However, for our purposes, we can think of it an the idea of extending the factorial.  We will just state the definition of the gamma function, and begin using it.  Late on, we may return to this definition to extend its use to other kinds of problems.  The gamma function is defined by 

\begin{equation}
    \Gamma(n) = \int_{x=0}^{x=\infty} x^{n-1}e^{-x}\, dx\,\,, \text{ where} n\in \mathbb{R}
\end{equation}
%
{
\centering\fbox{\includegraphics[scale=.4]{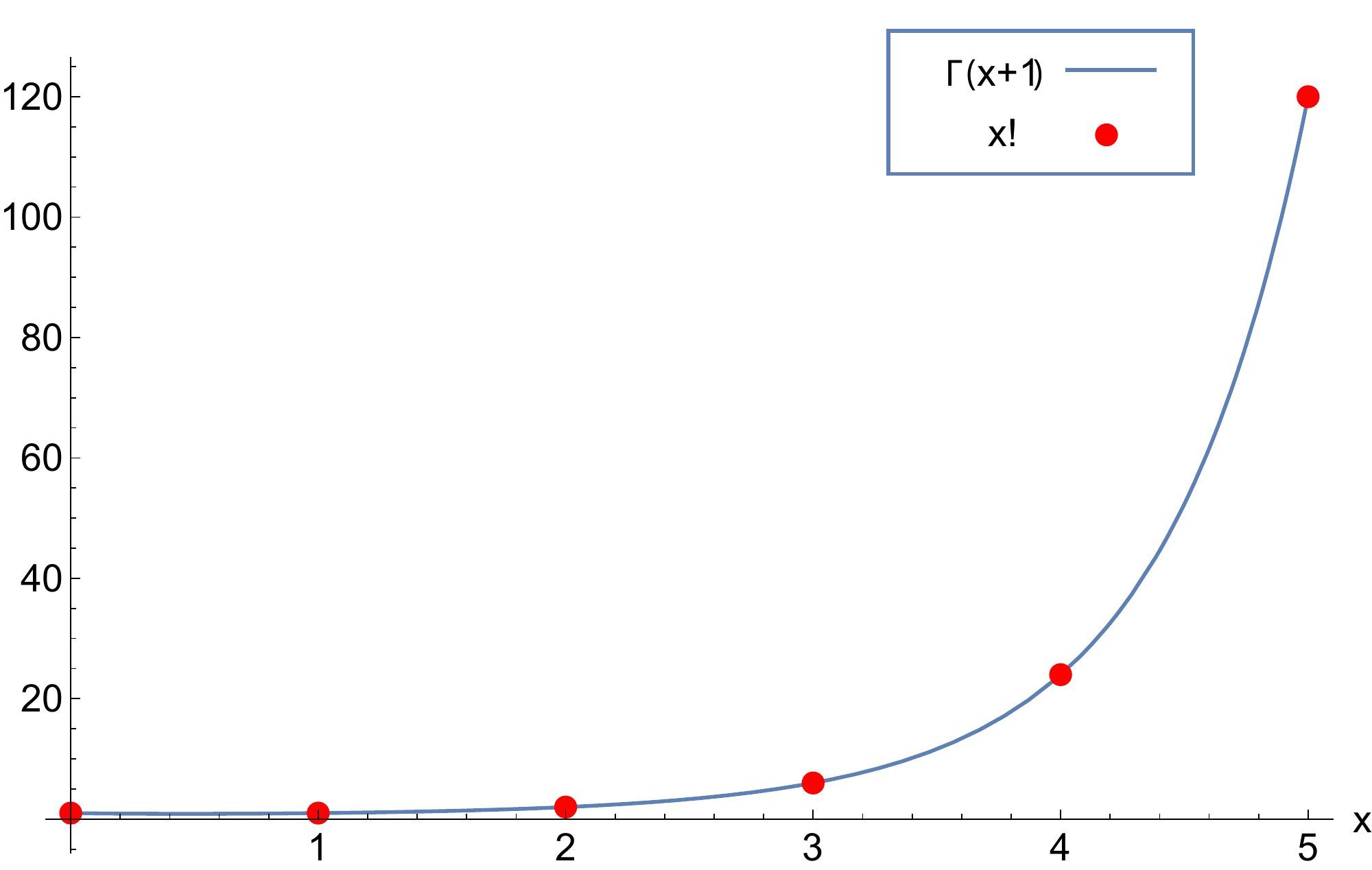}}
\vspace{2mm}
\captionof{figure}{The gamma function plotted over continuous $x$; for comparison, the factorial computed at the integer values of $x$ is shown in red.}
\label{gammafact}  
}
Here, we are using $n$ to indicate any positive real number, rather than only an integer.  While it is not conventional to use $n$ to indicate a real number, it helps us remember that when $n$ is equal to an integer value, it is equal to the factorial.

Now that we have defined the gamma function, we can turn to the area under the curve for the function $f(x)=(1-x^2)^n$.  It turns out that

\begin{equation} 
a(n) = \int_{x=-\tfrac{1}{2}}^{x=\tfrac{1}{2}}(1-x^2)^n \, dx = 
\frac{\sqrt{\pi } \Gamma (n+1)}{2 \Gamma \left(n+\frac{3}{2}\right)}
\end{equation}
which gives us the necessary result for the area under the curve, $a(n)$.
\\

\noindent{\it Some additional features of the gamma function}\\

We have the following identity (which will not be proven here)

\begin{equation}
    \Gamma(n) = (n-1)!
\end{equation}
There are many interesting properties for the gamma function.  Possibly the most relevant one is as follows

\begin{align*}
    \gamma(1) &= 1 \\
    \gamma(n+1)&=n\gamma(n)\\
    n!&=\gamma(n+1)
\end{align*}
The last of these is the statement that, for integer values, the gamma function  is equivalent to the factorial.  One can derive this by using the identity above, and building up starting from $n=1$.
\end{example}
\end{svgraybox}

So now, let's look at what happens to this function as $n$ becomes large.  In Fig.~\ref{deltaa}, the function is plotted for various values of $n$.  What starts out as a very flat function for $n=1$ rapidly becomes very highly peaked as the exponent $n$ increases.  To make this easier to discuss, define $\delta_n(x) = \frac{1}{a(n)} f_1(x)$ so that

\begin{equation}
    f_1(x)=(1-x^2)^n
\end{equation}
Note also the following

\begin{equation}
    \mathop {\lim }\limits_{n \to \infty } a(n) = \mathop {\lim }\limits_{n \to \infty } \frac{\sqrt{\pi } \Gamma (n+1)}{2 \Gamma \left(n+\frac{3}{2}\right)}= 0
\end{equation}
Thus, as $n$ increases, the area used to normalize the curve also tends toward zero.  

Now, note the following.  The maximum of the function $f_1(x)=(1-x^2)^n$ is at $x=0$, and its value is always $f_1(0)=1$.  However, every other point in the function has a value that is less than 1.  Recall, for a positive number, $m_0$ less than 1, we have

\begin{equation}
    m_0 > m^2_0 > m^3_0 > m^4_0\ldots
\end{equation}
This means that as $n$ grows, all of the values of the function become smaller, except the value at $x=0$, which remains 1 regardless of how large $n$ is.  Thus, once normalized by $a(n)$, the only option for the function is to become more and more peaked as $n$ becomes arbitrarily large, with the value at $x=0$ being equal to $1/a(n)$.  

Recalling that $a(n)$ is tending toward zero as $n$ becomes arbitrarily large, and we find that in the limit our function $f_\delta$ behaves quite unusually in the limit.  It apparently has zero width (as measured by the variance), and infinite height.  And, because of its construction, it also has a total area of 1; in other words

\[
 \mathop {\lim }\limits_{n \to \infty }\left(
\int_{-\infty}^{\infty} f_\delta(x,n) \, dx \right) = 1
\]

While it may seem like a peaked function with zero width could not somehow also have a finite area, note that the width (variance) times height (roughly, the area) of this function is, in the limit,  a $0\cdot \infty$ form.  Recall from Chapter 1, we had the example of the function $g(x) = \tfrac{1}{x}\sin(x)$.  In the limit, as $x\rightarrow 0$, we can show using L'H\^opital's rule that

\begin{equation}
    \mathop {\lim }\limits_{n \to \infty } \frac{\sin(x)}{x} = \mathop {\lim }\limits_{n \to \infty } \frac{\cos(x)}{1} = 1
\end{equation}
Therefore, it is not that unusual that a $0\cdot \infty$ form leads to a result that is neither 0 nor $\infty$.  This is essentially the case that we have for our function $f_\delta$.

It is also interesting to look at the $integral$ of the function $\delta_n(x)$ as $n$ increases.  Although that function is given mathematically by the expression given by Eq.~\eqref{stepfun}, it may be more helpful just to look at Fig.~\ref{stepa} to see how this function behaves.  When the values of $n$ are very small, the function $\delta_n(x)$ is somewhat wide and smooth.  Thus, as we integrate up to values of $x$ in the range $-1/2 < x <1/2$, the result is a smooth ``S" shaped function with a maximum of 1.  The function $F_\Phi$ has a maximum of 1 because the function $f_\delta$ is normalized to always have an area of 1.  

For $F_\Phi$, the function becomes more like a \emph{step} as the value of $n$ becomes arbitrarily large.  This is consistent with what we know for the function $f_\delta$; we know that as $n$ becomes arbitrarily large, the function $f_\delta$ becomes increasingly narrow.  Thus, the parts that contribute most to the integral of that function become concentrated around $x=0$; therefore, the function $F_\Phi$ resembles a sharp step up from zero to 1 at $x=0$ as $n$ becomes large.

\section{Delta Sequences}\indexme{delta function!delta sequence}

\begin{figure}[!ht]
\sidecaption[t]
\centering
\includegraphics[scale=.45]{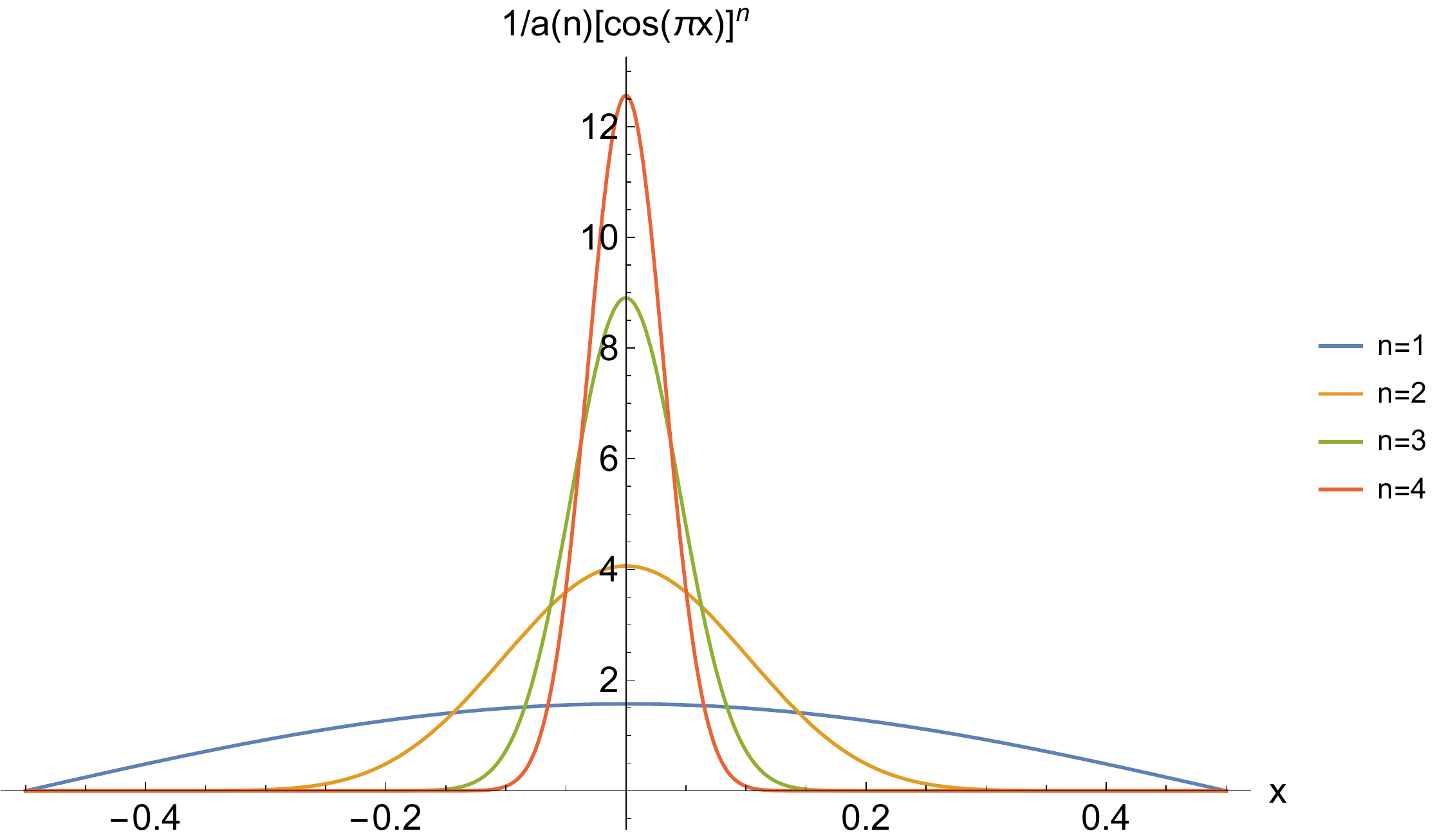}
\caption{{\bf Another delta sequence.}  The cosine function raised to the power $n$ and normalized by its area (as a function of $n$) is a delta sequence. }
\label{cosinesdelta}       
\end{figure}

In the previous section, we developed a \emph{sequence} of functions indexed by the integers, $n$ that became increasingly peaked and narrow as $n$ increased.  These functions were defined by

\[
    f_1(x)=(1-x^2)^n
\]
It turns out that there are any number of sequences of functions that one can make that, in the appropriate limit, the sequence defines a delta function. For example, the function 
\begin{equation}
    f(x)=
    \begin{cases}
    0& \textrm{for }x<-\tfrac{1}{2}\\
    \frac{\sqrt{\pi } \Gamma \left(\frac{n}{2}+1\right)}{\Gamma \left(\frac{n+1}{2}\right)}
    \left[\cos(n \pi x)  \right]^n & \textrm{for } -\tfrac{1}{2} \leq x \leq \tfrac{1}{2}\\
    0& \textrm{for } x>\tfrac{1}{2}
    \end{cases}
\end{equation}
is also a delta sequence (this is plotted in Fig.~\ref{cosinesdelta}).

%
The fact that there are multiple kinds of delta sequences that can lead to the same result indicates that the delta function is not a single entity; it is actually an equivalence class of sequences, any of which can be used equally well to define the (generalized) delta function.  A delta sequence, then, can be defined as follows.

\begin{definition}[delta convergent sequence.]  A delta (convergent) sequence, $\delta_n(x)$, is any sequence of functions such that

\begin{equation}
 \mathop {\lim }\limits_{n \to \infty }\left(
\int_{-\infty}^{\infty} \delta_n(x) f(x) \, dx \right) = f(0) \label{deltadef}
\end{equation}
\end{definition}
There are an infinite number of such functions.  For example, if g(x) is any non-negative function such that 

\begin{equation}
    \int_{-\infty}^{\infty} g(x)\,dx =1
\end{equation}
then rescaling this function in the form
\begin{equation}
    \delta_n(x) = n g(nx), ~~n=1,2,3,\ldots
\end{equation}
is a delta convergent sequence.

\section{Properties of the Delta Function}\indexme{delta function!properties}

Ultimately, the purpose for defining the delta function is because of its unique analytical properties, which make the representation of point sources in physical system very convenient.  There are two primary analytical properties of the delta function that make it convenient in applications.  Let $I$ be either a compact ($\{I:x\in [a,b]\}$) or non-compact ($\{I:x\in [-\infty,\infty]\}$) interval containing the point $x_0$.  Furthermore, we must assume that \emph{the function is at least continuous} at $x_0$ (i.e., there cannot be a jump discontinuity there).  Then, the following two properties are true for the delta function.

\begin{align}
&\textrm{Property 1: Unit integral.} &&\nonumber\\
 &&   &\int_{x\in I} \delta(x)\, dx = 1\\
 &\textrm{Property 2: The sifting property.} &&\nonumber\\
 &&  &\int_{x\in I} f(x)\delta(x-x_0) \, dx = f(x_0)
\end{align}
The first property arises by construction; the delta sequences are always structured such that they have unit area, so the limiting function must also have unit area.  The second of these properties is the so-called \emph{sifting} property of the delta function.  It is just an application of the definition given by Eq.~\ref{deltadef}, but one where the delta function is shifted by the amount $x_0$. The simplicity of the integration of the delta function
is one of its most attractive features. \\

\begin{svgraybox}
\begin{example}[The sifting property]

The sifting property of the delta function is called this because it offers a method to extract  (i.e., to \emph{sift out}) a single point of a function by integration.  Consider the following example.
\begin{equation*}
\int_{3}^5 y^3 \ln(y-2) \sin(y^2)e^{-y} \delta(y-4)\, dy
\end{equation*}
Ordinarily, this integral would be impossible to compute analytically (without the delta function, a symbolic integration software cannot resolve the integrand of this integral).  However, because of the sifting property of the delta function, this integral is simple.  Essentially, the integral is zero everywhere, except at the single point where the independent variable of integration is equal to 4.  At that point, we have

\[\delta(y-4) = \delta(4-4) = \delta(0)\]
Thus, the single point that contributes to the integral is the one where the argument of the delta function is zero, that is, the point $y=4$.  Thus, the result is

\begin{equation*}
\int_{3}^5 y^3 \ln(y-2) \sin(y^2)e^{-y} \delta(y-4)\, dy=4^3 \ln(4-2) \sin(4^2)e^{-4}
\end{equation*}
Which is just the integrand evaluated at $y=4$.

As a second, more abstract example, consider the following integration.  As above, assume that $x_0$ is any point defined within the domain of the integration.  Then, even the Gaussian function can be easily evaluated when integrated against a delta function.

\begin{equation*}
    \int_{-\infty}^{\infty} \frac{1}{\sqrt{2 a^2 \pi}} \exp\left(-\tfrac{1}{2}\frac{x^2}{a^2} \right) \delta(x-x_0) \, dx
    =\frac{1}{\sqrt{2 a^2 \pi}}\exp\left(-\tfrac{1}{2}\frac{x_0^2}{a^2} \right) 
\end{equation*}
We will see examples similar to this one in the course of studying the solution to partial differential equations.  There, the Gaussian is known as the heat kernel, and the delta function would represent the distribution of heat specified as an initial condition. 
\end{example}
\end{svgraybox}

\section{An Explanation of the Sifting Property}\indexme{delta function!sifting property}

For piecewise smooth functions, understanding how and why the delta function behaves the way it does is reasonably straightforwards.  We need only use the ideas of limits, and the expansion defined by the Taylor series.

 In the following development, we assume that the function $f$ is analytic near zero.  For concreteness, we will examine the following function, $f$

\begin{equation}
    f(x) = (8x+5)^2
\end{equation}

\begin{figure}[t]
\sidecaption[t]
\centering
\includegraphics[scale=.4]{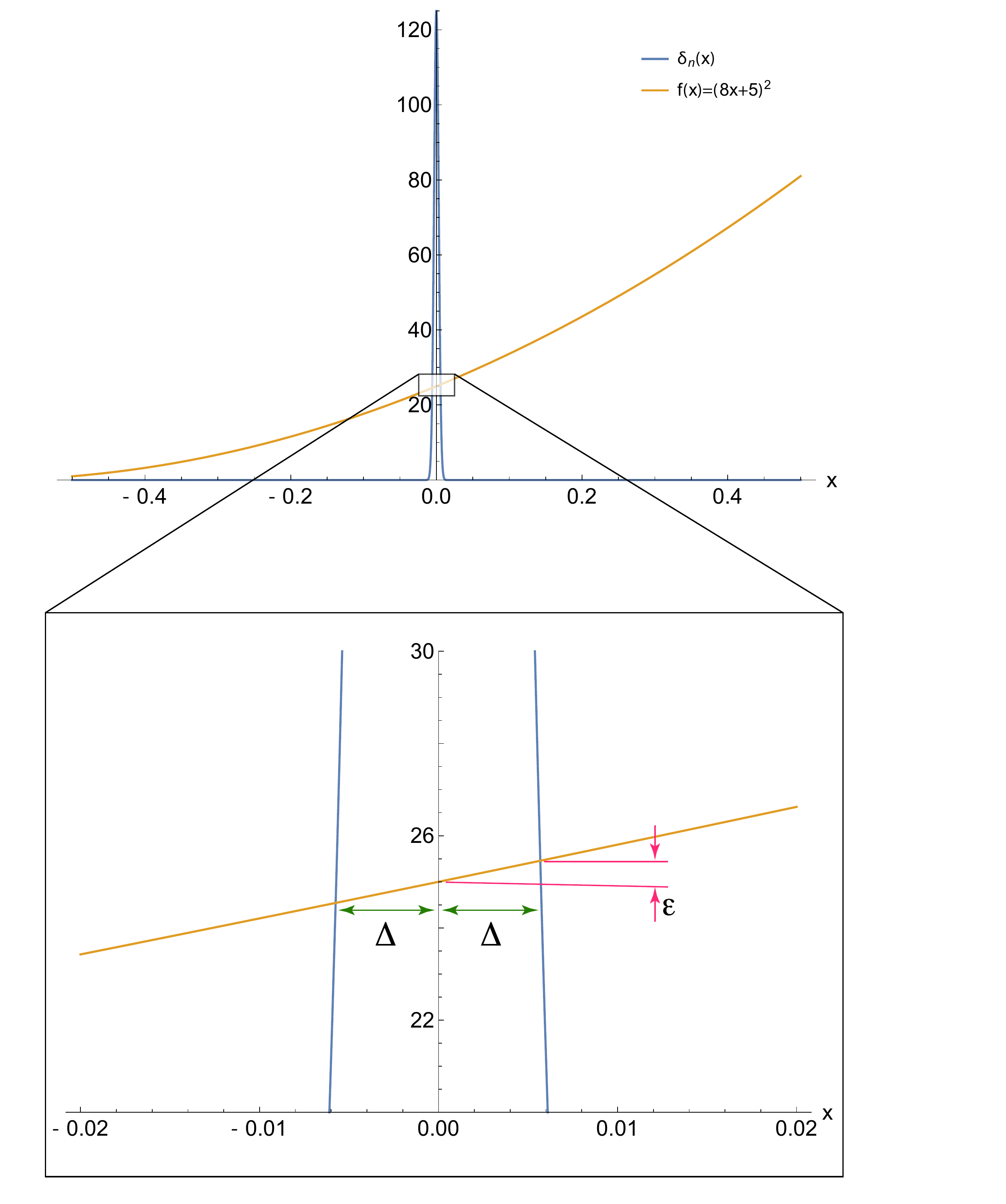}
\caption{{\bf A function, $f$ in the vicinity of a delta sequence function.}  The function changes by less than or equal to the amount $2\epsilon$ over the interval given by $\{I: -\Delta < x < \Delta\}$.  The value of $epsilon$ can be made as small as we like by making $\Delta$ smaller.  By definition, there is always a delta sequence with width less than $2\Delta$ around zero; we need only make the integer $n$ large enough.}
\label{deltataylor}       
\end{figure}

In Fig.~\ref{deltataylor}, the function is plotted, and one member of the of delta sequence functions is plotted for reference.

For the interval around zero, we can write the following Taylor series expansion for $f$ because $f$ is assumed to be analytic there.

\begin{equation}
    f(\epsilon)= f(0)+\Delta f'(0) +\frac{\Delta^2}{2!} f''(0) + \ldots
\end{equation}

Note that we can take $\Delta$ to be as small as we like because we can always find a delta sequence that has width less than $2\Delta$ simply by making $n$ large enough.  From chapter 1, we know that none of the derivatives of $f$ must be finite in the region around zero (this is a fact arising from the fact that $f$ is analytic).  Therefore, we can always find a $Delta$ sufficiently small such that 

\begin{align*}
    f(0) &\gg \Delta f'(0)\\
    f(0) &\gg \tfrac{1}{2!}\Delta f''(0)\\
    \ldots
\end{align*}

This completes the construction.  What we have shown is that we can always find a member of the delta sequence of functions such that the width of the function is small enough such that we can estimate $f$ by its value at $x=0$, with an error that can be made as small as we like.  This validates the use of the delta function as having the property, in the limit, of

\begin{equation}
\int_{x\in I} f(x)\delta(x) \, dx = f(0)
\end{equation}
The extension of the development to applications of the delta function at other points requires only that delta function be shifted (an affine transformation) appropriately.

\section{$^\star$A Computable Example of the Limit of a Delta Sequence}

One of the difficulties of working with delta sequences is that most of the results are understood (in the general cases, anyway) though abstract constructions.  In other words, it is not frequent that one can find a computable example that shows how a delta sequence converges in the sense that we have been discussing.  

There are some examples that are computable, with some effort.  In the example following, one example of a delta sequence is provided.  The goal for the example is to show that the sifting property of the delta function is valid when the sequence of functions are integrated against a continuous function $f(x)$.  This is done in two steps.  In the first step, the delta sequence is integrated against $f(x)$.  In the second step, the limit as $n\rightarrow\infty$ is evaluated after integration to show that the result is indeed the sifting property of the delta function.  It is important to note that the process of limits and integration do not necessarily commute!  In fact, in this case they distinctly do not commute.  If we first take the limit of the delta, and then attempt to to integrate, we find a situation where the delta function is supported at only a single point; regardless of what kind of integration one proposes, such an integral is necessarily zero!

\begin{svgraybox}
\begin{example}[A computable example of the sifting property]

There are many functions, both compact (i.e., defined on a closed interval) and non-compact (defined for the entire real number line).  For those who might be interested, compilations of such functions are available in the literature (e.g., \citet{dang2012dirac}).  The following function is one that is interesting, and also one that can be handled analytically.

\begin{equation*}
\delta_n(x) = 
\begin{cases}
\frac{1}{a_n} (1-4x^2)^n & \textrm{ for } \tfrac{1}{2} \le x \le \tfrac{1}{2} \\
 0 & \textrm{otherwise}
\end{cases}
\end{equation*}
where
\begin{equation*}
    a_n = \frac{\sqrt{\pi } \Gamma (n+1)}{2 \Gamma \left(n+\frac{3}{2}\right)}
\end{equation*}
This function is, in fact, the one provided at the start of this chapter.  Examples of the shape of this function with increasing $n$ are given in Fig.~\ref{deltaa}.  

One of the nice properties of this function is that it can be integrated analytically.  Thus, we can consider the following.\\

{
\centering\fbox{\includegraphics[scale=.5]{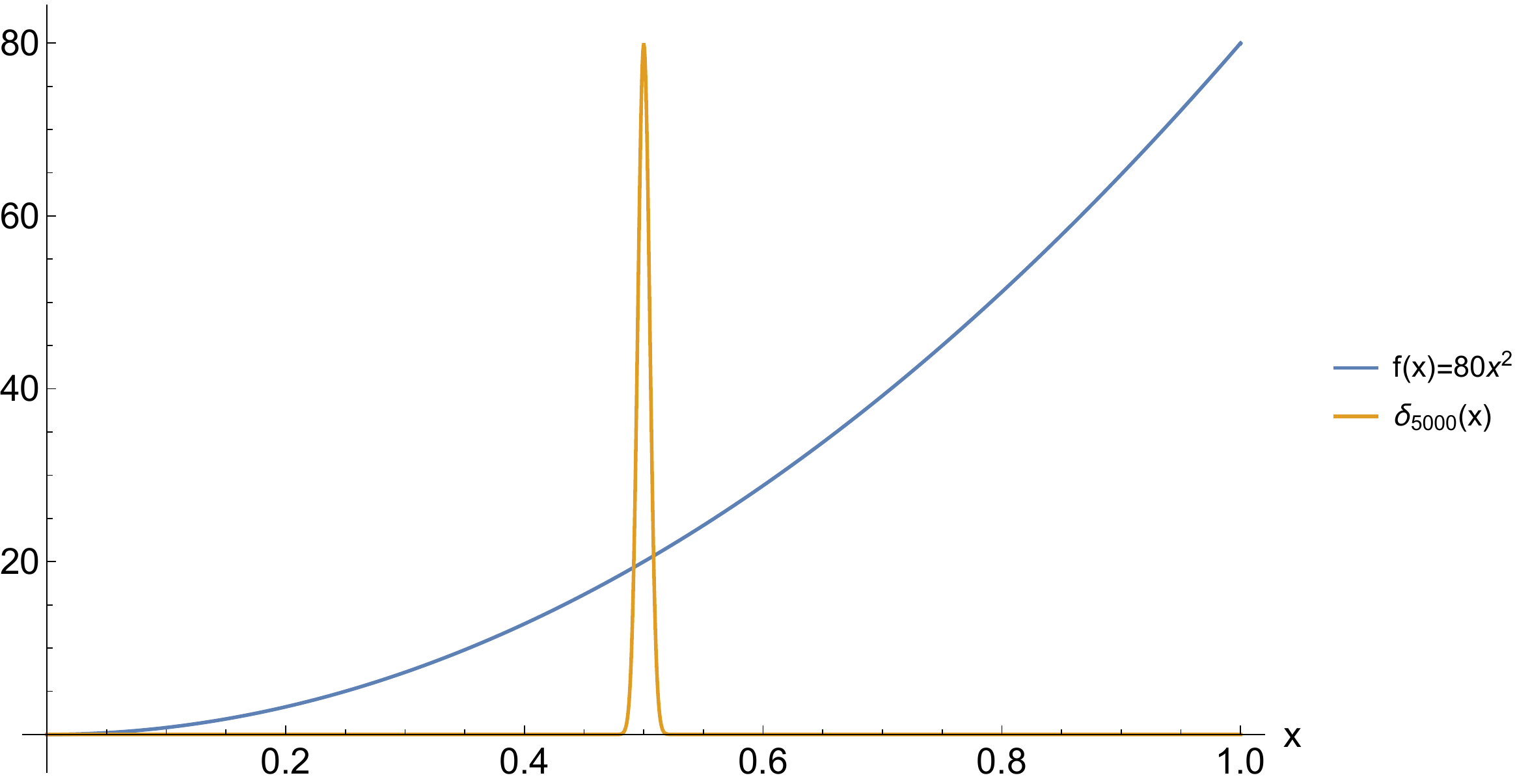}}
\vspace{2mm}
\captionof{figure}{The function $f(x)=80x^2$ plotted against the delta function $\delta_n(x-\tfrac{1}{2})$ defined above, with $n=5000$.}
\label{deltass}  
}
%
If we integrate our pre-delta function against an explicit function, say $f(x) = 80x^2$, then it may be possible to evaluate this integral.  In this example, we will shift the pre-delta function by one-half; an example of the pre-delta function and the function $f(x)$ are plotted in Fig.~\ref{deltass}.  Mathematically, the problme is represented as follows
\[ \int_{x=0}^{x=1} f(x) \delta_n(x-\tfrac{1}{2}) \, dx  
=
\int_{x=0}^{x=1}  \underbrace{80 x^2}_{f(x)}  \underbrace{ \frac{1}{a_n}(1-4(x-\tfrac{1}{2})^2)^n}_{\delta_n(x-\tfrac{1}{2})}  \, dx \]
It turns out that this definite integral can be computed for any value of $n$.  The result is

\[ \int_{x=0}^{x=1} f(x) \delta_n(x-\tfrac{1}{2}) \, dx= 80\frac{2+n}{2(3+2n)} \]

This is the result of the integration of any of the infinity of pre-delta functions, $\delta_n(x-\tfrac{1}{2})$,  against the function $f(x)=80 x^2$ for any value of $n$.  To find the final result, we need only take the limit of $n\rightarrow\infty$

\[ \mathop {\lim }\limits_{n \to \infty } 80\frac{2+n}{2(3+2n)} = 80\cdot\frac{1}{4}= 20
\]
This means that we have for our particular choice of pre-delta function, we have the result

\[\mathop {\lim }\limits_{n \to \infty } \left( \int_{x=0}^{x=1}80 x^2\left[\frac{1}{a_n}(1-4(x-\tfrac{1}{2})^2)^n \right] \, dx\right) = 20\]
and this result is exactly the result for $f(\tfrac{1}{2})=80(\tfrac{1}{2})^2$.  Thus, this provides at least one explicit example showing how the limit of delta functions yields exactly the result specified by the sifting property of the delta function.
\end{example}
\end{svgraybox}

\section{The Step Function and the (Generalized) Derivative of the Step Function}\indexme{step function!derivative}\indexme{Heaviside function!derivative}

\begin{figure}[t]
\sidecaption[t]
\centering
\includegraphics[scale=.8]{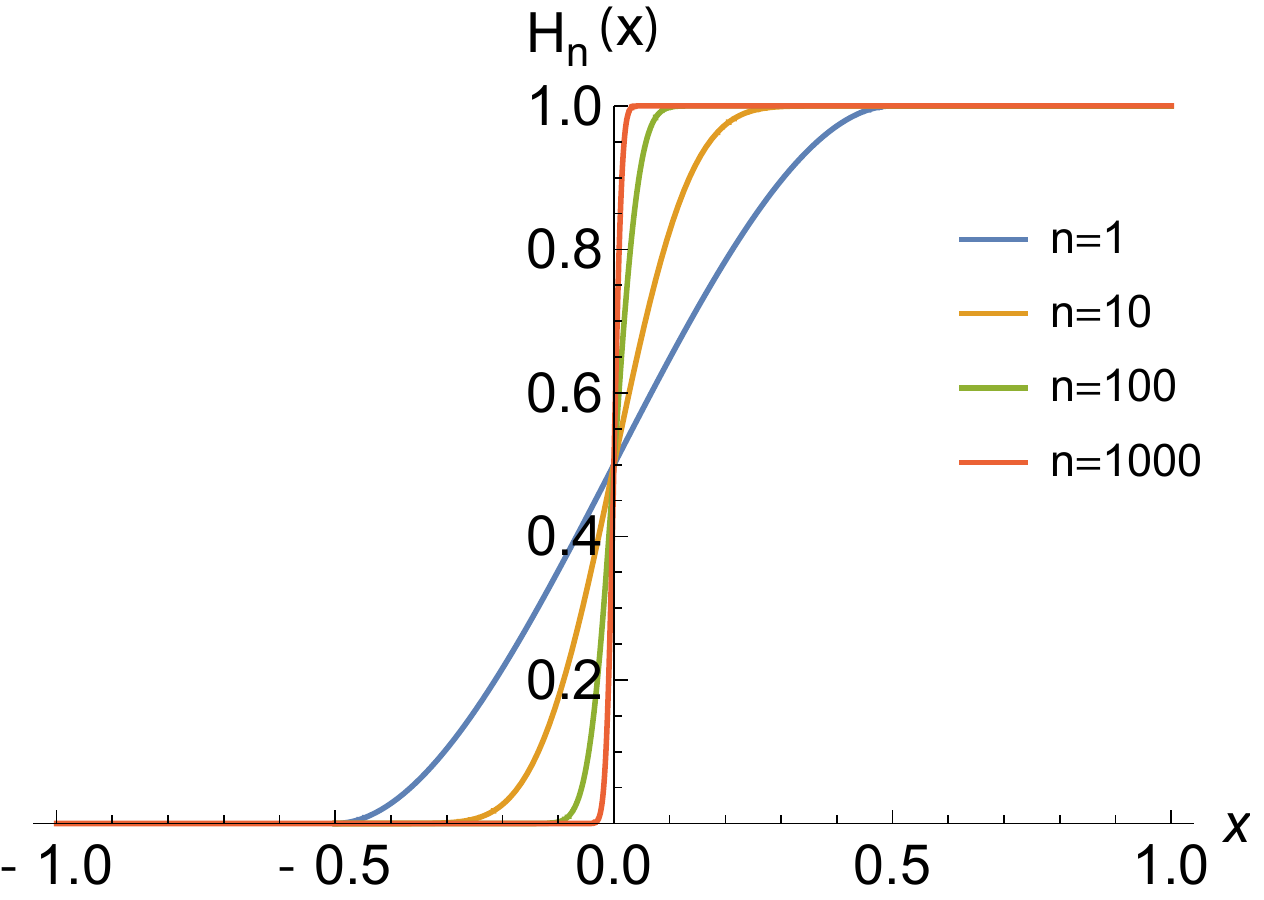}
\caption{{\bf Step function sequences.}  A \emph{sequence} of functions that approaches the \emph{integral} of the delta function in the limit.  This function converges to the step function as $n\rightarrow\infty$.}
\label{stepa2}       
\end{figure}

There is a reason that the delta function and the step function are covered together in this chapter.  It turns out that there is an interesting relationship between the two.  Consider the following representation.

\begin{equation}
H_n(x-x_0) = \int_{y=-\infty}^{y=x} \delta_n(y-x_0) \,dy
\label{stepdelt0}
\end{equation}
where $H_n(x-x_0)$ is a ``pre" step function; that is, it is similar to one of the S-shaped functions show in Fig.~\ref{stepa2} (where the functions are shown for $x_0=0$). As the value of $n$ increases, two things happen.  First, the function $\delta_n$ becomes more peaked and narrow (as, by now, we are accustomed to), and second, the S-shaped function becomes steeper and closer to a step.  It is clear from Eq.~\eqref{stepdelt0} that $H_n$ is a smooth and differentiable function.  The derivative is given by

\begin{equation}
\frac{d}{dx} H_n(x-x_0) = \delta_n(x-x_0) 
\label{stepdelt2}
\end{equation}

So we can consider the delta function to be the derivative of the step function.  Ordinarily, we should only attempt to interpret the delta function when it is under an integral; in other words, it is not a regular function but a generalized one, and technically it is not defined outside of an integral.  However, it is a common and even customary abuse of the notation to use the delta function as a symbol outside of an integral, particularly in physics, engineering, and the mathematics of partial differential equations.  However, it must always be kept in mind that the \emph{meaning} of this is that the delta function represents a sequence of functions whose limit must be interpreted only during the process of integration.

\section{Does the Delta Function Have a Fourier Series Representation?}\indexme{delta function!Fourier series}

This simultaneously deep and yet simple question.  The delta function is not a function at all, so asking whether or not it has a Fourier series representation may seem a bit odd.  However, the ``pre" delta functions, $\delta_n$ are in fact regular functions.  Depending upon which delta sequence is chosen, these pre-delta functions may even be $C^\infty$ smooth.

In the problems from the previous chapter, the Fourier series for a Gaussian function was examined.  It turns out that we can adopt this function to be a delta sequence as follows. First, we define the sequences as a function of $m$ (here, we have switched from $n$ to $m$ to prevent confusion when examining the Fourier series) as follows.

\begin{equation}
    \delta_m(x) =\frac{\alpha m}{\sqrt{\pi}\textrm{erf}(\alpha m)} \exp\left(-\alpha^2 m^2 x^2\right)
\end{equation}

\[A'_0= \frac{1}{2} \int_{-1}^1 \delta_m(x) \, dx = \frac{1}{2}\]
\[A_n= \frac{1}{1} \int_{-1}^1 \delta_m(x) \cos(n \pi x) \, dx = 
\exp\left(-\tfrac{\pi^2 n^2}{4 \alpha^2 m^2}  \right)
\]
Now note that in the limit of $m\rightarrow\infty$, we obtain the delta function rather than the delta sequences.  The Fourier coefficients take the values

\[A'_0= \frac{1}{2} \int_{-1}^1 \delta_m(x) \, dx = \frac{1}{2}\]
\[A_n= 1 \]
This leads to a series where all of the amplitudes are the same value - they are all unity (except $A'_0$).  What this says is that the delta function contains all possible frequencies, and the value of the amplitude does not, in general, decrease.  This is a very unusual result!  For more regular functions, it is possible to show that the amplitudes of the Fourier coefficients must decrease in value (this is related to Parseval's theorem), and, in fact, even yield a finite value when summed.  However, for the delta function, we do not have this relation; the sum of the amplitudes tends toward infinity.  While it is true that the Fourier series for the delta function violates Parseval's theorem, we should not be too surprised.  The delta function is not a regular function, and therefore the conditions required for proving Pareseval's theorem are not met by the delta function.

\begin{figure}[t]
\sidecaption[t]
\centering
\includegraphics[scale=.5]{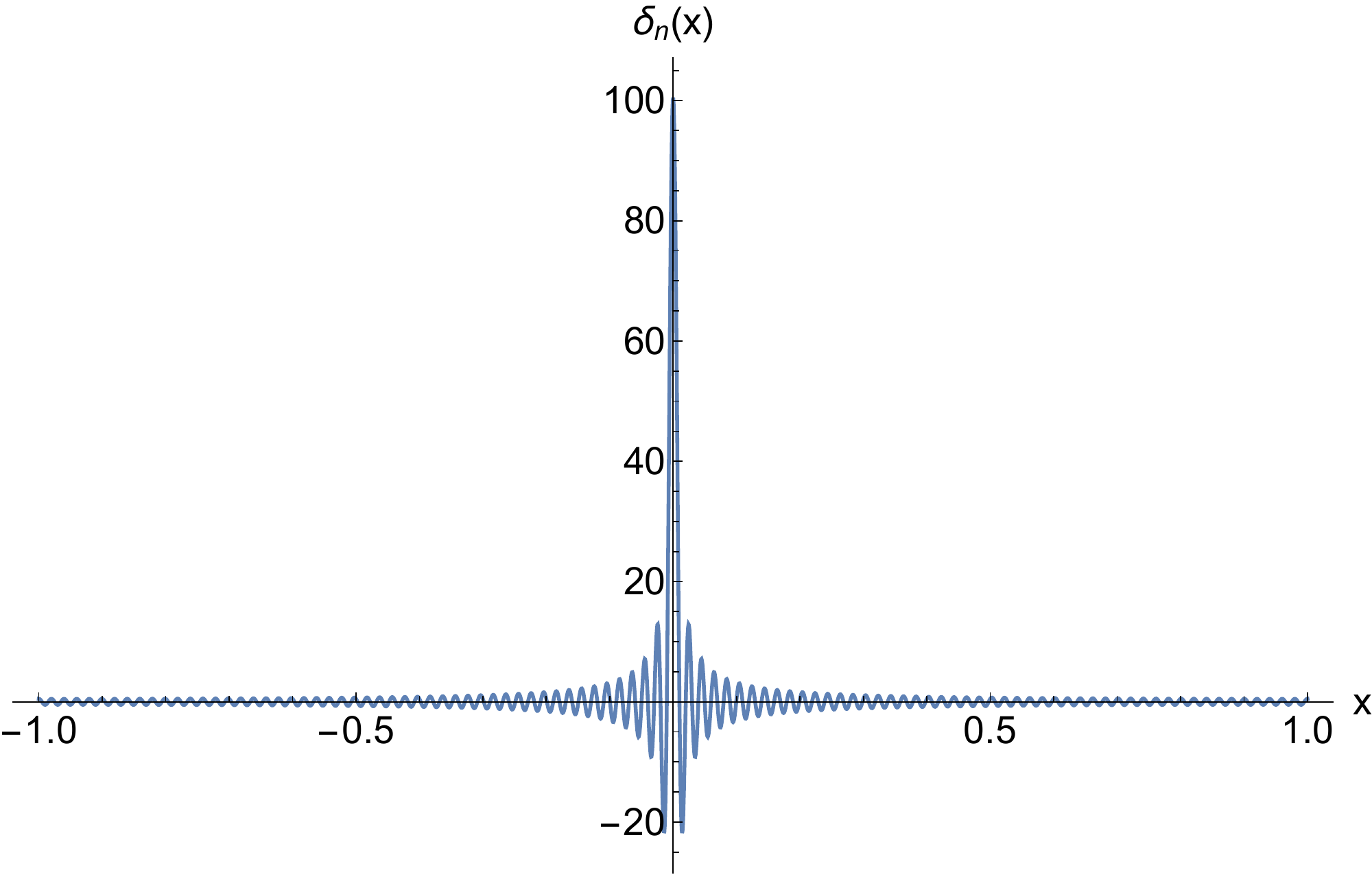}
\caption{{\bf The Fourier series sum for the delta function.}  The delta function has an unusual Fourier series expansion.  In this figure, the series is approximated by the first 101 terms in the sum.  Note that the resulting sum does actually become peaked and narrow at $x=0$, as expected for approximations to the delta function.}
\label{delta2}       
\end{figure}

It is interesting to look at the plot of the function predicted by summing up the Fourier series for the delta function; this plot is given in Fig.~\ref{delta2}.  The Fourier transform of the delta function is used in applications, such as in signal processing or modeling of point sources.  We will see more about Fourier analysis of the delta function when we study the Fourier transform later on.

\section{Some Identities for the Step and Delta Functions}\indexme{delta function!Identities}\indexme{Step function!Identities}

There are a few identities that can be useful when working with step and delta functions.  Some of these have been presented above, but are repeated here for reference.  

\begin{align}
    \delta(-t)& = \delta(t) \\
    \delta(t-a)&=\delta(a-t) \\
    \int_0^b \delta(x-a) f(x) \, dx &= f(a), ~~ b>a\\
    H(x-a)&=[1-H(a-x)]=[1-H(-(x-a))] \vspace{4mm}\\
    H(x)-H(-x)&=\textrm{sgn}(x),  { \large{\substack{\textrm{\transparent{0.0}\color{white}{!} }\\ \textrm{\hspace{15mm}sgn($x$) is the called the \emph{sign function},}\\  \textrm{\hspace{23mm}because it returns the sign of the argument}}}}\vspace{4mm}\\
    B(x;a,b)&= [H(x-a) - H(x-b)], ~~~b>a ~~~~\textrm{\hspace{5mm}The \emph{boxcar} function}\\
    \int_0^b \frac{d}{dx} H(x) \, dx& =  \int \delta(x) \, dx =1, ~~b>0\\
    \int_0^b \frac{d}{dx} H(x-a) \, dx &=  \int \delta(x-a) \, dx =1, ~~ b>a 
\end{align}

Two new functions are introduced here. The first is the \emph{boxcar} function, named because of its resemblance to a railroad boxcar.  The second is called the sign function, which has the useful property of returning the sign ($\pm 1$) of the argument of the function.

\begin{figure}[t]
\sidecaption[t]
\centering
\includegraphics[scale=.7]{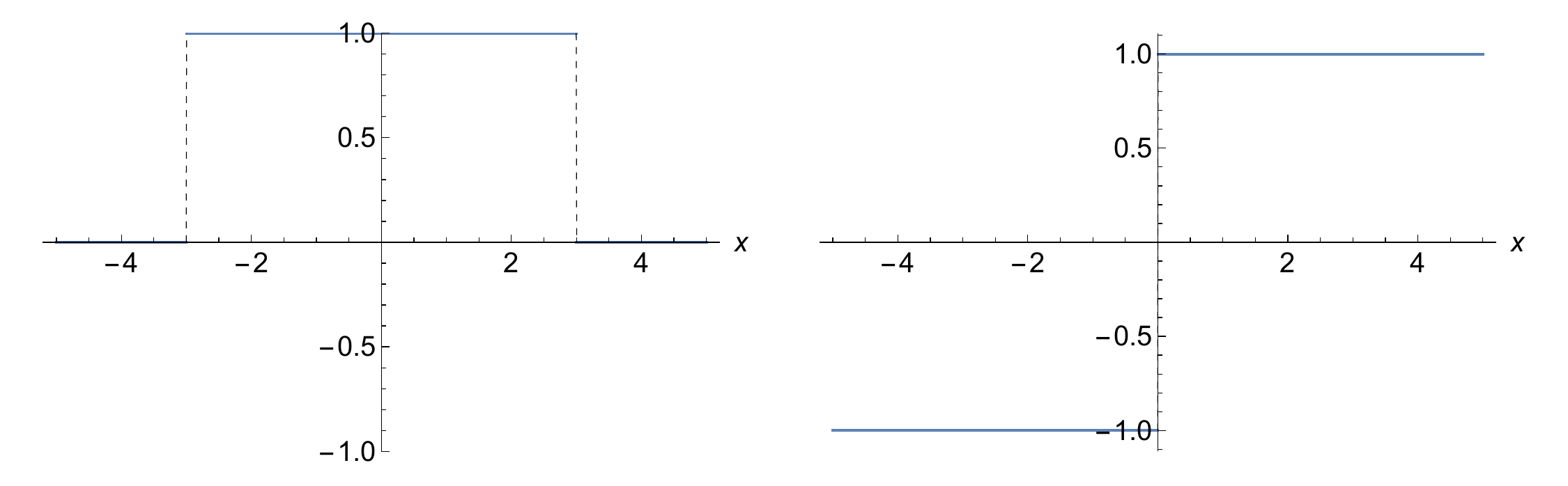}
\caption{{ (Left) The boxcar function with $a=-3$ and $b=3$. (Right) The sign function.}}
\label{boxcarsgn}       
\end{figure}
\newpage
\section*{Problems}
\subsection*{Practice Problems}

\noindent For the following problems, evaluate the integral.
\setlength{\columnsep}{2cm}

\begin{multicols}{2}
\begin{enumerate}
    \item \[ \int_{-1}^1 \delta(x) \, dx \]
    \item \[ \int_{-\infty}^{\infty} f(y) \delta(y)  \, dy \]
    \item \[ \int_{-\infty}^{\infty} y^3 \delta(y)  \, dy \]
    \item \begin{align*} \int_{-\infty}^{\infty} y^3 \delta(x-y)  \, dy \end{align*}
    \item \begin{align*} \int_{z=-1}^{z=1} \delta(z)f(x-z) \, dz,\\ 
    ~~\textrm{ where } 0<x<1\end{align*}
    \item \begin{align*}\int_{z=0}^{z=1} \delta(x-z)f(z) \, dz,\\ ~~\textrm{ where } 0<x<1 \end{align*}
    \item \begin{align*} \int_{-\infty}^{\infty}  \frac{\sin(y)}{\sin(x)}\delta(x-y) \, dy \end{align*}
    \item  \[ \int_{-\infty}^{\infty}  \frac{y^3-x^3}{x^2}\delta(x-y) \, dy \]
    \item \begin{align*} \int_{a}^{b} \sin({e^x}) \delta(x-x_0)\,dx,\\~\textrm{ where } a<x_0 <b\end{align*}
    \saveenumerate
\end{enumerate}
\end{multicols}

\setlength{\columnsep}{2cm}
\begin{enumerate}[topsep=8pt,itemsep=4pt,partopsep=4pt, parsep=4pt]
\restoreenumerate
    \item There are sequences that are the complement of delta sequences.  These are the zero sequences.  A delta sequence an extreme transformations of a function; when integrating a function against a delta sequence, the only information that passes through comes from a single point.  The zero sequences act in a complementary way; although its integral is always unity, when integrating a function against a zero sequence, in the limit the result is zero.  As an example, the following is an zero sequence on the real line.
    \begin{equation*}
   \iota_n(x) = 
        \begin{cases}
         \frac{1}{2 n}  &\textrm{ for } -n \le x \le n, ~~n\in\mathbb{N} \\
          0 & \textrm{ otherwise}
        \end{cases}
    \end{equation*}
    
    \item Show that the function 
    \[ (a n) Z_n\left(\frac{x}{a n}\right) \]
    Is an zero function as $n\rightarrow \infty$
    
    \item Show that every delta sequence can be transformed into an zero sequence by the substitution $n\rightarrow 1/n$.
        \saveenumerate
\end{enumerate}

\subsection*{Applied and More Challenging Problems}
\begin{enumerate}
    \item Heating of a surface at point can be represented by using the delta function.  For this problem, suppose you are interested in finding the steady state temperature distribution for an insulated wire (so it does not loose heat to the environment) of length $L$ and constant cross-sectional area, $A$.  The wire is heated at its center point, $x=L/2$, by a powerful laser (where a small amount of insulation has been removed so that the wire can be heated).  Suppose also that the left end of the wire is held in liquid nitrogen at $T_0=77$ Kelvin (K) and the right side is held at $T_1=194.5$ K in liquid $\textrm{CO}_2$.  Assume $L=1$ m.
    \begin{align*}
        K_T \frac{d^2 T}{d x^2} &= -q_0 \delta(x-L/2) \\
        T(0) &= T_0\\
        T(L) &= T_1\\
    \end{align*}
Note!  The delta function (curiously) must have units of $1/Length$ if it is to have an integral of unity (without units).  Thus, the equation above is dimensionally correct.
For this problem, assume the heat flux, $q_0$ is $q_0=$ 200,000 W m$^{-2}$, and the thermal conductivity is $K$= 385 W m$^{-1}$K$^{-1}$.  Please plot the final solution over the whole domain.

\item One curious feature about delta functions is that they allow discrete data to be given a functional representation.  Such features are part of why delta functions are not true ``functions".  To expand on this comment, first we just need to realize that the area under a point is identically zero when ``measured" by an integral.  As an example, think about the function $f(x)$, $x\in[0,1]$ containing a single point, as follows
\begin{equation*}
f_1(x)=
    \begin{cases}
    2 & x=\frac{1}{2}\\
    0 & \textrm{otherwise}
    \end{cases}
\end{equation*}
To integrate this, we might write
\begin{equation*}
    I_1(x) =\int_0^1 f_1(x)\,dx
\end{equation*}
But, noting that $f(x)$ is zero at every point \emph{except} at $x=a$, we must have the equivalent result
\begin{equation*}
    I_1(x) =\int_a^a f_1(x)\,dx \, = 0
\end{equation*}
In other words, by a classical analysis, the area under a single point is zero.  Now suppose we \emph{allow} the delta representation, which you may recall is only technically defined when integrated.   However, because of this understanding, we frequently write the delta (generalized) function as if it were a regular function.  Thus, the function $f(x)$ can be written as follows

\begin{equation}
    f_2(x) = 2\delta\left(x-\frac{1}{2}\right)
\end{equation}
So that upon integration we would find

\begin{equation*}
    I_2(x) =\int_0^1 \delta(x-\tfrac{1}{2})\,dx = 2
\end{equation*}
In other words, the discrete data point $(1/2,2)$ has an integral interpretation if we consider the delta function.

The \emph{meaning} is different between $f_1(x)$ and $f_2(x)$.  But, one curious thing that the delta function allows is to provide a way of representing discrete data in the same form as a regular function.  To be clear about this, let's consider the following discrete data set, $X$, given by the $(x,y)$ data pairs
$X=\{(1,1), (2,2.1), (3,2.9), (4,4.2), (5,5.05)\}$.  Based on the discussion above, we could write this data set in a form suitable for integration as follows.  To avoid any problems due to end point evaluations, suppose that the domain for this function is $x\in[0,6]$

\begin{equation}
    g(x) = 1\delta(x-1)+2.1\delta(x-2)+2.9\delta(x-3)+4.2\delta(x-4)+5.05\delta(x-5)
\end{equation}

With this definition, please try the following.
\begin{enumerate}
    \item Compute the average value of $X$ by using a conventional sum divided by the number of data points.  
    
    \item The average of a continuous function, say $f(x), x\in[0,L]$ is usually given by 
    \begin{equation*}
        A = \frac{1}{L} \int_0^1 f(x),\ dx
    \end{equation*}
    Where $L$ measures the magnitude of the domain.  One way of thinking about $L$ is that it is also defined by an integral that measures distance

    \begin{equation*}
        L = \int_0^L 1 \, dx = \left. x\right|_0^L = L
    \end{equation*}
    While this is a bit tautological, it is helpful in making the following correspondence.  When dealing with delta functions, we define the size of the domain differently.  For example, consider $g(x),~x\in[0,6]$

    \begin{equation*}
    L_\delta = \int_{x=0}^{x=6}  \delta(x-1)+\delta(x-2)+\delta(x-3)+\delta(x-4)+\delta(x-5) \, dx = 1+1+1+1+1 = 5
    \end{equation*}
    
    When we have delta functions involved, the measure we want for the domain is the \emph{number} of delta functions involved.  But, this arises (conceptually at least) from an analogue with continuous distance.  With this in mind, compute the average of $g(x)$ using an appropriate integral.

    \item The \emph{center of mass} of a discrete data set with $N$ data points of the form $(x_i,y_i),~i=1\ldots N$ is given by 
    \begin{equation*}
        \bar{y} = \dfrac{\sum_{i=1}^N x_i y_i}{\sum_{i=1}^N  y_i}
    \end{equation*}
    Compute this value for the data set $X$.  Then, illustrate how the same value can be computed using integral quantities of the function $g(x)$.
\end{enumerate}

\end{enumerate}
\abstract*{This is the abstract for chapter 00}

\begin{savequote}[0.55\linewidth]
`` Science is a partial differential equation.  Religion is a boundary condition."

\qauthor{Alan Turing, genius developer of modern computing theory. }
\end{savequote}

\chapter{Introduction to Partial Differential Equations, Conservation Laws, and Constitutive Equations}\label{introPDEs}
%
\def\CHAP {chapter05_general_structure_PDEs}
\section{Introduction}

Partial differential equations (PDEs) are the mathematical workhorses of continuum systems.  The term \emph{partial} differential equation indicates that the dependent variable depends on two or more independent variables.  This is in contrast to an ordinary differential equation (Chapter 2), where there is always one dependent and one independent variable.  For this chapter, we will use the symbol $u$ to indicate the independent variable, and the conventional symbols $x$ and $t$ to represent the independent variables corresponding to one spatial dimensions and one time dimension.  Examples of the three canonical types of PDEs (hyperbolic, parabolic, and elliptic) will be presented.  

\section{Terminology}

\begin{itemize}

\item {\bf Extensive variable.}  An extensive variable (or property) is one that depends upon the size of the domain of interest.  As an example, the mass of material in a domain is an extensive property of that domain.  \indexme{variables!extensive}\\

\item {\bf Intensive variable.}  An intensive variable (or property) is one that is independent size of the domain of interest.  While mass is an \emph{extensive} property of a domain, density is an \emph{intensive} property, because it depends in no way upon the total size of the domain.\indexme{variables!intensive}  \\

\item {\bf Flux.}  A flux is the amount of some extensive property passing through a unit area per unit time.  Technically, a flux may be more generally defined in terms of vector quantities; under these circumstances the flux is defined as the amount of the extensive property passing through a unit area per unit time, where the unit area is defined with respect to a specific unit vector indicating the direction normal to the surface; both normal and tangential fluxes can be defined in multiple dimensions.  In one space dimension these additional distinctions are not material. \indexme{flux}\\

\item {\bf Partial derivative.}  Partial derivatives are defined whenever there is more than one independent variable defined for function.  To be a little more concrete, consider derivatives of a function of one space dimension, $x$, and one time dimension, $t$.  The partial derivative in space would be formulated as for an ordinary derivative in space, where the time component of the function were simply held constant.  \\

\item {\bf Partial differential equation (PDE).} An differential equation that contains more than one independent variable.  The equations themselves are composed of algebraic combinations of functions multiplied by partial derivatives.   Often, such equations arise from physical considerations of continuum systems, although they may arise in any number of ways from any number of disciplines.  For such equations, partial derivatives do not necessarily have to appear for each independent variable. \\

\item {\bf Order of a PDE.}  The order of a PDE is given by the highest order of derivative that appears in the equation.\\

\item {\bf Transient PDE.}  Assume we have a partial differential equation with one space and one time variable.  A transient PDE is one where an initial condition evolves to new configurations in the spatial dimension as time changes. \\

\item {\bf Steady state conditions.} The steady state for a PDE is the condition where the time partial derivative is identically zero for all other combinations of the remaining independent variables.  For example, in one space and one time dimension with domain $x\in[a,b]$, the steady state is defined by $\partial u(x,t)/\partial t = 0$ for all points in $[a,b]$; note that if this occurs at some finite time, $t_1$ then it must also be true that $\partial u(x,t)/\partial t = 0$ for all $t>t_1$.  In other words, once the system achieves steady state, it remains in that state (unless the system is perturbed by external forces). \\

\item {\bf Ancillary conditions.}\indexme{partial differential equation!ancillary condition}  Like for ODEs, ancillary conditions are the additional information required to determine a \emph{particular} solution from the set of all possible solutions to a PDE.  Unlike ODEs, because more than one independent variable is involved, the ancillary conditions technically specify \emph{whole functions} rather than just constants (although constant functions are very common!)  \\

\item {\bf Parabolic PDE.} In one time and one space dimension, a parabolic PDE is one characterized by an initial configuration that evolves smoothly over time.  The classical example of a parabolic PDE is the heat  or diffusion equation.  For that equation, an initial distribution spreads out over time, becoming smoother and more uniform over time. \\

\item {\bf Hyperbolic PDE.} In one time and one space dimension, a hyperbolic PDE is one that represents wave-like behavior.  An initial perturbation is translated through the domain either without changing shape, or with shape changes (which can result either from self-sharpening behavior in nonlinear wave equations, to dissipating behavior in equations that represent energy loss from waveforms as they translate). The first-order wave equation is the primary hyperbolic example presented in this chapter.\\

\item {\bf Elliptic PDE.}  There is no exact analogue for an elliptic PDE in one time and one space dimension.  However, in one space dimension alone, and elliptic PDE can be considered to be the \emph{steady state} version of a corresponding parabolic PDE.  As an example, a steady-state diffusion-reaction problem can generate a solution in which the solution is not uniform in space, but at each spatial point the solution does not change in time.  Physically, this would indicate that the processes of diffusion and reaction balanced one another.  \\

\item {\bf Well-posed PDE problem}.  A well posed problem in the context of PDEs means the following \indexme{partial differential equation!well posed}
    \begin{enumerate}
        \item The solution exists.  In other words, there is sufficient information provided by the ancillary conditions such that the problem can be solved.
        \item The solution is unique.  This requirement indicates that not only must solutions be found, but one must have enough information to generate a \emph{particular} solution to the problem.  Another way of stating this is that the solution should not represent a \emph{class} of functions, but, rather, a specific function that is computable given all of the parameters involved.
        \item The solution must depend continuously upon the initial data.  This concept is a bit less obvious than the other two, but, in short, it means only that the solution to the problem cannot behave chaotically.  Here, not behaving \emph{chaotically} means that very small changes (or \emph{perturbations}) to the initial conditions do not lead to dramatically different results.  This is a feature of well-posed linear equations: small perturbations to the system lead only  to small changes in the solution.
    \end{enumerate}
\end{itemize}

\section{Partial Derivatives and the Types of Partial Differential Equations}

Many processes in science and engineering are governed by partial differential equations.  Essentially, most systems that can be approximated as being a \emph{continuum} and where one needs to understand the time-space behavior of the system (as opposed to the integrated qualities of the system as appears in material balances) are best described by partial differential equations.  Before continuing with examples, however, it is useful to recall the definition of a \emph{partial derivative}. 

\subsection{Partial derivatives}\indexme{derivative!partial}

A partial derivative describes the slope of a curve that is a function of more than one independent variable.  For example, suppose a function $u$ depends upon the spatial variables $x$ an $y$.  Then, the derivative in the $x-$ and $y-$directions are given by

\begin{align}
    \mathop {\lim }\limits_{\Delta x \to 0} \frac{\partial u(x,y)}{\partial x} & \equiv \frac{[u(x+\Delta x,y)-u(x,y)]}{\Delta x}\\
    \mathop {\lim }\limits_{\Delta y \to 0} \frac{\partial u(x,y)}{\partial y} & \equiv \frac{[u(x,y+\Delta y)-u(x,y)]}{\Delta y}\\
\end{align}
Here it is understood that $\Delta x$ is a positive value.  
Similarly, the second derivative with respect to $x$ and $y$ can be specified by

\begin{align}
     \mathop {\lim }\limits_{\Delta x \to 0} \frac{\partial^2 u}{\partial x^2} & \equiv \frac{[\tfrac{\partial u}{\partial x}(x+\Delta x,y)-u(x,y)]}{\Delta x}\\
    \mathop {\lim }\limits_{\Delta y \to 0} \frac{\partial u(x,y)}{\partial y} & \equiv \frac{[u(x,y+\Delta y)-u(x,y)]}{\Delta y}
\end{align}
These expressions can be iterated to derive even higher-order derivatives.  As a final note, there are a few common methods for indicating partial derivatives in PDEs.  The other conventional method is to use subscripts to indicate the independent variable of differentiation.  The subscript is repeated $n$ times to indicate an $n^{th}$ order derivative.  Thus, the following are equivalent

\begin{align}
    \frac{\partial u}{\partial x}& ~~\Leftrightarrow~~ u_x&~~  \frac{\partial u}{\partial t}& ~~\Leftrightarrow~~ u_t\\
    \frac{\partial^2 u}{\partial x^2}& ~~\Leftrightarrow~~ u_{xx}&~~  \frac{\partial^2 u}{\partial t^2}& ~~\Leftrightarrow~~ u_{tt}
\end{align}

\subsection{Characterization of Linear PDE Types}
Partial differential equations can be characterized primarily by three major features.

\begin{enumerate}
    \item {\bf Equation order.}  Partial differential equations are characterized in part by the order of the highest derivative (in any independent variable).  This is known as the \emph{order} of the PDE.  The vast majority of PDEs seen in science and engineering are of first or second order, although equations up to fourth order are sometimes encountered.  Equations above fourth order are seldom encountered, but they occasionally do arise as the description of some physical process.  All first order equations in 2 independent variables describe wave-like behavior; these are listed in Table~\ref{firstorder}.  Some examples of second-order equations are given in Table~\ref{secondorder}.  Relevant examples of linear equations of order higher than two are given in Table~\ref{higherorder}.\indexme{partial differential equation!order}\\
    
    \item {\bf Equation type.}  There are several well-studied behaviors seen in linear partial differential equations.  In particular, for linear second-order partial differential equations \emph{in two independent variables} of the form
    
    \begin{equation} a \frac{\partial^2 u}{\partial t^2}+ 2b \frac{\partial^2 u}{\partial t \partial x}+ c \frac{\partial^2 u}{\partial x^2}+ d\frac{\partial u}{\partial x} + ed\frac{\partial u}{\partial t} + f u + g = 0 \label{gen2ndPDE} \end{equation}
    
    every such can be described as exhibiting three types of behavior. Note that more than one type of behavior may be exhibited by the same equation if the coefficients are functions; that is, the equation type can change as a function of time or space. For equations with constant coefficients, only one equation type is defined, and it is does not depend on any of the independent variables. \\
    
    With the appropriate change of variables, every equation of the form of Eq.~\eqref{gen2ndPDE} can be put in its \emph{canonical form}.  While determining the equation type is a useful exercise, it does not necessarily tell us all that we need to know about a particular problem.  For example, the fact that the general equation given by Eq.~\eqref{gen2ndPDE} can always be put into its canonical form by a coordinate transformation completely ignores the fact that if such an transformation may make any associated boundary and initial conditions so complex as to render the problem nearly unsolvable.  Thus, the idea of the canonical form must be taken with the proverbial grain of salt.  While it is useful information, in more practical problems this information may be of minor importance.
    
    \begin{enumerate}
        \item {\bf Parabolic.}  The canonical form for the parabolic equation with one time and one space variable is \indexme{partial differential equation!parabolic}
        
        \[ \boxed{~~~~\frac{\partial u}{\partial t} -  \frac{\partial^2 u}{\partial x^2}=0~~~~} \]
        
       These equations have generally smooth, dissipative solutions.  The heat and diffusion conservation equations are examples.   
        \item {\bf Hyperbolic.} The canonical form for the parabolic equation with one time and one space variable is\indexme{partial differential equation!hyperbolic}
        
         \[\boxed{~~~~ \frac{\partial^2 u}{\partial t^2} -  \frac{\partial^2 u}{\partial x^2}=0~~~~} \]
        
        Hyperbolic equations are characterized by representing wave-like behavior.  These waves can be smooth and periodic (like ocean waves), or sharp, like the wave front of a compressible gas that creates a sonic boom.
        
        \item {\bf Elliptic.}  The canonical form for the parabolic equation with two space variables is \indexme{partial differential equation!elliptic}
        
        \[ \boxed{~~~~\frac{\partial^2 u}{\partial y^2} +  \frac{\partial^2 u}{\partial x^2}=0~~~~} \]
        
        While mathematically, there can be time-dependent elliptic equations, most physical applications come from steady-state descriptions of problems that are determined by their boundary conditions.  The steady, 2-dimensional distribution of heat in a plate would be an example.  Steady groundwater flow is also an example of an elliptic equation.  
    \end{enumerate}
    
    There are some problems with formal classification systems for PDEs that prevent them from being universally useful.  For example, there is a classification system for all linear second-order PDEs with two (or more) independent variables.  While such classifications can provide some insight, they still have weaknesses.  For example, in these classification systems, the role of derivatives less than order 2 are not considered at all.  This can create some difficulties, especially when these first-order terms dominate the problem.  As an example, the classical convection-diffusion equation in one spatial dimensions can be given by
    
    \[ \frac{\partial u}{\partial t} = - v_0 \frac{\partial u}{\partial x}+ D \frac{\partial^2 u}{\partial x^2}  \]
    
    While the traditional classification schemes presented in most texts on partial differential equations would qualify this as a \emph{parabolic} equation, this is not an entirely satisfactory answer.  In particular, when the convective term (the first-order term) is much larger than the diffusive term (the second-order term), then this equation behaves as a \emph{hyperbolic} equation.  This is of some significance, because, for example, numerical methods for solving parabolic (smooth) and hyperbolic (sometimes abrupt) equations can be quite different.    As a second example, consider the following elliptic equaiton in two spatial variables
    
    \[ K_x \frac{\partial^2 u}{\partial x^2} + K_y \frac{\partial^2 u}{\partial y^2}=0 \]
    
    while this is a perfectly acceptable elliptic PDE, we encounter a problem in the limit $K_y \rightarrow 0$.  In that limit, our equation is suddenly no longer elliptic by the conventional criterion.  While it is true that it can be then considered to be an ordinary differential equation (since there is only one independent variable), it does raise some problems in categorization.  \\
    
    Therefore, in this text while we will use the terms \emph{parabolic}, \emph{hyperbolic}, and \emph{elliptic}, to be descriptive, there will be no special effort made to cover the formal methods of characterizing linear PDEs by invoking a transformation of the coordinates.
    For higher-order PDEs, technically there is no conventional classification system; however, occasionally one of the categorizations listed above (parabolic, hyperbolic, elliptic) will be used if the solutions to the higher-order equation share similar features as the behavior for a second order equation. \\

    \item {\bf Linear or nonlinear.}  Linear equations describe many phenomena in engineering and science, often as an approximating behavior of a more general nonlinear equation.  The linearity of a PDE can be determined using the same techniques that are used to establish linearity for ODEs.  For some phenomena (e.g., turbulence in the Navier-Stokes equations), nonlinearity is essential to phenomenon being represented.  A few examples of nonlinear problems are given in Table~\ref{nonlinear}.  While nonlinear problems are both relevant and significant, their solutions generally require specialized methods that are beyond the scope of this introductory text.
\end{enumerate}

\begin{table}
\centering
\caption{Examples of first-order linear PDEs with two independent variables.  All first order PDEs represent waves of one form or another.} 
    \begin{tabular}{|c|c|c|c|c|}
        \hline
        \makecell{~\\Type\\~}   & \makecell{~Independent Variables~\\ in 2 Dimensions} & Name  & Example Application & Example Equation     \\ 
        \hline\hline
        \makecell{~\\~One-dimensional waves~\\~}  & \makecell{space, $x$ \\ time, $t$}              & \makecell{convection equation} &  Solute transport &
        $\frac{\partial {\bf u}}{\partial x}=-v_0 \frac{\partial u}{\partial x}, ~~~v_0>0$        \\ 
         \hline
       \makecell{~\\~One-dimensional waves~\\~\\~}   
       & \makecell{space, $x$ \\ time, $t$}   & \makecell{The Maxwell equations}   
       & Electrodynamics 
       & \makecell{$
       \frac{\partial E}{\partial x}=\frac{1}{c}\frac{\partial B}{\partial t}$ \\ ~\\
     $ \frac{\partial B}{\partial x}=\frac{1}{c}\frac{\partial E}{\partial t}+J $}\\
        \hline
    \end{tabular}
    \label{firstorder}
\end{table}

\begin{table}
\centering
\caption{Examples of second-order linear PDEs with two independent variables}
    \begin{tabular}{|c|c|c|c|c|}
        \hline
        \makecell{~\\Type\\~}   & \makecell{~Independent Variables~\\in 2 Dimensions }& Name  & Example Application  &Example Equation     \\ 
        \hline\hline
        \makecell{~\\Parabolic\\~}  
        & \makecell{space, $x$ \\ time, $t$}              
        & ~~The heat equation~~ 
        & Heat or mass transport 
        & $\frac{\partial u}{\partial t}=K\frac{\partial^2 u}{\partial x^2}$        \\ 
         \hline
        \makecell{~\\Elliptic\\~}   & \makecell{space, $x$ \\ space, $y$}   
        & ~The Laplace equation~    
        & Steady groundwater flow 
        &$~~K\frac{\partial^2 u}{\partial x^2}+K\frac{\partial^2 u}{\partial y^2}=0~~$ \\ 
        \hline
        \makecell{~\\~~Hyperbolic~~\\~} & \makecell{space, $x$ \\ time, $t$}   
        & The wave equation
        & \makecell{Movement of deep \\water waves}
        & $~~\frac{\partial^2 u}{\partial t^2}=K\frac{\partial^2 u}{\partial x^2}~~$  \\
        \hline\hline
         \makecell{~\\~~(none)~~\\~} & \makecell{space, $x$ \\ time, $t$}   
        & \makecell{The linear Korteweg--\\de Vries equation}
        & \makecell{Movement of shallow \\water waves}
        & $~~\frac{\partial u}{\partial t}=K\frac{\partial^3 u}{\partial x^3}~~$  \\
        \hline
        \makecell{~\\~~(none)~~\\~} & \makecell{space, $x$ \\ time, $t$}   
        & ~The vibration equation~
        & \makecell{Vibrations in beams \\or thin plates}
        & $~~\frac{\partial^2 u}{\partial t^2}=K\frac{\partial^4 u}{\partial x^4}~~$  \\
        \hline
    \end{tabular}
    \label{secondorder}
\end{table}

\begin{table}
\centering
\caption{Examples of a few linear, higher-order PDEs with two independent variables}
    \begin{tabular}{|c|c|c|c|c|}
          \hline
        \makecell{~\\Name\\~}  &   \makecell{~Independent Variables~\\ in 2 Dimensions} &   Application &  Equation     \\ 
        \hline\hline
        \makecell{~\\~~The linear Korteweg--~~\\de Vries equation\\~} & \makecell{space, $x$ \\ time, $t$}   
        & \makecell{~~~Movement of shallow~~~ \\water waves}
        & $~~\frac{\partial u}{\partial t}=K\frac{\partial^3 u}{\partial x^3}~~$  \\
        \hline
        \makecell{~\\~The vibration equation~\\~} & \makecell{space, $x$ \\ time, $t$}   
        & \makecell{Vibrations in beams \\or thin plates}
        & $~~~~~~~~\frac{\partial^2 u}{\partial t^2}=K\frac{\partial^4 u}{\partial x^4}~~~~~~~~$  \\
        \hline
    \end{tabular}
    \label{higherorder}
\end{table}

\begin{table}
\centering
\caption{Examples of some important nonlinear PDEs in two independent variables.}
    \begin{tabular}{|c|c|c|c|}
        \hline
        \makecell{~\\Name\\~}  &   \makecell{~Independent Variables~\\ in 2 Dimensions} &   Application &  Equation     \\ 
        \hline\hline
        \makecell{~\\Nonlinear diffusion\\~}  & \makecell{space, $x$ \\ time, $t$}              & \makecell{~Solute transport~\\~at high concentrations~} &
        $\frac{\partial { u}}{\partial t}=D(u) \frac{\partial^2 u}{\partial x^2}$      \\ 
         \hline
       \makecell{~\\Burgers' equation\\~}  & \makecell{space, $x$ \\ time, $t$}              &  \makecell{~one-dimensional analogue\\to Navier-Stokes~} &
        $~\frac{\partial {u}}{\partial t}=u \frac{\partial u}{\partial x}+\nu \frac{\partial^2 u}{\partial x^2}~$      \\ 
        \hline
         \makecell{~\\~The nonlinear Korteweg--~\\de Vries equation\\~}  & \makecell{space, $x$ \\ time, $t$}    &  \makecell{~shallow water waves with~\\~nonlinear restoring function~} &
        $~~\frac{\partial {u}}{\partial t}=\frac{\partial^3 u}{\partial x^3}-6 u \frac{\partial u}{\partial x}~~$      \\ 
        \hline
    \end{tabular}
    \label{nonlinear}
\end{table}

\section{The Origins of Partial Differential Equations in Science and Engineering: The Axioms of Conservation} 

There is a nearly inexhaustible supply of applications of PDEs in engineering and science.  Essentially, any systems that a subject to the laws of continuum mechanics are usually expressed as partial differential equations.  Usually, these partial differential equations arise by considering conservation principles.  Conservation principles are \emph{axiomatic} statements (or \emph{laws}) \indexme{axiom} that cannot be proven from more fundamental ideas, but they are consistent with all other information and measurements that exist.  The primary axiomatic statements regarding conservation principles in continuum mechanics can be listed as follows.

\begin{enumerate}
    \item {\bf Conservation of mass.}  This is encapsulated in the (approximate) idea that matter can be neither created nor destroyed. In reality, it is possible to convert some matter to energy (and vice versa), but this is not a usual situation in continuum mechanics, so the conservation of mass can be assumed to be valid.  An example of a conservation of mass expression would be the continuity equation that is encountered in fluid mechanics. \indexme{conservation!axiom of conservation of mass}
    
    \item {\bf Conservation of momentum.} In continuum systems, the conservation of momentum is analogous to Newtons laws for discrete bodies.  In short, conservation of momentum for continuum systems states that Newton's laws apply not only to discrete bodies, but to any portion one can imagine being cut out of a discrete body.  Conservation of momentum is where we get the important equations describing the flow of fluids (the Navier-Stokes equations), or the forces distributed in a solid container under pressure.  Technically, there are two parts of the axiom of conservation of momentum: Conservation of linear momentum, and conservation of angular momentum.  Most of the material in this text will focus only on problems of  linear momentum. \index{conservation!axiom of conservation of momentum}
    
    \item {\bf Conservation of energy.}  Conservation of energy is a complementary axiom to of the conservation of mass, and is usually stated that energy can be neither created nor destroyed.  As mentioned above, however, we know that mass and energy can be interchanged under somewhat extreme conditions (e.g., when objects are moving very fast, are very small, or quantum mechanical interactions are relevant).  Generally, continuum mechanics does not consider such systems, so the conservation of energy can safely be taken to be a truthful axiomatic statement.  \indexme{conservation!axiom of conservation of energy}
\end{enumerate}
As a related note, each of the quantities above represents a thermodynamically \emph{extensive} quantity; that is, the value of the mass, momentum, or energy computed depends upon the physical size of the system.  Conversely, one may consider any of these three quantities on a per volume basis (i.e., the mass density, the momentum density, or the energy density), which makes these densities intensive quantities.   For the most part, in continuum mechanics it is best to discuss the dependent variables in terms of densities rather than extensive quantities.  For later use, we will use lower case letters (in particular, the letter $u$) to indicate intensive variables, and upper case letters (e.g., $U$) to indicate extensive variables.

Usually, the development of partial differential equations come from application of these conservation laws to systems that can be considered to be a \emph{continuum}.  A continuum is simply any region of space filled with a material in which can be assumed to behave as if it varied smoothly from point to point.  We know that, fundamentally, this cannot be true because all matter is composed of atoms, which are distinctly not continuous at the molecular scale.  However, when we consider large numbers of atoms together, their behavior can often be well approximated as if it were a continuous system.  In fact, this is the basic idea behind the field of continuum mechanics.  In continuum mechanics, a material is represented by \emph{representative volumes} of atoms that are numerous enough that they behave as if they were a continuum.  
A coherent example of a representative volume can is given by the application of the ideal gas law to noble gases.  Noble gases tend to obey the ideal gas law rather closely at near standard temperature and pressures.  Clearly, a gas is made up of a collection of individual atoms or molecules, but if there are a large enough number of them, then they behave as if they can be described by as a single bulk medium.  Hence, a large enough molar concentration of a noble gas at nearly earth-surface temperatures and pressures will be well described by the ideal gas law.  

However, if one could (at least conceptually) maintain a small number of gas atoms (say, $n=10$ atoms) at room temperature in, say, a 1 liter volume, the ideal gas law would \emph{fail} to predict the pressure very accurately.  In fact, the very definition of pressure (force per unit volume) would be challenging for a such a small number of atoms.  The problem is that the amount of momentum transfer between a small number of atoms and the walls of the container would be a rapidly varying function of time.  Because we think about the molecules as having more-or-less random thermal motion (i.e., they have kinetic energy that is proportional to the temperature), with only 10 atoms in a 1 liter container, the atoms collide with the walls in a way that appears somewhat random. 

To make the is clearer, suppose we fix a time interval $\Delta t$ such that, on average, 5 of the molecules contact the walls of the container over the time period.  Of course, with so few molecules involved, sometimes this number would be higher (say, 7 or 8 contact the wall in $\Delta t$) and sometimes lower (say 2 or 3 contact the wall in $\Delta t$).  Viewed this way, the number of atoms having collided with the walls of the container would be a random variable.  Thus, the pressure that one might measure would vary wildly from near zero (when no molecules contact the container walls) to a maximum value when all ten molecules contact the surface over the time interval $\Delta t$.  Clearly, this situation would lead to a pressure that varied in time, proportional to the number of molecules that contacted the surface over the time interval.  Therefore the ideal gas law would only be met in an average sense.  If the molecules are assumed to behave independently, then the variance of the number of collisions with the walls would be inversely proportional to the number of molecules involved, consistent with the Law of Large Numbers that you may have studied in statistics.  Thus, one can always drive down the variance of the behavior of our example gas by increasing the number of molecules.  With a very large number of molecules, the variance rapidly becomes incredibly small.  Under these conditions, the gas can be treated as a continuum, and the ideal gas law would hold with very high fidelity.

\section{Derivation of Differential Conservation Equations: An Introduction to Continuum Mechanics}

Assuming that a system can be treated as a continuum, then one can develop conservation equations on the basis of differential balances for the quantity of interest.  In this section, we will focus on generating a general balance equation for mass, linear momentum, or energy, regardless of which quantity is of interest.  Each of these quantities subscribe to the same fundamental axiom for conservation.  Following the discussion above, this axiom can be stated as follows.\indexme{conservation!axiom of conservation}

\begin{axiom}[Axiom of Conservation for Mass, Momentum, or Energy-- Form 1]

For any \emph{isolated system}, the amount of mass, momentum, or energy contained in that system is a constant.

\end{axiom}
Here, it is necessary to define the word \emph{isolated}.\index{system!isolated}\indexme{isolated system}  An \emph{ideal system} \indexme{system!ideal} is one that can be isolated from the rest of the universe using perfect boundaries that prevent the transport of mass, momentum, and energy.  Thus the terms \emph{isolated} and \emph{ideal} are synonymous when applied to systems.  While this makes for a reasonable modeling approximation in many cases, it should be noted that no real system can be ideal.

While the axiom of conservation of mass, momentum, and energy is probably familiar as presented, it is not necessarily the most useful form of the axiomatic statement.  The problem is that we are often interested in sub-volumes of an isolated system.  Fortunately, there are other statements of the axiom of conservation of mass, momentum, and energy that are equivalent to the version given above, but are more useful for applications to sub-volumes of isolated systems.  

\begin{axiom}[Axiom of Conservation for Mass, Momentum, or Energy-- Form 2]
 For any sub-volume cut from an isolated system, the following balance is maintained in accordance with  the laws of conservation of mass, momentum, and energy. 
For any volume, $V$ cut out of a continuum, the following balance law is valid for the extensive variable $U$ (where $U$ represents the mass, momentum, or energy in the volume $V$)
 
 \begin{equation}
  \underbrace{ \begin{Bmatrix}rate\;of\;\\accumulation\;\\of\;$U$\;in\;$V$\end{Bmatrix}}_{\numcircledtikz{1}}
  =
  \underbrace{\begin{Bmatrix}rate\;of\;$U$;\\entering\;the\;\\volume\;$V$\end{Bmatrix}}_{\numcircledtikz{2}}
  -
  \underbrace{\begin{Bmatrix}rate\;of\;$U$\;\\leaving\;the\;\\volume\;$V$\end{Bmatrix}}_{\numcircledtikz{3}}
  +
  \underbrace{\begin{Bmatrix}rate\;that\;sources\;or\;reactions\\\;change\;the\;amount\\\;of\;$U$\;in\;\;the\;volume\;$V$\end{Bmatrix}}_{\numcircledtikz{4}}
  \label{wordbal}
 \end{equation}
 
\end{axiom}

This form of the axiom of conservation is equivalent to the first, but it more useful for developing a differential balance.  The axiom gives us a direct way to setting up the differential balance; we need only convert our statement given in words into an equivalent statement using the language of mathematics.  

\subsection{Development of the Word Statement of the Conservation Axiom} \indexme{conservation!word statement}

To begin, we define a the sub-volume of the isolated system as illustrated in Fig.~\ref{repvol}; keep in mind that this volume is meant to be one that will become a differential volume in the limit as $\Delta x \rightarrow 0$.  In this introduction to computing balances, we will impose some constraints so that the problem can be considered to vary only in one spatial direction (the $x-$direction).  Thus, we make the following, rather strong assumption. 

\subsubsection*{\bf Assumptions for a 1-Dimensional Process Geometry.} \indexme{geometry!1-dimensional}\indexme{geometry!process geometry} 
All intensive properties of interest change only in the $x-$direction; thus, the intensive property of interest is constant over any cross-sectional area perpendicular to the $x-$direction.  We will further require that the cross-sectional area, $A$, be constant for 1-dimensional systems.  While all real systems exist in three spatial dimensions, the \emph{processes} involved can often be idealized as involving fewer dimensions.  To make this more concrete, assume that the cylinder illustrated in Fig.~\ref{repvol} is an insulated wire and that we are interested in heat transport though the system.  If the boundary conditions are uniform over the ends of the wire, then we can make reasonable symmetry arguments to suggest that transport of heat is primarily in the axial direction (that is, aligned with the $x$-axis).  For such systems, we can \emph{neglect by assumption} transport of heat in the $y$- and $z$-directions (that is, in the planes perpendicular to the $x$-axis.  So, while the \emph{physical geometry} of the system is 3-dimensional, the \emph{process geometry} for heat transport is for all practical purposes 1-dimensional for this system. \\

To start defining the terms for the mathematical statement of the balance, we note the following

\begin{enumerate}
    \item The intensive variable of interest (mass density (concentration), momentum density, or energy density) is denoted by the dependent variable $u(x,t)$.  The variable $u$ depends only on the independent variables $x$ and $t$. \indexme{variables!intensive}
    
    \item  The \emph{flux}, $j(x,t)$, is defined as the amount of mass, momentum, or energy passing through a unit area in a unit time.  Thus, the fluxes described here must have units of mass per unit area per unit time, momentum per unit area per unit time, or energy per unit area per unit time.  Because the property of interest has a 1-dimensional process geometry, then there are non-zero fluxes only in the $x$-direction.  The general expression above can be re-written for the 1-dimensional process geometry as follows

        \begin{align}
        \begin{Bmatrix}rate\;of\;\\U\;\\entering\;V\end{Bmatrix}
        &=\int_{y,z\in A(x)} {{j}_x}(x,t) \, dy \, dz \\
        &= A j_x(x)
    \end{align}
A similar expression is easily generated for the rate of $U$ leaving the volume $V$
        \begin{align}
        \begin{Bmatrix}rate\;of\;\\U\;\\leaving\;V\end{Bmatrix}
        &=\int_{y,z\in A(x+\Delta x)} {{j}_x}(x+\Delta x,t) \, dy \, dz \\
        &= A j_x(x+\Delta x)
    \end{align}
    
\begin{figure}[t]
\sidecaption[t]
\centering
\includegraphics[scale=.65]{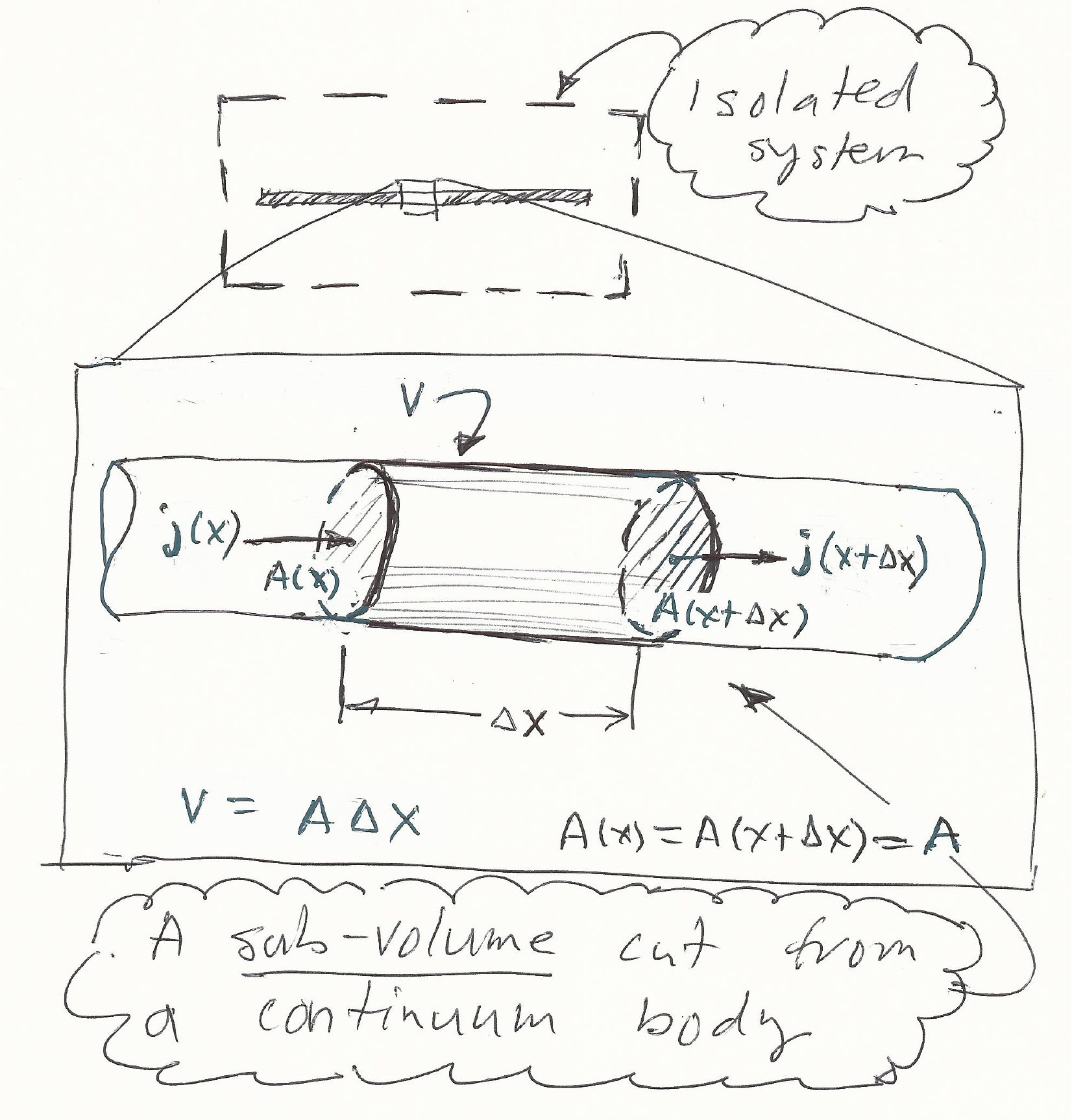}
\caption{{\bf A representative volume extracted from a continuum body.}  For our purposes, we will assume various symmetries that allow us to consider only one spatial dimension.  In particular, the system shown has a constant cross-sectional area, $A$, and variables of the system change only in the $x-$direction (thus, they are constant on every cross-sectional area.)}
\label{repvol}       
\end{figure}

\item The total amount of the extensive property within the volume $V$ is given by the integral of its density over the volume.  Density is the \emph{intensive} counterpart to  mass. Thus, integrating density over the control volume, $V$ generates the mass  within the control volume.  In other words

 \begin{align}
        \begin{Bmatrix}total\;\\amount\;of\;\\U\;\\in\; V\end{Bmatrix}
        &=\int_{x,y,z\in V(x)} u(x,y,z,t) dx\, dy \,dz \nonumber\\
        &= A \int_{x}^{x+\Delta x} u(x,t) \,dx 
    \end{align}
We are slightly abusing notation here; the variable of integration within the integral should, for absolute clarity, be different than the variable used for the bounds.  However, for this development, the abuse is not likely to confuse and is likely to improve understanding, so we will not stand on propriety and instead opt for clarity.  

\end{enumerate}

With these definitions in place, we are in a position to evaluate each of the terms denoted ${\numcircledtikz{1}}$ -- ${\numcircledtikz{4}}$ in Eq.~\eqref{wordbal}.  Each term is taken in sequence in the material following.

\subsection*{Item $\numcircledtikz{1}$ : The accumulation term.}

The accumulation term is defined as the ``rate of accumulation of $U$ in $V$."  Given that we have an expression for the \emph{total} amount of $U$ in $V$, rate can be found by taking the time derivative.   Note that because we have \emph{two independent variables} ($x$ and $t$), the derivatives must be defined as \emph{partial derivatives}.  We have as a result

 \begin{align}
        \begin{Bmatrix}rate\;of\;\\accumulation\;\\of\;U\;in\;V\end{Bmatrix}
        &= A\frac{\partial}{\partial t} \int_{x}^{x+\Delta x} u(x,t) \,dx 
    \end{align}
Note that, although our word form is in terms of the total amount of $U$ in $V$, this has been converted to the intensive quantity $u$ in the mathematical representation.  It should be clear from the steps above why this is the case!

\subsection*{Item $\numcircledtikz{2}$ : The rate of $U$ entering $V$}
The groundwork done above makes this somewhat simple to evaluate.  Recalling that we have a 1-dimensional process geometry, the rate of flux entering the volume $V$ is given by 

\begin{equation}
      \begin{Bmatrix}rate\;of\;$U$;\\entering\;the\;\\volume\;$V$\end{Bmatrix} = A j_x(x)
\end{equation}

\subsection*{Item $\numcircledtikz{3}$ : The rate of $U$ leaving $V$}
Similarly, we find that the rate that mass is \emph{leaving} the volume $V$ is

\begin{equation}
      \begin{Bmatrix}rate\;of\;$U$;\\entering\;the\;\\volume\;$V$\end{Bmatrix} = A j_x(x+\Delta x)
\end{equation}

\subsection*{Item $\numcircledtikz{4}$ : Source or Sink Term in $V$}
Finally, we consider the amount of internal reaction or other source or sink within the volume to be given by a function $s$. Recalling that we have assumed that the process geometry is 1-dimensional, then this means that the source/sink term can depend only on $x$;  the total amount of reaction or source or sink is then given by 
    
    \begin{align}
        \begin{Bmatrix}rate\;that\;sources\;or\;reactions\\\;change\;the\;amount\\\;of\;u\;in\;\;the\;volume\;V\end{Bmatrix}
        &=\int_{x,y,z\in V(x)} s(x,y,z,t) dx dy dz \nonumber\\
        &= A \int_{x}^{x+\Delta x} s(x,t) \,dx
    \end{align}
    
Now we are in a position to consider converting our word statement given above into a mathematical statement.  Substituting the mathematical statements of the word quantities developed above, the axiomatic statement given by Eq.~\ref{wordbal} is given in mathematical terms by

\begin{equation}
    A\frac{\partial}{\partial t} \int_{x}^{x+\Delta x} u(x,t) \,dx =
     A j(x)- A j(x+\Delta x) +
    A \int_{x}^{x+\Delta x} s(x,t) \,dx
\end{equation}
Clearly, the constant $A$ plays no role in this expression, so we can simplify this to

\begin{equation}
    \frac{\partial}{\partial t} \int_{x}^{x+\Delta x} u(x,t) \,dx = j(x)- j(x+\Delta x) +
    \int_{x}^{x+\Delta x} s(x,t) \,dx 
\end{equation}

At this juncture, this does not quite look like a differential balance equation.  However, there is one important piece of information that we have not yet used.  Because we are considering our system to be a \emph{continuum}, we are free to let the volume $V$ be as small as we like.  More specifically, for the one-dimensional representation, we can let $\Delta x \rightarrow 0$.  This will simplify our expression significantly.  To start, note that the variable of integration $x$ is such that $x_0 < x < x+\Delta x$.  Thus, as 
$\Delta x \rightarrow 0$, we must also have $|x -x_0| \rightarrow 0$.  Now consider the following two Taylor series expansions for $u(x,t)$ and $s(x,t)$ around the point $x_0$. These are given by

\begin{align}
    u(x,t) &= u(x_0,t) + (x -x_0)\left.\frac{\partial u}{\partial x}\right|_{(x_0,t)}+\ldots\\
    s(x,t) &= s(x_0,t) + (x -x_0)\left.\frac{\partial s}{\partial x}\right|_{(x_0,t)}+\ldots
\end{align}
Recall, we have the situation where $|x -x_0| \rightarrow 0$.  In Chapter 1, we noted that all derivatives of analytic functions must remain bounded within their domain.  This means that we can \emph{always} find a small enough value for $\Delta$ such that the second term in each of the Taylor series above is as small as we like. In the limit as $\Delta x$ tends toward zero (but is not zero), the Taylor series results in the following approximations

\begin{align}
    u(x,t) &= [u(x_0,t) + \alpha_1 \Delta x] \\
    s(x,t) &= [s(x_0,t) + \alpha_2\Delta x]
\end{align}
where the final term represents a measure of the error involved, with $\alpha_1$ and $\alpha_2$ being two finite constants.  This means the two integrals can be simplified as follows

\begin{align}
    \frac{\partial}{\partial t} \int_{x}^{x+\Delta x} [u(x_0,t) + \alpha_1 \Delta x] \,d\xi &= j(x,t)- j(x+\Delta x,t) +
    \int_{x}^{x+\Delta x} [s(x_0,t) + \alpha_2\Delta x] \,d\xi \nonumber\\ \vspace{2mm}\nonumber\\
    \frac{\partial}{\partial t} [u(x_0,t) + \alpha_1 \Delta x]\int_{x}^{x+\Delta x}  \,d\xi& = j(x,t)- j(x+\Delta x,t) +
   [s(x_0,t) + \alpha_2\Delta x] \int_{x}^{x+\Delta x}  \,d\xi \nonumber\\
   \vspace{2mm}\nonumber\\
   \frac{\partial}{\partial t} [u(x_0,t) + \alpha_1 \Delta x]\Delta x& = j(x,t)- j(x+\Delta x,t) +
   s(x_0,t) + \alpha_2\Delta x] \Delta x\nonumber \\
    \vspace{2mm}\nonumber\\
     \frac{\partial}{\partial t} [u(x_0,t) + \alpha_1 \Delta x]& = -\frac{[j(x+\Delta x,t)-j(x,t)]}{\Delta x}+
   [s(x_0,t) + \alpha_2\Delta x]
\end{align}
Note in the final line of this expression, we have rearranged the two flux terms so that the resulting expression is in the conventional form for the definition of the derivative.  Now, if we take the limit as $\Delta x\rightarrow 0$, note that $u(x_0,t)\rightarrow u(x,t)$.  We find the result known as the general conservation equation

\begin{equation}
\boxed{
        ~ ~~\frac{\partial u(x,t)}{\partial t}   = -\frac{\partial j(x,t)}{\partial x}+ s(x,t) ~~~
         }
         \label{conservationeq}
\end{equation}

\section{Constitutive Equations and Flux Laws}

It is not an understatement to say that that Eq.~\eqref{conservationeq} is the most important equation in continuum mechanics. In one very simple equation, a relatively universal statement about the conservation of mass, momentum, and energy is expressed.  Plus the statement is mathematically interpretable in terms of physical intuition: accumulation in a control volume is equal to the difference between mass in and mass out (the flux term) plus any generation or reaction within the volume.

One thing that the conservation equation \emph{does not} do is express conservation of $u$ in terms of only the dependent variable $u$.  Instead, we are faced with needing to know both the fluxes, $j$, and the source terms, $s$, in order to solve the conservation equation.  The source terms are generally straightforward: they are either specified by somewhat familiar terms representing kinetic reactions (or equivalent processes for energy or momentum), or they are nonhomogeneous terms (rendering the equation itself nonhomogeneous).  The flux terms, however, are somewhat different.  We know on the basis of the physics of the problem that the flux of $u$ should somehow involve the variable $u$ itself. The expressions that relate the flux $j$ to a function of the dependent variable, $u$ are called \emph{constitutive equations}.  They play a central roll in both the thermodynamics of nonequilibrium systems and in continuum mechanics (in fact, these two fields overlap significantly).  There are some deep results in the theory of constitutive equations, and they continue to be an active area of research.  In the material that follows, we will list some of the more important constitutive relationships from the perspective of mathematical modeling of continuum systems, but, aside from a few side notes, we will not dwell too heavily on the details of their genesis.  While many constitutive equations were initially determined empirically (and were often called 
\emph{laws} given that they were otherwise axiomatic statements at the time of their discovery), many of them have been given more fundamental interpretations though their formulation in statistical mechanics. 

It can be useful to think of PDEs as having characteristics of the three fundamental ``types" that we have discussed (parabolic, hyperbolic, and elliptic).  In the material that follows, we will combine the general conservation equation derived above (Eq.~\ref{conservationeq}) with various flux laws to obtain representative equations of each of these three types.  As an example of a \emph{hyperbolic equation}, we will examine the first-order ``pure convection" equation; such an equation would describe, for example, how a solute plume might move along a river.  For \emph{parabolic equations}, we will illustrate examples where fluxes are given by gradient laws.  Representatives illustrating the conservation of mass, momentum, and energy will be presented.  Parabolic PDEs are one of the most useful representations for describing the continuum mechanical behavior of quantities that are specified by gradient flux laws.  Finally, for \emph{elliptic equations} we will examine a reaction-diffusion problem at steady state.   

There is an important note to be made at this juncture.  The second-order wave equation (also a hyperbolic equation) is not one that arises naturally in the context of combining conservation and flux laws.  While it arises from considerations of \emph{conservation of momentum}, its derivation has some key differences, the foremost of these being that there is no simple flux law for the momentum in that case.  Instead, we will adopt the following approach.  The \emph{first-order wave equation} (a hyperbolic equation) is one that does, in fact, arise from the combination of Eq.~\ref{conservationeq} and a flux law; in this chapter we focus on deriving it.  A derivation of the second-order wave equation, using the more-or-less traditional approach, is given in the Appendix to this chapter.    Later, during the study of the Fourier transform, we will have some new tools that allow us to re-visit this problem.  In particular, we will be able to make a link between the first-order and second-order wave equations.

For each of the examples that appear below, for now we restrict the discussion to \emph{one space dimension} and one dimension representing time.  The spatial restriction will be relaxed in material that appears later on in the text where we consider the solution to multiple spatial dimensions.  

\subsection{First-Order (Hyperbolic) PDEs via Convective Fluxes}

The most fundamental kind of flux is referred to by either the term \emph{convective flux} or \emph{advective flux}.  There is some disagreement in the literature (and, occasionally, even debates) about the terminology, usually with strong statements that one or the other of these two terms is the ``correct" terminology followed by some historical fact presented as proof.   Regardless of these details, meaning of language is created via usage, and both terms are used to describe the same phenomenon widely enough that one should not attempt to distinguish between the two words unless the author has carefully defined the meaning.  From here forward, only the term \emph{convective flux} will be used.

\begin{figure}[t]
\sidecaption[t]
\centering
\includegraphics[scale=.6]{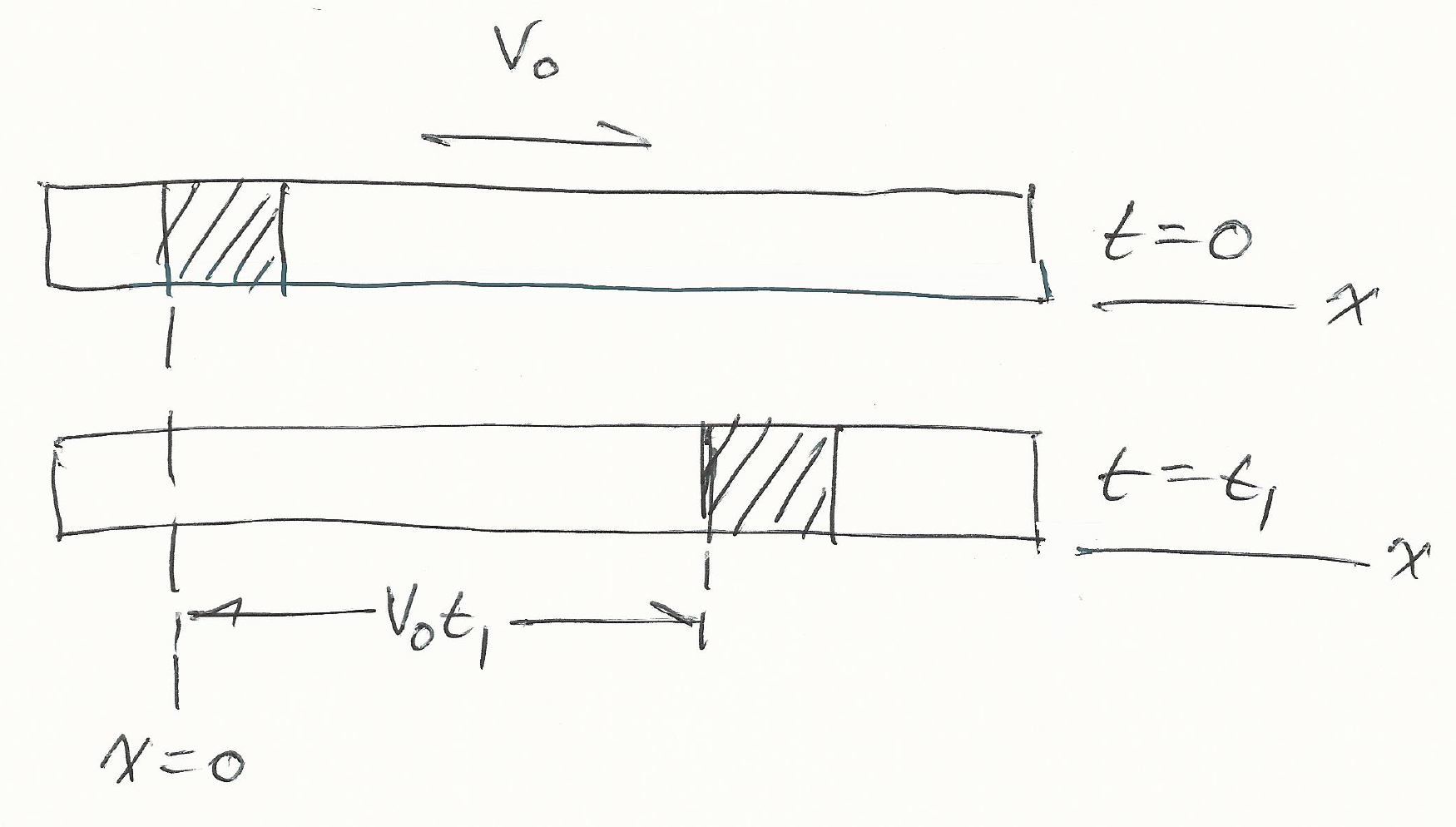}
\caption{{\bf Convective flux.}  The effect of a purely convective flux in an ideal plug flow reactor with a uniform velocity profile.}
\label{convectpfr}       
\end{figure}

A convective flux is one that transports the intensive quantity, $u$ by virtue of motion of the whole medium.  As an explicit example, if one were to release 1,000 rubber ducks on  the surface of a river, the motion of these ducks could be said to define the convective flux of ducks (assuming they were normalized by the appropriate definition of cross sectional area).  Another familiar example is the motion of a pulse of a chemical species in an ideal plug flow reactor, which is illustrated in Fig.~\ref{convectpfr}.  The idea is that the motion of $u$ follows the bulk velocity of the medium in which $u$ is embedded; obviously this makes the most sense if that medium is a fluid (gas or liquid), although there are forms of matter that are not clearly one or the other or are mixtures of matter types, but still exhibit bulk motion.  Magma, sand flows, and plasmas are examples of matter that can exhibit liquid-like convective fluxes.  Mathematically, a convective flux is stated by the relationship

\begin{equation}
    j_c(x,t) = u(x,t) v(x,t) 
\end{equation}
where $u$ is the dependent variable (concentration, momentum per unit volume, or specific energy), and $v$ is the velocity field that describes the motion of the material in which $u$ is embedded.

\begin{svgraybox}
\begin{example}[Mass conservation: The continuity equation. A hyperbolic equation.]\label{continuity0}

A simple and important example of convective fluxes combined with the conservation equation is the \emph{continuity0} equation.  For now, assume that we interpret $u$ as the density of a gas, $u \Leftrightarrow \rho$ in a steady flow field in one dimension, with no source or sink terms ($s=0$).  This system is composed of only one chemical species -- the fluid.  The convective flux is then $j_c(x,t) = u(x,t) v(x)$, and the conservation equation gives

\begin{equation}
    \frac{\partial u(x,t)}{\partial t}   = -\frac{\partial (u v) }{\partial x}
\end{equation}
\end{example}
\end{svgraybox}

\begin{svgraybox}
\begin{example}[Mass conservation: The first order wave equation. A hyperbolic equation.]\label{1stwave}

A first order wave equation occurs under conditions such as the ideal plug flow reactor illustrated in Fig.~\ref{convectpfr}.  For this example, we can think of $u$ as representing the concentration of a species dissolved in a fluid (liquid or gas) in steady flow in a uniform flow field.  Thus, the system is composed of at least two chemical species: one is the fluid, and the other the chemical species that is dissolved in the fluid.  It is assumed that there are no reactions or other sources, so that $s=0$.  The flux is given by the same relation as above: $j_c(x,t) = u(x,t) v(x)$.  The conservation equation gives

\begin{equation}
    \frac{\partial u(x,t)}{\partial t}   = -\frac{\partial (u v) }{\partial x}
\end{equation}
For many liquids, the concept of incompressibility is applicable.  While all matter is \emph{compressible} to some extent, liquids are often well approximated as being incompressible.  In short, this means simply that no more liquid can enter a control volume than that which exits the volume (i.e., there can be no liquid accumulation in a volume).  Another way of stating this is that the density of the liquid is a constant.  From the continuity equation above, we immediately find that $\partial v/\partial x = 0$ for an incompressible fluid.  This also implies, then, that the fluid velocity must be a constant, $v(x)=v_0$.  Under these conditions, the conservation equation can be written in the form

\begin{equation}
    \frac{\partial u(x,t)}{\partial t}   = - v_0\frac{\partial u}{\partial x}
\end{equation}
Which is a first-order wave equation.  Although we are not technically \emph{solving} PDEs at this juncture, a little thought will indicate that in 1D, if the initial condition is simply being translated through space at a constant velocity, $v_0$.  If the initial condition is given by the function $f(x)$, then the solution to the PDE is given by the translated version of the initial condition for time $t$, i.e., 

\begin{equation}
    u(x,t) = f(x-v_0 t)
\end{equation}

\end{example}
\end{svgraybox}

As a reminder, second-order hyperbolic equations do not arise readily from the combination of a conservation law and a flux law.  This is not to say that such equations are not important, only that their structure is distinct from the others in this chapter (which all arise from the the combination of Eq.~\eqref{conservationeq} with a flux law).  Second-order wave (hyperbolic) PDEs will be revisited later under the topic of Fourier transforms.  There, it will be possible to show a useful correspondence between first-order and second-order hyperbolic equations.  

\subsection{Second-Order (Parabolic) PDEs via Gradient Flux Laws}

Our general equation for conservation given by Eq.~\ref{conservationeq} can generate a wide array of interesting PDEs that describe many important physical processes.  An interesting symmetry exists between this conservation equation and flux laws.  Above, we saw how the convective flux law $j= u(x,t) v(x,t)$ generated a first-order wave equation that describes, for example, the motion of a solute in a plug flow reactor.  In this section, we examine a second category of flux laws-- these are the linear gradient laws.  Such flux laws are useful for describing the dispersive transport of mass, momentum, and energy under the appropriate conditions.

A gradient law is any constitutive expression that describes the flux as being proportional to the spatial gradient of $u$.  In one-dimension, the gradient is represented simply by the derivative, $\partial u/\partial x$.  As mentioned above, gradient laws have a long history in the field of continuum mechanics.  A few are given in the following.

\begin{align}
   & \textrm{Fick's Law (mass):} &j_g &= -D \frac{\partial u} {\partial x}&&u = \textrm{concentration} \\
    & \textrm{Fourier's Law (energy):} &j_g &= -K \frac{\partial u} {\partial x}&&u = \textrm{temperature}  \\
     & \textrm{Newton's Law of Viscosity (momentum):} &j_g &= -\mu \frac{\partial u} {\partial x}&&u = \textrm{velocity}
\end{align}
While gradient laws were in fact determined experimentally originally, in 1905 Albert Einstein illustrated one of the first \emph{derivations} of a gradient law from what was to become known as statistical mechanics.  Einstein was able to show that if large molecules could be conceived of as being represented by a very large number of non-interacting particles moving in random directions that changed direction over every fixed increment of time, then such motion naturally led to Fick's law as a continuum approximation to the motion.  While his derivations might be viewed as being somewhat rough by modern standards, one also has to remember that the very notion of the existence of molecules was not entirely accepted by the scientific community in 1905.  Not only did Einstein develop one of the first derivations of a constitutive equation from a more fundamental theory (statistical mechanics), but his work was one of the final indications that the molecular theory of matter must in fact be valid.  Einstein's theory was tested just a few years later (1908) by a brilliant French experimentalist named Jean-Baptiste Perrin.  After Perrin's work was published, there were virtually no longer any valid arguments that could be effectively waged against the molecular nature of matter. 

\begin{svgraybox}
\begin{example}[Mass conservation: Fick's second law.  A parabolic equation.]\label{1stpar}

Sometimes Fick's gradient law combined with the conservation equation is called ``Fick's second law".  Given the information above, the result is straightforward to develop.  Substituting the gradient law for mass for the flux, $j$, the conservation equation gives (again, assuming that $s=0$, so there are no reactions or sources in this formulation)

\begin{align}
    \frac{\partial u}{\partial t}   &= -\frac{\partial}{\partial x}\left(-D\frac{\partial  u}{\partial x}\right) \nonumber \\
    \intertext{or}
     \frac{\partial u}{\partial t}   &= D\frac{\partial^2 u  }{\partial x^2} 
     \label{diffusioneq}
\end{align}
Clearly, the use of Fourier's law in the conservation equation for energy will yield analogous results.  Thus, Eq.~\eqref{diffusioneq} is known both as the \emph{diffusion} equation and the \emph{heat} equation.

\end{example}
\end{svgraybox}

\begin{svgraybox}
\begin{example}[Momentum conservation: Burgers' Equation. Mixed parabolic and hyperbolic.]\label{burgersexample}

Burgers' equation is a one-dimensional analogue to the well-known Navier-Stokes equations describing fluid motion.  While Burgers' equation was initially derived as a tool for studying the Navier-Stokes equations (and, thus, one without particular physical relevance), it was later understood to also describe certain kinds of waves in fluids.  The flux for Burgers' equation consists of the combination of a convective flux for momentum, and Newton's gradient law for viscosity.  For conservation of momentum, the variable of interest is $q=\rho u$; in other words, it is the momentum per unit volume.  Thus, the flux takes the form 

\begin{equation}
    j = q - \mu \frac{\partial q}{\partial x} 
\end{equation}
Substituting this result into the general conservation equation (and assuming that $s=0$) yields

\begin{align}
    \frac{\partial q}{\partial t}   &= -\frac{\partial}{\partial x}\left(q - \mu \frac{\partial q}{\partial x} \right) \nonumber \\
    \intertext{or, upon substituting $q=\rho u$ and assuming constant density }
     \frac{\partial u}{\partial t}   &= -u\frac{\partial u}{\partial x} +\frac{\mu}{\rho}\frac{\partial^2 u  }{\partial x^2} 
     \label{burgerseq}
\end{align}
This is the first example of a conservation equation that contains both \emph{convective} and \emph{gradient} fluxes.

\end{example}
\end{svgraybox}

\begin{svgraybox}
\begin{example}[Mass conservation: Convection-diffusion-reaction.  Mixed parabolic and hyperbolic. ]\label{cdr}

The convection-diffusion-reaction equation is a familiar one science and engineering, and it has been used to describe everything from the motion of dissolved contaminants in groundwater to the spread of infectious diseases.  Technically, systems that have convective motion experience \emph{dispersion} which is the combination of spreading due to molecular diffusion plus spreading due to fluctuations (or velocity gradients) in the velocity field.  Regardless, both diffusion and dispersion have the same mathematical form (assuming that neither is a strong function of space), and mathematically they are identical problems.   The reaction term is expressed in the source or sink term, $s$, appearing in the conservation equation.  This is the first example to contain a non-zero source term.  For this effort, consider $u$ to be the concentration of a reacting chemical species, $v_0$ to be a spatially constant (uniform) one-dimensional velocity field, and the reaction rate to be first-order in $u$, that is $s= - k u$.  The resulting expression is very similar to the previous result for Burgers' equation.

\begin{align}
    j &= u v_0 - D \frac{\partial u}{\partial x} \\
    s &= - k u
\end{align}

Substituting this result into the general conservation equation yields
\begin{align}
    \frac{\partial u}{\partial t}   &= -\frac{\partial}{\partial x}\left(uv_0 - D \frac{\partial u}{\partial x} \right)-k u \nonumber \\
    \intertext{or, upon using the fact that $v_0$ is a constant }
     \frac{\partial u}{\partial t}   &= -v_0\frac{\partial u}{\partial x} +D\frac{\partial^2 u  }{\partial x^2} - k u
\end{align}
This is the first example of a conservation equation that contains both \emph{convective} and \emph{gradient} fluxes and a source term.

\end{example}
\end{svgraybox}
%
\begin{svgraybox}
\begin{example}[Mass conservation: The Groundwater Flow Equation in 1-dimension.  Elliptic. ]\label{laplace}

Examples of \emph{elliptic} problems with two spatial variables present the challenges noted above (specifically, the fluxes become vectors, which is a topic that will be covered in a separate chapter). Presenting examples of elliptic equations with one space dimension is also somewhat challenging.  However, the following is a reasonable example.  To start, recall the one-dimensional conservation equation on $x\in[0,L]$ of the form

\[
\frac{\partial u(x,t)}{\partial t}   = -\frac{\partial j(x,t)}{\partial x}+ s(x,t)
\]

In groundwater flow, Darcy's law is the gradient constitutive equation that relates hydraulic head, $h$ to the rate of groundwater flow per unit area $j$.  In one-dimension, the expression takes the form

\[ j = -{K} \frac{\partial u}{\partial x} \]

where here $u=$hydraulic head in meters, $K$ is the hydraulic conductivity in $m/s$.  We are ignoring a amplification term known as the storativity for ease in the presentation.  Substituting the flux into the general conservation equation yields

\[
\frac{\partial u(x,t)}{\partial t}   = K\frac{\partial^2 h(x,t)}{\partial x^2}+ s(x,t)
\]
While this result is technically a parabolic equation, we consider specifically the case where the groundwater flow is at steady state.  Steady state is defined by the state in which, for all locations $x$ in the domain, $\partial u/\partial t = 0$.  Thus, the balance equation is given by 

\[ 0 = K\frac{\partial^2 h(x)}{\partial x^2}+ s(x)\]
Technically, this equation is \emph{elliptic}.  More conventionally, this equation would appear with two spatial variables in the form of \emph{Poisson's} equation.

 \[ K_x \frac{\partial^2 u}{\partial x^2} + K_y \frac{\partial^2 u}{\partial y^2}= -s(x) \]

If we had conditions of symmetry on the boundaries of the $x-y$ domain, we can imagine conditions in which the flow in the $y$-direction is zero.  Under those conditions, we would recover an equation that is identical to the one-dimensional form that was derived above (even though the problem applied to a 2-dimensional domain).
\end{example}
\end{svgraybox}

Hopefully at this juncture, you can see that the simple derivation of the conservation equation is indeed a powerful tool.  Combined with the constitutive expressions for fluxes and representations for the source term, a variety of meaningful continuum mechanical expressions can be derived.

\section{Ancillary Conditions}

The problem of determining how many ancillary conditions are needed for a PDE is very much like the analogous problems for ODEs.  Here, for concreteness, we will continue to discuss problems with at most two independent variables, one for space and one for time.  The case of two independent variables in space adds a complication; it introduces fluxes as vector quantities.  Thus, a fuller analysis of the the elliptic problems discussed earlier in this chapter will be delayed until the chapter on multidimensional conservation equations.

\subsection{Order of the PDE and Determining Necessary Ancillary Conditions}

For problems in one space dimension, it is relatively straightforward to determine the necessary number of ancillary conditions.  For clarity, we consider the diffusion / heat equation in the form

\begin{align}
    \frac{\partial u} {\partial t} = \alpha^2 \frac{\partial^2 u}{\partial x^2}
\end{align}
where $\alpha^2$ is a parameter that is real and (necessarily) positive, and represents either the diffusion coefficient or the thermal conductivity.  The way to think about these problem is to consider how many \emph{integrations} are needed to eliminate the derivatives that appear.  Here we are using the word ``integrations" somewhat informally (although, when we explore separation of variables, we will find that the solutions are indeed found by integration in the conventional sense that we have used the term in the study of ODEs).

In this example, we have the following derivatives to resolve
\begin{enumerate}
    \item One time derivative of first order.
    \item One space derivative of second order.
\end{enumerate}
Thus, to resolve the problem, we must conduct one ``integration" in time, and two ``integrations" in space.  This generates a set of three unknown constants, one from each integration.  In reality, the unknowns may be whole \emph{functions} because we are dealing with two dimensions; while the unknown values are fixed in one of the two dimensions, they remain potentially free in the other dimension.  This will be more clear when we consider an example.

\subsection{Types of Ancillary Conditions}

There is an additional detail that needs some attention regarding the ancillary conditions.  Generally, we should think of ancillary conditions as a hierarchy of conditions which usually start with specifying the value of the function, and may (depending on the order of the equation) include specifying derivatives of the function.  In other words, suppose we consider the two ancillary conditions in space needed for the diffusion equation given above.  Before continuing on, we will make note of the three fundamental types of ancillary conditions \emph{for second-order PDEs composed of one space and one time variable}, which we assume to be represented by the dependent variable $u$, and independent variables $x$ and $t$ (thus, the function being specified by the ancillary conditions is $u(x,t)$).

\begin{enumerate}
    \item {\bf Specified value conditions.} In PDEs, specifying the value of the solution at its time-space boundaries is the most fundamental kind of ancillary condition.  As will be discussed below, for a problem to be well-posed it often must have some portion of its boundary for each independent variable identified by a specified value condition.  The convention in the study of PDEs is that space and time variables are treated somewhat differently (and, in fact, they are physically quite different-- we can think of traveling backwards or forwards in even one space dimension, but only forwards in time!)  Thus, there are separate names given to the two different kinds of specified value conditions.
    \begin{enumerate}
        \item {\bf Specified value in the space dimension}.\indexme{ancillary conditions!Dirichlet boundary condition}\indexme{ancillary conditions!specified boundary value}\indexme{ancillary conditions!first-type}  These conditions are called \emph{Dirichlet} or \emph{first-type} boundary conditions.  The term ``boundary" is used to emphasize that they are spatial conditions.  Thus, these conditions specify the value of the dependent variable, $u$, at one or more spatial locations for all values of the independent variable, $t$, defining the time domain.  This boundary condition is named after the German mathematician Johann Peter Gustav Lejeune Dirichlet (1805--1859).
        
        \item {\bf Specified value in the time dimension}.\indexme{ancillary conditions!initial condition}  When applied to the time dimension, the specified value condition is called an \emph{initial condition}.  Thus, and initial condition specifies the value of the dependent variable, $u$, at $t=0$ (or some other time considered to be the start of the clock associated with the problem) for all values of the independent variable, $x$, defining the spatial domain.  
    \end{enumerate}
    
    \item {\bf Specified derivative conditions.}  Specified derivative conditions are, in a sense, ``weaker" than specifying the value of the dependent variable.  Thus, specifying the derivative of the dependent variable gives the problem more freedom to adapt to the requirements imposed by the PDE itself.  In other words, one can meet a specified derivative condition at a boundary, but the value of the function at that boundary is otherwise unconstrained by the boundary condition itself.  Because of this, some care has to be taken when applying specified derivative conditions.  It is possible to generate problems that have no unique solution if only specified derivative conditions are imposed.  Like specified value conditions, the convention in the study of PDEs is that different terminologies are used for space and time conditions of this type.
    
        \begin{enumerate}
         \item {\bf Specified derivative in the space dimension.}\indexme{ancillary conditions!Neumann boundary condition}\indexme{ancillary conditions!specified boundary derivative}\indexme{ancillary conditions!second-type} These conditions are known as \emph{Neumann} or \emph{second-type} boundary conditions.  These conditions specify the derivative of the dependent variable, $u$, at one or more spatial locations for all values of the independent variable, $t$, defining the time domain.  The condition is named after the German mathematician Carl Gottfried Neumann (1832--1925). 
    
        \item {\bf Specified derivative in the time dimension.} \indexme{ancillary conditions!Cauchy condition}  When applied to the time dimension, the specified derivative condition is still referred to as an \emph{initial condition}; there is no special reason for this other than convention.  Thus, this initial condition specifies the derivative of the dependent variable, $u$, at $t=0$ (or some other time considered to be the start of the clock associated with the problem) for all values of the independent variable, $x$, defining the spatial domain.  In more mathematical contexts, the specified value plus specified derivative initial condition is sometimes also called a \emph{Cauchy condition} in honor of the French mathematician Baron Augustin-Louis Cauchy (1789--1857).
        \end{enumerate}
    
    \item {\bf Linear combinations of specified value and specified derivative conditions.}\indexme{ancillary conditions!Robin condition}\indexme{ancillary conditions!third-type}  There are instances where a linear combination of the two, more fundamental, types of ancillary conditions is useful and contains and often necessary physical content.  These kinds of conditions will not be brought up in any detail for the remainder of this chapter, but an example of a useful linear combination of ancillary conditions as a model for a physical boundary condition will be discussed in the material related to \emph{solving} PDEs.  When applied to spatial boundaries, the linear combination condition is called a \emph{Robin} condition after the French mathematician Victor Gustave Robin (1855--1897).
    
    \item {\bf Periodic conditions.} \indexme{ancillary conditions!periodic} There are some applications where the boundary conditions are also required to behave periodically on the domain.  This is not actually a separate boundary condition \emph{per se}, so much as a constraint on the kinds of solutions that can be attained.  The use of conditions to amend or constrain the boundary conditions is a topic that is more advanced than the material being covered in this text, so we will not discuss the in detail.
\end{enumerate}

\subsection{Applications of Ancillary Conditions with One Space and One Time Dimension }

Let's return to the diffusion problem introduced above. 
\begin{align}
    \frac{\partial u} {\partial t} = \alpha^2 \frac{\partial^2 u}{\partial x^2}
\end{align}
We know that we need three ancillary conditions total to account for the undetermined information (functions) created by integrating the equation. For the diffusion problem, the integrations involved necessarily require two \emph{boundary} conditions, and one \emph{initial} condition to be specified. 

We could specify these conditions by any of the following options

\begin{align*}
    \intertext{Option 1 (two Dirichlet conditions, one initial condition)}
    & u(0,t) = u_0(t) & u(L,t) = u_L(t) & & u(x,0)=f(x) \\
    \intertext{Option 2 (one Dirichlet, one Neumann condition, one initial condition)}
    & \frac{\partial u(0,t)}{\partial x}= g_0(t) & u(L,t)  = u_L(t) & & u(x,0)=f(x) \\
    \intertext{Option 3 (one Dirichlet, one Neumann condition, one initial condition)}
    & u(0,t) = u_0(t) & \frac{\partial u(L,t)}{\partial x}= g_L(t)& & u(x,0)=f(x)  
\end{align*}
In other words, we should generally work up through the sequence of functions in the order $u, u_x, u_{xx}, \cdots$ when specifying the ancillary conditions.  Another way of stating this is that, for problems in one space and one time dimension, both space and time should have at least one Dirichlet (specified value) condition specified at some point in the problem.\\

Missing from the list above is the possible combination 

\begin{align*}
    \intertext{Option 4 (two Neumann conditions)}
    & \frac{\partial u(0,t)}{\partial x}= g_0(t)  & \frac{\partial u(L,t)}{\partial x}= g_L(t)&  & u(x,0)=f(x) 
\end{align*}

The reason that this problem is not listed is that it requires special handling.  This is best shown by using a direct example.

\begin{svgraybox}
\begin{example}[The problem of two Neumann conditions. ]\label{2neumann}

Consider the elliptic problem above for groundwater flow.  Suppose we have the following PDE and two boundary conditions on the interval $x\in[0,L]$

\begin{align*}
    K_x \frac{\partial^2 u}{\partial x^2} &= 0 \\
    \frac{\partial u(0,t)}{\partial x} &= g_0 \\
    \frac{\partial u(L,t)}{\partial x} &= g_L\\
\end{align*} 
This problem is actually an easy one to solve.  It requires just two direct integrations of the first equation.  An outline of the process is as follows

\begin{align*}
    \frac{\partial}{\partial x}\left(\frac{\partial u}{\partial x}\right) &= 0~~\textrm{ Divide by $K_x \ne 0$; rewrite derivatives in equivalent form.}\\
     \frac{\partial u}{\partial x} &= C_1~~~~\textrm{  Integrate once, as indefinite integration}\\
       u(x) &= C_1 x + C_2~~~~\textrm{  Integrate a second time, as indefinite integration}\\
\end{align*}
Now, note that the derivative of the solution implies

\[ \frac{\partial u}{\partial x} = C_1 \]
Indicating that the derivative is, everywhere, a constant.  Thus, for the original statement to make sense, we must have $g_0=g_L=G$, for $G$ a constant.  In other words, the solution does not allow the flux at the two ends of the domain to be different from one another.  The solution at this juncture is 

\[ u(x) = G x + C_2 \]

Unfortunately, we have run out of unique boundary conditions with which to solve for $C_2$!  This situation may be repairable if one can generate some other constraint on the system (usually considering some physical conservation principle) to replace the lack of a second, unique boundary condition.  For example, if the average velocity were known, i.e., 

\[\frac{1}{L}\int_{0}^{L}u(x)\, dx = U_0\]
Then, this would be sufficient to replace the lack of information that we have by having a repeated Neumann condition.  Specifically, note

\begin{align*}
    \frac{1}{L}\int_{0}^{L}u(x)\, dx &= U_0 \\
    \frac{1}{L}\int_{0}^{L} (G x + C_2)\, dx &= U_0\\
    \frac{1}{L}\left( \frac{GL^2}{2} + C_2L\right) &= U_0 \\
    \implies C_2 &= U_0 -\frac{GL}{2}
\end{align*}

This is an indication that the application of two Neumann conditions must be considered very carefully in general to assure that the following are met:
\begin{enumerate}
    \item The problem must make \emph{physical} sense, including in the limits as time becomes arbitrarily large.
    \item The problem must have sufficient content so that the problem is \emph{well posed}.  This just means that the solution exists, the solution is unique, and the solution does not exhibit ``chaotic" behavior.
\end{enumerate}

\end{example}
\end{svgraybox}

An problem combining all of the ideas in this section is given in the following example illustrating a set of boundary and initial conditions for a reaction-diffusion equation that is \emph{well posed}.

\begin{svgraybox}
\begin{example}[A well-posed reaction-diffusion equation. ]\label{rxndiff0}

A \emph{well-posed} problem is one in which the following conditions are met.
\begin{enumerate}
    \item The problem has a solution.  This means that the PDE itself is constructed so as to be internally consistent (i.e., it makes mathematical or physical sense), and that it has sufficient ancillary information so that solutions may be found.
    
    \item The problem is unique.  This means that not only can we find solutions, but we can find particular solutions that represent specific conditions for which we would like solutions.  This means that we do not find an equivalence class of solutions (where each solution is valid, but no two solutions are linear combinations of the others), but rather a single solution that is fully specified by appropriate ancillary conditions.
    
    \item The problem does not behave chaotically.  This is usually stated as requiring that the problem depend \emph{continuously upon the initial data}, although that is not the most intuitive statement that can be made.  In short, this requirement indicated that if the ancillary conditions of the problem are perturbed by a small amount, the problem solution is changed in a small amount.  While this statement can be made mathematically more rigorous, the essential idea is sufficient for our purposes.
\end{enumerate}

Consider the parabolic problem for a diffusion-reaction system on the interval $x\in[0,L]$.

\begin{align*}
  \frac{\partial u}{\partial t}=  D\frac{\partial^2 u}{\partial x^2}-k u  \\
\end{align*} 
This problem has one spatial derivative of order 2, and one time derivative of order 1.  Thus, we need two ancillary conditions for the spatial derivatives, and one for the time derivative.  Now, consider the following statement

\begin{align*}
&&  \frac{\partial u}{\partial t}&=  D\frac{\partial^2 u}{\partial x^2}-k u  \\
&B.C.~1&  u(0,t) &= u_0 \\
&B.C.~2&  \frac{\partial{u(L,t)}}{\partial x} &= 0 \\
&I.C.~1&  u(x,0) &= f(x) 
\end{align*} 
According to the guidelines above, the boundary conditions include one Dirichlet condition, and one other condition (in this case a Neumann condition, although a second Dirichlet conditions would be acceptable also if it represented the physical context of the problem).  Thus, the spatial ancillary conditions appear to meet our requirements.  Note also that the ancillary conditions are technically functions.  Even though they are constant functions, both boundary conditions assign a particular value of $u(0,t)$ or $\partial u(L,t)/\partial x$ to each boundary for every possible value of $t$; thus, these two statements represent functions of time.  A little though indicates why this must be true; if the functions were not defined for all times of relevance, then there would exist some times for which we did not have the appropriate ancillary conditions, and the problem would no longer be well posed.  An analogous statement can be made about the initial condition. The initial condition indicates what the state of the system is at $t=0$ (or, in principle, any time, $t_0$, that was assigned to define the ``start" of the problem).  Of course, this initial conditions must state what happens for every location $x$ in the domain, or else there would exist locations for which no ancillary space data existed, and the solution in time could not be found (the problem, again, would not be well posed).  

\end{example}
\end{svgraybox}

\section{Transformations to Obtain Solutions to PDEs}

Transformations of PDEs are an application of the idea that ``if you have a problem you do not know how to solve, try to convert it into a problem you do know how to solve."  While, \emph{in principle} this is a great idea, one of the problems with the general use of transformations in PDEs it is difficult to state definite rules to find such transformations, especially for nonlinear problems. 

\subsection{Symmetries and Transformations}

Some transformations for linear and nonlinear problems can be learned by understanding what are called the differential symmetries of the problem (i.e., the symmetries associated with the infinitesimal geometry of the system).  While that topic is well beyond the scope of this text, the basic results for the approach were found by a German Amalie Emmy Noether (1882--1935).  Not only is it notable that Emmy Noether was a female mathematician at a time where there were very few, but her work led to some deep theories and even entire discplines in mathematics that did not exist before her work.  
In short, Noether’s theorem proclaims that every conservation law has a continuous \emph{symmetry}, and vice versa — for every continuous symmetry, there’s an associated conservation law.  

What is meant by a symmetry here can be roughly thought about in terms of our conventional notion of the word.  Take, for example, a square drawn of a piece of paper.  If the square is rotated by any integer multiple of $\tfrac{\pi}{2}$ radians, then it appears to be unchanged.  Thus, it has a rotational symmetry with respect to rotations of $\tfrac{\pi}{2}$ radians.  Similarly, a square is unchanged by reflections across its horizontal or vertical axes; thus, it has two kinds of reflection symmetry.  Many of the familiar and unfamiliar kinds of symmetries are discussed in a text by the famous physicist Herman Weyl \citep{weyl2015symmetry}.  While most familiar examples represent discrete symmetries (i.e., the rotational symmetry for a square happens for only certain angles of rotation), continuous symmetries are a bit more complex.  While the details are not important here, one can think of continuous symmetries as symmetries that exist for a sequence of infinitesimal actions (e.g., rotations, translations).  

Regardless, the importance of symmetries and of Emmy Noether's theorem have been critical to the study of PDEs, and to the study of various physical pheneomena (often expressed through PDEs).  A few deep results in quantum mechanics have even been discovered by looking for certain kinds of physical symmetries in the mathematical laws that describe it.  

\subsection{An Example Transformation}

For our purposes, we might best understand how transformations can lead to a simplified problem by an explicit example.  In the following example, we propose a continuous symmetry in time; that is, we find that if we exponentially re-scale the dependent variable in time, we capitalize on a particular kind of symmetry (a dilation-of-time symmetry) that puts the problem in a new, much simpler form.  

\begin{svgraybox}
\begin{example}[Transformation of a reaction-diffusion equation. ]\label{transdiff}

Consider a parabolic problem similar to the one described in the previous example for a diffusion-reaction system on the interval $x\in[0,1]$.  This problem meets our (somewhat limited) definition of being well-posed. 

\begin{align*}
&&  \frac{\partial u}{\partial t}&=  D\frac{\partial^2 u}{\partial x^2}-k u  \\
&B.C.~1&  u(0,t) &= 0 \\
&B.C.~2&  u(1,t) &= 0 \\
&I.C.~1&  u(x,0) &= f(x) 
\end{align*} 
In the next chapter, we will examine finding solutions to this problem using separation of variables.  Actually, the solution to the pure diffusion problem is well known, and is reasonably straightforward.  While the solution to the reaction-diffusion problem is possible using separation of variables, it does contain a number of complications that make the solution unwieldy.  

Thus, we consider the following well-known transformation by \citet{danckwerts1951} (who was a chemical engineer) for transforming this problem into a pure diffusion problem.  The transformation is accomplished by defining the following new dependent variable $w$

\[ u(x,t) = w(x,t) \exp(-kt) \]

To adopt this transformation, we simply substitute the right-hand side for $u$ in both the PDE and ancillary conditions.  The following results

\begin{align*}
&&  \frac{\partial}{\partial t}  (w(x,t)e^{-kt})&=  D\frac{\partial^2}{\partial x^2} (w(x,t) e^{-kt})-k (w e^{-kt}) \\
&B.C.~1&  w(0,t) e^{-kt} &= 0 \\
&B.C.~2&  w(1,t) e^{-kt} &= 0 \\
&I.C.~1&  w(x,0) e^{0} &= f(x) 
\end{align*}
Applying the product rule for derivatives and extracting the exponential from the space derivative yields.  

\begin{align*}
&&  \frac{\partial w(x,t)}{\partial t}e^{-kt}  -k w(x,t)e^{-kt}&=  e^{-kt}D\frac{\partial^2}{\partial x^2} w(x,t)-k w(x,t) e^{-kt} \\
&B.C.~1&  w(0,t) e^{-kt} &= 0 \\
&B.C.~2&  w(1,t) e^{-kt} &= 0 \\
&I.C.~1&  w(x,0) e^{0} &= f(x) 
\end{align*}

Noting that the exponentials are not zero, and that for the initial condition the exponential is identically 1, we have the final result

\begin{align*}
&&  \frac{\partial w(x,t)}{\partial t}&= D\frac{\partial^2 w(x,t)}{\partial x^2}  \\
&B.C.~1&  w(0,t)  &= 0 \\
&B.C.~2&  w(1,t)  &= 0 \\
&I.C.~1&  w(x,0) &= f(x) 
\end{align*}
Which is a diffusion-like equation.  
\end{example}
\end{svgraybox}

\section{Elliptic PDEs Arising as the Steady State of Parabolic PDEs}

So far, we have discussed hyperbolic (the first-order wave equation), and parabolic (the heat/diffusion equation) as they are developed from applying various flux laws.  We have not yet discussed \emph{elliptic} equations.  Here, we discuss how such equations often arise as the steady-state version of a parabolic equation.  Because we have restricted ourselves to one space and one time dimension in this chapter, the resulting equation is a somewhat simple elliptic problem (it ends up being an ODE rather than a PDE).  However, the important points to be made here are (1) how elliptic equations often arise as the steady-state version of a parabolic equation, and (2) the development of some understanding of how elliptic problems behave physically.    

The steady state for a transient problem (i.e., one that is a function of time, and the solution evolves over time from an initial configuration) is simply the state where the partial derivatives in time approach zero.  The word \emph{approach} is used here because for transient second-order PDEs, the steady state usually occurs as time becomes arbitrarily large (i.e., as $t\rightarrow \infty$); thus, it represents a \emph{limiting} behavior.

When transient problems are composed of one space and one time dimension, then the steady state solutions are identical in form to ordinary differential equations.  There are some philosophically delicate points here to be made about the use of partial versus total derivatives for such problems (technically, as $t\rightarrow \infty$, they are no longer functions of time, thus there is only one independent variable, so one can certainly argue that they are formally, in fact, ODEs in the limiting case), but we will not dwell on these.  Instead, we will treat such problems as being \emph{de facto} ODEs, while continuing to keep in mind that the equation itself came from a more general form.  An example of handling  the determination of a steady state solution is given below.

\begin{svgraybox}
\begin{example}[A steady-state reaction-diffusion equation. ]
Let's start with the reaction-diffusion problem given above.  Specifically, we have

\begin{align*}
&&  \frac{\partial u}{\partial t}&=  D\frac{\partial^2 u}{\partial x^2}-k u  \\
&B.C.~1&  u(0,t) &= u_0 \\
&B.C.~2&  \frac{\partial{u(L,t)}}{\partial x} &= 0 \\
&I.C.~1&  u(x,0) &= f(x) 
\end{align*} 
Now, as time becomes arbitrarily large, we expect this equation to reach some steady state.  To determine if this is, in fact, true, we can examine what happens if the time derivative is set to zero.  If a solution exists, then the problem has a steady state.  Specifically, we have the problem (after dividing the balance equation through by $D$).  

\begin{align*}
&&   \frac{d^2 u}{d x^2}-\frac{k}{D} u &=0 \\
&B.C.~1&  u(0,t) &= u_0 \\
&B.C.~2&  \frac{d{u(L,t)}}{d x}  &= 0 \\
\end{align*} 
Note that we no longer have need for the initial condition, because there is no integration in time required; the partial derivatives have also been replaced by ordinary derivatives, since time is no longer an independent variable in this equation.  

Solving this equation is the same as described in Chapter 2 on ODEs.  It is a homogeneous equation, and solving it requires only that we first determine which case it represents by computing the roots of the characteristic equation.  Noting $a=1$, $b=0$, and $c=-k/D$, we have the roots given by

\[
    s_1 = +\sqrt{k/D}~~~~~~~~~~~~~s_2 = -\sqrt{k/D}
\]
For ease in notation, set $\beta =|\sqrt{k/D}|$.  The solution is now

\[  u(x) = C_1 \exp(\beta x) + C_2 \exp(-\beta x)  \] 
\[  du(x)/dx = C_1 \beta \exp(\beta x) - C_2\beta \exp(-\beta x) \]
\end{example}
The second boundary condition gives us
\[ 0 = C_1 \beta \exp(\beta L) - C_2 \beta \exp(-\beta L) \]
and the first boundary condition
\[ u_0 = C_1  + C_2 \implies -C_2 = C_1 -u_0\]
Substituting this last relation into the second boundary condition expression gives
\[ C_1 \left[\exp(\beta L)+\exp(-\beta L) \right] = u_0 \exp(-\beta L)  \]
or, recalling the definition of the hyperbolic cosine $\cosh(x) = 1/2(e^x+e^{-x})$, we have
\[ C_1 = \frac{u_0 \exp(-\beta L)}{2 \cosh(\beta L)} \]
\[ C_2 = u_0-\frac{u_0 \exp(-\beta L)}{2 \cosh(\beta L)} \]

The final solution is

\[ u(x) = \frac{u_0 \exp(-\beta L)}{2 \cosh(\beta L)}\exp(\beta x)
+\left(u_0 -\frac{u_0 \exp(-\beta L)}{2 \cosh(\beta L)}\right)
\exp(-\beta x) \]

For reference, the solution to this problem is given in Fig.~\ref{diffrxn}.

{
\centering\fbox{\includegraphics[scale=.4]{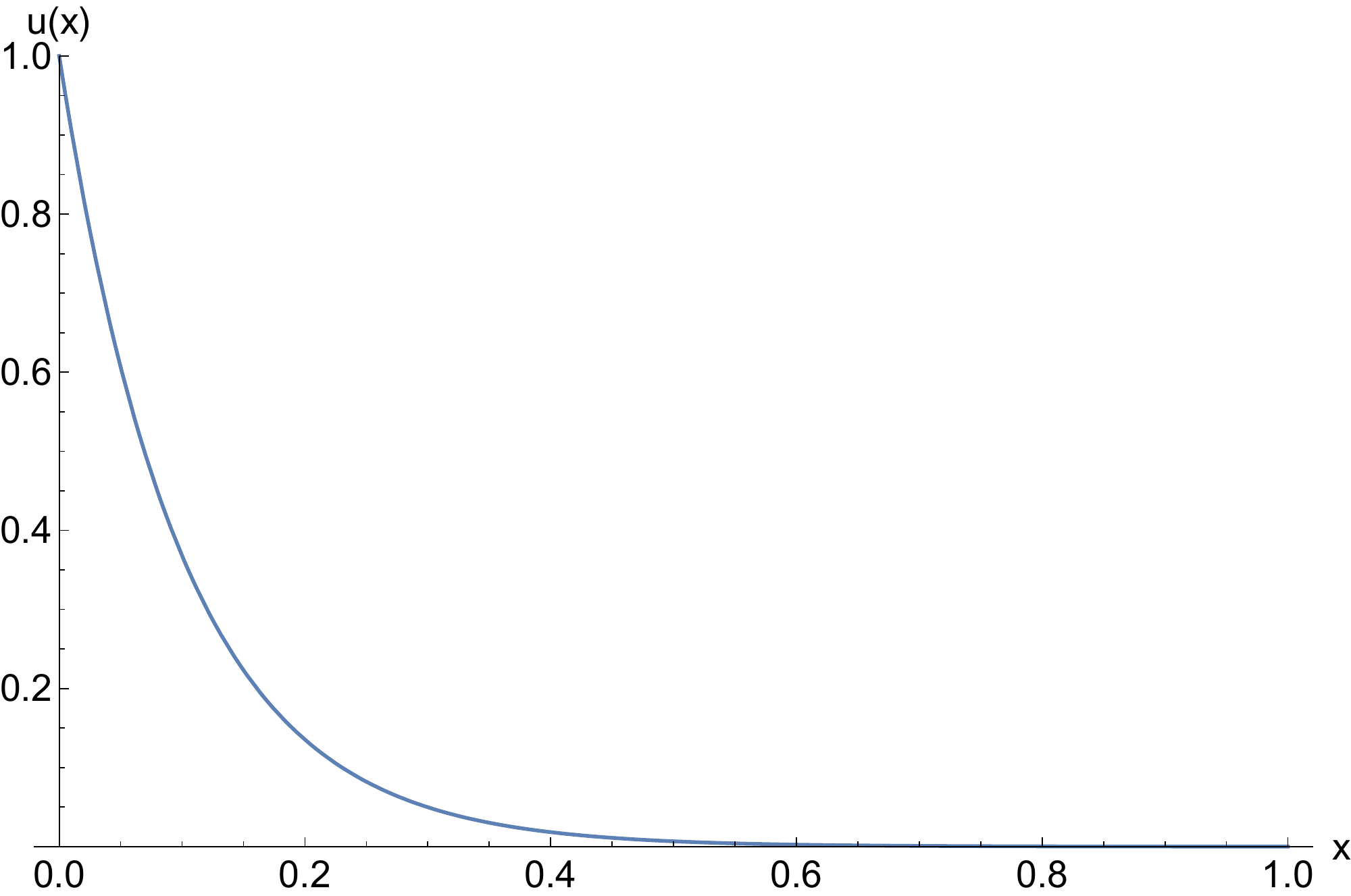} }
\vspace{2mm}
\captionof{figure}{The solution to the steady-state diffusion-reaction problem.}\label{diffrxn}
 
}

\end{svgraybox}

\section*{$^\star$Appendix: Derivation of the Wave Equation}
\addcontentsline{toc}{section}{~~~~~~~~$^\star$Appendix: Derivation of the Wave Equation}

There are a number of different \emph{kinds} of waves, but for our purposes it is sufficient to consider the following types
\begin{enumerate}
    \item \textbf{Longitudinal waves}. For example, the first-order wave equation describing a purely translating signal; or the compressive waves that comprise sound.
    \item \textbf{Transverse waves}. An example here might be waves on a body of water, where the vertical displacement from the equilibrium position defines the wave.
    \item \textbf{Combinations of Longitudinal and Transverse}.  One example of this kind of motion is observable in mechanical springs.  If you have ever played with a slinky \citep{gluck2010project}, you have probably observed both longitudinal (compressive) waves, and transverse waves.  When one periodically drives one end of a slinky, the toy spring develops a sine-function like shape, but also experiences regions of compression and extension.  Thus, both wave types are present.  
\end{enumerate}

\begin{figure}[t]
\sidecaption[t]
\centering
\includegraphics[scale=.35]{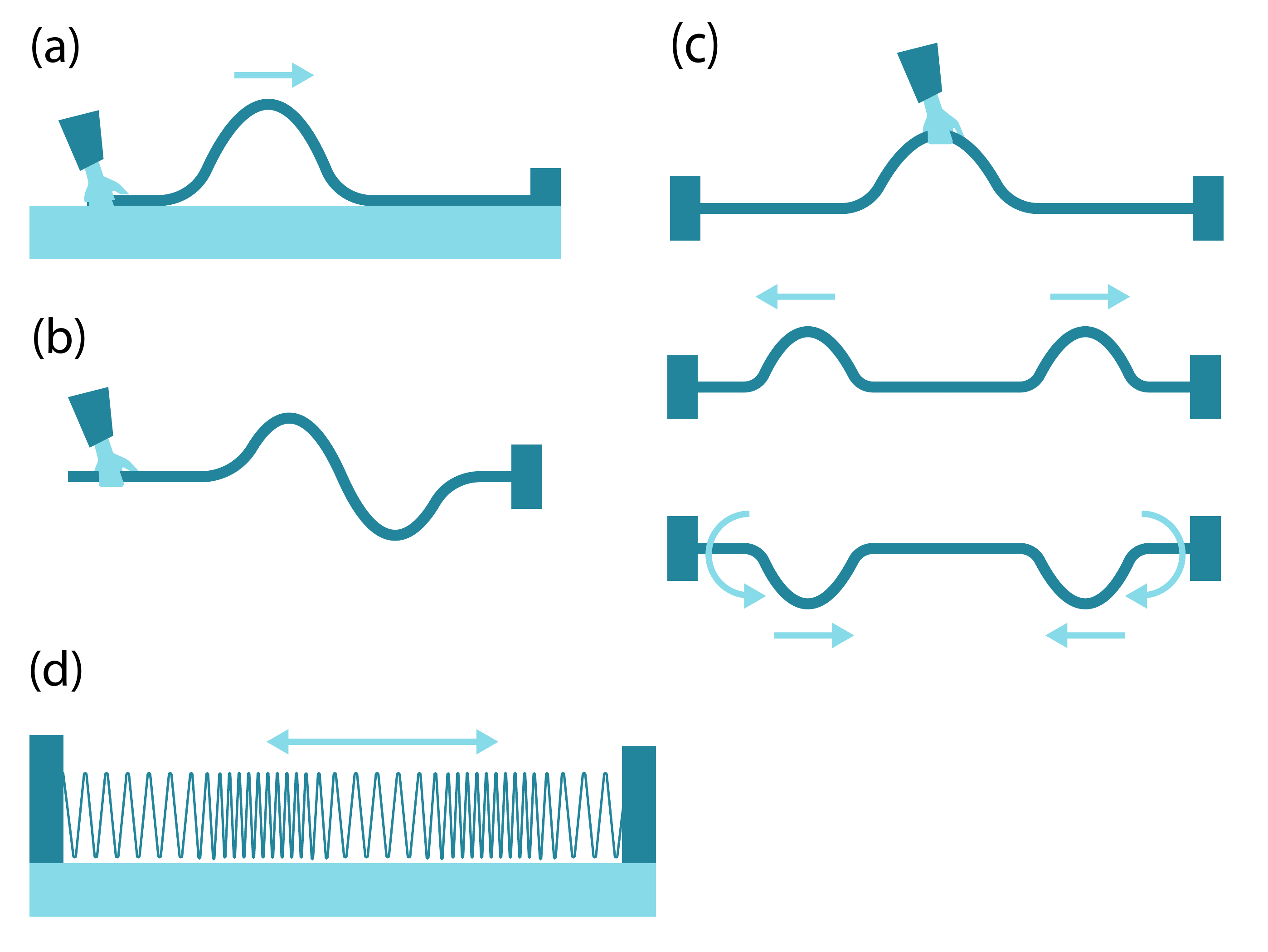}
\caption{Wave forms that can be observed in a slinky.  (a) A supported slinky disturbed to make a transverse wave that can be described by a first-order wave equation.  (b) A suspended slinky disturbed to create a transverse wave that requires a second-order wave equation to describe. (c) A suspended slinky with an initial condition created by perturbing the center; two waves are formed travelling in opposite directions. (d) A longitudinal wave created by compressing the center of a supported slinky; this also requires a second-order wave equation to describe. }
\label{fig:slinky}       
\end{figure}

In Fig.~\ref{fig:slinky}, several cases for waves are presented.  In the first (Fig.~\ref{fig:slinky}(a)) a slinky supported by the floor has only one direction to displace; that is, it can show displacements vertically above the floor.   A perturbation at one end of the slinky  will propagate along the slinky without changing shape.  This kind of wave is a primarily transverse wave that would be described by a \emph{first-order wave equation}.  As a second kind of wave, consider a slinky that is suspended and anchored at a wall. A suspended slinky can displace in both the positive and negative directions, as shown in Fig.~\ref{fig:slinky}(b).  As a different case, suppose the slinky is suspended between two walls (Fig.~\ref{fig:slinky}(c)). If one drives the toy spring with a sinusoidal initial displacement (i.e., the initial condition) in the middle, \emph{two waves are formed}, each of half the height of the initial displacement, and with the two waves traveling in opposite directions.  Upon contact with the wall, they are reflected, and their displacements are inverted.  This kind of motion is described by a second-order wave equation.  As a final mode of wave motion, we return to the case where the slinky is supported on the floor (Fig.~\ref{fig:slinky}(d)).  A longitudinal wave can be made by compressing part of the slinky and releasing it.  The longitudinal wave then travels by compressing and extending the slinky parallel to its center axis, but without creating vertical displacements. Done in the center of the slinky, again two waves will be formed, and these longitudinal waves will travel in both directions, and will also be described by a second-order wave equation.  Many motions of an actual slinky will be a combination of both transverse and longitudinal waves.

While the slinky is a toy, it is a very useful way to think about wave motions.  It also gives us a good way to begin thinking about these waves.  To start, recall that a first-order wave equation in an infinite medium is given by 

\begin{align}
    \frac{\partial u}{\partial t} &= -c  \frac{\partial u}{\partial x} ~~x\in[0,\infty) \label{eq:1wave}\\
    u(x,0) &= f(x)
\end{align}
Where here $u$ represents the transverse displacement, and  $a$ represents the width of an initial \emph{single} disturbance near the left boundary.  As an example, one might use the function $f(x)=\sin(\pi x)[H(x)-H(x-a)]$ to represent the initial condition of a \emph{single sine pulse} that travels along the system.  A solution to this problem can be computed by a number of methods, although we have not yet learned about these.  Here, we take a somewhat inverted approach by providing the solution, and then checking to validate it.  The solution to the problem above is

\begin{equation}
    u(x,t) = f(x-c t)
\end{equation}
Regardless of what function $f(x)$ is chosen.  In fact, we can codify this as a theorem, with following proof.

\begin{theorem}[Solution to the first-order wave equation with initial condition]
The solution to the first order wave equation with translation speed $c$ and initial condition $u(x,0)=f(x)$ is $u(x,t)=f(x-c t)$.
\end{theorem}

\noindent\textbf{Proof}.\\
The proof is not too difficult, but does require a change of variables, and appropriate differentiation of those variables.  Let

\begin{equation}
    y(x,t) = x- c t
\end{equation}
with this definition, note the following.
\begin{align}
    \frac{\partial y}{\partial t}& = -c \\
    \frac{\partial y}{\partial x} & = 1
\end{align}
Now, technically $y$ is a function.  However, using the normal rules of differentiation (the chain rule), we can compute derivatives with respect to $y$.  These derivatives can be computed as follows.

\begin{align}
    \frac{\partial u}{\partial t} &= \frac{\partial u}{\partial y} \frac{\partial y}{\partial t} \nonumber \\
    & = -c \frac{\partial u}{\partial y}\\
       \frac{\partial u}{\partial x} &= \frac{\partial u}{\partial y} \frac{\partial y}{\partial x} \nonumber \\
    & =  \frac{\partial u}{\partial y}
\end{align}
Now, substituting these results into Eq.~\eqref{eq:1wave}, we find that both sides are indeed equal, indicating that $f(x-c t)$ is indeed a solution.  Because the problem is a linear one of first order, we also know that the solution is unique. \qed\\

\begin{figure}[t]
\sidecaption[t]
\centering
\includegraphics[scale=.6]{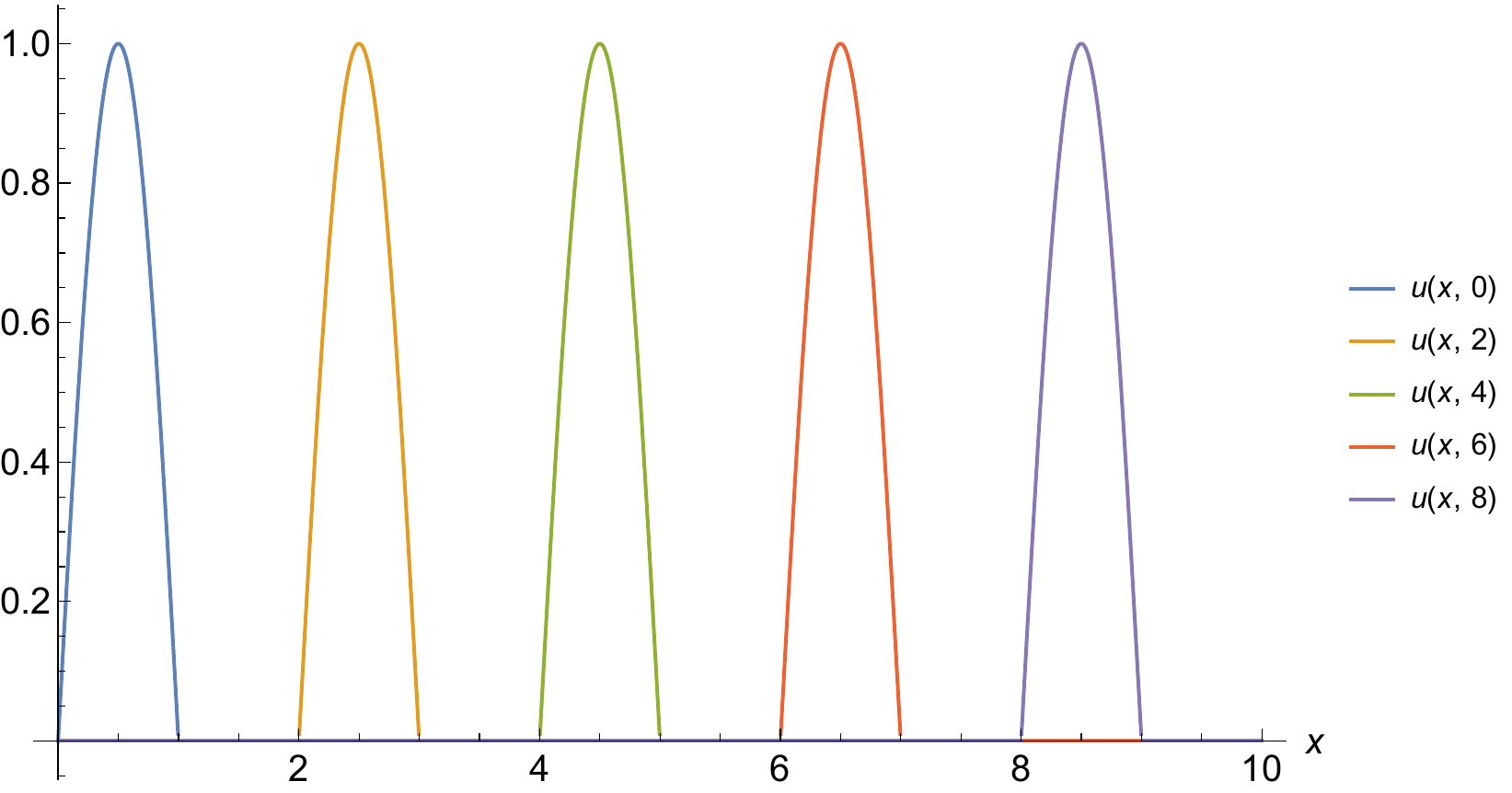}
\caption{Solution for a first-order wave equation, with $c=1~\si{m/s}$ and $a=1~\si{m}$. The transverse displacements, assumed to be non-dimensional, are translations of the function $f(x)=\sin(\pi x)[H(x)-H(x-1)]$ by the amount $c t$.  Times are at $t=0, 2, 4, 6,$ and $8$ seconds, and distance along the $x$-axis is in meters.}
\label{fig:firstwave}       
\end{figure}

The first-order wave equation describes the motion of a wave that translates along the $+x$-direction at velocity $c$, but does not change shape.  However, there are some deficiencies in this description.  First of all, it does not describe the physics of what happens if one disturbs the center of a medium rather than one end.  For that case (e.g., Fig.~\ref{fig:slinky}(b)), the translations propagate in both directions away from the initial disturbance.  Second, because the first-order wave equation involves only first-order derivatives in time, a single initial condition (the initial \emph{position}) is needed to determine any unknown constants of integration.  However, in some physical situations, we might want to specify both an initial configuration (shape) in space \emph{and} and initial velocity.  We cannot do this with a first-order wave equation.

One solution is to use what we know to think about what appropriate solutions might look like for the case where initial disturbances occur in the middle of a wave medium.  For now, let's consider only solutions where the initial velocity is identically zero, and the initial displacement, $f(x)$ is a symmetric function around $x=0$.    We expect waves to propagate in both directions ($+x$ and $-x$) from the disturbance.  Thinking about the conservation of energy, the potential energy of the initial disturbance must be converted into kinetic energy of the two waves.  Thus, each wave will have the same shape as the initial disturbance $f(x)$, but be only one-half the height (here, height is a proxy for mass).  From what we know about the first-order wave equation, we might propose a solution of the form

\begin{equation}
    u(x,t) = \tfrac{1}{2}[f(x-c t)] + \tfrac{1}{2}[f(x+c t)]
    \label{2ndwavesolution}
\end{equation}
This solution is a \emph{supposition} based on qualitative notions of the behavior of the first-order wave equation and a very rough constraint  arising from the conservation of energy (we suggested that height is proportional to energy, but we did not prove this in detail!).  Now the question arises-- ``if this represents a solution, to what PDE is it a solution?".  The proposition is this: ``This is a solution to the second-order wave equation."  Now, the problem is to show that this is true.  Again, we take a somewhat inverted approach by providing the solution, and then checking to validate it.  Below, we have a theorem and a proof that follows much along the lines of the proof for the first-order wave equation.

\begin{theorem}[The Second-Order Wave Equation]
    Suppose we have a wave-like solution in an 1-dimensional infinite medium with the form $u(x,t)=\tfrac{1}{2}[f(x-c t)] + \tfrac{1}{2}[f(x+c t)]$.   Here, $f(x)$ is the initial condition where $f(x)$ is a symmetric function around $x=0$, and it is assumed that the initial velocity is identically zero.  Then, the associated PDE describing this motion is given by 
    \begin{equation}
        \frac{\partial^2 u}{\partial t^2}=c^2 \frac{\partial^2 u}{\partial x^2}
        \label{eq:secondorderwave}
    \end{equation}
\end{theorem}

\noindent\textbf{Proof.}\\

The proof that the PDE given above (the classical ``wave equation") is indeed solved by $u(x,t)=\tfrac{1}{2}[f(x-c t)] + \tfrac{1}{2}[f(x+c t)]$ follows very closely the proof associated with the first-order wave equation.  Here, we have two translation variables, so we need two definitions as follows.

\begin{align}
    y(x,t) &= x- c t\\
    z(x,t) &= x+ c t
\end{align}

with this definition, note the following.
\begin{align}
    \frac{\partial y}{\partial t}& = -c,~~~~~~
    \frac{\partial y}{\partial x}  = 1\\ 
    \frac{\partial z}{\partial t}& = c,~~~~~~~~~
    \frac{\partial z}{\partial x}  = 1
\end{align}
Here, both $y$ and $z$ are a functions.  Again, we can use the normal rules of differentiation (the chain rule), to compute derivatives with respect to $y$ and $z$.  These derivatives can be computed as follows.

\begin{align}
    \dfrac{\partial u}{\partial t} &= \dfrac{\partial u}{\partial y} \dfrac{\partial y}{\partial t}+\dfrac{\partial u}{\partial z} \dfrac{\partial z}{\partial t} = \Cline{-c \dfrac{\partial u}{\partial y}+c\dfrac{\partial u}{\partial z}}\\
       \dfrac{\partial u}{\partial x} &= \dfrac{\partial u}{\partial y} \dfrac{\partial y}{\partial x}+\dfrac{\partial u}{\partial z} \dfrac{\partial z}{\partial x} =  \Cline{~~~~\dfrac{\partial u}{\partial y}+\dfrac{\partial u}{\partial z}}\\
        \dfrac{\partial^2 u}{\partial t^2} 
        &=  \dfrac{\partial }{\partial t}\left( -c \dfrac{\partial u}{\partial y}+c\dfrac{\partial u}{\partial z}\right) 
        =  -c\dfrac{\partial }{\partial y}\left( \dfrac{\partial u}{\partial t}\right) 
        =  -c\dfrac{\partial }{\partial y}\left( \dfrac{\partial u}{\partial t}\right) +c\dfrac{\partial }{\partial z}\left( \dfrac{\partial u}{\partial t}\right) = \Cline{c^2\dfrac{\partial^2 u}{\partial y^2}+c^2\dfrac{\partial^2 u}{\partial z^2}}\\
           \dfrac{\partial^2 u}{\partial x^2} 
        &=  \dfrac{\partial }{\partial t}\left(  \dfrac{\partial u}{\partial y}+c\dfrac{\partial u}{\partial z}\right) 
        =  \dfrac{\partial }{\partial y}\left( \dfrac{\partial u}{\partial t}\right) 
        =  \dfrac{\partial }{\partial y}\left( \dfrac{\partial u}{\partial t}\right) +\dfrac{\partial }{\partial z}\left( \dfrac{\partial u}{\partial t}\right) = \Cline{\dfrac{\partial^2 u}{\partial y^2}+\dfrac{\partial^2 u}{\partial z^2}}
\end{align}
Finally, substituting the underlined results for $\tfrac{\partial^2 u}{\partial t^2}$ and $\tfrac{\partial^2 u}{\partial x^2}$ into Eq.~\eqref{eq:secondorderwave}, we find that the second-order wave equation is indeed consistent with the proposed solution. \qed\\

This development was not a derivation of the wave equation per se.  It was actually (1) the proposal of a solution, and then (2) the illustration that the solution was \emph{consistent} with possible solutions to the second-order (classical) wave equation.  Recall that it was predicated by several requirements (symmetry, zero initial velocity).  It turns out that the second-order wave equation allows even more general solutions than the one proposed.  Because there are two derivatives in time, both and initial condition and initial \emph{velocity} are allowed as ancillary conditions that determine the two unknown constants of integration in time.  It is also not necessary that the initial wave by symmetric; non-symmetric initial conditions are also allowed by the second-order wave equation.  

Conventional derivations of the second-order wave equation (directly from the conservation of momentum) are relatively easy to come by-- most introductory texts on PDEs include some kind of a derivation.  These derivations, however, are at best, poorly motivated, and often impose many assumptions that are either incorrect or unclear.  The careful derivation of wave equations based on sounds physics is available in the literature \citep{antman1980equations,yong2006strings}.  For those interested, these two references are a good source for a sound derivation.  While the derivation of the second-order wave equation itself is not particularly difficult, it uses concepts that are so different from those presented elsewhere in this chapter that the derivation represents a substantial tangent from the main material in this text. Thus, the derivation of the wave equation from first principles is purposefully excluded.  However, we will return to this problem for more exploration in our studies of Fourier transforms.  There, some additional insight as to the nature of the second-order wave equation can be found. 

Some comments on the restrictions, on the \emph{kinds} of functions $f$ that can serve as solutions to the first- or second-order wave equation are needed.   To this point, we have not mentioned specifically if there are restrictions.  By the nature of the PDE itself, one might assume that the functions must have at least  continuous first-order derivatives (in both $x$ and $t$) for the first-order wave equation, and second-order derivatives (in both $x$ and $t$) for the second-order wave equation (i.e., the function $f$ must be of the differentiability class $C^1$ or  $C^2$). While this is certainly a \emph{sufficient} condition for solutions to the wave, it is not a necessary one.  While we cannot explore the details in this text, it is still worth noting that the solutions to the wave equation can be interpreted in the sense of the \emph{generalized functions} which have been mentioned previously in Chap.~\ref{deltachap}.  When viewed this way, it turns out that the function $f$ can even be a function with discontinuities (e.g., such as the \emph{boxcar} or \emph{step} function defined by $B(x)=H(x-a)-H(x-b), ~b>a$).  

\newpage
\section*{Problems}



\subsection*{Applied and More Challenging Problems}

\begin{enumerate}

\item The conventional diffusion equation takes the form (as derived earlier the chapter) 

\begin{equation}
 \frac{\partial u}{\partial t} = D \frac{\partial^2 u}{\partial x^2}
 \label{diffusioneq2}
\end{equation}
This derivation assumed that the diffusion coefficient was a constant, and did not depend on space.  Assume now that we have the case that the diffusion constant depends upon space; that is, it is given by a function of $x$ in the form $D(x)$.  Starting by re-stating Fick's law (i.e., the diffusive mass flux constitutive equation) for this situation.  Then, substitute this result into the general conservation equation.  How does the result for the case of a spatially-dependent $D(x)$ change from the result given by Eq.~\eqref{diffusioneq2}?  Note, this last question is asking for more of a qualitative description rather that an extensive discussion.\\

\item Find the steady-state solution to the following first-order hyperbolic problem on the domain $x\in[0,L]$.

\begin{align*}
 &&   \frac{\partial u}{\partial t} &= -v_0\frac{\partial u}{\partial x} - k u \\
 &B.C.~1& u(0,t) &= u_0 \\  
 &I.C.~1 & u(x,0)& = 0
\end{align*}
Does the solution to this problem match your intuition?  This last question is a qualitative one, so a qualitative answer is all that is expected.

\item  A theme in this text is one in which, when faced with a problem that one does not know how to solve, change the problem into one that can be solved.  While this is easy to state, it is not always easy to know how to transform such problems, however! \\

One interesting transformation that was invented by the famous chemical engineer Peter Danckwerts (1916--1984) is one that transforms a diffusion-reaction problem into one of only diffusion, greatly simplifying the solution process.  To start, consider the following problem with Neumann boundary conditions.

\begin{align}
 &\hspace{-20mm}&   \frac{\partial u}{\partial t} &= D \frac{\partial^2 u}{\partial x^2} - k u \label{eqdiffusereact}\\
 &\hspace{-20mm}B.C.1 & \left. \frac{\partial u}{\partial x}\right|_{(0,t)} &= 0\\
 &\hspace{-20mm}B.C.2 & \left. \frac{\partial u}{\partial x}\right|_{(L,t)} &= 0 \\
 &\hspace{-20mm} I.C. & u(x,0) &= f(x) \label{diffusereactic}
\end{align}

Now, consider the following transformation of variables

\[ u(x,t) = \exp(-k t) w(x,t) \]

Show that substituting this transformation into the diffusion-reaction problem given by Eqs.~(\ref{eqdiffusereact})-(\ref{eqdiffusereact}) will transform that set of equations into the following equations in terms of the variable $w$

\begin{align}
&\hspace{-10mm}& \frac{\partial w}{\partial t} &= D \frac{\partial^2 w}{\partial x^2}  \\
&\hspace{-10mm}B.C.1 & \left. \frac{\partial w}{\partial x}\right|_{(0,t)} &= 0\\
&\hspace{-10mm}B.C.2 & \left. \frac{\partial w}{\partial x}\right|_{(L,t)} &= 0 \\
&\hspace{-10mm}I.C. & w(x,0)&= f(x) 
\end{align}

\item Burgers' equation is a nonlinear equation that is used as a prototype for the Navier-Stokes equation in 1-dimension, and in applications to dissipative wave transport.  The equaiton is nonlinear, taking the form

\begin{equation}
u_t + u u_x - \nu u_{xx} = 0
\end{equation}
Nonlinear equations are generally difficult to solve analytically (if they can be solved at all!)  For this problem, there is a famous two-step transformation that makes the equation linear.  This transformation, called the Cole-Hopf transformation begins as follows.  First, we define the following relationship among variables.

\[ u = w_x \]
Now, substituting this into the original equation, we find an equation in $w$ of the form

\[ w_{xt} + w_x w_{xx} - \nu w_{xxx}=0 \]

It is not clear at this juncture that anything has been made better, but we will press on.  First note the following identity

\[ \frac{\partial} {\partial x}\left(\frac{1}{2} w^2_x \right) = w_x w_{xx} \]

Note that we have mixed notation for derivatives here, which is a perfectly acceptable thing to do as long as you keep everything straight.  This result is easy to check by a simple application of the chain rule.  It is also useful to recall that compositions of derivatives combine as follows

\[ \frac{\partial}{\partial x} w_x = \frac{\partial}{\partial x} \left(\frac{\partial w}{\partial x}\right) = \frac{\partial^2 w }{\partial x^2} = w_{xx} \]\\

Now, this means that the transformed Burgers' equation can be rewritten as follows

\[ w_{xt} + \frac{\partial}{\partial x}\left(\frac{1}{2} w^2_x\right)  - \nu w_{xxx}=0 \]
All of this work has paid off a bit, because now we can extract one partial derivative with respect to $x$ outside the expression as a whole, giving

\[\frac{\partial}{\partial x}\left( w_{t} + \frac{1}{2} w^2_x  - \nu w_{xx}\right)=0 \]

Both sides of this can be integrated once with respect to $x$ using an indefinite integration.  The result is 

\[ w_{t} + \frac{1}{2} w^2_x  - \nu w_{xx}=C_1 \]
While this equation seems perhaps \emph{better} than the original one, it is still nonlinear because of the term $1/2 w^2_x$.  Now, for the second transformation, which is the basis for this problem.  Show that the transformation of variables given by 

\[ w(x,t) = -2\nu \ln[z(x,t)]  \]

Will lead to a \emph{linear} second-order PDE of the diffusion equation form, i.e., 

\[ z_t = \nu z_{xx} \]
Where we have assumed that $C_1=0$.

\item The Cahn--Hilliard (CH) equation is a balance equation that predicts how two phases in a fluid will separate starting from a highly mixed initial configuration. In this case, $u$ is actually a dependent variable that encodes the internal \emph{configuration} of the geometry of the two phases.   The variable $u$ is initially viewed as a variable that takes on the values $u \in [0,1]$.  The value of $u$, then, indicates one of the two phases.  A value of 1 indicates the first phase, and a value of 0 indicates the second phase.  Variables that take on values of 1 or 0 to represent phases are sometimes called \emph{phase indicator} variables in the study of continuum mechanics.\\

The nonlinear balance equation for the phase-indicator variable $u$ takes the form
    
    \begin{equation}
        \frac{\partial u}{\partial t} = \frac{\partial^2}{\partial x^2}\left(Du^3 -Du -\gamma 
        D\frac{\partial^2 u}{\partial x^2} \right)
        \label{CH}
    \end{equation}

If you are interested in the topic, an example of the application of the CH equation in one dimension is reported in the following paper.\\
 \begin{quote}
\bibentry{argentina2005}
\end{quote}
The behavior of the equation is very interesting.  Starting from a (relatively) uniform initial condition representing a fully mixed fluid, the CH equation predicts how the fluid will separate.  As an analogue, think about how an oil-water emulsion can form, over time, large oil droplets as the small droplets coalesce (this may be an observation you have had in something as mundane as vinegar-and-oil salad dressing...)  While the applications of the CH are most interesting in 2- and 3-dimensional systems, there are 1-dimensional solutions also.  As an example, Fig.~\ref{chfig} provides the time-space value of the indicator function for one such solution.

\begin{figure}[t]
\sidecaption[t]
\centering
\includegraphics[scale=.35]{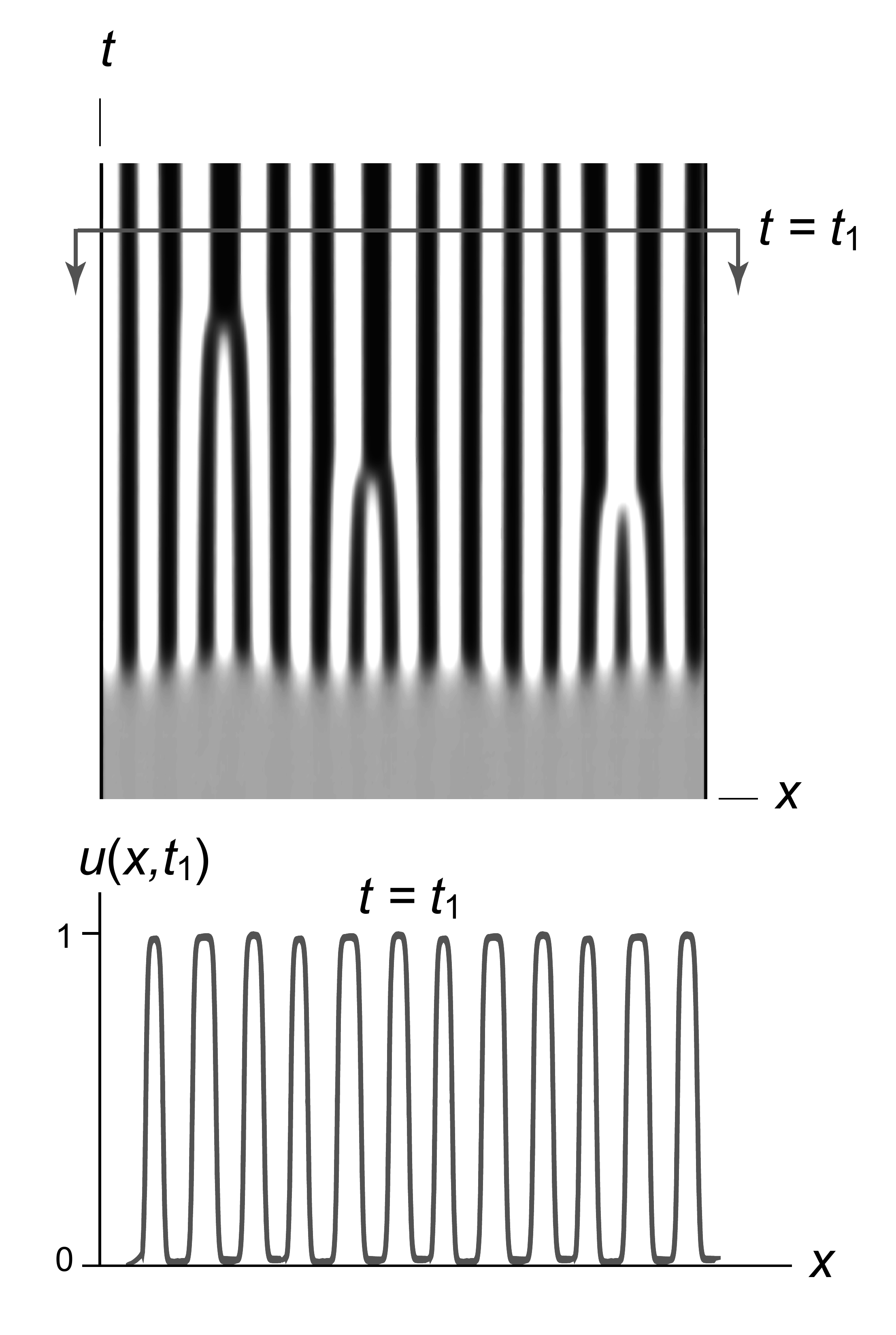}
\caption{{\bf Solution to the Cahn-Hilliard Equation.} This equation transform an initial mixture of phases (denoted by the gray region near the time origin) to one in which the two phases separate.  Note that this is a case illustrating the complex behavior known as ``emergent behavior" (see Chapter 1).  The upper figure shows the value of $u(x,t)$ as a function of both $x$ and $t$; the values in the domain represent an intensity plot (white=0, black=1, with gray representing values in between) for each time space location.  The bottom figure illustrates the spatial variation of $u(x,t_1)$ for a fixed value of time; it corresponds to the horizontal line labeled ``$t=t_1$" in the upper figure.  It is remarkable that this nonlinear equation generates such extreme behavior.  Figures based on data presented by \citet{argentina2005}. }
\label{chfig}       
\end{figure}

Recall, the general conservation equation (with no source) for 1-dimension takes the form

\begin{equation}
\frac{\partial u(x,t)}{\partial t}   = -\frac{\partial }{\partial x}j(x,t)
\label{simplecons}
\end{equation}

Comparing Eq.~\eqref{simplecons} with the expression given by Eq.~\eqref{CH}, determine what the flux, $j(x,t)$ must be as a function of $u$.  To be clear-- you need only compare these two equations to determine what the the flux must be interpreted as in Eq.~\eqref{CH}; this involves a bit of algebra and application of the derivative operators, but nothing else.  To start, think about how you can pull out the operator $-\partial/ \partial x$ from the right-hand side of Eq.~\eqref{CH}.  Once you do that, what this operator acts upon will then be the flux, $j(x,t)$.

\item Consider a parabolic problem, similar to the one described in example \ref{transdiff}, on the interval $x\in[0,1]$.  This problem is slightly different from the example because the reaction rate coefficient is a function of time.  The balance equation, boundary, and initial conditions are 

\begin{align*}
&&  \frac{\partial u}{\partial t}&=  D\frac{\partial^2 u}{\partial x^2}-k(t) u  \\
&B.C.~1&  u(0,t) &= 0 \\
&B.C.~2&  u(1,t) &= 0 \\
&I.C.~1&  u(x,0) &= f(x) 
\end{align*} 
For this problem, show that the transformation

\[ u(x,t) = \exp\left(-\int_{\tau=0}^{\tau = t} k(\tau) \, d\tau \right)  \]

can be used to transform the reaction-diffusion problem into a purely diffusion problem.

\item On the interval $x\in(-\infty, \infty)$ the heat/diffusion equation

\begin{align*}
    \frac{\partial u}{\partial t} &= D \frac{\partial^2 u}{\partial x^2}\\
    u(x,0) &= \delta(x)
\end{align*}
has a solution known as the \emph{fundamental solution} of the heat/diffusion equation.  This solution is 

\begin{equation*}
    u(x,t) = \frac{1}{\sqrt{4 \pi D t}}\exp\left(-\frac{x^2}{4 D t}\right) 
\end{equation*}
Taking derivatives and substituting into the heat/diffusion PDE, show that this is indeed a solution to the PDE.

\item We have discussed first-order wave equations. In the Appendix we have shown that for the first-order wave equation

\begin{align*}
    \frac{\partial u}{\partial t} &= -c \frac{\partial u}{\partial x} \\
    u(x,0) = f(x)
\end{align*}
has the general solution $u(x,t)=f(x-c t)$.  That derivation involved a change of variables, and was slightly complex.  For this problem, show that the for the specific initial condition $f(x)=\exp(-x^2)$, that the solution $u(x,t) = \exp[-(x-c t)^2]$ is a solution to the first-order wave equation, and that it meets the proposed initial condition.  Do this \emph{directly} by taking derivatives of the given form for $u(x,t)$.

\item We have discussed first-order wave equations, and, in the Appendix, the relationship between first- and second-order wave equations.  In that derivation, we showed that the second-order wave equation
\begin{align*}
    \frac{\partial^2 u}{\partial t^2} &= -c^2 \frac{\partial^2 u}{\partial x^2} \\
    u(x,0) &= f(x)\\
\end{align*}
has the general solution $u(x,t)=\tfrac{1}{2}f(x-c t)+\tfrac{1}{2}f(x+c t)$.  That derivation involved a change of variables, and was moderately complex.  For this problem, show that the for the specific initial condition $f(x)=\exp(-x^2)$, that the solution $u(x,t) = \tfrac{1}{2}\exp[-(x-c t)^2]+\tfrac{1}{2}\exp[-(x+c t)^2]$ is a solution to the seconds-order wave equation, and that it meets the proposed initial condition.  Do this \emph{directly} by taking derivatives of the given form for $u(x,t)$.
\end{enumerate}
\abstract*{This is the abstract for chapter 00}

\begin{savequote}[0.55\linewidth]
``The basic relations between Lie groups, special
functions; and the method of separation of variables have recently been
clarified... 
The main ideas relating the symmetry group of a linear partial differential
equation and the coordinate systems in which the equation admits
separable solutions are most easily understood through examples."

\qauthor{ Willard Miller, Mathematician who made advances on separation of variables as a general technique.}
\end{savequote}

\chapter{Separation of Variables (SOV)}
%
\def\CHAP {chapter06_separation_of_variables}
%
\section{Introduction}

The method of separation of variables (SOV) \indexme{separation of variables} has been studied for hundreds of years, with the first known result given by \citet{liouville1846quelques}.  While the idea is seemingly simple, it has actually been studies intensively up to the present day in an effort to better understand the conditions under which the method provides solutions \citep{bluman2010applications, kalnins2018separation}.  It turns out that the method of separation of variables is tied in with the structure (more specifically, the differential geometry) of the underlying partial differential equations and their boundary conditions.  The description of \emph{when and why} the separation of variables method works relies on understanding a wide array of mathematical concepts ranging from topology to algebraic structures known as groups.  Exploration of these ideas as applied to partial differential equations has led to some of the most celebrated results in applied mathematics.  Interested readers can find a good introduction to the topic in the text by \citet{miller1977symmetry}, who also provides some history of the topic.

All of this may make it sound like the method is difficult.  In fact, the method itself poses no particular challenges in applying.  In fact, it was originally developed by mathematicians such as Liouville (and later, John Bernoulli) who used physical motivation and observations from what they knew of existing solutions to PDEs to propose the method empirically.  In other words, it was easy to see that, for many interesting cases, the method of separation of variables could successfully generate solutions.  It is only when attempting to understand fully as to what conditions the method works (and why it does when it is successful) that one encounters the more advanced notions of symmetries and differential geometry.  

For our purposes, we can take a more empirical tactic with the method.  We will \emph{propose} that the method be tried, and then \emph{check} to see if this is indeed the case.  There is nothing inherently unscientific or wrong with such an empirical approach!  In fact, it is used very frequently in the study of mathematics.  Recall, in our study of second order ODEs, we proposed that useful solutions to a homogeneous ODE of second order with constant coefficients could be represented by $u(x) = \exp(-s x)$; this proposition, of course, led ultimately to developing all possible solutions for such second order ODEs.   Thus, this will be the approach here, and, as we will ultimately determine, the method works for a number of PDEs on domains that are sufficiently regular in some sense.  

As is frequently the case when first learning new material, it is best to begin to understand the approach by starting out with an example.  In fact, we will start by illustrating the solution for the three canonical equations (the parabolic diffusion equation, the hyperbolic wave equation, and the elliptic Laplace equation) in two variables.  After we see that it \emph{can} work in at least the case of this example, we can begin to push our understanding of the method a bit further.

Before continuing, note it is not easy to give \emph{general conditions} for separability of PDEs, even for linear ones.  However, for the three canonical examples (the parabolic diffusion equation, the hyperbolic wave equation, and the elliptic Laplace equation) separation it will always be possible to use separation of variables as a solution approach if \indexme{separation of variables!conditions for application}

\begin{enumerate}
    \item The equations have constant coefficients
    \item The equations are homogeneous
    \item The equations contain no mixed derivatives
    \item The boundary conditions are linear, homogeneous, and are aligned with the coordinate system axes
    \item The initial condition is defined by a function that is assumed to have a valid Fourier series representation
\end{enumerate}
Note that these conditions are not \emph{necessary} ones for separation of variables to be used.  However, for our purposes, they are \emph{sufficient ones}, meaning that if these conditions are met, then separation of variables can be used adopting the approach described in this chapter.  This also implies (by the definition of \emph{necessary} and \emph{sufficient}) that there may be (and, in fact, there are!) solutions that can be found via SOV where not all of the four conditions above are met.

In the remainder of this chapter, we will explore the general solutions for the three canonical equations (parabolic heat/diffusion equations, hyperbolic wave equations, and the elliptic Laplace-like equations).  For each of these, our focus will be on systems with only two independent variables (one space, one time for both the parabolic and hyperbolic equations; two space variables for the Laplace-like equations), and for homogeneous specified value (Dirichlet) boundary conditions.  An important component of this analysis will be the continuation of building skills in the \emph{modeling} component of the problems that we address.  In particular, we will show how the use of various constraints based on the \emph{physical} (rather than \emph{mathematical}) context will help us to develop solutions that represent the specific quantities that are of interest to us.  Once we have developed the basic ideas of SOV for the three canonical equations, we will turn out attention to cases that involve other kinds of boundary conditions.  Boundary conditions are an essential part of mathematical modeling, and they deserve some special attention to more fully understand how the various kinds of boundary conditions affect the models that we develop.

\section{Terminology}

\begin{itemize}
\item \textbf{Separation of variables.} A solution process for some linear PDEs.  For a PDE in $n$ independent variables, ${\bf x} = x_i,~i=1\ldots n$, it is assumed that the solution is given by $u({\bf{x}})=f_1(x_1)f_2(x_2)\cdots f_3(x_3)$.\\  

\item \textbf{Transient.}  The word \emph{transient} in PDEs is used to indicate an process that changes in time. As an example, consider a drop of dye added to a transparent container of water.  Starting from this initial state (the drop), the dye continues to spread in time; in more formal terms, if $u({\bf x},t)$ is the dye concentration for all locations at all possible times, then initially $\tfrac{\partial u}{\partial t} \neq 0$.\\

\item \textbf{Steady state.}  The converse of transient.  A process that is at steady state does not change in time.  If $u({\bf x},t)$ represents the dependent variable of interest, then  at all locations one has $\tfrac{\partial u}{\partial t} = 0$.  It is possible for a transient process to \emph{reach} steady state.  In most of our models, this will require $t\rightarrow \infty$, but there are processes modeled by PDEs that can reach steady state in finite time.  For example, the diffusion-reaction problem with a zero-order reaction (meaning that the reaction rate is independent of the concentration), the process goes to completion in some finite time.  Many transient problems tend toward a \emph{steady state} as time grows arbitrarily large; this will depend on the particular form of the PDE involved.  As an example, many heat/diffusion equations have a well-defined steady state.  Conversely, many wave-type equations have no steady state as we have defined it. \\

\item \textbf{The heat/diffusion equation}.  One of the three classical equations (heat/diffusion, wave equation, Laplace equation) usually studied in an introduction to PDEs.  It is a \emph{parabolic} equation, and describes phenomena that have \emph{diffusive} characteristics.  In one space and one time dimension, $(x,t)$ it takes the form
\begin{equation}
\frac{\partial u}{\partial t} = \alpha^2 \frac{\partial^2 u}{\partial x^2} 
\end{equation}
The fact that the coefficient is $\alpha^2$ is to emphasize that it \emph{must} be a positive quantity.  However, this is usually understood, and this formalism is often not maintained in practice.  Unlike the wave equation (below), $\alpha^2$ has physical significance (it is the heat or mass diffusion constant), but $\alpha$ (not squared) does not have a direct physical interpretation.\\

\item \textbf{The wave equation}.  Another of the classical equations studied in introductory PDEs. It is a hyperbolic equation, and represents phenomena that have \emph{wave like} behavior.  There are actually two wave equations studied in introductory PDEs.  The first is the first-order wave equation.  In one space and one time dimension, $(x,t)$ it takes the form  
\begin{equation}
    \frac{\partial u}{\partial t} = -v_0 \frac{\partial u}{\partial x}
\end{equation}
This equation describes the translation of an initial disturbance (initial condition) rectilinearly at velocity $v_0$.  The initial condition translates without changing shape.

The second form of the wave equation is one that is more commonly thought of as \emph{the} wave equation.  In one space and one time dimension, $(x,t)$ it takes the form  
\begin{equation}
     \frac{\partial^2 u}{\partial t^2} = c^2 \frac{\partial^2 u}{\partial x^2}
\end{equation}
The fact that the coefficient is $c^2$ is to emphasize that it \emph{must} be a positive quantity.  Strangely, unlike the heat/diffusion equation, the squared term is usually maintained in practice.  In part this is because of the physical significance of the quantity $c$, which can be interpreted as a measure of the wave velocity.\\

\item \textbf{The Laplace equation}.  The third of the classical PDEs studied in introductory courses, the Laplace equation is the the archetype for \emph{elliptic} equations in two or more independent variables (spatial dimensions).  Often, the steady-state for a heat/diffusion type equation results in a Laplace equation.  The conventional form for the Laplace equation in two spatial dimensions is

\begin{equation}
   K \frac{\partial^2 u}{\partial x^2} + K \frac{\partial^2 u}{\partial x^2} = 0
\end{equation}
The Laplace equation is named after 
Pierre-Simon, marquis de Laplace (23 March 1749 --5 March 1827) a famous French scholar who contributed to an enormous number of fields including engineering, mathematics, statistics, physics, astronomy, and philosophy.
If the right-hand side of the Laplace equation is specified by some function of the independent variables, then the equation is called the Poisson equation.  This equation is named after  Baron Sim\'eon Denis Poisson (21 June 1781 – 25 April 1840).  It turns out that Poisson was an ardent student of Laplace, and Laplace considered Poisson to be practically like a son to him.\\

\item \textbf{Harmonic function.}  Solutions to the Laplace equation are called \emph{harmonic functions}.  Harmonic functions have some unique properties that make them worthy of study in their own right.\\

\item \textbf{Principle of superposition.}  Linear problems have a unique status in mathematics-- because of the linearity, one can sum two solutions to problems, and the result is another solution.  This can be exceptionally useful, for example, for breaking a complex problem up into two simpler problems, where the sum of the two simpler problems is equivalent to the more complex one.  Then, the solution to the two simpler problems can be combined lineally to generate a solution that is a valid solution to the more complex problem.   This will be discussed more explicitly toward the end of the chapter.
\end{itemize}

\section{The Heat/Diffusion Equation}\indexme{second-order PDEs!heat/diffusion equation}

For the first example, we will consider the ubiquitous example of the heat/diffusion equation.  The parabolic heat/diffusion equation is possibly one of the most useful equations in engineering and physics, and it can be used to describe an enormous array of phenomena ranging from (obviously!) diffusion of mass to the spreading of momentum in a fluid by viscous shear (Newton's law of viscosity combined with the conservation of momentum equation).  In one space dimension, the balance of mass, $u(x,t)$, due only to (linear) diffusion takes the form

\begin{equation}
    \frac{\partial u}{\partial t} = \alpha^2 \frac{\partial^2 u}{\partial x^2}, ~~x\in[0,\ell]
\end{equation}
Here, we have used $\alpha^2$ as the coefficient to remind us that this equation is a heat/diffusion equation only for cases where the coefficient is \emph{positive}.  While one can formulate an ``backward" heat/diffusion equation, such an equation has some technical challenges, not the least of which is that it violates the second law of thermodynamics!  (It is interesting to note, however, that the backward heat equation does have uses in inverse modeling of data).  The domain length is $\ell$, although later we will use $\ell=1$ to simplify the initial examples.
 
As discussed in the previous chapter, this is an equation of second order, and there are two independent variables.  Thus, we will need two ``integrations" in space, and one ``integration" in time to resolve the problem.  This, then, corresponds to the requirement of supplying two boundary (fixed space) conditions, and one initial (fixed time) condition in order to specify a unique solution.  For this example, we will investigate solutions on $x\in [0,1]$ and for $t\ge 0$, with the time being measured initially from $t=0$.   The ancillary conditions will be homogeneous Dirichlet conditions at either end of the domain, and general function for the initial condition is specified.  On physical grounds, we impose two additional constraints on the problem: (1) The solutions for $u(x,t)$ must remain bounded for all time, and (2) all solutions must represent a finite and constant amount of total mass (which we impose by requiring the slightly stronger condition that the Fourier series of the initial condition exists).  The whole problem, including the physical constraints, can then be written out as


\begin{align}
&&    \frac{\partial u}{\partial t} &= \alpha^2 \frac{\partial^2 u}{\partial x^2}, \textrm{  for $x\in[0,\ell]$} \\
&B.C.~1& u(0,t) &= 0 \label{bc1}\\
&B.C.~2& u(\ell,t) &= 0 \label{bc2}\\
&I.C.~1& u(0,x) &= \varphi(x), \textrm{  for $x\in[0,1]$} \\
&Constraint~1 & &\substack{\hspace{-9mm}\textrm{\normalsize Solutions stay bounded as}\\  \hbox{\hspace{-35mm} \fontsize{10}{25}\selectfont\( t\rightarrow \infty\)}}\\
&Constraint~2 & &\substack{\hspace{-9mm}\textrm{\normalsize The initial condition has a valid} \vspace{1mm}\\  \hbox{\hspace{-14mm} Fourier series representation}}
\end{align}
Here ``B.C." stands for \emph{boundary condition}, and ``I.C." stands for \emph{initial condition}.  Note, the last condition represents a constraint imposed by the \emph{physics} of the problem.  We do not expect a diffusion problem to grow arbitrarily large; in fact, we expect the opposite.  Diffusion-like problems tend to spread mass out in space, creating a mass distribution that is as uniform as possible in the long-time limit (while still respecting the boundary conditions imposed).  While it may not be obvious that this constraint is needed, it is certainly a constraint that matches our understanding of the physical processes involved.  Later, we will find that this constraint allows us to extract the physically relevant solutions from the set of all mathematically possible solutions.  Mathematically possible solutions do not always represent solutions that are physically relevant!

Before continuing to the separation of variables method, we will take a short aside to determine a \emph{physical} property of heat equations that will be useful later on.  This aside is a property of heat equations that is general to heat equations with homogeneous boundary conditions.

\subsection{A Property of the Heat/Diffusion Equation: Decreasing Variance}\label{variance}
\indexme{heat/diffusion equation!decreasing variance property}
In thermodynamics, we learn that the entropy of a sequence of states for isolated system must always increase (or not change).  We can consider a 1-space dimension heat problem to be isolated in the sense that zero Dirichlet or Neumann conditions represent the appropriate conditions.  In the case of the Dirichlet conditions, the assumption is that the system has heat reservoirs that are contain substantial more heat energy than does the medium of interest (e.g., typically for heat problems in 1-dimension, one considers the transport of heat along an insulated rod between two reservoirs fixed at constant temperature.)  Now, consider the case where we have homogeneous Dirichlet or Neumann boundary conditions (i.e., regardless of which condition we choose, the value is zero at the boundaries $x=0$ and $x=\ell$).  Now, consider multiplying the heat/diffusion equation by the independent variable $u$. This yields

\begin{align}
&&    u\frac{\partial u}{\partial t} &= \alpha^2 u\frac{\partial^2 u}{\partial x^2}, \textrm{  for $x\in[0,\ell]$} 
\end{align}
With a little manipulation, this can easily be put in the form

\begin{align}
&&    \frac{1}{2}\frac{\partial u^2}{\partial t} &= \alpha^2 \frac{\partial}{\partial x}\left(u\frac{\partial u}{\partial x}\right)- \alpha^2\frac{\partial u}{\partial x}\frac{\partial u}{\partial x}, \textrm{  for $x\in[0,\ell]$} 
\end{align}
(You can confirm this by expanding the derivatives on the left and right sides of the expression).  Integrating this result over $x\in[0,1]$ gives

\begin{align}
&&   \int_{x=0}^{x=\ell} \frac{1}{2}\frac{\partial u^2}{\partial t} \, dx&= \alpha^2  \int_{x=0}^{x=\ell} \frac{\partial}{\partial x}\left(u\frac{\partial u}{\partial x}\right)\, dx- \alpha^2  \int_{x=0}^{x=\ell} \frac{\partial u}{\partial x}\frac{\partial u}{\partial x}\,, \textrm{  for $x\in[0,\ell]$} 
\end{align} 
Note that the first term on the right hand side is easily integrated to give 

\begin{equation}
    \left. \alpha^2  \left(u\frac{\partial u}{\partial x}\right) \right|_{x=0}^{x=\ell}
\end{equation}
which is identically zero because the boundary conditions are assumed to be homogeneous.  This gives the result

\begin{align}
&&  \frac{\partial }{\partial t}\int_{x=0}^{x=\ell} u^2 \, dx&= - 2\alpha^2  \int_{x=0}^{x=\ell} \frac{\partial u}{\partial x}\frac{\partial u}{\partial x}\,, \textrm{  for $x\in[0,\ell]$} \label{maxprinciple}
\end{align} 
This is an interesting expression; it indicates that the time-rate-of-change of the quantity $u^2$ is always negative (i.e., decreasing).  In more mathematical explorations of the heat/diffusion equation, this kind of analysis is the first step in proving a \emph{maximum principle} which indicates that the maximum value of the independent variable, $u$, always occurs either (1) internal to the domain at the initial condition, or (2) on the boundaries. Suppose that we had the case that the average temperature in the system was $\overline{u} = 0$.  If this were the case, then $u^2$ would represent the variance of the temperature (or concentration) in the system.  Thus, one way of thinking about Eq.~\ref{maxprinciple} is that it indicates that the \emph{variance} of the system is always decreasing; or, stating this differently, the system is becoming closer to uniform as time increases.  This is exactly in line with the concept of increasing entropy as the system progresses from state to state in time.  It is also consistent with our observations and intuition about heat or mass transport phenomena.  We know that heat tends to spread out in a body, or that chemicals tend to diffuse until their spatial gradients are relaxed.  Thus, generally speaking, the short proof that we outlined above states that heat and mass tend to spontaneously \emph{spread out} in space rather than spontaneously {concentrating} in space.  This is a useful property to have regarding the heat/diffusion equation.   Not only does it tell us something useful about its overall behavior (and that this behavior is consistent with our observations about the universe), but it will also serve as a useful constraint to help determine the correct mathematical solution to the problem.  We will see this in the next section.

\subsection{Separation of Variables for the Heat/Diffusion Equation}\indexme{heat/diffusion equation!separation of variables}

To begin the SOV method, we start by positing the following form for the solution.  To make things easier, we will take $\ell=1$; converting to other domain lengths is straightforward using the principles established in Chapter 3 (Fourier Series).

\begin{equation}
    u(x,t) = X(x)T(t)
    \label{SOV}
\end{equation}
Viewing this empirically, we will simply \emph{try} this solution to see if it leads to a useful result.  Toward that goal, we substitute Eq.~\eqref{SOV} into the governing differential equation to get

\begin{align}
    \frac{\partial (X(x)T(t))}{\partial t} = \alpha^2 \frac{\partial^2 (X(x)T(t))}{\partial x^2}
\end{align}
Because $X(x)$ and $T(t)$ are each functions of only one variable, we can simplify the notation by adopting the following $X'\Leftrightarrow d X/d x$, $X''\Leftrightarrow d^2 X/d x^2$, $T'\Leftrightarrow d T/dt$.  This way, the original PDE takes the form

\begin{align}
    X(x)T'(t) = \alpha^2  X''(x)T(t)
    \label{SOV2}
\end{align}
Dividing both sides by $\alpha^2 X(x)T(t)$ gives

\begin{align}
 \underbrace{\frac{T'(t)}{\alpha^2 T(t)}}_{\substack{\textrm{Function of}\\\textrm{only }t\\G(t)}}  =  \underbrace{ \frac{X''(x)}{X(x)}}_{\substack{\textrm{Function of}\\\textrm{only }x\\F(x)}}
\end{align}

Notice that an unusual thing has occurred here: for the equation in this form, apparently the left-hand side is a function of only $t$ (which we might call $G(t)$), and the right-hand side is a function of only $x$ (which we might call $F(x)$).  In short, this means that for all possible combinations of $x\in(0,1)$ and $t > 0$, we have the condition that

\begin{equation}
    G(t) = F(x), \textrm{ for any possible combinations of } x\in(0,1) \textrm{ and } t>0
\end{equation}
While this may not seem particularly astounding on first look, it is actually a remarkably strong condition that is being imposed.  Whatever the functions $X(x)$ and $T(t)$ are, they are such that the two sides of Eq.~\eqref{SOV2} are equal to each other.  Because the left-hand side depends only on $t$, and the right-hand side depends only on $x$, and because this equation must be satisfied for any pair that we select in the domain $(x,t)$, the \emph{only} possible choice is that 

\begin{equation}
    G(t)= F(x) = Constant
\end{equation}
In the material that follows, we will call this constant $\pm \lambda^2$.  Note that again we have chosen the square of our constant value to impose the condition that the constant itself be positive.  At this juncture, we have not assigned a sign to the constant, hence we list it with the $\pm$ operator.

A little though at this juncture indicates that we have the following condition

\begin{align}
    \frac{T'(t)}{\alpha^2 T(t)} = \frac{X''(x)}{X(x)} = \pm \lambda^2
    \label{SOV3}
\end{align}
And thus, what has happened in actuality is that the time and space components of this problem have been \emph{separated} up into two, ODEs expressed in variables that are orthogonal to one another.  To see this more clearly, we can rewrite the expression above as

\begin{align}
    \frac{T'(t)}{\alpha^2 T(t)} &=\pm \lambda^2& \frac{X''(x)}{X(x)} &= \pm \lambda^2
    \intertext{or, equivalently, as the two separate ODEs of the form}
    \frac{dT}{dt} &=\pm \lambda^2 \alpha^2 T(t)& \frac{d^2 X}{dx^2} &= \pm \lambda^2
    \label{SOV4}
\end{align}

These are two, homogeneous ODEs of first and second-order respectively.  We studied these cases back in Chapter 2, and their solutions are relatively straightforward.  A comment about notation is warranted here.  While the original equation involves \emph{partial} derivatives, note that here we have adopted the conventional notation of ordinary derivatives.  This is, in part, due to a long history of this convention, and it occasionally causes some confusion.  However, a some reflection on the definition of these derivatives can help clarify the situation.  Because $X(x)$, by construction, depends only on the independent variable $x$, there is no other variable for which one can consider a derivative for this function.  Thus, technically, this defines an ordinary derivative.  A similar argument can be made for the ODE in the time variable.  Thus, we have managed to take a single PDE, with partial derivatives in $x$ and $t$, and express it as the combination of two ODEs, where the time and space influences evolve independently from one another.  The good news is ODEs are something that we already have familiarity with.

To obtain a solution to our PDE, we now need to find the solutions to these two ODEs.  let's take the easiest of the two-- the first-order ODE in time.  This problem is solvable using a simple separation of variables.

\begin{align}
    \frac{dT}{dt} &=\pm \lambda^2 \alpha^2 T(t)\\
    \intertext{or,}
    \frac{dT}{T} &= \pm \lambda^2 \alpha^2 \,dt \\
    \intertext{Integrating both sides over $t\in[0, t] \Rightarrow T\in [T(0),T(t)]$ gives}
    \int_{T(0)}^{T(t)} \frac{dT}{T} &= \int_{0}^{t} \pm \lambda^2 \alpha^2 \, dt
\end{align}
The solution is the conventional exponential one of the form

\begin{align}
    \ln\left( \frac{T(t)}{T(0)}\right) &= \pm \lambda^2 \alpha^2 t \\
    \intertext{Or, solving for $T(t)$}
    T(t) &= T(0)\exp[\pm \lambda^2 \alpha^2 t]
\end{align}
Technically, both the positive and negative components of this equation represent \emph{solutions} to the original ODE.  However, we know from physical reasoning that the concentrations in a pure diffusion problem have their maximum somewhere in the initial state of the system ($t=0$) (including the boundaries, assuming that the boundaries are steady).  In fact, this is exactly what we proved to ourselves in Section \ref{variance}.  This is important,  because it indicates that the positive exponential root of the time component of the ODE is not a valid one because it leads to temperatures or concentrations that grow arbitrarily large as time increases, which is not consistent with the proof that the variance must decrease with increasing time.  Thus, we have explicitly that only the negative exponential is \emph{physically} relevant

\begin{align}
    T(t) &= T(0)\exp[- \lambda^2 \alpha^2 t]
\end{align}

This represents the solution to one of the two ODEs we have generated.  The second is given by 

\begin{equation}
    \frac{d^2 X}{dx^2}  + \lambda^2 = 0
\end{equation}
where we have adopted the sign imposed by the solution to the ODE in time. This second-order ODE has the characteristic function

\begin{equation}
    s^2 +\lambda^2 = 0
\end{equation}
Substituting this into the quadratic formula immediately gives the two possible roots for the equation

\begin{align}
    & s_1 = -\lambda ~i  && s_2 = + \lambda ~i
\end{align}
As a reminder, this represents Case 3 for the solution to a homogeneous second order ODE with constant coefficients; thus, the two roots are complex conjugates.  In chapter 2, we found that the solution to this case is given by

\begin{equation}
    X(x) = C \cos(\lambda x) + D \sin(\lambda x)  \label{solutionODE2}
\end{equation}
Where $C$ and $D$ must be determined by two ancillary conditions.  Before proceeding, note that the general solution (before any anciallry conditions are imposed) is given by

\begin{equation}
    u(x,t)=T(t)X(x) = T(0)\exp[- \lambda^2 \alpha^2 t]\left[ C \cos(\lambda x) + D \sin(\lambda x) \right]
\end{equation}
Each of $T(0)$, $C$ and $D$ are undetermined constants; because $T(0)$ multiplies both $C$ and $D$, it is apparent that these three constants are not independent.  In fact, there are really only two independent constants.  Thus, the solution can be written as

\begin{equation}
    u(x,t)=T(t)X(x) = \exp[- \lambda^2 \alpha^2 t]\left[ A \cos(\lambda x) + B \sin(\lambda x) \right]
\end{equation}
where $A= C \times T(0)$ and $B= D\times T(0)$.  In the future, we will use this form as the general solution for the heat/diffusion PDE.  Note that this kind of reduction in the number of constants will appear in applications of the SOV method to other kinds of PDEs.  For convenience, we will often reduce the equations in the future by \emph{specifying} one of the constants to be unity, which leaves the second constant to be set via the ancillary conditions.  Thus, for the discussion here, we can make the statement $T(0) = 1$, and this is equivalent to removing the extra (unnecessary) constant.

Continuing forward, we note that for this particular case the two ancillary conditions that we need are given by Eqs.~\eqref{bc1}-\eqref{bc2}, which we write here as

\begin{align}
    &B.C.~1& T(t)X(0) &= 0 \\
    &B.C.~2& T(t)X(1) &= 0 
\end{align}
Because these boundary conditions hold for all time, and because $T(t)$ is not, in general, identically zero, we must then have the following two equivalent conditions on the spatial boundaries

\begin{align}
    &B.C.~1& X(0) &= 0 \\
    &B.C.~2& X(1) &= 0 
\end{align}
Imposing the first of these on the solution given by Eqs.~\eqref{solutionODE2} gives

\begin{align}
     A \cos(0) + B \sin(0) =0
     \intertext{or}
     A = 0
\end{align}
Thus, we at least know now that, whatever the solution to the problem, it involves a decaying exponential in space, and a sine function.  To eliminate the second unknown, we impose the second of the two boundary conditions

\begin{align}
     B \sin(\lambda) =0 \label{part2}
\end{align}
While it is possible that the solution to this is $B=0$, we note that this would give $A=B=0 \Leftrightarrow u(x,t) \equiv 0$.  This solution is known as the \emph{trivial} or \emph{null} solution.  While it does represent a valid solution to the problem, it corresponds \emph{only} to the problem where the initial conditions is zero.  In other words, it is the solution to the problem where the initial and boundary conditions are all zero, thus there are no spatial derivatives.  In other words, the initial and boundary conditions are zero, and then nothing happens ever for all time.  One can see why this might be called the trivial solution.

Thus, in general we \emph{cannot have} both $A$ and $B$ be equal to zero.  We conclude that we must have $B\ne 0$, and thus we can divide both sides of Eq.~\eqref{part2} by $B$ to give

\begin{equation}
     \sin(\lambda) =0 
\end{equation}
Curiously, while this may initially seem like it is not much of an improvement, a little thought indicates that this has a non-trivial solution.  In fact, it has an infinite number of them!  We must have the solution

\begin{equation}
    \lambda = 0,~\pi,~2\pi,~3\pi \ldots
\end{equation}
Which indicates that $\lambda$ can be \emph{any} integer multiple of $\pi$.  Note: we have not considered \emph{negative} values of $\lambda$ because of the odd symmetry of the sine function $\sin(-\lambda)= -\sin(\lambda)$.  The minus sign that arises here can be adsorbed ultimately into the constant $B$, which remains undetermined at this juncture.  A similar argument can be made for solutions where the cosine function was selected, so we will not repeat this argument in the future.  Thus, the application of the two boundary conditions gives us not just a single solution, but actually an \emph{infinite number} of possible solutions, one for each value of $n$.  Because there exists an infinite number of solutions, then there must also be a corresponding infinite number of multiplicative constants, $B_n$.  This is easy to see by noting that if $\sin(n \pi x)$ is a solution, then $B_n\sin(n \pi x)$ is also a solution.  

At this juncture, we have the following:

\begin{enumerate}
    \item We have decomposed the solution using $u(x,t) = T(t) X(x)$.
    \item We have shown that this decomposition can be put back into the governing PDE, and the result is two ODEs that involve the constant $-\lambda^2$.
    \item The two solutions $T(t)$ and $X(x)$ can be found, are represented by $T(t) = \exp[-\lambda^2 \alpha^2 t]$ and $X(x) = B\sin(\lambda x)$, where we must have $\lambda = n \pi$, where $n = 0, ~1,~2,~3 \ldots$.  The function $X(x)$ has been determined using the two spatial ancillary (boundary) conditions. The initial condition has not yet been used.
    \item the solution to the problem for $X(x)$ indicates that we have an infinite number of possible solutions, one for every value of $n$; each such solution has (potentially) a different associated value for the remaining constant, $B_n$.  
\end{enumerate}

Using this information, apparently our solution at this point in the analysis is given by the infinite set of possible independent solutions, $u_n(x,t)$

\begin{align}
  u_n(x,t) &= T(t)X(x)\nonumber\\
  &=\exp[-n^2\pi^2 \alpha^2 t]B_n\sin(n \pi x)
\end{align}
Clearly, if the $u_n$ represent independent solutions, then so does its sum.  Thus, the most general solution we can write at this point in our analysis is as follows

\begin{align}
  u_n(x,t) &=\sum_{n=0}^{\infty}\exp[-n^2\pi^2 \alpha^2 t]B_n\sin(n \pi x)
\end{align}

We have an unusual seeming situation here.  We have a single remaining condition (the initial condition), but we have an infinite number of constants, $B_n$, to determine.  While initially this may seem intractable, we should note that, unlike the constant value boundary conditions, the initial condition contains in general a \emph{function}.  Thus, in some sense, the initial condition contains more information than the two boundary conditions.  To make more sense of this, we can apply the initial condition to see what happens.  Recall, the initial condition states

\begin{align}
&I.C.~1& u(0,x) &= \varphi(x) 
\end{align}
Thus, at the initial time, $t=0$, we must have the following 

\begin{align}
  u(x,0) &=\sum_{n=0}^{\infty}\exp[-n^2\pi^2 \alpha^2(0)]B_n\sin(n \pi x) \label{solnsum}\\
  \intertext{or,}
    \varphi(x) &=\sum_{n=0}^{\infty} B_n\sin(n \pi x)
\end{align}
We assumed at the start of the analysis that $\phi(x)$ is a function that meets the requirement that its Fourier series exists (see Chapter 3); thus, this result indicates that the $B_n$ are, for this particular case, \emph{exactly the amplitude constants for the Fourier sine series representation of} $\phi(x)$!  In some ways this result should not be surprising.  However, it is interesting to note that the SOV method automatically generates the Fourier sine series for the initial condition.  In fact, this analysis is close to how Fourier originally found the Fourier series in the first place (he was interested in finding solutions to the heat equation).  We have already studied Fourier series in Chapter 3, and using that information, we can simply write down the value for $B_n$ by

\begin{equation}
    B_n = 2\int_{x=0}^{x=1} \varphi(x) \sin(n \pi x) \, dx
\end{equation}
This, then, formally concludes the solution to the problem.   Note that for $n=0$, the resulting contribution to the sum in Eq.~\eqref{solnsum} is zero.  For this reason, it is customary to begin the numbering for the sine series at $n=1$.  The final result for the heat/diffusion problem with two homogeneous Dirichlet conditions on the interval $x\in[0,1]$ is

\begin{align}
\Aboxed{
  u_n(x,t) &=\sum_{n=1}^{\infty}\exp[-n^2\pi^2 \alpha^2 t]B_n\sin(n \pi x)}\\
  \intertext{where}
\Aboxed{B_n &= 2\int_{x=0}^{x=1} \varphi(x) \sin(n \pi x) \, dx }
\end{align}

Now that we have seen how the method works, it will be useful to see an example worked out in detail.  Note that the general method for SOV outlined above also works for the linear wave and Laplace equations in two variables, with very minor modifications.  We will explore those solutions after an example.

\begin{svgraybox}
\begin{example}[Separation of Variables: The Heat/Diffusion Equation with Dirichlet Boundary Conditions]\label{diffusion_example}

Consider the following equation on the interval $x\in [0,1]$.

\begin{align*}
&&    \frac{\partial u}{\partial t} &= \alpha^2 \frac{\partial^2 u}{\partial x^2},~ x \in [0,1]\\
&B.C.~1& u(0,t) &= 0 \\
&B.C.~2& u(1,t) &= 0 \\
&I.C.~1& u(0,x) &= \delta(x-\tfrac{1}{2}) \\
&Constraint & &\substack{\hspace{-9mm}\textrm{\normalsize Solutions stay bounded as}\\  \hbox{\hspace{-35mm} \fontsize{10}{25}\selectfont\( t\rightarrow \infty\)}}\\
&Constraint & &\substack{\hspace{-9mm}\textrm{\normalsize The initial condition has a valid} \vspace{1mm}\\  \hbox{\hspace{-14mm} Fourier series representation}}
\end{align*}
Rather than reproducing all of the details above, the SOV method will be conducted here in a summary fashion.  This is done to help build familiarity with the method; eventually, once the details are well understood, one can ``jump" to later stages of the solution with confidence, and then complete the solution from there.  To get us started, we recall that the substitution of $u(x,t) = T(t)X(x)$ yields two, separated equations of the form

\begin{align*}
    T(t) &= 1\cdot \exp(-\lambda^2 \alpha^2 t) \\
    X(x) &= A \cos(\lambda x) + B \sin(\lambda x) 
\end{align*}
Where, here, recall that we are arbitrarily setting the constant in the first ODE to unity because it is essentially a redundant constant.  Thus, the general solution, before applying any other conditions, is

\begin{equation*}
    u(x,t)= \exp(-\lambda^2 \alpha^2 t)[ A \cos(\lambda x) + B \sin(\lambda x)]
\end{equation*}
Application of the first boundary condition at $x=0$ rapidly shows us that we must have

\begin{equation*}
   0= \exp(-\lambda^2 \alpha^2 t)[ A \cos(0) + B \sin(0)]
\end{equation*}
Whatever the value of the exponential is, we know it is not zero for any finite time.  Thus, we can divide both sides by the exponential, and also evaluate $\cos(0)$ and $\sin(0)$.  This gives

\begin{equation*}
    0 = A + 0\cdot B \rightarrow A = 0
\end{equation*}
The cosine functions do not contribute to the solutions of this problem.  Thus, the solutions must be of the form

\begin{equation*}
    u(x,t)= \exp(-\lambda^2 \alpha^2 t) B \sin(\lambda x)
\end{equation*}
Applying the second boundary condition gives

\begin{equation*}
   0= B \sin(\lambda)
\end{equation*}
The constant $B$ cannot be zero unless the initial condition is also zero (which it is not).  Thus we must have 

\begin{align*}
   0&= \sin(\lambda)\\
   \intertext{or,}
   \lambda &= \pi,~2\pi,~3\pi\ldots =n\pi, ~\textrm{for $n\in \mathbb{N}$}
\end{align*}
The solution is now (recognizing that the most general solution must be a linear combination of all possible independent solutions)

\begin{equation*}
    u(x,t)=\sum_{n=1}^\infty \exp(-n^2 \pi^2 \alpha^2 t) B_n \sin(n \pi x)
\end{equation*}
Using the initial condition, we have

\begin{align*}
    \delta(x-\tfrac{1}{2}) = \sum_{n=1}^\infty  B_n \sin(n \pi x)
\end{align*}
This indicates that the coefficients $B_n$ are just the corresponding amplitudes for the Fourier series for the delta function.  Recall, while the delta function is not a \emph{proper} function, we can work with it as though its Fourier series makes sense (this was covered in Chapter 4).  The coefficients for the delta function are particularly easy to determine

\begin{align*}
 B_n &= 2\int_{x=0}^{x=1}  \delta(x-\tfrac{1}{2}) \sin(n \pi x) \, dx \\
 &= 2\sin\left(\frac{n \pi}{2}\right)
\end{align*}
Therefore, our final solution takes the form

\begin{equation*}
    u(x,t)=\sum_{n=1}^\infty \exp(-n^2 \pi^2 \alpha^2 t) \left[2\sin\left(\frac{n \pi}{2}\right)\right] \sin(n \pi x)
\end{equation*}
There is an important features of the solution that is revealed here.  Each value of $n$ corresponds to a particular frequency of the sine function.  Each such function is modulated by the decaying exponential term, which is a function of $n^2$.  What this indicates is that, the high-frequency components of the solution decay much more quickly than the low-frequency components of the solution.  From a physical perspective, this does make some sense.  A high-frequency wave in $x\in[0,1]$ is a function where the peaks are many, and each peak is close to its neighbors.  With just a small amount of diffusion, these peaks will fully mix, and be essentially eliminated from the solution.  For low-frequency components of the solution, there are many fewer peaks on the same interval, and each peak is, relatively speaking, far from its neighbors.  Thus, it takes much more time for the diffusion process to mix these peaks, because they have a much larger characteristic distance to diffuse!   We will examine this behavior more in a separate example.

For reference, the plot of the solution for various times is given in Fig.~\ref{deltasoln}.\\

{
\centering\fbox{\includegraphics[scale=.5]{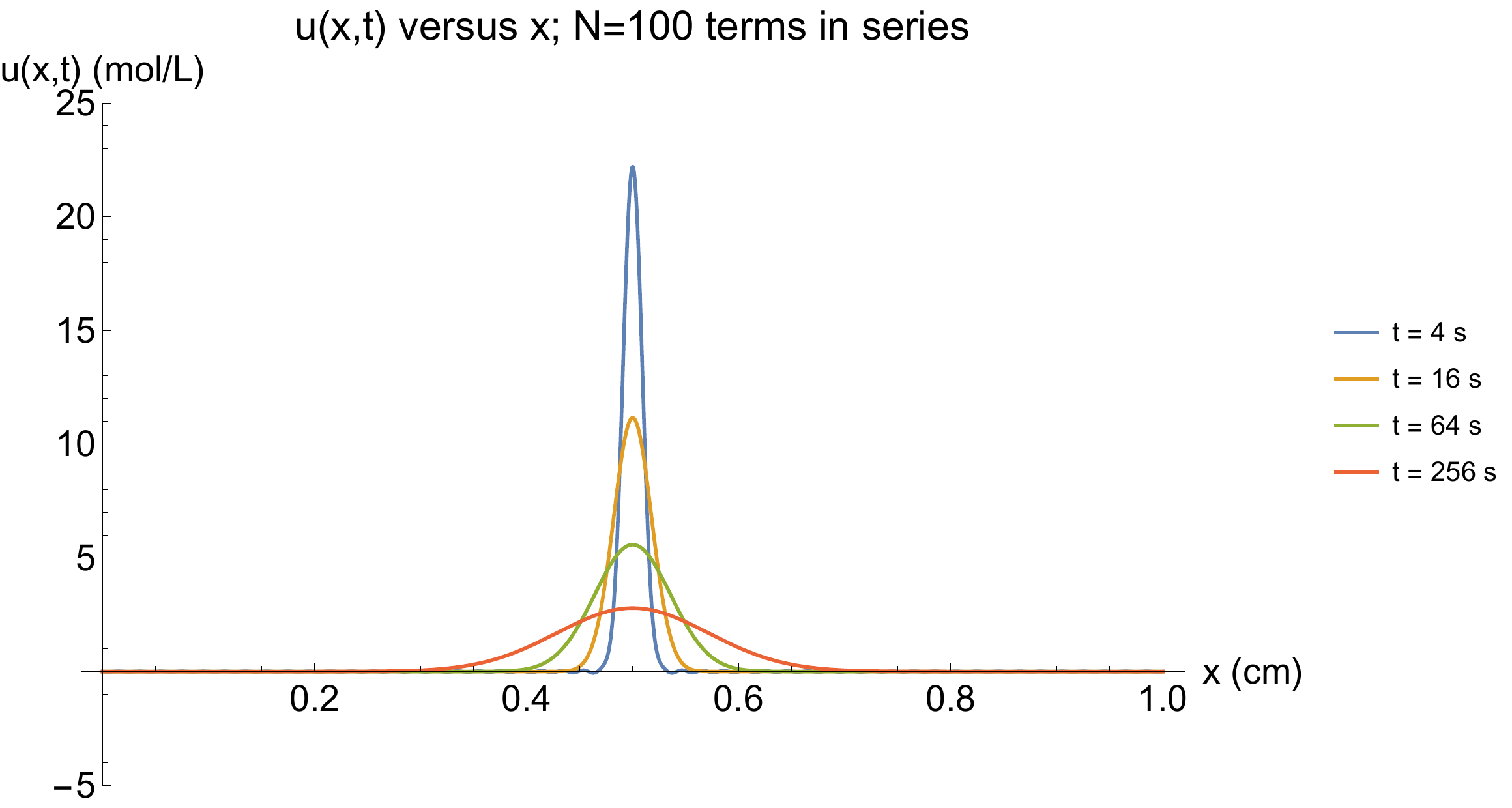} }
\vspace{2mm}
\captionof{figure}{The solution to the diffusion problem for several selected times.  For this plot, $x\in [0,1]$ cm, $0< t < 256$ s, $\alpha^2 = 1\times 10^{-5}$ cm$^2$/s. For this figure, the number of terms in the series was truncated at $N=100$.}\label{deltasoln}
 
}

\end{example}
\end{svgraybox}

\section{The Wave Equation}\indexme{second-order PDEs!wave equation}

The wave equation, as discussed in the previous chapter, represents any number of phenomena including water waves (e.g., gravity waves in the ocean), vibrations, and compression waves (e.g., the transmission of sound). Recall, the wave equation in one space dimension requires two space and two time ancillary conditions.  The example given here with Dirichlet conditions would represent, for example, the wave motion of a plucked string.  Again, for simplicity we assume that the interval is $x\in[0,1]$, and that time is measured starting from $t=0$.  For concreteness, we can think of this example as one where $u$ measures the \emph{displacement} of a string fixed at two ends as function of time and space.  Wave equations for other phenomena (e.g., for compression waves, $u$ represents the density of the compressible medium) have equally easily interpretable, but physically distinct, meanings for $u$.   Note that in the form we are considering, there is no \emph{damping} term, and thus the classical wave equation (as we will see) does not allow for the loss of energy from the system.  This represents an approximation that will be discussed in additional detail.

\begin{align}
&&    \frac{\partial^2 u}{\partial t^2}& = c^2 \frac{\partial^2 u}{\partial x^2}, \textrm{  for $x\in[0,1]$}\\
&B.C.~1& u(0,t) &= 0 \label{bc1w}\\
&B.C.~2& u(1,t) &= 0 \label{bc2w}\\
&I.C.~1& u(0,x) &= \varphi_0(x), \textrm{  for $x\in[0,1]$} \\
&I.C.~2& \frac{\partial u(0,x)}{\partial t} &= \varphi_1(x), \textrm{  for $x\in[0,1]$} \\
&Constraint~1 & &\substack{\hspace{-9mm}\textrm{\normalsize The initial condition has a valid} \vspace{1mm}\\  \hbox{\hspace{-14mm} Fourier series representation}}\\
&Constraint~2 & &\substack{\hspace{-9mm}\textrm{\normalsize The initial condition is at least} \vspace{1mm}\\  \hbox{\hspace{-21mm} $C^1$ continuous in time}}
\end{align}
This represents our starting point for our analysis of the wave equation using SOV.  Note that, for this case, $u$ represents the displacement from zero for the string. Thus the two initial conditions (required because of the second-order derivative in time) represent the initial displacement and the initial velocity for the string. 

\subsection{A Property of the Wave Equation: Conservation of Energy}\indexme{wave equation!conservation of energy}

Before continuing, it is helpful to consider the physical context of this equation.  If we think of $u$ as being the \emph{displacement} of a vibrating string fixed at both ends (such as a guitar string), then this equation describes how the displacements of such a string evolve in time.  If one imagines, for example, using a strobe light to capture instantaneous images at times $t_i, ~i=1,~2,~3,~4,~5$ of the vibrating string, the solutions $u(x,t_i)$ would represent the state of the displacement field (from zero displacement) at the the five indexed times.  

Our derivation of the wave equation did not account for any kinds of energy loss; thus, we would expect the resulting model to be an approximation, valid for early times, that did not include any energy dissipation mechanism (since we did not model one!)  Thus, the total energy in the string remains constant at all times.  If $u$ is the displacement field, and the one-dimensional density is $\rho = M/\ell$ (where $\ell$ is the length of the string, equal to unity for this example), then it is not hard to show that the following equation defines the energy per unit length in the string for all $(x,t)$ (this is not shown here, but is available as a problem to work out)

\begin{equation}
   e(x,t) = \rho \left(\frac{\partial u}{\partial t}\right)^2 + \rho c^2\left(\frac{\partial u}{\partial x}\right)^2
\end{equation}
Correspondingly, the total energy, $E_0$ of the string is given by

\begin{equation}
   E_0 = \rho \int_{x=0}^{x=\ell} \left(\frac{\partial u}{\partial t}\right)^2 + c^2\left(\frac{\partial u}{\partial x}\right)^2\, dx \label{tenergy}
\end{equation}
By definition, we must have that $E_0$ is a constant.  Thus, the energy in the string is conserved for all times.  This provides a third, important constraint for our problem.  Whatever the solution of the problem is, it must at least obey the following conservation of energy principle (imposed independently by the physical interpretation of the system) 

\begin{align}
&Constraint~3 & &\substack{\hspace{-9mm}\textrm{\normalsize The energy $E_0$ remains} \vspace{1mm}\\  \hbox{\hspace{-13mm} constant for all time}}
\end{align}
With this property established, we can continue the analysis.  We will find that this bit of \emph{physical} analysis constructed ahead of time will make the solution to the problem substantially less complex!  Note that all real physical systems have some kind of loss of energy (in accordance with the second law of thermodynamics), so the wave equation must represent an appropriate \emph{idealization} of systems that represents an approximation.

The conservation of energy for the undamped wave equation has some significant physical ramifications for the solutions.  In particular, it indicates that there are no \emph{steady state} solutions for the wave equation (except the trivial solution, where the initial condition is zero everywhere, and stays zero forever).  While there are no steady-state solutions, the wave equation does have \emph{periodic} ones.  This means that, if one fixes themselves at a point in the domain and observes a value for the dependent variable (say, the wave height), that value will reoccur again, and will actually reoccur an infinite number of times!  Periodicity is an interesting kind of behavior that appears in many physical systems ranging from pendulum motion, to the motion of the tides.  From this perspective, the wave equation aligns with our physical reasoning (up to a point!)   

In all real systems, there is a gradual loss of energy (damping) in media that support waves, so ultimately there can be a steady-state solution for systems where energy is lost.

\subsection{Separation of Variables for the Wave Equation}\indexme{wave equation!separation of variables}

As with the case for the heat/diffusion equation, we start by positing that a solution may exist of the form

\begin{equation}
    u(x,t) = T(t)X(x)
\end{equation}
and then substitute this form into the governing PDE to see if solutions of this form are, in fact, possible (and useful!).  Following the example of the heat/diffusion equation, we can very quickly get to the following expression

\begin{equation}
    \frac{T''(t)}{c^2 T(t)} = \frac{X''(x)}{X(x)} 
\end{equation}
As previously, one now makes the argument that independent functions of $t$ and $x$ can only be equal to one another for all possible pairs, $(x,t)$, if they are equal to the same constant.  By convention (and for distinction among the three basic problems that we will solve by SOV), we set the contant equal to $\omega^2$.  Thus, we have 

\begin{equation}
    \frac{T''(t)}{c^2 T(t)} = \frac{X''(x)}{X(x)} = \pm \omega^2
\end{equation}
Or, equivalently

\begin{align}
   T''(t) &=  \pm c^2\omega^2 T(t)\\
   X''(x)&= \pm \omega^2 X(x)
\end{align}

As with the case of the heat/diffusion equation, we can start our analysis with the time variable.  Because we have not yet assigned a sign to the separation constant, $\omega$, we have three possible solutions.  Thus, we need to examine the cases where the separation constant is equal to $-\omega^2$, $\omega^2 = 0$, and $+\omega^2$.  In other words, we need to examine the following ODEs.  To make these easier, we have computed the roots to the characteristic equation in the right-hand column

\begin{align}
&(a)&-&\omega^2 \Rightarrow& T''(t)  + c^2\omega^2 T(t)&= 0 & s_1=c\omega, ~~s_2 = -c\omega \\
&(b)& &~0~~ \Rightarrow& T''(t) &= 0 & s_1=s_2 = 0 \\
&(c)&+&\omega^2 \Rightarrow& T''(t)  - c^2\omega^2 T(t)&= 0 & s_1=i c\omega, ~~s_2 = -i c\omega 
\end{align}
The corresponding solutions are, then, as follows

\begin{align}
&(a)&-&\omega^2 \Rightarrow& T(t) = \alpha \exp(c \omega t)+\beta \exp(-c \omega t) \label{a}\\
&(b)& &~0~~ \Rightarrow& T(t) = \alpha + \beta t \label{b}\\
&(c)&+&\omega^2 \Rightarrow& T(t) = \alpha \cos(c \omega t) + \beta \sin(c \omega t) \label{c}
\end{align}

While all of these represent mathematically possible solutions to the wave equation, only one of them is consistent with our desire for the solution to have a constant, finite energy.  Solution (a) (Eq.~\eqref{a}) can be excluded because, regardless of what parameters we pick for $\alpha$ and $\beta$ (unless they are both zero, giving the trivial solution), the resulting solution will either grow exponentially or decay exponentially fast as $t\rightarrow \infty$.  Thus, energy cannot be maintained as a constant.  For solution (b) (Eq.~\ref{b}), the total energy given by Eq.~\eqref{tenergy} would grow (or decay) linearly in time for $\beta \ne 0$.  However, if $\beta=0$, then the solution would have no time dependence at all, which again corresponds only to the trivial solution for the given boundary conditions.  This leaves only Eq.~\eqref{c} as a possible solution.  While the solution is oscillatory in time, it is at least bounded regardless of how large time grows.  Thus, at this juncture, we adopt this as the only feasible solution, and turn our attention to the second ODE in space.  Note that this means that we are adopting the case where the separation constant is given explicitly by $+\omega^2$.  

For the space solution, we now have a fixed value for the separation constant.  Thus, we have only one possibility for the form of the solution

\begin{align}
    X''(x)&= + \omega^2 X(x) \\
    \intertext{or,}
    X''(x) - \omega^2 X(x)&=0
\end{align}
From the analysis above, we can immediately write down the solution as

\begin{equation}
    X(x) = A \cos(\omega x) + B \sin(\omega x)
\end{equation}
The solution as a whole to this stage is then given by

\begin{align}
    u(x,t) & = T(t)X(x) \\
    u(x,t) &= [\alpha \cos(c \omega t) + \beta \sin(c \omega t)][A \cos(\omega x) + B \sin(\omega x)]
\end{align}
This solution at least seems consistent with the solution explored previously for the heat equation, except now instead of having spatial waves that decay in time, we have waves in both time and space.  To make further sense out of this result, we can, as before, begin applying the ancillary conditions.  To start, we can apply the two Dirichlet boundary conditions.  The first of these gives

\begin{equation}
    0 =  [\alpha \cos(c \omega t) + \beta \sin(c \omega t)][A \cos(0) + B \sin(0)]
\end{equation}
We have already determined that both $\alpha$ and $\beta$ cannot be zero.  Also, while the remainder of the time component of the solution may be zero for some times (since it is the sum of two oscillating functions that are each sometimes zero), they are not zero for \emph{all} times, which is what is required by the boundary condition.  This means that we must have 

\begin{align}
    [A \cos(0) + B \sin(0)] &= 0 \\
    \intertext{or,}
    A= 0
\end{align}
The complete solution at this juncture is now

\begin{align}
    u(x,t) &= [\alpha \cos(c \omega t) + \beta \sin(c \omega t)][B \sin(\omega x)]
\end{align}
indicating that the spatial component of the solution will ultimately be represented by only sine functions.  Applying the second boundary condition (and applying the same logic regarding the time-component of the solution) gives the result

\begin{align}
   [B \sin(\omega)]=0
\end{align}
This is identical to the problem for the heat equation, and we find that the only possibility is that $\omega$ take on the values

\begin{equation}
    \omega = \pi, ~1\pi, ~3\pi \ldots = n \pi,~~ n=1,~2,~3, \ldots
\end{equation}
Here, remember that the4 case $n=0$ does lead to a solution for the problem, but it corresponds to the zero solution ($\sin(0\pi x)=0$), so by convention we do not include this in the final result.  

We have made some significant progress here.  If we substitute the possible values for $\omega$ into our solution now, we have 

\begin{align}
    u_n(x,t) &= [\alpha_n \cos(c n \pi t) + \beta_n \sin(c n \pi t)][B_n \sin(n \pi x)], ~~n=1,~2,~3,\ldots
\end{align}
Note that we have included a subscript $n$ to the solution variable, $u_n$, because we actually have an infinite number of them now, one for each value of $n$.  Note that we have again a situation where the time and space constants are multiplied, and it should be clear that we can set $B_n=1$ with no loss of generality (because $\alpha_n B_n$ forms a single new constant, and the same is true for $\beta_n B_n$; rather than renaming constants, we can just set $B_n=1$, and still have two undetermined constants in the end.)  As with the heat equation, we also conclude that the most general possible solution at this point is just the linear superposition of all of the possible solutions

\begin{align}
    u(x,t) &= \sum_{n=1}^\infty [\alpha_n \cos(c n \pi t) + \beta_n \sin(c n \pi t)] \sin(n \pi x)
\end{align}

From our previous experience with heat/diffusion problem, it is possible to guess at how things will evolve from here.  We have \emph{two} ancillary initial conditions for this problem, and imposing those will generate series that allow us to reproduce the two initial conditions at the time $t=0$.  Imposing the second initial condition first, we find

\begin{align}
    u(x,0) &= \sum_{n=1}^\infty \alpha_n \sin(n \pi x)
    \intertext{or,}
    \varphi_0(x) &= \sum_{n=1}^\infty \alpha_n \sin(n \pi x)
\end{align}
Without doing any additional work, we note that this is exactly the sine series expansion of the function $\varphi_0(x)$.  The constants $\alpha_n$ are thus known, and are given explicitly by

\begin{equation}
    \alpha_n = 2 \int_{x=0}^{x=1} \varphi_0(x) \sin(n \pi x) \, dx
\end{equation}
For applying the second boundary condition, we first need the time derivative of our proposed solution

\begin{align}
   \frac{ \partial u}{\partial t} &= \sum_{n=1}^\infty c \pi n [-\alpha_n \sin(c n \pi t) + \beta_n \cos(c n \pi t)] \sin(n \pi x)
\end{align}
And, for $t=0$ we must have

\begin{align}
   \varphi_1(x) &= \sum_{n=1}^\infty c \pi n  \beta_n  \sin(n \pi x)
\end{align}

If we let $\gamma_n = c \pi n  \beta_n$, we find 

\begin{align}
   \varphi_1(x) &= \sum_{n=1}^\infty \gamma_n  \sin(n \pi x)
\end{align}
so that this result is the Fourier series for $\varphi_1(x)$ with

\begin{align}
    \gamma_n &= 2 \int_{x=0}^{x=1} \varphi_1(x) \sin(n \pi x) \, dx  \\
    \intertext{or, converting back to $\beta_n$}
    \beta_n &= \frac{2}{c \pi n} \int_{x=0}^{x=1} \varphi_1(x) \sin(n \pi x) \, dx
\end{align}

Our final solution for this problem is then given by 

\begin{align}
\Aboxed{
    u(x,t) &= \sum_{n=1}^\infty [\alpha_n \cos(c n \pi t) + \beta_n \sin(c n \pi t)] \sin(n \pi x)}\\
\Aboxed{ \alpha_n &= 2 \int_{x=0}^{x=1} \varphi_0(x) \sin(n \pi x) \, dx}\\
\Aboxed{\beta_n &= \frac{2}{c \pi n} \int_{x=0}^{x=1} \varphi_1(x) \sin(n \pi x) \, dx}
\end{align}

\begin{svgraybox}
\begin{example}[Separation of Variables: The Wave Equation with Dirichlet Boundary Conditions]

Guitars are at least a familiar instrument to most people (although they do have fewer direct applications to engineering problems than would be hoped for in a perfect world), so they do make a good platform to understand the intuition behind the wave equation.  The length between the two bridges of a guitar (the scale length) varies between about 62.9 to 65.1 cm, but $L=64$ cm would be a good choice representing most acoustic and electric guitars currently made.  It turns out that the frequencies of guitar stings plucked on an open standard tuning (EADGBE) are reasonably well known, since the musical scale is independent of the instrument that it is played on!  One reference using a specific set of strings measured the following values for the frequencies, $f$, of a standard guitar tuning (all values in 1/s=Hertz=Hz) : E- 82.41, A- 110.00, D- 146.83, G- 196.00, B- 246.94, E(high)- 329.63.  For this example, we will examine the motion of the G string, which travels with a frequency of 196.00 Hz.  Recall the relationship between wavelength, $\lambda$ and length, $L$ is

\begin{equation}
    \lambda = 2 L
\end{equation}

This is because the fundamental solution to the problem is represented by half of a full sine series (i.e., the sine function goes from zero, to a maximum, and back to zero in a distance L; this is only \emph{half} of the full sine function, which has a negative component that has to be traversed before one full period of the sine function is expressed.)  To compute the value of $c$, then, we need only recall the expression

\begin{align}
    c &= f\lambda \\
    &= f 2L= 196~s^{-1}\times 2\times 0.64~m = 250.88~m/s
\end{align}

\begin{align*}
&&    \frac{\partial^2 u}{\partial t^2} &= c^2 \frac{\partial^2 u}{\partial x^2},~ x \in [0,L]\\
&B.C.~1& u(0,t) &= 0 \\
&B.C.~2& u(L,t) &= 0 \\
&I.C.~1& u(0,x) &= 4 \delta_0 \tfrac{x}{L}(1-\tfrac{x}{L}) \\
&I.C.~2& \frac{\partial u}{\partial t} = 0
\end{align*}
Where here $\delta_0$ is the initial maximum displacement (with units of length) during a plucking of the string.  

Rather than reproducing all of the details above, the SOV method will be conducted here in a summary fashion, as we did for the heat/diffusion problem.  

\begin{align*}
    T(t) &=  \alpha \cos(c \omega t) + \beta \sin(c \omega t) \\
    X(x) &= A \cos(\omega x) + B \sin(\omega x)
\end{align*}
Where, here, recall that we are arbitrarily setting the constant in the first ODE to unity because it is essentially a redundant constant.  Thus, the general solution, before applying any other conditions, is

\begin{equation*}
    u(x,t) = [\alpha \cos(c \omega t) + \beta \sin(c \omega t)][A \cos(\omega x) + B \sin(\omega x)]
\end{equation*}
Application of the first boundary condition at $x=0$ rapidly shows us that we must have

\begin{equation*}
   0=  A \cos(0) + B \sin(0) 
\end{equation*}
Which gives us the result $A=0$.
The cosine functions do not contribute to the solutions of this particular problem. Note however -- for other boundary conditions, this may change!  The solution to this point must be of the form

\begin{equation*}
    u(x,t)= [ \alpha \cos(c \omega t) + \beta \sin(c \omega t)] \sin(\omega x)
\end{equation*}
Applying the second boundary condition (again, noting that the time component of the solution is generally non-zero, thus we can focus on the space component as the portion that must be zero) gives

\begin{equation*}
   0= B \sin(\omega L)
\end{equation*}
The constant $B$ cannot be zero unless the initial condition is also zero (which it is not).  Thus we must have 

\begin{align*}
   0&= \sin(\omega L)\\
   \intertext{or,}
   \omega L &= \pi,~2\pi,~3\pi\ldots =n\pi, ~\textrm{for $n\in \mathbb{N}$}\\
   \omega &= \pi/L,~2\pi/L,~3\pi/L\ldots =\frac{n\pi}{L}, ~\textrm{for $n\in \mathbb{N}$}
\end{align*}
Recognizing that the most general solution must be a linear combination of all possible independent solutions, we form the sum to represent this.  Recall, we have the multiplication of $B$ with $\alpha$ and $\beta$; this generates two new constants.  Rather than rename them, we can use the original constant names, and think of $B=1$ as a condition we can impose without eliminating possible solutions.

\begin{equation*}
    u(x,t)=\sum_{n=1}^\infty [ \alpha_n \cos(c n \pi/L  t) + \beta_n \sin(c n \pi/L  t)]\sin(n \pi x/L)
\end{equation*}
Using the initial condition, we have

\begin{align*}
    4 \delta_0 \tfrac{x}{L}(1-\tfrac{x}{L}) = \sum_{n=1}^\infty  \alpha_n \sin(n \pi x/L)
\end{align*}
This indicates that the coefficients $\alpha_n$ are just the corresponding amplitudes for the Fourier series for the initial condition.  This is given by

\begin{align*}
\alpha_n &= \frac{2}{L} \int_{x=0}^{x=L} 4 \delta_0 \tfrac{x}{L}(1-\tfrac{x}{L})\sin(n \pi x/L) \, dx \\
 &= 16\delta_0 \frac{1-(-1)^n}{\pi ^3 n^3}
\end{align*}
The second boundary condition indicates that the initial velocity of the string is zero everywhere.  Thus, we need the Fourier series for the zero function, which is just zero (in other words, for this case, the initial velocity does not impact the solution because it is zero, and we have $\beta_n=0$ for all values of $n$).  Therefore, our final solution takes the form

\begin{equation*}
    u(x,t)=\sum_{n=1}^\infty  \cos(c n \pi  t) \left[ 16\delta_0 \frac{1-(-1)^n}{\pi ^3 n^3}\right] \sin(n \pi x/L)
\end{equation*}
There is some interesting interpretation that can be done here as we did with the heat/diffusion equation.

For reference, the plot of the solution for various times is given in Fig.~\ref{deltasoln}.  For this plot, we have taken $c=250.88$ m/s.\\

{
\centering\fbox{\includegraphics[scale=.45]{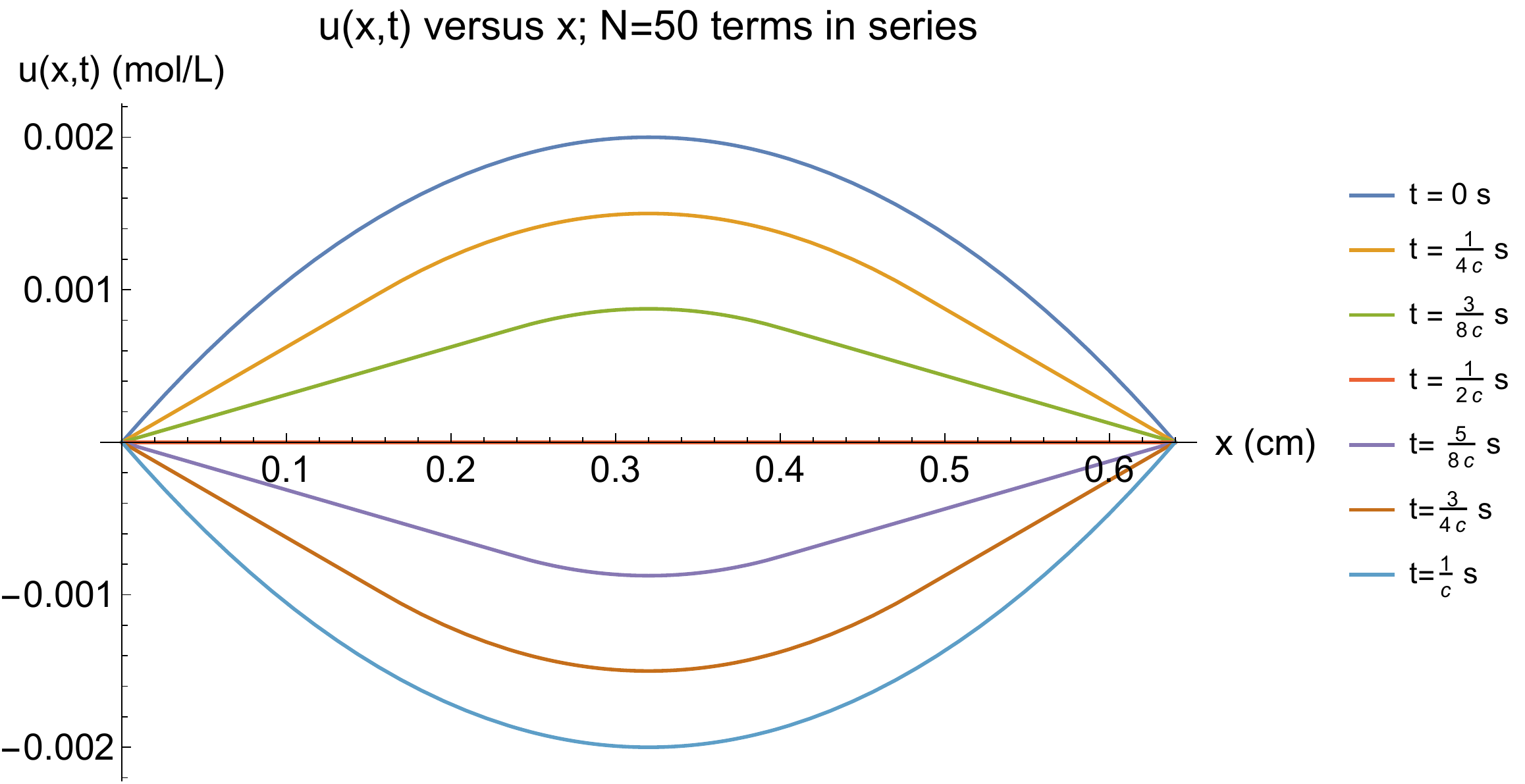} }
\vspace{-2mm}
\captionof{figure}{The solution to the wave problem of a plucked guitar string; spatial profiles are given for for several selected times.  For this plot, $x\in [0,64]$ cm, $0< t < 1/c = 0.004$ s. For this figure, the number of terms in the series was truncated at $N=50$.}\label{wavesoln}
}
~\\

{
\centering\fbox{\includegraphics[scale=.5]{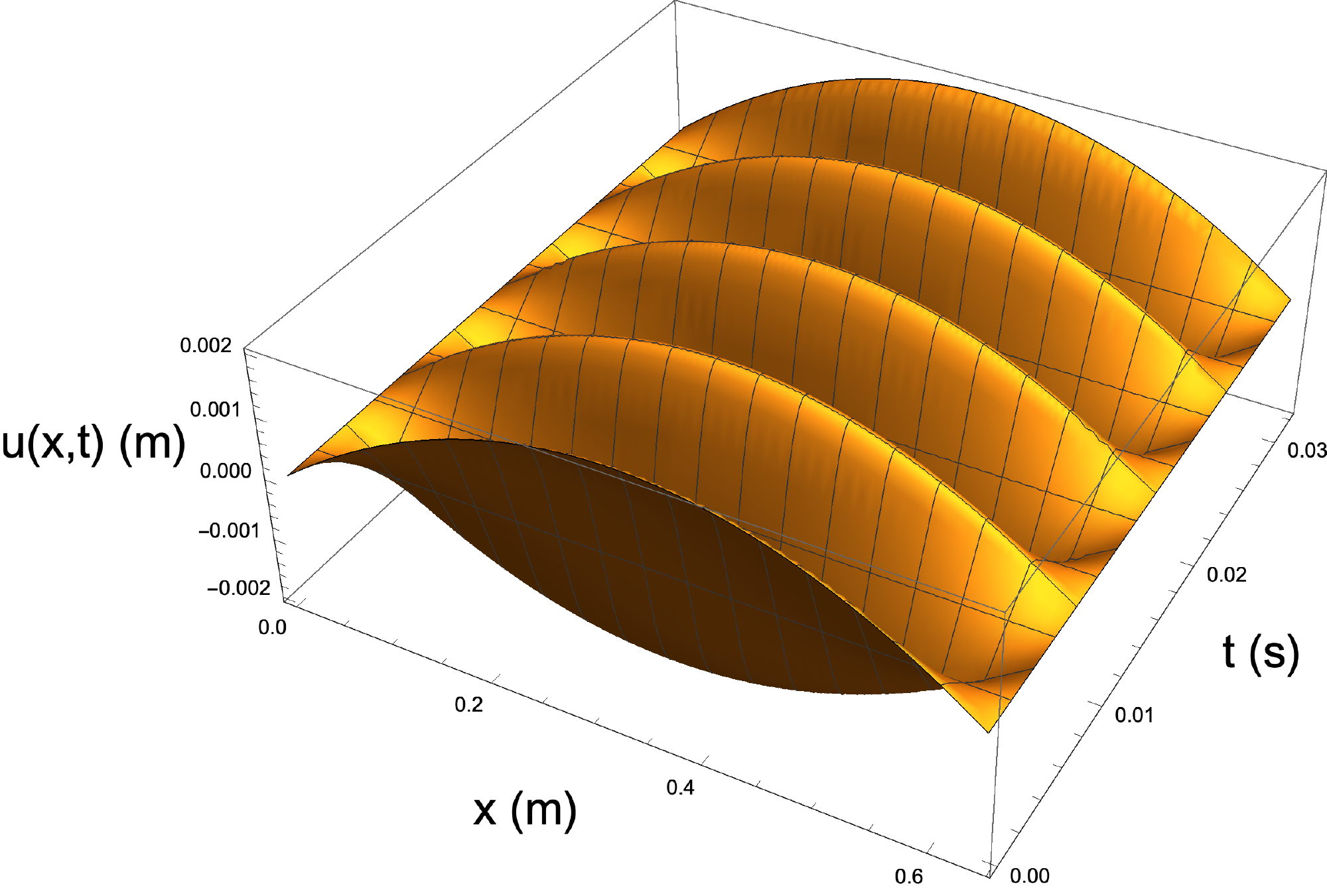} }
\vspace{2mm}
\captionof{figure}{The solution to the wave problem of a plucked guitar string. Here, the full time-space solution is shown as a surface.  Parameters are the same as for the previous figure.  Note that a total of four full periods in time are illustrated.  This solution assumes that energy is not lost to heat, which is an approximation that does not hold for long times (otherwise guitar strings would never stop vibrating after being plucked!)}
 
}

\end{example}
\end{svgraybox}

\section{The Laplace Equation}\indexme{second-order PDEs!Laplace equation}

The Laplace equation is an \emph{elliptic} equation in space; it has been such a fundamental part of applied mathematics that it has it is  a true archetype in PDEs.  In practice, the Laplace equations often arises as the \emph{steady state} form of the heat/diffusion equation; that is to say, the steady-state version of the heat equation in multiple dimensions yields a Laplace equation.  Technically, we have already investigated a one-dimensional Laplace-like equation in Section \ref{rxndiff0}, where we examined a steady-state diffusion-reaction problem.  A more general analysis of PDEs with multiple space dimensions is the subject of a separate chapter.  However, for illustrating the features of elliptic equations, the 2-dimensional Laplace equation is the simplest PDE that can be studied.  The two-dimensional Laplace equation takes the form

\begin{equation}
    \frac{\partial^2 u}{\partial x^2} +  \frac{\partial^2 u}{\partial y^2} = 0
\end{equation}
where $x$ and $y$ represent the conventional orthogonal axes in two dimensions.  When the right-hand side of this equation is a function of $x$ and $y$, then the equation is called a Poisson equation.  We will not study the Poisson equation in this chapter on separation of variables.

To more thoroughly understand this equation, first note the following heat/diffusion equation in two space dimensions.

\begin{align}
    &heat/diffusion& \frac{\partial u}{\partial t} & = K\left(\frac{\partial^2 u}{\partial x^2} + \frac{\partial^2 u}{\partial x^2} \right) &
\end{align}
Now note, the two \emph{steady state} versions of these problems are identical, and both are of the form of the Laplace equation.

\begin{align}
    &heat/diffusion~(steady)& 0 & = \frac{\partial^2 u}{\partial x^2} + \frac{\partial^2 u}{\partial y^2} & & &
\end{align}
Recall that the heat/diffusion equation is a parabolic equation that changes in time.  In the long-time limit of this equation, the result is mathematically an elliptic equation.  This is interesting in its own right.  Elliptic equations apply to a wide variety of physical processes, and many of them can be thought of as being the steady-state of some more complex \emph{transient} problem. 

\subsection{A Property of the Laplace Equation: The Maximum Principle}\indexme{Laplace equation!Maximum principle}

As mentioned above, the 2-dimensional Laplace equation can be thought of as the steady-state solution to the heat/diffusion equation in two spatial dimensions (and one time dimension!)  We have already illustrated that the heat equation with homogeneous boundary conditions is characterized by having a variance for $u$ that decreases in time. Thus, one can think of the Laplace equation as being one that has the smallest possible variance of all possible solutions to the heat equation.  It turns out that this smallest possible variance happens as time tends toward infinity, and these solutions are equivalent to the solutions of the Laplace equation with the same boundary conditions.  

A little thought will indicate that for homogeneous boundary conditions, the only steady-state solution to the heat equation is given by the solution

\begin{equation}
    u(x,t\rightarrow \infty) = C_1
\end{equation}
where $C_1$ is a constant that is possibly zero.  So, in general The steady state solutions to the heat/diffusion equation are not tremendously interesting at long times when the boundary conditions are homogeneous.  We will handle nonhomogeneous for the heat/diffusion equation (and wave equation) later.  For now, we can consider the Laplace equation in two space dimensions by thinking about the $y-$ spatial direction as being ``time like" (for the purposes of pondering the qualities of the problem and its solution), and the $x-$direction as being space like.  Under those circumstances, we might suspect that our overall scheme of separation of variables will still work, just by analogy with other equations in two independent variables that we have examined.  In the material below, we will see that SOV does in fact work for the Laplace equation considered in this way.  However, before continuing on, it is useful to establish a few other properties for solutions to the Laplace equation.

For largely historical reasons, solutions to the Laplace equation are called \emph{harmonic functions}.  This bit of vocabulary is useful primarily because one finds the term harmonic function used everywhere in the mathematics and engineering literature.  The fact that harmonic functions have been given a special name is also an indication that they must have interesting properties (we rarely give uninteresting things special names!)  There are a number of interesting properties regarding harmonic functions that are worth being noted.  The following are two theorems regarding harmonic functions in two dimensions that are easy enough to prove \citep[][Chp.~4]{olver2014}, but we will not do so here. 

\begin{theorem}[Average of a Harmonic Function]
Suppose $u(x,y)$ is the solution to the Laplace equation (with, in general, non-homogeneous boundary conditions)in the two-dimensional domain, $D(x,y)$.  Then, for every point $(x_0,y_0)$ within the domain $D$, the following is true

\begin{equation}
    u(x_0,y_0)= \frac{1}{2 \pi}\int_{\theta = 0}^{\theta=2 \pi} u(x_0+a \cos \theta, y_0+a \sin \theta)\, d\theta
\end{equation}
where $a$ is the radius of any circle that does not cross the boundary $D$ (although it can be tangent to the boundary).  In other words, the value at the center of every circle in the domain is the \emph{average} value of the function $u$ within the circle.  
\end{theorem}

This is a very interesting statement that allows many other things to be determined about harmonic functions.  Together with the next theorem (again, offered without proof), these properties tell us almost everything that we might want to know about harmonic functions.

\begin{theorem}[Harmonic Functions are Analytic]
Recall that a function is analytic in a domain, $D(x,y)$, when it has a convergent Taylor series for each point $(x,y)$ contained in the domain.  All harmonic functions are analytic at each point in the interior of the domain $D$ (this emphasizes that they may not be analytic on the boundary of the domain).
\end{theorem}

Using the first of these two theorems, we can prove one additional thing about harmonic functions.  While the proof is not difficult, we report only the result here.

\begin{theorem}[Maximum-Minimum Theorem for Harmonic Functions]
The solution to the Laplace equation has both its maximum and its minimum values of $u$ on the boundary of the domain $D$.  This means that if the boundary conditions are specified by functions of $u$ that are finite, then the solution inside the domain is also finite.
\end{theorem}

\subsection{SOV for the Laplace Equation}\indexme{Laplace equation!separation of variables}

 To be specific, consider the following steady-state heat transport problem on a unit square (where $x\in[0,1]$ and $y\in[0,1]$).  For this problem, think of the dependent variable, $u$ as being the temperature in a 1 meter by 1 meter square with a unit value for the thermal conductivity.  The problem statement is

\begin{align}
 &&   \frac{\partial^2 u}{\partial x^2} +  \frac{\partial^2 u}{\partial y^2} &= 0 \\
 &B.C.~1& u(x,0)&=\sin(\pi x) \\
 &B.C.~2& u(x,1) &=0 \\
 &B.C.~3& u(0,y) &= 0 \\
 &B.C.~4&  u(1,y) &=0 \\
 &Constraint& &\textrm{\hspace{-10mm}The solution and its derivatives remain bounded.} \nonumber
\end{align}
The last of these constraints is motivated again by physical conditions.  Recalling that we are thinking about these kinds of problems as being the steady-state version of some parabolic (heat/diffusion) type equation, then the physical constraints that we have already assumed for those must apply (thus, the solutions must be bounded).  The second part of the constraint can be considered to be a condition that requires that all \emph{fluxes} remain finite.  Recall that fluxes are often given by linear flux relations (i.e., a constant times the first derivative of the independent variable); thus, demanding that the derivatives remain bounded means that no fluxes are allowed to approach infinite (positive or negative) values. This is certainly consistent with the physics of the problem!

Now, if at least conceptually we think of the $y-$direction as being ``time like", then in some sense, this problems seems like one in which the first tow boundary conditions give information that propagates into the domain, much like an initial condition would.  While the physics of the problem represents, in fact, a steady state condition, thinking about the problem in this way is somewhat helpful.  It at least suggests that our conventional SOV approach has some hope of working, since it has been successful with problems that have an analogous structure.  

To check to see if our SOV approach will work for this problem, we can start with the conventional supposition: We want to determine if a solution of the form 

\begin{equation}
    u(x,y) = X(x)Y(y)
\end{equation}
can be found.  Substituting directly into the Laplace equation, we find

\begin{equation}
    X''Y + X Y'' = 0
\end{equation}
and, after dividing both sides by $XY$ we have

\begin{equation}
    \frac{X''}{X} = - \frac{ Y''}{Y} 
\end{equation}
As with the previous solutions for the heat/diffusion and wave equations, the only possibility for this result to be valid is if the two sides (both of which are independent from one another) are equal to the same constant.  For this case, we still have to determine the sign of the constant, even though it is clear that we will have one positive and one negative value of the separation constant in the final solution.  Recall, for the heat/diffusion problem, we already established that the appropriate choice for generating a solution with homogeneous Dirichlet conditions in space was that the constant should be negative.  Here, we have same kind of problem.  In the $x$-direction, we have two homogeneous boundary conditions.  Thus, the proper sign for the separation constant, by analogy, is the one that leads to the separated ODE $X''+\lambda X = 0$.  This requires a negative separation constant for this case.  Note, if we were to switch the roles of $x$ and $y$ (i.e., if we were to apply the nonhomogenous boundary condition $y(0,y) = \sin(\pi x)$), we would have reached the opposite conclusion!

Choosing constant to be ``-$\lambda^2$", then we have the two ODEs

\begin{align}
    X'' +\lambda^2 X &= 0 \\
    Y'' - \lambda^2 Y & = 0
\end{align}
As before, we can take these equations one at a time to determine whether or not our proposed solution is actually a reasonable one. Because the $x-$variable has two homogeneous Dirichlet conditions [$u(x=0,y) = 0$ and $u(x=1,y) = 0$], we propose starting with that one.  The solution is the familiar one that we also saw in the analysis of the heat/diffusion equation.

\begin{equation}
    X(x) = A \cos(\lambda x) + B \sin(\lambda x) 
\end{equation}
The two homogeneous boundary conditions demand that, whatever the solution for $X(x)$, it must be zero at the two boundaries in the $x-$direction.  This gives again a familiar result

\begin{align}
    X(0)& = 0 \nonumber\\
    0 & = A \cos(0) + B \sin(0) \nonumber\\
    \intertext{and, as we have seen previously, this implies}
    A &= 0
\end{align}
Thus, in the $x-$direction, the functions we use to ``build" our solution will be the sine functions.  Applying the second boundary condition gives

\begin{align}
    X(1) &= 0 \nonumber \\
    0 &= B \sin(\lambda) \nonumber \\
    \lambda & = n \pi, ~~ n=1,~2,~3\ldots~~\textrm{ (or, equivalently } n\in \mathbb{N}) \nonumber \\
    \intertext{This implies then}
    X_n(x) &= B_n \sin(n \pi x), ~~n\in \mathbb{N}
\end{align}
As with our previous solutions, we find that the solution demands that there be an infinite number of possible sine functions involved, each with a different frequency as modulated by the quantity $n\pi$.  

Turning our attention to the second ODE, we have 

\begin{equation}
    Y'' - \lambda^2 Y  = 0
\end{equation}

The characteristic equation for this problem is (noting $a=1$, $b=0$, and $c=-\lambda^2$)

\begin{equation}
    s^2 -\lambda^2 = 0
\end{equation}
and this has roots
    
\begin{align}
    s&= \frac{0\pm\sqrt{4\lambda^2}}{2}
    \intertext{or,}
    s_1 &= \lambda, ~~ s_2 = -\lambda
\end{align}
The solutions that correspond to these roots are the exponential ones.  

\begin{align}
    Y(y) &= C_1 \exp(\lambda y) +C_2 \exp(-\lambda y)
\end{align}
The two remaining boundary conditions will help us work out the final solution for the problem.  At this juncture, we have 

\begin{align}
    u_n(x,y) &= X_n(x)Y_n(y) \nonumber \\
    u_n(x,y) &= B_n \sin(n \pi x) [C_1 \exp(n \pi y) +C_2 \exp(-n \pi y)]
\end{align}
A good rule of thumb is to always apply homogeneous boundary conditions \emph{first}; it is also helpful to apply boundary conditions where one of the independent variables is zero.  To that end, we note that ``B.C. 2" is homogeneous, and evaluated where $y=0$.  
\begin{align}
        u_n(x,y) &= 0\nonumber\\
        B_n \sin(n \pi x) [C_1 \exp(\lambda y) +C_2 \exp(-\lambda y)] = 0
\end{align}

Solving this, we find

\begin{equation}
    C_1 = -C_2 \exp(-2 \lambda)
\end{equation}
Substituting this into the solution gives us

\begin{equation}
   u(x,y) = B_n \sin(n \pi x)C_2[-\exp(-2\lambda)\exp(\lambda y) +\exp(-\lambda y)]
\end{equation}
As we have seen in previous examples, we can set $B_n C_2 \Rightarrow B_n$ (i.e., we can take $C_2 = 1$) because the product of two undetermined constants mathematically represents a single undetermined constant.  Thus, we have 

\begin{equation}
   u_n(x,y) = B_n \sin(n \pi x)[-\exp(-2\lambda)\exp( \lambda y) +\exp(-\lambda y)]
\end{equation}

Our final step is to apply the remaining boundary condition at $y=0$.  Before proceeding, we note that the most general condition for the solution is a linear combination of all possible solutions to this point.  While we could have imposed that condition at several points in the previous analysis, it is most convenient to impose it it now.  This gives the solution in the form (where we have also substituted $\lambda = n\pi$)

\begin{equation}
   u(x,y) = \sum_{n=1}^{\infty} B_n \sin(n \pi x)[-\exp(-2 \pi)\exp(-n \pi y) +\exp(-n \pi y)]
\end{equation}
Now, imposing the initial condition, we find

\begin{align}
    \sin(\pi x) = \sum_{n=1}^{\infty} B_n \sin(n \pi x)[1-\exp(-2 n \pi)]
\end{align}
The term on the right-hand side is just a Fourier series for $\sin(\pi x)$, but one that has a rescaled value for $B_n$ (note that the rescaling function, $[1-\exp(-2n\pi)$, does not depend on $x$).  We can multiply both sides by $\sin(m \pi x)$, integrate, and use the conventional properties of orthogonality to complete the problem.

\begin{align}
    \int_{x=0}^{x=1} \sin(\pi x) \sin(m \pi x) \, dx = \sum_{n=1}^{\infty} B_n \int_{x=0}^{x=1}\sin(n \pi x)\sin(m \pi x)\,dx [1-\exp(-2 n \pi)]
\end{align}
On the left-hand side of this expression, we must have that $m=1$; all of the other possible integer values for $m$ lead to an integral that is zero by orthogonality.  On the right hand side, because $m=1$, we must also have only a single term in the sum corresponding to $n=1$.  Putting this all together, we find

\begin{align}
   \frac{1}{2} =  B_1 \int_{x=0}^{x=1}\sin( \pi x)\sin( \pi x)\,dx [1-\exp(-2  \pi)]
\end{align}

Or, we find that 

\begin{equation}
    B_1 = \frac{1}{1-\exp(-2 \pi)}, ~~B_n=0~~\textrm{for } n>1
\end{equation}
The final solution to our problem is

\begin{equation}
   u(x,y) = \frac{1}{1-\exp(-2 \pi)} \sin(\pi x)[-\exp(-2 \pi)\exp(- \pi y) +\exp(-\pi y)]
\end{equation}
Upon additional algebra, the result is the somewhat simple looking result

\begin{equation}
   u(x,y) = \sin(\pi x)\exp(-\pi y)
\end{equation}
Solutions to the Laplace equation are interesting from a number of perspectives.  In many applications, they are the surfaces that \emph{minimize} the energy associated with the surface.  This would include both the steady-state versions of the heat/diffusion equation, and the steady-state version of the wave equation.   In both cases, we have discussed the notion of the energy associated with the solution $u$.  Even though the underlying equations have different physical interpretations, it is interesting to note that they are both governed by the same steady-state solution.  In both cases, this particular steady state solution is the \emph{minimal surface}, which is just a general term used to describe a surface that minimizes some ancillary function (in this case, the energy).\\

\begin{figure}[t]
\sidecaption[t]
\centering
\includegraphics[scale=.35]{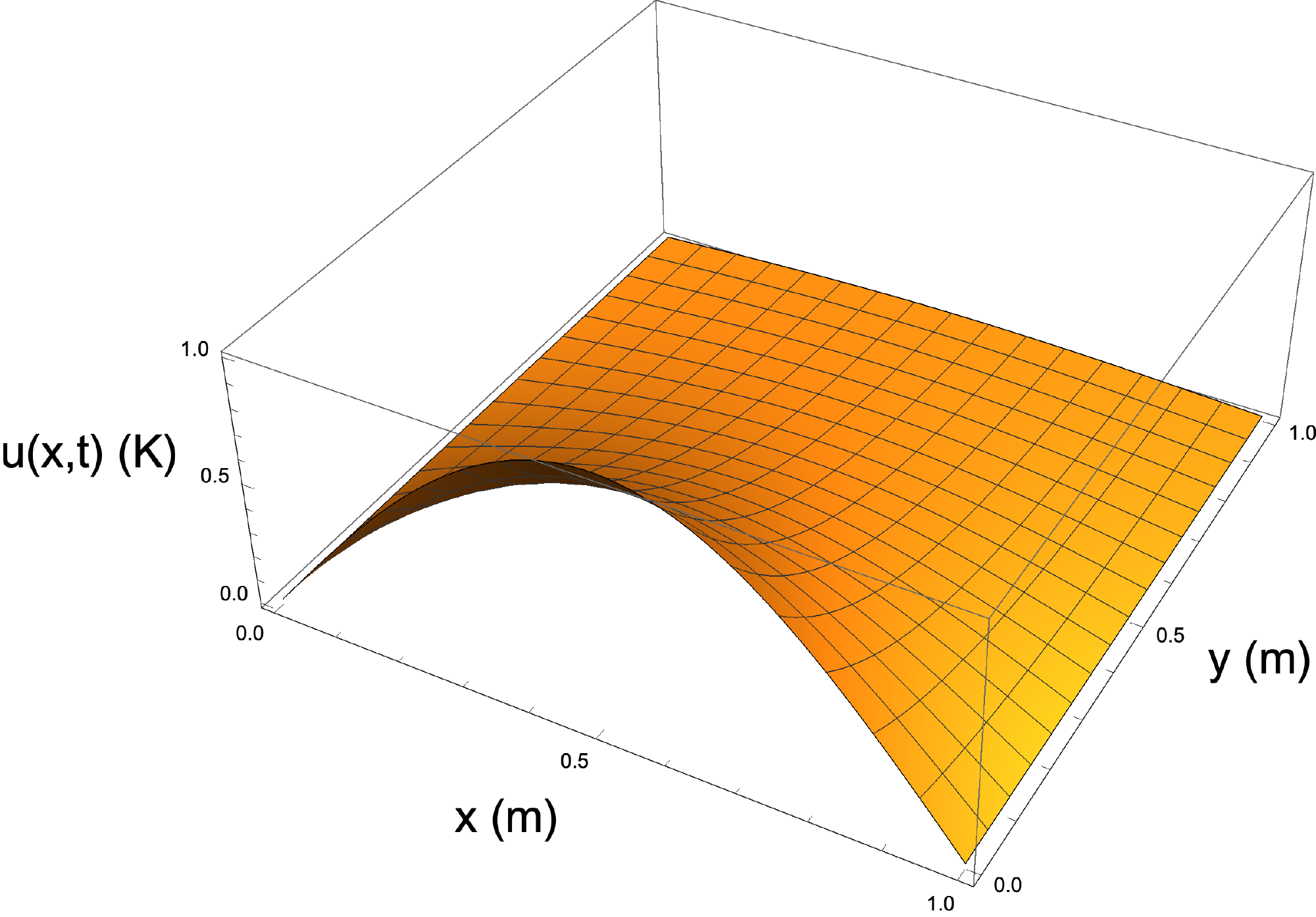}
\includegraphics[scale=.35]{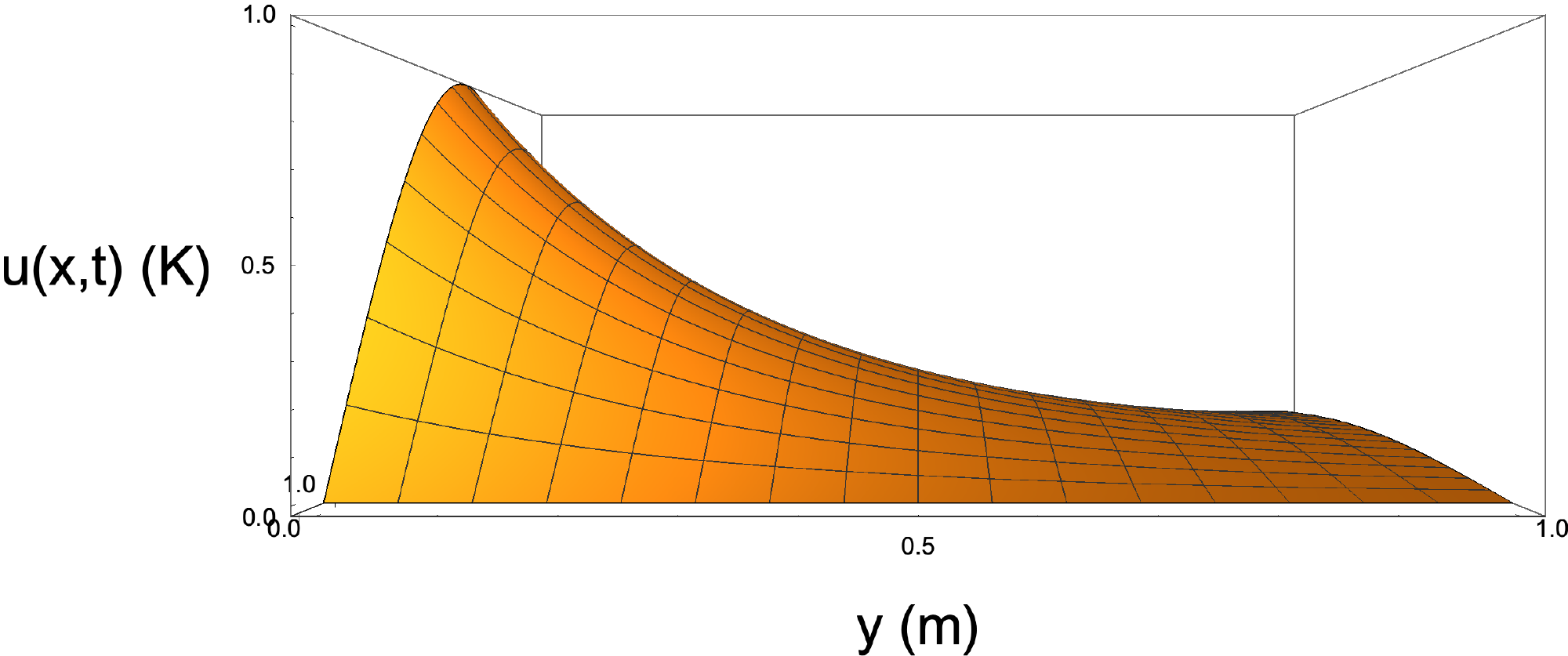}
\caption{{\bf Solutions to Laplace Equation.} The Laplace equation solved via the SOV method.  This approach works when the boundaries in one direction are homogeneous, and in the other direction there is one homogeneous boundary, and one (non-zero) specified value boundary.  }
\label{Laplacefigs}       
\end{figure}

%
\section{Linear Superposition}\indexme{second-order PDEs!linear superposition}\indexme{superposition, linear}
%
Linearity is a very powerful property when it exists.  While it is true that \emph{many if not most} differential equations used in physics and engineering are most generally nonlinear, almost all of them have linear regimes (or approximations) that can be very useful.

One of the most useful concepts in the theory of linear differential equations is the concept of \emph{superposition}.  This concept is not arcane or complex; it is just a natural consequence of linearity.  There are a number of situations in which one may want to consider decomposing a complex linear problem into two or more simpler linear problems.  The critical feature here is that the \emph{sum} of the simpler problems must be, in some sense, equivalent to the original complex problem.

\begin{definition}[Linear Superposition.]\indexme{Superposition!linear}\indexme{Superposition} For any two solutions to a linear and homogeneous differential equation, the linear combination of those solutions is also a solution to the differential equation.  This is sometimes called the \emph{principle of superposition}, and it applied equally well to ODEs and PDEs. 
\end{definition}

The best way to understand the principle of linear superposition is, as is normal, through an example.

%
%
\begin{svgraybox}
\begin{example}[Separation of Variables: The use of the linear superposition principle]\label{linsup}

\begin{align}
 &&   \frac{\partial^2 u}{\partial x^2} +  \frac{\partial^2 u}{\partial y^2} &= 0 \\
 &B.C.~1& u(x,0)&=\sin(\pi x) \\
 &B.C.~2& u(x,1) &=\sin(\pi x) \\
 &B.C.~3& u(0,y) &= 0 \\
 &B.C.~4&  u(1,y) &=0 \\
 &Constraint& &\textrm{\hspace{-10mm}The solution and its derivatives remain bounded.}\nonumber
\end{align}

\end{example}
\end{svgraybox}

%
%
\begin{svgraybox}
\begin{example}[Separation of Variables: The Laplace equation as a model for soap film geometry]

{
\centering\fbox{\includegraphics[scale=.45]{\CHAP/wave_times.pdf} }
\vspace{-2mm}
\captionof{figure}{The solution to the wave problem of a plucked guitar string; spatial profiles are given for for several selected times.  For this plot, $x\in [0,64]$ cm, $0< t < 1/c = 0.004$ s. For this figure, the number of terms in the series was truncated at $N=50$.}
}
~\\

{
\centering\fbox{\includegraphics[scale=.5]{\CHAP/wave_spacetime.pdf} }
\vspace{2mm}
\captionof{figure}{The solution to the wave problem of a plucked guitar string. Here, the full time-space solution is shown as a surface.  Parameters are the same as for the previous figure.  Note that a total of four full periods in time are illustrated.  This solution assumes that energy is not lost to heat, which is an approximation that does not hold for long times (otherwise guitar strings would never stop vibrating after being plucked!)}
 
}

\end{example}
\end{svgraybox}

\abstract*{This is the abstract for chapter 00}

\begin{savequote}[0.55\linewidth]
``Roughly speaking what Fourier developed was a mathematical way of converting any pattern, no matter how complex, into a language of simple waves. He also showed how these wave forms could be converted back into the original pattern. In other words, just as a television camera converts an image into electromagnetic frequencies and a television set converts those frequencies back into the original image, Fourier showed how a similar process could be achieved mathematically. The equations he developed to convert images into wave forms and back again are known as Fourier transforms."

\qauthor{ Michael Talbot, in ``The Holographic Universe", 1991 }
\end{savequote}


\chapter{Fourier Transforms}\label{Fouriertransforms}
%
\def\CHAP {chapter10_Fourier_Transforms}
%
\section{Introduction}

Throughout the material presented so far in this text, we have seen that Fourier series (and related series based on functions defined by Sturm-Liouville theory) have broad utility.  Fourier series have applications in a range of disciplines, from computer science (where they can be a component of various compression algorithms) to time-series analysis of natural phenomena.  It is not an overstatement to say that the development of the Fourier series was one of the most important mathematical results in modern history.

One class of problems that is \emph{not} addressable by Fourier series analysis is the class of functions on the entire real line.  For a (sufficiently well-behaved) function on any \emph{finite} interval, a combination of shifting the function and periodic extension lead to the ability to represent the function as a Fourier series.  However, for functions defined on the real line (i.e., $\lim_{L\rightarrow \infty} x\in (-L,L)$), there is no corresponding definition for the Fourier series.

It is actually possible to extend the basic notion of the Fourier series to the entire real line.  In short, this conversion will lead to the definition of the continuous Fourier Transform.  Although the Fourier Transform is not particularly complicated, there are many details that often frustrate easy understanding of the transform.  In addition, there is a second method known as the finite Fourier transform that applies to periodic systems. In fact, it will turn out that this latter transform is just an application of the Fourier series method described in chapter \ref{FS_1}.  However, viewing the approach as a \emph{transform} rather than as an infinite trigonometric series for representing a function involves a change in perspective, and a useful one at that.

In this chapter, the intent is to break the understanding of the Fourier transform down into more readily digestible pieces.  While there are many uses for the Fourier transform, we will use the approach to continue to study solution method for PDEs.

\section{Terminology}

\begin{itemize}

\item {\bf Integral Transform.}  A method where a function is integrated against some \emph{kernel} over its entire domain; the result is that its independent variable is replaced by a new independent variable contained in the kernel function.  The purposes is to map the use the properties of integration to (1) bring out specific physical properties of the function (e.g., as Fourier transforms do by converting to wave space), or (2) to simplify some of the operations of the original equations (e.g., derivatives often become algebraic when integral transformed). A key component of the transforms is that they are both linear and invertible.  In other words, a transformed function can be recovered by computing the appropriate inverse transform.   \indexme{integral transforms}\\

\item \textbf{Fourier Transform.}  The Fourier transform is a generalization of the Fourier series to the infinite real line.  It is an integral transform.  If $x$ is the original independent variable, and $\xi$ is the kernel variable, then the kernel function is  $\exp(-i \xi x)$.  Upon transforming a function of $x$ a new function is recovered that is a function of $\xi$ and no longer a function of $x$.  This is sometimes represented by $f(x) \xrightarrow{\mathscr{F}[f]} \hat{f}{\xi}$.

\item \textbf{Finite Fourier Transform}.  A representation of the conventional Fourier series as a transform method.  Even though a finite Fourier transform occurs on a bounded interval $x\in[a,b]$, it can still be cast as an integral transform.  In this case, the inverse transform is the classical infinite sum that defines the Fourier series.

\item \textbf{Euler's Identity.}  We have seen this identity previously, but as a reminder it is the identity linking sine, cosine, and the complex exponential.  Two frequently encountered forms of the identity are
\begin{align*}
    \exp({i x}) &= \cos(x) + i \sin(x) \\
     \exp({i \pi x}) &= \cos(\pi x) + i \sin(\pi x) 
\end{align*}
The second of these is useful because the two trigonometric functions go through a full period as $x$ goes between 0 and 1.

\item \textbf{The Spectrum}.  For Fourier series and transforms, the transformation of the independent variable, $x$, is to a new variable, $\xi$, that measures the frequency of the component of the sine and cosine functions used to represent it.  The relationship between the \emph{amplitude} of the trigonometric functions versus the corresponding frequency variable, $\xi$.  

\item \textbf{Convolution}.\indexme{convolution}  An integral operation that integrates one function against another by a continuous process of shifting one of the two functions.  For two functions $f(x)$ and $g(x)$, the convolution would be defined by on interval $x\in[a,b]$ by
\begin{equation*}
(f*g)(x) = \int_{a}^{b} f(x)g(x-z) \, dz
\end{equation*}
Although convolutions on finite intervals do occur \citep{ljubarskiui1977convolution}, they are more routinely found on the semi-infinite or infinite lines. For Fourier series, the convolution is on the infinite real line, and we have the following definitions
\begin{equation*}
(f*g)(x) = \int_{-\infty}^{\infty} f(z)g(x-z) \, dz = \int_{-\infty}^{\infty} f(x-z)g(z) \, dz
\end{equation*}

\item \textbf{Square-Integrable.}  A function is said to be square-integrable on $x\in[a,b]$ if 
\begin{equation*}
    F(x) = \int_a^b f^2(x) \, dx < \infty
\end{equation*}
Functions that are square-integrable belong to a space called the Hilbert space, $L^2$.  Square-integrable functions are an important class of functions in many areas of analysis.

\item \textbf{Tempered Distribution.}  A version of distribution theory that allows the Fourier transform to be generalized to more kinds of functions.  While the conventional requirement for a Fourier transform to exist is that the function be \emph{square integrable} on $x\in(-\infty,\infty)$, the theory of tempered distributions allows the extension of Fourier Transforms to (1) generalized functions (such as the delta function), and (2) functions that do not grow \emph{too fast} (and this means only that they must generally grow slower than an exponential function).  Although the theory of tempered distributions is not covered in this text, their existence is needed to explain the use of several functions that arise naturally in applications of the Fourier transform.

\end{itemize}

\section{Return to the Fourier Series: The Complex Fourier Series}\indexme{Fourier Series!Complex}

Recall that every function defined over a symmetric interval around zero, $x\in[-L,L]$, can be expressed as the sum of an \emph{even} plus an \emph{odd} function

\begin{equation}
    f(x)=f_{even}(x)+f_{odd}(x)
\end{equation}
Also recall, that over the same symmetric interval, the Fourier sine series is able to represent only odd functions, whereas the Fourier cosine series is able to represent only even functions.  Thus, in general, we need both the Fourier sine and cosine series to represent a function (one to represent $f_{even}$ and one to represent $f_{odd}$).  The expressions developed previously were of the form

\begin{subequations}
\begin{align}
    f(x) &= A'_0+\sum_{n=1}^{n=\infty} A_n \cos\left(n \pi \tfrac{x}{L}\right) +\sum_{n=1}^{n=\infty} B_n \sin\left(n \pi \tfrac{x}{L}\right) \label{fs1}\\
    A'_0 &= \frac{1}{2L} \int_{x=-L}^{x=L} f(x)\, dx~~(\textrm{the average})\\
    A_n &= \frac{1}{L} \int_{x=-L}^{x=L} f(x) \cos\left(n \pi \tfrac{x}{L}\right)\, dx \\
    B_n &= \frac{1}{L} \int_{x=-L}^{x=L} f(x) \sin\left(n \pi \tfrac{x}{L}\right)\, dx \label{fs4}
    \end{align}
\end{subequations}
These results summarize what we have developed for the general Fourier series on a symmetric interval around $x=0$.  

Although generally there is no need to do so, we can also express these results entirely using the complex exponentials defined by Euler's equation.  Recall, Euler's equation stated

\begin{equation}
    e^{i n \pi x/L} = \cos\left(n \pi \tfrac{x}{L}\right)+ i\sin\left(n \pi \tfrac{x}{L}\right)
\end{equation}
Making linear combinations of this result allows us to define the following relationships

\begin{align}
    \cos\left(n \pi \tfrac{x}{L}\right) &= \frac{e^{i n \pi x/L}+e^{-i n \pi x/L}}{2} \\
    \sin\left(n \pi \tfrac{x}{L}\right) &=\frac{e^{i n \pi x/L}-e^{-i n \pi x/L}}{2i}
\end{align}
Substituting these into Eqs.~\eqref{fs1}-\eqref{fs4} is algebraically messy, but otherwise not difficult.  The result is

\begin{subequations}
\begin{align}
    f(x) &= A'_0+\sum_{n=1}^{n=\infty} A_n \frac{e^{i n \pi x/L}+e^{-i n \pi x/L}}{2} +\sum_{n=1}^{n=\infty} B_n \frac{e^{i n \pi x/L}-e^{-i n \pi x/L}}{2i} \label{fse1}\\
    A'_0 &= \frac{1}{2L} \int_{x=-L}^{x=L} f(x)\, dx~~(\textrm{the average})\\
    A_n &= \frac{1}{L} \int_{x=-L}^{x=L} f(x) \frac{e^{i n \pi x/L}+e^{-i n \pi x/L}}{2}\, dx \\
    B_n &= \frac{1}{L} \int_{x=-L}^{x=L} f(x) \frac{e^{i n \pi x/L}-e^{-i n \pi x/L}}{2i}\, dx \label{fse4}
\end{align}
Or, collecting like exponential terms

\begin{align}
    f(x) &= A'_0+\sum_{n=1}^{n=\infty} \frac{A_n-i B_n}{2} e^{i n \pi x/L}+\sum_{n=1}^{n=\infty} \frac{A_n+i B_n}{2} e^{-i n \pi x/L}\\
    A'_0 &= \frac{1}{2L} \int_{x=-L}^{x=L} f(x)\, dx~~(\textrm{the average})\\
 \frac{A_n-i B_n}{2}& =   \frac{1}{L} \int_{x=-L}^{x=L} f(x) e^{-i n \pi x/L} \, dx  \\
  \frac{A_n+i B_n}{2} & =   \frac{1}{L} \int_{x=-L}^{x=L} f(x) e^{i n \pi x/L} \, dx 
\end{align}
\end{subequations}
Finally, note that it is convenient to rename the coefficients being used as follows

\begin{align}
    f(x) &= C_0+\sum_{n=1}^{n=\infty} C_n e^{i n \pi x/L}+\sum_{n=1}^{n=\infty} D_n e^{-i n \pi x/L}\\
    C_0 &= \frac{1}{2L} \int_{x=-L}^{x=L} f(x)\, dx~~(\textrm{the average})\label{fsC0}\\
    C_n& =   \frac{1}{2L} \int_{x=-L}^{x=L} f(x) e^{-i n \pi x/L} \, dx \label{fsCn}  \\
  D_n & =   \frac{1}{2L} \int_{x=-L}^{x=L} f(x) e^{i n \pi x/L} \, dx \label{fsDn}
\end{align}
As a last step, note that $C_n$ and $D_n$ are defined in the same way, with the exception that the exponents are of opposite sign.  However, if we sum the expression for $C_n$ from $n=-1$ to $n=-\infty$ we have the following handy relationship

\begin{equation}
    \sum_{n=1}^{n=\infty} D_n e^{-i n \pi x/L}=\sum_{n=-1}^{n=-\infty} C_n e^{i n \pi x/L}\label{fsCD}
\end{equation}
Combining Eqs.~\eqref{fsC0}, \eqref{fsCn}, \eqref{fsDn}, and \eqref{fsCD}, the result is a single sum expressing the entire series

\begin{align}
    f(x) &=\sum_{n=-\infty}^{n=-1} C_n e^{i n \pi x/L}+C_0 +\sum_{n=1}^{n=\infty} C_n e^{i n \pi x/L}\\
    C_0 &= \frac{1}{2L} \int_{x=-L}^{x=L} f(x)\, dx\\
    C_n &= \frac{1}{2L} \int_{x=-L}^{x=L} f(x)e^{-i n \pi x/L}\, dx
\end{align}
Note: As with the Fourier series, it is not always true that the expression for $C_n$ evaluates to the proper quantity when $n=0$.  However, when it does, we can write the series in an even more compact form as follows

\begin{subequations}
\begin{align}
    f(x) &=\sum_{n=-\infty}^{n=\infty} C_n e^{i n \pi x/L}\\
    C_n &= \frac{1}{2L} \int_{x=-L}^{x=L} f(x)e^{-i n \pi x/L}\, dx~~n=\ldots -2, -1, 0, 1, 2 \ldots
\end{align}
\end{subequations}

As a reminder, the results above are simply the results for the general Fourier series for a symmetric interval around zero, $x\in[-L,L]$, where we have expressed the sine and cosine functions by the complex exponential.  While the forms of the expressions are simpler looking, they mean exactly the same thing as the interpretation in terms of sine and cosine series.  Upon summation, the result for the Fourier series expressed in this form should be identical to that computed using the sine and cosine series directly.  Particularly important to note here is that the resulting series should contain no complex numbers.  Even though complex numbers arise in intermediate steps, in the resulting final expression only real numbers should appear (assuming that the original function was a real one to begin with!)  Like the conventional Fourier series studied in Chapter \ref{FS_1}, the functions must meet the Dirichlet conditions outlined 

Note that because of the change of bounds for the summation, we still have to compute exactly the same number of coefficients for any finite approximation to the series (that number being $2N+1$, where $N$ is the number of terms that we decide to compute, and the plus one comes from the $C_0$ constant).  In addition, some find this form of the expression to be less intuitive than the one involving sine and cosine functions (and not involving any complex numbers).  

\begin{svgraybox}
\begin{example}[A Fourier Series Using the Complex Series Form]
For the function $f(x)=e^x$ on $[-1,1]$, the conventional Fourier series is found from
\begin{align}
    A'_0&=\frac{1}{2}\int_{x=-1}^{x=1} e^x \,dx  = \sinh[1] \\
    A_n&=\frac{1}{2}\int_{x=-1}^{x=1} e^x \sin[n \pi x]\, dx = \frac{2  (-1)^n \sinh (1)}{1+\pi ^2 n^2}\\
    B_n&=\frac{1}{2}\int_{x=-1}^{x=1} e^x \cos[n \pi x] \,dx -\frac{2 \pi  (-1)^n n \sinh (1)}{1+\pi ^2 n^2}\\
\intertext{with }
e^x & = A'_0+\sum_{n=1}^{n=\infty} A_n \cos(n \pi x) +\sum_{n=1}^{n=\infty} B_n \sin(n \pi x)
\end{align}
The series solution for this is plotted in Fig.~\ref{ex}.  The corresponding complex exponential series can be found as follows.  

{
\centering\includegraphics[scale=.7]{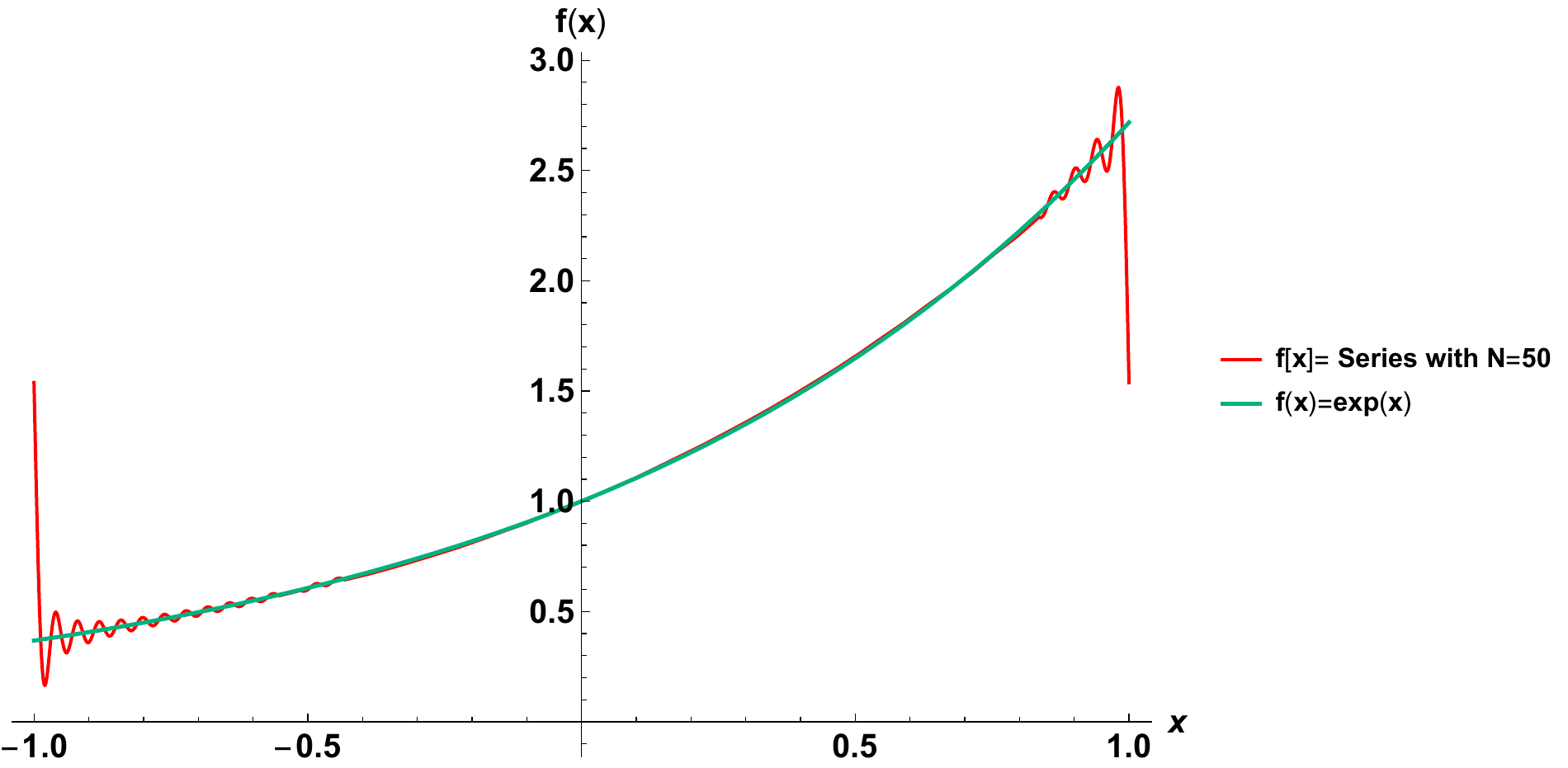}
\captionof{figure}{Fourier series for $f(x)=e^x$ with $N=50$ terms.}
\label{ex}  
}

\begin{align*}
C_n &= \frac{1}{2} \int_{x=-1}^{x=1} e^x e^{-i n \pi x}\, dx~~n=\ldots -2, -1, 0, 1, 2 \ldots\\
&= \frac{(-1)^n (L+i \pi  n) \sinh (L)}{L^2+\pi ^2 n^2}
\intertext{Thus, the resulting series is}
e^x &= \sum_{n=-\infty}^{n=\infty} \frac{(-1)^n (L+i \pi  n) \sinh (L)}{L^2+\pi ^2 n^2} e^{i \pi x}
\end{align*}
For plotting purposes, we can re-express this result as

\begin{align*}
e^x &= \sum_{n=-\infty}^{n=\infty} \frac{(-1)^n  \sinh (L)}{L^2+\pi ^2 n^2} (L+i \pi  n)\left[\cos(n \pi x) + i \sin(n \pi x)\right]
\end{align*}
Note that the final product in this sum involves the multiplication of two complex numbers.  This multiplication when carried out involves many canceling terms.  With substantial algebra, one can show that all of the imaginary components cancel entirely (as they must).  The result is

\begin{align*}
e^x &= \sum_{n=-\infty}^{n=\infty} \frac{(-1)^n  \sinh (L)}{L^2+\pi ^2 n^2} (L+i \pi  n)\left[\cos(n \pi x) + i \sin(n \pi x)\right]\\
&=\sinh (1)+\frac{2 \pi  \sinh (1) \sin (\pi  x)}{1+\pi ^2}-\frac{4 \pi  \sinh (1) \sin (2 \pi 
   x)}{1+4 \pi ^2}
   +\frac{6 \pi  \sinh (1) \sin (3 \pi  x)}{1+9 \pi ^2}-\frac{8 \pi 
   \sinh (1) \sin (4 \pi  x)}{1+16 \pi ^2}\\
   &\hspace{15mm}+\frac{10 \pi  \sinh (1) \sin (5 \pi  x)}{1+25
   \pi ^2}\\
   &\hspace{15mm}-\frac{2 \sinh (1) \cos (\pi  x)}{1+\pi ^2}+\frac{2 \sinh (1) \cos (2 \pi 
   x)}{1+4 \pi ^2}
   -\frac{2 \sinh (1) \cos (3 \pi  x)}{1+9 \pi ^2}+\frac{2 \sinh (1) \cos
   (4 \pi  x)}{1+16 \pi ^2}\\
   &\hspace{15mm}-\frac{2 \sinh (1) \cos (5 \pi  x)}{1+25 \pi ^2}+\ldots
\end{align*}
The plot of this series is identical to that shown previously by Fig.~\ref{ex}.  In summary, while the complex exponential representation of the Fourier series is \emph{always} equivalent to the series expressed in trigonometric functions, it is not necessarily more convenient for all purposes.  In particular, if one wants to actually sum the series, the complex exponential form is generally not more convenient.  However, for developing the extension of the Fourier series to the real line (generating the Fourier transform), the complex form will have some advantages. 
\end{example}
\end{svgraybox}

\section{The Fourier Series as a Transform}\indexme{Fourier transform!finite Fourier transform}
Before developing the Fourier transform on the real line, it is instructive to first revisit the definition of the Fourier series, but viewed in the context of a transform.  In this context, the Fourier series is sometimes called the \emph{Finite Fourier Transform} (cf. \citet{churchill1972}).  This renaming represents only a change in perspective of the operations involved; functionally, nothing about the details of the Fourier series are changed when thinking of it as a transformation applied to a finite interval.  To start, we first review the concept of a linear transform.

\subsection{Linear transforms}\indexme{Fourier transform!linearity}

{\noindent \bf Definition.} In the context of function spaces, a linear transformation, $\mathscr{L}(\cdot)$, takes a function as an argument (its domain), and generates a new function as its output (its range).  Explicitly, this is defined by
\begin{equation}
    \mathscr{L}[f] = g
\end{equation}
Note that the independent variables associated with the $f$ and $g$ may be different.

\begin{svgraybox}
\begin{example}[Linear Transformations on Functions]
Consider the following two functions
\begin{align*}
    f(x) &= 4x\\
    g(y) &= 6y
\end{align*}
Then the linear operation $f(g)$ is defined by
\begin{equation*}
    f[g(y)]= 4(6y) = 24 y
\end{equation*}
This can be considered a linear transformation between one function space to another.

As a second example, consider the operation of integration of a function.  The integration process can be thought of as a linear transform (or, equivalently, linear operator).  Note the following

\begin{equation*}
    \int (\alpha f(x)+\beta g(x))\, dx = \alpha \int f(x)\, dx+\beta \int  g(x)\, dx
\end{equation*}
The property that a transformation preserves the operations of addition and multiplication, as illustrated above, define a linear transformation.  Thus, integration can be thought of as is a linear transformation.  As a concrete example consider the following.
\begin{align*}
    f(x) &= 4x^2+2\\
    \mathscr{L}(f) &= \int_{x=0}^{x=1} f(x) \, dx\\
    \intertext{Therefore}
    \mathscr{L}(4x^2+2)&=  \int_{x=0}^{x=1} (4x^2+2) \, dx\\
    &= 4\int_{x=0}^{x=1} x^2+ \, dx+  2\int_{x=0}^{x=1} 1 \, dx\\
    &= \frac{10}{3}
\end{align*}
The process of applying the integration as an operation is clearly a linear one.
\end{example}
\end{svgraybox}

\subsection{The Finite Fourier Transform}

Now that we have had some reminders about the concept of a linear transform, we can revisit the Fourier series, but viewing it as a transform of one function (in its function space) to another function (in a different function space).  Because this mapping between functions spaces is reversible, we will be examining the Fourier series as a set of two operations: (1) a forward transform that take a function in a real variable and sends it to a new function defined by the sequence of amplitudes and frequencies, and (2) an inverse function that takes the transformed function of amplitudes and frequencies, and sends it back to being a conventional function in the original independent variable.

\begin{figure}[t]
\sidecaption[t]
\centering
\includegraphics[scale=.4]{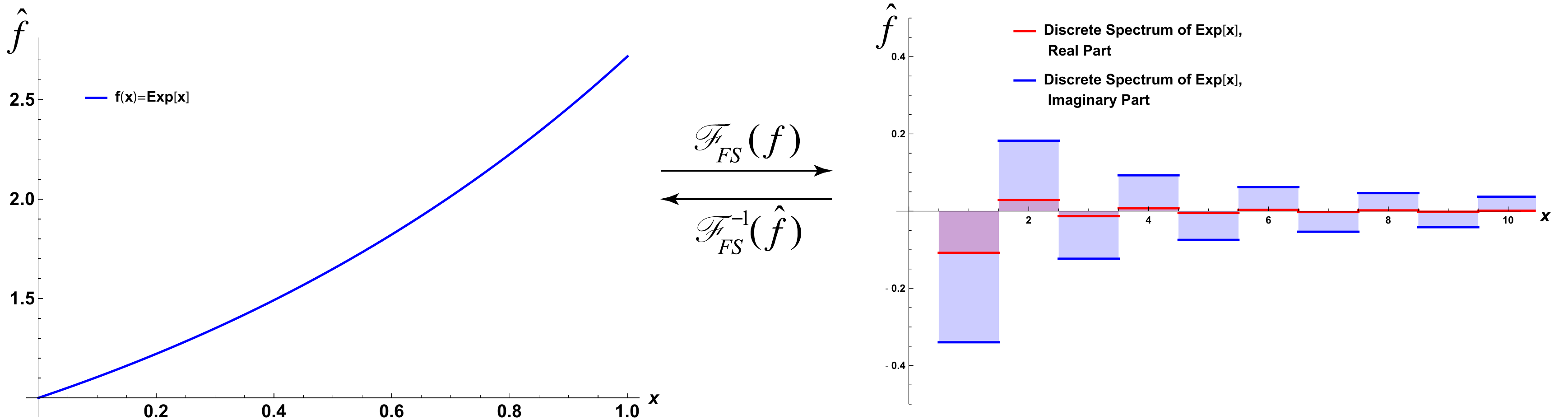}
\caption{The function $f(x)=\exp(x)$ (left) and its finite Fourier transform spectrum (right).}
\label{finiteFT}       
\end{figure}

Recall, the definition of the complex Fourier series given previously is
\begin{subequations}
\begin{align}
    f(x) &=\sum_{n=-\infty}^{n=\infty} C_n e^{i n \pi x/L}\\
    C_n &= \frac{1}{2L} \int_{x=-L}^{x=L} f(x)e^{-i n \pi x/L}\, dx~~n=\ldots -2, -1, 0, 1, 2 \ldots
\end{align}
\end{subequations}
We can actually write this as a single expression, as follows
\begin{align}
    f(x) &=\sum_{n=-\infty}^{n=\infty} \left(\frac{1}{2L} \int_{x=-L}^{x=L} f(x)e^{-i n \pi x/L}\, dx \right)e^{i n \pi x/L}\label{fourierint}
\end{align}
Written in this way, we can interpret the Fourier series as a sequence of two operations.
\begin{enumerate}
    \item First, there is an operation that takes a function, and generates a transformed version of the function.  Noting that the function can be thought of as the pair $(x, f(x))$, then the process of finding the Fourier amplitudes generates a new pair $(n, C_n)$ that represents the original function, but in a transformed space (with domain $n$, and range $C_n$).  Symbolically, we suppose we use $\mathscr{F}_{FS}(\cdot)$ to represent this transformation.  Then, we could write either of the following to represent the transformation process
    \begin{align*}
        &(x,f(x))\xrightarrow{{{\mathscr{F}_{FS}}}}(n,{C_n})\\
        &\mathscr{F}_{FS}\left[f(x)\right] = C_n(n)
    \end{align*}
    By tradition, when describing the Fourier series as a linear transform, the following notation (or any one of many similar notations) is adopted: $C_n(n)  \Leftrightarrow \hat{f}(n)$.  This is so that the form of the transformed function, $\hat{f}(n)$ , provides is a reminder that it came originally from the function $f(x)$.  Thus, the two transforms given above are written in the traditional manner by
        \begin{align*}
        &(x,f(x))\xrightarrow{{{\mathscr{F}_{FS}}}}(n,\hat{f}(n))\\
        &\hspace{15mm}\textrm{or}\\
        &\mathscr{F}_{FS}\left[f(x)\right] = \hat{f}(n)
    \end{align*}
    With these definitions in place, we can think of the finite Fourier transform as being the inner integral in Eq.~\eqref{fourierint}.  This yields the definition of the forward transform
     
    \begin{equation}
    \boxed{
      ~\mathscr{F}_{FS}\left[f(x)\right] = \hat{f}(n) = \frac{1}{2L}
        \int_{x=-L}^{x=L} f(x)e^{-i n \pi x/L}\, dx~
        }
    \end{equation}

    \item Second, there is an \emph{inverse} operation that returns a function, $\hat{f}(n)$,  in Fourier space back to a real function of the original variable, $f(x)$.  
    
   \begin{align*}
        &(n,\hat{f}(n)) \xrightarrow{{{\mathscr{F}^{-1}_{FS}}}} (x,f(x))\\
         &\hspace{15mm}\textrm{or}\\
        &\mathscr{F}^{-1}_{FS}\left[\hat{f}(n)\right] = f(x)
    \end{align*}
    
    Again, in reference to Eq.~\eqref{fourierint}, note that if we adopt the notation given directly above for $\hat{f}$, then we have a definition for the \emph{inverse transform} $\mathscr{F}_{FS}^{-1}$
    
    \begin{equation}
    \boxed{
       ~\mathscr{F}_{FS}^{-1}\left[ \hat{f}(n)\right]= f(x) =\sum_{n=-\infty}^{n=\infty}  \hat{f}(n) e^{i n \pi x/L}~
       }
    \end{equation}
    
\end{enumerate}
An example showing the function $f(x)=e^x$ and its Finite Fourier transform spectrum is given in Fig.~\ref{finiteFT}.

\section{The Fourier Transform}\indexme{Fourier transform}\label{fourtrans}

We are now in a position where we can more carefully examine the idea of how to define the Fourier series on the real line, i.e., for the open interval $~I= \lim_{L\rightarrow\infty}\left[-L < x < L\right]$.  As mentioned above, the Fourier transform is an extension of the idea of the Fourier series, but extended to the entire real line.  The Fourier transform was an original invention of Fourier; before him, the notion of integral transforms was not developed.  While the Fourier series opened entirely new concepts regarding the very definition of the function, the Fourier integral extended these ideas, and had enormous impact across mathematics ranging from integration theory to the solutions to partial differential equations.

In the following,  the Fourier transform is derived from the Fourier series.  This derivation proceeds in a way that reasonably matches our intuition (i.e., if the Fourier transform is the extension of the Fourier series to the entire real line, then we should be able to develop it from the Fourier series by letting $L\rightarrow\infty$).  However, there are many technical issues in the \emph{formal} derivation of the Fourier transform that we will not examine in detail.  Therefore, the derivation below should be considered the outline to a more careful, detailed proof.  Regardless, the development has substantial intuitive value, and the presentation is frequently offered as the de facto derivation of the Fourier transform

\subsection{Derivation of the Fourier Integral}\index{Fourier transform!derivation}\indexme{Fourier transform!Fourier integral}

To start, note that by definition of the complex Fourier series (or, equivalently, the finite Fourier transform) for any admissible function $f$ on $x\in[-L,L]$ is given by

\begin{align}
    f(x) &= \sum_{n=-\infty}^{n=\infty} C_n e^{i n \pi x/L}  \\
    C_n &= \frac{1}{2L} \int_{x=-L}^{x=L} f(x) e^{-i n \pi x/L}\, dx
\end{align}
Combining these into a single expression (as we did for the finite Fourier transform), we find

\begin{align}
        f(x)& = \sum_{n=-\infty}^{n=\infty} \left[ \frac{1}{2L} \int_{x=-L}^{x=L} f(x) e^{-i n \pi x/L} \, dx\right] e^{i n \pi x/L} \label{FT001}\\
    \intertext{Now, let }
    \xi_n &= \frac{n \pi}{L} \,\, \Leftrightarrow \,\, n= \frac{\xi_n L}{\pi}\nonumber\\
     \Delta \xi_n  &=  \left(\frac{(n+1) \pi}{L} - \frac{n \pi}{L}\right)=\frac{\pi}{L}
\end{align}
Substituting these into Eq.~\eqref{FT001} gives the following sum
\begin{align}
        f(x)& = \sum_{n=-\infty}^{n=\infty} \left[ \frac{1}{2\pi} \int_{x=-L}^{x=L} f(x) e^{-i \xi x} \, dx\right] e^{i \xi x} \Delta \xi_n
\end{align}
This sum is now in the form of a conventional Riemann sum (i.e., the sum that was covered in introductory calculus and led to the definition of the integral).  The integral in this case is obtained in the limit as $L\rightarrow\infty$ (or, equivalently, as $\Delta \xi_n \rightarrow 0$)
\begin{align}
        f(x)& ={\mathop {\lim }\limits_{\begin{subarray}{l} 
  \Delta \xi  \to 0 \\ 
  L \to \infty  
\end{subarray} }} \left( \sum_{\xi L/\pi=-\infty}^{\xi L/\pi=\infty} \left[ \frac{1}{2\pi} \int_{x=-L}^{x=L} f(x) e^{-i \xi x} \, dx\right] e^{i \xi x} \Delta \xi_n\right)
\end{align}
The final result is
\begin{align}
        f(x)& = \frac{1}{\sqrt{2 \pi}}\int_{\xi=-\infty}^{\xi=\infty} \left[ \frac{1}{\sqrt{2 \pi}} \int_{x=-\infty}^{x=\infty} f(x) e^{-i \xi x} \, dx\right] e^{i \xi x} \, d\xi\\
\end{align}
Note that we have distributed the factor $1/(2 \pi)$ in two parts, so that the forward and reverse transforms have symmetry in the constants in front of them.  While the derivation above is compelling, there is at least one detail that was not described fully.  We have apparently let the index, $n$, be transformed from an index set over the integers to a real number defined on the entirety of the real line.  This was somewhat hidden by the transformation of variables from $n$ to $\xi_n$, but this transformation did occur.  While we did not justify this step rigorously, such justifications are available (see \citet{lanczos1966}).

As with the finite Fourier transform, we define both a forward and inverse transform as follows.

\begin{equation}
\boxed{
         ~~~~~~\hat{f}(\xi)=\mathscr{F}[f(x)]=\frac{1}{\sqrt{2 \pi}} \int_{x=-\infty}^{x=\infty} f(x) e^{-i \xi x} \, dx~~~~~~}\label{fouriertransformdef}\\
\end{equation}
\begin{equation}
          \boxed{~~~~
        f(x) =\mathscr{F}^{-1}[f(x)]= \frac{1}{\sqrt{2 \pi}}\int_{\xi=-\infty}^{\xi=\infty} \hat{f}(\xi) e^{i \xi x} \, d\xi~~~~}
\end{equation}

While this transform has direct correspondence with the finite Fourier transform, there are some additional considerations that must be addressed.  These are covered briefly as follows.

\begin{enumerate}
    \item For the Fourier transform to exist, we must have the condition that the function $f$ be \emph{square-integrable} on $x\in(-\infty,\infty)$.  This means only that
    \begin{equation}
    \boxed{~~~~
       \int_{x=-\infty}^{x=\infty} \lvert f(x)\rvert^2 \, dx < \infty~~~~}
    \end{equation}
    \emph{and} that this finite number is unique.  
    While we may not have occasion to have to use this directly (since most of our Fourier transforms will be available in tables), it is good to understand that only functions that are finite in this sense are transformable.  Later in the chapter, we will illustrate how we can relax this requirement for a few special cases such as the delta function, the Heaviside function, and the constant function.  
    
    \item The smoothness of the resulting transformations is governed by the smoothness of the function being transformed, in a manner that is completely analogous to that for Fourier series.  In other words, the representation
    
    \begin{align}
        f(x)& =\mathscr{F}^{-1}[f(x)]= \frac{1}{\sqrt{2 \pi}}\int_{\xi=-\infty}^{\xi=\infty} \hat{f}(\xi) e^{i \xi x} \, d\xi
\end{align}
will converge uniformly for functions that are analytic.  For functions that are continuous, the Fourier transform will also converge pointwise.  For functions with jump discontinuities, the Fourier transform will converge in the same way that the Fourier series do; that is, it converges to the value

    \begin{equation}
        f(x) = \frac{f(x^+)+f(x^-)}{2}
    \end{equation}

\item Because the delta function has a representation as a Fourier series, it also has a representation in the transformed space.  In other words, the transform of generalized functions like the delta function can be given meaning.  In this sense, the generalized functions act the same as any other function in the transformed space.
 \end{enumerate}

\subsection{Notation and Restrictions on the Fourier Transform}\indexme{Fourier transform!notation}\indexme{Fourier transform!restrictions}

There is probably no single area in all of applied mathematics where the adoption of multiple kinds of notation have created as much confusion as for the case of the Fourier transform.  There are several items regarding notation that users of Fourier transform tables need to be aware of.  These are as follows.

\begin{enumerate}
    \item First, the coefficient in front of the integral defining the Fourier transform, $\mathscr{F}$ and the inverse, $\mathscr{F}^{-1}$, is given for our case as the symmetric quantity $1/\sqrt{2 \pi}$ (i.e., the same constant appears in front of the transform and its inverse).  Some texts prefer to leave the factor $1/(2 \pi)$ in front of the transform defining $\mathscr{F}[f(x)]$, and a coefficient of 1 for the inverse transform. 
    \item  In our definition, the forward transform is defined by the set of exponentials $e^{-i \xi x}$, and the inverse by $e^{i \xi x}$.  Some texts switch the definition of which transform is defined by the basis functions $e^{-i \xi x}$.
    \item There are multiple variables used for the Fourier variable, $\xi$.  This variable is frequently called the \emph{wave space variable} (or \emph{wave vector} in multiple dimensions).  The symbols $k$, $\alpha$, and $\omega$ are also frequently adopted for the variable that we have denoted by $\xi$.
    \item The notation indicating the Fourier transform of a function can vary. While the most frequent notation is to indicate the transform by a ``hat" (i.e., $\hat{f}=\mathscr{F}[f]$), other notations frequently occur.  Also widely used is the convention of using a lower case letter and corresponding capital letter to indicate the transform pair, (i.e., $F=\mathscr{F}[f]$).  Yet others use additional notations, such as an overbar (i.e., $\overline{f}=\mathscr{F}[f]$; we will reserve this notation specifically for the \emph{Laplace transform}) to indicate the transform.  
\end{enumerate}
The primary difficulty with the use of different formalisms for the Fourier transform is primarily related to the use of tables for looking up the appropriate transform of a function.  In short, the results listed in any table will be a function of the definitions for the transforms themselves and for the notations associated with the transforms.  While for cases where the fundamental transform definitions are identical, the primary problem is only a change in variable names (e.g., if a table were to use $\alpha$ instead of $\xi$).  However, for tables that adopt a different definition for the transform-inverse transform pair, the differences are more substantial.  The most important point is to be sure that if one uses a table (or software, such as Mathematica) to determine a Fourier transform, then the definitions of the transform pairs must match those that the definitions used in the table.

\textbf{Restrictions to the Fourier Transform.}  There are few restrictions on the Fourier transform, at least ones that we need to be concerned with.  The conventional requirement for the Fourier transform to exist is that the function being transformed is \emph{square integrable}.  Thus, for the transform of $f(x)$ to exist, we need to assure that \indexme{square integrable}

\begin{equation}
    I = \int_{-\infty}^{\infty} f^2(x) \, dx < \infty
\end{equation}
A few functions that occur routinely in applications (e.g., the Heaviside function) are not square integrable.  For those, there is an extension to the theory of Fourier transforms known as the theory of tempered distributions.  We will not discuss this theory in any detail, but it is important to at least note that such a framework exists.  In this case, noting this framework, we can make sensible intepretations for a few functions (such as the Heaviside) that would not otherwise be transformable.

\section{Fourier Transform Identities}\indexme{Fourier transform!identities}

The Fourier transform has some interesting properties (more correctly, identities) that occur between the real and transformed spaces.  It is, in part, these properties that make the Fourier transform useful in applications.  In summary, the Fourier transform of the following functional forms has a simple representation in transform space: (1) functions where the independent variable is shifted by a constant, (2) functions where the independent variable is multiplied by a constant scaling factor, (3) functions that are defined by the derivative operation.  Each of these cases is briefly discussed below.

\subsection{Fourier Transform of Shifted Functions}\indexme{Fourier transform!shifted functions}

\begin{figure}[t]
\sidecaption[t]
\centering
\includegraphics[scale=.4]{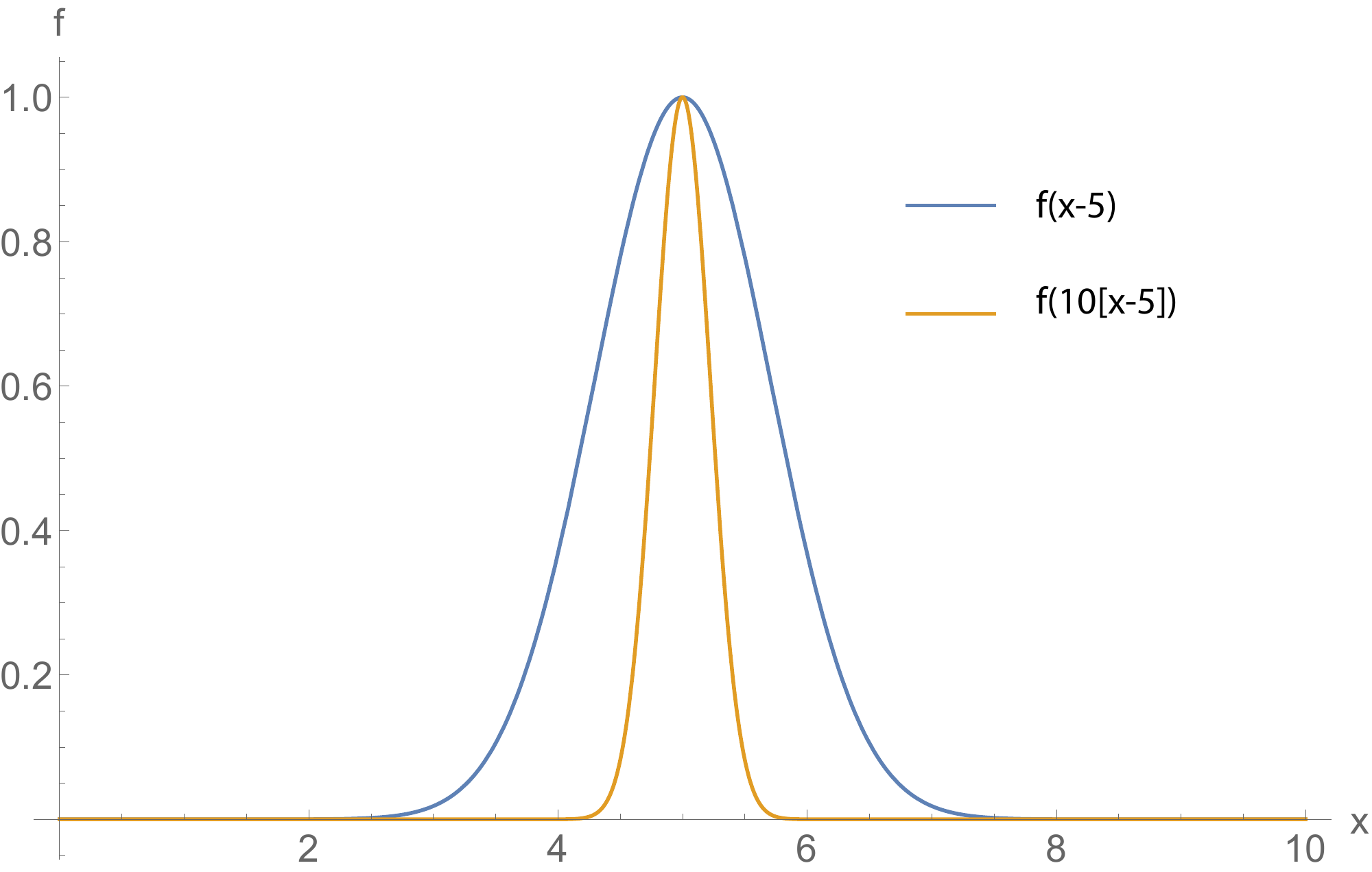}
\caption{The function $f$ upon both shifting and scaling so that it takes the form $f(a[x-c])$.  Here, the constant $a$ controls the scale of the function (the larger the value of $a$, the narrower the function), and $c$ controls the shift (subtracting a constant shifts functions to the right; adding a constant shifts functions to the left).}
\label{shiftscale}       
\end{figure}

Recall, a function $f(x)$ is shifted to the right when a constant $c$ is subtracted from the independent variable to form $f(x-c)$ (see Fig.~\ref{shiftscale}).  
When one transforms a function of the form $f(x-a)$, the following identity holds

\begin{equation}
    \mathscr{F}[f(x-a)] = e^{-i a\xi}\hat{f}(\xi)
\end{equation}
This identity is easily proven by making the variable transform $z = x-c$, and substituting this into the integral definition of the Fourier transform.

\subsection{Fourier Transform of Scaled Functions}\indexme{Fourier transform!scaled functions}

Scaling the independent variable of a function by a constant expands or compresses the functional relationship to the domain.  In Fig.~\ref{shiftscale}, the effect of multiplying the independent variable by a constant $a$ is illustrated.  Intuitively, one might expect the compression of the domain that such a scaling (with $a>1)$ creates to affect the Fourier representation by increasing the high-frequency components (i.e., the narrower function requires higher frequency components to resolve).  This is essentially what happens.  The following relationship holds for Fourier transforms of a scaled function

\begin{equation}
            \mathscr{F}[f(ax)] = \frac{1}{a}\hat{f}\left(\frac{x}{|a|}\right)
\end{equation}
This identity can be proved by the variable transformation $z = ax$.and substituting this into the integral definition of the Fourier transform.

\subsection{Fourier Transform of Derivatives}\indexme{Fourier transform!derivatives}

From the perspective of solving partial differential equations, the Fourier transform has \emph{one primary useful feature}: it turns derivatives (of the transformed independent variable) into algebra.  This is very easy to prove using integration by parts; this derivation is left as an exercise.  Assuming that we \emph{require} the following two conditions

\begin{align*}
&B.C.1&    u(x,t) &\rightarrow 0 ~~as~ |x|\rightarrow \infty \\
&B.C.2&    \left.\frac{\partial u}{\partial x}\right|_{(x,t)}&\rightarrow 0 ~~as~ |x|\rightarrow \infty
\end{align*}

Then, we have the result
\begin{align}
    \mathscr{F}\left[ \frac{\partial u}{\partial x}\right] &=\frac{1}{\sqrt{2 \pi}} \int_{x=-\infty}^{x=\infty} \left.\frac{\partial u}{\partial x}\right|_{(x,t)} e^{-i \xi x} \, dx \nonumber\\
    &=i \xi \hat{u}
\end{align}
This result can be iterated for higher-order derivatives (assuming that for each such higher derivative, the derivative is constrained so that it goes to zero as $|x|\rightarrow\infty$)

\begin{equation}
    \mathscr{F}\left[ \frac{\partial^2 u}{\partial x^2}\right] = (i \xi)^2 \hat{u} = -\xi^2 \hat{u}
\end{equation}
Note that the following action of the transform on derivatives that are not with respect to a transformed variable

\begin{align}
    \mathscr{F}\left[ \frac{\partial u}{\partial t}\right] &=\frac{1}{\sqrt{2 \pi}} \int_{x=-\infty}^{x=\infty} \left.\frac{\partial u}{\partial t}\right|_{(x,t)} e^{-i \xi x} \, dx \nonumber\\
    &=\frac{\partial }{\partial t}\frac{1}{\sqrt{2 \pi}} \int_{x=-\infty}^{x=\infty} u(x,t) e^{-i \xi x} \, dx \nonumber\\
    &= \frac{\partial \hat{u}}{\partial t}
\end{align}

While these results may not necessarily seem significant, they are actually is exceptionally useful.  The following examples helps illustrate the transform of the derivative.

\begin{svgraybox}
\begin{example}[Solution to the first-order convection equation using the Fourier transform]\label{firstorderwaveFT}\indexme{Fourier transforms!First-order wave equation}
Consider the following problem, representing the pure convection of a pulse of solute at velocity $c$ with concentration given by $u$.

\begin{align*}
    &&\frac{\partial u}{\partial t}&= -c \frac{\partial u}{\partial x}\qquad -\infty < x < \infty\\
&B.C.1&    u(x,t) & \textrm{~~remains bounded for all $x$ and $t$} &&\\
&I.C.& u(x,0)&=\exp(-x^2)
\end{align*}
This problem represents the translation in the $+x$ direction of the initial pulse, unchanged in shape, at velocity $c$.  From what we know about this problem from our study of separation of variables, we expect the solution to be of the form $u(x,t)=\exp[-(x-c t)^2]$.  However, our goal here is to \emph{derive} this result.  Toward that end, we take the Fourier transform of both sides of the PDE and of the initial condition.  This gives

\begin{align*}
    &&\frac{d \hat{u}}{d t}&= -i \xi c \hat{u}\qquad -\infty < x < \infty\\
&B.C.1&    \hat{u}(\xi,t) & \textrm{~~remains bounded for all $\xi$ and $t$} &&\\
&I.C.& u(\xi,0)&=\frac{1}{\sqrt{2}}\exp(-\xi^2/4)
\end{align*}
The solution to this problem reduces to a first order, integrable ODE.  The solution is

\begin{equation*}
    \hat{u}(\xi,t) =\frac{1}{\sqrt{2}}\exp(\frac{-\xi^2}{4}) \exp\left[ - i c t \xi \right]
\end{equation*}
To complete the problem, we need to invert the transform.  Taking the inverse of both sides yields

\begin{equation*}
    \mathscr{F}^{-1}[\hat{u}(\xi,t)] = \mathscr{F}^{-1}\left[\exp\left( - i c t \xi \right) \frac{1}{\sqrt{2}}\exp(\frac{-\xi^2}{4}) \right]
\end{equation*}
Here, the easiest way to proceed is to recognize that the first exponential is one that shifts the following transformed function.  Thus, the way forward is to note that the second function, $\frac{1}{\sqrt{2}}\exp(\frac{-\xi^2}{4})$, has inverse transform $\exp(-x^2)$.  The presence of the first exponential indicates that this function should be shifted by the amount $-c t $ in the independent variable $x$.  Thus, the result is 

\begin{equation*}
    u(x,t) = \exp\left[-(x-c t)^2 \right]
\end{equation*}
\end{example}
\end{svgraybox}

\begin{svgraybox}
\begin{example}[Solution to the diffusion equation using the Fourier transform]\indexme{Fourier transform!heat/diffusion problem}
Consider the following diffusion / heat transport problem 
\begin{align*}
    &&\frac{\partial u}{\partial t}&= D \frac{\partial^2 u}{\partial x^2}\qquad -\infty < x < \infty\\
&B.C.1&    u(x,t) &\rightarrow 0 ~~as~ |x|\rightarrow \infty \\
&B.C.2&    \left.\frac{\partial u}{\partial x}\right|_{(x,t)} &\rightarrow 0 ~~as~ |x|\rightarrow \infty\\
&I.C.& u(x,0)=\delta(x)
\end{align*}
This problem represents the diffusion of a ``point" concentration or initial distribution of heat at the center of an infinite domain.  Typically in Fourier transform problems, the two boundary conditions are assumed to be true even when they are not stated; in short, this just indicates that the behavior of most physical systems as one gets arbitrarily distant from the initial condition.  Although, technically speaking, periodic initial conditions might also be possible to \emph{solve}, they often do not make physical sense (e.g., it would require an infinite amount of energy to heat up an infinitely large system!).

Although this problem may look somewhat daunting, the Fourier transform makes it pretty simple.  Transforming both sides of the equation and the initial condition, we have

\begin{align*}
    &&\frac{\partial \hat{u}}{\partial t}&= -\xi^2 D \hat{u}\qquad -\infty < x < \infty\\
&I.C.& \hat{u}(x,0)&=\mathscr{F}[\delta(x)]
\end{align*}
Note that we do not need the two boundary conditions at this juncture, since they were accounted for in the integration by parts that lead to the transformation of the derivative.

The Fourier transform of the initial condition is listed in tables; however, for this particular transform, the definition is also sufficient

\begin{align*}
   \mathscr{F}[\delta(x)]&=\frac{1}{\sqrt{2 \pi}} \int_{x=-\infty}^{x=\infty} \delta(x) e^{-i \xi x} \, dx \\
   &= \frac{1}{\sqrt{2 \pi}}
\end{align*}
Alternatively, from a table of transforms it is possible to find the result
\begin{equation}
    \mathscr{F}[\delta(x-a)] = \frac{1}{\sqrt{2 \pi}} e^{-i a \xi}
\end{equation}
which, evaluated at $a=0$, generates the same result for the transform.  Now, we have a simple first-order ODE of the form
\begin{align*}
    &&\frac{d \hat{u}}{d t}&= -\xi^2 D \hat{u}\qquad -\infty < x < \infty\\
&I.C.& \hat{u}(x,0)&=\frac{1}{\sqrt{2 \pi}}
\end{align*}
This is easily solved using separation of variables to compute the integral

\begin{align*}
    &&\frac{d \hat{u}}{\hat{u}}&= -\xi^2 D dt\qquad -\infty < x < \infty\\
&I.C.& \hat{u}(x,0)&=\frac{1}{\sqrt{2 \pi}}
\end{align*}
This can be recognized to give the conventional exponential solution

\begin{equation}
    \hat{u}(\xi,t) = \frac{1}{\sqrt{2 \pi}}\exp\left( -\xi^2 D t \right)
\end{equation}
This provides a \emph{solution} to the problem, but the solution is still in the transformed space.  In other words, we have the solution for $\hat{u}(\xi,t)$, but we want the solution for ${u}(x,t)$.  To get this solution, we need only use the inverse transform on both sides of the last equation

\begin{equation}
    \mathscr{F}^{-1}\left[\hat{u}(\xi,t)\right] =   \mathscr{F}^{-1}\left[\frac{1}{\sqrt{2 \pi}}\exp\left( -\xi^2 D t \right) \right]
\end{equation}
Obtaining the inverse transform of the left-hand side of this equation is easy; by definition $\mathscr{F}^{-1}\left[\hat{u}(\xi,t)\right]= u(x,t)$.   For the right-hand side, we can at least extract the constant (in other words, you can always extract a constant from an integral) giving the result

\begin{equation}
    u(x,t) =  \frac{1}{\sqrt{2 \pi}} \mathscr{F}^{-1}\left[\exp\left( -\xi^2 D t \right) \right]
    \label{Fsolution}
\end{equation}
At this juncture, we could try inverting the expression using the definition.  However, this will prove to be very difficult in most cases.  Most work with Fourier transforms involves using tables (or software like Mathematica) to compute the inverse Fourier transforms.  To use a table, the process is to look at the transform side of the table (these tables always pair functions and their transforms in two columns) to find a transform that has (except for constants) the same form as our result does.  To that end, most tables will list the following pair
\begin{equation*}
    e^{-a^2 x^2} \Leftrightarrow \frac{1}{a\sqrt{2}}e^{-\xi^2/(4 a^2)}
\end{equation*}
We just need to make our transformed result match this transform (and thereby figure out what the associated constants are).  To do so, we first note that we can operate on transform pairs with constants (or variables that are not part of the transform) as if they were equations.  In other words, if we multiply a function by a constant, its transform will also be multiplied by the same constant.  To see how this works, we can work out the example above in detail. We start with the exponential, and note that the two transformed functions (our function and the one listed in the table) will be identical for the case that 

\begin{align*}
    -\xi^2 D t &= \frac{-\xi^2}{4 a^2}\\
    \intertext{Or, solving this for $a^2$ and $a$}
     a^2&= \frac{1}{4 D t}\\
     a &= \frac{1}{\sqrt{D t/4}}
\end{align*}
Substituting this back into our transform pair, this gives us

\begin{equation*}
    e^{-\tfrac{x^2}{4 D t} } \Leftrightarrow \sqrt{2D t}e^{-\xi^2 Dt}
\end{equation*}
Recall, we are trying to find the inverse transform of just $\mathscr{F}^{-1}\left[\exp\left( -\xi^2 D t \right) \right]$ that appears in Eq.~\eqref{Fsolution}.  We are very close with the transform pair listed above.  Now, we consider multiplying both sides by $1/\sqrt{2D t}$. This gives the result

\begin{equation*}
   \frac{1}{\sqrt{2D t}} e^{-\tfrac{x^2}{4 D t} } \Leftrightarrow e^{-\xi^2 Dt}
\end{equation*}
which is exactly the transform result that we need.  Now, we can invert the transform given in Eq.~\eqref{Fsolution}

\begin{align*}
    u(x,t) &=  \frac{1}{\sqrt{2 \pi}} \mathscr{F}^{-1}\left[\exp\left( -\xi^2 D t \right) \right] \\
    &=  \frac{1}{\sqrt{2 \pi}} \frac{1}{\sqrt{2D t}} e^{-\tfrac{x^2}{4 D t} }\\
    &= \frac{1}{\sqrt{\pi}\sqrt{4  D t}} e^{-\tfrac{x^2}{4 D t} }
    \label{Fsolution1}
\end{align*}
This result was not too difficult to obtain.  The result is an important one; the function that we derived is called the \emph{heat kernel}, and it appears widely in solutions to the heat (or diffusion) equation.
\end{example}
\end{svgraybox}


\section{Convolutions}\indexme{Fourier transform!convolutions}\indexme{convolutions}

In this section, the general concept of convolutions will be introduced.  While the goal for this material is ultimately for applications to the Fourier transform (and, later, for the Laplace transform), here some general properties of convolutions are discussed and illustrated.  Before we explore the use of convolutions in Fourier transforms, this is a good juncture to discuss what convolutions \emph{do}.  One of most physically direct way to see how convolutions ``work" is to explore their use in smoothing operations.  Now, not all convolutions are necessarily smoothing operations; but the example of smoothing using a convolution is a nice one that provides a sense of what they do both functionally and geometrically.  

\subsection{Convolution Definition and Examples}

Convolutions are general mathematical processes that have applications in a wide variety of areas.  For example, they can be used to smooth noisy functions, and their discrete versions are a central part of artificial neural networks called convolutional neural networks (CNNs).  In signal processing, convolutions are used to create filters to either highlight or decrease specific features of the signal.  If you have ever used a \emph{box} filter on a set of data in, say, a spreadsheet program, this is actually a discrete convolution.  The process generally works as follows.  First, a width for the box filter is chosen in terms of numbers of cells that it will span in, say, one row of a spreadsheet.  Suppose we pick that number to be three; this would mean that for any cell $i$, the filter would incorporate the cell above it (cell $i-1$) and the cell below it (cell $i+1$).  Now, suppose one takes the average over these three cells (i.e., the sum of the three cells multiplied by $1/3$).  The value in cell $i$ is now replaced by this average value, and the process is continued by proceeding to the next cell in the column and using the same relative method for determining the average.  In the end, one obtains the  moving box average, or simply the box filtered result for the data in the column.  In the new data set, the values are smoothed by accounting for the value before and after each data point.  This is the essential notion of a convolution, and it can be extended to \emph{continuous} functions.

For continuous functions $f$ and $g$, on the real line, the convolution is defined by either of the following two expressions

\begin{equation}
 ~~~~~~~~~ (f*g)(x) =  \gamma \int_{-\infty}^{\infty} f(x-\xi)g(\xi) \, d\xi \label{convo1}
  \vspace{-4mm}
\end{equation}
\begin{align}
  (g*f)(x) &=  \gamma \int_{-\infty}^{\infty} g(x-\xi)f(\xi) \, d\xi \label{convo2}\\
  \textrm{Note:} (f*g)(x)&=(g*f)(x) \nonumber
\end{align}
The constant $\gamma$ is a rescaling value that is usually set to $\gamma=1$.  Describing what is going on during a convolution is a bit tricky, and here some images will help.  To start, lets define some compact function that is also a density function (so that its area is 1). It is not at all necessary that this function be compact or a density, but it does help with the visualization and interpretation.  Here is the function we will use

\begin{equation}
    w(x;\alpha) = \frac{\pi  a}{4}  \left[H\left(x+\tfrac{1}{\alpha}\right)-H\left(x-\tfrac{1}{\alpha}\right)\right] \cos \left(\frac{\pi
    \alpha x}{2}\right)\label{density0}
\end{equation}
where here $\alpha$ controls the width of the density function (larger $\alpha$ corresponds to narrower $w$.

\begin{figure}[t]
\sidecaption[t]
\centering
\includegraphics[scale=.8]{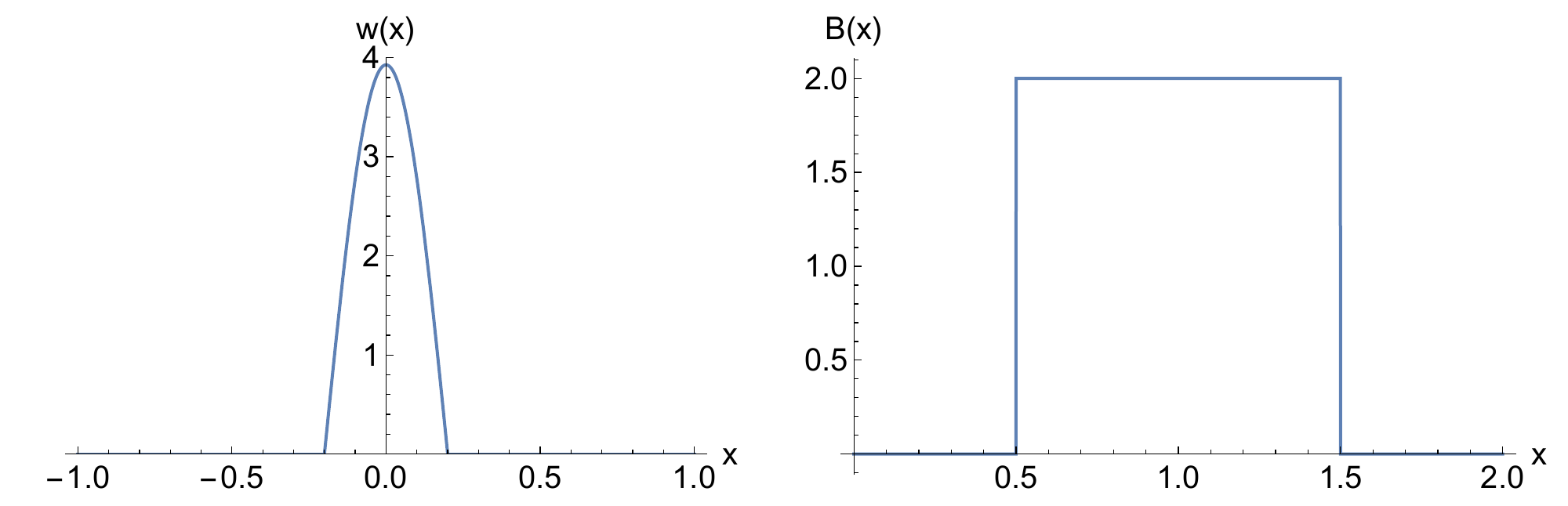}
\caption{(Left) The density function given by Eq.~\eqref{density0}, with $a=5$. (Right) The boxcar function, $2B$ with $a=\tfrac{1}{2}$ and $a+\Delta=\tfrac{3}{2}$.}
\label{densitybox}       
\end{figure}

As a second function, we define the \emph{boxcar} function (presumably because it can take the form of a railroad boxcar under the correct conditions).  This function is the difference of two Heaviside functions, and is defined by

\begin{equation}
   B(x;a;a+\Delta) = \left[ H(x-a)-H(x-(a+\Delta))\right]
\end{equation}
These two functions $w$ and $B$ are plotted in Fig.~\ref{densitybox}.  Note that the density function $w$ is both symmetric around $x=0$.  This is not necessary for convolutions in general; it is, however, necessary for smoothing operations.  The symmetry of the density function around $x=0$ creates a result where the boxcar function is smoothed, but its center of mass does not change.   The fact that $w$ is a \emph{density} function with unit area assures that the total area of the smoothed boxcar function is the same as the original.  When conducting smoothing operations using convolutions, these are often desirable properties.

Now, to compute the convolution, we need only compute either of Eqs.~\eqref{convo1} or \eqref{convo1}.  The results are identical, and it makes no difference in the solution which function is selected to represent $f$ or $g$ in evaluating these integrals.  Note however, the choice of which function contains the displacement $(x-\xi)$ can influence the difficulty of the resulting integral.  Suppose we set up our convolution as follows.  Note that for $B$ we have $a=\tfrac{1}{2}$ and $a+\Delta=\tfrac{3}{2}$, and for $w$ we have taken $\alpha=5$

\begin{align}
  B'(x; \tfrac{1}{2},\tfrac{3}{2})& =  \int_{-\infty}^{\infty} w(x-\xi)B(\xi) \, d\xi \nonumber \\
   &= \int_{-\infty}^{\infty}\left[ \frac{5 \pi }{4}  \left[H\left((x-\xi)+\tfrac{1}{5}\right)-H\left((x-\xi)-\tfrac{1}{5}\right)\right] \cos \left(\frac{5\pi
     (x-\xi)}{2}\right) \right] \nonumber\\
    &\times \left[ H(\xi-\tfrac{1}{2})-H(\xi-\tfrac{3}{2})\right]\,d \xi
    \label{theconvolution}
\end{align}
While this expression is not easy to evaluate, it is possible to do so.  Later on, we will see how Fourier transforms can be used to \emph{evaluate} this integral.  The result is the following somewhat lengthy result

\begin{figure}[t]
\sidecaption[t]
\centering
\includegraphics[scale=.5]{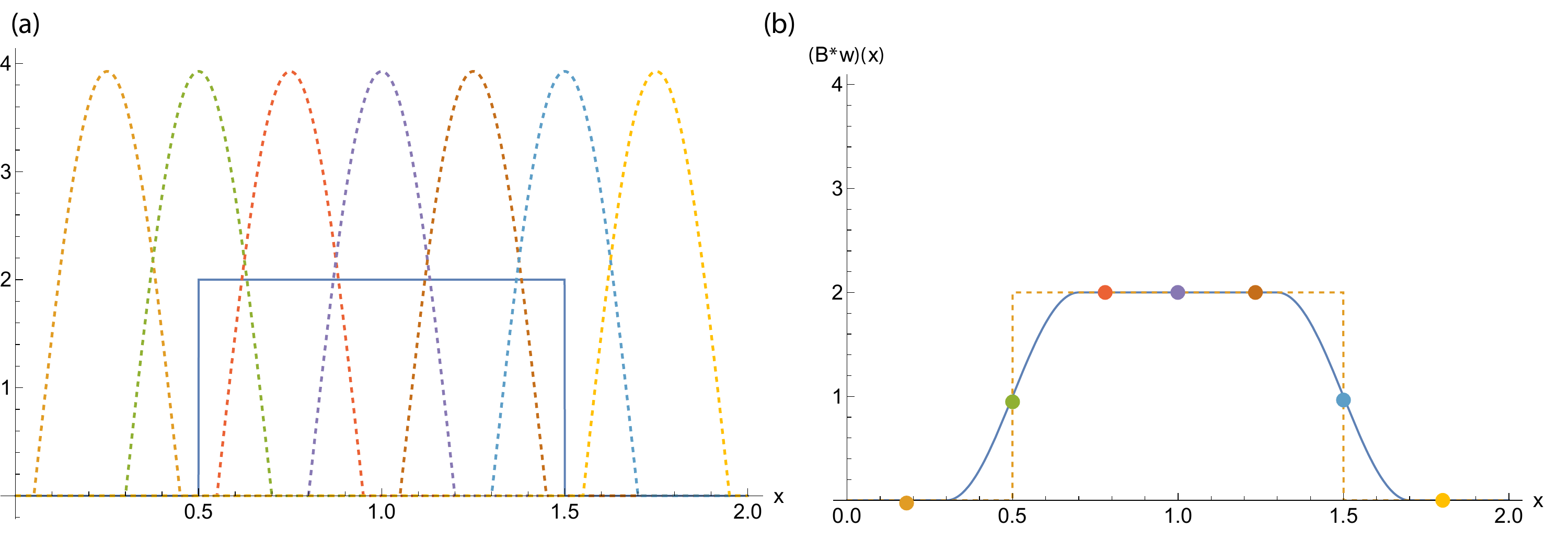}
\caption{(a) The idea in the convolution given by Eq.~\eqref{theconvolution} is that shifted copies of the density function $w$ are integrated against the boxcar function.  We can thing of each such integration as occurring for a fixed value of $x$.  Then, each integration over the $\xi$ variable yields a single point on the new curve shown in (b).  (b) The function resulting from the convolution.  Each colored circle on the curve corresponds to the integration of a shifted copy of $w$ times $B$; colors correspond to curves shown in (a).  Both the center of mass and the area of $B$ are maintained at their original values when $w$ is a symmetric density function centered at $x=0$.}
\label{densitybox2}       
\end{figure}

\begin{align}
B'(x; \tfrac{1}{2},\tfrac{3}{2}) &=
    \frac{1}{2 \sqrt{2}}\textrm{sgn}\left(\frac{17}{10}-x\right) \left[-\sin \left(\frac{5 \pi  x}{2}\right)-\cos \left(\frac{5 \pi 
   x}{2}\right)+\sqrt{2}\right]\nonumber \\
   -&\textrm{sgn}\left(\frac{3}{10}-x\right) \left(-\sin \left(\frac{5 \pi  x}{2}\right)+\cos
   \left(\frac{5 \pi  x}{2}\right)+\sqrt{2}\right)\nonumber \\
   -&\textrm{sgn}\left(\frac{7}{10}-x\right) \Bigg[\sin \left(\frac{5 \pi 
   x}{2}\right) -\cos \left(\frac{5 \pi  x}{2}\right)
   +\sqrt{2}\Bigg]\nonumber \\
   +&\textrm{sgn}\left(\frac{13}{10}-x\right) \left[\sin
   \left(\frac{5 \pi  x}{2}\right)+\cos \left(\frac{5 \pi  x}{2}\right)+\sqrt{2}\right]
\label{bigconvosol}
\end{align}
This result is plotted in Fig.~\ref{densitybox2}.  
\newpage
Here, the sign function $sgn(x)$ is defined by 
\begin{align}
    \textrm{sgn}(x)&=\begin{cases}
        -1 & x<0 \\
        0 & x=0 \\
        1 & x>0
    \end{cases}
\end{align}

\begin{figure}[t]
\includegraphics[scale=.5]{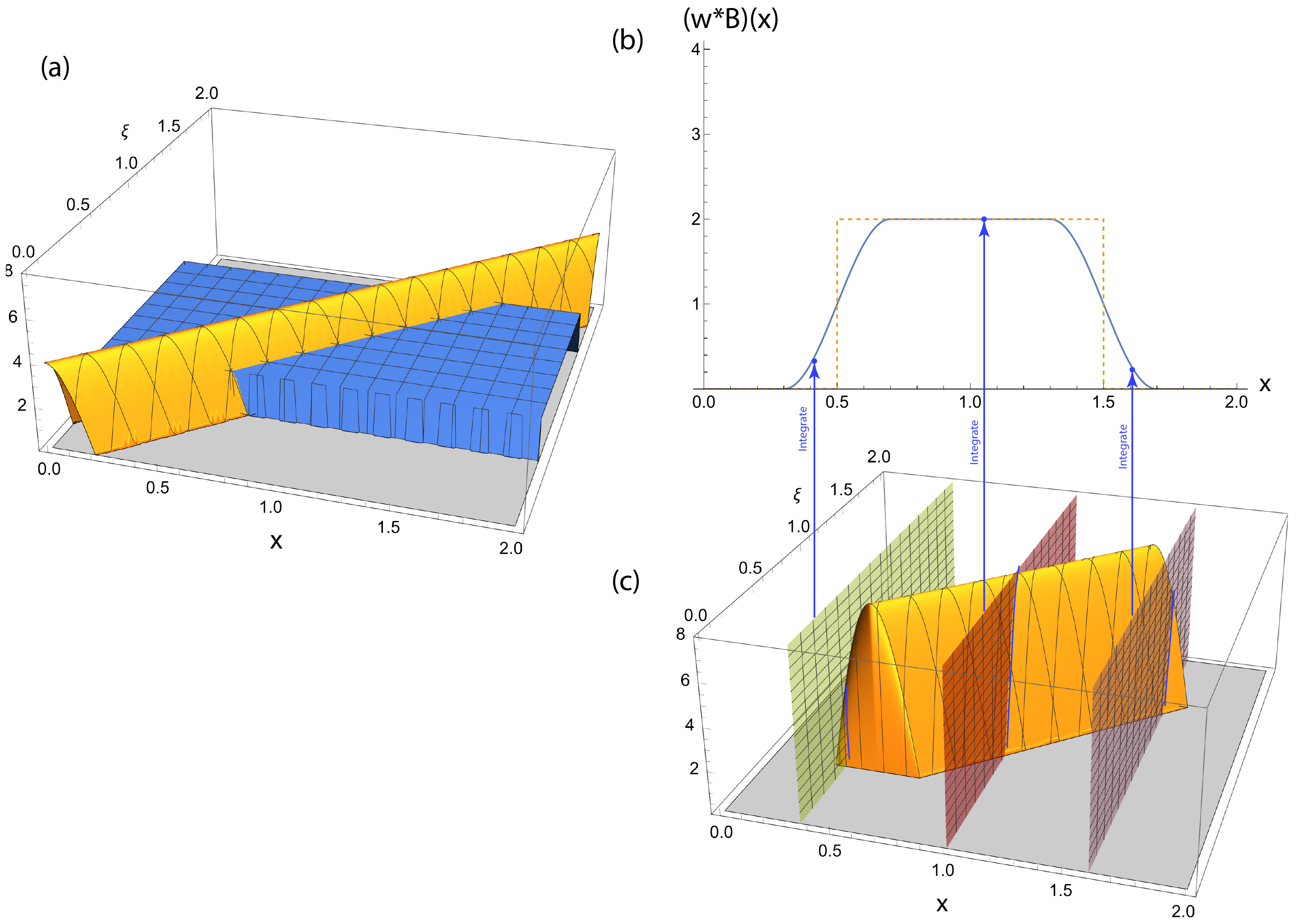}
\caption{(Left) A graphical interpretation of the convolution $w*B=\int_{-\infty}^{\infty} w(x-\xi)B(\xi)\, d\xi$. (a) The two functions $B(\xi)$ and $w(x-\xi)$ shown as a function of the two independent variables $x$ and $\xi$.  (b) The result of the convolution is a smoothed version of the original function $B$.  Note that here, the dashed line is $B(x)$, whereas the solid blue line is the convolution.  (c) The convolution, $(w*B)(x)$, is found by integrating along the $\xi$-direction for fixed values of $x$; each point on the convolution is given by one such integration.  The function shown in orange is $w(x-\xi)\times B(\xi)$.  The blue line indicates the intersection of the function $w(x-\xi)\times B(\xi)$ with vertical planes (constant $x$) perpendicular to the $x$-$\xi$ coordinate plane.}
\label{integration2}    
\end{figure}
%
In Fig.~\ref{integration2}, there is a second presentation of the convolution process.  Here, the idea is to think of the convolution as happening in a 2-dimensional space defined by the coordinates $x$ and $\xi$.  Recall that the integrand in the convolution is given by $w(x-\xi)\times B(\xi)$.  Thus, when plotted in $(x,\xi)$ space, the function $B(\xi)$ is independent of $x$, as shown in Fig.~\ref{integration2}(a) (blue curve).  However, the weighting function \emph{is} a function of both $x$ and $\xi$.  Recall that the weighting function was symmetric around $x=0$.  Thus, $w(x)=w(-x)$. If we like, we could rewrite the function $w$ as $w(-(\xi-x))$ which is identical to $w(\xi-x)$.  So, while the function $B$ stays stationary in $(x,\xi)$ space, the function $w$ we can consider as one that is continually being translated by the amount $x$ in $(x,\xi)$ space.  This results in the function illustrated in Fig.~\ref{integration2}(a) (orange curve).  

Finally, the convolution itself involves integrating the multiplication of these two curves as $\xi$ goes from $-\infty$ to $\infty$.  Because both $w$ and $B$ are compact, the actual range of integration is the interval $\xi \in[\tfrac{3}{10},\tfrac{17}{10}]$ (note that these two fractions appear in the solution given by Eq.~\ref{bigconvosol}.  This region is plotted in Fig.~\ref{integration2}(c).  In that figure, three representative planes with fixed values of $x$ are illustrated.  The curves formed by the intersection of any such plane and the function $w(x-\xi)\times B(\xi)$ (the orange surface in Fig.~\ref{integration2}(c)) are integrated during the convolution.  The area of each such curve gives a single point on the solution curve which is shown in Fig.~\ref{integration2}(b).

\section{Convolutions in Fourier Transforms}\indexme{Fourier transform!convolutions}\indexme{convolutions}

There are times when a Fourier transform yields a function that does not have a known inverse function, but it can be expressed as the \emph{product} of two functions whose inverse is known.  While this might seem on the surface to be an unusual event, it turns out that it happens quite frequently because of the structure of most initial-boundary value problems. 

The following is a fairly remarkable feature of Fourier transforms, further illustrating the utility of the method.  When one would like to inverse transform the product of two Fourier transformed functions, the result can be expressed

\begin{align}
    \mathscr{F}^{-1}[\hat{g}(\xi)\cdot \hat{f}(\xi)] &= \frac{1}{\sqrt{2 \pi}} \int_{\xi =-\infty}^{\xi=\infty} f(x-\xi)g(\xi)\, d\xi\\
    &= \frac{1}{\sqrt{2 \pi}} \int_{\xi =-\infty}^{\xi=\infty} f(\xi)g(x-\xi)\, d\xi
\end{align}
Thus, the result states that the inverse transform of the product of two Fourier transformed functions $\hat{f}$ and $\hat{g}$ is just the convolution of the functions $f$ and $g$ in their \emph{non-transformed} state.  Note that for this definition, the scaling factor $\gamma$ is not zero.  This scaling factor shows up in the convolution directly from the definition of the Fourier transform.

We can also define the inverse of the convolution above, where one is taking the Fourier transform of the convolution of functions $f$ and $g$.

\begin{align}
    \mathscr{F}[{g}(x)\cdot {f}(x)] &= \frac{1}{\sqrt{2 \pi}} \int_{x =-\infty}^{x=\infty} \hat{f}(x-\xi)\hat{g}(\xi)\, dx\\
    &= \frac{1}{\sqrt{2 \pi}} \int_{x =-\infty}^{x=\infty} \hat{f}(\xi)\hat{g}(x-\xi)\, dx
\end{align}
While it might not be obvious at this point how this relationship might be useful, it can be easily seen through an explicit example.  In this example, we consider the problem that we examined earlier, but with a more complicated (or, at least, in principle so) initial condition.

As an example of the power of the Fourier transform of convolutions, the following example is provided.
\begin{svgraybox}
\begin{example}[Computing a convolution using the Fourier transform]
One way that the convolution identities can be used is to actually compute the value of a convolution using the Fourier transform.  Suppose we take the following two functions

\begin{align}
    f(x) &= 2 \sqrt{2} \exp(-x^2) \\
    g(x) &= H(x+a)-H(x-a)
\end{align}
The convolution of these two quantities is (using the weighting factor $\gamma = 1/\sqrt{2 \pi}$)

\begin{equation}
    (f*g)(x) = \dfrac{1}{\sqrt{2 \pi}} \int_{-\infty}^\infty 2 \sqrt{2} \exp[-(x-\xi)^2][H(x+a)-H(x-a)]\, d\xi
\end{equation}
This integral would present a serious challenge to compute directly.  However, using the properties of the Fourier transform of a convolution, we know that the transform will be just the product of $\hat{f}$ and $\hat{g}$.  Note that the Fourier transform tables (see the appendix for this chapter) provide the following
\begin{align}
    \hat{f} &= 2 \exp\left(-\frac{\xi^2}{4}\right) && (\textrm{entry \ref{expx2} in the Fourier transform table)}  \\
    \hat{g} & =\dfrac{1}{\xi} \sqrt{\dfrac{2}{\pi}} \sin(a \xi)&& (\textrm{entry \ref{boxcartable} in the Fourier transform table)}
\end{align}
Thus, we have
\begin{equation}
\mathscr{F}\left[ (f*g)(x)\right] = 2 \exp\left(-\frac{\xi^2}{4}\right)\dfrac{1}{\xi} \sqrt{\dfrac{2}{\pi}} \sin(a \xi)\label{fgtrans}
\end{equation}
Now, examining the Fourier transform tables again, we find tht Eq.~\ref{fgtrans} is on the table as entry number \ref{erferf}. Thus, taking the inverse transform of both sides of Eq.~\ref{fgtrans} gives the result

\begin{equation}
    (f*g)(x) = \textrm{erf}(a+x)-\textrm{erf}(a-x)
\end{equation}
This result was obtained without too much effort, but with careful use of the properties of the Fourier transform of convolutions.
\end{example}
\end{svgraybox}

In general, determining the value of convolutions is not the primary use of the Fourier transform.  On the real line, the Fourier transform can provide a very effective method for solving linear partial differential equations.  This requires, of course, that all of the terms in the equation have well-known Fourier transforms.  The process is one in which a PDE is changed, via the Fourier transform, into an ODE which is simpler to solve.  Once the problem is solved, the inverse transform can be used to obtain the form of the function back in the original (real space) variables.  There are some subtleties involved in using the Fourier transform and the tables of transforms; it takes some practice to become familiar with the process.  While not always the best solution, the Fourier transform is nonetheless an additional tool in the arsenal of resources available to solve PDEs.  The following example illustrates the process on the familiar heat/diffusion equation on the real line.

\begin{svgraybox}
\begin{example}[A solution to the heat/diffusion equation using the Fourier transform- General initial conditions]

Consider the following diffusion / heat transport problem 
\begin{align*}
    &&\frac{\partial u}{\partial t}&= D \frac{\partial^2 u}{\partial x^2}\qquad -\infty < x < \infty\\
&B.C.1&    u(x,t) &\rightarrow 0 ~~as~ |x|\rightarrow \infty \\
&B.C.2&    \left.\frac{\partial u}{\partial x}\right|_{(x,t)} &\rightarrow 0 ~~as~ |x|\rightarrow \infty\\
&I.C.& u(x,0)=\varphi(x)
\end{align*}
This problem represents the diffusion of a general initial condition in an infinite domain.  It is an extension of the simpler initial condition (delta condition) that we examined earlier.

Transforming both sides of the equation and the initial condition, we have

\begin{align*}
    &&\frac{\partial \hat{u}}{\partial t}&= -\xi^2 D \hat{u}\qquad -\infty < x < \infty\\
&I.C.& \hat{u}(x,0)&=\mathscr{F}[\varphi(x)]
\end{align*}
Upon transforming the PDE, we have a simple first-order ODE of the form
\begin{align*}
    &&\frac{d \hat{u}}{d t}&= -\xi^2 D \hat{u}\qquad -\infty < x < \infty\\
&I.C.& \hat{u}(x,0)&=\hat{\varphi}(x)
\end{align*}
This is easily solved using separation of variables to compute the integral

\begin{align*}
    &&\frac{d \hat{u}}{\hat{u}}&= -\xi^2 D dt\qquad -\infty < x < \infty\\
&I.C.& \hat{u}(x,0)&=\hat{\varphi}(x)
\end{align*}
This can be recognized to give the conventional exponential solution

\begin{equation}
    \hat{u}(\xi,t) = \hat{\varphi}(x)\exp\left( -\xi^2 D t \right)
\end{equation}
This provides a \emph{solution} to the problem, but the solution is still in the transformed space.  In other words, we have the solution for $\hat{u}(\xi,t)$, but we want the solution for ${u}(x,t)$.  To get this solution, we need only use the inverse transform on both sides of the last equation

\begin{equation}
    \mathscr{F}^{-1}\left[\hat{u}(\xi,t)\right] =   \mathscr{F}^{-1}\left[\hat{\varphi}(x)\exp\left( -\xi^2 D t \right) \right]
\end{equation}
Obtaining the inverse transform of the left-hand side of this equation is easy; by definition $\mathscr{F}^{-1}\left[\hat{u}(\xi,t)\right]= u(x,t)$.   For the right-hand side, we are a bit stuck, since we have not yet specified what $\hat{\varphi}$ is.  However, the convolution theorem can come to the rescue here.

\begin{equation}
    u(x,t) =  \mathscr{F}^{-1}\left[\hat{\varphi}(x)\exp\left( -\xi^2 D t \right) \right]
\end{equation}
which, using the convolution theorem, is

\begin{align}
    u(x,t)& = \frac{1}{\sqrt{2 \pi}} \int_{x=-\infty}^{x=\infty}\varphi(x-\xi)\frac{1}{\sqrt{2D t}} e^{-\tfrac{x^2}{4 D t} }\\
    \intertext{or, equivalently}
     u(x,t)& = \frac{1}{\sqrt{2 \pi}} \int_{x=-\infty}^{x=\infty}\varphi(\xi)\frac{1}{\sqrt{2D t}} e^{-\tfrac{(x-\xi)^2}{4 D t} }
\end{align}

\end{example}
\end{svgraybox}

This last result is particularly relevant.  The term $1/\sqrt{2 D t} \exp(-(x-\xi)^2/(4 D t)$ is a somewhat famous expression known generally as the \emph{heat kernel} \indexme{kernel!heat}.  The equation is a general one that defines the solution for \emph{any} initial condition $\varphi$ by forming the appropriate convolution with the heat kernel.  If the resulting integral can be computed, then one has an explicit solution to this otherwise very difficult problem.  Even if the solution cannot be computed analytically, the result is at least in the form of an integral.  Such integrals can often be computed very accurately using specific numerical methods designed to determine the value of finite integrals on infinite domains.

\section{Extensions of the Fourier Transform}\indexme{Fourier transforms!Extensions to}

There are a number of functions that occur regularly in applications, but fail to have a Fourier transform, at least by the conventions that we have identified.  Recall, the primary requirement that we stated for the Fourier transform to exist is that the function being transformed must have a unique finite square integral (i.e., it is square-integrable).  However, a number of common and useful functions are not square-integrable.  For example, both $\sin(x)$ and $\cos(x)$ are not square integrable their square integrals tend toward infinity.   Similarly, the delta function is not square integrable because there the delta is technically a \emph{distribution} not a function; it turns out that there is no general notion of the square of the delta function.  This puts us in a somewhat difficult position, because these functions are all fairly fundamental, and are the kinds of things that appear routinely in applications.

Fortunately, there is a more general framework that lends a sound interpretation of the Fourier transform of the such functions.  This theory involves \emph{distribution theory} as was briefly introduced in Chp.~\ref{deltachap}, and, combined with the definition of the Fourier transform leads to a more general framework known as the Fourier transform of \emph{tempered distributions}.  Very roughly, tempered distributions are simply functions that ``do not grow too fast".  Here not growing too fast means that the functions grow at some finite polynomial rate (and not at some exponential rate).    We will not be investigating tempered distributions per se, but it is important for us to know that the developments above do have a secure basis for understanding.  In the developments that follow, we will develop Fourier transforms for a few functions that are not transformable according to our initial requirements of being square-integrable.  However, each of them belong to the class of functions that are the tempered distributions.  Because this more sophisticated theory allows us to determine the transform of such functions, then if we can find the transform (and its inverse) then we are in a sense ``allowed" to expand our library of transformable functions.  In the material below, we focus only on the transforms of (1) the Heaviside and delta functions, (2) the constant function, and (3) the sine and cosine functions.  There are a few additional functions (such as the natural logarithm) where similar concepts apply, but the results are not given in detail here (however, their transforms are given in the table of Fourier transforms).

\subsection{Fourier Transform of the Delta Function}

As discussed in Chapter \ref{deltachap}, the delta is not a true \emph{function} (even though we use the name \emph{delta function}); rather, it is a mathematical construct that represents a impulse (with infinite magnitude) applied to a single point.  The delta function can be given clear mathematical meaning when it is integrated against other (appropriately constrained) functions.  

It makes some intuitive sense, then, that the delta function can be represented by a Fourier transform; a Fourier transform simply integrates the input function, $f$, against the integral kernel $e^{-i\xi x}$ to determine the transformed function.  Using the conventional properties of the delta function and the shifting property of Fourier transforms, we can easily derive the following

\begin{align}
    \mathscr{F}[\delta(x)] &=\frac{1}{\sqrt{2 \pi}} \int_{-\infty}^{\infty} \delta(x) e^{-i \xi x}\, dx \nonumber\\
    &= \frac{1}{\sqrt{2 \pi}}
\end{align}

\begin{align}
    \mathscr{F}[\delta(x-a)] &=\frac{1}{\sqrt{2 \pi}} \int_{-\infty}^{\infty} \delta(x-a) e^{-i \xi x}\, dx \nonumber\\
    &= \frac{1}{\sqrt{2 \pi}}e^{-i a \xi} 
\end{align}

Note that this last identity implies the following

\begin{align}
    \mathscr{F}[\delta(x+a)] &= \frac{1}{\sqrt{2 \pi}}e^{i a \xi} 
\end{align}

and using the Fourier inversion formula gives

\begin{align}
    \mathscr{F}[e^{-i a \xi}] = \sqrt{2 \pi} \delta(x)
\end{align}

\subsection{Fourier Transform of the Heaviside Function}

To develop the Fourier transform of the Heaviside function, we take a bit of an indirect path.  Recall, the Heaviside function (or the \emph{unit step function}) is given by 
\begin{equation}
    H(x) =
    \begin{cases}
       0 & x < 0 \\
        1 & x \ge 0
    \end{cases}
\end{equation}
From our discussion of admissible functions for Fourier transform in \S \ref{fourtrans}, this function does not qualify because it is not square integrable.  Nonetheless, a sensible notion of the Fourier transform of the constant function can be developed. 
To start, consider the following transformed function

\begin{equation}
    \hat{f}(\xi) = \frac{1}{\sqrt{2 \pi}} \frac{1}{i \xi}
\end{equation}
While this proposed transform has simply proposed without explanation, we can proceed to see where the proposition leads us.  We can consider the inverse Fourier transform for $\hat{f}$ by computing the following integral.

\begin{align}
    \mathscr{F}^{-1}\left[ \frac{1}{\sqrt{2 \pi}} \frac{1}{i \xi} \right] &= \frac{1}{\sqrt{2 \pi}} \int_{-\infty}^\infty  \frac{1}{\sqrt{2 \pi}} \frac{1}{i \xi} e^{i \xi x} d\xi \nonumber\\
    &= \frac{1}{2 \pi} \int_{-\infty}^\infty  \frac{e^{i \xi x}}{i \xi}  d\xi\nonumber \\
    &= \frac{1}{2 \pi i} \int_{-\infty}^\infty  \frac{\cos(\xi x)+i \sin(\xi x)}{\xi}  d\xi \nonumber \\
    \intertext{Now, noting that $\sin(\xi)/\xi$ is an \emph{even} function, and $\cos(\xi)/\xi$ is an \emph{odd} function, it is clear that the cosine term will evaluate to zero.  Further, the sine term is symmetric about $x=0$, leading to}
    &= \frac{1}{2 \pi} \int_{-\infty}^\infty  \frac{\sin(\xi x)}{\xi}  d\xi \nonumber\\
    &= \frac{1}{\pi} \int_{0}^\infty  \frac{\sin(\xi x)}{\xi} d\xi  \nonumber
\end{align}
As a final step, we will express this integral as follows

\begin{equation}
   \frac{1}{\pi}  {\mathop {\lim }\limits_{L \to \infty}}\left[ \int_{0}^L  \frac{\sin(\xi x)}{\xi} d\xi \right]
\end{equation}

This integral in the square brackets is called the sine integral.  Many symbolic mathematical packages can evaluate this integral; in particular, as $L\rightarrow \infty$, we have 

For $x>0$
\begin{equation}
  \frac{1}{\pi}   {\mathop {\lim }\limits_{L \to \infty}}\left[  \int_{0}^L  \frac{\sin(\xi x)}{\xi} d\xi \right] = \frac{1}{2}
\end{equation}

and for $x<0$

\begin{equation}
  \frac{1}{\pi}   {\mathop {\lim }\limits_{L \to \infty}}\left[  \int_{0}^L  \frac{\sin(\xi x)}{\xi} d\xi \right] = -\frac{1}{2}
\end{equation}

At this point it is helpful to develop a new type of step-like function which we might call the zero-average Heaviside function.  Suppose that we define the zero mean Heaviside, $H_0$ by 

\begin{equation}
    H_0(x) =
    \begin{cases}
        -\frac{1}{2} & x<0 \\
        ~~~\frac{1}{2} & x>0
    \end{cases}\label{h0transform}
\end{equation}

With this definition and the result above, we have the following 
\begin{equation}
    \mathscr{F}[H_0(x)]= \frac{1}{\sqrt{2 \pi}} \frac{1}{i \xi}
\end{equation}
With a little thought, we can relate the standard Heaviside function to the zero-mean Heaviside function simply by shifting $H(x)$ vertically by $1/2$.

\begin{equation}
H_0(x) = H(x) - \frac{1}{2}
\end{equation}
Thus, we can write Eq.~\eqref{h0transform} as follows

\begin{equation}
   \hat{f}=\mathscr{F}\left[ H(x)-\frac{1}{2} \right]  =  \frac{1}{\sqrt{2 \pi}} \frac{1}{i \xi}
\end{equation}
So, we have the Fourier transform for the function $H(x)-\frac{1}{2}$, but not the one for $H(x)$.  Note that this is in line with our initial requirement that the functions amenable to Fourier transform must be integrable; $H(x)-1/2$ has, by symmetry properties, a zero integral and thus qualifies (N.B., one must be careful with such arguments on infinite domains-- sometimes the infinite does things that are quite counterintuitive). 

We can proceed forward by exploiting the linearity of the Fourier transform.   

\begin{align}
    \mathscr{F}\left[ H(x)-\frac{1}{2} \right]  &=  \frac{1}{\sqrt{2 \pi}} \frac{1}{i \xi}
    \end{align}

 \begin{align}   
     \mathscr{F}\left[ H(x)\right]- \mathscr{F}\left[\frac{1}{2} \right]&=  \frac{1}{\sqrt{2 \pi}} \frac{1}{i \xi} 
\end{align}
From the material regarding the Fourier transform of the delta function, we have
 \begin{align}
   \mathscr{F}\left[ H(x)\right]- \sqrt{\pi} \delta(\xi)&=  \frac{1}{\sqrt{2 \pi}} \frac{1}{i \xi}
\end{align}
Adding $\sqrt{\pi} \delta(\xi)$ to both sides yields the final result
\begin{align}
   \mathscr{F}\left[ H(x)\right]&=  \sqrt{\frac{\pi}{2}} \delta(\xi)+ \frac{1}{\sqrt{2 \pi}} \frac{1}{i \xi} 
\end{align}

The developments above could be repeated to define the \emph{negative Heaviside function}. 

\begin{equation}
    H_-(x) =
    \begin{cases}
        1 & x < 0 \\
        0 & x \ge 0
    \end{cases}
\end{equation}
So that 
\begin{align}
   \mathscr{F}\left[ H(x)\right]&=  \sqrt{\frac{\pi}{2}} \delta(\xi)- \frac{1}{\sqrt{2 \pi}} \frac{1}{i \xi} 
\end{align}

\subsection{The Fourier Transform of the Constant Function}

Once we have defined the Fourier transform of the delta function and the unit step function, some somewhat odd things happen.  For example, we have found the Fourier transform for the unit step function, which is a decidedly non-integrable function on the real line.  The \emph{reason} that these odd things happen is because we have allowed the delta function into our world of acceptable functions.  While we can make sense of the Fourier transform of the delta function, we must remember that the delta function is actually a \emph{distribution}, and not a function at all!  While we have managed to illustrate the properties of the Fourier transform of the delta and Heaviside functions, it has used some arguments that are not necessarily obvious ones.  

We have another ``trick up our sleeves", so to speak, to develop the Fourier transform of a tempered distribution but using our existing knowledge of the Fourier transform.  The constant function is a tempered distribution (it is the polynomial defined by  $c x^0$, and grows at rate zero).  To achieve this, we can use the results for the Heaviside and negative Heaviside functions above. Note the following:

\begin{equation}
    c(H(x)+H_-(x)) = c, ~~x\in(-\infty,\infty)
\end{equation}
However, we already have the Fourier transform for these two functions from the discussion in the previous section.  Thus

\begin{equation}
    \mathscr{F}[c(H(x)+H_-(x))] = \mathscr{F}[c]
\end{equation}
and from the Fourier transforms for $H$ and $H_-$ we have 

\begin{equation}
c\left(\sqrt{\frac{\pi}{2}} \delta(\xi)+ \frac{1}{\sqrt{2 \pi}} \frac{1}{i \xi}\right)  +  c\left( \sqrt{\frac{\pi}{2}} \delta(\xi)- \frac{1}{\sqrt{2 \pi}} \frac{1}{i \xi} \right) = \mathscr{F}[c]
\end{equation}
Upon simplifying and rearranging, we have

\begin{equation}
\mathscr{F}[c]=c\sqrt{\frac{\pi}{2}} \delta(\xi)   +  c \sqrt{\frac{\pi}{2}} \delta(\xi) 
\end{equation}
Recall that $\mathscr{F}[c]$ is just the integral 

\begin{equation}
    \mathscr{F}[c] = \frac{1}{\sqrt{2 \pi}} \int_{-\infty}^\infty c \exp(-i x \xi) \, dx
\end{equation}
Substituting this on the left-hand side gives

\begin{equation}
\frac{1}{\sqrt{2 \pi}} \int_{-\infty}^\infty c \exp(-i x \xi) \, dx=c\sqrt{\frac{\pi}{2}} \delta(\xi)   +  c \sqrt{\frac{\pi}{2}} \delta(\xi) 
\end{equation}
And simplifying this result gives us 

\begin{equation}
\frac{1}{\sqrt{2 \pi}} \int_{-\infty}^\infty 1\times \exp(-i x \xi) \, dx=\sqrt{2 \pi} \delta(\xi) 
\end{equation}
Rewriting the left-hand side as $\mathscr{F}[1]$ gives us a final result of
\begin{equation}
\mathscr{F}[1]=\sqrt{2 \pi} \delta(\xi) 
\end{equation}
or, equivalently
\begin{equation}
\mathscr{F}\left[ \frac{1}{\sqrt{2 \pi}}\right]= \delta(\xi) 
\end{equation}

\subsection{The Fourier Transform of $\sin(a x)$ and $\cos(a x)$}

The two functions $\sin( a x)$ and $\cos(a x)$ are interesting functions to consider transforming.  Again, they do not meet the requirements we set out for functions that have a Fourier transform, because they are not square integrable.  However, clearly we have that $sin(x)^2  \le 1$ and $cos(x)^2  \le 1$ for all values of $x$, and we now know that the function $f(x)=1$ has a Fourier transform.  Thus, it is reasonable to wonder if both the sine and cosine functions can transformed in some way that makes sense to us without resorting to the theory of tempered distributions.  We can use the following construction based on our existing notions of the Fourier transform, but again relaxing the assumption that we can relax the condition that functions must be square-integrable.  

First, we note the following identities

\begin{align}
    \cos(a x) &=\dfrac{1}{2}(\exp[i a x]+\exp[-i a x]) \\
    \sin(a x) &=-i \dfrac{1}{2}(\exp[i a x]-\exp[-i a x]) 
\end{align}
Examining the table of Fourier transforms, note that transform number \ref{shiftdelt} gives the needed transform, i.e., 

\begin{equation}
    \mathscr{F}\left[\sqrt{\dfrac{1}{2\pi}}\exp(-i a x)\right] =  \delta(\xi-a)
\end{equation}
In a fairly straightforward manner, we can use this information to prove 

\begin{align}
    \mathscr{F}[\cos(a x)] &=  \sqrt{\dfrac{\pi}{2}}\left[\delta(\xi +a)+\delta(\xi-a)\right]\\
    \mathscr{F}[\sin(a x)] &= i \sqrt{\dfrac{\pi}{2}}\left[\delta(\xi +a)-\delta(\xi-a)\right]
\end{align}
Note that for $a=0$, this reduces to the appropriate Fourier transform of the constant function $f(x)=1$ and of the zero function $f(x)=0$, respectively.

\begin{svgraybox}
\begin{example}[A solution to the Poisson equation using the Fourier transform]
In the chapter on separation of variables, we explored some solutions to the (homogeneous) Laplace equation.  The Laplace equation still involved only two independent variables, but for that case the two variables are both spatial.  We were able to show that separation of variables for this case operated essentially the same way as it did for other problems in two independent variables.  

We can use the Fourier transform to solve the Laplace equation on the infinite plane in two spatial variables, but we will see that this requires the use of the transform of each variable independently.  In fact, we are able to illustrate another strength of the Fourier transform in this example, and that is the ease in which one can develop solutions for nonhomogeneous problems when using transform methods.

To start, recall that when the Laplace equation has a source term it is usually called the Poisson equation.  Consider the following Poisson equation on the infinite 2-dimensional plane.

\begin{align*}
    &&\frac{\partial^2 u}{\partial x^2}+ \frac{\partial^2 u}{\partial y^2}&=-\delta(x)\delta(y)\qquad -\infty < x < \infty, ~-\infty < x < \infty\\
&B.C.1&    u(x,y) &\rightarrow 0 ~~as~ x^2+y^2\rightarrow \infty \\
&B.C.2a&    \left.\frac{\partial u}{\partial x}\right|_{(x,y)} &\rightarrow 0 ~~as~ x^2+y^2\rightarrow \infty\\
&B.C.2b&    \left.\frac{\partial u}{\partial y}\right|_{(x,y)} &\rightarrow 0 ~~as~ x^2+y^2\rightarrow \infty\\
\end{align*}
Here, we have used $d^2=x^2+y^2$ as the square of the distance between any point and the origin as the appropriate measure on the plane.  

Now, consider transforming both sides of the equation with respect to the variable $x$.  Recall, the Fourier transform is given by

\begin{equation*}
    \mathscr{F}[f(x)] = \hat{f}(\xi) = \int_{-\infty}^\infty f(x) \, e^{-i x \xi} \, dx
\end{equation*}
Transforming the $x-$variable, we obtain the result
\begin{align*}
    -\xi^2{\hat{u}}(\xi,y)+ \dfrac{\partial^2 \hat{u}}{\partial y^2}=\sqrt{\dfrac{1}{2\pi}}\delta(y)
\end{align*}
Take careful note that here we have transformed only the variable $x$ by applying the Fourier transform in the variables $x$ and $\xi$, as has been the convention thus far.  At this juncture, we have eliminated one of the two derivative terms via the transform.  However, we are left with a second-order nonhomogeneous equations.  You may recall from our review of ODEs that such equations are solvable using the methods of variation of parameters, but some effort is required in obtaining solutions.  

As an alternative, consider conducting a second Fourier transform, but this time in variable $y$.  Note that we have $-\infty < y < \infty$, and these are exactly the conditions where the Fourier transform is useful.  There is no reason that we cannot use two independent transforms to transform the quantities in both $x$ and in $y$, but we will have to take extra effort to keep track of the variables involved.  Suppose we use the pair $(x,\xi)$ to indicate the Fourier transform pair of independent variables on $x$, and $(y,\eta)$ for the independent variable pair on $y$.   The transform of the $y$ coordinate, then, is given by 

\begin{equation*}
    \mathscr{F}_y[f(y)] = \bar{f}(\eta) = \int_{-\infty}^\infty f(y) \, e^{-i y \eta} \, dy
\end{equation*}
Now, with this definition in place, we can conduct a second transform of our Poisson equation.  This gives

\begin{align*}
    \mathscr{F}_y\left[ -\xi^2{\hat{u}}(\xi,y)\right]+ \mathscr{F}_y\left[\dfrac{\partial^2 \hat{u}}{\partial y^2} \right]=\mathscr{F}_y\left[-\sqrt{\dfrac{1}{2\pi}}\delta(y)\right]  
\end{align*}

Recalling that constants pass through integrals, and that $\xi$ is independent of $y$, then we have the result

\begin{align*}
    -\xi^2{\bar{\hat{u}}}(\xi,\eta)-\eta^2 {\bar{\hat{u}}}(\xi,\eta)=-\sqrt{\dfrac{1}{2\pi}}\sqrt{\dfrac{1}{2\pi}}
\end{align*}
where $\bar{\hat{u}}$ represents the function $u$ Fourier transformed in both the $x$ and $y$ variables.  Note that we now have a function of $\xi$ and $\eta$, indicating that both spatial independent variables have been transformed.  We are also left with a purely algebraic equation, which is very convenient.  Solving this for $\bar{\hat{u}}$ gives

\begin{align*}
    {\bar{\hat{u}}}&=-\sqrt{\dfrac{1}{2\pi}}\sqrt{\dfrac{1}{2\pi}}\dfrac{1}{\xi^2 + \eta^2} \\
    &=-\dfrac{1}{2}\sqrt{\dfrac{1}{2\pi}}
    \left\{
    \sqrt{\dfrac{2}{\pi}}\dfrac{1}{\xi^2 + \eta^2}\right\}
\end{align*}
Inverting this requires only that we keep careful track of which variables we are working on, and the fact that $\xi$ and $\eta$ are independent.  First, we can find the inverse transform with respect to $\eta$.  This is given symbolically by

\begin{align*}
    \mathscr{F}_y^{-1}\left[{\bar{\hat{u}}}\right]&=-\dfrac{1}{2}\sqrt{\dfrac{1}{2\pi}}
    \mathscr{F}_y^{-1}\left[\left\{
    \sqrt{\dfrac{2}{\pi}}\dfrac{1}{\xi^2 + \eta^2}\right\}\right]
\end{align*}

Note that the term on the right-hand side in braces matches entry \ref{laplaceinv} from the table of Fourier transforms.   To be clear, here we set $a=\xi$, and note that the requirement that $a>0$ implies $\xi \Rightarrow \lvert \xi \rvert$.  This yields the inverse transform (where  $\eta \rightarrow y$) as follows

\begin{align*}
    {\hat{u}}=-\dfrac{1}{2}\sqrt{\dfrac{1}{2\pi}}\left\{\dfrac{1}{\lvert \xi \rvert}\exp(-\xi\lvert y \rvert)\right\}
\end{align*}
Now we have only the $\xi$ variable to invert.  Looking at the expression, we can see that entry number \ref{lntrans} is very close to our needed inverse.  To assist in inverting, we can rearrange the result as follows

\begin{align*}
    {\hat{u}}(\xi,y)=-\dfrac{1}{2}\sqrt{\dfrac{1}{2\pi}}\dfrac{1}{\sqrt{2 \pi}}\left\{ \dfrac{\sqrt{2 \pi}}{\xi}\exp(-\xi\lvert y \rvert)\right\}
\end{align*}
Note that now the shorthand notation for the transform in the $y$-variable is given by an overbar so that we can distinguish between transforms in the two variables.  

In the expression above, the term in braces matches the table exactly, so we find the final inverse transform to give us the result

\begin{align*}
    u(x,y)=-\dfrac{1}{4 \pi} \ln(x^2 + y^2)
\end{align*}
This is a well-known result, and is often called the \emph{fundamental} solution for the Poisson problem on the infinite plane.  Noting that we have defined $d(x,y) = \sqrt{x^2+y^2}$, we can use the properties of logarithms to express this result by
\begin{align*}
    u(x,y)&=-\dfrac{1}{4 \pi} \ln(x^2 + y^2)\\
    &=-\dfrac{1}{2 \pi} \ln\left[ (x^2 + y^2)^{\frac{1}{2}} \right] \\
    &=-\dfrac{1}{4 \pi} \ln[d(x,y)]
\end{align*}\\
As a final note, the solution above is still a solution if we add \emph{any} constant to the problem.  Because we have relaxed the necessity that the functions be square integrable, we have no additional constraint to fix this constant for the infinite plane.  Thus, we must think of our solution as an \emph{equivalence class of solutions}, each differing only by some constant.  This also indicates that the Fourier transform of $\ln(a^2+x^2)$ can, for some cases, involve the addition of a second part containing a delta function $c \delta(\xi)$, with some as yet undetermined constant $c$.  This additional delta function would, upon inversion, generate the constant term that, technically, is part of our potential set of solutions. A plot of this somewhat unusual-looking function appears below.\\
\newpage

{\vspace{4mm}
\centering\includegraphics[scale=1.1]{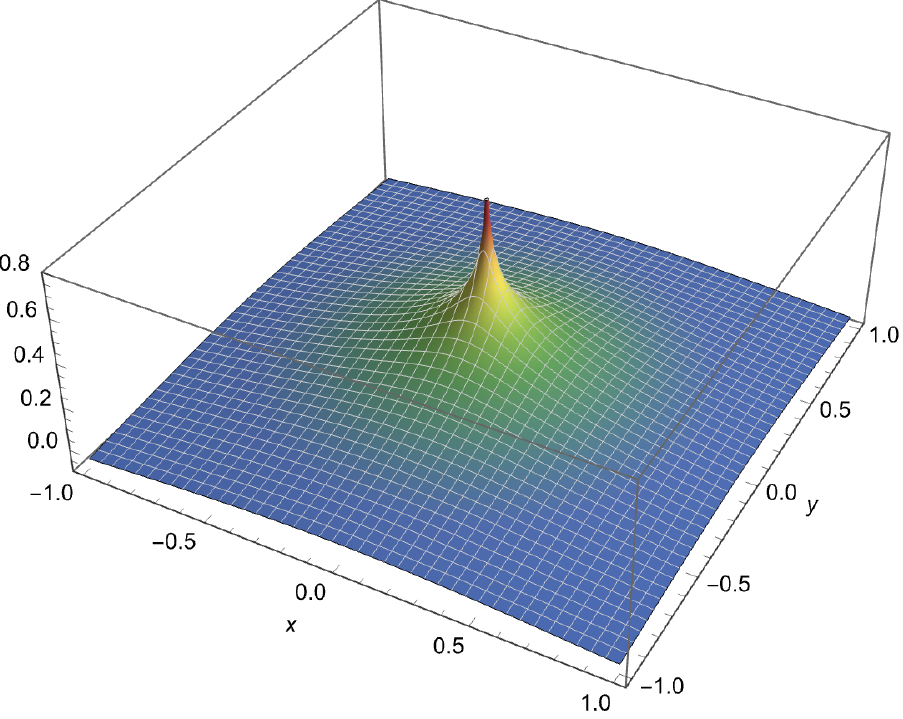}
\captionof{figure}{The solution for the 2-dimensional Poisson equation. Note that the lower boundary of the vertical axis is set to $-1/10$.} 
}
\end{example}
\end{svgraybox}

\section{$^\star$Hyperbolic Equations Redux}\indexme{Fourier transform!wave equations}

As mentioned in Chap.~\ref{introPDEs}, second-order hyperbolic equations can arise from consideration of conservation of momentum applied to a continuum material.  However, they do not lend themselves to the combination of a simple conservation law and a flux law, as was the case for parabolic equations like the heat/diffusion equation.
With the help of the Fourier transform, we are in a position to develop a derivation of the second-order hyperbolic equation from the first-order description.  This is not conventional; the two equations are not frequently related this way in textbooks.  However, approaching the problem this way provides an alternative perspective, and perhaps more physical content than it typically available in the more conventional derivations.

To start, lets take a look again at the Fourier transform of the first-order wave equation.  Recall, we have 

\begin{align}
  &&  \frac{\partial u}{\partial t}& = -c  \frac{\partial u}{\partial x}&&\\
    &I.C.& u(x,0) &= f(x) &&
\end{align}
with the conventional requirements that the functions $u$ and $f$ be square integrable (or, recalling the extensions, are at least tempered distributions).  The Fourier transform of this expression is

\begin{align}
    \frac{d \hat{u}}{dt} &= - i \xi c &&\\
    \hat{u}(\xi,0)&= \hat{f}(\xi) &&
\end{align}
The solution to this problem is reasonably straightforward, and gives an exponential of the form

\begin{equation}
    \hat{u}(\xi,t) = \Cline{\hat{f}(\xi) \exp(-i \xi c t)}\label{firstFT}
\end{equation}

Now, for the second-order equation, recall we have 

\begin{align}
  &&  \frac{\partial^2 u}{\partial t^2}& = c^2  \frac{\partial^2 u}{\partial x^2} &&\\
    &I.C.1& u(x,0) &= f(x) &&\\
    &I.C.2& \left.\frac{\partial u}{\partial t}\right|_{(\xi,0)} &= 0 &&
\end{align}
And the transform is 
\begin{align}
  &&  \frac{d^2 \hat{u}}{d t^2} -c^2 \xi^2 \hat{u}&=0 &&\\
    &I.C.1& \hat{u}(\xi,0) &= \hat{f}(\xi) &&\\
    &I.C.2& -i \xi \hat{u}{(\xi,0)} &= 0 &&
\end{align}
The solution to this second-order homogeneous ODE is fairly straightforward.  It is

\begin{equation}
    u(\xi,t) = \tfrac{1}{2} \Cline{\hat{f}(\xi)\exp(-i c \xi t)} +\tfrac{1}{2} \hat{f}(\xi)\exp(i c \xi t)\label{secondFT}
\end{equation}
At this point, it is possible to see a close relationship between the solution to the first-order wave equation and the second-order wave equation in Fourier space.   Note that the first order wave equation has a solution with the transform of the initial conditions modified by a single exponential.  The solution to the second-order wave equation has a similar solution, except the initial condition is modified by two different exponentials.  You may recall that exponentials of the form $\exp(-i \xi ct)$ or $\exp(i \xi ct)$ simply \emph{shift} their inverse functions back in real space by the amount $\pm c t$, depending on the sign of the argument of the exponential.  Note also, that the two underlined terms in Eqs.~\eqref{firstFT} and \eqref{secondFT} are identical!

In Fourier space, it is easy to see that these two solutions are related.  The first-order solution consists of the initial condition \emph{translated} forward in space by the amount $c t$.  The second-order equation includes two terms, both of which are scaled (by $\tfrac{1}{2}$) versions of the initial condition, and both of which are translated by $c t$, but in opposite directions.

The inversion of these two equations is relatively easy to do when we note that the role of the exponentials is only one of translating the inverted functions.  Thus, for both solutions we obtain the following

\begin{align}
    &\textrm{first-order}& u(x,t) &= f(x-ct) && \\
    & \textrm{second-order}& u(x,t) &= \tfrac{1}{2}f(x-ct)+\tfrac{1}{2}f(x+ct) && 
\end{align}
The analysis of the first- and second-order wave equations on the infinite line using the Fourier transform make the correspondence between the two equations reasonably easy to see.  Both expressions lead to translations in space by the same quantity $c t$.  As discussed previously, the first-order wave equation applied to problems where there is a single wave propagating in a single direction (e.g., a pulse of solute moving along with the fluid in a plug-flow reactor).  The second-order wave equation corresponds to two waves, each with half the original height, and moving in opposite directions.  Here, one might think of a disturbance in a narrow canal (e.g., suppose one throws a large stone into an irrigation canal); the waves start at the location of the initial disturbance, but propagate away in both directions.

All of this was predicated on the requirement that the initial velocity $\partial u/\partial t$ was identically zero in our analysis.  One of the additional feature that make the second-order wave equation  useful is that it contains more unknown constants of integration (thus, it requires one additional condition on $\partial u/\partial t$ in order to determine an unique solutions).  When $\partial u/\partial t \ne 0$, other interesting phenomena are manifest, and the similarity between the first-order and second-order wave equations begins to be less obvious.  However, the examples given above should provide some substantial evidence that, fundamentally, the two equations are indeed related!  The second-order wave equation can be thought of, quite accurately and correctly, as a generalization of the first-order wave equation.

\newpage

\section*{Appendix. Fourier Transform Table \vspace{-10mm}} \indexme{Fourier transform!table of transforms}
\begin{table}[H]
\begin{centering}
\captionsetup{singlelinecheck=off}
\def\arraystretch{2.4}
\caption{Fourier transform pairs.}
{\centering \large
\begin{tabular}{rccccc}
&$f(x)$ & ~~~~~~~~~~~~~ & $\hat{f}(\xi)$  &~~~~~~~~~~~~~Constraints &  \\ 
\cline{1-6}
\rownumber~~~&$\dfrac{\partial f(x)}{\partial x}$&   & $i \xi \hat{f}(\xi)$ \\
\rownumber~~~&$\dfrac{\partial^2 f(x)}{\partial x^2}$&   & $- \xi^2 \hat{f}(\xi)$ \\
\rownumber~~~&$f(ax)$&   & $\dfrac{1}{a} \hat{f}(\xi/a)$ &  $a>0$ &\\
\rownumber~~~&$f(x-a)$&   & $\exp(-i a \xi) \hat{f}(\xi)$\\
\rownumber~~~& $\sqrt{\dfrac{1}{2\pi}}$ &   & $\delta(\xi)$  &\\
\rownumber~~~\newtag{\arabic{magicrownumbers}}{zerotrans}& $0$ &   & $0$  &\\
\rownumber~~~& $1$ &   & $\sqrt{2\pi}\delta(\xi)$  &\\
\rownumber~~~& $\delta(x)$ &   & $\sqrt{\dfrac{1}{2\pi}}$\\
\rownumber~~~\newtag{\arabic{magicrownumbers}}{shiftdelt}& $\delta(x-a)$ &   & $\sqrt{\dfrac{1}{2\pi}}\exp(-i a \xi) $\\
\rownumber~~~\newtag{\arabic{magicrownumbers}}{shiftdelttrans}& $\sqrt{\dfrac{1}{2\pi}}\exp(-i a x) $ &   & $\delta(\xi-a)$\\
\rownumber~~~& $H_0(x)$ &  &$\sqrt{\dfrac{\pi}{2}} \delta(\xi)$\\
\rownumber~~~& $H(x)$ &  &$\sqrt{\dfrac{\pi}{2}} \delta(\xi)+ \dfrac{1}{\sqrt{2 \pi}} \dfrac{1}{i \xi} $\\
\end{tabular}\par 
}
\end{centering}
\end{table}
\begin{table}[H]
\def\arraystretch{3.5}
{\centering \large
\begin{tabular}{rccccc}
&$f(x)$ & ~~~~~~~~~~~~~~ & $\hat{f}(\xi)$  &~~~~~~Constraints~~~~~~ &  \\ 
\cline{1-6}
%
\rownumber~~~& $\begin{cases}H(x-a)\\1-H(a-x)\end{cases}$ &  &$\sqrt{\dfrac{\pi}{2}} \delta(\xi)+ \dfrac{1}{\sqrt{2 \pi}} \dfrac{\exp(-i a \xi)}{i \xi} $\\
\rownumber~~~& $\begin{cases}H(a-x)\\1-H(x-a)\end{cases}$ &  &$\sqrt{\dfrac{\pi}{2}} \delta(\xi)- \dfrac{1}{\sqrt{2 \pi}} \dfrac{\exp(-i a \xi)}{i \xi} $\\
\rownumber~~~\newtag{\arabic{magicrownumbers}}{expx2}&$e^{-a^2 x^2}$&   & $\dfrac{1}{a\sqrt{2}}e^{-\tfrac{\xi^2}{4 a^2}}$ & $a>0$\\
\rownumber~~~&$e^{-a |x|}$&   & $\sqrt{\dfrac{2}{\pi}}\dfrac{a}{a^2+\xi^2}$ & $a>0$\\
\rownumber~~~\newtag{\arabic{magicrownumbers}}{laplaceinv}&$\dfrac{1}{a}e^{-a |x|}$&   & $\sqrt{\dfrac{2}{\pi}}\dfrac{1}{a^2+\xi^2}$ & $a>0$\\
\rownumber~~~& $x e^{-a|x|}$ &   & $-2\sqrt{\dfrac{2}{\pi}} \dfrac{i a \xi}{(\xi^2 + a^2)^2}$ & $a>0$\\
\rownumber ~~~&$\frac{1}{\sqrt{2 a}} \exp\left(-\frac{x^2}{4a} \right)$ & & $\exp(-a \xi^2)$ & $a>0$\\
\rownumber ~~~\newtag{\arabic{magicrownumbers}}{lntrans}&$\ln\left(x^2+a^2\right)$ & & $
\dfrac{\sqrt{2 \pi}}{\lvert \xi \rvert}\exp(-a \lvert \xi \rvert)$ &  $a>0$\\
\rownumber~~~& $x e^{-a x^2}$ &   & $-\dfrac{i \xi  }{2 \sqrt{2}\,\, a^{3/2}} \exp\left(-\dfrac{\xi ^2}{4 a}\right)$ & $a>0$\\
\rownumber~~~& $x^2 e^{-a x^2}$ &   & $-\dfrac{2 a-\xi^2  }{4 \sqrt{2}\,\, a^{5/2}} \exp\left(-\dfrac{\xi ^2}{4 a}\right)$ & $a>0$\\
\rownumber~~~& $\cos(ax)$ &   & $\sqrt{\dfrac{\pi}{2}} [\delta(\xi+a)+\delta(\xi-a)]$\\
\rownumber~~~& $\sin(ax)$ &   & $i\sqrt{\dfrac{\pi}{2}} [\delta(\xi+a)-\delta(\xi-a)]$\\
\end{tabular}
}
\end{table}

\begin{table}[H]
\def\arraystretch{3.5}
{\centering \large
\begin{tabular}{rccccc}
&$f(x)$ & ~~~~~~ & $\hat{f}(\xi)$  &~~~~Constraints~~~~~~ &  \\ 
\cline{1-6}
%
\rownumber~~~& $\textrm{arctan}(x)$ & & $-\dfrac{i}{\xi} \sqrt{\dfrac{\pi }{2}} \exp\left(-| \xi | \right)$ & \\
\rownumber~~~& $\cos(\pi x)[H(x+1)-H(x-1)]$ &  &$\sqrt{\dfrac{2}{\pi }}\dfrac{ \xi  \sin (\xi )}{\pi ^2-\xi ^2} $\\
\rownumber~~~& $\sin(\pi x)[H(x+1)-H(x-1)]$ &  &$-\sqrt{2 \pi }\dfrac{i  \sin (\xi )}{\pi ^2-\xi ^2} $\\
\rownumber~~~& $\textrm{erf}(x)$ & &  $-\dfrac{1}{\xi}\sqrt{\dfrac{2}{\pi }} i \exp\left(-\frac{\xi ^2}{4}\right)$&\\
\rownumber~~~\newtag{\arabic{magicrownumbers}}{erferf}& $\textrm{erf}(a+x)-\textrm{erf}(a-x)$ & &  $\dfrac{2}{\xi}\sqrt{\dfrac{2}{\pi }} \exp\left(-\frac{\xi ^2}{4}\right) \sin (a \xi )$&\\
\rownumber~~~& $\textrm{erf}(x)-\textrm{erf}(x-a)$ &  & $\begin{aligned}
 \dfrac{1}{\xi }  \sqrt{\dfrac{2}{\pi }}& i \exp\left[-\dfrac{1}{4} \xi  (\xi +4 i a)\right] \\
 &\times\left[1-\exp\left(i a \xi \right) \right]\end{aligned}$\vspace{5mm}\\
\protect\raisebox{8mm} \rownumber~~~& $\begin{aligned}T(x;&-a,a)=\\
&\begin{cases}
(1+\frac{1}{a}x),& \hspace{-2mm}-a<x<0 \vspace{2mm} \\
(1-\frac{1}{a}x),& \hspace{-2mm}~0<x<a
\end{cases}\end{aligned}$ &
& \hspace{-5mm} $-\sqrt{\dfrac{2}{\pi }}\dfrac{ (\cos (a \xi )-1)}{a \xi ^2}$ & \hspace{0mm}$a>0$\vspace{3mm}\\
\rownumber~~~\newtag{\arabic{magicrownumbers}}{boxcartable}& $B(x;-a,a)=H(x+a)-H(x-a)$ &   & $\sqrt{\dfrac{2}{\pi}}\dfrac{\sin(a \xi)}{\xi}$ & \hspace{-10mm}$a>0$\\
\end{tabular}
}
\end{table}

\begin{table}[H]
\def\arraystretch{3.5}
{\centering \large
\begin{tabular}{rccccc}
&$f(x)$ & ~~~~~~ & $\hat{f}(\xi)$  &~~~~Constraints~~~~~~ &  \\ 
\cline{1-6}
%
\rownumber~~~\newtag{\arabic{magicrownumbers}}{boxcartable2}& $\begin{aligned}B&(x+b_0;-a,a)\\&=H(x+a+b_0)-H(x-a+b_0)\end{aligned}$ &   & $\sqrt{\dfrac{2}{\pi}}\dfrac{\sin(a \xi)}{\xi}\exp\left( i b_0 \xi\right)$ & \hspace{-10mm}$\begin{aligned}a&>0\\b_0&>0\end{aligned}$\\
\rownumber~~~\newtag{\arabic{magicrownumbers}}{boxcartable3}& $\begin{aligned}B&(x-b_0;-a,a)\\&=H(x+a-b_0)-H(x-a-b_0)\end{aligned}$ &   & $\sqrt{\dfrac{2}{\pi}}\dfrac{\sin(a \xi)}{\xi}\exp\left( -i b_0 \xi\right)$ & \hspace{-10mm}$\begin{aligned}a&>0\\b_0&>0\end{aligned}$\\
\rownumber~~~& $\dfrac{1}{a^2+x^2}$ &   & $\dfrac{1}{a}\sqrt{\dfrac{\pi}{2}}\exp(-a\lvert \xi \rvert)$  &\\
\rownumber~~~& $\dfrac{a}{x^2+a^2}$ &   & $\sqrt{\dfrac{\pi}{2}} \exp\left(-a\lvert\xi\rvert\right)$ & $a>0$\\
\rownumber~~~& $\dfrac{x}{a^2+x^2}$ &   & $-\textrm{sgn}(\xi) i a \sqrt{\dfrac{\pi}{2}}\exp(-a\lvert  \xi \rvert)$  & $a>0$\vspace{4mm}\\
\protect\raisebox{12mm} {\rownumber}~~~& $\begin{aligned} &\dfrac{a^4}{a^4+x^4}\\
 &~\\
 &~\\
 &~\\
 &~\end{aligned}$ 
 &   & $\begin{aligned}\left(\frac{1}{4}+\frac{i}{4}\right)& \sqrt{\pi } a \exp\left(-\frac{(1+i) a \xi }{\sqrt{2}}\right)\times \\
 &\Bigg[e^{\sqrt{2} a \xi }
   \left(e^{i \sqrt{2} a \xi }-i\right) H(-\xi )\\
   &-i e^{i \sqrt{2} a \xi } H(\xi )+H(\xi )\Bigg]\end{aligned}$ &$a>0$\vspace{4mm}\\
\end{tabular}
}
\vspace{-6mm}
\newline

Note.  If you use Mathematica to conduct your Fourier transforms, you will need to set the following as an option for the \texttt{FourierTransform} and \texttt{InverseFourierTransform} commands: \texttt{FourierParameters -> \{0, -1\}}.  This will assure that the \emph{definition} of the Fourier transform is identical to the one adopted in this text.  See the Mathematica help pages for more details.\vspace{-4mm}\\

Function Notes:
\begin{align*}
   &\textrm{erf}& &\textrm{The error function,}~~~  \textrm{erf}(x)=\frac{2}{\sqrt{\pi}}\int_0^{x} \exp(-\eta^2)\,d\eta. \\
     &\textrm{erfc}& &\textrm{The complimentary error function},~~~ 1-\textrm{erf}(t)\\
     &\delta& & \textrm{The Dirac delta function} \\
     &H&& \textrm{The Heaviside function}\\
     &B(x;-a,a) && \textrm{The boxcar function on $x\in[-a,a]$}\\
     &T(t;-a,a)&& \textrm{The symmetric triangle function on $x\in[-a,a]$, with maximum of 1 located at $t=0$}\\
     &\textrm{sgn}(x) &&\textrm{The sign function; sgn$(x)=1,~x>0$, sgn$(x)=-1,~x<0$, sgn$(x)=0,~x=0$}\\
     &H_0(x) = 
    \begin{cases}
        -\frac{1}{2} & x<0 \\
        ~~~\frac{1}{2} & x>0
    \end{cases} && \textrm{The zero-mean Heaviside function}
\end{align*}
\end{table}

\newpage
\section*{Problems}
\subsection*{Practice Problems}
\begin{enumerate}
\item Prove that $\mathscr{F}[f(x-a)] = e^{-i a x} \hat{f}(\xi)$.
\item Prove that $\mathscr{F}[f(ax)] = \frac{1}{|a|} \hat{f}\left(\frac{\xi}{a}\right)$
\item Prove that $(f*g)(x) = (g*f)(x)$.  Hint: this can be done using a transform of variables of the form $\eta = (x-\xi)$.

\item Work out the following transforms in detail.  
\begin{align}
    \mathscr{F}[\cos(x)] &=  \sqrt{\dfrac{\pi}{2}}\left[\delta(\xi +1)+\delta(\xi-1)\right]\\
    \mathscr{F}[\sin(x)] &= i \sqrt{\dfrac{\pi}{2}}\left[\delta(\xi +1)-\delta(\xi-1)\right]
\end{align}
Begin by using the Euler identity to transform the cosine and sine functions to equivalent functions expressed in terms of the complex exponential.  Then, directly integrate the resulting expressions using the definition of the Fourier transform given by Eq.~\ref{fouriertransformdef}.  Finally, you can use the definition of the transform given in the table under entry number \ref{shiftdelt} to identify the delta function in the result.

\item Note the following two transforms (entries number \ref{shiftdelt} and \ref{shiftdelttrans} in the table of Fourier transforms).
\begin{align*}
    \mathscr{F}\left[\delta(x-a)\right]&=\sqrt{\dfrac{1}{2\pi}}\exp(-i a \xi) \\
    \mathscr{F}\left[\sqrt{\dfrac{1}{2\pi}}\exp(-i a x)\right]&=\delta(\xi-a)
\end{align*}
There appears to be some symmetry between the forward and inverse transforms.  To see why this is true, start with the following inverse transform

\begin{equation*}
 \left[\sqrt{\dfrac{1}{2\pi}}\exp(-i a x)\right]=\mathscr{F}^{-1}\left[\delta(\xi-a)\right]
\end{equation*}
Note that using the definition of the inverse, this is
\begin{equation*}
\sqrt{\dfrac{1}{2\pi}}\exp(-i a x)=\dfrac{1}{\sqrt{2 \pi} }\int_{-\infty}^\infty \delta(\xi-a) \exp[i \xi x] \, d\xi
\end{equation*}
Propose the change of variables $\xi = -\eta$ to transform this integral.  Take extra care in transforming the bounds of integration!  After the transformation is complete, switch back to the original variable $\xi$ by simply renaming $\eta \rightarrow \xi$.

Finally, note that $\delta(x-a)=\delta(a-x)$.  This symmetry relationship allows us to develop a result that shows that the Fourier transform of the delta function can be used in reverse.  That is, if you know the transform of the delta function, then you also have the transform of the function $1/\sqrt{2 \pi} \exp[-i a x]$ by reversing the role of the function and its transform, and switching the place of the variables.  To be absolutely clear, we want to show from the development above that the transform given by \ref{shiftdelt} \emph{implies} that the transform given by \ref{shiftdelttrans} is also true.

There was one key feature about the function $\delta(x-a)$ that allows this to be generalized-- what is that?  Can you make a general statement about the conditions for function $f(x)$ such that the Fourier transform can be used in this ``backwards" manner?

\saveenumerate\end{enumerate}

\subsection*{Applied and More Challenging Problems}
\begin{enumerate}
\restoreenumerate

\item Solve the following heat/diffusion equation problem using the Fourier transform.
\begin{align}
 && \frac{\partial u}{\partial t}& = D \frac{\partial^2 u}{\partial x^2}, ~~x\in(-\infty,\infty)& \\
&condition~1& u(x,t)&\rightarrow 0 \textrm{~as~} |x|\rightarrow \infty ~\textrm{for all $x$ and $t$}&\\
&condition~2& \frac{\partial u}{\partial x}&\rightarrow 0 \textrm{~as~} |x|\rightarrow \infty  ~\textrm{for all $x$ and $t$}&\\
 &I.C. & u(x,0)& = \exp(-x^2)~~x\in(-\infty,\infty)& 
\end{align}

\item Solve the following heat/diffusion equation problem using the Fourier transform.
\begin{align}
 &&   \frac{\partial u}{\partial t}& = D \frac{\partial^2 u}{\partial x^2}, ~~x\in(-\infty,\infty) &\\
&condition~1& u(x,t)&\rightarrow 0 \textrm{~as~} |x|\rightarrow \infty ~\textrm{for all $x$ and $t$}&\\
&condition~2& \frac{\partial u}{\partial x}&\rightarrow 0 \textrm{~as~} |x|\rightarrow \infty  ~\textrm{for all $x$ and $t$}&\\
& I.C.& u(x,0) &= \delta(x) &
\end{align}

\item Solve the following heat/diffusion equation problem using the Fourier transform.  Solve this problem two ways. (1) First, try the transformation $u(x,t)=w(x,t)\exp(-k t)$ to develop a new problem in $w$ that does not contain the reaction term, and then solve that using the Fourier transform.  (2) Second, try solving the problem directly using the Fourier transform.  Which method seemed easier to you?
\begin{align}
 &&   \frac{\partial u}{\partial t}& = D \frac{\partial^2 u}{\partial x^2} -k u, ~~x\in(-\infty,\infty)& \\
&condition~1& u(x,t)&\rightarrow 0 \textrm{~as~} |x|\rightarrow \infty ~\textrm{for all $x$ and $t$}&\\
&condition~2& \frac{\partial u}{\partial x}&\rightarrow 0 \textrm{~as~} |x|\rightarrow \infty  ~\textrm{for all $x$ and $t$}&\\
& I.C.& u(x,0) &= \delta(x)&
\end{align}

\item \textbf{A problem of pure convection.}  In some instances, such as in packed beds or in rivers, it is permissible to make the approximation that the transport is primarily \emph{convective} (or sometimes the word \emph{advective} is used), meaning that dissolved species in the fluid (gas or liquid) phase move with the fluid with very little spreading (or dispersion).  If one has, for example, a first-order reaction representing the conversion of a chemical compound in a river to some end product, one possible model for the reactant given by the following first-order PDE

\begin{align}
 &&   \frac{\partial u}{\partial t}& = -v_0 \frac{\partial u}{\partial x}, ~~x\in(-\infty,\infty) &\\
&condition~1& u(x,t)&\rightarrow 0 \textrm{~as~} |x|\rightarrow \infty ~\textrm{for all $x$ and $t$}&\\
& I.C.& u(x,0) &= \beta \exp\left(- \alpha x^2\right)
\end{align}
This model assumes that the initial condition is a Gaussian distribution of concentration (changing in distance, $x$ along the river path, and constant in planes perpendicular to $x$).  This assumes that the river can be approximated as having a 1-dimensional process geometry.

Solve the problem above using Fourier transforms.  Then, plot the solution to the problem in the form $(1/b) u(x,t)$  for $x\in[-10,125]$ \si{m} and $u\in[0,1]$. Use $v_0=10$ \si{km \cdot h^{-1}}, $\alpha=0.1$, and $\beta=1$.  Plot the function at times $t=0.1, 1, 5$ and $10$ hours.  Note that we have not specified the \emph{units} for $\beta$, but they must be the same as those for $u$; thus $u/b$ is dimensionless.

\item \textbf{A problem of pure convection and reaction.}  Assume that for the pure convection in a river problem given above, we now also have the presence of a first-order transformation reaction.  The concentration for an initially Gaussian concentration would be specified by the following first-order PDE

\begin{align}
 &&   \frac{\partial u}{\partial t}& = -v_0 \frac{\partial u}{\partial x}-\gamma u, ~~x\in(-\infty,\infty) &\\
&condition~1& u(x,t)&\rightarrow 0 \textrm{~as~} |x|\rightarrow \infty ~\textrm{for all $x$ and $t$}&\\
& I.C.& u(x,0) &= \beta \exp\left(- \alpha x^2\right) &
\end{align}
This model assumes that the initial condition is a Gaussian distribution of concentration (changing in distance, $x$ along the river path, and constant in planes perpendicular to $x$).  This assumes that the river can be approximated as having a 1-dimensional process geometry.  Here, $\gamma$ is a first-order kinetic reaction rate; set $\gamma = 0.1$ \si{h^{-1}}.

Solve the problem above using Fourier transforms.  Then, plot the solution to the problem in the form $(1/b) u(x,t)$  for $x\in[-10,125]$ and $u\in[0,1]$. Use $v_0=10$  \si{km \cdot h^{-1}}, $\alpha=0.1$, and $\beta=1$.  Plot the function at times $t=0.1, 1, 5$ and $10$ hours.  Note that we have not specified the \emph{units} for $\beta$, but they must be the same as those for $u$; thus $u/b$ is dimensionless.

\item \textbf{A convection-diffusion problem.}  Suppose we would like to solve the following convection-diffusion problem using the Fourier transform.  

\begin{align}
 &&   \frac{\partial u}{\partial t}& = D \frac{\partial^2 u}{\partial x^2}- v_0 \frac{\partial u}{\partial x}, ~~x\in(-\infty,\infty) &\\
&condition~1& u(x,t)&\rightarrow 0 \textrm{~as~} |x|\rightarrow \infty ~\textrm{for all $x$ and $t$}&\\
&condition~2& \frac{\partial u}{\partial x}&\rightarrow 0 \textrm{~as~} |x|\rightarrow \infty  ~\textrm{for all $x$ and $t$}&\\
& I.C.& u(x,0) &= \beta \exp(-x^2) &
\end{align}
Note: you will want to use the transform given by \#4 in the table of transforms.

\item \textbf{A convection-diffusion-reaction problem.}  Suppose we would like to solve the following convection-diffusion problem using the Fourier transform.  

\begin{align}
 &&   \frac{\partial u}{\partial t}& = D \frac{\partial^2 u}{\partial x^2}- v_0 \frac{\partial u}{\partial x}- k u, ~~x\in(-\infty,\infty) &\\
&condition~1& u(x,t)&\rightarrow 0 \textrm{~as~} |x|\rightarrow \infty ~\textrm{for all $x$ and $t$}&\\
&condition~2& \frac{\partial u}{\partial x}&\rightarrow 0 \textrm{~as~} |x|\rightarrow \infty  ~\textrm{for all $x$ and $t$}&\\
& I.C.& u(x,0) &= \beta \delta(x-x_0) &
\end{align}
Note: you will want to use the transform given by \#4 in the table of transforms.

\item\textbf{The wave equation in an infinite medium}.\label{fourierwave}  The wave equation behaves very differently from the heat/diffusion equation. In the versions that we have developed so far, waves travel without damping.  Thus, in an infinite medium (i.e., mathematically, the real line in 1-D), an initial disturbance located near $x=0$ at $t=0$ will travel outward from the disturbance toward infinity for all $t>0$.  If the wave equation is the second-order one (which will be written out below), then an initial disturbance travels out in both the $+x$ and $-x$ directions towards infinity with increasing time, $t$.  We can actually use the Fourier transform to solve this equation in an infinite medium just as we would for the heat/diffusion equation.  The steps below will walk us through such a solution.

\begin{enumerate}
    \item Begin by transforming the wave equation \emph{and} its initial condition by conducting a Fourier transform with respect to the variable $x$.  For the wave equation, assume the following balance equation and ancillary conditions
    \begin{align*}
        \dfrac{\partial^2 u}{\partial t^2} & = c^2 \dfrac{\partial^2 u}{\partial u^2} \\
        u(x,0)&=H(x+1)-H(x-1)\\
       \left. \dfrac{\partial u}{\partial t}\right|_{(x,0)}& = 0\\
       u(x,t) &\textrm{~and~} \partial u/\partial x ~\rightarrow 0 \textrm{~as~} \lvert x\rvert \rightarrow \infty
    \end{align*}
\noindent Note that here we have both and initial condition and an initial derivative for $u$ with respect to time.  As we have discussed, we need one ancillary condition for each derivative in time and in space.  Because we have a second-order derivative in time for this equation, we need \emph{two} conditions to eliminate the unknown constants of integration that will result.  Also recall that the conditions on $u$ in space as $\lvert x \rvert\rightarrow \infty$ serve as the two necessary ancillary conditions in space.

\item The result from your previous step should take the form
\begin{align*}
    \dfrac{\partial^2 \hat{u}}{\partial t^2} & + c^2 \xi^2 \hat{u}=0 \\
    \left. \dfrac{d \hat{u}}{dt}\right|_{(\xi,0)} &= 0\\
    \hat{u}(\xi,0) & = \dfrac{1}{\xi}\sqrt{\dfrac{2}{\pi}} \sin(\xi)
\end{align*}

This equation may initially look challenging, but it is in fact just a linear second-order ODE.  We have studied the solutions to such equations in Chap.~\ref{ODEChap}.  Recall that the solution for an equation like this is found by first developing the characteristic equation, and then solving for the roots of that equation.  The form of the roots tell you which of three cases the solution belongs to.  Do not forget to use the two ancillary conditions to eliminate the two constants of integration!

\item Your result for the solution is not on the Fourier transform table as such.  However, we can generate the solution as follows.  First, obtain the transforms given on the table as entries \ref{boxcartable2} and \ref{boxcartable3}.  Note that each of those involve an imaginary exponential function that can be converted to trigonometric functions by using the Euler identity: $\exp(i b_0 \xi) = \cos(b0 \xi)+i\sin(b_0 \xi)$ and $\exp(-i b_0 \xi) = \cos(b0 \xi)-i\sin(b_0 \xi)$.  So, make this substitution to these to transform pairs.  Then, average the two results (i.e., take $\tfrac{1}{2}$ of the sum of the two transforms).  You should now have a result that matches your solution in Fourier space, and this will allow you to compute the inverse transform.  Your final result should be as follows
\begin{equation}
    u(x,t) = \underbrace{\dfrac{1}{2}[H(x+1+ct)-H(x-1+ct)]}_\textrm{half-height wave traveling in negative direction}+ \underbrace{\dfrac{1}{2}[H(x+1-ct)-H(x-1-ct)]}_\textrm{half-height wave traveling in positive direction}
\end{equation}
So the initial single perturbation (the boxcar function) is split into two boxcar functions of half the height, and traveling with velocity $c$ in opposite directions.

\item Suppose we think of the area under the wave times a unit thickness as being proportional to mass. At $t=0$ we can think of the initial wave as as having two parts moving at velocity $c$ (for the positive half), and $-c$ (for the negative half); this is what leads to its eventual splitting for $t>0$.  We have discussed previously that this wave equation conserves energy.  If we use the kinetic energy alone to represent the total, then each wave has an energy equal to $A c^2$.  Compute the energy in the initial condition, and show that this amount is conserved for all $t>0$.
\end{enumerate}

\item \textbf{The Poisson equation in an infinite medium}.  Solve the following problem using Fourier transforms for $a>0$, $b>0$.  Your solution should be a shifted version of the solution for the Poisson equation given in the example problems.

\begin{align*}
    &&\frac{\partial^2 u}{\partial x^2}+ \frac{\partial^2 u}{\partial y^2}&=-\delta(x-a)\delta(y-b)\qquad -\infty < x < \infty, ~-\infty < x < \infty\\
&B.C.1&    u(x,y) &\rightarrow 0 ~~as~ x^2+y^2\rightarrow \infty \\
&B.C.2a&    \left.\frac{\partial u}{\partial x}\right|_{(x,y)} &\rightarrow 0 ~~as~ x^2+y^2\rightarrow \infty\\
&B.C.2b&    \left.\frac{\partial u}{\partial y}\right|_{(x,y)} &\rightarrow 0 ~~as~ x^2+y^2\rightarrow \infty\\
\end{align*}

\saveenumerate
\end{enumerate}
\vspace{5mm}

\textbf{Some problems involving hypergeometric functions.}  The hypergeometric functions \indexme{hypergeometric function} are series solutions to differential equations.  There existence has been known about since at least the mid-1600s.  The famous mathematicians  Leonhard Euler and  Carl Friedrich Gauss studied the properties of these mathematical functions in Late 1700s and early 1800s.  The \emph{generalized} hypergeometric functions \indexme{hypergeometric function!generalized} appear to have been introduced by \citet{Clausen1828} in 1828.  

There are many ways of representing the hypergeometric functions.  We will not detail this here, but mention only a few of them and provide their definitions.  The following two \emph{generalized hypergeometric} functions are relevant to the material that will follow.

\begin{align}
    _0F_1(x;a) &= \sum_{n=0}^\infty\frac{1}{n!} \frac{\Gamma(a)}{\Gamma(a+n)} x^n\\
    _0F_2(x;a;b) &= \sum_{n=0}^\infty\frac{1}{n!}  \frac{\Gamma(a)}{\Gamma(a+n)}\frac{\Gamma(b)}{\Gamma(b+n)} x^n
\end{align}
Note that these are presented in a non-standard form in terms of the independent parameters $a$ and $b$; this is done to make the presentation clearer.  If you study generalized hypergeometric functions in the future, the notation will be slightly different from this.  These functions have well-defined infinite series in general, and many computer mathematics packages (such as Mathematica) have definitions for these functions.

In addition to these, we introduce a function that is related to the hypergeometric functions called the Airy $Ai(x)$ function.  This special function was discovered by  British astronomer George Biddell Airy (1801--1892).  It is defined in terms of the hypergeometric functions as follows.

\begin{equation}
    Ai(x) = \frac{1}{3^{2/3}\Gamma\left(\frac{2}{3}\right)} {\,\,_0F_1}\left(\tfrac{x^3}{9},\tfrac{2}{3}\right)
    -
    x \frac{1}{3^{1/3}\Gamma\left(\frac{1}{3}\right)} {\,\,_0F_1}\left(\tfrac{x^3}{9},\tfrac{4}{3}\right)
\end{equation}
There are many equivalent expressions for the Airy $Ai(x)$ function, but we this form has been chosen to help illustrate how the generalized hypergeometric functions provide the framework for many solutions to higher-order PDEs.   Because the Airy function shows up in applications to physical systems (e.g., optics,    quantum mechanics, waves in fluids), it has been given its own symbol, and is defined in many computer mathematics packages (such as Mathematica).  

Now, the reason that these functions are important is that they allow us to describe the solutions to PDEs with higher-order derivatives in space.  Toward that goal, we now note the following Fourier transform pairs.

\begin{table}[H]
    \centering
    \def\arraystretch{2.4}
    \begin{tabular}{c|c|c}
    ~ &function  & transform \\
    \hline
    A.  & $~~\frac{1}{\sqrt{2\pi}} \exp\left(\alpha_3 i k^3 t \right)~~$  & 
    $\frac{1}{(3\alpha_3 t)^{1/2}}\textrm{Ai}\left(\frac{x}{(3\alpha_3 t)^{1/2}} \right) $\\
    B.  &  $~~\frac{1}{\sqrt{2\pi}} \exp\left(-\alpha_4 k^4 t \right)~~$ & 
    $\frac{1}{8 \pi (\alpha_4 t)^{3/4}}\left[   8 \sqrt{\alpha_4 t} \Gamma\left(\tfrac{5}{4} \right) {_0F_2}\left(\tfrac{x^4}{256 \alpha_4 t},\tfrac{1}{2},\tfrac{3}{4}\right) -x^2\Gamma\left(\tfrac{3}{4} \right) {_0F_2}\left(\tfrac{x^4}{256 \alpha_4 t},\tfrac{5}{4},\tfrac{3}{2}\right)    \right]   $
    \end{tabular}
    \caption{Fourier transforms for two functions that arise from (spatial) higher-order PDEs.}
\end{table}

With this information available, it is possible for us to examine two higher-order (in space) PDEs.  This is done in the following two questions.

\begin{enumerate}
\restoreenumerate

\item \textbf{Third-order PDE.}  \indexme{high-order PDEs!third order}  As mentioned above, third-order derivatives in space arise in several interesting physical applications, such as in optics, the description of water waves, or in modeling vapor deposition for the epitaxial growth of crystal films \citep{lam2001multiscale}.  Solve the following third-order PDE with delta initial condition using Fourier transforms.

\begin{align}
 &&   \frac{\partial u}{\partial t}& = -\alpha_3 \frac{\partial^3 u}{\partial x^3}, ~~x\in(-\infty,\infty) &\\
&condition~1& u(x,t)&\rightarrow 0 \textrm{~as~} |x|\rightarrow \infty ~\textrm{for all $x$ and $t$}&\\
&condition~2& \frac{\partial u}{\partial x}&\rightarrow 0 \textrm{~as~} |x|\rightarrow \infty  ~\textrm{for all $x$ and $t$}&\\
&condition~3& \frac{\partial^2 u}{\partial x^2}&\rightarrow 0 \textrm{~as~} |x|\rightarrow \infty  ~\textrm{for all $x$ and $t$}&\\
& I.C.& u(x,0) &= \beta \delta(x) &
\end{align}
Note here that $\alpha_3$ has the units of $L^3/T$ (length cubed per unit time).  After finding the solution, plot the solution for $\alpha_3 = 2$ \si{cm^3/hr} and $\beta=10$.  If this is done in Mathematica, note that the Airy $Ai(x)$ function is given by
 \begin{equation*}
   Ai(x)  = \texttt{ AiryAi[x]}
\end{equation*} 
For the plot use $x\in[-10,10]$, and for the range use the interval $u\in[-10,10]$.  Use three plot times equal to $t=0.1, 1,$ and $10$.

\item \textbf{Fourth-order PDE.}  \indexme{high-order PDEs!fourth order} Fourth-order derivatives in space arise in a number of interesting problems.  While they are not generally well known (even in those that study PDEs), they do have applications including (a) describing ice formation, (b) modeling liquid flows in the lungs, (c) and smoothing of raster images in computational image processing \citep{greer2006fourth,you2000fourth}.

Find the solution to this problem using Fourier transforms.

\begin{align}
 &&   \frac{\partial u}{\partial t}& = -\alpha_4 \frac{\partial^4 u}{\partial x^4}, ~~x\in(-\infty,\infty) &\\
&condition~1& u(x,t)&\rightarrow 0 \textrm{~as~} |x|\rightarrow \infty ~\textrm{for all $x$ and $t$}&\\
&condition~2& \frac{\partial u}{\partial x}&\rightarrow 0 \textrm{~as~} |x|\rightarrow \infty  ~\textrm{for all $x$ and $t$}&\\
&condition~3& \frac{\partial^2 u}{\partial x^2}&\rightarrow 0 \textrm{~as~} |x|\rightarrow \infty  ~\textrm{for all $x$ and $t$}&\\
&condition~4& \frac{\partial^3 u}{\partial x^3}&\rightarrow 0 \textrm{~as~} |x|\rightarrow \infty  ~\textrm{for all $x$ and $t$}&\\
& I.C.& u(x,0) &= \beta\delta(x) &
\end{align}
Note here that $\alpha_3$ has the units of $L^4/T$ (length to the fourth per unit time).

After finding the solution, plot the solution for $\alpha_4 = 10$ \si{cm^4/hr} and $\beta=10$.  If this is done in Mathematica, note that the hypergeometric function $0_0F_2$ is given by the command

\begin{equation*}
   _0F_2(x;a;b)  = \texttt{ HypergeometricPFQ[\{\}, {a, b}, x]}
\end{equation*}
For the plot use $x\in[-10,10]$, and for the range use the interval $u\in[-2,5]$.  Use three plot times equal to $t=0.1, 1,$ and $10$.

\item \textbf{Fifth order PDE.}  While it may or may not be a surprise at this juncture, it turns our that there are applications for fifth-order PDEs.  In it simplest possible form, one can write \indexme{high-order PDEs!fifth order}

\begin{align}
 &&   \frac{\partial u}{\partial t} &= -\alpha_5 \frac{\partial^5 u}{\partial x^5},~~ x\in(-\infty, \infty)&\\
  &I.C. & u(x,0) & = \delta(x)  &
 \intertext{plus the conventional conditions on the function and its derivatives.}\nonumber
\end{align}
Such equations have been studied in the context of linearized versions of the (generally nonlinear) fifth-order Korteweg-De Vries (KDV) Equation.\indexme{Korteweg-De Vries equation}  The KDV equation first became popular because its solutions involved traveling waves called \emph{solitons}.  The study of such solutions occurred in the mid 1800s, and then again a resurgence of research occurred in the 1950s and onward.  The fifth-order version of this equation has been used to describe, among other phenomena, the motion of plasma waves or capillary-gravity water waves \citep{bridges2002stability}.  The solution to this problem has the somewhat formidable-looking solution involving generalized hypergeometric functions with seven parameters.

\begin{align}
    u(x,t)& = \frac{\sqrt{\frac{1}{2} \left(5+\sqrt{5}\right)}}{1171875000\, \pi \,
   (\alpha_5\,t)^{9/5}}\Bigg\{
    16500 \,x^3 \Gamma
   \left(-\frac{21}{5}\right)\, \Bigg[x^5 \,\,
   _0F_7\left(\{\};\frac{11}{10},\frac{6}{5},\frac{13}{10},\frac{3}{2},\frac{8}{5},\frac{17}
   {10},\frac{9}{5};\frac{x^{10}}{(50000\,\alpha_5\,t)^2}\right) \nonumber\\
   -&8400 \,(\alpha_5\,t) \,\,
   _0F_7\left(\{\};\frac{1}{2},\frac{3}{5},\frac{7}{10},\frac{4}{5},\frac{11}{10},\frac{6}{5
   },\frac{13}{10};\frac{x^{10}}{(50000\,\alpha_5\,t)^2}\right)\Bigg] \nonumber\\
   %
   +&9375
   \left(\sqrt{5}-1\right) (\alpha_5\,t)^{1/5}\, x^2\, \Gamma \left(-\frac{12}{5}\right) \Bigg[4200
   (\alpha_5\,t) \,
   _0F_7\left(\{\};\frac{2}{5},\frac{1}{2},\frac{3}{5},\frac{7}{10},\frac{9}{10},\frac{11}{1
   0},\frac{6}{5};\frac{x^{10}}{(50000\,\alpha_5\,t)^2}\right) \nonumber
   \end{align}
   \begin{align}
   -&x^5 \,
   _0F_7\left(\{\};\frac{9}{10},\frac{11}{10},\frac{6}{5},\frac{7}{5},\frac{3}{2},\frac{8}{5
   },\frac{17}{10};\frac{x^{10}}{(50000\,\alpha_5\,t)^2}\right)\Bigg]\nonumber \\
   -&31250
   \left(\sqrt{5}-1\right) (\alpha_5\,t)^{2/5} \,x\, \Gamma \left(-\frac{8}{5}\right) \Bigg[1800\, (\alpha_5\,t) \,\,
   _0F_7\left(\{\};\frac{3}{10},\frac{2}{5},\frac{1}{2},\frac{3}{5},\frac{4}{5},\frac{9}{10}
   ,\frac{11}{10};\frac{x^{10}}{(50000\,\alpha_5\,t)^2}\right)\nonumber \\
   -&x^5 \,
   _0F_7\left(\{\};\frac{4}{5},\frac{9}{10},\frac{11}{10},\frac{13}{10},\frac{7}{5},\frac{3}
   {2},\frac{8}{5};\frac{x^{10}}{(50000\,\alpha_5\,t)^2}\right)\Bigg]\nonumber \\
   +&4418826048\,\, (\alpha_5\,t)^{3/5}\, \Gamma
   \left(-\frac{39}{5}\right) \Bigg[600\, (\alpha_5\,t) \,
   _0F_7\left(\{\};\frac{1}{5},\frac{3}{10},\frac{2}{5},\frac{1}{2},\frac{7}{10},\frac{4}{5}
   ,\frac{9}{10};\frac{x^{10}}{(50000\,\alpha_5\,t)^2}\right)\nonumber\\
   -&x^5 \,
   _0F_7\left(\{\};\frac{7}{10},\frac{4}{5},\frac{9}{10},\frac{6}{5},\frac{13}{10},\frac{7}{
   5},\frac{3}{2};\frac{x^{10}}{(50000\,\alpha_5\,t)^2}\right)\Bigg]\Bigg\}
\end{align}

Using Mathematica (or the symbolic mathematics program of your choice) do the following.  (a) Illustrate that this equation is indeed a solution to the fifth-order PDE given above.  (b) For the conditions $\alpha_5=1$ \si{m^5/h}, plot the solution on $x\in[-10,10]$ with $u\in[-1,1]$.  Plot these solutions for times $t=0.1, 1$ and $10$ hours.  What happens to this solution as $t\rightarrow \infty$?  Just examining the terms involving time (those outside of the hypergeometric functions), can you \emph{prove} what the solution tends to as $t\rightarrow \infty$?
\end{enumerate}
\abstract*{This is the abstract for chapter 00}

\begin{savequote}[0.55\linewidth]
``Man follows only phantoms."

\qauthor{ Pierre-Simon Laplace, discoverer of the Laplace transform }
\end{savequote}


\chapter{Laplace Transforms}
%
\def\CHAP {chapter11_Laplace_transforms}\indexme{Laplace transform}
%
\section{Introduction}

The Laplace and Fourier transforms have many similarities.\indexme{Laplace transform}  In some texts, the Laplace transform is discussed first; in others the Fourier transform is examined first.  Because the Fourier transform has a connection with the Fourier series, and because it's inverse is easily understood, it has been investigated first in this text.  Although the Laplace transform is often viewed as being less \emph{intuitive} than the Fourier transform, it still contains much that appeals to the physical aspects of the problem.  Here, we approach the Laplace transform as an intuitive extension of the Fourier transform; this extension is motivated by increasing the space of functions for which the transform is defined.

\section{Terminology}

\begin{itemize}
\item \textbf{Laplace transform.}  An integral transform for on the interval $[0,\infty)$, usually used to transform the time variable.  Laplace transforms can be though of as a generalization of the Fourier transform, where the kernel of the transform is augmented by a real-valued decaying exponential.  Thus, the kernel is $\exp(-\alpha t)\exp(-i \xi t)$.  Usually, this is written by defining $s=\alpha + i \xi$, so that the kernel is $\exp(-s t)$.   In this chapter, the Laplace transform is of a function $f$ then defined by  

\begin{equation*}
    F(s) = \int_{\tau=0}^{\tau=t} f(\tau) \exp(-s \tau) \, d\tau
\end{equation*}
Because there is very little chance of misunderstandings, in the presentation that follows we forgo the formal use of the variable of integration ($\tau$) in favor of using simply $t$ inside integrals.  This is technically poor form, but it is consistent with the notation of most textbooks.
\\

\item \textbf{Complex plane.}  Recall that complex numbers involve both a real and an imaginary component.  Thus, $s=\alpha + i \xi$ is a complex number.  Complex numbers can also be thought of as a pair, $(\alpha,\xi)$, where $\alpha$ gives the real component, and $\xi$ the imaginary component. Such numbers can be plotted on the complex plane in the same way that one might plot a pair of real numbers on the plane.  One can also define a \emph{function} of complex numbers, $f(\alpha+i \xi)=f(\alpha,\xi)$.  Such functions can be represented as a surface with height $f$ and coordinates on the complex plane of $(\alpha, \xi)$.  \\

\item \textbf{Convolution in Laplace transforms.}  The convolution has been discussed a number of times, and it arises again in the study of Laplace transforms.  Because the Laplace transform is on the half-line, the convolution is similarly restricted.  Thus, for two functions $f$ and $g$, the half-line convolution can be given by 
\begin{equation*}
    (f*g)(t)= \int_0^\infty f(t-\tau)g(\tau)\, d\tau = \int_0^\infty f(t)g(t-\tau)\, d\tau
\end{equation*}\\

\item \textbf{Contour integral.}  We have defined the complex plane, and functions that might be defined on the complex plane.   There is an analogue to \emph{line integrals} from calculus that generalizes integration to the complex plane.  While such integrals (and associated theorems) are not difficult, they do require significant study of complex analysis.  This topic is outside of the intended breadth of this text.  However, even though we do not describe the process of how to compute contour integrals, it is worthwhile to define them.  The inverse Laplace transform is most commonly understood as a contour integration process.  
\end{itemize}

\section{The Laplace Transform}

One of the motivations for the Laplace transform is to increase the types of functions that can be analyzed with the method.  For example, recall that for a Fourier transform, a function had to at least have a finite square integral on $x \in (-\infty, \infty)$.  Thus, for example, the simple-looking function $f(x) = x$ has no Fourier transform.

One way of thinking about the Laplace transform is that it extends the kinds of functions that can be transformed by including a conditioning function within the transform itself.  The role of this conditioning function is to eradicate ``fast growing" components of a function that prevent it from having a defined square integral.  Transforms, like the Fourier or Laplace transforms, are qualified as \emph{functionals} that map entire \emph{functions} from one space to another.  Recall, the Fourier transform of the function $f(x)$, assuming that it exists, would be given in the complex exponential form.  Laplace transforms are often used to transform functions that are dependent on time.  Thus, by convention we will adopt $t$ as the independent variable for this chapter.  For easy correspondence with the Laplace transform, we will write out the Fourier transform in the independent variable $t$.  The result is as follows.
\begin{equation}
    \hat{f}(\xi) = \mathscr{F}[f(t)] = \frac{1}{\sqrt{2 \pi}}\int_{-\infty}^{\infty} f(t)  e^{-i \xi t}\, dt
\end{equation}
In the case of Fourier transforms, the functional \emph{uniquely} maps functions with independent variable $t$ (here, noting our change in notation for this chapter only) to a new function, $F$, entirely in terms of the variable $\xi$, i.e., 
\begin{equation}
    \mathscr{F}[f(t)]\rightarrow \hat{f}(\xi)
\end{equation}
As a matter of vocabulary, the function $e^{-i \xi t}$ is called the \emph{kernel} of the Fourier transform.\index{Fourier transform!kernel}\indexme{kernel}

Now, there is nothing  inherently sacred about this particular transformation (other than its long history and relatively direct physical interpretation!)  In  fact, were we to modify the definition slightly to some other well-defined functional (by changing the kernel function), the transformation would still be potentially valuable.  There are actually a large  number of integral  transforms,  each adopting a different kernel function that emphasizes particular features of the transformed function.  For example, one may find on Wikipedia that, at this writing, there are nearly 20 different kinds of integral transforms defined there.  Thus, the exploration of different kernels for conducting the transforms has a place in applied mathematical analysis, so long as the proposed kernel can be shown to have some particular value (such as allowing the solution to a problem that is not otherwise solvable).

For our discussion, let's consider one of the \emph{simplest} modifications of the Fourier transform that we can devise that helps expand the space of transformable functions.  As a concrete example, let's consider the function $f(t)=x$ on $-\infty < t <  \infty$.  Clearly, the Fourier transform for this function is not well-defined.  The  Fourier transform would give
\begin{equation}
    \mathscr{F}[f(t)] = \frac{1}{\sqrt{2 \pi}}\int_{-\infty}^{\infty} t \cdot e^{-i \xi t}\, dt
\end{equation}
Because $e^{-i\xi t}$ is an oscillatory function, this integral is unbounded as $\lvert t \rvert \rightarrow \infty$.  However, suppose we decide that instead of $e^{-i \xi t}$, we use the kernel $e^{-\alpha t}e^{-i \xi t}$,\indexme{Laplace transform!kernel} where $\alpha>0$ is some real number.  What would such an addition do to the behavior of the kernel function?  For one thing, we know that (for positive values of $t$ at least), $e^{-\alpha t}$ is a function that decays quickly relative to, for example, \emph{every} polynomial.  In this case, what we mean by \emph{quickly} is that for any polynomial in $t$, $P_n(t)$ (where $n$ is the degree or order of the polynomial) then we always have

\begin{equation}
   \underset{t\rightarrow\infty}{\lim} \,{P_n(t)}{e^{-\alpha t}} = 0
\end{equation}
for any $\alpha>0$.  (To see this, one can simply expand the exponential as its power series, which converges everywhere; the exponential has polynomial terms of all orders).  The practical outcome of this addition is that now functions that were previously not integrable are integrable with this new kernel.  However, we do have one detail to attend to.  Again, let's return to our definition for the Fourier transform, but adopting this modification.  For now, let's refer to our modified Fourier transform as the $\mathscr{F}'$ transform.  For $f(t)=t$, we would have the result

\begin{equation}
    \mathscr{F}'[f(t)] = \frac{1}{\sqrt{2 \pi}}\int_{-\infty}^{\infty} t \cdot e^{-\alpha t} e^{-i \xi t}\, dt
\end{equation}
The exponential $e^{-\alpha t}$ accomplishes our goal for $t\ge 0$, but actually makes the problem worse for $t<0$.  In other words, for $t<0$, the sign of the exponential is \emph{positive}, and the exponential grows without bounds.  There are some functions that still converge for this transform on the bounds $-\infty < t < \infty$; ~in this application, the transformation is known as the \emph{bilateral} or \emph{two-sided} Laplace transform.  This is a somewhat specialized application that will not be further discussed in this text.  Instead of insisting on convergence over  $-\infty < t < \infty$, we will adopt the convention that our $\mathscr{F}'$ transform is used only for functions defined on $t\in[0,\infty)$.  Under these circumstances the increasing exponential is no longer a problem.  Restricting the domain of integration, we can now compute the integral 

\begin{equation}
    \mathscr{F}'[t] =\dfrac{1}{\sqrt{2 Pi}} \int_{0}^{\infty} t \cdot e^{-\alpha t} e^{-i \xi t}\, dt
\end{equation}
this integral, although an improper one, can be computed. 
\begin{equation}
    \mathscr{F}'[t] = \dfrac{1}{\sqrt{2 \pi}}\int_{0}^{\infty} t \cdot e^{-\alpha t} e^{-i \xi t}\, dt = \dfrac{1}{\sqrt{2 \pi}}\dfrac{1}{(a+i\xi)^2}
\end{equation}
So, our new transform seems to have an the benefit of having a convergent (and thus computable) integral for functions like $f(t)=t$ that previously could not be transformed.

This new transform is known as the \emph{Laplace} transform.  However, before we complete our definition, we must make note of two additional changes that are conventional for the Laplace transform.  First, the constant $1/\sqrt{2 \pi}$ is not included in the forward transform; rather, it is combined with the same constant in the \emph{inverse} transform.  This changes nothing about the transform pairs, only how the constant $1/(2\pi)$ is divided between the forward and inverse transforms.  Second, again by convention, we set $s=\alpha+i \xi$.   These two changes provide the conventional form of the Laplace transform, denoted by $\mathscr{L}$. 

\begin{equation}
    \mathscr{L}[f(t)] = \frac{1}{\sqrt{2 \pi}}\int_{0}^{\infty} f(t) \cdot e^{-s t} \, dt
\end{equation}
Now note that making these changes to our previous analysis of the function $f(t)=t$ we have

\begin{align}
    \mathscr{L}[f(t)] &= \int_{0}^{\infty} f(t) \cdot e^{-s t} \, dt  = \dfrac{1}{s^2}
\end{align}
So, this represents a rather significant finding!  With the Fourier transform, we were limited to functions that decayed away quickly as $|t|$ (or $|x|$ in the notation of the last chapter) became large, so that the Fourier integral remained finite.  Now, we have significantly expanded the potential functions we can transform, as long as they grow slower than exponentially (or even exponentially fast, as long as the resulting integral converges; there will be more discussion on this later).  We have paid a small price for this: Our transform now applies only to the positive portion of the real line.  However, given the correspondence with initial value problems where the positive portion of the real line is the only part that corresponds to physical reality, this suggests that the Laplace transform will be potentially useful for transforming problems in the time domain.

We are now have a clear definition for the Laplace transform.

\begin{definition}[Laplace Transform]\indexme{Laplace transform!definition}
The Laplace transform, $\mathscr{F}[f(t)]=F(s)$, is an integral transform primarily used to transform functions on the half-real line, $t\in[0,\infty)$. It is defined by

\begin{center}
\boxed{
\begin{aligned}
&\\
~~\mathscr{L}[f(t)] &= F(s) = \int_{0}^{\infty} f(t) \cdot  e^{-st}\, dt~~\\
&
\end{aligned}
}
\end{center}
where $s$ is a complex number. We restrict our analysis to functions, $f(t)$, where the following must be true

\begin{enumerate}
    \item $f$ is piecewise continuous on $0\le t \le T$, for all $T>0$.  This means only that $f$ has a finite number of points of discontinuity.
    \item We must have 
    \begin{equation}
   \underset{t\rightarrow\infty}{\lim} \,{f(t)}{e^{-s t}} = 0
\end{equation}
so that the transform is given by a finite integral.  

\item Sometimes, a more practical version of this last statement is given by the following:  There are \emph{real numbers} $M$ and $c$ such that
\begin{equation}
    |f(t)| \le M e^{c t} \textrm{ for all values of } t
\end{equation}
Then, the Laplace transform exists for $Re(s)>c$, which is to say that the \emph{real} part of $s$ (the real number $\alpha$ in $s=\alpha+i \xi$) is greater than $c$.  This, then, means that the decaying exponential defined as part of the Laplace transform \emph{faster} than the increasing exponential $e^{c t}$; the result is a finite integral.
\end{enumerate}
\end{definition}

\begin{svgraybox}
\begin{example}[Computation of the Laplace transform.]

The power of the Laplace transform can start  be seen as we investigate functions that could not be examined with the Fourier transform.  For example, consider the following function that diverges to infinity as $t\rightarrow\infty$

\begin{equation}
    f(t) = a t^5
\end{equation}
Clearly this funtion grows rapidly as $t$ becomes large.  The Laplace transform, however, is still defined because this polynomial grows slower than the decaying exponential that we added to the transform.  Thus, the result is finite. The transform is formally given by

\begin{equation}
    \mathscr{L}[a t] = \int_{0}^{\infty} a t^5  e^{-st}\, dt
\end{equation}
This problem is solvable by a straightforward (and tedious) use of integration by parts five times, reducing $t^5$ by one order of the exponent each time integration by parts is used.  The steps are not shown here; however, the result is

\begin{align}
    \mathscr{F}[at^5] &=a\left( \int_{0}^{\infty} t^5   e^{-st}\, dt \right)  \\
    \mathscr{F}[at^5] &= \dfrac{120 a}{s^6}
\end{align}

\end{example}
\end{svgraybox}

\section{Some Notes About Laplace Transformed Functions}

Unlike the Fourier transform, it is a bit harder to relate the \emph{shape} of a Laplace transformed function to the properties of the original function.  For example, the Fourier transform of a pure sine function is a delta function shifted to the appropriate wave number; intuitively, this makes sense, because the sine funciton contains a single wavelength, and this is represented in its transform.  For a Laplace transform, however, the transform of a sine function is a hyperbolic function in the transform space -- not a result that immediately appeals to an easy physical interpretation.

A second issue is that the transformed functions are a function of $s=\alpha + i \xi$-- a  complex number.  Recall from Chp.~\ref{chaprev} that a complex variable can be though of as a two variables, one representing the real-line component, and one representing the imaginary component.  While it this is also true for the Fourier transform, the transform variable for that case is given by $0+i \xi$  Thus, a Fourier-transformed function can be plotted as a function of only $\xi$.   To properly interpret a Laplace-transformed function, however, we must plot the function as a complex number $F(s)$ $s=\alpha+i \xi = (\alpha,\xi)$; that is, we must plot it as if it were a function in 2-dimensions!

It can be useful to examine the shapes of a few functions and their transforms, which is done in the following example.

\begin{svgraybox}
\begin{example}[A few examples of functions and their transforms.]

Recall from Chapter \ref{chaprev}, a complex number $s=\alpha+i \xi$ can be represented in polar coordinates by specifying its magnitude, $r$ (or absolute value), and its argument, $\theta$.  These are defined by 

\begin{align*}
    r &= \sqrt{\alpha^2+\xi^2} \\
   Arg(s)&= \theta = \tan^{-1}\left( \frac{\xi}{\alpha}\right)
\end{align*}
One complicating feature of the Laplace transform is that the \emph{result}, $\mathscr{L}[f(t)]$ of the transform is the complex function F(s). Recall that a function of a complex number can be thought of as being a function in a two-dimensional plane, where the vertical axis represents the \emph{imaginary} component of $s$, and the horizontal axis represents the \emph{real} part of $s$.  Further complicating things is that the output of the function is, in general, itself a complex number.  Thus, in some very reasonable ways, the Laplace transform is a vector function, that returns a vector as its output, i.e., 

\begin{equation*}
    F(s) = F_r(\alpha+i \xi) +i F_I(\alpha + i \xi)
\end{equation*}
where here $F_r$ represents the real component of $F(s)$ ($F_r = Re[F(s)]$), and $F_I$ represents the imaginary part of $F(s)$ ($F_I = Im[F(s)]$).  This complicates the matter because we now have a result with four dimensions $(\alpha, \xi, F_r, F_I)$!  There are a number of methods for representing such functions, however.  Here, we will map the functions as follows.  First, the horizontal and vertical components will represent the real and imaginary components of $s$ as is customary.  Second, the magnitude (or modulus) of $F(s)$ will then be plotted on the vertical axis.  Finally, the argument ($\theta$) will be overlaid on the surface as a color field to represent the fourth dimension of our plot.  While this is a lot of information to take in, the plots of our transformed functions are still interesting to examine.  A few examples appear below.
\end{example}
\end{svgraybox}
\vspace{-7mm}
\begin{svgraybox}
\begin{enumerate}
\item The function $f(t)=1$ for $t>0$ is has the transform $\mathscr{L}(1)=\tfrac{1}{s}$.   Note that at $s=0+i 0$, the transform function has singularity ($\tfrac{1}{s}$ tends toward infinity there) which in complex analysis is known as a \emph{pole}.  This imposes no particular problems in complex analysis.  There is an analogue to the Stokes theorem (called the \emph{residual theorem}\indexme{complex number!residual theorem} in the complex plane that allows an integral over a complex domain to be replaced by a closed line integral.  Converting to the complex coordinates $s= \alpha +i \xi$ (as we defined above) we can plot the resulting transform in the complex plane.  This plot is given in Fig.~\ref{g1p0}.
  {
\begin{centering}
\includegraphics[scale=.4]{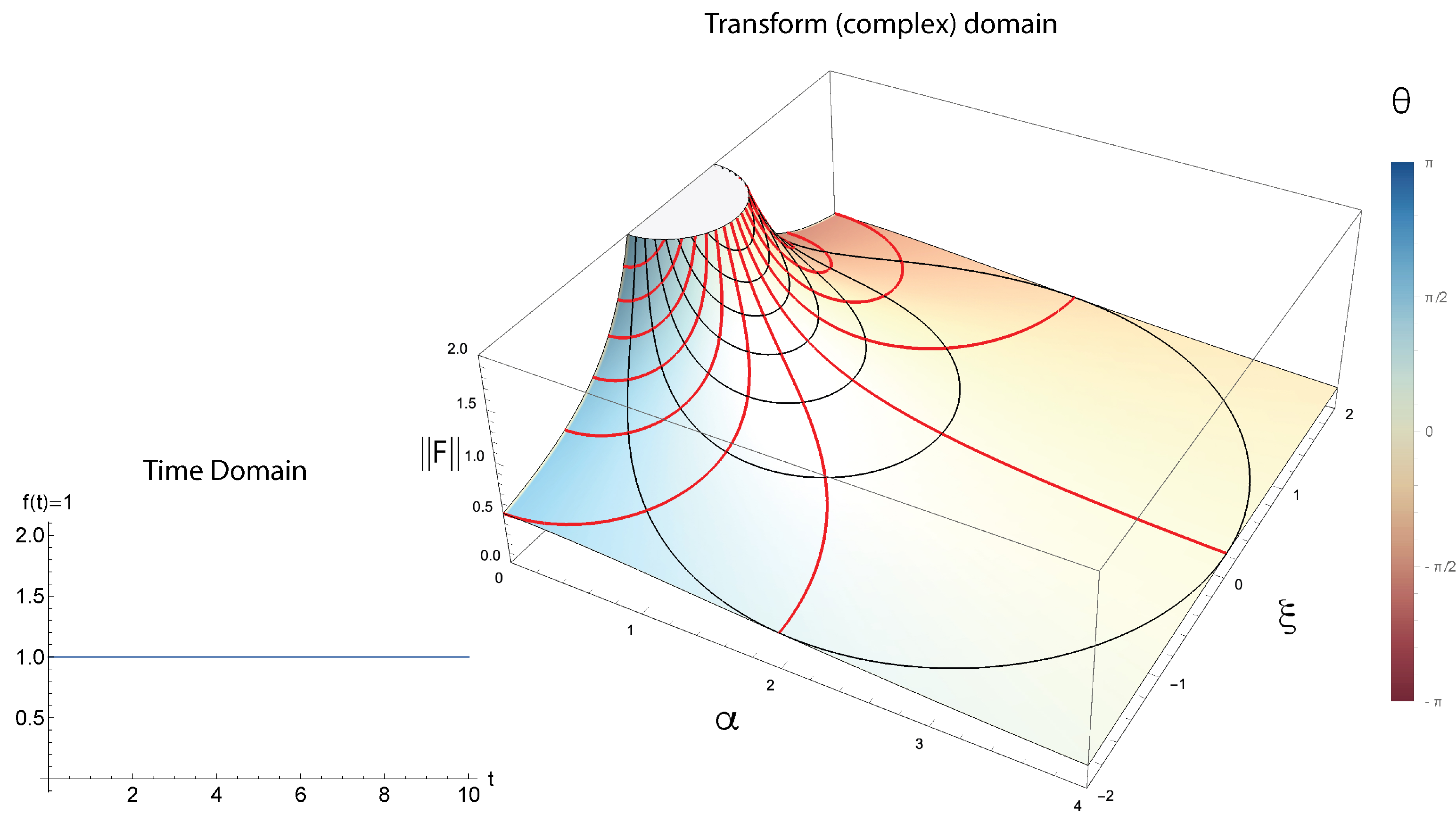}
\end{centering}
\vspace{-2mm}
\captionof{figure}{Laplace transform of $f(t)=1$. The vertical component gives the magnitude ($r$) of the function, and the color indicates the argument, $\theta$.  The black and red lines indicate lines of constant imaginary (black) and real (red) components of $s$.\vspace{3mm}}
\label{g1p0}  
}  
\saveenumerate
\end{enumerate}
\end{svgraybox}
%
%
\vspace{-14mm}
\begin{svgraybox}
\begin{enumerate}
\restoreenumerate
\item  The function $f(t)=sin(t)$ for $t>0$ has the transform $\mathscr{L}[\sin(t)]=\tfrac{1}{1+s^2}$.  Here, there are singularities as $s$ tends towards $s=1+0i$ or $s=-1+0i$.  The plot of the transform is given in complex coordinates $s= \alpha +i \xi$ in Fig.~\ref{g1p2}.

     {
\begin{centering}
\includegraphics[scale=.3]{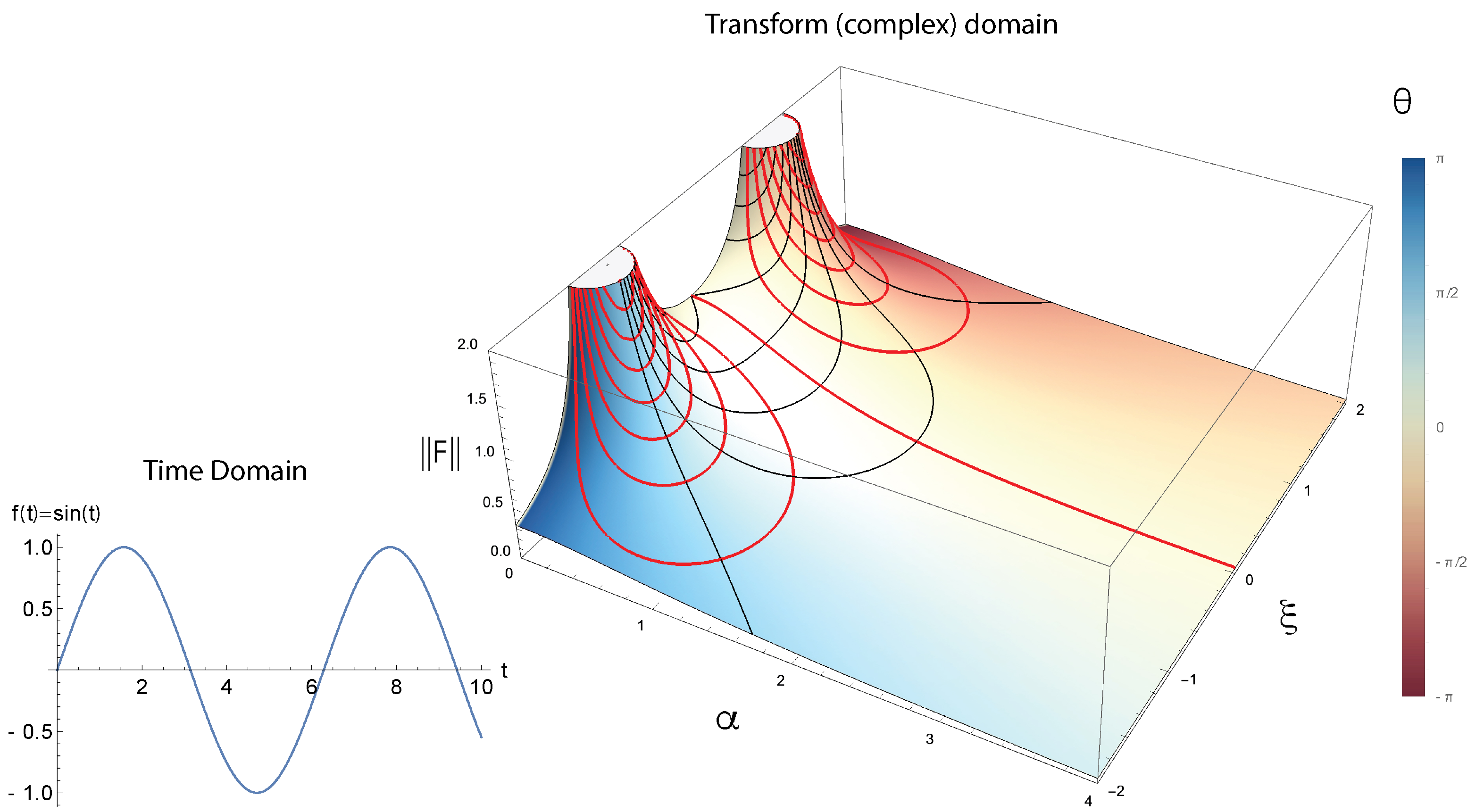}
\end{centering}
\vspace{-4mm}
\captionof{figure}{Laplace transform of $f(t)=\sin(t)$.The vertical component gives the magnitude ($r$) of the function, and the color indicates the argument, $\theta$.  The black and red lines indicate lines of constant imaginary (black) and real (red) components of $s$\vspace{3mm}}
\label{g1p2}  
}  
\saveenumerate
\end{enumerate}
\end{svgraybox}
\vspace{-10mm}
\begin{svgraybox}
\begin{enumerate}
\restoreenumerate
\item  The function $f(t)=e^{2t}$ for $t>0$,  with transform $\mathscr{L}[e^{2t}] = \frac{1}{s-2}$.  Note that for this function, the transform $\mathscr{L}[e^{2t}]$ is not defined for $Re(s)\le 2$.  Recall, this happens because the exponential $e^{-st}$ must decay faster than the exponential $e^{2t}$ grows; thus we must have $Re(s)> 2$.  This result means that there is a discontinuity in the function that spans the entire complex plane along the line $s=2+i \xi$ for all values of $\xi$.  As mentioned above, this does not present significant problems in the inversion integral.  In the complex plane, one can use the residual theorem (analogous to the Stokes theorem of vector calculus) to compute the inverse.

     {
\begin{centering}
\includegraphics[scale=.35]{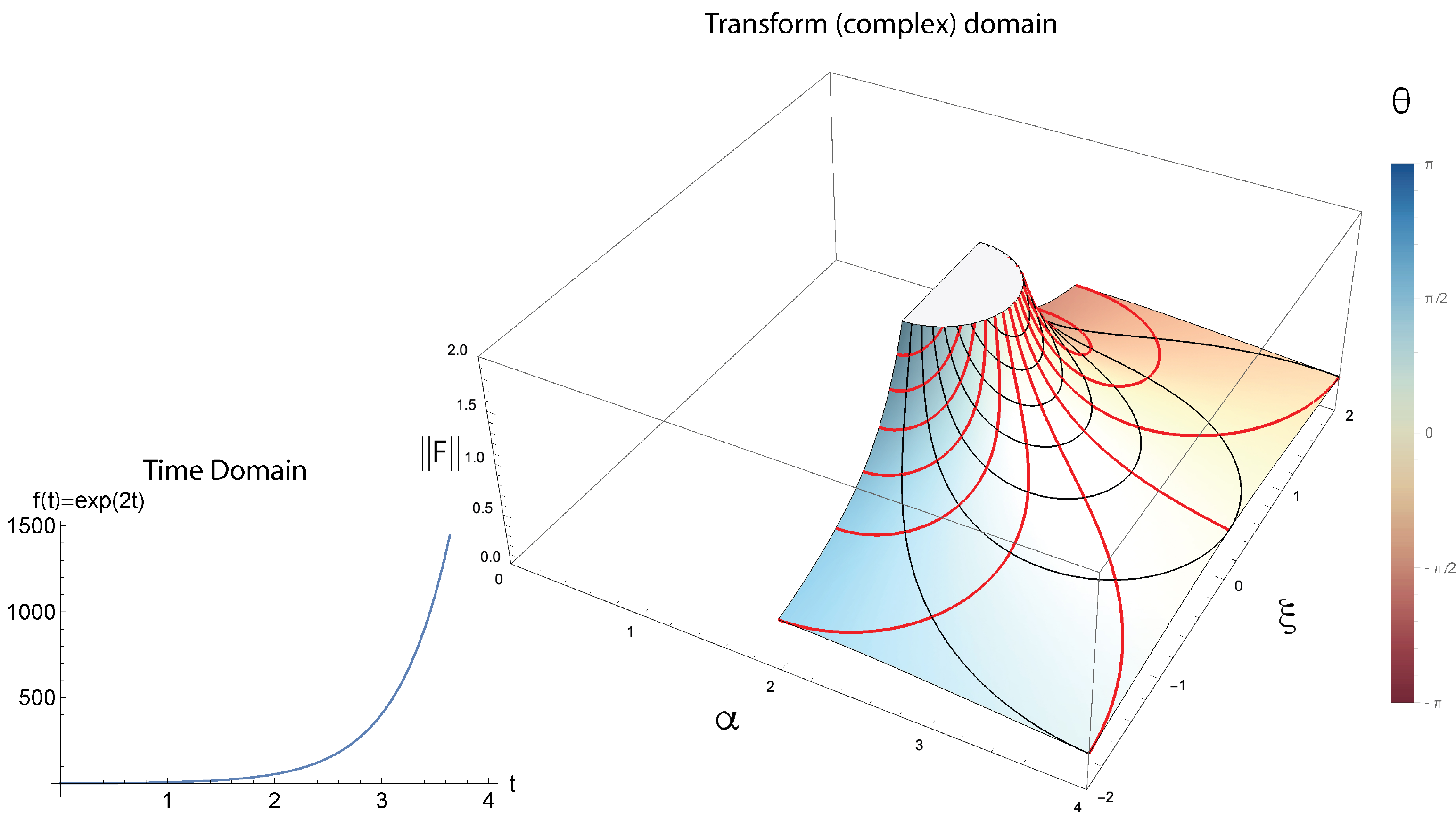}
\end{centering}
\vspace{-2mm}
\captionof{figure}{Laplace transform of $f(t)=e^{2t}$, with transform $\mathscr{L}[e^{2t}] = \frac{1}{s-2}$.  The vertical component gives the magnitude ($r$) of the function, and the color indicates the argument, $\theta$.  The black and red lines indicate lines of constant imaginary (black) and real (red) components of $s$\vspace{3mm}}
\label{g1p}  
}  
\end{enumerate}
\end{svgraybox}

\section{The Inverse Laplace Transform}\indexme{Laplace transform!inverse}

There is a subtlety in the definition of the Laplace transform integral that has a significant effect on the inverse transform.  In analogy with the Fourier transform, we would expect the \emph{inverse} Laplace transform to involve integrating the transformed function with respect to the transform variable, $s$.  On the surface, this is both familiar (by comparison with the Fourier transform), and reasonable.  However, recall that now $s$ is a complex variable.  Depending on how much you have previously studied complex variable analysis, you may know that functions of complex variables can be viewed as if they were functions of two independent variables-- the real variable, and the imaginary variable (this was briefly covered in Chapter \ref{chaprev}).  Thus, integration with respect to a complex variable can be viewed as being an integration in the complex \emph{plane} (where, conventionally, the vertical axis is called the \emph{imaginary} axis, and horizontal axis is the \emph{real} axis.)

\subsection{$^\star$The Bromwich Integral}

The inverse Laplace transform, $\mathscr{L}^{-1}[F(s)]$ can formally be given in terms of what is called the Bromwich integral\indexme{complex number!Bromwich integral}\indexme{Laplace transform!Bromwich integral}.  While we will not pursue this formalism further (it requires an understanding of complex variable analysis), it is useful to think about this integral conceptually because it corresponds to the inverse of the Fourier transform.  The inverse transform in terms of the complex Bromwich integral is given by

\begin{equation}
    \mathscr{L}^{-1}[F(s)] = \underset{\xi \rightarrow\infty}{\lim}\, \frac{1}{2 \pi i} \int_{s=b-i \xi}^{s=b+i \xi} F(s) e^{st}\, ds
\end{equation}
where here the ``lim" has been added to emphasize that it is only the imaginary component of the variable that is evaluated at infinity in the integration step.  Also note that the normalization constant, analogous to the constants appearing in the Fourier transform, appears only in the inverse transform.  Clearly, the net result is the same, so this choice is done only for convenience in the definition of the transform (and creates a slight asymmetry in the definition of the inverse transform).

\subsection{Post's Inversion Theorem}\indexme{Laplace transform! Post's inversion theorem}

The Bromwich integral allows one to compute the \emph{inverse} Laplace transform, but it requires that one also have expertise in complex analysis, which is not a common area of study in more modern instruction in applied mathematics.  

However, there are alternatives.  In 1930, a mathematician named Emil Post \citep{post1930generalized} presented the following Laplace inversion formula as a the limit of a sequence.  Assume that we have the Laplace transform pair $f(t) \rightleftharpoons F(s)$.

\begin{theorem}[Post's inversion theorem]
Suppose we have the Laplace transform pair $f(t) \rightleftharpoons F(s)$.  Then, the inversion of $F(s)$ can be determined by computing the following limit

\begin{equation}
f(t) = \mathop {\lim }\limits_{k \to \infty }  \frac{(-1)^k}{k!} \left( \frac{k}{t}\right)^{k+1} F^{(k)}\left( \frac{k}{t}\right)
\end{equation}
where $F^{(k)}$ is the $k^{th}$ derivative of $F(s)$ with respect to $s$ [N.B., the independent variable is changed after the differentiation process so that $F^{(k)}(s) \rightarrow F^{(k)}(k/t)$].
\end{theorem}

It is somewhat frequently stated that this inversion formula is not of practical value, but this statement should be viewed dubiously.  There are several reasons to view this inversion formula as a useful tool, including the following.

\begin{enumerate}
    \item Post's inversion formula is expressed entirely in terms of real variables, and is thus does not involve complex analysis
    \item There are now symbolic software packages that can be used to help find a recurrence for the derivative  evaluate the limits of terms such functions as the one specified in Post's inversion theorem. 
\end{enumerate}

\begin{svgraybox}
\begin{example}[Use of Posts Inversion Theorem.]

In the example above, we saw that $\mathscr{L}[a t] = \frac{a}{s^2}$.  Show that the inverse transform of $F(s)=1/s^2$ is in fact, $f(t)=t$ using Post's inversion theorem.

To start, we need the derivatives of $F(s)$ with respect to $s$.  A little work allows us to find the general form for this derivative.

\begin{align*}
    F^{0}(s) &= \frac{1}{s^2} \\
      F^{1}(s) &= -2\frac{1}{s^3} \\
        F^{2}(s) &= 6\frac{1}{s^4} \\
          F^{3}(s) &= -24\frac{1}{^5} \\
            F^{4}(s) &= 120\frac{1}{s^6} \\
            &\ldots \\
            F^{(k)}(s) &=(-1)^k (k+1)! \frac{1}{s^{k+2}}
\end{align*}
 Now, making the change $s\rightarrow k/t$, and substituting into Post's inversion formula, we need to evaluate
 
 \begin{equation*}
f(t) = \mathop {\lim }\limits_{k \to \infty }  \frac{(-1)^k}{k!} \left( \frac{k}{t}\right)^{k+1} (k+1)! (-1)^k\frac{1}{\left(\frac{k}{t}\right)^{k+2}}
\end{equation*}
Noting that $(-1)^k(-1)^k = 1$ and  that $(k+1)!/k! = (k+1)$, we have 
 
\begin{equation*}
f(t) = \mathop {\lim }\limits_{k \to \infty }  \frac{(k+1)!}{k!}  \frac{k^{k+1}}{t^{k+1}}  \frac{t^{k+2}}{k^{k+2}}
\end{equation*} 
While this looks somewhat messy, we can simplify a little bit by making the substitution $m = k+1$.  Then we have

\begin{equation*}
f(t) = \mathop {\lim }\limits_{k \to \infty }  \frac{m!}{(m-1)!}  \frac{(m+1)^{m}}{t^{m}}  \frac{t^{m+1}}{(m+1)^{m+1}}
\end{equation*} 
noting that $m!/(m-1)! = m$, then this reduces to

\begin{align*}
f(t)& = \mathop {\lim }\limits_{m \to \infty }  \frac{m}{(m+1)} t\\
&= t
\end{align*} 

\end{example}
\end{svgraybox}
%
Fortunately, most of the Laplace transforms and their inverses have been computed and tabulated previously.  In fact, this is the very thing that makes the transform methods so useful -- much of the effort for computing the transforms and inverse transforms is removed from the user because most of the ``interesting" functions have been computed previously.  In the appendix to this chapter, there is a somewhat lengthy table of transforms and inverse transforms.  Regardless of the presence of such tables, however, it is important to understand the transforms themselves.  This means understanding what the transform is doing when it is computed, why it works, and how one might compute its inverse.  The material above provides one essentially all that one might need in this respect, and, in particular, the Post inversion theorem allows one to compute the inverse Laplace transform without having to have previously studied the theory of complex variable analysis.

\section{The Laplace Transform of the Heaviside and Delta Functions}\indexme{Laplace transform!delta function}\indexme{Laplace transform!Heaviside function}

Laplace transforms are used more widely than Fourier transforms, in part because of their applications to signal processing.  Recall that the Laplace transform is used on the half-line $t\in[0,\infty)$, thus there is a natural zero time built into the transform.  This corresponds to many processes that evolve from initial condition.   For example, the Laplace transform is routinely used to predict the behavior in time of a circuit containing resistors, inductors, and capacitors.  Considering that these are the building blocks of classical (pre-transistor) circuits, solutions to such problems are an important application.  Another example might be the analysis of power grids.  Suppose that one wanted to understand what might happen on a grid if there were a sudden ``spike" load (delta-function like load) or sudden load added continuously at time $t=a$ (a Heaviside load).  These correspond to real physical phenomena that are of interest to those who plan and run power grid infrastructure.  Again, the Lapalce transform would be a potential tool for analysis of these problems.

It turns our that the Heaviside and delta functions are used routinely in modeling many kinds of dynamic systems because of their obvious interpretations as spike or sudden constant sources or sinks.  Unlike the Fourier transform, there is no real difficulty in determining their transforms (once one has adopted the formalisms associated with the delta function).  The transform of the Heaviside function is given by the following

\begin{align}
   \mathscr{L}[ H(t-a)]& = \int_0^\infty H(t-a)\exp(-t s) \, dt \nonumber\\
   &= \int_a^\infty 1\times\exp(-t s) \, dt \nonumber\\
   &= -\left. \dfrac{1}{s} \exp(-t s) \right|_{a}^{\infty} \nonumber\\
   &= \dfrac{1}{s} \exp(-a s) 
\end{align}
And for the delta function

\begin{align}
   \mathscr{L}[ \delta(t-a)]& = \int_0^\infty \delta(t-a)\exp(-t s) \, dt \nonumber\\
   &= \exp(-a s) 
\end{align}
%
\subsubsection*{A Warning Note About the Application of the Delta Function at a Discontinuity}

It is necessary to make an important note about the Laplace transform of the delta function.  The casual use of generalized functions without carefully understanding their deeper theory can sometimes lead to trouble.  One problem that has arisen classically with the Laplace transform is the interpretation of a delta function placed at the origin.  Why might this create a problem?  To start, lets consider the following

\begin{equation}
    \int_0^\infty \delta(t) \exp(-s t) \, dt
\end{equation}
Now, this integral looks innocent enough.  But, recall, we defined the delta function as a sequence of functions, indexed by $n$, that ``converge" in some sense to the delta function as $n\rightarrow\infty$.  We proposed a number of functions and these were usually symmetric around the origin (e.g., One example we used in the chapter on delta functions involved powers of the cosine function, $(cos[\pi t])^n$, which is symmetric around $t=0$).  Suppose we think of our delta function as the limit of a delta sequence $\delta_n(t)$.  Considering this problem further, we might write for the integral above as following

\begin{equation}
    \int_0^\infty \delta(t) \exp(-s t) \, dt = \int_{-\infty}^\infty \left[\lim_{n\to\infty}\delta_n(t)H(t) \exp(-s t)\right] \, dt
    \label{impossibleint}
\end{equation}
If our delta function were indeed a sequence of functions that are symmetric around zero, we might expect such a series to converge to yield an integral of 1/2.  Why would this be the case?  For each member of the delta sequence (for all $n>0$), exactly half of the delta sequence is on $(-\infty,0)$, and the other half is on $(0, \infty)$.  When only the positive portion is taken (which is assured by the Heaviside cutoff), then the result should be 1/2 rather than the more conventional value of 1.  

In fact, the problem is somewhat unsolvable, at least in the context of conventional distribution theory.  Depending upon the delta sequence one chooses, the integral given by Eq.~\eqref{impossibleint} can be made to converge to any number between 0 and 1.  Another way of stating this is that the multiplication of a delta function and a Heaviside function does not have a unique result, and should be avoided (at least by us!)  Rather than attempting to patch up the notion of how to define the delta function applied at $t=0$, we will simply state that this is not an allowable operation in our framework for Laplace transforms.  This would apply anywhere that a delta and Heaviside function might overlap-- for example, the Laplace transform of $\delta(t-a)H(t-a)$ is similarly not well defined.  Thus, we are free to use delta functions, but we cannot apply them anywhere in the domain where there is already a discontinuity!  This does not actually cause us much trouble as the circumstances where this might occur physically are few (although, this problem has been noted in signal processing theory and in applications to physics).  While there are extensions to the theory of distributions that can make sensible interpretations about such multiplications of the delta and Heaviside, for our purposes the best route forward is simply to indicate that these quantities are not defined in our framework.

\section{The Laplace Transform of Functions with a Discontinuity}\indexme{Laplace transform!discontinuous functions}

\begin{figure}[t]
\sidecaption[t]
\centering
\includegraphics[scale=.5]{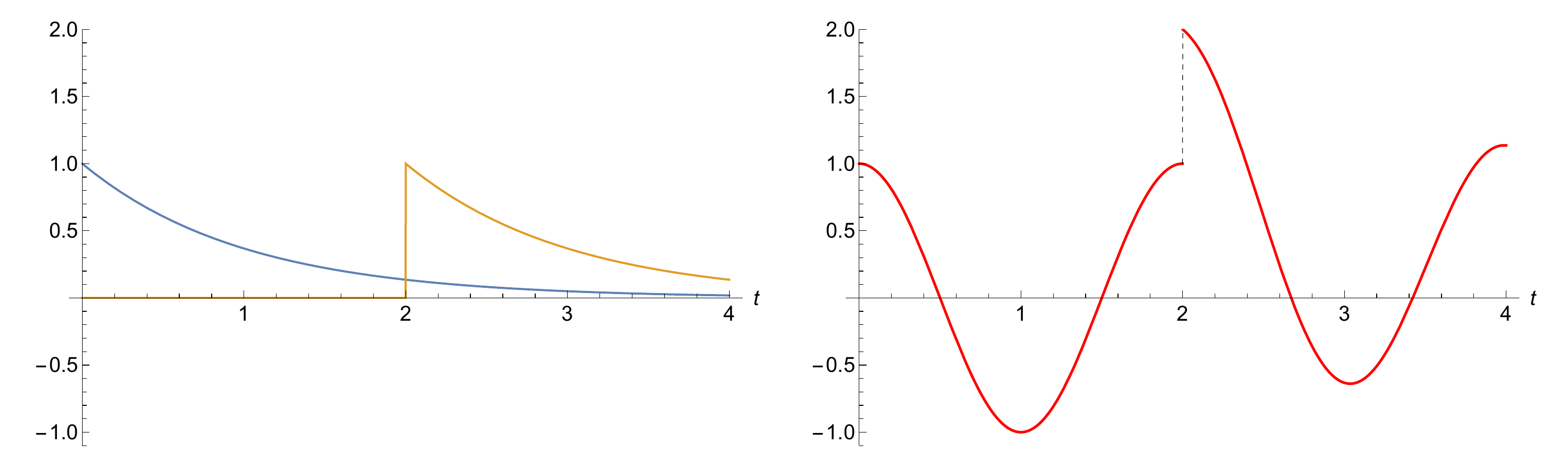}
\caption{(Left) The exponential function $\exp(-t)$ and the shifted and cutoff exponential function $\exp[-(t-2)]H(t-2)$. (Right) The function $w(t) =  \cos(\pi t) + \exp[-(t-2)] H(t-2)$.  The dashed line indicates the jump. }
\label{convolvefig}       
\end{figure}

The Heaviside function, $H(t-a)$ is an example of a function with a single discontinuity located at $t-a$.  Above, we were able to compute the Laplace transform of this discontinuous function, suggesting that more general functions with discontinuities might be be transformable.  In fact, this will prove to be the case.  To start thinking about this problem, let's consider the following example.  Consider the following two functions $f$ and $g$ given by

\begin{align}
f(t) & = \cos(\pi t) \\
g(t)= & \exp(-t)
\end{align}
We can construct a discontinuous function by a number of methods, but let's try the following one.  Between $0 < t < a$, the function will be specified by $f(t)$; after time $t=a$, the function is specified by the sum of $f(t)$ and $g(t-a)$.  Note that this latter function is \emph{shifted} before being added to $f(t)$. Suppose we let this shifted time be given by $t'=t-a$.  Thus the function $g(t')$ is identical to $g(t)$, but it is translated in time forward in time by the amount $a$ (see Fig.~\ref{convolvefig}).   Note that the shift is \emph{forward} in time simply because at $t=a$, $t'=a-a=0$; thus, the value $g(t=0)$ gets mapped in the new coordinate system to the time $g(t'=0)$, or, equivalently, to $g(t-a)$ in the original time coordinates.  Thus, we have the following discontinuous function, which we will call $w(t)$, specified by 

\begin{equation}
    w(t) = \begin{cases}
        \cos(\pi t) & t< a \\
        \cos(\pi t) + \exp[-(t-a)] & t \ge a
    \end{cases}
\end{equation}
This is not the most useful form for us to work with.  We can make an equivalent statement in one equation by using the Heaviside function as follows. 

\begin{equation}
    w(t) =
        \cos(\pi t) + \exp[-(t-a)] H(t-a), ~~~~ ~\textrm{for all }t>0
\end{equation}
To be very clear, this kind of problem is perfectly well defined generally.  However, note that by the warning in the previous section, we cannot consider the function $g(t)=\delta(t)$ because then we would have a delta function multiplying a Heaviside function at the point of discontinuity (which we have disallowed).

At this juncture, we can proceed formally by applying the Laplace transform to both sides of the expression.

\begin{equation}
    W(s) =
       \mathscr{L}[ \cos(\pi t)] + \mathscr{L}[ \exp[-(t-a)] H(t-a)], ~~~~ ~\textrm{for all }t>0
\end{equation}
Where $W(s)$ is the Laplace transform of $w(t)$.  The first of these two functions is simply the Laplace transform of $\cos(\pi t)$.  One can compute this transform by converting to complex exponentials using Euler's identity (i.e., $\cos(\pi t) = \tfrac{1}{2}[\exp(i \pi t)-\exp(-i \pi t)]$), and then integrating.  Much like for the Fourier transform, however, we rarely will have to compute the transform of common functions ourselves, because the transforms are well-known.  The table at the end of this chapter contains a large number of such transforms; in particular, it gives us

\begin{equation}
    \mathscr{L}[ \cos(\pi t)] = \dfrac{s}{\pi^2+s^2}
\end{equation}
Now, what about $\mathscr{L}[ \exp[-(t-a)] H(t-a)]$?  Well, this one is not on the table per se, but it is also an exponential, so it is not to difficult to integrate.  We have

\begin{align}
    \mathscr{L}[ \exp[-(t-a)] H(t-a)]& = \int_{t=0}^\infty \exp[-(t-a)] H(t-a)\exp(-t s)\, dt \nonumber\\
    & = \int_{t=a}^\infty \exp[-(t-a)]\exp(-t s)\, dt
\end{align}
While this integral can be done directly, for the purposes of broader applications we will approach the problem as follows.  First, let $\tau = t-a$ so that $dt = d\tau$ and $t=\tau+a$.  Now note that the integral converts to the following

\begin{align}
    \mathscr{L}[ \exp[-(t-a)] H(t-a)] & = \int_{\tau=0}^\infty \exp[-(\tau)]\exp(-(\tau+a) s)\, d\tau \nonumber \\
    & = \int_{\tau=0}^\infty \exp[-(\tau)]\exp(-\tau s) \exp(-a s)\, d\tau \nonumber \\
    & = \exp(-a s)\int_{\tau=0}^\infty \exp[-\tau]\exp(-\tau s) \, d\tau 
\end{align}
This last integral is just the Laplace transform of the negative exponential multiplied by $\exp(-a s)$.  Note that the name we use for the variable of integration is not material to the transform of the function. Referring to the table of transforms, we find that

\begin{equation}
   \exp(-a s)  \mathscr{L}[\exp(-\tau)] =\exp(-a s) \dfrac{1}{s+1}
\end{equation}
where here we need only note that in the integration, $\tau$ plays exactly the same role as does $t$ (and, in fact, we can simply relabel our independent variable $t$ to prevent any confusion). Our final result can now be given as follows.

\begin{align}
    W(s) &=  \mathscr{L}[ \cos(\pi t)] + \mathscr{L}[ \exp[-(t-a)] H(t-a)] \nonumber\\
    &=  \dfrac{s}{\pi^2+s^2}+\exp(-a s) \dfrac{1}{s+1}
\end{align}
This result can be generalized as follows.

\begin{theorem}[Laplace transform of a shifted and cut off function]
Suppose we have a function $g$, shifted forward in time by the amount $a$, whose value is identically zero for $t<a$ (i.e., cut off for $t<a$).  This function can be represented by  $w(t) = g(t-a)H(t-a)$.  The Laplace transform of $w(t)$ is given by

\begin{equation*}
    W(s) = \exp(-a s) G(s)
\end{equation*}
Where $G(s)=\mathscr{L}[g(t)]$; i.e., $G(s)$ is the Laplace transform of the function in the absence of shifting.
\end{theorem}
One useful way to look at this result is that the transform of a shifted and cutoff function, $g(t-a)H(t-a)$, is represented in the Laplace domain as the Laplace transform of $g(t)$ but with an exponentially-decaying weight of $\exp(-a s)$.    

\begin{svgraybox}
\begin{example}[Laplace transform of a jump.]  Suppose we have a power system with load, $\ell(t)$ equal to a constant value of $\ell_1$ (in appropriate units).  At time $t=a$ a new additional load is added to give a total load of of $\ell_2$ units.  The function defining this situation can be written as follows

\begin{equation}
    \ell(t)=\begin{cases}
    \ell_1(t) & t < a \\
    \ell_1(t)+\ell_2(t-a) & t \ge a
    \end{cases}
\end{equation}
Find the Laplace transform of this function.

\textbf{Solution.}
We can convert this piecewise function using the formalism of the Heaviside as a cutoff filter.  The result is
\begin{equation*}
    \ell(t)=\ell_1(t) +\ell_2(t-a)H(t-a)
\end{equation*}
The Laplace transform is given by

\begin{align*}
 L(t)& =    \mathscr{L}[\ell_1(t)] + \mathscr{L}[\ell_2(t-a)H(t-a)] \\
    & =  L_1(t) + \exp(-a s) L_2(t) 
\end{align*}
Note- even though we do not have explicit formulas for $\ell_1(t)$ and $\ell_2(t)$, we were still able to determine the Laplace transform symbolically.  At this point, we could be given \emph{any} admissible functions (i.e., that have a defined Laplace transform), and the result could be found by simply substituting the appropriate transformed functions.
\end{example}
\end{svgraybox}

\section{Laplace Transform of the Derivative}

One of the most familiar applications of the Laplace transform is for solving first- and second-order ODEs on the real half-line (i.e., on $t\in[0,\infty)$).  However, they are  also \emph{very useful} for solving PDEs.  In both cases, the real power of the Laplace transform arises because of what it does to derivatives.  Much like the Fourier transform, the Laplace transform changes derivatives in real space into algebraic quantities in transformed space.  However, unlike the Fourier transform, the Laplace transform comes with a tiny bit of overhead.  Because the integration defining the transform has a lower bound of zero, the integration leads to capturing the initial conditions of the problem when the integration is completed.   Thus, we have the following identities for the Laplace transform of derivatives.  

\begin{theorem}[Laplace trasnform of derivatives]
For any function $f(t)$, where $f(t)$, $f'(t)$ and $f''(t)$ are continuous, then the following transformations are valid.
\begin{align}
    \mathscr{L}\left[\frac{d f(t)}{d t} \right]& =\mathscr{L}[f'(t)]= s F(s) -f(0) \\
    \mathscr{L}\left[\frac{d^2 f(t)}{d t^2} \right]& =\mathscr{L}[f''(t)]= s^2 F(s) -sf(0)-f'(0)
\end{align}
\end{theorem}
\noindent where here $f(0)$ and $f'(0)$ represent known ancillary conditions for $f$.  The proof for this can be done by a simple application of integration by parts, and will be left for an exercise.

There are a few notes to make here.  First, higher-order derivatives can be similarly defined if needed (e.g., see \citet{spiegel1965laplace}).  Second, for the Laplace transform of derivative of order $n$ to exist, we need the function $f$ to be differentiable $n$ times such that the resulting derivative is a continuous function.  However, $n^{th}$ derivative need only be \emph{sectionally} continuous (i.e., the $n^{th}$ derivative function might have non-differentiable points).  Finally, even the condition of differentiability can be relaxed somewhat if needed.  The following theorem defines the Laplace transform for functions that fail to be continuous.

\begin{theorem}
Suppose $w(t)$ is not continuous, but instead it contains a jump discontinuity at $t=a$, and is represented by $w(t) = f(t) + g(t-a)H(t-a)$.  Then, the expression for the derivative is modified as follows
\begin{align}
    \mathscr{L}\left[\frac{d w(t)}{d t} \right]& =\mathscr{L}[w'(t)]=  s F(s)- w(0)+s \exp(-s a) G(s).
\end{align}
Note that for the case where $a=0$, this expression holds with $w(t)=f(t)+g(t)$ by definition.
\end{theorem}
\noindent\textbf{Proof.}\\
By definition, we have

\begin{align}
    \mathscr{L}\left[\frac{d w(t)}{d t} \right] = \mathscr{L}\left[\frac{d f(t)}{d t} \right]+ \mathscr{L}\left[\frac{d }{d t}\left(g(t-a)H(t-a) \right)\right]
\end{align}
Here, $f(t)$ is a regular, continuous function ($w(t)$ is discontinuous).  Thus,  $\mathscr{L}[f]=s F(s)- f(0)$.  Now the trick is to determine the second of the two transforms on the right-hand side; to be explicit, this is given by

\begin{equation}
    \mathscr{L}\left[\frac{\partial }{\partial t}\left(g(t-a)H(t-a) \right)\right] = \int_0^\infty \frac{d }{d t}\left(g(t-a)H(t-a) \right) \exp(-s t) \, dt   
\end{equation}
Here, we can use integration by parts.  Suppose we integrate by parts choosing

\begin{align*}
    u=\exp(-st) && dv= \dfrac{d}{dt}[g(t-a)H(t-a)] \\
    du = -s \exp(-st) && v= g(t-a)H(t-a)
\end{align*}
Integration by parts gives us

\begin{align}
    \mathscr{L}\left[\frac{\partial }{\partial t}\left(g(t-a)H(t-a) \right)\right] 
    &= \left. \exp(-s t) g(t-a)H(t-a)\right|_0^\infty - (-s) \int_0^\infty g(t-a)H(t-a)\exp(-s t) \, dt \nonumber \\
    &= 0 + s\int_0^\infty g(t-a)H(t-a)\exp(-s t) \, dt \nonumber \\
    &= s \exp(-s a) \int_0^\infty g(t)H(t)\exp(-s t) \, dt ~~\textrm{(using the shifting property)}\nonumber \\
    &= s \exp(-s a) \int_0^\infty g(t)\exp(-s t) \, dt \nonumber \\
    &= s \exp(-s a) G(s) 
\end{align}
Putting this all together, we obtain 

\begin{equation}
    W(s) = s F(s)- f(0)+s \exp(-s a) G(s)
\end{equation}
\qed
\begin{svgraybox}
\begin{example}[Laplace transform of derivatives.]  Find the Laplace transform of the following 
\begin{enumerate}
    \item Find $\mathscr{L}[f'(t)]$ for $f(t) = t^2+f_0$. The ancillary condition is $f(0) = f_0$, where $f_0$ is some constant.
    \item Find $\mathscr{L}[f''(t)]$ for $f(t) = \sin(t)$ if the ancillary conditions are $f(0) = 5$, and $f'(0)=2$.  
    \item Find Find $\mathscr{L}[w'(t)]$ for $w(t) =(t-5)^2 H(t-5)$ if the ancillary condition is $w(0) = 0$.
\end{enumerate}
\textbf{Solutions}.

\begin{enumerate}
    \item Recall, for this case we have the theorem $[f'(t)]= s F(s) -f(0)$; thus, we need the transform of $f$ and the value of $f(0)$.  The transform is given by $\mathscr{L}[t^2+f_0]=2/s^3+f_0/s$.  Note that from the expression for $f$, we have $f(0)=f_0$.  According to the theorem, then, we have 
    \begin{align*}
        \mathscr{L}[f'(t)]&= s(2/s^3+ f_0/s) -f_0 \\
        &= 2/s^2 +f_0-f_0 \\
        &=2/s^2
    \end{align*}
 
\item According to  the theorem above, the Laplace transform of the second derivative is given by $\mathscr{L}[f''(t)]= s^2 F(s) -sf(0)-f'(0)$.  Thus, we need the transform of $f$, and the values of $f(0)$ and $f'(0)$.  Collecting these, we have
\begin{align*}
\mathscr{L}[f(t)]&=\mathscr{L}[\sin(t)] = \dfrac{1}{1+s^2} \\
f(0)&=\sin(0)=0 \\
f'(0)&=\cos(0)=1
\end{align*}
Using the theorem, the result is
\begin{equation*}
\mathscr{L}[f''(t)]= s^2\dfrac{1}{1+s^2}-1
\end{equation*}
Or, with a little algebra, this can be put in the form
\begin{equation*}
\mathscr{L}[f''(t)]= -\dfrac{1}{1+s^2}
\end{equation*}

\item  The theorem for the transform of the derivative of discontinuous functions requires that we identify $w(t)=f(t)+g(t-a)H(t-a)$.  Here, it is clear that $g(t-a)H(t-a)=(t-a)^2H(t-a)$, but what should we do about the function $f(t)$?  One solution is to take $f(t)=0$.  Then, we find

\begin{align*}
   W(s) &= s\times 0- 0+s \exp(-5 s) \dfrac{2}{s^3}\\
   &=  \dfrac{2 \exp(-5 s)}{s^2}
\end{align*}

\end{enumerate}

\end{example}
\end{svgraybox}

\subsection{Convolutions in Time}

In the chapter on Fourier transforms, the convolution (which we ordinarily think of as exising in space) was defined by the integral of a product of two functions where one function was shifted relative to the other.  The following useful properties of the convolution were shown to hold true: (1) the \emph{transform} of the convolution of two functions was equal to the product of the transforms of the functions, and (2) the \emph{inverse transform} of the product of two transformed functions was equal to the convolution of the functions (in their non-transformed state).  These same properties hold for convolutions in the Laplace transform.  One major difference is that the convolutions are defined only on the half-line for the Laplace transform.  To be specific, the convolution of two functions $f$ and $g$ is given by

\begin{equation}
    (f*g)(t) =\int_{0}^t f(t-\tau)g(\tau) \, d\tau
\end{equation}

\begin{figure}[t]
\sidecaption[t]
\centering
\includegraphics[scale=.7]{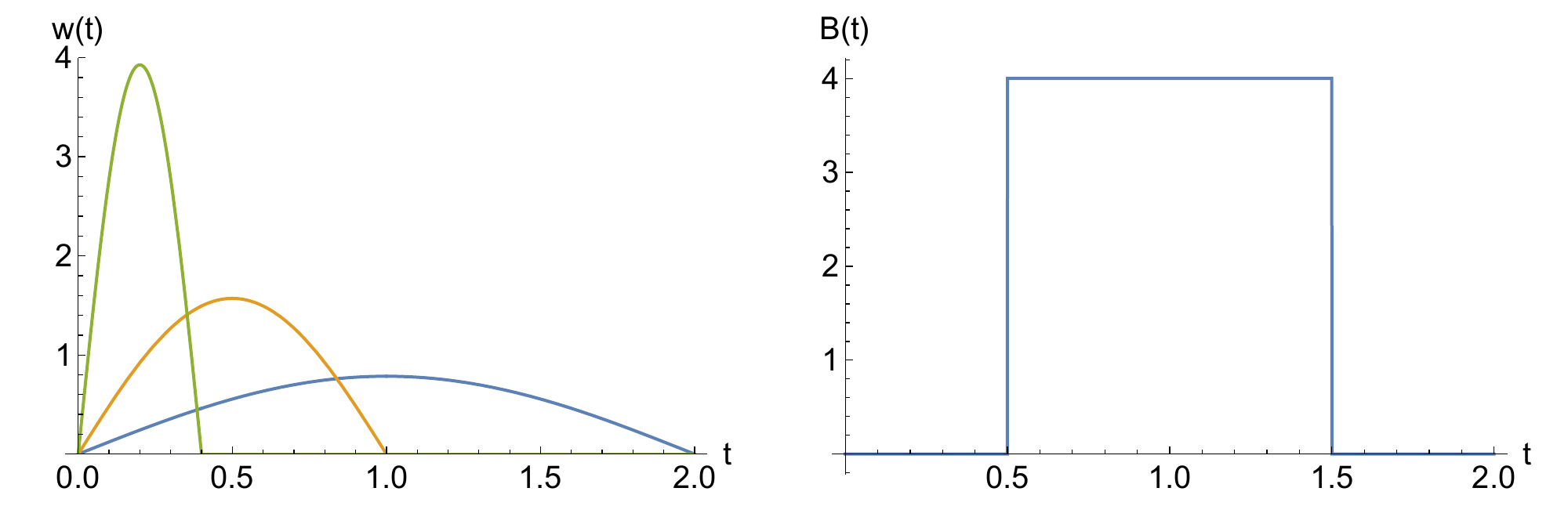}
\caption{(Left) The weighting function $w(t)$ with parameters $a=1$ (blue), $a=2$ (orange) and $a=5$ (green). (Right) The boxcar function $B(t;t_1,t_2)$ with parameters $h_0 = 4$, $t_1=\tfrac{1}{2}$, and $t_2=\tfrac{3}{2}$.}
\label{wandB}       
\end{figure}

To see the effect of the ``one sidedness" of the transform on the resulting convolution, we can repeat an analysis similar to the one we computed as an example of convolutions with the Fourier transform.  Here, again we consider two compact functions. First, we have a \emph{compact} but continuous function defined by 

\begin{equation}
    w(t) = \frac{a \pi}{4}  \left[H(t)-H\left(t-\tfrac{2}{a}\right)\right] \cos
   \left(\tfrac{\pi}{2}   \left(at-1\right)\right)
\end{equation}
where $a$ is any positive number.  Recall that a compact function is nonzero on some finite interval, and zero everywhere else.  For this function, it is compact on the interval $[0,2/a]$.  It is also a technically \emph{probability density} function or normalized \emph{weighting function} because the area under the curve is $1$ for all values of $a$.  The second function we will look at is the \emph{boxcar} function, which is the difference between two Heaviside functions.   It is defined as it was in the chapter on Fourier transforms.

\begin{equation}
    B(t;t_1,t_2) = h_0 [H(t-t_1)-H(t-t_2)]
\end{equation}
where here $t_1$, $t_2$, and $h_0$ are parameters.  The first two are the starting and ending points for the function interval (that is, the function is non-zero in $t\in[t_1,t_2]$), and $h_0$ sets the height of the function.  
Plots of $w(t)$ and $B(t)$ are given in Fig.~\ref{wandB}.

\begin{figure}[t]
\sidecaption[t]
\centering
\includegraphics[scale=.7]{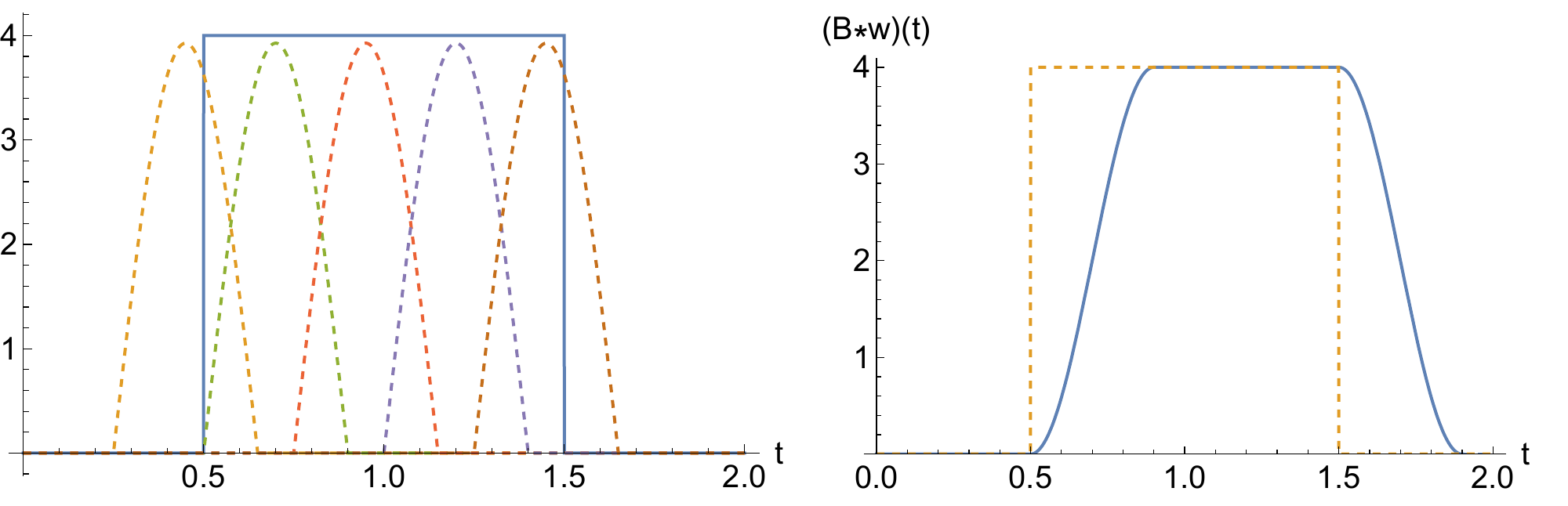}
\caption{(Left) The boxcar function with five copies of the weighting function superimposed at discrete positions.  The convolution computes the integral of the multiplication of these two functions as they are \emph{continuously} overlapped over the bounds of the integral.  (Right) The boxcar function $B(t;\tfrac{1}{2},\tfrac{3}{2})$ with parameters $h_0 = 4$ after convolution with the weighting function $w(t)$ with $a=5$.  The result is the smoothed function shown in blue.  The original function is shown in the dashed orange line for comparison. }
\label{convolvefig2}       
\end{figure}

The convolution defined on the half-line looks different than the convolution defined on the whole line (see Fig.~\ref{densitybox2} for comparison), even thought the two functions used in these two examples have the same shape.  The difference is that, for the Fourier transform, the weighting function $w$ was centered around zero.  For the Laplace transform, this is not possible because the transform is defined for only $t>0$.  Thus, the weighting function must be shifted forward in time by one-half the width of $w(t)$ so that it is defined over a positive domain.

The practical result of this is that the convolution behaves qualitatively differently than for the zero-centered weighting function used in the previous example given in the Chapter on Fourier transforms (Fig.~\ref{densitybox2}).  If one examines Fig.~\ref{convolvefig2}, it is possible to see this difference.  Note that the \emph{shape} of the resulting function is the same as for the previous example, but instead of the center of mass of the convolution being conserved (i.e., the same before and after the operation), it has moved to the right.  This results exclusively because the weighting function has itself been moved to the right so that it is defined on the domain $t>0$.

These differences are not merely an interesting side note.  Convolutions in time are often used to represent the \emph{response} of a system to some signal or process in time.  Thus, in Fig.~\ref{convolvefig2} we would not expect the convolution to allow a \emph{response} to the weighting function at times earlier than $t=1/2$ which is the time that the the boxcar function begins.  Allowing the convolution to do so would, in essence, allow the weighting function to weight components of the boxcar function that had not arrived yet!  This feature of time convolutions is sometimes known as \emph{obeying causality} or as the convolution being \emph{non-anticipative}.  In short, it means simply that time convolutions do not allow their weighting functions to act on portions of the signal (in this example, the boxcar function) that have not yet arrived in time.  No output can precede the input that led to it.

\section{Convolutions and the Laplace Transform}

Laplace transforms, like Fourier transforms, have nice properties when transforming convolution integrals.  There are two ways that the convolution integrals can arise.  First, when attempting to invert a Laplace transform, the result might be expressible as a convolution in the time variable.  Second, one might find that they have an integro-differential equation in time which involves a convolution; the Laplace transform allows consideration of such complex problems.  The following two theorems are stated without proof.

\begin{theorem}[Inverse convolution theorem.]
Suppose $\mathscr{L}[f(t)]=F(s)$ and $\mathscr{L}[g(t)] = G(s)$.  Let $w(t)=\mathscr{L}^{-1}[F(s)G(s)]$.  Then we have 
\begin{equation}
    w(t) = \int_0^t f(t-\tau)g(\tau) \, d\tau = \int_0^t f(\tau)g(t-\tau)\, d\tau
\end{equation}
\end{theorem}

\begin{svgraybox}
\begin{example}[Inverse convolution example]
Suppose we were working on a problem using Laplace transforms, and came to the following result
\begin{equation*}
    W(s) = \dfrac{a}{s^2(a^2+s^2)}
    \label{invertmelap}
\end{equation*}
Inspecting the table of Laplace transforms, we see that there is no such transform listed. However, if we let $F(s)=1/s^2$ and $G(s)=a/(a^2+s^2)$, we can note that each of these transforms are listed on the table.  We have 
\begin{align*}
    F(s)&=\frac{1}{s^2} \rightarrow f(t) = t \\
    G(s)&=\frac{a}{a^2+s^2} \rightarrow g(t) = \sin(a t) 
\end{align*}
The inversion of Eq.~\eqref{invertmelap} is then given by
\begin{equation*}
    w(t) = \int_0^t \tau \sin[a(t-\tau)]\, d\tau
\end{equation*}
Such convolutions are not always easy to solve.  In this case, however, the integral can be done by parts.  The result is 

\begin{equation*}
    w(t) = \dfrac{1}{a^2}[a t - \sin(at)]
\end{equation*}
Technically, we have discovered a new Laplace transform pair this way!  We could put the following pair in the table if we were so inclined.
\begin{equation*}
    \dfrac{1}{a^2}[a t - \sin(at)] \xrightarrow{\mathscr{L}} \dfrac{a}{s^2(a^2+s^2)}
\end{equation*}
\end{example}
\end{svgraybox}

The second convolution theorem allows us to convert convolutions in real space to multiplications in transform space.  

\begin{theorem}[Forward convolution theorem.]
Suppose $\mathscr{L}[f(t)]=F(s)$ and $\mathscr{L}[g(t)] = G(s)$.  Let $w(t)$ be given by the convolution
\begin{equation}
    w(t) = \int_0^t f(t-\tau)g(\tau) \, d\tau = \int_0^t f(\tau)g(t-\tau)\, d\tau
\end{equation}

Then 
\begin{equation}
    \mathscr{L}[w(t)]= W(s) = F(s)G(s)
\end{equation}
\end{theorem}

    \begin{svgraybox}
\begin{example}[Laplace transform of a convolution]

Integro-differential equations are notoriously difficult to solve.  The Laplace transform provides one method for potential solution.  Take the following example

\begin{align*}
    \frac{du}{dt}& = \int_0^t \exp(t-\tau) u(\tau) \, d\tau \\
    u(0)&=1
\end{align*}
For this problem, if we take the Laplace transform of both sides, we find
\begin{equation*}
    s U(s) - 1 = \dfrac{1}{s-1} U(s)
\end{equation*}
Solving this for $U(s)$ gives

\begin{equation*}
    U(s) = \dfrac{s-1}{s(s-1)-1}
\end{equation*}
Using the table of transforms, we find the solution

\begin{equation*}
    u(t) = \tfrac{1}{5} e^{t/2} \Bigg[5 \cosh \left(\tfrac{\sqrt{5} t}{2}\right)
-\sqrt{5} \sinh \left(\tfrac{\sqrt{5} t}{2}\right)\Bigg]
\end{equation*}
We can easily check this solution.  While it is tedious to work out, the derivative of the solution is

\begin{equation*}
    \dfrac{du}{dt} = \frac{2 e^{t/2} \sinh \left(\frac{\sqrt{5} t}{2}\right)}{\sqrt{5}}
\end{equation*}

We can also compute, in this case, the integral

\begin{align*}
    \dfrac{du}{dt} &= \int_0^t \tfrac{1}{5} e^{\tau/2} \Bigg[5 \cosh \left(\tfrac{\sqrt{5} \tau}{2}\right)
-\sqrt{5} \sinh \left(\tfrac{\sqrt{5} \tau}{2}\right)\Bigg]\exp(t-\tau)\, d\tau \\
&= \frac{2 e^{t/2} \sinh \left(\frac{\sqrt{5} t}{2}\right)}{\sqrt{5}}
\end{align*}
So, the solution checks out.  A plot of the solution is provided below.

 \vspace{2mm}
{
\centering
\includegraphics[scale=.75]{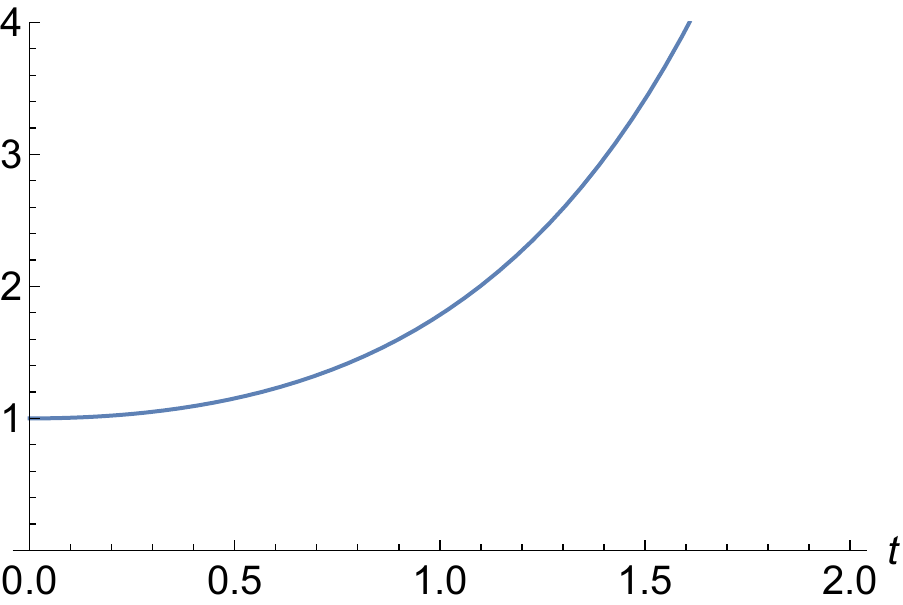}
\vspace{-2mm}
\captionof{figure}{The solution to the integro-differential equation specified in this example.}  
} 

\end{example}
    \end{svgraybox}

Note that this feature of Laplace transforms allows us to transform general integrals of time.  The trick here is to realize that a generic integral in time is just a convolution with the function $\mu(t) = 1$, $t\in(-\infty,\infty)$.  Note that \emph{this function is distinct from the Heaviside function} because it contains no jumps anywhere in the domain.  Additional note that all translations of this function are identical-- they are all equal to for all values of $t$; thus $\mu(t) = \mu(t-a)$ for all values of $a$.  Now, suppose $f(t)$ is a function that has a well-defined Laplace transform.  Then, consider the following integral.

\begin{equation}
    w(t) = \int_{\tau=0}^{\tau=t} f(\tau) \, d\tau
\end{equation}
where $\tau$ is the variable of integration, and we have been more explicit about the variable of integration on the integral bounds.  Noting that we have defined the unit function $\mu(t) = 1$, then we can write this function as the convolution

\begin{align}
     w(t) &= \int_{\tau=0}^{\tau=t} 1\times f(\tau) \, d\tau \nonumber\\
     &=  \int_{\tau=0}^{\tau=t} \mu(t-\tau) f(\tau) \, d\tau
\end{align}
Noting that $\mathscr{L}[\mu(t)]=\mathscr{L}[1]=\tfrac{1}{s}$, then the result for \emph{any} integral of time is the transform

\begin{equation}
    \mathscr{L}\left[ \int_{\tau=0}^{\tau=t} f(\tau) \, d\tau\right]
    = \dfrac{1}{s} F(s)
\end{equation}
\section{Solutions to ODEs Using Laplace Transforms}

One of the primary uses of Laplace transforms is to solve linear ordinary differential equations.   There are a few advantages to finding solutions this way.  First, complicated problems, such as nonhomogeneous equations or integro-differential equations, can be handled with Laplace transforms in a fairly straightforward manner.  The primary difficulty that arises in solving ODEs this way is to determine the appropriate \emph{inverse} transform once the problem is completed.

There is not much more for us to learn about the applications to ODEs.  Note that by convention, Laplace transforms are used to transform the \emph{time} variable.  However, the primary feature defining the Laplace transform is not the name of the independent variable (i.e., time, indicated by $t$), but the fact that the transform is defined on the half-line.  So, any ODE defined on the half-line is amenable to Laplace transform, regardless of whether the independent variable is time or space (or something else).  A few examples will help solidify the methods for solving linear ODEs using Laplace transforms.

\begin{svgraybox}
    \begin{example}[A steady-state diffusion-reaction problem]
As mentioned above, the application of the Laplace transform does not require that the independent variable be time, but, rather, that the independent variable be defined on the half-line.  The following is an example where we look at the solution to a steady-state reaction-diffusion problem defined on $x\in[0,\infty)$.  Recall, the heat/diffusion equation with a first-order reaction at steady state is, by definition, no longer a function of time (i.e., by definition, $\partial u/\partial t \equiv 0$ for steady state to exist).  Suppose we have a steady reaction-diffusion problem on the half real line as follows

\begin{align*}
    D\frac{d^2 u}{dx^2} - k u &= 0, ~~x\in[0,\infty)\\
    u(x=0) &= u_0 \\
   \left. \frac{d u}{dx}\right|_{(x\rightarrow \infty)}& = 0
\end{align*}
For simplicity, let $\alpha = k/D$.  Noting that the domain for space is $x\in[0,\infty)$, the Laplace transform is the correct one to use.  We have one bit of trouble here, however.  We are not given the derivative of the function at $x=0$; instead, we are given a \emph{requirement} that must be respected as $x$ grows large.  We can still proceed, however.  Assume for now that we specify

\begin{equation*}
    \left.\frac{d u}{dx}\right|_(x=0) = c_1
\end{equation*}
where $c_1$ is some currently unknown constant.  Then, we can proceed as follows.
Transforming the ODE yields

\begin{align*}
    [s^2 U(s)-s u_0-c_1] -\alpha U(s) &= 0\\
\end{align*}
Note that the Laplace transform \emph{always} automatically incorporates the ancillary information into the solution.  This is easy to solve, with the result being

\begin{align*}
   U(s) &= \dfrac{s u_0 +c_1}{s^2-\alpha}\\
\end{align*}
Referring to the Laplace transform tables, we find that the inverse transform is given by 

\begin{align*}
   \mathscr{L}^{-1}[U(s)] &= u_0\mathscr{L}^{-1}\left[ \dfrac{s }{s^2-\alpha}\right]\\
   &=u_0 \cosh[\sqrt{\alpha} x]
   +c_1\mathscr{L}^{-1}\left[ \dfrac{1 }{s^2-\alpha}\right]\\
   &=u_0 \cosh[\sqrt{\alpha} x]+ c_1  \sinh[\sqrt{\alpha} x]\\
   &=u_0 \cosh\left[\sqrt{\dfrac{k}{D}} x\right] +c_1 \sinh \left[\sqrt{\dfrac{k}{D}} x\right]\\
\end{align*}
To evaluate the constant, we can take derivatives and then specify that these must be zero.  Alternatively, we note $\cosh(x)/sinh(x) \rightarrow 1$ as $x \rightarrow \infty$.  Therefore, if we set $c_1 = -1$, then our ancillary condition is met.  If we set $k/D = 0.1$ and $u_0=1$, the the solution is given by the plot shown below.
{
\vspace{2mm}
\centering
\includegraphics[scale=.85]{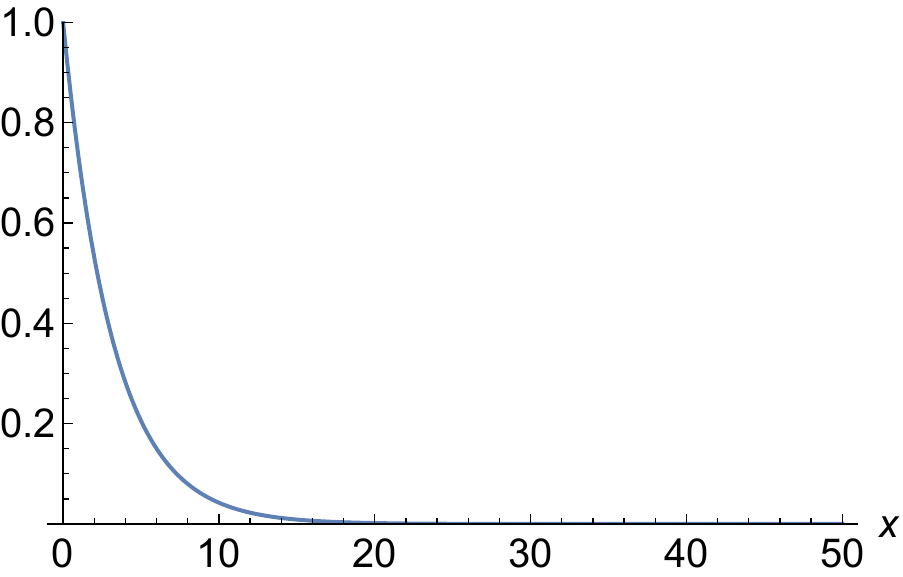}
\vspace{-2mm}
\captionof{figure}{The solution for the oscillating mass problem described above.  Here, the vertical axis defines the displacement relative to the equilibrium (zero) position.}  
}
    \end{example}
\end{svgraybox}

\begin{svgraybox}
    \begin{example}[Solution of a second-order ODE using Laplace transforms]
        The following problem arises in mechanics when analyzing the motion of a dynamic mass-spring system with friction that damps the motion.  In this expression, the variable $y(t)$ represents the displacement from the equilibrium position, and it can be positive or negative.

        \begin{equation}
           \underbrace{m y''(t)}_{mass\times accelleration} =~~  -\hspace{-3mm}\underbrace{\beta_0 y'(t)}_{friction~damping} - \hspace{-0mm}\underbrace{k_0 y(t)}_{spring~force} 
        \end{equation}
        If the mass is also subjected to an external force $F(t)$ that drives the motion (e.g., one could imagine a magnetic driver for a ferromagnetic mass), then this forcing function is added to the right-hand side, yielding a nonhomogeneous ODE.  Suppose that the system is driven periodically by an oscillating and decaying electromagnetic force, $f(t) =\exp(-t)\sin(t)$.   The equation, and its initial position and velocity could be given by an expression of the form

        \begin{align*}
            y''(t) + \beta y'(t) + k y(t)& = \exp(-t)\sin(t) \\
            y'(0)&= 0\\
            y(0)&= 0
        \end{align*}
        where here, to simplify, $\beta = \beta_0/m$ and $k=k_0/m$.  Ordinarily, this would pose a very difficult problem to solve.  Using Laplace transforms will make this somewhat easier, but we will still encounter some challenges in inverting the solution.  For this problem, assume that $\beta=2$ and that $k=2$.  Proceeding, it is fairly direct to find the transform of the ODE.  

        \begin{align*}
           \left[ s^2 Y(s) - sy(0) -y'(0)\right] +  \beta \left[ s Y(s) - y(0)\right]+ k Y(s) = \dfrac{1}{1+(1+s)^2}
        \end{align*}
        Using the ancillary conditions and the values for $\beta$ and $k$, we have
        \begin{equation}
            (s^2+  2 s + 2)Y(s) = \dfrac{1}{1+(1+s)^2}
        \end{equation}
        Note that the left-hand side can be written as
            \begin{equation}
        [(s+1)^2+1]Y(s) = \dfrac{1}{1+(1+s)^2}
        \end{equation}
        Solving for $Y(s)$ gives
        \begin{equation}
            Y(s) = \dfrac{1}{1+(1+s)^2}\dfrac{1}{1+(1+s)^2}
        \end{equation}
        The good news is that we have a solution in Laplace space.  However, we still need to invert this back to functions of time.  Here, note that we know the inverse transform of $ \mathscr{L}^{-1}\left[1/(1+(1+s)^2)\right]= \exp(-t)\sin(t)$ from the tables of Laplace transforms.  In this case, we have the multiplication of the transform by itself.  Setting $F(s) =1/[1+(1+s)^2]$ and $G(s)=1/[1+(1+s)^2]$, from the convolution theorem we have

        \begin{align*}
            \mathscr{L}^{-1}[Y(s)]&= \mathscr{L}^{-1}\left[\dfrac{1}{1+(1+s)^2}\dfrac{1}{1+(1+s)^2}\right] \\
            y(t) &= \int_0^t \exp(-\tau)\sin(\tau)\exp[-(t-\tau)]\sin(t-\tau)\,d\tau
        \end{align*}
        This integral can be computed by using Euler's identity to convert the sine functions to exponentials.  The integral can then be computed (it is the integral of an exponential function), and Euler's identity can be used in the inverse fashion to recover sine and cosine functions from the result.  The computations are algebraically complicated, but otherwise straightforward.  The result is

        \begin{equation}
            y(t) = \dfrac{1}{2} \exp(-t)\left[ -t \cos(t) +
            \sin(t)\right]
        \end{equation}
        As always, it is a good idea to check our answers.  Toward that goal, we have
        \begin{align*}
            y'(t) &= \dfrac{1}{2} \exp(-t) t \sin (t)-\dfrac{1}{2} \exp(-t) (\sin (t)-t \cos (t))\\
            y''(t) &= \dfrac{1}{2} \exp(-t) \sin (t)-\exp(-t) t \sin (t)+\dfrac{1}{2} \exp(-t) t \cos (t)+\dfrac{1}{2}
   e^{-t} (\sin (t)-t \cos (t))
        \end{align*}
        Computing the left-hand side sum $y''(t)+2 y'(t)+2 y(t)$ gives the result $\exp(-t)\sin(t)$.  This matches the right-hand side of the original ODE, so we have verified that our solution is correct.  A plot of the solution appears below.

        {
        \vspace{2mm}
\centering
\includegraphics[scale=.7]{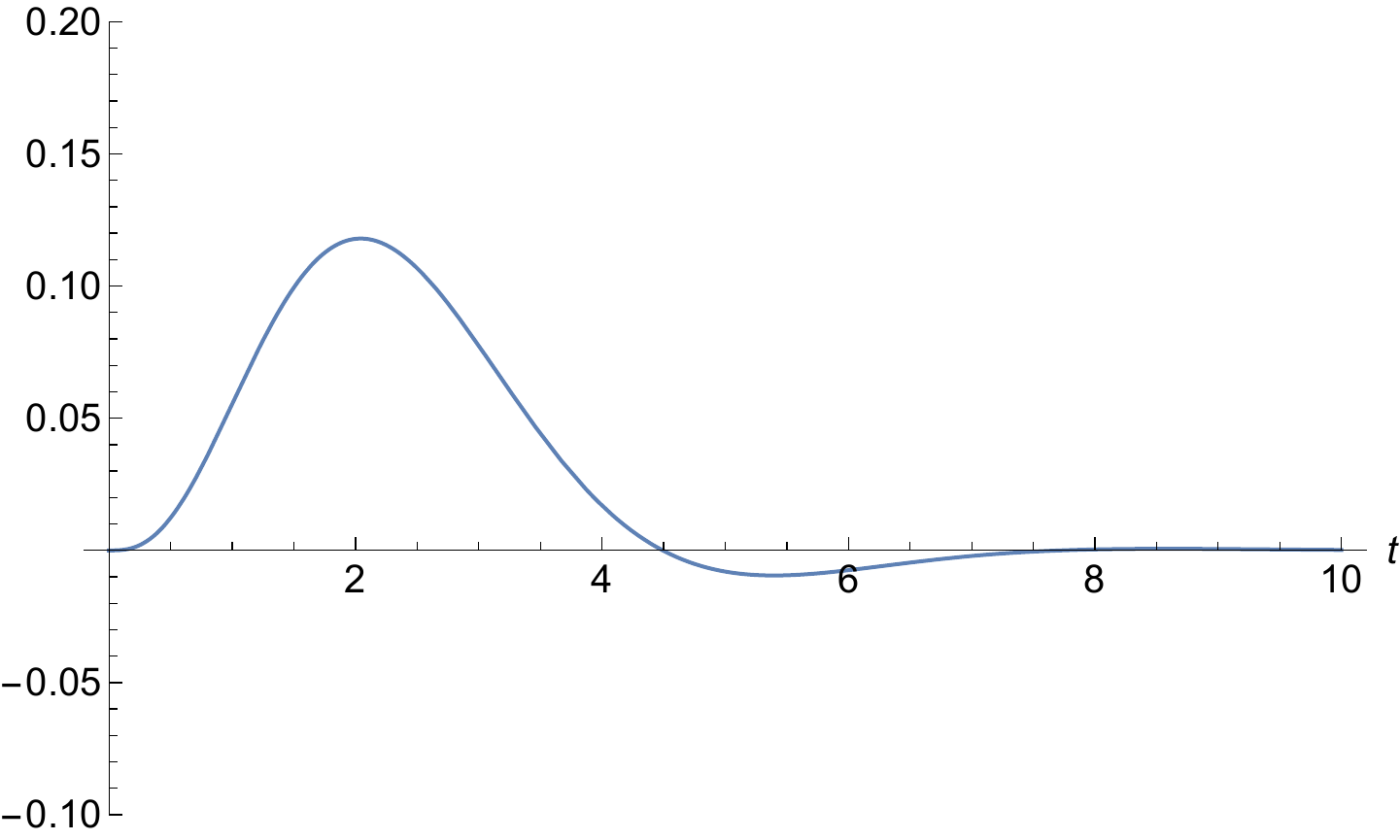}
\vspace{-2mm}
\captionof{figure}{The solution for the oscillating mass problem described above.  Here, the vertical axis defines the displacement relative to the equilibrium (zero) position.}  
} 
    \end{example}
\end{svgraybox}

\section{Solutions to PDEs Using Laplace Transforms}

The solution to PDEs using Laplace transforms is very similar to the use of Fourier transforms for the same purpose.  The primary benefit realized by using the Laplace transform is, as for the Fourier transform, the conversion of derivatives into algebraic quantities.  As might be expected, the Laplace transform is a tool that to transform independent variables on the half-line.  

There are essentially no new concepts to learn in applications of the Laplace transform to PDEs.  However, the practice of solving such problems does, at times, require a little creativity.  To illustrate, some examples of solving PDEs using the Laplace transform are given in the following.

\begin{svgraybox}
    \begin{example}[All roads lead to Rome: Solution to the transient first-order wave equation in an infinite domain]\label{wave_laplace1}
    In Chap. \ref{introPDEs} we introduced the first-order wave equation.  That equation was revisited in Chap. \ref{Fouriertransforms}, where we used Fourier transforms to solve the problem on the real line in space for a specified initial condition.  In particular, in Example \ref{firstorderwaveFT}, we found the solution to the following problem.

\begin{align*}
    &&\frac{\partial u}{\partial t}&= -c \frac{\partial u}{\partial x}\qquad -\infty < x < \infty\\
&B.C.1&    u(x,t) & \textrm{~~remains bounded for all $x$ and $t$} &&\\
&I.C.& u(x,0)&=\exp(-x^2)
\end{align*}
Using Fourier transforms, we found that the solution was
\begin{equation*}
    u(x,t) = \exp\left[-(x-c t)^2 \right]
\end{equation*}
This solution should, mathematically speaking, be entirely independent of the method of solution (assuming that the method is a valid one!)  Thus, we \emph{should be able} to find this same solution via the Laplace transform.  To do so, we can, as usual, take the transform of the PDE (which will automatically incorporate the initial condition when using the Laplace transform).  The result of that step is as follows

\begin{align*}
    &&\mathscr{L}\left[\frac{\partial u}{\partial t}\right]&=\mathscr{L}\left[ -c \frac{\partial u}{\partial x}\right]\qquad -\infty < x < \infty\\
&I.C.& u(x,0)&=\exp(-x^2)
\end{align*}
with the understanding that the solution must remain bounded.  Taking the Laplace transform, we find (after some rearrangement) the following first-order nonhomogeneous ODE.

\begin{align*}
   \frac{\partial U(s)}{\partial x}+ \dfrac{s}{c} U(s)=\frac{1}{c}\exp(-x^2)   \\
\end{align*}
The solution to linear first-order ODEs in two variables is always possible (assuming that they are well-posed), at least as an integral solution.  Recall, we need an \emph{integrating factor} of $\exp(s x/c)$ to solve this particular problem.  Multiplying both sides of the equation by the integrating factor yields

\begin{align*}
   \exp\left(\frac{s x}{c}\right)\frac{d U(x,s)}{d x}+ \exp\left(\frac{s x}{c}\right)\dfrac{s}{c} U(x,s)=\frac{1}{c}\exp\left(\frac{s x}{c}\right)\exp(-x^2) 
\end{align*}
The first two terms are the derivative of $\exp\left(\frac{s x}{c}\right) U(x,s)$, thus we have

\begin{align*}
  \dfrac{d}{dx}\left[ \exp\left(\frac{s x}{c}\right)U(x,s)\right]=\frac{1}{c}\exp\left(\frac{s x}{c}\right)\exp(-x^2)  
\end{align*}
Integrating both sides with respect to $x$ gives us

\begin{align*}
 \int_{z=-\infty}^{z=x} \dfrac{d}{dz}\left[ \exp\left(\frac{s z}{c}\right)U(x,s)\right]\, dz
& =
  \frac{1}{c}\int_{z=-\infty}^{z=x} \exp\left(\frac{s z}{c}\right)\exp(-z^2)\, dz   
\end{align*}
where here we have been careful to use a variable of integration ($z$) because the independent variable $x$ appears in the bounds of the integral.
The left-hand side of this equation is easy to evaluate.  The right-hand side is actually somewhat \emph{difficult} to evaluate because it integrates to some kind of an error-function.  Rather than evaluate the right-hand side, we will leave it in integral form for now.  This is \emph{not} usual.  However, it will ultimately serve us well (and the author, of course, has already solved this problem, so this choice comes from knowledge of at least one method that leads to a solution!)  Thus, we have

\begin{align*}
\left.\left[ \exp\left(\frac{s z}{c}\right)U(x,s)\right]\right|_{z=-\infty}^{z=x}
& =
  \frac{1}{c}\int_{z=-\infty}^{z=x }\exp\left(\frac{s z}{c}\right)\exp(-z^2)\, dz   
\end{align*}
and evaluating, we have
\begin{align*}
 \exp\left(\frac{s x}{c}\right)U(x,s)
& =
  \frac{1}{c}\int_{z=-\infty}^{z=x} \exp\left(\frac{s z}{c}\right)\exp(-z^2)\, dx   \\
  \intertext{or, simplifying}
  U(x,s)&=  \frac{1}{c}\int_{z=-\infty}^{z=x} \exp\left(\frac{s (z-x)}{c}\right)\exp(-z^2)\, dz
\end{align*}
Recall that the Laplace transform is an integration.  You may recall from calculus that we can change the order of integration of two integrals if the functions involved have defined integrals, and the region that the double integral covers is \emph{simple} (that is, the region formed on the $s-x$ plane is bounded by two functions).  For improper integrals, the situation is a bit more complicated, and the general result is known as Fubini's theorem.\indexme{Fubini's theorem}  In short, this theorem says that we can change the order of integration if the absolute value of the functions converge for each integral.  This is true for our integrals, so changing the order of integration is allowed. (N.B., for the purposes of comparison of Fourier and Laplace transforms, this example goes beyond what are covering in this text; you would not be expected to know when one might apply Fubini's theorem!)  In our case, we have the following result

\begin{align*}
\mathscr{L}^{-1}[U(x,s)]&=  \mathscr{L}^{-1}\left[ \frac{1}{c}\int_{z=-\infty}^{z=x} \exp\left(\frac{s (z-x)}{c}\right)\exp(-z^2)\, dz\right]\\
&=  \frac{1}{c}\int_{z=-\infty}^{z=x}  \mathscr{L}^{-1}\left[\exp\left(\frac{s (z-x)}{c}\right)\right]\exp(-z^2)\, dz
\end{align*}
Consulting the Laplace transform tables, this gives

\begin{align*}
u(x,t)&=  \mathscr{L}^{-1}\left[ \frac{1}{c}\int_{z=-\infty}^{z=x} \exp\left(\frac{s (z-x)}{c}\right)\exp(-z^2)\, dz\right]\\
&=  \frac{1}{c}\int_{z=-\infty}^{z=x}  \mathscr{L}^{-1}\left[\exp\left(\frac{s (z-x)}{c}\right)\right]\exp(-z^2)\, dz\\
&=\frac{1}{c}\int_{z=-\infty}^{z=x}  \delta\left( \frac{t c - x + z}{c}\right)]\exp(-z^2)\, dz
\end{align*}
Here, the delta function is nonzero only for $t c - x + z=0$ or $z=x-c t$.  Thus the integration yields the result

\begin{align*}
u(x,t)&=  \exp[-(x-c t)^2]
\end{align*}
Comparing with the result via the Fourier transform, we see that we have arrived at the same answer.  However, this case represents one where the Laplace transform is perhaps \emph{not} the most convenient route to a solution.  This is a good concept to keep in mind: while all valid mathematical methods applied to a particular problem should lead to the \emph{same} solution, it is not necessarily true that each of them require the same effort!  In this case, the Fourier transform solved the problem for the first-order wave equation in just a few lines of mathematics.  For the Laplace transform, the results are the same, but the method was substantially more complicated.  This is not true in general -- for some problems, the Laplace transform will yield more easily obtained solutions than the Fourier transform.  When both methods are suitable, one can try both methods.  It should become reasonably clear after a bit of analysis if one method will yield results more easily than the other.
    \end{example}
\end{svgraybox}

\begin{svgraybox}
    \begin{example}[Solution to the transient second-order wave equation in an infinite domain]\label{wave_laplace}
    Solutions to the wave equation are to many less intuitive than solutions to, for example, the heat/diffusion equation.  This is a result of the second-order derivative in time, which does not generate the same intuitive interpretation as does the first derivative in time (which can be interpreted an accumulation at fixed location term)  

   Laplace transforms allow solutions to the wave equation by transforming the time variable; for second-order derivatives in time, this is especially helpful.  Let's specify a wave equation in an half-infinite domain.  Here, we will have the boundary condition at $x=0$ to ``drive" the problem. 
   \begin{align*}
     &&   \dfrac{\partial^2 u}{\partial t^2} & = c^2 \dfrac{\partial^2 u}{\partial x^2} \\
     &B.C.1&   u(0,t) &=\sin(t)\\
     &B.C.2& u(x,t)& ~\textrm{ remains bounded for all $x>0$ and $t>0$}\\
      &I.C.1&  u(x,0)&=0\\
      &I.C.2& \left. \dfrac{\partial u}{\partial t}\right|_{(x,0)}& = 0\\
    \end{align*}
Taking the Laplace transform of the PDE gives us the following result

\begin{align*}
 &&  \frac{d^2 U}{d x^2} -s^2 U(s) &=0&& \\
 & B.C.1.& U(x,s) &= \dfrac{1}{s^2+1}&&\\
  &B.C.2& U(x,s)& ~\textrm{ remains bounded for all $x>0$ and $s>0$}
\end{align*}
 This is a homogeneous, second-order ODE in the variable $x$ with a well-specified boundary condition at $x=0$.  Note that here the two boundary conditions have been transformed as well as the PDE itself.  Also note, as it typical for Laplace transforms, the ancillary conditions associated with the variable being transformed (in this example, the initial conditions, because the transformed variable is $t$) are incorporated directly into the transformation via the rules for transforms of derivatives.  

 To solve this problem, we use the standard methods for second-order ODEs with constant coefficients.  The characteristic equation is $r^2 - s^2/c^2 r = 0$.  Using the quadratic formula (noting $a=1$, $b=0$, and $c=-s^2$, this yields two real roots, $r_1=-s/c$ and $r_2=s/c$.  The solution is

 \begin{equation*}
     U(x,s) = c_1 \exp\left(-\frac{s x}{c} \right)+ c_2  \exp\left(\frac{s x}{c} \right)
 \end{equation*}
Noting that the solution must remain bounded for all $x>0$ and $t>0$, then we must have that $c_2=0$.  The constant $c_1$ can be found using the boundary condition.  At $x=0$ we find that $c_1=1/(s^2+1)$.  The solution is, then
 \begin{equation*}
     U(x,s) = \dfrac{1}{s^2+1} \exp\left(-\frac{s x}{c} \right)
 \end{equation*}
The final step is to convert this solution back into real time rather than the transform variable $s$.  Taking the inverse transform of both sides of the solution gives us

 \begin{equation*}
     \mathscr{L}^{-1}[U(x,s)] =\mathscr{L}^{-1}\left[ \dfrac{1}{s^2+1} \exp\left(-\frac{ x}{c}s \right)\right]
 \end{equation*}
There are two ways that the inversion of the right-hand side can be done.  First, we might note that $1/(s^2+1)$ and $\exp(-xs/c)$ are both transforms that are available on the table of Laplace transforms.  Thus, the convolution method could be used here.  A simpler approach is to recall the \emph{shifting} property of the Laplace transform.  The exponential term in the result can also be interpreted as a shifting operator, whose role is to shift the independent variable $t$ by the amount $-x/c$ during the inversion of $1/(s^2+1)$.  Taking this latter approach, we find the solution to be 

\begin{equation*}
    u(x,t) = \sin\left(t-\tfrac{x}{c}\right)H\left(t-\tfrac{x}{c}\right)
\end{equation*}

The solution for a few representative times are provided in the plot below.

{
\vspace{2mm}
\centering
\includegraphics[scale=.65]{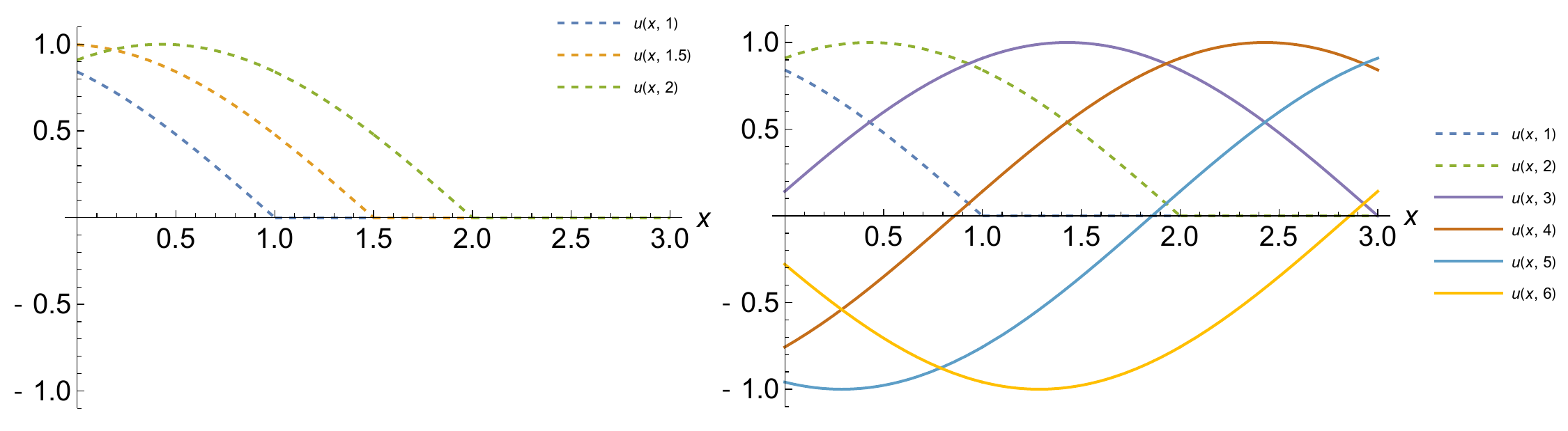}
\vspace{-2mm}
\captionof{figure}{The solution to the wave equation on a half-infinite domain.  For this plot, the vertical axis is wave displacement (distance), $x$ is measured in \si{m}, and  $c=1~\si{m/s}$. The times are between $t=1$ and $t=6$ \si{s}.}  
}
    \end{example}
\end{svgraybox}

\begin{svgraybox}
    \begin{example}[solution to the transient heat/diffusion equation on the half-line]\label{heat_laplace}
    The heating of an insulated rod by the application of a constant temperature at one end is a problem we have examined previously using separation of variables.  If the rod is very long compared to the characteristic time scale of heat transport and only ``early" times (several multiples of the characteristic time scale) are of interest, it can be useful to treat the rod as being infinitely long.

    Our heat transport problem will take the following form.  Here, assume that the temperature, $u$ is given in non-dimensional form (i.e., it has been normalized such that the boundary temperature is 1 and the initial temperature is zero).
    \begin{align*}
      &&  \frac{\partial u}{\partial t}& = K \frac{\partial^2 u}{\partial x^2}&&  \\
      &B.C.1& u(0,t) &= u_1 &&\\  
       &B.C.2& u(x,t)&< \infty\\  
       & I.C.& u(x,0) &= 0
    \end{align*}
    
    The interpretation here is that the rod is initially of uniform temperature, $u=0$.  At $t=0$, the left end is suddenly put in contact with a large reservoir of heat at temperature $u(0,t) = u_1$.  The second boundary condition is set to require that the temperatures are not allowed to grow without bound, which would be physically unrealistic.

    First, note that in principle, we could conduct the Laplace transform on either $t$ or on $x$ (or both!), since both variables are defined on the interval $[0,\infty)$.  To keep the problem from getting complicated, we will transform only the variable $t$ (which is the conventional choice).  
    
    The solution can be found by determining Laplace transform of the PDE and the two boundary conditions.  The boundary conditions are assumed to depend on time, even if they are constant functions of time (and this is clear in the list of independent variables for the two boundary conditions).  The transform is given by the following

        \begin{align*}
      &&  sU(s)-0& = K \frac{\partial^2 U}{\partial x^2}&& 0\le x <\infty,~ t>0 \\
      &B.C.1& U(0,s) &= \frac{u_1}{s} &&\\  
       &B.C.2& U(x,s)  & < \infty&&\\   
       & I.C.& u(x,0) &= 0
    \end{align*}
Note that, because we are taking the transform of the time variable, the Laplace transform commutes with the spatial derivative.  In other words, \emph{assuming that $t$ is the transformed variable}, then $\mathscr{L}[\partial u/\partial x] = \partial/\partial x \mathscr{L}[u]= \partial U/\partial x$.  Obviously this argument can be repeated for higher-order derivatives.  Now we have a second-order, non-homogeneous ODE with constant coefficients.  The ODE takes the form

\begin{align*}
   &&    \frac{\partial^2 U}{\partial x^2} - \frac{s}{K}U(s) &= 0\\
       &B.C.1& U(0,s) &= \frac{u_1}{s} &&\\  
       &B.C.2& U(x,s)  & < \infty&&\\   
\end{align*}
The solution to this ODE is straightforward.  The characteristic equation is given by $r^2 -s/K r =0$, and this indicates that $a=1$, $b=0$, and $c=-s/K$.  Thus the roots are $r_1 = +\sqrt{4s/K}$ and $r_2 = -\sqrt{4s/K}$.  The solution must be of the general form

\begin{equation*}
    U(x,s) = c_1 \exp\left(-2x\sqrt{\frac{s}{K}}\right) +c_2\exp\left(2x\sqrt{\frac{s}{K}}\right) 
\end{equation*}
And, without much effort, we can immediately see that $c_2=0$ because of the boundedness imposed by the second boundary condition.
\begin{equation*}
    U(x,s) = c_1 \exp\left(-2 x\sqrt{\frac{s}{K}} \right)
\end{equation*}
Using the first boundary condition, we find that $c_1= \frac{u_1}{s}$, which gets us to the point (noting that the quantity $2x$ is moved to inside the root)  
\begin{equation*}
    U(x,s) = \dfrac{u_1}{s} \exp\left(-\sqrt{\frac{4x^2}{K}}\sqrt{s} \right)
\end{equation*}
From the table of Laplace transforms, the entry number \ref{laplacetable1} is the one needed here.  Taking $a=\sqrt{\frac{4x^2}{K}}$, then we have the Laplace inversion of the form

\begin{align*}
    u(x,t) &= u_1 \textrm{erfc}\left( \frac{\sqrt{\frac{4x^2}{K}}}{2\sqrt{t}}\right) \\
    &= u_1 \textrm{erfc}\left( \sqrt{\frac{x^2}{K t}}\right)
\end{align*}

The solution for a few representative times are provided in the plot below.
{
\vspace{2mm}
\centering
\includegraphics[scale=.7]{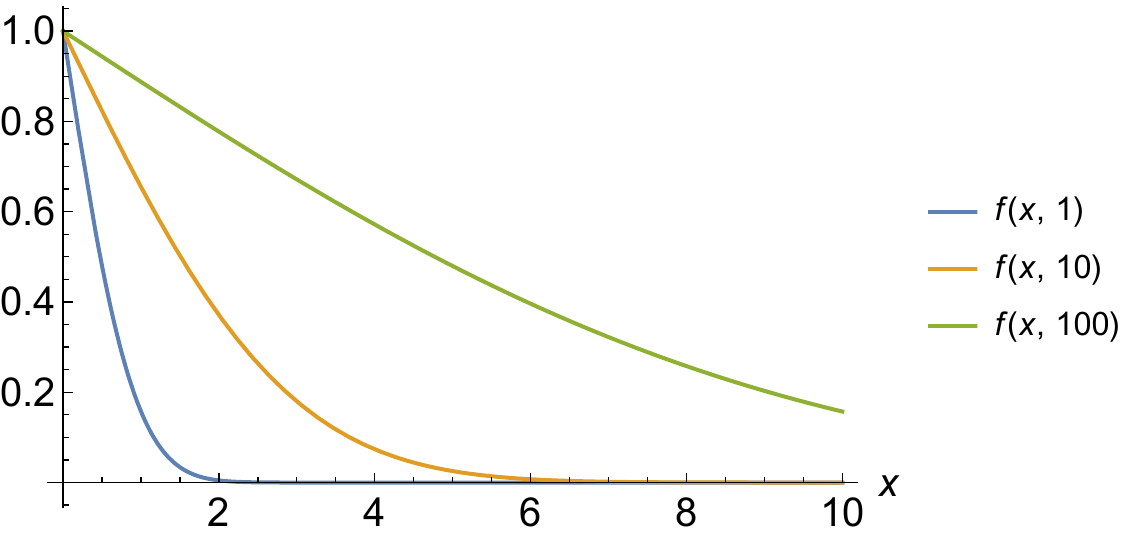}
\vspace{-2mm}
\captionof{figure}{The solution to the heat equation on a half-infinite domain.  For this plot, the vertical axis is temperature (difference), the temperature, $u$ has been non-dimensionalized such that $u_1=1$, $x$ is measured in \si{mm}, and  $K=1~\si{mm^2/min}$. The times are taken to be $t=1$, 10, and 100 minutes.}  
}
    \end{example}
\end{svgraybox}




\newpage
\section*{Appendix 2. Laplace Transform Table} 
\begin{table}
 \centering
\captionsetup{singlelinecheck=off}
\def\arraystretch{2.3}
\caption{Laplace transform pairs.}
 $\textrm{It is assumed that $a$ is a real number ($a\in\mathbb{R})$, and $Re(s)>0$ unless other constraints are provided}$.\\
 $n$ is assumed to be a positive integer ($n\in \mathbb{N}$)\\
{\large
\setcounter{magicrownumbers}{0} 
\begin{tabular}{rccccc}
&$f(t)$ & ~~~~~~~~~ & $F(s)$   &~~~ & \hspace{-30mm}Constraint ~~~~~~~~~\\ 
\cline{1-6}
\rownumber~~~&$\dfrac{\partial f(t)}{\partial t}$&   & $s F(s)-f(0)$ &\\
\rownumber~~~&$\dfrac{\partial^2 f(t)}{\partial t^2}$&   & $s^2 F(s)-s f(0)-f'(0)$ &\\
\rownumber~~~&$\dfrac{\partial^n f(t)}{\partial t^n}$&   & $\begin{aligned}s^n F(s)-s^{(n-1)} f(0)&-s^{(n-2)}f'(0)\\&-\ldots - f^{(n-1)}(0)\end{aligned}$ &\\
\rownumber~~~&$\exp(-a t)$&   & $\dfrac{1}{s+a}$ & $a>0$\\
\rownumber~~~&$\exp(a t)$&   & $\dfrac{1}{s-a}$ & ~~$s>a \textrm{~and~} a>0$\\
\rownumber~~~&$\exp(a t)(1+a t)$&   & $\dfrac{s}{(s-a)(s-a)}$ & ~~$s>a \textrm{~and~} a>0$\\
\rownumber~~~&$\dfrac{\exp(-a t)-\exp(-b t)}{a+b}$&   & $\dfrac{1}{(s+a)(s+b)}$ & ~~$ a,b>0$\\
\rownumber~~~&$\dfrac{\exp(a t)-\exp(b t)}{a-b}$&   & $\dfrac{1}{(s-a)(s-b)}$ & ~~$s>a,b \textrm{~and~} a,b>0$\\
\rownumber~~~&$\dfrac{\exp(a t)-\exp(-b t)}{a+b}$&   & $\dfrac{1}{(s-a)(s+b)}$ & ~~$s>a,b \textrm{~and~} a,b>0$\\
\rownumber~~~&$\frac{t e^{-a t}}{(a-b)^2}+\frac{t e^{-b t}}{(a-b)^2}+\frac{2 e^{-a t}}{(a-b)^3}-\frac{2
   e^{-b t}}{(a-b)^3}$&   & $\dfrac{1}{(s+a)^2(s+b)^2}$ & ~~$a,b\in\mathbb{R}$\\

\protect\raisebox{8mm} \rownumber~~~&$\begin{aligned}
\frac{1}{3} e^{-t/2}& \Bigg[\sqrt{3} \sin \left(\tfrac{\sqrt{3} t}{2}\right)\\
&+3 \cos
   \left(\tfrac{\sqrt{3} t}{2}\right)\Bigg]\end{aligned}$&   & $\dfrac{s+1}{s(s+1)+1}$ & \\
\end{tabular}
}
\end{table}


\begin{table}
\centering
\def\arraystretch{3}
{\large
\begin{tabular}{rccccc}
&$f(t)$ & ~~~~~~~~~ & $F(s)$   &~ & \hspace{-18mm}Constraint ~~~~~~~~~\\ 
\cline{1-6}
%
\protect\raisebox{8mm} \rownumber~~~&$\begin{aligned}
\tfrac{1}{5} e^{t/2} &\Bigg[5 \cosh \left(\tfrac{\sqrt{5} t}{2}\right)\\
&-\sqrt{5} \sinh \left(\tfrac{\sqrt{5} t}{2}\right)\Bigg]\end{aligned}$&   & $\dfrac{s-1}{s(s-1)-1}$ & ~~$s>1$\\
\rownumber~~~&$t^0=1$&   & $\dfrac{1}{s}$ & \\
\rownumber~~~&$t$&   & $\dfrac{1}{s^2}$ &\\
\rownumber~~~&$t^2$&   & $\dfrac{2}{s^3}$& \\
\rownumber~~~&$t^{\tfrac{1}{2}}$&   & $\dfrac{\sqrt{\pi}}{2s}\sqrt{\dfrac{\pi}{s}}$ &\\
\rownumber~~~&$t^{-\tfrac{1}{2}}$&   & $\sqrt{\dfrac{\pi}{s}}$ &\\
\rownumber~~~&$t^n$&   & $\dfrac{n!}{s^{(n+1)}}$ &~~~$n \in \mathbb{N}$\\
\rownumber~~~&$t^{a}$&   & $s^{-(1+a)}\Gamma\left(1+a \right)$ & ~~~$a > -1$ \\
\rownumber~~~&$t^{-a}$&   & $s^{-(1-a)}\Gamma\left(1-a \right)$ & ~~~~~~~~~~~~~$a \ge 0$ \\
\rownumber~~~&$t^{\tfrac{1}{n}}$&   & $s^{-\left(1+\tfrac{1}{n}\right)}\Gamma\left(1+\dfrac{1}{n} \right)$ & ~~~~~~~~~~~~~$n \in \mathbb{N}$ \\
\rownumber~~~&$t^{-\tfrac{1}{n}}$&   & $s^{-(1-\tfrac{1}{n})}\Gamma\left(1-\dfrac{1}{n} \right)$ & ~~~~~~~~~~~~~$n \in \mathbb{N}$ \vspace{3mm}\\
\end{tabular}
}
\end{table}
\begin{table}
\centering
\def\arraystretch{2.9}
{\large
\begin{tabular}{rccccc}
&$f(t)$ & ~~~~~~~~~ & $F(s)$   &~~  Constraint ~\\ 
\cline{1-5}
&~~~ & ~~~~~~~~~ & ~~~  &~~~~~~ & ~~~~~~~~~\vspace{-15mm}\\ 
\protect\raisebox{5mm} \rownumber~~~&$\dfrac{1}{a^2+t^2}$ & &$\begin{aligned}-\dfrac{1}{a}&\text{Ci}\left({a s}\right) \sin\left({a s}\right)\\
+&\dfrac{1}{2 a} \left[\pi -2 \, \text{Si}\left({a s}\right)\right]\cos\left({a s}\right)\end{aligned}$&  ~~~~~~~~~~~~~$a>0$\vspace{2mm}\\
\protect\raisebox{5mm} \rownumber~~~&$\dfrac{t}{a^2+t^2}$ & &$\begin{aligned}~&\text{Ci}\left({a s}\right) \cos\left({a s}\right)\\
+&\tfrac{1}{2} \left[\pi -2 \, \text{Si}\left({a s}\right)\right]\sin\left({a s}\right)\end{aligned}$&  ~~~~~~~~~~~~~$a>0$\\
\rownumber~~~&$J_0(at)$&   & $\dfrac{1}{\sqrt{s^2+a^2}}$ &~~~~~~~~~~~~~$a>0$\\
\rownumber~~~&$\textrm{erf}\left(\frac{a}{2\sqrt{t}}\right)$&   & $\dfrac{1-\exp(-a\sqrt{s})}{s}$ &~~~~~~~~~~~~~$a>0$\\
\rownumber\newtag{\arabic{magicrownumbers}}{laplacetable1}~~~&$\textrm{erfc}\left(\frac{a}{2\sqrt{t}}\right)$&   & $\dfrac{\exp(-a\sqrt{s})}{s}$ &~~~~~~~~~~~~~$a>0$\\
\rownumber~~~&$\delta(t)$&   & $1$ &\\
\rownumber~~~&$\delta(t-a)$&   & $\exp(-a s)$ &~~~~~~~~~~~~~$a>0$\\
\rownumber~~~&$\begin{cases}H(t-a)\\1-H(a-t)\end{cases}$&   & $\dfrac{\exp(-a s)}{s}$ & ~~~~~~~~~~~~~$a>0$\\
\rownumber~~~&$\begin{cases}H(a-t)\\1-H(t-a)\end{cases}$&   & $\dfrac{1-\exp(-a s)}{s}$ & ~~~~~~~~~~~~~$a>0$\\
\protect\raisebox{3mm} \rownumber~~~&
$\begin{aligned}
B(t;a,a+\Delta)=& \\ 
[H(t-a)-&H(t-(a+\Delta)] 
\end{aligned}$ &   &$\dfrac{\exp(-a s)[1-\exp(-\Delta s)]}{s}$&$~~~~~~~~~~~~~\Delta>0 \textrm{~and~} a>0$\vspace{4mm}\\
\protect\raisebox{9mm} \rownumber~~~& $\begin{aligned}
T(t;a)=~~ &\\
\dfrac{t}{a}[H(t)&-H(t-a)]\\
+\Big(2-&\frac{t}{a}\Big) [H(t-a)- H(t-2 a)] 
 \end{aligned}$&  &$\dfrac{\exp(-2 a s) [\exp(a s)-1]^2}{s^2}$ &~~~~~~~~~~~~~$a>0$\\
\end{tabular}
}
\end{table}
\begin{table}
\centering
\def\arraystretch{3}
{\large
\begin{tabular}{rccccc}
&$f(t)$ & ~~~~~~~~~ & $~~~~~~~~~~~~~~~F(s)~~~~~~~~~~~~~$   &~~~~~~~~~~~~~~~~~  Constraint ~~~~~~~~~~\\ 
\cline{1-6}
&~~~ & ~~~~~~~~~ & ~~~  &~~~~~~~~~~~~~~~~~ & ~~~~~~~~~~~~~~~\vspace{-15mm}\\ 
\rownumber~~~&$\sin(a t)$&   & $\dfrac{a}{s^2+a^2}$ & ~~$a\in \mathbb{R}$\\
\rownumber~~~&$\cos(a t)$&   & $\dfrac{s}{s^2+a^2}$ & ~~$a\in \mathbb{R}$\\
\rownumber~~~&$\exp(-b t)\sin(a t)$&   & $\dfrac{a}{a^2+(b+s)^2}$ & ~~$a\in\mathbb{R},~ b\in \mathbb{R}$\\
\rownumber~~~&$\exp(-b t)\cos(a t)$&   & $\dfrac{b+s}{a^2+(b+s)^2}$ & ~~$a\in\mathbb{R},~ b\in \mathbb{R}$\\
\rownumber~~~&$\sin\left(\tfrac{\pi}{a} t\right)[H(t)-H(t-a)]$ && $\dfrac{\pi a [1-\exp(-a s)]}{a^2s^2+\pi^2}$ &~~~~~~~~~~~~~$a>0$\\
\rownumber~~~&$\cos\left(\tfrac{\pi}{a} t\right)[H(t)-H(t-a)]$ && $\dfrac{s a^2[1 + \exp(-a s)]}{a^2 s^2+\pi ^2}$ &~~~~~~~~~~~~~$a>0$\\
\rownumber~~~&$\sin(a t+b_0)$&   & $\dfrac{s\sin(b_0)+a\cos(b_0)}{s^2+a^2}$ & ~~$a<0$ or $a=0$  or $a>0$\\
\rownumber~~~&$\cos(a t+b_0)$&   & $\dfrac{s\cos(b_0)-a\sin(b_0)}{s^2+a^2}$ & ~~$a<0$ or $a=0$  or $a>0$\\
\rownumber~~~&$t\sin(a t)$&   & $\dfrac{2 a s}{(s^2+a^2)^2}$ \\
\rownumber~~~&$t\cos(a t)$&   & $\dfrac{s^2-a^2}{(s^2+a^2)^2}$ & \\
\end{tabular}
}
\end{table}
%
\begin{table}
\centering
\def\arraystretch{3}
{\large
\begin{tabular}{rccccc}
&$f(t)$ & ~~~~~~~~~ & $~~~~~~~~~~~~~~~F(s)~~~~~~~~~~~~~$   &~~~~~~~~~~~~~~~~~  Constraint ~~~~~~~~~~\\ 
\cline{1-6}
&~~~ & ~~~~~~~~~ & ~~~  &~~~~~~~~~~~~~~~~~ & ~~~~~~~~~~~~~~~\vspace{-15mm}\\ 
\rownumber~~~&$\sinh(a t)$&   & $\dfrac{a}{s^2-a^2}$ &  $\begin{aligned}~~&a<0 \textrm{~or~} a=0\textrm{~or~} a>0\\ &\textrm{~and~} s>\lvert a \rvert\end{aligned}$\\
\rownumber~~~&$\cosh(a t)$&   & $\dfrac{s}{s^2-a^2}$ & $\begin{aligned}~~&a<0 \textrm{~or~} a=0\textrm{~or~} a>0\\ &\textrm{~and~} s>\lvert a \rvert\end{aligned}$\\
\end{tabular}
}
%
\begin{align*}
\textrm{Notes:}&&\\
    &J_0& &\textrm{Bessel function of the first kind, zero order}\\
     &\textrm{erf}& &\textrm{The error function,}~~~  \textrm{erf}(t)=\frac{2}{\sqrt{\pi}}\int_0^{t} \exp(-\tau^2)\,d\tau. \\
     &\textrm{erfc}& &\textrm{The complimentary error function},~~~ 1-\textrm{erf}(t)\\
     &\delta& & \textrm{The Dirac delta function} \\
     &H&& \textrm{The Heaviside function}\\
     &B(x;a,a+\Delta) && \textrm{The boxcar function on $t\in[a,a+\Delta]$}\\
     &T(t;a)&& \textrm{The symmetric triangle function on $t\in[0,2a]$, with maximum of 1 located at $t=a$}\\
     & \text{Ci} && \textrm{The cosine integral function, ~~~} \text{Ci}(z) = -\int_{\tau=z}^\infty \frac{\cos(\tau)}{\tau} \,d\tau\\
   & \text{Si} && \textrm{The sine integral function, ~~~}\text{Si}(z) = -\int_{\tau=0}^{\tau =z} \frac{\sin(\tau)}{\tau} \,d\tau \\
   & \Gamma &&\textrm{The gamma function}
\end{align*}
\end{table}
\clearpage

\section*{Problems}
\subsection*{Practice Problems}
Solve the following problems using Laplace transforms.
\begin{enumerate}
    \item $\begin{aligned}[t] ~~u''(t) + 2 u'(t) + 5 u(t) &= 0 & u(0)&=1 &u'(0)=1\end{aligned}$
    \item $\begin{aligned}[t]~~ u''(t) - 2 u'(t) + 2 u(t) &= 0 & u(0)&=0 &u'(0)=1\end{aligned}$
    \item $\begin{aligned}[t]~~ u''(t) - 2 u'(t) + 2 u(t) &= 0 & u(0)&=5 &u'(0)=0\end{aligned}$
    \item $\begin{aligned}[t]~~ u''(t) - 3 u'(t) + 2 u(t) &= 0 & u(0)&=0 &u'(0)=-1\end{aligned}$
    \item $\begin{aligned}[t]~~ u''(t)-3 u'(t)+2 u(t)&=10 \cos(t) & u(0)&=0 &u'(0)=-1\end{aligned}$
      \item $\begin{aligned}[t]~~ u''(t)-2 u'(t)+2 u(t)&=\exp(-t)  & u(0)&=0 &u'(0)=0\end{aligned}$
    \item $\begin{aligned}[t]~~ u''(t)-4 u'(t)+4 u(t)&=\cos (t)  & u(0)&=0 &u'(0)=0\end{aligned}$
    \item $\begin{aligned}[t]~~ u''(t)-
    u'(t)+4 u(t)&=\cos (t)  & u(0)&=0 &u'(0)=0\end{aligned}$
    \item $\begin{aligned}[t]~~ u''(t)-4 u'(t)&=H(t-5)  & u(0)&=0 &u'(0)=0\end{aligned}$
    \item $\begin{aligned}[t]~~ u''(t)- u'(t)&= 3t^2 H(t-1)  & u(0)&=0 &u'(0)=1\end{aligned}$
    
\end{enumerate}



\subsection*{Applied and More Challenging Problems}

\begin{enumerate}
    \item Prove that $\mathscr{L}\left[ \dfrac{d f}{d t} \right] = s F(s)-f(0)$.  Use integration by parts to prove this assertion.

    \item Solve
    \begin{align*}
        \dfrac{du}{dt} &= -\int_0^t \exp[-(t-\tau)] u(\tau) \, d\tau \\
        u(0)&=1
    \end{align*}

    \item Redo the analysis of Example \ref{heat_laplace}, but this time assume that the initial condition is $u(x,0)=u_0$, a constant.  You will find that you generate a nonhomogeneous equation for $U$ upon transformation.  Although frequently nonhomogeneous problems create substantial extra work, sometimes when the nonhomogeneous term is a constant, we can propose a change of variables.  Please try the change of variables $V= (U + K/s u_0)$.  You should end up with the ODE and boundary conditions
    \begin{align*}
      &&    \frac{\partial^2 V}{\partial x^2} - \frac{s}{K}V(s) &= 0\\
      &B.C.1& V(0,s) &= \frac{u_1}{s}-\frac{K}{s} u_0 &&\\  
       &B.C.2& V(x,s)  & < \infty&&\\  
    \end{align*}
    Solve this problem, and then substitute back to the original variable $U$.  You should be able to invert the final result.

    \item Post's inversion theorem is not particularly handy to use, but it is one of only a few methods for avoiding contour integration in the complex plane (which is not something covered in this text).  To prove to yourself that you can derive some of the inversions, however, please try the following problem.  Using Post's inversion theorem, show that the inverse of $F(s) = 1/s$ is the function $f(t) = 1$.  Hint: Try to first reduce all powers of $k$ to simpler forms.  For example 
    \begin{align*}
        t^{k+1} t^{-(k+2)} & = (t^k) t (t^{-k}) t^{-2} = t^{-1}
    \end{align*}
    These kinds of algebraic manipulations will allow you to successfully compute the appropriate quantity using Post's inversion theorem.

    \item Solve the following problem using the Laplace transform.  Note that, unlike Example \ref{wave_laplace1}, this problem will give an integral (arising from an integration factor) that can be (and should be!) done directly before attempting to invert the result.  

    \begin{align*}
    &&\frac{\partial u}{\partial t}&= -c \frac{\partial u}{\partial x}\qquad -\infty < x < \infty\\
&I.C.& u(x,0)&=H(x)
\end{align*}
Plot your results for $c=1~\si{m/s}$ and $t=1, 2$ and $3~\si{s}$ (all on the same plot).  Plot on the domain $0~\si{m} < x < 5 ~\si{m}$.

\item Solve the following first-order wave equation using the Laplace transform.  

    \begin{align*}
    &&\frac{\partial u}{\partial t}&= -c \frac{\partial u}{\partial x}\qquad -\infty < x < \infty\\
&I.C.& u(x,0)&=H(x)\exp(-x)
\end{align*}
Upon being transformed, this problem will lead to a nonhomogeneous first-order ODE that can be solved using an integrating factor.  As you work through this problem, be sure to keep clear the distinctions between the variables $s$ and $x$.  The following are a few hints to help assure that you are successful in solving this problem.  First, note that the integrations in $x$ will have bounds between $-\infty$ and $x$.  Second, it is \emph{extremely helpful} in this case to formally use a variable of integration for integrating in space (the variable $z$ is a handy one that helps remind you that the integration is still a spatial one).  Finally, recall that the Heaviside function  when integrated can be accounted for by restricting the bounds of integration (i.e., $H(z)$ is zero for$z<0$).  As a final hint, after using the integration factor, simplification of the integral by accounting for the Heaviside function, and moving all exponential functions to the right-hand side of the result, you should find a solution to the ODE of the form

\begin{equation*}
    U(x,s) = \frac{1}{c}  \int_{z=0}^{z=x} \exp\left[-\frac{s}{c}(x-z)\exp(-z) \right] \, dz
\end{equation*}
For this integral, note that $z>0$ and that $x>z$. So, in all possible cases, this integral is composed of \emph{two decaying exponentials in $z$}, so it must be convergent.  To actually solve the integral, however, it is much more convenient to put it in the form

\begin{equation*}
    U(x,s) = \frac{1}{c} \exp\left[-\frac{s}{c}x \right] \left\{\int_{z=0}^{z=x}  \exp\left[z(\frac{s}{c}-1) \right] \, dz\right\}
\end{equation*}

This result is offered as a way to \emph{check} your results, of your analysis.  Do make sure that you show all of your work to get to this point.  Note that the integral in this result has the variable of integration $z$, so the term $(1-s/c)$ can be treated as essentially a constant for the purposes of integration.  After integrating, be sure to separate all functions involving $s$ to the extent possible.   You should find a result that can be inverted using the table of Laplace transforms.  Plot your results for $c=1~\si{m/s}$ and $t=1, 2$ and $3~\si{s}$ (all on the same plot).  Plot on the domain $0~\si{m} < x < 5 ~\si{m}$.

\item  Solve the following heat/diffusion problem using Laplace transforms.
\begin{align*}
    && \frac{\partial u}{\partial t} &= \alpha^2 \frac{\partial^2 u}{\partial x^2}&& x\in(-\infty, \infty)\\
    &B.C.1 & u(x,t)&~~\textrm{remains bounded for all $(x,t)$}&&\\
     &B.C.2 & \left.\frac{\partial u}{\partial x}\right|_{(x,t)}&~~\textrm{remains bounded for all $(x,t)$} &&\\
     &I.C.&& u(x,0)=\sin(x)+1 &&
\end{align*}
Upon transforming this problem, you will find that the result is an inhomogeneous second-order ODE in the variable $x$.  You will need to use the method variation of parameters (Chp.~\ref{variationofparams}) to complete the problem. Plot solutions for $\alpha^2=1\times 10^5~\si{cm^2}{s}$ and times $t=1,100, 1000$ \si{s} over the range $0<x<100$ \si{cm}.

\item  Solve the following heat/diffusion problem using Laplace transforms.
\begin{align*}
    && \frac{\partial u}{\partial t} &= \alpha^2 \frac{\partial^2 u}{\partial x^2}&& x\in(0, \infty)\\
    &B.C.1 & u(x,t)&=\sin(t)+1\\
     &B.C.2 & \left.\frac{\partial u}{\partial x}\right|_{(x,t)}&~~\textrm{remains bounded for all $(x,t)$} &&\\
     &I.C.&& u(x,0)=0 &&
\end{align*}
Your solution for this problem will not be explicit.  You should end up with a convolution integral as your final result.

\item  Solve the following heat/diffusion problem using Laplace transforms.
\begin{align*}
    && \frac{\partial^2 u}{\partial t^2} &= c^2 \frac{\partial^2 u}{\partial x^2}&& x\in(0, \infty)\\
    &B.C.1 & u(0,t)&=\sin(t)\\
     &B.C.2 & \left.\frac{\partial u}{\partial x}\right|_{(x,t)}&~~\textrm{remains bounded for all $(x,t)$} &&\\
     &I.C.&& u(x,0)&=0 &&
\end{align*}
Plot solutions for $\alpha^2=1$ and times $t=1,5, 10$.
\end{enumerate}

\include{chapter16_introduction_to_wavelets}
\abstract*{This is the abstract for chapter 3.}

\begin{savequote}[0.55\linewidth]
``Tesler's Theorem: AI is whatever hasn't been done yet."
\qauthor{Lawrence Gordon Tesler, American computer scientist}
\end{savequote}
\begin{savequote}[0.55\linewidth]
``To fight an exponential, it seems reasonable to arm oneself with other exponentials."
\qauthor{From \emph{On the Expressive Power of Deep Architectures} by Bengio and Delalleau (2011)}
\end{savequote}


\chapter{Primer on Feedforward Neural Networks: An Analytical Approach}\label{FNN}
%
\def\CHAP{chapter17_intro_to_the_mathematics_of_deep_neural_networks}
%

\section{Introduction}
Artificial neural networks (ANNs), or simply \emph{neural networks},  have become a regularly used tool in many areas of applied mathematics.  A succinct but not entirely general definition of an ANN is that an ANN consists of  weighted sums of compositions of functions.  Such networks are typically used to either (1) fit a function to a set of example data (regression), or (2) to segment a set of example data into categories (classification) where the boundaries of the categories are functions predicted by the network.  The basic ideas underlying neural networks are the same regardless of whether they are used for regression or classification.  To help keep this introduction simple, the primary focus of the material will be on \emph{fitting} data sets by generating approximation functions using feedforward neural networks.  From here forward we will assume that the kinds of ANNs that we discuss have the express purpose of generating approximate functions to fit a data set; however, in the background, we should keep in mind that the basics of ANNs discussed here are the for data classification as well.  The objective of deep learning is to reconstruct a relationship between input and output. We  assume that there exists an unknown function $\hat{y}$ that approximates the know data examples (which associate a collection of proposed independent variables with a set of observations (the independent variables).  Each such function $\hat{y}$ contains a number of adjustable coefficients.  To make such an approximation, we propose a metric to measure the errors between the $\hat{y}$ and the observed data, and then attempt to optimize the adjustable coefficients to provide the least error.  Training a network on the dataset should then return
a deep neural network $\hat{y}$ that is close to the observed data by the chosen metric.

For now it is enough to think of ANNs as a hierarchical network of compositions of transformation functions.  In their simplest form (the only form we will consider here), the networks are arranged in a hierarchy called, suggestively, \emph{layers}.   Within each layer, there are one or more \emph{nodes}.  A node represents a linear or nonlinear operation; we will describe these operations in more detail later.    A fully connected feedforward neural network (FNN) is one where each node within a layer receives input from all of the nodes the previous layer, and outputs information to all of the nodes in the subsequent layer.  In this material of this chapter, our focus will be on FNNs exclusively.  
            \begin{figure}[t!]
             \sidecaption[t]
           \centering
            \includegraphics[scale=0.5]{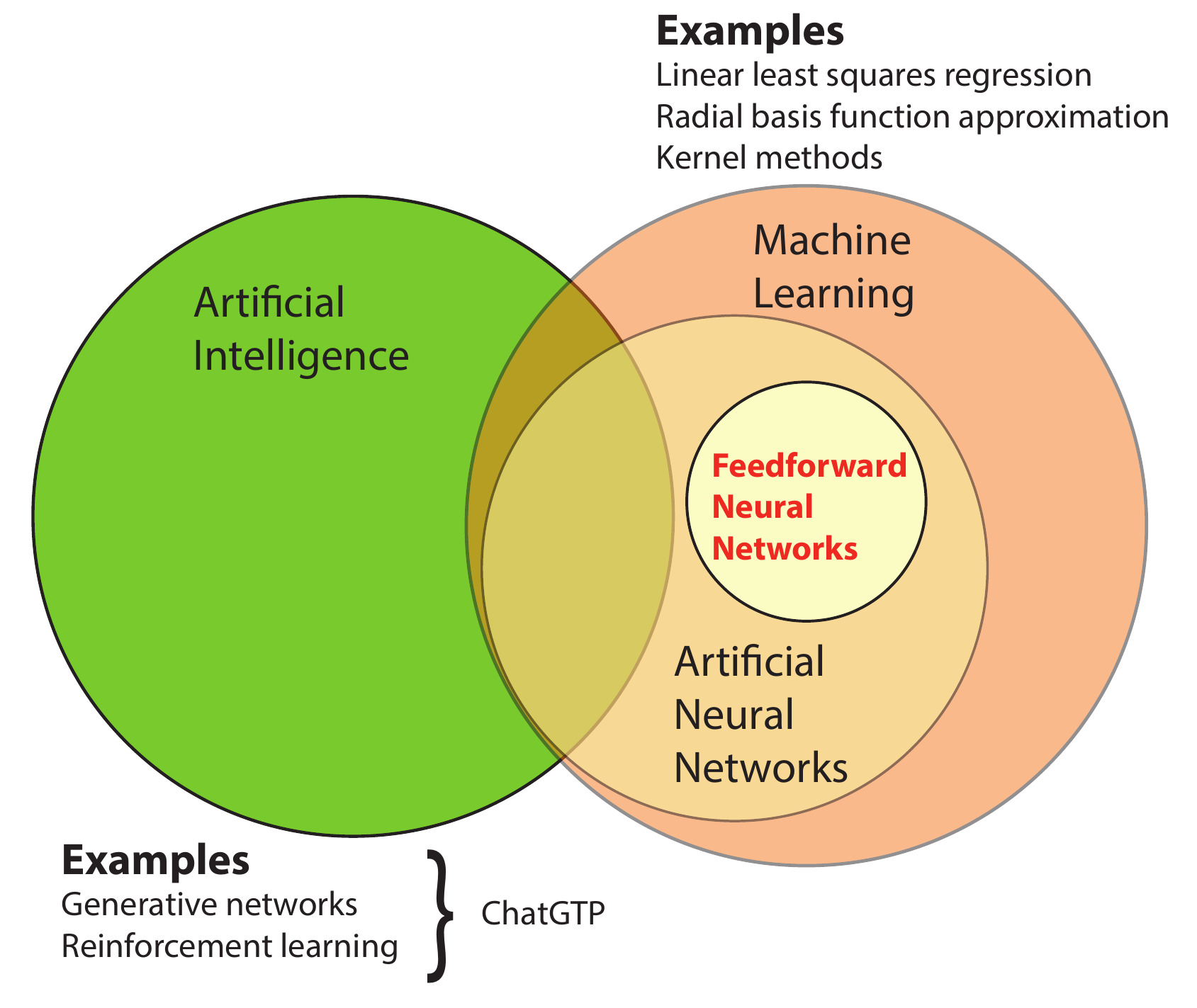}
            \vspace{-5mm}
            \caption{Machine learning is not necessarily a subset of artificial intelligence.  As an example, few would consider the method of linear least-squares fitting to represent ``artificial intelligence". However, it can be represented by a simple feedforward network, and, hence, is a subset of machine \emph{learning} methods.}
            \label{fig:types}       
            \end{figure}

The \emph{depth} of FNNs is related to how many \emph{layers} of compositions of functions there are between the input and output.  In the past there was significant discussion about the function of network depth.  It is now realized that adding layers to a network increases its expressive power (essentially, how many regions of the independent variables the functions to be estimated is broken up into) roughly exponentially as the network width and depth increases \citep{bengio2011expressive}. The increasing depth of FNNs has been critical in the development of neural networks that can approximate the behavior of thousands to millions (or more) of independent variables (as of 2023, the GTP-3 chatbot uses over 175 billion independent variables).  

In computational sciences (and other disciplines) there is a concept known as the \emph{curse of dimensionality} \citep{bellman1957dynamic} that relates to the problem of describing high-dimensional data using naive methods.  In short, this concept recognizes that as the number of independent and dependent variables (the dimension) increases, the volume of the space needed to represent them increases.  This volume expansion can occur so rapidly with increasing dimension that the available data become sparse (i.e., if one represented the data in a high-dimensional matrix in the independent variables, nearly all the entries would be zero). Thus, in order to obtain a reliable result, the amount of data needed would grow exponentially with the dimensionality.  Among the properties that make neural networks useful, and in particular deep neural networks with many layers, is their ability to overcome the curse of dimensionality through function compositions.

\section{Terminology}

In the study of FNNs (and ANNs in general), terminology is an important component.  It is probably not an overstatement to say that the differences in terminology between computational sciences and applied mathematics has hampered understanding of what FNNs are (in a mathematical sense), and how they operate.  In the glossary of terms below, the terminology established is intended to be both clear and consistent with existing practice.  The definitions presented are specific to fully-connected feedforward networks (sometimes called multilayer perceptron networks).  When it comes to the presentation of FNNs in graph form, there is not currently a single unified approach. In the definitions below, a particular effort has been made to present a systematic and reproducible method to both form and interpret the graph representations of feedforward networks.  The glossary following defines many of the terms needed to describe FNNs.  The glossary is organized hierarchically rather than alphabetically because this helps organizes the information in a more practical way.

\begin{itemize}

\item \textbf{Artificial intelligence}.\indexme{artificial intelligence} There is no one accepted definition for term \emph{artificial intelligence}.  Rather than attempting to parse all possible attempts and options to date, for our purposes we can define the words are meaning the following: development of computer systems able to perform tasks normally requiring human intelligence. Examples include image recognition, generating sensible and logically correct responses to spoken input, and generation of new data (that did not exist previously) based on appropriate training (e.g., generation of images from word cues).  While this definition is somewhat circular (the notion of  "normally requiring human intelligence" is frequently revised after a new computational algorithm achieves this goal!), it is still a useful one. \\

\item \textbf{Machine learning}.\indexme{machine learning}  Much like the term artificial intelligence, the term \emph{machine learning} has a no single unambiguous definition.  For our purposes, the term will be used to indicate any algorithm (regardless of whether this algorithm is considered artificial intelligence or not) that can \emph{use a measure of the current error to automatically generate an new model that has smaller error}. In applications to ANNs, the \emph{learning loop} uses some variation of \emph{gradient descent} (defined below) to optimize the cost function residual, and \emph{backpropagation} to update the weights so that a better approximation can be represented.  While machine learning is often presented as a \emph{subset} of \emph{artificial intelligence}, that distinction is not imposed here.  For example, linear least-squares fitting is expressible as an neural network that has a known algorithm, and can be successfully used on high-dimensional data sets with many data points as the examples for input.  However, few would consider such methods to represent ``artificial intelligence".  \\

\item \textbf{Artificial neural network}.\indexme{network! artificial neural}\indexme{artificial neural network}  An artificial neural network (FNN) is a hierarchical network constructed from weighted compositions of functions.  The networks are arranged in a hierarchy called  \emph{layers}. Each layer contains a number of \emph{nodes}, where each node represents a linear or nonlinear function.  The vector of independent variables for the network is known as the \emph{input} or \emph{feature set}.  The result of the network is a vector or scalar (or other mathematical object) containing the dependent variable, and is called the \emph{output} or \emph{target}.  The input data is processed through the network by making a sequence of transformations of the input.  These transformations are done layer-wise, with the output from one layer forming the input to the next.   The input to any node in an internal layer of the network is a weighted sum of the output from the nodes from the previous layer.  The functions used to transform the variables are frequently called \emph{activation functions}, but for the problem of function fitting, might better be thought of as \emph{basis} functions.
\\

\item \textbf{(Fully connected) feedforward network (FNN)} (also called a \textbf{multilayer perceptron network or MLP}).\indexme{network!fully connected feedforward (FNN)}\indexme{artificial neural network!feedforward network (FNN)} \indexme{feedforward network (FNN)!definition}    An artificial neural network that has connections among layers in only one direction: from the input toward the output.  It is assumed that each node within a layer receives input from \emph{all} of the nodes the previous layer, and outputs information to \emph{all} of the nodes in the subsequent layer.  The feedforward structure is distinct from artificial neural networks that contain feedback loops, or parallel loops that .  This definition is also distinct from any \emph{algorithmic} process used to optimize the weights of the neural network.  Sometimes simple feedforward networks are called multilayer perceptrons, a terminology left over from work done on primarily linear networks (i.e., networks where are transformations conducted are linear ones).  The term \emph{multilayer perceptron} has come to include, however, simple feedforward networks that utilize nonlinear transformations; thus the two terms are interchangeable in modern usage.  While FNNs (MLPs) are frequently lumped in with artificial intelligence, there are good reasons not to do so.  While FNNs use nonlinear activation functions, they do not use nonlinear network structures (e.g., loops or reinforcement feedback) to find solutions.   \\

\item \textbf{Deep neural network}.\indexme{network!deep}  There is some disagreement about what constitutes a \emph{deep} neural network.  For the purposes of this work, any neural network that contains \emph{at least two hidden layers} (defined below) that invokes nonlinear transformations of the input data is a deep neural network (cf. \citep{chui1994neural}).  The reason that two layers has been selected is that \emph{nonlinearity} alone is not sufficient; in principle, even the early binary perceptron models were nonlinear because they contained a thresholding (step) activation function.  One of the key features about adding \emph{layers} to a neural network (which increases the system size linearly) is that the nonlinearity of the activation functions creates an exponentially more complex solutions to the problem.  To do something equivalent with a single-hidden-layer network would require increasing the width exponentially.   This is an important concept. Deep neural networks can be thought of the same expressive power as shallow ones, but with many fewer total nodes.  For example, it has been shown that for rectified linear unit (ReLU) activation functions, an exponentially greater number of nodes are needed for a single hidden-layer network to have the same expressive power (the ability to approximate functions) as a deep network \citep{bengio2011expressive, montufar2014number}. \\

\item \textbf{Shallow neural network}.\indexme{network! shallow}  Like many terminologies regarding neural networks, the word \emph{shallow} does not appear to have a uniform meaning.  In this work, the word \emph{shallow} will mean a ANN that contains zero or one hidden layers.  Networks with a single hidden layers are known to have good capacity for reproducing continuous functions (as measured by, say, the square of deviations, or $L_2$, metric).  However, such networks may require a number of units that increase exponentially compared with the number of units in a deep neural network (which increase only polynomially with depth for the same expressive power).    \\

\item \textbf{Network architecture}.\indexme{network!architecture}\indexme{artificial neural network!architecture}\indexme{architecture}  The network architecture is the specification of the of operations put together to form a network; it is often graphically represented.  For fully connected FNNs, the network architecture $(L, {\bf p})$ is fully specified by 
a positive integer $L$ indicating the number of hidden layers, and the vector ${\bf p} =
(p_1,...,p_L )$ indicating the width of each hidden layer.\\

\item \textbf{Node}\indexme{network!node}\indexme{node} (or sometimes \textbf{unit} or \textbf{neuron}; in graph theory sometimes \textbf{vertex}).  In the definition of FNNs above, the concept of a node was introduced.  A node at layer $k+1$ is simply an \emph{operator} (the word \emph{function} is usually used in the computational sciences literature) that transforms its input layer $k$ (see Fig.~\ref{fig:genericann}.  These transformations may be linear or nonlinear.  However, the success of FNNs in a general context for function estimation (or, equivalently, classification of data) requires that some nonlinear transformation be conducted at some of the layers of nodes.  In the material of this text, \textbf{every node in a network represents some operator} that when combined with the input, generates an \emph{activation function}.  Nodes can be \emph{input} nodes (located at the boundary where the independent variables enter the network, indicated by a square), \emph{output} nodes (located at the boundary where the prediction of the dependent variable is formed, indicated by a square), or \emph{interior} nodes (indicated by circles); these definitions should be apparent from usage.\\

\item \textbf{Node number}. \indexme{node!number}  In this work, \emph{structure preserving} notation indicates a notation in which the network layer structure is embedded in the parameter and function representations.  For this notation, each \emph{node number} is denoted by a pair indicating the layer, and the node number within that layer.  Thus, $(1,2)$ indicates ``layer 1, node number 2".  In small networks, we will use \emph{simplified notation}, where each node is given a unique integer.  While this notation does not indicate connectivity well, it does make the presentation simpler, which is ideal for small networks.\\

\item \textbf{Source node}.\indexme{node!source}  Nodes are connected by links (defined below).  A node in layer $k$ with a link whose output is directed to node $n$ in layer $k+1$ is called a \emph{source node} for node $n$.  The node $n$ in layer $k+1$ may have as many source nodes as there are total nodes in layer $k$.\\

\item \textbf{Destination node}.\indexme{node!destination}  A node in layer $k+1$ with a link whose input is directed from a source node in layer $k$.\\

\item \textbf{Bias node}.\indexme{node!bias}  Homogeneous functions are such that when each of their independent variables is zero, the returned dependent variable is also zero.  To allow for shifts in the functions defined for an FNN, one adds \emph{bias nodes} that represent the addition of a constant (bias).  Thus the magnitude of the shift is determined by its weight, $b$.  The weight variable is uniquely identified by a subscript using the convention described under the definition of  \emph{link}. See Fig.~\ref{fig:genericann} for an example of the implementation a bias term in graph form.\\

\item \textbf{Weight}.\indexme{weight}  \indexme{network!weight}  As described above, each link (or edge) of the network will have an associated weight.  There is no standardized notation for weights.  In this chapter, in \emph{structure preserving} notation, the connections among nodes will be indicated by a subscript for the source node, and a superscript for the destination node; color will be added when it can be to increase clarity.  Thus $w_{1,2}^{2,2}$ or $\wubs{1,2}{2,2}$ indicates the weight linking node $(1,2)$ (layer 1, node 2) with node $(2,2)$ (layer 2, node 2).  For \emph{simplified} notation, each node is given a unique integer index.  Thus, a weight between nodes $j$ and $k$ is specified by the simpler notation $w_{j,k}$. \\

\item \textbf{Link}\indexme{network!link}\indexme{link}(in graph theory, sometimes \textbf{edge}).  Each node in an FNN is connected to at least one other node by a link.  A link is a \emph{directed} line, usually indicated by an arrowhead, showing the direction of information flow.  Every node in a system will have at least one input link (indicated by an arrow pointing to the node) and one output link (indicated by an arrow pointing away from the node).  Links show more than simple connectivity in an FNN.  Each link in an FNN is associated with a particular \textbf{weight}, $w$.  Implicit in this organization is the concept that the input links associated with a node in layer $k+1$ define a linear combination of the output from the source nodes at layer $k$.  This linear combination is the sum of the output from each source node, multiplied by its respective weight.  Input and output links will be shown as dashed lines to distinguish them, and they can be though of as carrying a weight of unity.  \textbf{In short, every link represents a weighted term originating at a source node, and ending at a destination node.  The formation of a linear combination of the output from all source links with associated weights is \emph{implied} at every destination node.}  N.B.  Because in the \emph{simplified notation} each node in the network is given a unique integer label, each link in the network can be specified by these two integers.  For labeling the links, the convention is to list as subscripts the \emph{source node} first, and the \emph{destination node} second.  Thus $w_{1,2}$ would indicate the weight associated with the link for the \emph{output} from node 1, and directed to the \emph{input} of node 2.  If the source node involved is a bias node, then the convention is to use the symbol $b$ instead of $w$, but the convention for subscripts is the same as for regular weights.\\

\item \textbf{Layer}.  As mentioned above, and FNN is composed of nodes arranged in layers.  An example of a generic FNN is given in Fig.~\ref{fig:genericann}.  The concept of layers is a convenient organizing structure for fully-connected simple feedforward networks because node connectivity is conveniently illustrated with a layered network topology.  For networks that are not fully connected or involve feedback loops, the concept of layers becomes somewhat less important and is replaced by the concept of simple connectivity among the nodes. \\

\item \textbf{Hidden layer}.  A hidden layer is, technically, any layer that is not an input or output layer.  Thus, any network with three or more layers will contain at least one hidden layer.  Note that hidden layers are generally nonlinear (although there may be instances of deep networks with many layers that impose a layer of linear transformations for specific purposes).  All hidden layers are made up of \emph{interior} nodes.\\

\item \textbf{Features}. (also known as \textbf{inputs} or \textbf{independent variables}).  The discipline of machine learning uses terminology that is common in computational science.  In neural networks, sometimes the vector of \emph{independent variables} are frequently called \emph{features} to indicate their role as the descriptors that allow one to predict the outcome of a complex function (N.B., here ``complex" indicates degrees of freedom, not the presence of imaginary components).   As an example, observations of $N$ temperatures, $T_i$, of a resistor at discrete times, $t_i$, would generate $N$ pairs of the form $(t_i, T_i)$.  The $N$ values of $t_i$ form the features for the data set.  In machine learning, the feature sets can often become quite large.  The fact that neural networks (coupled with appropriate optimization and weight-adjustment algorithms) are able to generate optimized fits to data sets with such large numbers of features is one of the primary motivations for using them. \\

\item \textbf{Target}. (also known as \textbf{outputs} or \textbf{dependent variables}).  The collection of dependent variables.  Often the target is a scalar value that is specified by a real number associated with a \emph{specific} vector of the independent variables.  As an example, observations of $N$ temperatures, $T_i$, of a resistor at discrete times, $t_i$, would generate $N$ pairs of the form $(t_i, T_i)$.  The values $T_i$ form the target values for the single independent variable (feature) $t_i$.  In principle, one may have targets that are specified by a vector or tensor (matrix) output. By convention, the  target values predicted by a neural network is indicated by typesetting a circumflex (hat) over the output variable to indicate that the resulting function is an estimator.\\

            \begin{figure}[t!]
            \sidecaption[t]
            \centering
            \includegraphics[scale=0.4]{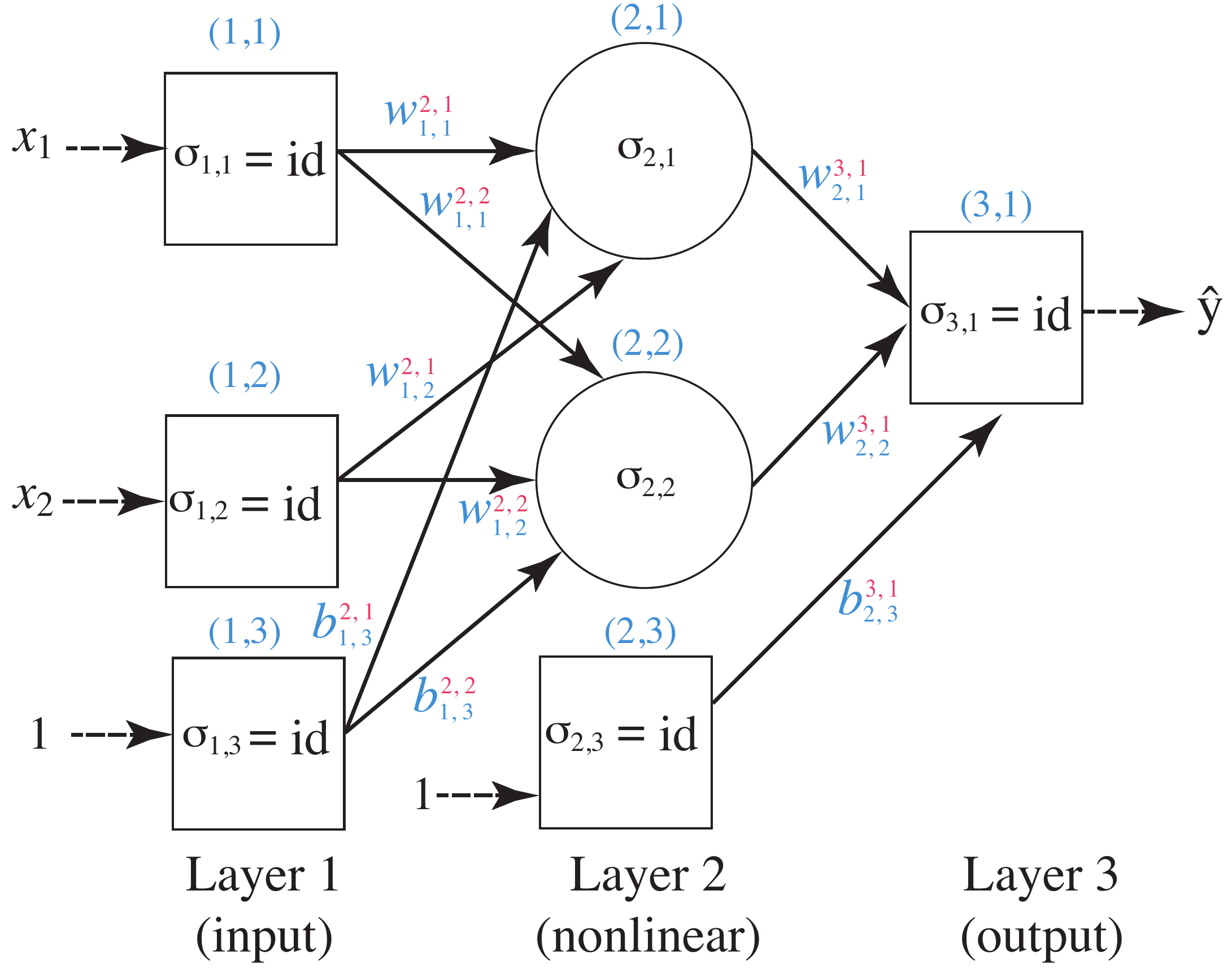}
            \vspace{-5mm}
            \caption{A generic feedforward artificial neural network, illustrated in graph form with structure-preserving notation.  In this example, there are two independent (feature) variables, $x_1$ and $x_2$, and one dependent (target) variable, $y$.  The example contains three layers, and one of these is a \emph{hidden layer} (i.e., it is neither the first nor the last layer of the network).}
            \label{fig:genericann}       
            \end{figure}
%

\item \textbf{Activation function} (or \textbf{basis function}).  Each node in a neural network represents an \emph{operator} that transforms the input.  By convention, these operators are called activation functions.  The technical difference between an \emph{operator} and a \emph{function} is not likely to cause confusion, so we will tend to stay with the convention.   This function may be linear, or nonlinear.  Because compositions of linear functions map back to a (new) linear function, an essential feature in modern neural networks with one or more hidden layers is that the activation functions are \emph{nonlinear}.  It is not an overstatement to say that the use of such nonlinear activation functions is one of the reasons that FNNs have become useful for function approximation.  While the term \emph{activation function} is a popular one, in some applications it may be better to think of these functions as being \emph{basis} functions, from which new functions can be derived by rescaling, shifting, and linearly combining the basis functions. \\

Several examples of activation functions are given in Fig.~\ref{fig:activations}.  The selection of an activation function depends on the particular application (e.g., if the result needs to have continuous derivatives; if the networks is very deep).  Although it does not have continuous derivatives, the rectified linear unit ($relu$) activation function has become the most widely used function in deep neural networks.  These activation functions will be the primary ones discussed in this chapter, and are described in additional detail in the material following. \\

            \begin{figure}[t!]
            \sidecaption[t]
           \centering
            \includegraphics[scale=0.45]{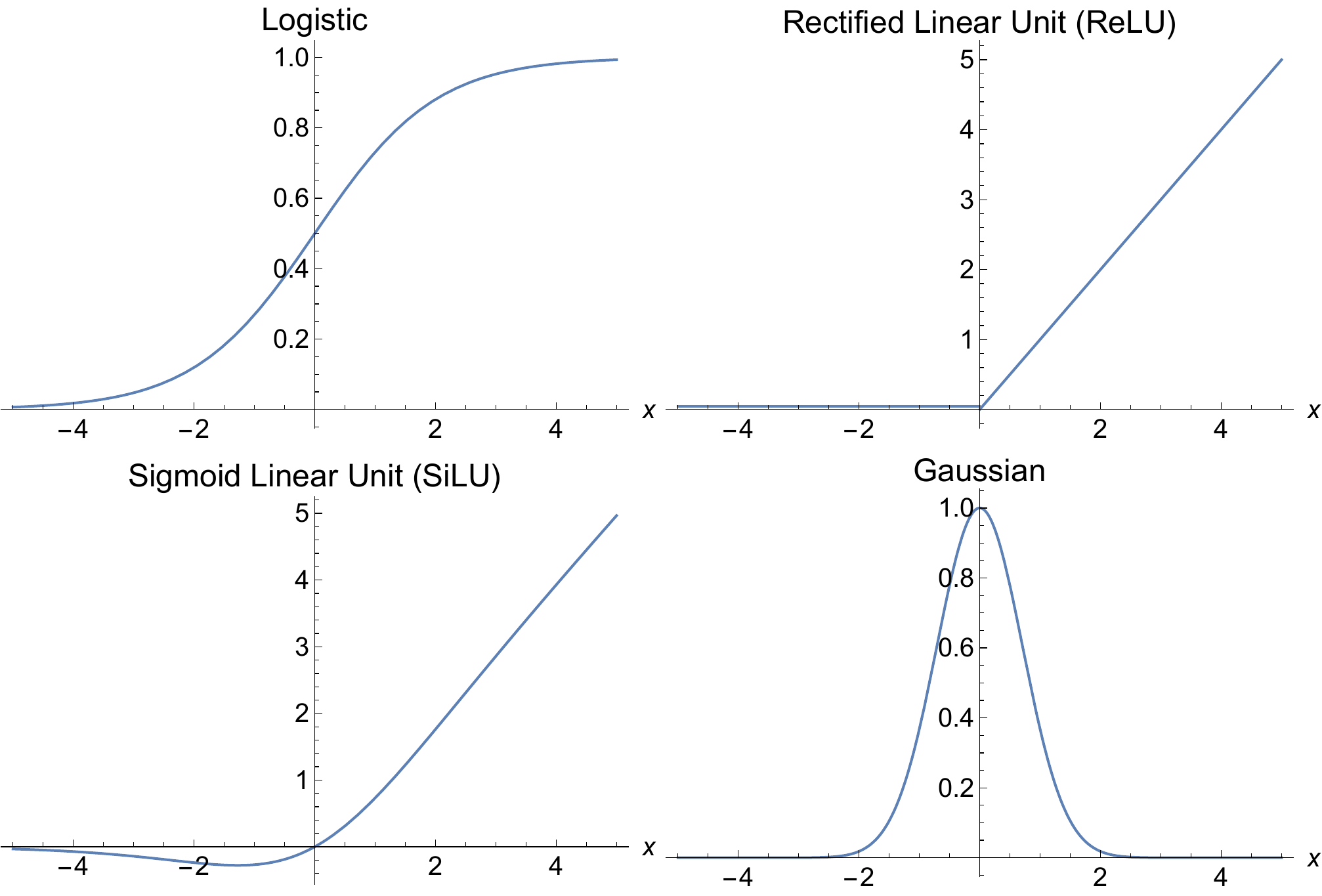}
            \vspace{-5mm}
            \caption{Examples of four activation functions.  Activation functions a selected on the basis of network purpose, and the depth of the network.  Functions whose derivative approaches zero as $|x|\rightarrow \infty$ often do not perform well in deep networks, because the (numerical) propagation of information back through the network is proportional to the derivatives of the function.}
            \label{fig:activations}       
            \end{figure}

\item \textbf{Loss function}.\indexme{network!loss function}  In order to optimize a network, some differentiable metric of the objective (goodness-of-fit for regression; appropriate data separation for classification) is required.  For analytic purposes, all metrics of goodness of fit (the absolute value of the error, $L_1$ norm; or the sum of squared errors, $L_2$ norm) for a finite-dimensional vector space are all equivalent.  However, in practice, computations may be facilitated by using one error metric over another.  For our purposes, we will adopt the sum-of-squared-errors as the norm exclusively.  Some of our examples will be on infinite-dimensional (functional) vectors spaces where different norms are not necessarily equivalent.  This will not create problems that are noteworthy.\\

\item \textbf{Identity operator}.  There is occasionally use for an identity operator as one of the possible operators (functions), $\sigma$, that can be used at a node.  For our purposes, the identity operator is defined by

\begin{equation*}
    \textrm{id}(x) = x
\end{equation*}
That is, the operator leaves the input unchanged.  \\

\item \textbf{Optimization}.\indexme{network!optimization}  The process by which one determines the set of parameters (treated as independent variables) that returns the minimum of an objective function formulated with specified error metric.  The error metric is usually the absolute difference between a model and the data (the $L_1$ metric) or the least-squares difference (the $L_2$ metric), although others are used.  Optimization can be done by analytic or numerical means.  For a finite-dimensional problem (e.g., minimizing a function based on a finite number of data), then in principle all error metrics should yield equivalent results.  In practice, because of the approximation methods used in numerical schemes, sometimes one error metric may provide superior properties compared with another.  Sometimes the optimization process is indicated by specifying an ``$\textrm{arg min}$" operator.  Thus one might write $\textrm{arg min}[f(x; w_1, w_2)]$ to return the values $w_1$ and $w_2$ such that the specified error were minimized.\\

\item \textbf{Mulitiresolution}.  The word \emph{multiresolution} came into common use when the theory of wavelets became popular, and at some point became synonomous with wavelets.  However, it has since taken on a more general meaning.  In this chapter, \emph{multiresolution} means only that \emph{local refinements} can be made to a function without disturbing the remainder of the function.  This is a recurring theme in the analysis of ANNs, so this generalized definition of multiresolution is a useful one.\\

\item \textbf{Backpropagation}.\indexme{network!backpropagation}  Backpropagation is the \emph{numerical scheme} used to adjust the weights when optimizing an ANN in response to the current state of the error in the solution.  This method relies on an application of the chain rule to assign adjustments to the weights.  Backpropagation methods were discovered in the early 1970s, but it took some time for the algorithms to be assimilated by the machine learning community.  Because this is a scheme used for solving ANNs numerically, it will not be a component of the material in this chapter.  However, the scheme is so central to the development and solution methods for ANNs, it is necessary to be aware of the scheme and its purpose. \\

\item \textbf{Supervised machine learning (in the context of neural networks)}.  Suppose we have $N$ examples of data in the form $(\bf{x}_i,y_i)$, where $\bf{x}_i$ is a vector containing the independent variables associated with observation $y_i$.  The problem  of parameterizing an FNN to conduct regression or classification using data sets with known independent and dependent variables is called \emph{supervised learning}.  Here, the word supervised indicates that the example data are explicitly labeled with their appropriate independent variables.   Large neural networks generally cannot be (or, would not be practical to be) solved using strictly analytical methods.  Instead, such networks are solved computationally.  The \emph{primary} significant difference in solving simple feedforward neural networks via computation (as opposed to analytically) is the method used to minimize the loss function.  In machine learning for FNNs, the method of optimization is usually done via one of several versions of \emph{gradient descent} methods coupled with a method of \emph{backpropagation} which is a scheme for accurately updating the weights in the network.  Because only analytically tractable networks are introduced here, neither gradient descent nor backpropagation methods are discussed further.

\end{itemize}

\section{Interpretation of the Graphical Representation of FNNs}\indexme{network!graphs}

It is helpful to outline the way that ANNs are described graphically before discussing them further; this way we have the appropriate ``structural vocabulary" to better understand examples. 
As we know, the basic operation of a neural network is to take input, transform the input using a linear or nonlinear operator, and then create output.  This process is done at nodes in a sequence of nonlinear operations arranged in layers (as defined above).  Between each layer, the output from the nodes of the previous layer is transformed by creating a linear composition of those outputs.  To make the process more transparent, these operations are frequently presented as connected graphs.  While there appears to be little uniformity in the literature, the following rules are adopted in this text to describe ANNs graphically.

\begin{enumerate}

    \item Dashed lines indicate input or output.  The direction of the associated arrow allows one to determine whether it represents input or output.

    \item In the notation adopted here, all nodes contain activation functions, are indicated by the symbol $\sigma$ with a subscript indicating its node number.  Node numbers are specified by the ordered pair $(i,j)$, where $i$ represents the layer index, and $j$ indicates the node index within that layer.  Frequently, only one form for the activation function is used, in which case the activation function does not need to be indexed. The activation functions can be linear or nonlinear (particular forms will be discussed later).  As a special case for input and output nodes, the operator is the identity operator $\sigma(\cdot)$\raisebox{2mm}{$\frac{\textrm{def}}{\Large=}$} $\textrm{id}(\cdot)$, which returns as output the value of the input. 

    \item Squares represent input or output nodes.  These are also called \emph{boundary} nodes, indicating that they are nodes that communicate external input to or output from the network.  There is always at least one input and one output layer; the output layer may not perform any nonlinear transformations on the data.

    \item A circle on the graph represent a node with accompanying activation function; the activation function can be linear or nonlinear.  These operators transform the weighted input coming into them to generate a single output value.  All nodes indicated with circles are \emph{interior} nodes, indicating that they are not an input or output boundary of the graph.

    \item Solid lines directed between nodes (with an arrow indicating direction) represent a weighted output.\indexme{network!structure-preserving notation}  The output of the node at the tail of the line is weighted by a weight value, $\wubs{i,j}{k,l}$, given above the directed line.  In this notation, $i,j$ represents a node in layer $k-1$, and $k,l$ represents a node in layer $k$.  This is known as \emph{structure-preserving notation}.

    \item Each internal node will have one or more weighted inputs directed toward them.  The input to every internal node is a sum of the weighted inputs pointing to that node.
    
    \item Bias terms are inputs that may occur at any layer of the network. Such nodes take unity as the input, transform the input with the identity, and output the number ``1".  This output is then weighted by a weight value, $\wubs{i,j}{k,l}$, given above the directed line.  Frequently, bias terms are added to each internal node, and these are numbered the same way as other nodes. Thus, $\bubs{i,j}{k,l}$ represents a bias term for the link between node $(i,j)$ and node $(k,l)$ Because the bias nodes can complicate the graph, it is not uncommon to simply note that each internal node contains a bias term without explicitly representing them on the graph.
    
\end{enumerate}
%

            \begin{figure}[t!]
            \sidecaption[t]
           \centering
            \includegraphics[scale=0.4]{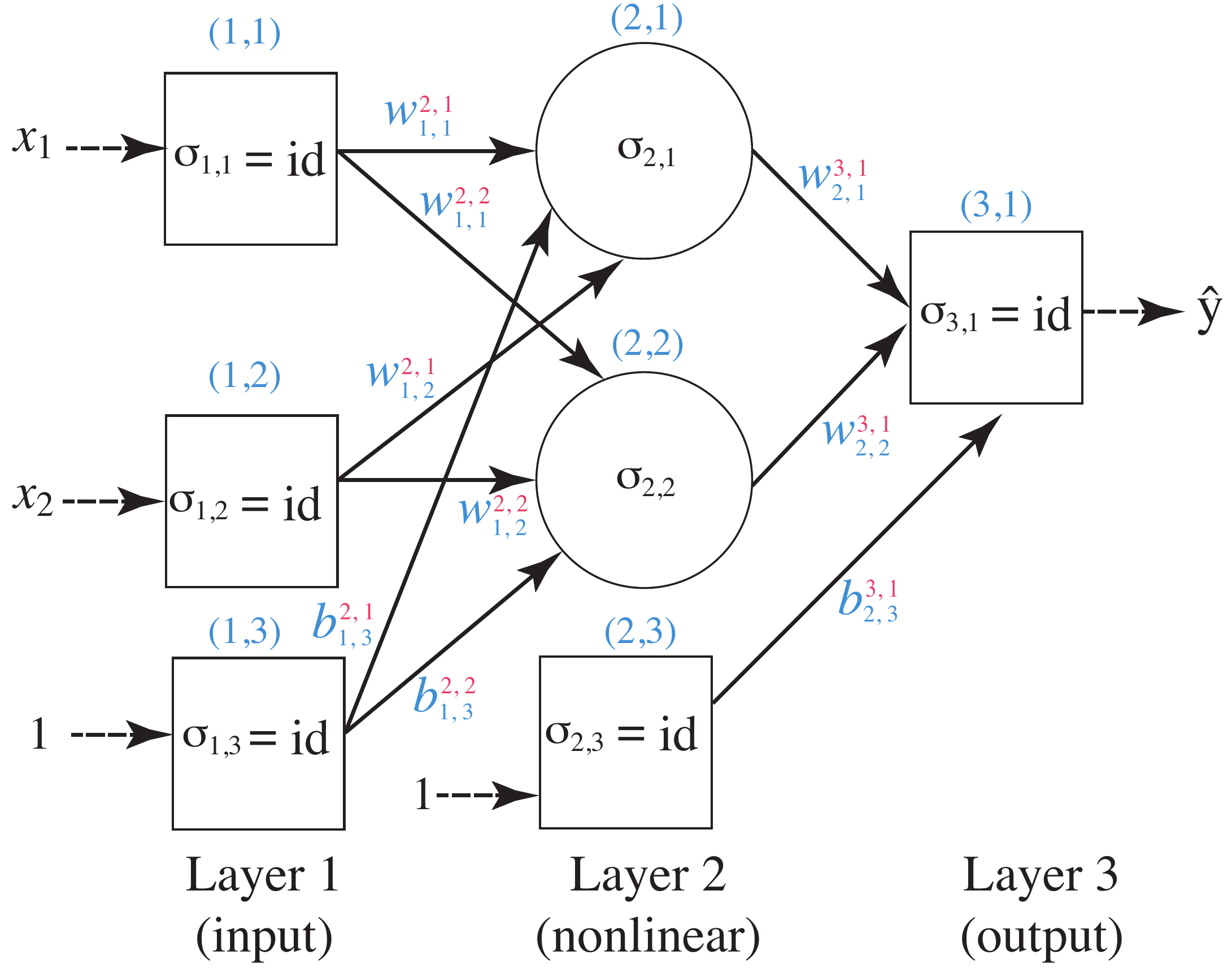}
            \caption{The feedforward artificial neural network from Fig.~\ref{fig:genericann}, repeated here for convenience.  This network is specified with structure-preserving notation.}
            \label{genericann2}       
            \end{figure}

To complete the connection between the graph and the mathematical representations, we consider the example of the network shown in Fig.~\ref{fig:genericann}, and repeated here as Fig.~\ref{genericann2} for convenience.  The output from any layer will is usually indicated by $\mathcal{N}$ (for ``network"), subscripted by its layer index.  
Starting with the first layer, we note that the linear composition rule indicated by the arrows connecting the first and second layer would provide the following \emph{output} from layer 1

\begin{align}
    \mathscr{N}_1(x_1,x_2,1) &= \begin{pmatrix}
   \wubs{1,1}{2,1}    & ~~\wubs{1,2}{2,1} \\
    \wubs{1,1}{2,2}    & ~~\wubs{1,2}{2,2} 
    \end{pmatrix}
    \begin{pmatrix}
        x_1\\
        x_2
    \end{pmatrix}
    +
    \begin{pmatrix}
       \bubs{1,3}{2,1}\\
        \bubs{1,3}{2,2}
    \end{pmatrix}\nonumber\\
    &=\begin{pmatrix}
    \wubs{1,1}{2,1} x_1 + \wubs{1,2}{2,1} x_2 +  \bubs{1,3}{2,1} \\
     \wubs{1,1}{2,2} x_1 + \wubs{1,2}{2,2} x_2 + \bubs{1,3}{2,2}
    \end{pmatrix}
\end{align}
The second layer takes the vector $\mathcal{N}_1$ as input, and forms the following output (as indicated by the arrows between layers 2 and 3)

\begin{align}
    \mathcal{N}_2(x_1,x_2,1)= \begin{pmatrix}
        \sigma_{2,1}\left(\wubs{1,1}{2,1} x_1 + \wubs{1,2}{2,1} x_2 +  \bubs{1,3}{2,1}\right) \\
        \sigma_{2,2}\left( \wubs{1,1}{2,2} x_1 + \wubs{1,2}{2,2} x_2 + \bubs{1,3}{2,2}\right)
        \end{pmatrix}
        \cdot
        \begin{pmatrix}
            \wubs{2,1}{3,1}\\
           \wubs{2,2}{3,1}
        \end{pmatrix}
        + \bubs{2,3}{3,1}
\end{align}
Note that this output is now scalar valued.  Again, this is very typical of ANNs.  A vector-valued input is transformed by a sequence of compositional vector operations.  In the last step, a contraction reduces the dimensionality of the output (usually) to a scalar, although other output forms (e.g., vector or tensor outputs) are possible.

Finally, the third layer operates on the output $\mathcal{N}_2$ with the identity.  The result is

\begin{align}
   \hat{y}(x_1,x_2,1)&= \begin{pmatrix}
        \sigma_{2,1}\left(\wubs{1,1}{2,1} x_1 + \wubs{1,2}{2,1} x_2 +  \bubs{1,3}{2,1}\right) \\
        \sigma_{2,2}\left( \wubs{1,1}{2,2} x_1 + \wubs{1,2}{2,2} x_2 + \bubs{1,3}{2,2}\right)
        \end{pmatrix}
        \cdot
        \begin{pmatrix}
            \wubs{2,1}{3,1}\\
           \wubs{2,2}{3,1}
        \end{pmatrix}
        + \bubs{2,3}{3,1}\nonumber\\
        &=  \wubs{2,1}{3,1}\sigma_{2,1}\left(\wubs{1,1}{2,1} x_1 + \wubs{1,2}{2,1} x_2 +  \bubs{1,3}{2,1}\right) +
        \wubs{2,2}{3,1} \sigma_{2,2}\left( \wubs{1,1}{2,2} x_1 + \wubs{1,2}{2,2} x_2 + \bubs{1,3}{2,2}\right)
        + \bubs{2,3}{3,1}
\end{align}
While in this example the application of the identity in layer 3 is a bit tautological, it does serve a specific purpose.  Viewed this way, the nodes of \emph{every} layer conduct \emph{some} operation.  It is not uncommon to use a filtering function as the final step in the network, in which case the final activation function would not be the identity.  

Even a small network can rapidly lead to long and complex compositional functions.  The utility of network representations begins to become quite apparent for networks of any substantial size.

\subsection{More on Notation}\indexme{network!simplified notation}

The notation above is useful because the network structure is embedded in the notation.  Thus, one can generate the mathematical form of the network from the graph; one can also easily generate the graph from the mathematical representation because the structure of the graph is fully represented in the notation.  

While this kind of notation can be useful (e.g., in coding, or with networks that are not strictly feedforward), it is not ideal for examining small networks such as we will use.  The alternative is to abandon the labeling of layers.  This means that each node receives a single integer index, and each weight requires only two integers to indicate its connection.  The downside is that the graph must accompany the mathematical representation or there can be a lack of clarity about the location of nodes (although with strictly feedforward networks this is less of a problem).   

            \begin{figure}[t!]
            \sidecaption[t]
           \centering
            \includegraphics[scale=0.4]{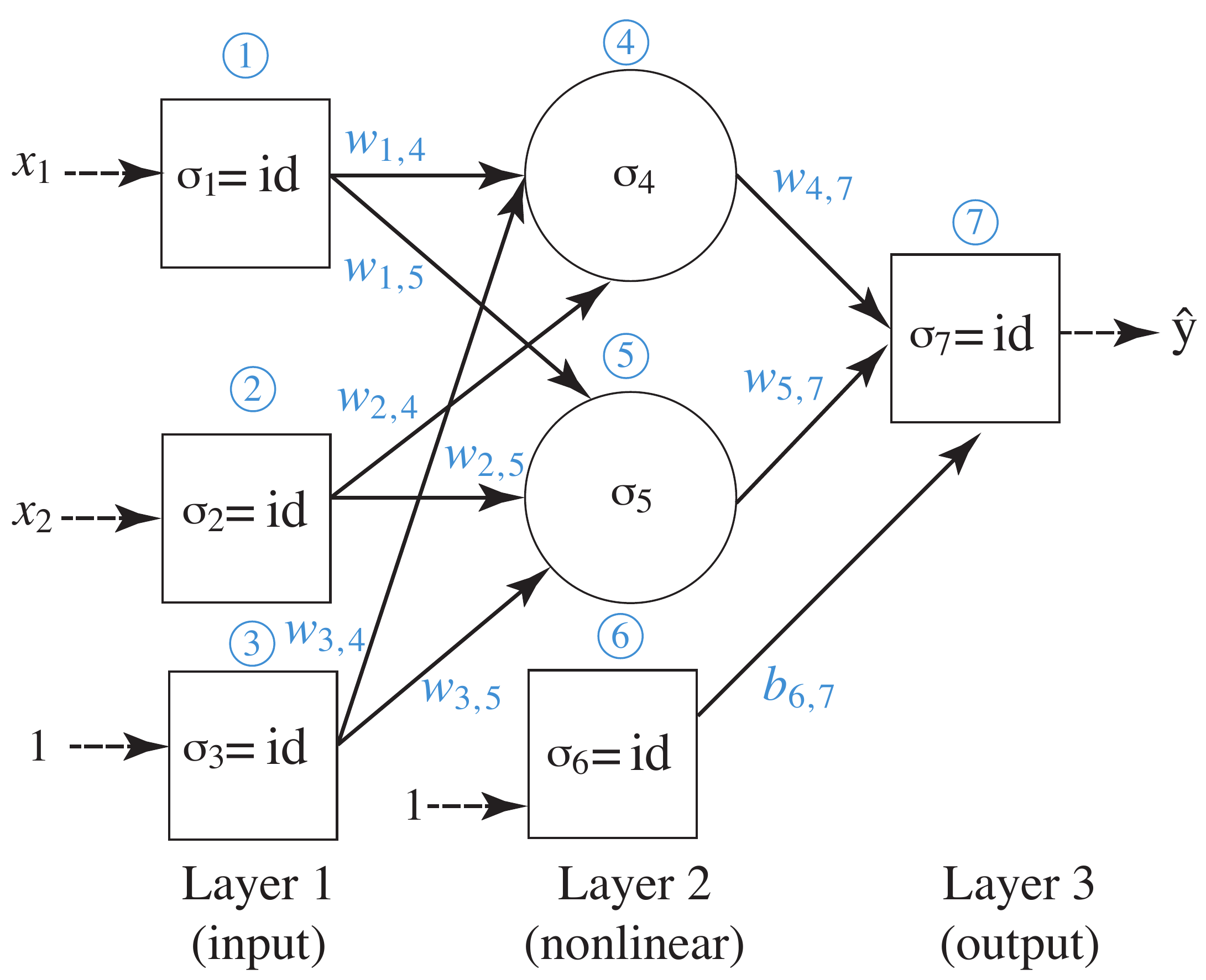}
            \caption{The feedforward artificial neural network from Fig.~\ref{fig:genericann}, repeated here for convenience.  This network is specified with simplified notation that does not preserve structure.}
            \label{genericann_simple}       
            \end{figure}

An example of the \emph{simplified} notation is given in Fig.~\ref{genericann_simple}.  For small networks, this notation is substantially easier to understand, and it will be the default option for the remainder of this chapter.  For comparison, the analysis above can be repeated.  The result is a bit easier to read, and is probably better for initial forays into understanding ANNs.  

\begin{align}
   \hat{y}(x_1,x_2,1)&= \begin{pmatrix}
        \sigma_{4}\left(w_{1,4} x_1 + w_{2,4} x_2 +  b_{3,4}\right) \\
        \sigma_{5}\left( w_{1,5} x_1 + w_{2,5} x_2 + b_{3,5}\right)
        \end{pmatrix}
        \cdot
        \begin{pmatrix}
            w_{4,7}\\
           w_{5,7}
        \end{pmatrix}
        + b_{6,7}\nonumber\\
        &=   w_{4,7}\sigma_{4}\left(w_{1,4} x_1 + w_{2,4} x_2 +  b_{3,4}\right) +
         w_{5,7} \sigma_{5}\left( w_{1,5} x_1 + w_{2,5} x_2 + b_{3,5}\right)
        + b_{6,7}
\end{align}

\section{History and Purpose of ANNs}\indexme{artificial neural network!history}

The development of modern artificial neural networks provides an interesting case history that provides and example of how scientific advances actually unfold.  While ANN-like networks had existed in the mathematical literature for some time, Warren McCulloch and Walter Pitts \citep{mccullough1943logical} published what is thought to be the first paper that linked the concept of \emph{cognition} with that of computation.  Their model was essentialy a binary one, where the sum of binary inputs are filtered by a threshold (or Heaviside) function; a binary output of 1 indicates that the sum is over the threshold (or, 1 indicates ``true"), and a 0 indicates that the sum is below the threshold (or, 0 indicatese ``false").  A history and analysisi of their work is available in the paper by \citet{piccinini2004first}.  While their work is no longer considered an archetype for cognition, it nonetheless started the era of bio-inspired methods for computation. 
The theory developed by McCulloch and Pitts represents the first foray into what have become modern neural networks.  However, their model did not contain a \emph{learning} step.  In other words, it was simply positioned at some state (i.e., the threshold was set to some specified value), and then the network responded to the inputs.

\begin{figure}[t!]
            \sidecaption[t]
           \centering
            \includegraphics[scale=0.4]{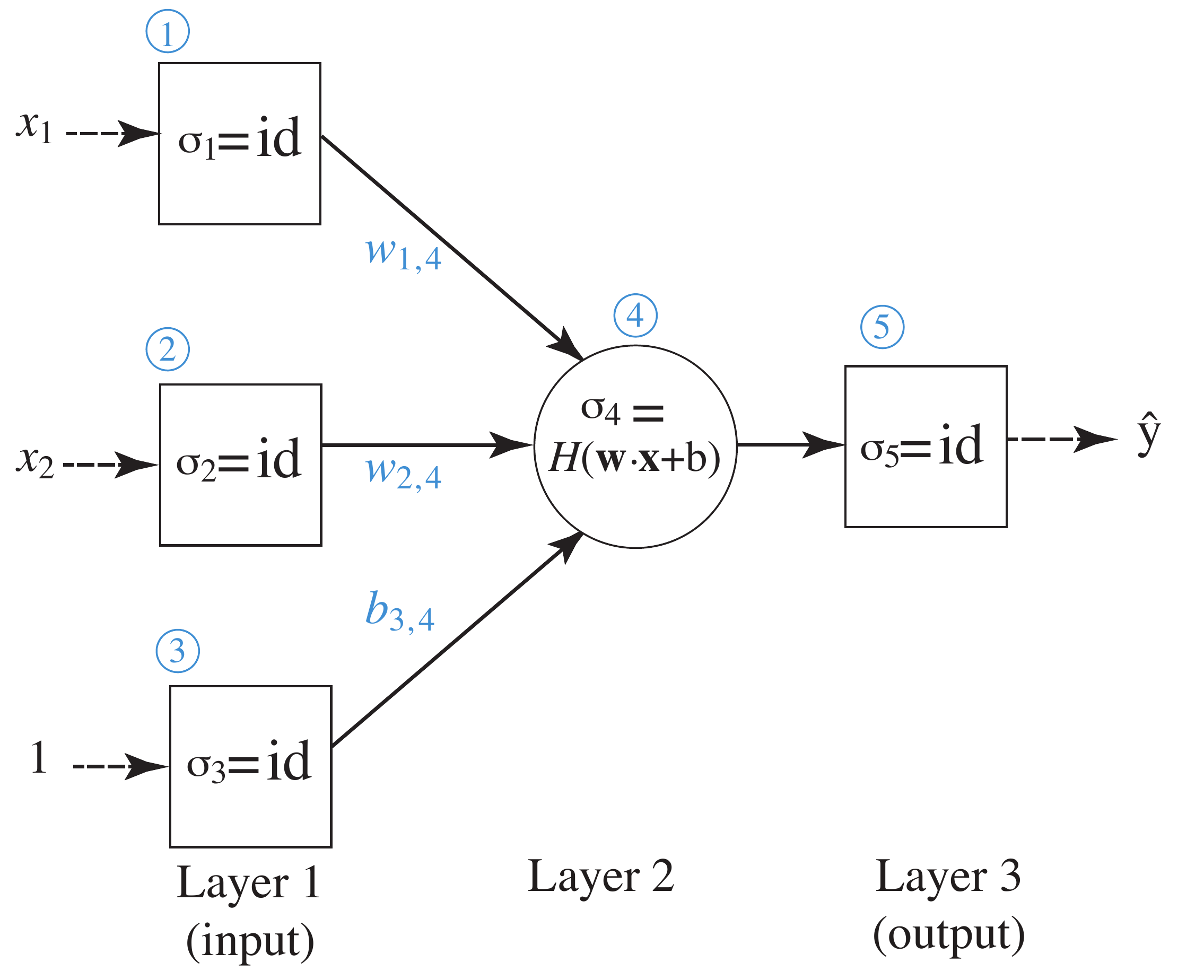}
            \vspace{-5mm}
            \caption{The classical perceptron model of Rosenblatt.}
            \label{fig:classicpercep}       
            \end{figure}

In the late 1950s, a psychologist named Fran Rosenblatt picked up the mantle of machine learning by proposing a network with a feedback step for adjusting the weights on the basis of the output error.  The basic unit of computation was called the \emph{perceptron}.   A graphical representation of the classical perceptron appears in Fig.~\ref{fig:classicpercep}.  In short, the perceptron is a binary classifier that takes a vector of input values ${\bf x}=(x_1, x_2, \ldots,1)$, and constructs a linear combination.  This combination is then sent through a Heaviside function which acts as a thresholding unit.  The output is either a 1 or 0 depending on the input function.  Critically, there was a simple algorithm that was able to adjust the weights in a systematic way so that if a solution were possible, the perceptron would converge to it.  Unfortunately, if classification was not possible (specifically, if the set to be classified was not separable using a linear function), then the learning method would never converge.  The limitations on the perceptron were that (1) it could learn only linear relationships in the data, and (2) the learning algorithm relied on a binary error.  This requires that the classification be error-free; that is, the data must be completely separable by a linear function.  Even a single point that prevents linear separation would be sufficient for the method to never converge.

\begin{figure}[t!]
           \centering
            \includegraphics[scale=0.9]{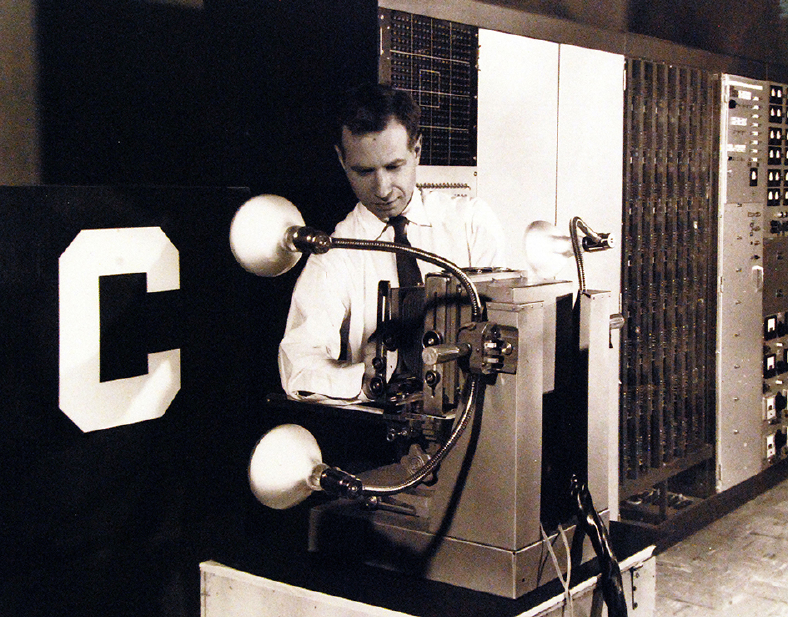} 
            \hspace{5mm}\includegraphics[scale=0.9]{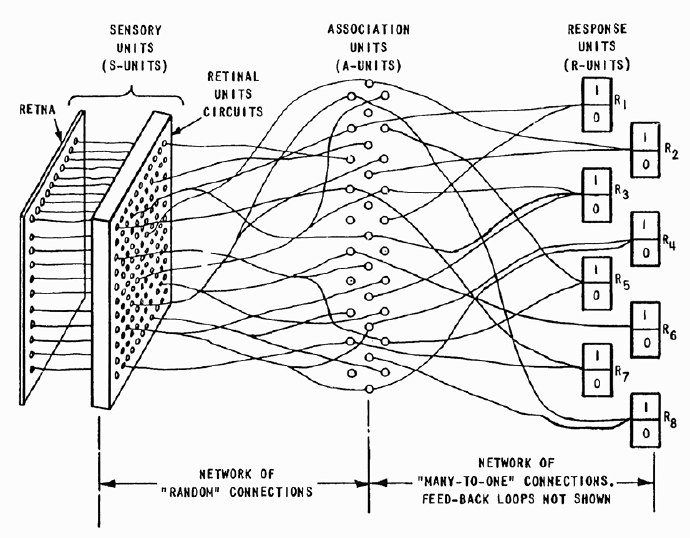}
            \caption{(Left).  Roseanblatt working with the Mark I Perceptron, circa 1960.  (Right) A schematic of network comprising the Mark I Perceptron from the operator's manual  (adapted from \citep{hay1960mark}).   Photograph by the National Museum of the U.S. Navy, licensed under the Public Domain Mark 1.0.}
            \label{fig:perceptimages}       
            \end{figure}

The work by Rosenblatt continued with the first successful construction of a bio-inspired computer. The Mark I perceptron was electro-mechanical device developed during 1957 and 1958 by Frank Rosenblatt, Charles Wightman, and others.  The primary purpose of the Mark I perceptron was to conduct pattern recognition for symbols.  Weights were encoded in current potentiometers, and weight updates during learning were performed by automatic adjustment of the potentiometers by electric motors.  A picture of Resenblatt with the Mark I and a network diagram of the apparatus appears in Fig.~\ref{fig:perceptimages}.

The next significant step in the theory for learning algorithms was developed by Bernard Widrow and his doctoral student Ted Hoff at Stanford university in 1960.  They developed a network using an improved learning rule; the method was dubbed the adaptive linear element (ADALINE) rule or the Widrow-Hoff rule.   
The main difference in this new approach was how the feedback error was defined and used to adapt the weights and bias of the two perceptrons. While Rosenblatt used the classification error (either a 1 or 0) to as the learning metric, ADALINE introduced the concept of a loss function (or objective function).  This allowed non-integer errors to be computed, which was necessary to allow weight modification to be proportional to the error metric.  Widrow and Hoff also developed an early form of a gradient descent algorithm known as the \emph{delta learning rule}.  The details of this algorithm are available elsewhere; however the important points were that it was a simple algorithm that allowed weight modifcations to be proportional to a non-integer error metric.

While there were many good ideas proposed in this early model, there was a perfect storm of occurrences that led to a dramatic decrease in interest in continued research on ANNs.  First, the original perceptron was dramatically oversold by the media and the scientists involved.  While there were early successes, the somewhat hyperbolic promises of a new world of ``thinking machines" never materialized.  Second, possibly in part as a response to the hype regarding the perceptron, Marvin Minsky and Seymour Papert of MIT release an unpublished technical manuscript (later, after revision, to be released as a book \citet{minsky1969perceptrons,minsky1988perceptrons}) illustrating that it was impossible for these linear, single-layer perceptron networks to learn simple \emph{nonlinear} binary functions; the ``exclusive or" function is the conventional example .  While the intent of Minsky and Papert did advance the field, it also resulted was ultimately a cooling of interest in artificial neural networks in general.  It did not help that Minsky and Papert were uncertain that multilayer perceptrons, and they prognosticated that they  would be unlikely to rectify the problem.   Simultaneously, Rosenblatt claimed that the generalization of perceptron capabilities could be improved considerably by including multiple layers.  However, with no know training method available, these technical difficulties prevented any direct implementations of multi-layer perceptrons.  The dynamics of those who supported versus those who disparaged single layer perceptrons (and the idea that multilayer perceptrons were not useful) fell fairly distinctly into two camps.  This is an example of a somewhat pathological diversion that can occur in science for somewhat non-scientific reasons.  While such pathological paths in science are generally rare (and cannot last indefinitely), it can be useful to understand why they occur so that they might be prevented in the future. This event in particular been examined in some detail in the fascinating analysis by \citet{olazaran1996sociological}.  

In the 1970s through the mid 1980s, work on ANNs was being conducted by only a handful of dedicated researchers.  However, during this time various forms of the backpropagation algorithm were generated; the invention seems to have happened simultaneously by several researchers \citep{linnainmaa1970representation, werbos1974beyond, lecun1985procedure,rumelhart1986learning}.  This algorithm was a crucial step in the modern development of ANNs.   In addition, some researchers continued to explore the use of neural networks with more than one hidden layer, despite the concerns levied by Minsky and Papert.  

In the early 1980s, \citet{hopfield1982neural,hopfield1984neurons} published several papers that helps restart the work on ANNs, and ushered what has more-or-less become the second revolution in deep learning.  As a final impetus for the reemergence of research on ANNs, Finally in 1986 and 1987, James McClelland and David Rumelhart published a set of textbooks that summarized much of what was known to date \citep{mcclelland1986parallel,mcclelland1987parallel}. In particular, these new texts published the backpropagation
learning algorithm, which spread knowledge of the algorithm widely.  These advances led to the birth of FNNs (also called multilayer perceptrons, MLPs), and a general revitalization of research on the topic of artificial neural networks.

\section{Some Principles of Feedforward Neural Networks}\indexme{feedforward network (FNN)!principles}
Before delving into the construction and solution for ANNs, it is useful to outline why one should study them at all.  While ANNs have been researched for many years, it has been only in more contemporary times that the computational algorithms and hardware that allow the solution to deep neural networks to become routine.  Recall, deep neural networks (DNNs) are neural networks that have more than one hidden layer.  The reason that the number of hidden layers is important is because the methods to optimize networks with more than one hidden layer is a relatively recent innovation.  However, once research on DNNs became practical, it was discovered that they were able to compute solutions \emph{substantially} more efficiently than previously known methods.

The reason that ANNs in general are so effective at fitting or classifying high-dimensional data is an area of intense research, both from the perspective of generating descriptive theory, and more empirical efforts to produce efficient algorithms.  While it has been somewhat common to characterize ANNs as \emph{black boxes}, this is not an entirely correct characterization.  While a complete theory of ANNs generally does not yet exist, substantial progress has been made to develop theory for many specific architectures.  The literature on this topic is enormous, and the mathematical methods that have been used involve very specific technical methods from function theory, topology, classical approximation theory, and statistics.  The large array of approaches that have been used is perhaps one of the best indicators that the topic is an interesting one.  While there is much yet to learn about ANNs, the progress on understanding the theory describing them has been significant in the past few years, and there is little evidence that such progress will not continue.

The following represent general principles that can be gleaned (as of 2023) from the literature regarding ANNs generally, and DNNs in particular.  Here the word ``principles" does not imply ``axiom", but is the colloquial understanding of the word (i.e., these represent summary observations about the current state of knowledge).  Presented in no particular order, the principles below are based on both the qualitative and quantitative knowledge regarding ANNs to date.  

\begin{enumerate}

\item \textbf{Artificial neural networks are unreasonably effective}.  One finds in the literature comments regarding the \emph{unreasonable} effectiveness of deep neural networks.  The purpose of the word ``unreasonable" here is to underscore that, until recently, methods that are as efficient as ANNs are has not been known to exist.  For example, the theory of approximating nonlinear functions with piecewise-linear splines is well understood, but comparatively inefficient.  A full understanding of \emph{why} ANNs can accomplish this has only recently begun to be understood.  To date, partial explanations appear to involve (a) the compositional way that activation functions are used (in deep neural networks) to reducing the size of the parameter space, (b) the self-similar nature of such compositional structures (and, presumably, the data that they represent), and (c) the existence of reduced-dimensional representations that capture ``hidden" structure within the data \citep{recanatesi2019dimensionality}.  Because most data sets include some kind of (perhaps unrecognized) structure, a number of researchers are pursuing the idea that compositional functions are particularly well suited to modeling such data.  The practical result is that if a basis function set is well-suited for representing certain kinds of data structures, then it might do so with fewer adjustable parameters.  As an analogy here, we can think of methods to represent the function $f(x)=\sin(\pi x)$ on $x\in[0,1]$.  While a Taylor series would represent one possible set of polynomial functions that approximate $f(x)$, it may take many terms to decrease the approximation error to the desired level.  However, a Fourier series converges to this function exactly with only a single term!  A paper by \citet{daubechies2022nonlinear} outlines some of the current thinking in regard to \emph{why} deep neural networks are so effective.

\item \textbf{Universal approximation theorems exist}.  In the early 1990s, a number of researchers began to explore the question of the \emph{kind of functions} that could be represented by ANNs.  This concept is sometimes called the \emph{expressiveness} of a neural network.  Most of the early work (e.g., \citet{cybenko1989approximation}) used classical approximation theory (here, we can think of this as a more general theory that incorporates the Fourier series as a special case) to describe what the capabilities of ANNs could be.  Others have pointed out that the work of Schwartz in the 1940s provides a basis for understanding the representation problem.   More recent efforts have focused on two themes: (1) generating results for the kinds of functions (e.g., the $\relu$ function) that are actually used in practice, (2) generating estimates of the accuracy of the approximations as a function of network size, and (3) accounting for the unique functional structures induced by deep networks (e.g., \citet{petersen2018optimal}).  These most recent results suggest that not only can $\relu$ feedforward networks approximate (the are \emph{dense in}) Lebesgue integrable functions (using the integral of the square-error ($L_2$) metric), but the convergence rates are improved when one uses \emph{deep} rather than shallow networks.  It is important to note, however, that approximation theorems for ANNs has existed since the 1980s, and do not depend on the depth of the network.  While the \emph{ability} to approximate very general functions is clearly important (ANNs would not be useful if they could not approximate wide classes of functions!), their unusual \emph{effectiveness} is not inherently tied to their ability to approximate functions.

\item \textbf{Multiresolution representations are involved}.  When we studied Fourier series, we learned that they are a representation that is global over the domain.  That is, each weighted member of the series contributes in to the solution in a periodic way that is the same throughout the domain.  Wavelets, on the other hand, are functions that have peaked non-periodic components (and are sometimes even compact functions) that can selectively emphasize features locally (in some subset of the domain).  Later on, we will see some examples where particular activation functions can be made to represent both local features (e.g., peaks) and global features (e.g., components easily represented by, say, a Fourier seires) of functions.   Thus, it appears that ANNs, at least in principle, can exhibit \emph{multiresolution} behavior.  Here, \emph{multiresolution} means only that \emph{local refinements} can be made to a function without disturbing the remainder of the function. This was recognized quite some time ago in the literature \citet{chui1994neural}, but only recently became a topic of broader interest in the community studying ANNs.  One might expect (or hope?) that the optimization process would select activation functions that represent the appropriate local or global representations depending upon which were most efficient.  Currently, it is not well understood how local versus global representations are internally parsed in applications \citep{raghu2021vision}.  Again, the paper by \citet{daubechies2022nonlinear} summarizes much of the current understanding of the multiresolution behavior expressed by ANNs.

\item \textbf{It is not necessary to find the global minimum during optimization}.  Generally, it is not possible to determine if a complex ANN with many nodes and layers has been optimized to find the global minimum for some specified loss function; the problem is generally not convex, thus local many minima can occur.  Instead, training methods seek to optimize neural networks using statistical sampling of the data.  The hope is that such methods predict some, possibly local, optimization that provides an approximation that is sufficient given an acceptable loss.  Under such conditions, whether or not the optimum that is found is the global optimum or not is immaterial; rather, it is sufficient \emph{a priori} by statement of an acceptable error in approximation.  Recent research \citep{haeffele2017global} has suggested that with appropriately configured ANNs, local minima are also global minima (i.e., there are many local minima, but they all have the same value and each represent the global minimum).

\item \textbf{Overcompleteness is the rule, not the exception}.  When attempting to decompose a known function, complete, orthogonal bases can be very efficient.  However, for inverse problems where we want to determine the best function to fit a set of data, the use of complete orthogonal bases can actually be a hindrance.  The problem is that there is that searching parameter space for a unique solution (such as those presented by complete orthogonal bases) presents is a difficult problem.  Instead, using an overcomplete set of basis function allows one to search for\emph{any acceptable} solution; with overcomplete bases, there may be many such solutions \citep{lewicki2000learning}.  This is, to some degree, coupled with the previous principle regarding global versus local minima.

\item \textbf{Regularization is a useful practice.}  Regularization is a catch-all term that indicates the addition of information (via constraints on the solution space, or by  by adding terms to the loss function that penalize certain kinds of results).  One of the reasons that ANNs are successful is the application of such methods to \emph{regularize} the resulting solutions.   The effects of regularization are to reduce the parameter space that is searched through, and to reduce the number of parameters that are needed to generate acceptable solutions.  Thus, when a large number of solutions are all consistent with the data, some kinds of regularization (e.g., parameter number penalties in the loss function) breaks ties in favor of the solution with smallest number of parameters.  Despite these seemingly obvious payoffs, exactly how effective regularization is for deep neural networks is still a topic of investigation \citep{zhang2021understanding}. The results of \citet{zhang2021understanding} (among others) suggest that regularization can be helpful in finding solutions efficiently, it is not always necessary.
\end{enumerate}

\addcontentsline{toc}{section}{{Part I: Linear Feedforward Networks}}
\section*{\Large PART I:  Linear Feedforward Networks}

\section{Feedforward Networks for \uline{Linear} Regression}\indexme{feedforward network (FNN)!linear regression}

Now that the necessary background about ANNs and FNNs has been presented, it is time to put some of these ideas to use.  Nearly everyone is familiar, at least in concept, with the idea of linear regression.  It turns out that the problem of fitting a set of data with the ``best" possible solution is one that makes a good first foray into the applications of feedforward neural networks.

Quite frequently, FNNs are used to conduct either simple regression or classification on a set of data.  As such, one usually assumes that the data (which may or may not have a component representing noise from uncertainty in the data collection process) are, in actuality, described by \emph{some} function that is unknown.  If the data are categorical, then this function describes the boundary between categories in the proposed space of independent variables.  If the data represent some (presumed) continuum process, then it is assumed that regression (with some goodness-of-fit metric) is an appropriate method to generate an approximation of the function as it depends on its independent variables.  There is a third use for FNNs that is not generally discussed in the literature.  This is the purposeful \emph{approximation} of known functions (e.g., as one would do with, for example, a Fourier series).  

In this section, we will consider only linear regression problems on sets of data.  Later, we will extend these ideas to nonlinear regression on data, and nonlinear approximation of specified functions.

\subsection{An example of FNNs for linear regression of data}

The use of examples an effective way to see how the concept of FNNs can be combined with ideas from other areas of applied mathematics.  In the present example, a two-layer FNN will be constructed to find the best-fit solution to a set of data.  Here, the assumption is that the data are well-represented by a \emph{linear function}.  The example being presented is somewhat trivial; we will fit a plane through three points.  This is a linear problem that has a single unique solution for a specified loss function (and this solution can be found by more conventional methods); however, examining simple versions of otherwise complex problems can be a very effective means for bringing clarity to the methods used.  

For the purposes of example, suppose that we have the following set of data in the independent variables (features) $x_1$ and $x_2$, with the goal of fitting these data with a linear model.

\begin{table}
\sidecaption[b]
    \begin{tabular}{|c|c|c|c|}
    Example number  & ~~~$x_1$ ~~~ &  ~~~$x_2$ ~~~ &  ~~~$y$ ~~~\\
    \hline
       1  &  1 & 1 & 5\\
       2  &  1 & 2 & 7\\
       3  &  2 & 1 & 8\\
    \end{tabular}
    \caption{Data for a linear perceptron  model.  This data set is equivalently specified by the vectorized variable ${\bf{R}}=\{(1,1,5),(2,3,7),(2,1,8)\}$}.
    \label{tab:data1}
\end{table}
\noindent  Notation for independent variables (and their associated data) can be challenging at times with high-dimensional data sets.  To help reduce this problem, we denote the following subsets of the data.  The symbol ${\bf{X}}$ is used to indicate the set of independent variable data.  Thus, for the example above, we have set of ordered pairs of data given in the table above by ${\bf{X}}(x_1,x_2)$ (i.e., ${\bf{X}}(x_1,x_2)=\{(1,1),(2,3),(2,1)\}$).  Thus, we can refer to the entire independent variable (feature) data set by referring to the variable ${\bf{X}}(x_1,x_2,y)$, or when we do not need an explicit list of the independent variables, simply ${\bf{X}}$.  With a hopefully clear meaning, ${\bf X}_i$ indicates the $i^{th}$ tuple of independent variables.  A similar notation is used to describe the set of dependent data. Thus, we have ${\bf{Y}}(y)=\{5, 7, 8\}$.  Finally, the entire data set is a relation among the independent and dependent variables, denoted by ${\bf{R}}=\{{\bf X}R{\bf Y}\}$.  With a slight abuse of notation, in the example above we represent this by ${\bf R}=\{(1,1,5),(2,3,7),(2,1,8)\}$

The appropriate linear model is given by

\begin{equation}
 \hat{y}(x_1,x_2) = w_{1,4} x_1 + w_{2,4} x_2 + b_{3,4}
 \label{simple_percept_model}
\end{equation}

            \begin{figure}[t!]
            \sidecaption[t]
           \centering
            \includegraphics[scale=0.45]{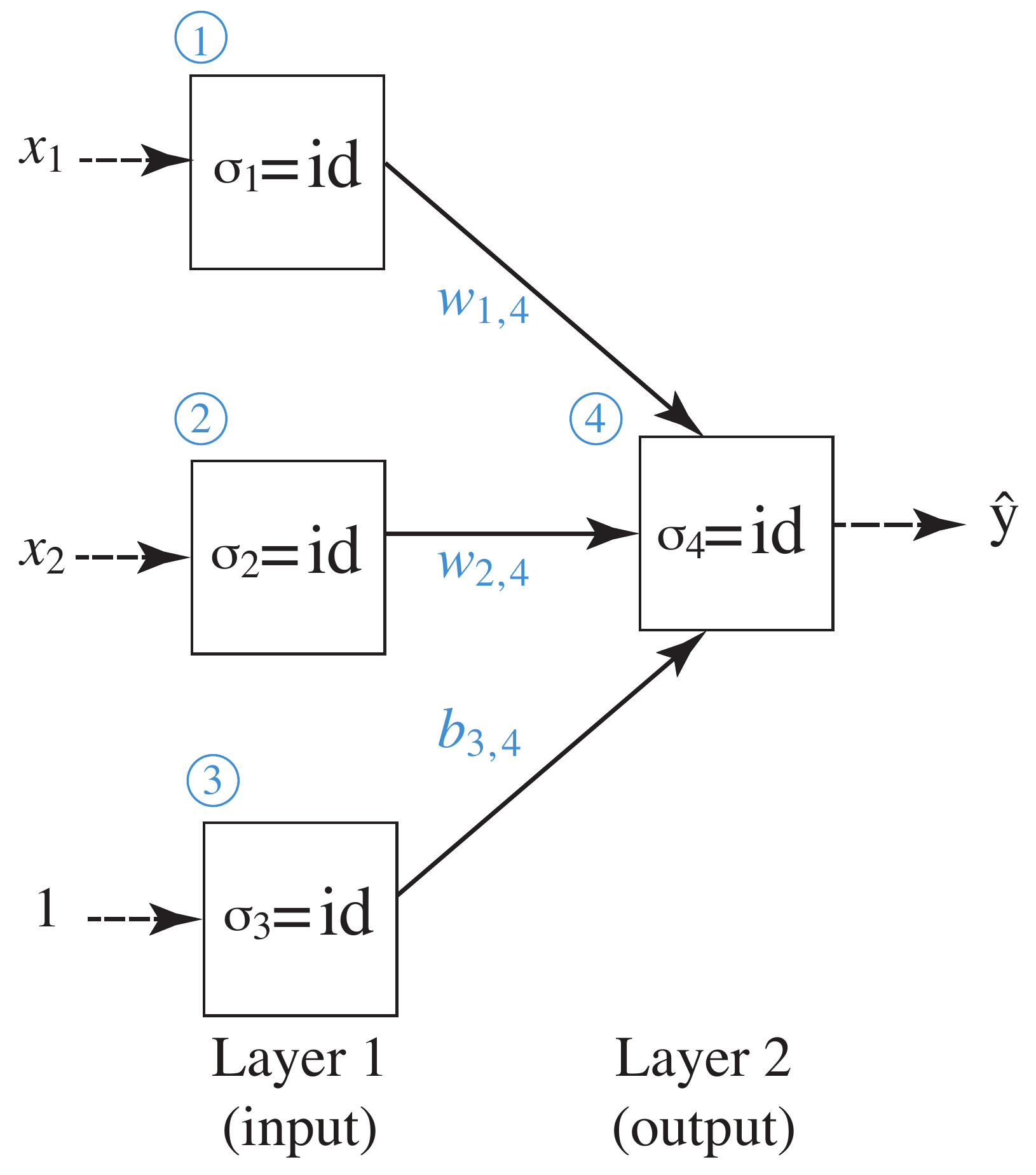}
            \vspace{-5mm}
            \caption{A two layer network.  This is a fully \emph{linear network} with a single input and a single output layer; there are no hidden layers (hence, no interior nodes).  While the original perceptrons represented binary functions, modern usage of the word frequently includes single-node networks with continuous activation functions (such as the linear one illustrated here).}
            \label{fig:simpleperc}       
            \end{figure}
%
Assuming that the solution we seek is indeed linear, then the network required to generate the parameters desired can be given by a \emph{simple perceptron}.  A simple perceptron is a network with no hidden layers.  Historically, the perceptron was the first physically-realized (i.e., constructed in the lab as an electro-mechanical systems) binary-classifying network that was shown to exhibit \emph{machine learning}; this was accomplished in the 1940s by Warren McCulloch and Walter Pitts \citep{mccullough1943logical}.  The initial goal for the simple perceptron was to model a biological neuron (hence, the common usage today of the term ``neural network").  Modern understanding of biological neural networks would suggest that this model is substantially over-simplified; but it still generated a useful archetype for network models.  While it is sometimes thought that simple perceptrons were strictly linear devices, this is not entirely true.  The output of the original perceptron models involved a third layer containing a nonlinear thresholding function that assigned a binary output.  Frequently, two layer models with the identity transformation in each node and a single collating layer (as layer 2) are referred to as perceptrons even though they are not binary classifiers.
A simple (linear) perceptron models provide a good platform for understanding how FNNs work.  For the problem at hand, a simple perceptron network is structured as in Fig.~\ref{fig:simpleperc}. 

At this juncture, we have all the components that we need to both formulate and solve the problem represented by this neural network.  The only remaining step is to develop an explicit form for the loss function, and then minimize it.  The least-square-error loss function is given by the square difference between the observed target values, $y_j$ and those predicted by the approximating function $\hat{y}_j$ at the same values of the independent variables as is associated with the targets, $y_j$.  This can be written compactly by 

\begin{equation}
    J({\bf R}; w_{1,4},w_{2,4},b_{3,4}) = \sum_{j\in\ell(1)} [y_j-\hat{y}(x_{1,j},x_{2,j})]^2
    \label{lossfunc1}
\end{equation}
where here  $\ell(1)=\{1,2,3\}$ is an index set for the number of data points.  The subscript ``$i$" on $x_1$ and $x_2$ is then used to index the associated value of $x_1$ and $x_2$ from the example data.  Sometimes the loss function is normalized by the number of data appearing in the function, but this is done primarily as a normalization in computational optimization.  If we optimize using analytical methods, the normalization is immaterial.

With such a simple network, the loss function is easily written out explicitly, and this will help to clarify the notation further.  The loss function given above can be more verbosely written as
\begin{align}
    J({\bf R}; w_{1,4},w_{2,4},b_{3,4}) &=  [5-w_{1,4}-w_{2,4}-b_{3,4}]^2 + [7-w_{1,4}-2 w_{2,4}-b_{3,4}]^2\nonumber \\
    &+ [8-2w_{1,4}-w_{2,4}-b_{3,4}]^2
\end{align}
It is necessary to tend to another notational detail.  As with the list of dependent and independent variables comprising the data, the list of parameters that are being optimized can grow quite long.  Thus, it is common to give them a single summary vector notation.  Here, we define $\boldsymbol{\theta}=\{w_{1,4},w_{2,4},b_{3,4}\}$.  This vector of parameters can be indexed similarly to those for ${\bf X}$ (i.e., $\theta_{1,3}=w_{1,4}$, etc.)  Now, the expression above can be more clearly written out by 

\begin{align}
    J({\bf R}; \boldsymbol{\theta}) &=  [5-w_{1,4}-w_{2,4}-b_{3,4}]^2 + [7-w_{1,4}-2 w_{2,4}-b_{3,4}]^2\nonumber \\
    &+ [8-2w_{1,4}-w_{2,4}-b_{3,4}]^2
\end{align}

\begin{figure}[t!]
            \sidecaption[t]
           \centering
            \includegraphics[scale=0.45]{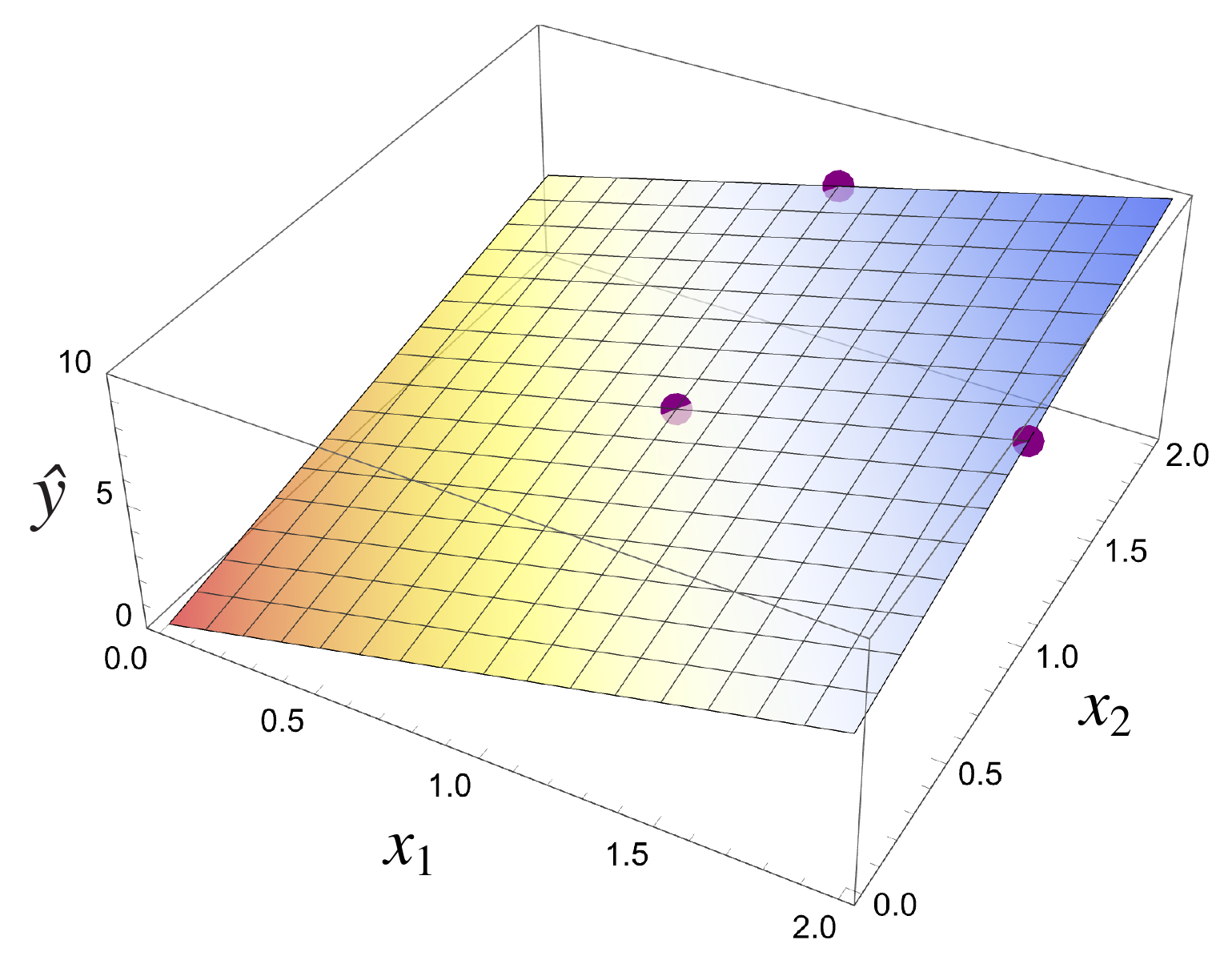}
            \vspace{-5mm}
            \caption{The optimized solution to Eq.~\eqref{simple_percept_model} with Eq.~\eqref{lossfunc1} as the loss function.  Points plotted represent the data used to \emph{train} the simple perceptron network.}
            \label{fig:simplepercplot}       
            \end{figure}

Now the goal is to optimize the weights of the loss function such that it reaches its minimum.  To do so, we can take derivatives of the loss function with respect to each adjustable parameter, and set each such derivative equal to zero.  If the resulting set of equations can be solved simultaneously, then we have a solution to the optimization problem.  Computing the appropriate derivatives leads to the following set of equations

\begin{align}
    12 w_{14}+ 10 w_{24}+8 b_{34}-56 &=0 \\
    10 w_{14}+12 w_{24}+8 b_{34}- 54 &=0\\
    8 w_{14}+8w_{24}+6 b_{34} -40 &=0
\end{align}
This linear system is easily solved to yield
\begin{align*}
    w_{14}&= 3\\
    w_{24}&= 2 \\
    b_{34} &= 0
\end{align*}
And this yields the \emph{approximating function}
\begin{equation}
 \hat{y} = 3 x_1 +2 x_2 
\end{equation}
Because there is only a single unique plane between any three (unique) points in $\mathbb{R}^3$, the error is identically zero for this case; this is easily verified by computing the value of Eq.~\eqref{lossfunc1} with this approximating function.  A plot of the fit to the data and the data used in fitting appear in Fig.~\ref{fig:simplepercplot}.

\begin{svgraybox}
\begin{example}[Fitting a linear expression to a data set]

In the previous example, the system was somewhat of a special case because there is only one plane through three lines in $\mathbb{R}^3$.  Thus, the residual error in this example was identically zero.  However, the process changes very little if there are more data point (or, indeed, even if there are more independent variables).  Suppose we had a fourth data point in the set of data given by the triplet $(2,2,13)$, so that ${\bf{R}}=\{(1,1,5),(2,3,7),(2,1,8),(2,2,13)\}$.  The analysis would be identical, except that the new loss function would be 

\begin{align}
J({\bf R}; w_{1,4},w_{2,4},b_{3,4}) &=  [5-w_{1,4}-w_{2,4}-b_{3,4}]^2 + [7-w_{1,4}-2 w_{2,4}-b_{3,4}]^2\nonumber \\
    &+ [8-2w_{1,4}-w_{2,4}-b_{3,4}]^2+[13-2 w_{1,4} - 2 w_{2,4}-b_{3,4}]^2
\end{align}

The values of the derivatives are slightly changed, yielding the set of equations

\begin{align}
   20 w_{14}  +  18 w_{24} + 12 b_{34} -108 &=0\\20 18 w_{14}  +  20 w_{24} + 12 b_{34} -106 &=0\\20 12 w_{14}  +  12 w_{24} + 8 b_{34} -66 &=0
\end{align}
The solutions is routine, and found to be
\begin{align}
    w_{14} & = \frac{9}{2} \\
    w_{24} & = \frac{7}{2} \\
    b_{34} &= -\frac{15}{4}
\end{align}

The loss function is no longer zero.  The best we can do is the minimum value for the loss function given by $J({\bf R}; 9/2, 7/2, -15/4) = 9/4$.
A revised plot based on this optimization is given below.

{\hspace{14mm}
\centering\includegraphics[scale=.4]{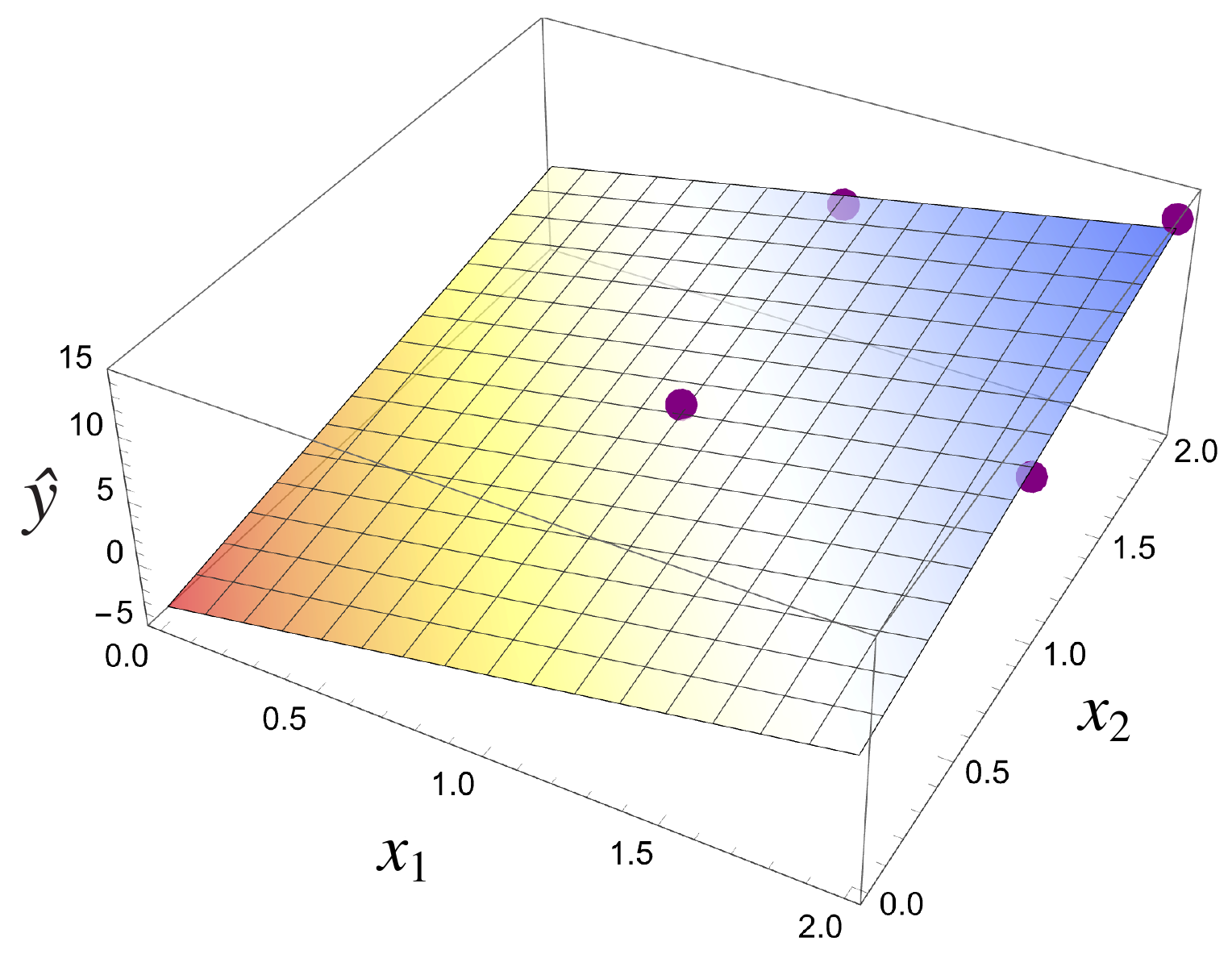}
\vspace{4mm}
}
\captionof{figure}{Example of fitting a plane to four data points by minimizing the sum of the squared errors.}
\end{example}
\end{svgraybox}

\addcontentsline{toc}{section}{{Part II: Nonlinear Feedforward Networks}}
\section*{\Large PART II:  Nonlinear Feedforward Networks}\indexme{feedforward network (FNN)!nonlinear}

In the examples previously, neural networks with simple linear functions were illustrated, and it was shown that the conventional linear least-squares fitting algorithm can be thought of as a two-layer neural network.  While this is a good example for illustrating the \emph{concept} of neural networks, it does not tap into many of the strengths that modern neural networks embody.  In particular the following are not possible to do with such simple networks: (1) Nonlinear functions are, by definition, not well represented by a linear map, and (2) the inherent \emph{efficiency} realized by deep neural networks is not achieved by such simple networks.  

In the material that follows, we relax both of these to the extent that is practical in an analytic setting; that is, we will examine problems which can be optimized \emph{without} the use of a computer.  One may wonder if direct analytical optimization represented \emph{machine learning} or not; reference to the definition of machine learning from the glossary would suggest that it does not!  This is indeed the case.  The algorithms involved in machine learning are \emph{numerical methods} that can be used to iteratively train large and complex ANNs.  However, ANNs exist and can be studied \emph{without} algorithms that qualify as machine learning.  One might argue (and this is the perspective of this chapter) that understanding the structure and behavior of ANNs analytically is an important first step in the study of ANNs in the more general context that machine learning algorithms allow.  

Our goal in this second part of the chapter will be to extend our understanding of fully connected feedforward neural networks (FNNs).  This will be done by (1) allowing nonlinearity by adopting nonlinear activation functions, and (2) begining to examine the role of depth by looking at neural networks with one or more hidden layers.  Our focus will be primarily networks with a single hidden layer which are technically not \emph{deep}.  However, the concepts associated with nonlinear representation and determination of the optimal parameters are not substantially different for deep networks, at least in the analytical context.  With neural networks with two or more layers, the differentiation of the network becomes more complex because of the need for the chain (or composition) rule for differentiation, and this mirrors one of the early challenges for numerical optimization of neural networks that was ultimately solved by backpropagation algorithms.  This second part of the chapter begins with a discussion of a particular nonlinear activation function known as the rectified linear unit ($\relu$) function.

\section{The $\relu$ Function}\indexme{ReLU function}

For computations to occur, we must choose an activation function to use on the network.  Ultimately the depth and purpose of the network will have some influence on the choice of activation function.

For deep neural networks, the choice of the activation functions makes a \emph{significant difference} in how well the network performs (i.e., how easy or difficult it is to optimize).  The rectified linear unit ($\relu$) function has been found to outperform other activation functions for both regression and classification tasks in terms of computational efficiency.  It has several advantages that have led to its wide use.  The primary reason are practical ones: (1) it is very economical to compute the function (and its derivative), and (2) the $\relu$ function does not suffer from the ``vanishing gradients" problem when $|x|\rightarrow \infty$. This latter property is important in deep networks for the \emph{backpropagation} of weight information through the network.  As a concrete counter example, of the sigmoidal functions (see Fig.~\ref{fig:activations}),  were very popular in the past.  However, vanishing gradients as $|x|\rightarrow \infty$ is a problem for these functions, and it somewhat stymied research on deep networks for some time \citet{glorot2010understanding}.

The $\relu$ function has a number of the other nice features.  It is a nonlinear and non-polynomial function, which is essential for its use as a basis function \citet{pinkus1999approximation}.  While it is nonlinear, it also has the property that any weighted sums of the $\relu$ function generate a piecewise-linear result.  This is, interestingly, also true for compositions: the composition of a $\relu$ function with any weighted sum of $\relu$ functions generates a piecewise linear function.  These qualities give the $\relu$ function an algebraic structure that makes analysis of large networks more tractable.

Among the various possible definitions, we adopt the following for the $\relu$ function defined on $x\in (-\infty,\infty)$

\begin{equation}
    \reluf{x} = 
    \begin{cases}
        0 & x<0 \\
        1 & x>0 \\
        0 & x = 0
    \end{cases}\label{reludef}
\end{equation}

While the $\relu$ function does not have continuous derivatives, its other features more than compensate for this difficulty.  This is a good example of a development that occurred for largely practical reasons that later is shown to have a particular theoretical structure that makes it especially effective.
with
\begin{equation}
    \frac{d}{dx}\reluf{x} = 
    \begin{cases}
        0 & x<0 \\
        1 & x>0 \\
        \textrm{undefined} & x = 0
    \end{cases}
\end{equation}
In the literature on FNNs, the derivative of the $\relu$ function is often assigned an arbitrary derivative (usually between 0 and 1) at $x=0$.  However, it turns out that the discontinuity in the derivative at $x=0$ posed neither analytical or numerical problems for regression problems.  In the section below regarding the calculus of $\relu$ functions, some additional explanation is provided.

Because it will be helpful notation in the material that follows, we define the function $H_0$ as the Heaviside function 

\begin{equation}
    H_0(x)=
    \begin{cases}
        0 & x\le 0 \\
        1 & x>0 \\
    \end{cases}
\end{equation}
Thus, and alternative expression for the $\relu$ function is 

\begin{equation}
    \reluf{x} = x H_0(x)
\end{equation}

\section{Algebraic Properties of $\relu$}\indexme{ReLU function!algebraic properties}

What is important in applications is (generally) to understand how one rescales, translates, and adds $\relu$ functions.  While $\relu$ functions are nonlinear, they do, nonetheless, have some properties that are nearly linear (and this is, in part, the reason that it has been so successful).  To begin this discussion, we focus on a single independent variable; extensions to multiple independent variables is discussed subsequently.
The rescaling (by multiplication by a scalar), translation, and addition of $\relu$ functions are described below.  Note that the ReLU function is \emph{nonlinear}.  The $\relu$ function passes neither the tests of homogeneity, nor additivity.  In other words, we have 
\begin{equation}
    \reluf{a x} \ne a \relu{x}
\end{equation}
and
\begin{equation}
    \reluf{x+y}\ne \reluf{x}+\reluf{y}
\end{equation}
The piecewise-linear form of the ReLU function does give it some advantages.  The quasi-linearity of the function helps tremendously in analytic efforts for solutions and interpretation of FNNs.  In addition, the ReLU function avoids some of the problems created by other candidates for activation functions that have been used in the past.  For example, the logistic (``s" shaped) function is smooth, but because its derivatives rapidly tend toward zero as $|x|\rightarrow \infty$, such functions can create difficulties when using numerical methods such as gradient decent and conventional backpropagation of information throughout the network (which both rely on some linearly proportional function of the \emph{derivatives}).  

An affine function in one variable is of the form $f(x)=a x$ \emph{plus all of its translations} (e.g., $f(x)=a(x+\beta) = ax + a\beta$ is a linear function plus all of its translations, $\beta$. The $\relu$ function is piecewise-affine; because there is much existing theory on such functions, this theory has been useful for better understanding the function of deep neural networks using $\relu$  as the activation function.

\subsection{Translation}\indexme{ReLU function!translation}

It is a bit unusual to \emph{start} the discussion of the algebra of $\relu$ functions with the definition of translation operations, but they are possibly the most important concept to understand for applications to ANNs.
Suppose in one dimension we define $T_\beta$ as an operation that \emph{translates} the function it operates on by the amount $\beta$ (which can a be positive or negative quantity), then  we have 

\begin{equation}
 T_\beta(\reluf{x}) = \reluf{x+\beta}
\end{equation}
The translation operator acts only on the variable $x$, so some care is needed in interpreting it correctly.  For example, if a scalar multiplies $x$, we interpret the translation as follows.

\begin{equation}
     T_\beta\reluf{ax} = \reluf{a(x+\beta)}
\end{equation}
This operation is correct regardless of the sign of $a$.  Sequential translations can be defined (although one does not usually encounter this)

\begin{align}
    T_\beta T_\gamma \reluf{ax} & = T_\beta \reluf{a(x+\gamma)} \nonumber\\
    &= \reluf{a([x+\beta]+\gamma)} \nonumber\\
    &= \reluf{a(x+\beta+\gamma)} 
\end{align}

\noindent Note that the translation operator is a linear one:

\begin{equation}
    T_\beta\left[\reluf{w_{1,3} x}+\reluf{a_2 x}\right] = \reluf{w_{1,3} (x+\beta)}+\reluf{a_2 (x+\beta)}
\end{equation}

While $\beta$ the the amount of the \emph{translation} of the $\relu$ function, it is not in general equivalent to the \emph{bias} term that is typically adopted in applications.  To make this clear, define  $b = a \beta$.  Then 
            \begin{figure}[b]
            \sidecaption[t]
           \centering
            \includegraphics[scale=0.35]{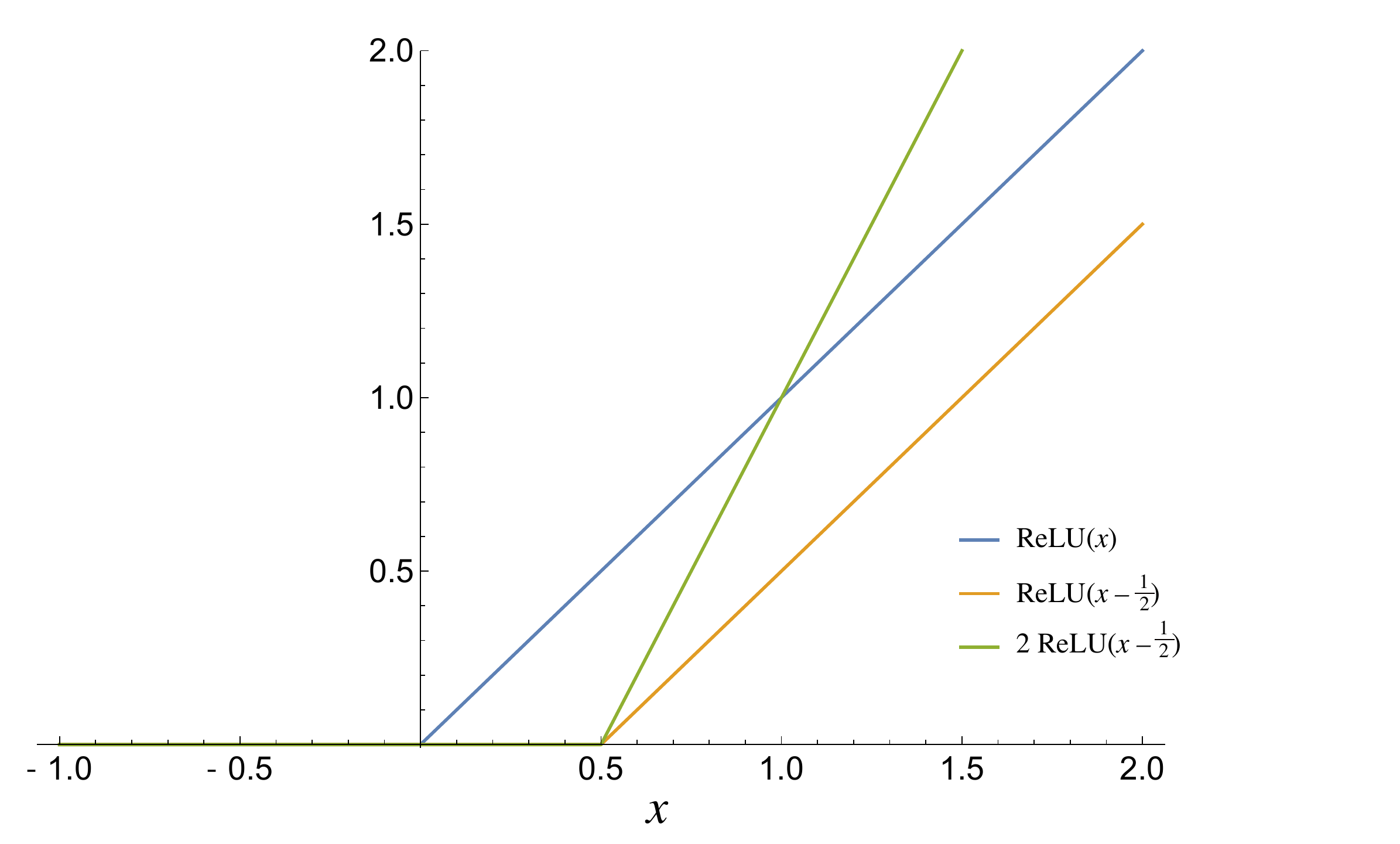}
            \vspace{-5mm}
            \caption{Relu functions, both shifted and scaled.}
            \label{fig:relushiffig}       
            \end{figure}
\begin{align}
 \reluf{a(x+\beta)}& =\reluf{a x +  b} \nonumber 
\end{align}
where $b$ is the \emph{bias} term. From here forward, the symbol $b$ will always imply the \emph{bias}, whereas the symbol $\beta$ will be used to indicate the translation.  In Fig.~\ref{fig:relushiffig}, shifted and scaled functions $relu$ are plotted for reference. \\

\subsection{Addition}

\noindent The addition of $\relu$ functions is straightforward.  Combining the operations above with that of addition gives 

\begin{align}
 T_{\beta_1}(a_1 \reluf{x}) + T_{\beta_2}(a_2 \reluf{x}) &= \reluf{a_1(x +\beta_1)}+\reluf{a_2(x + \beta_2)} \nonumber \\
 & = \reluf{a_1 x +b_{2,3}}+\reluf{a_2 x +b_2}
\end{align}
%
            \begin{figure}[t]
            \sidecaption[t]
            \centering
            \includegraphics[scale=0.4]{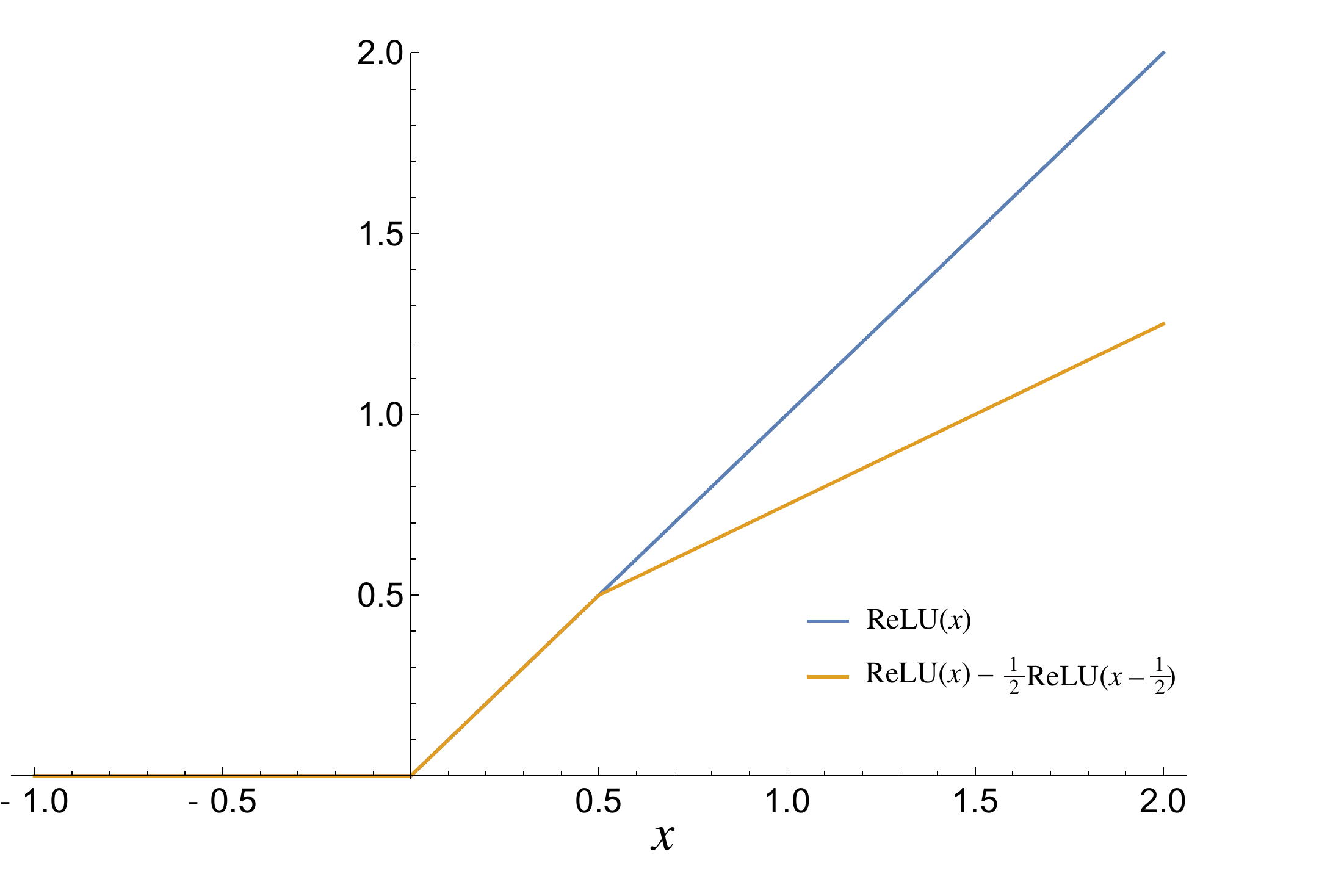}
            \vspace{-5mm}
            \caption{Addition of non-degenerate $\relu$ functions.}
            \label{fig:reluadd}       
            \end{figure}
%
An example of the sum of $relu$ functions appears in Fig.~\ref{fig:reluadd}.

While this definition is straightforward, note that when $\beta_1 = \beta_2=\beta$, the sum is a \emph{degenerate} case; that is, sum of the two functions can be replaced by a single $\relu$ function.

\begin{align}
 & = \reluf{a_1 (x +\beta)}+\reluf{a_2 (x + \beta)} \nonumber \\
&= T_{\beta}(a_1 \reluf{x}) + T_{\beta}(a_2 \reluf{x}) \nonumber \\
& = T_{\beta}(\reluf{(a_1 +a_2)x} ) \nonumber \\
\intertext{Setting $a_3 = a_1+a_2$ gives the result }
&= \reluf{a_3 (x + \beta)}
\end{align}
In applications, degenerate cases are inefficient (by unnecessarily increasing the number of parameters involved).  Depending upon the system being investigated using $\relu$, one might also seek to avoid cases that are \emph{nearly degenerate}, i.e., for two one-dimensional $\relu$ functions shifted by $\beta_1$ and $\beta_2$ respectively, then if $\epsilon_0$ represents a minimum (scaled) separation or ``resolution'' between any two $\relu$ functions, a nearly degenerate case would be indicated by

\begin{equation}
    \left|1- \frac{\beta_1}{\beta_2}\right|  < \epsilon_0
\end{equation}
The problem created by nearly degenerate cases is a practical one.  For such cases, the concept of nearly degenerate would be established \emph{a priori}; it means only that the expressiveness gained by adding two $\relu$ functions that are sufficiently close together (in terms of the shift variable, $\beta$) is not warranted by the quality of the data or other external factors.

\subsection{Multiplication by a scalar}\indexme{ReLU function!scalar multiplication}

Multiplication by a positive scalar is commutative with the $\relu$ function.  Thus, we can write for any real number $a$

\begin{equation}
 \reluf{a^2 x}= a^2 \, \reluf{x} 
\end{equation}
Unfortunately, this kind of linear behavior does not extend to cases where the constant might be negative.  Suppose $a$ is any real number, then we have

\begin{align}
    \reluf{a (x+\beta)} &= \reluf{\textrm{sng}(a)|a| (x+\beta)} \nonumber\\
    &= |a| \,\reluf{\textrm{sng}(a)(x +\beta)}
\end{align}
where recall $\beta=b/a$, and the signum function is defined by 
\begin{equation}
    \textrm{sgn}(a)=\frac{a}{|a|} = \begin{cases}
        +1 & \textrm{ for } a>0 \\
        ~~0 & \textrm{ for } a=0\\
        -1 & \textrm{ for } a<0 
    \end{cases}
\end{equation}

Conversely, we also have
\begin{align}
   c \reluf{a'(x+\beta)} &= |c|\textrm{sng}(c)\reluf{a(x+\beta)} \nonumber\\
    &=  \textrm{sng}(c)\,\reluf{|c|a'(x+\beta)}
\end{align}

Note that we can now define $a=a' |c|$ to write
\begin{align}
   c \reluf{a'(x+\beta)}  &=  \textrm{sng}(c)\,\reluf{a(x+\beta)}
\end{align}
This last relation does allow some simplification of the operation of scalar multiplication with a $\relu$ function. In the last layer of a network, we often have forms that are simple weighted sums of $\relu$ functions.  This last relation allows us to replace those weights by $\pm 1$, where now the magnitude is fixed, but the sign must be determined.  Later, we will encounter examples where this property allows us to make effective constraints that allow optimization to be done more efficiently.

            \begin{figure}[t]
            \sidecaption[t]
           \centering
            \includegraphics[scale=0.8]{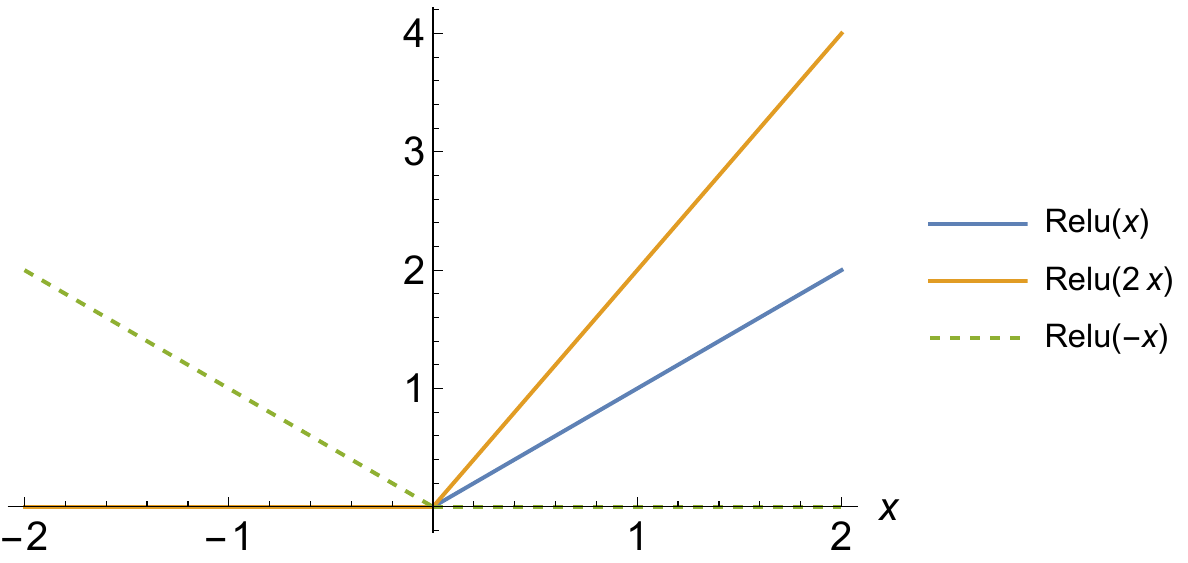}
            \vspace{-5mm}
            \caption{Three $\relu$ functions, differing only by the constant, $a$.}
            \label{relusaxfig}  
            \vspace{4mm}
            \end{figure}

\subsection{Multiplication of {$\relu$} functions}

Although not commonly arising in applications, the product of $\relu$ functions is nonetheless well-defined.  At the end of this chapter, we will provide an example of using a $\relu$ squared function as the activation.  To define the product of $\relu$ functions, it is helpful to adopt the following form for the function 

\begin{equation}
    \reluf{a x + b} = \reluf{a(x+\beta)}= a(x+\beta)H_0(a(x+\beta))
    \label{genrelu2}
\end{equation}
Here, $H_0(x)$ is the Heaviside function defined previously.

Assume we have the multiplication of $\reluf{a_1 (x - b_{2,3})}$ and $\reluf{a_2 (x - b_2)}$.
With Eq.~\eqref{genrelu2}, it is easy to verify (using the properties of the Heaviside function)

\begin{align}
    \reluf{a_1 x - b_{2,3}} \cdot \reluf{a_2 x - b_2} &= (a_1x-b_{2,3})H_0\left[a_1 x-b_{2,3}\right] \left( a_2 x-b_2\right)H_0\left[ a_2x-b_2\right] 
   \end{align}
From this point, there are four possible cases regarding the configuration of the Heaviside functions.  Thus, there is no \emph{simple} simplification.  The $\relu^2$ case is one that is used in applications, and has a simplification.  This case can be simplified as follows.

\begin{align}
    \textrm{ReLU}^2({a (x - b)})
     &= H(a x+b)\reluf{ax+b)} H(a x+b)\reluf{ax+b)}\nonumber\\
     &= H(a x+b)(a x-b)^2
\end{align}

\begin{figure}[!ht]
            \sidecaption[t]
            \centering
            \includegraphics[scale=0.4]{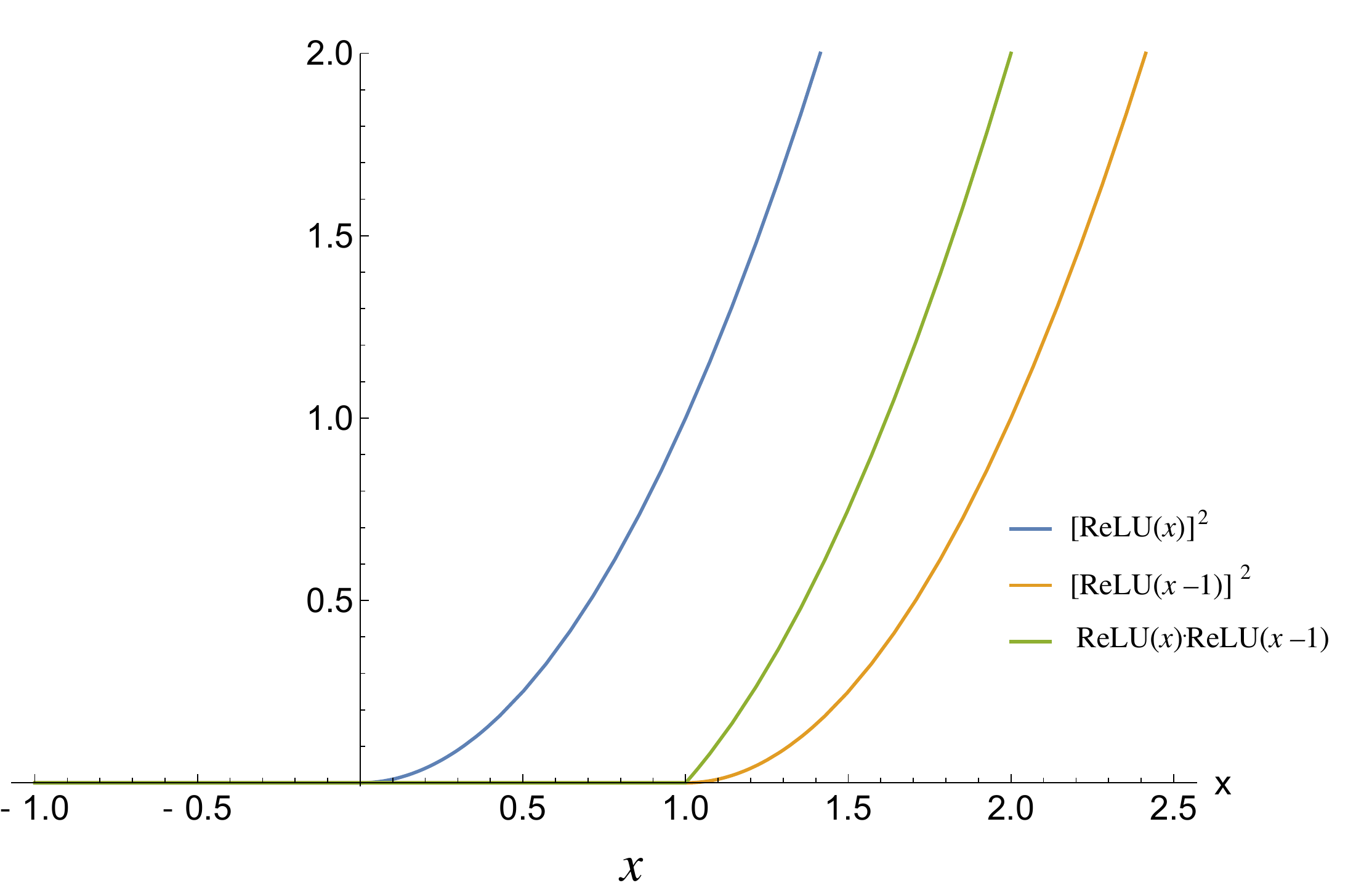}
            \vspace{-5mm}
            \caption{Multiplication of $\relu$ functions.}
            \label{relusmultfig}       
            \end{figure}

Example plots of the multiplication of $relu$ functions appears in Fig.~\ref{relusmultfig}.

\subsection{Composition of {$\relu$} Functions}\indexme{ReLU function!composition}

\begin{figure}[!ht]
            \sidecaption[t]
           \centering
            \includegraphics[scale=0.55]{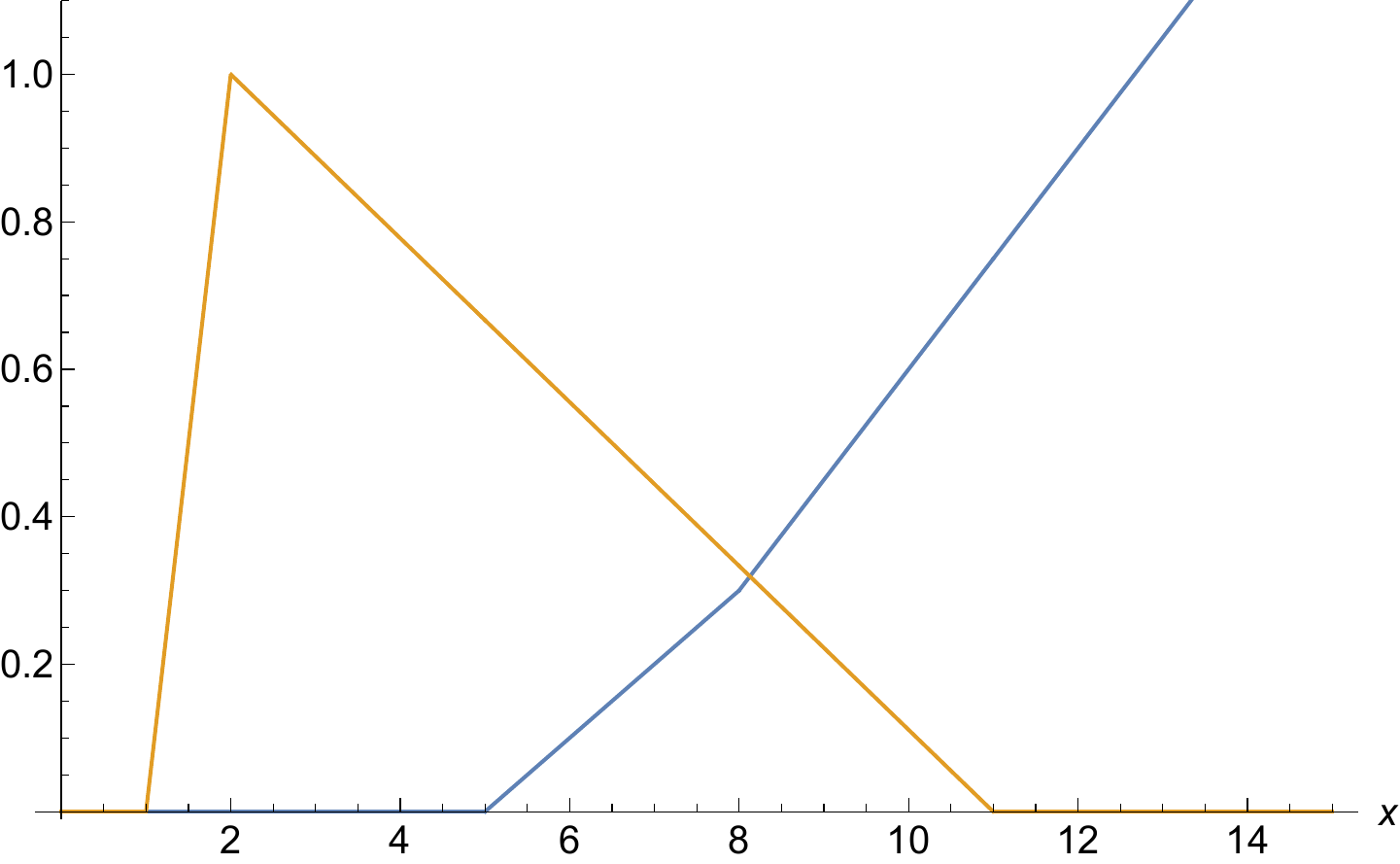}
            \vspace{-5mm}
            \caption{The composition of two ReLU functions. Orange line: $f(x)= \relu\left[{\reluf{x-1}-\frac{10}{9}\reluf{x-2}}\right]$.  Blue line:  $f(x)= \frac{1}{10}\relu\left[{\reluf{x-5}+\frac{1}{2}\reluf{x-8}}\right]$}
            \label{fig:relucomp}       
            \end{figure}

Recall that the composition of two functions $f$ and $g$ with independent variable $x$ is a given by

\begin{equation}
    F(x;f,g) = f(g(x))
\end{equation}
Simple compositions of ReLU functions have the property that they generate another ReLU function.  Here, by \emph{simple} we mean that only a single ReLU function is composed with another ReLU function.  For a simple compositions we have the following.

\begin{align}
   \reluf{ \reluf{ax-b}} &=\reluf{ a\reluf{x-b/a}} \nonumber\\
   & =\reluf{a\reluf{x}-b/a} \nonumber \\
   &= \reluf{x-b/a}
\end{align}
While the composition of two ReLU functions is somewhat uninteresting, things change dramatically when one considers the composition to a ReLU functions with a linear combination of two others.  Because of the non-additivity of ReLU functions, such compositions result in complicated relationships that cannot be easily expressed in simpler terms.  In other words, expressions like $\reluf{\reluf{a_1 x + b_{2,3}} + \reluf{a_2 x + b_2}}$ cannot be more simply expressed in general.  As a concrete examples, consider 

\begin{align}
   f(x)&= \relu\left[{\reluf{x-1}-\frac{10}{9}\reluf{x-2}}\right]\\
   f(x) &= \frac{1}{10}\relu\left[{\reluf{x-5}+\frac{1}{2}\reluf{x-8}}\right]
\end{align}
The plots of these compositions appear in Fig.~\ref{fig:relucomp}.  It is not difficult to show that \emph{all} such compositions of the sum of two $\relu$ functions are positive functions (which includes the zero function).  Without delving into details, it can be shown that the properties of compositions of $\relu$ functions are one part of the explanation of why neural networks have been so successful.  

As mentioned above (and this is true even for fixed-width neural networks) as the depth increases the underlying space of independent variables shows an exponential increase in the number of regions which are defined by linear functions (cf., \citet{hanin2019universal}). In other words, the $\relu$ functions create a piecewise-continuous function of the independent variables, and the number of pieces comprising the function increase at (asymptotically) exponential order.  This feature of compositional $\relu$ functions is quite successful.  In short, one could accomplish the same result using a \emph{single layer}, but an equivalent single layer network would contain exponentially more parameters than would a multi-layer network.  Because these results are determined by asymptotic analysis of estimates, one generally needs to have a sufficiently large independent variable dimension and sufficient number of layers in order to achieve nearly exponential scaling in practice.  However, this is how many deep neural networks are structured in conventional practice, so the power of compositions of functions can be very effective.

\section{Calculus of $\relu$ Functions}\indexme{ReLU function!calculus}

\subsection{Integration of the $\relu$ function}

Using the alternative expression for the $\relu$ function can be found by using the Heaviside (step) function.

\begin{equation}
    H_0(x) = 
    \begin{cases}
    0 & x \le 0 \\
    1 & x > 0
    \end{cases}
\end{equation}
With this definition, we have the alternative expression for the $\relu$ function

\begin{equation}
    \reluf{x} = x H_0(x)
\end{equation}

\begin{align}
   \frac{d}{dx}\reluf{x} &= H_0(x)  \\
   \reluf{x} &= \int_{z=0}^{z=x} H_0(z) \, dz \\
   \reluf{x} &= x H_0(x) \label{relu2}
\end{align}
of these, the last one provides a convenient way of integrating the $\relu$ function. Assuming $x>0$, then

\begin{align}
    \int_{-\infty}^{z=x} \reluf{z} \,dz & = \int_{z=0}^{z=x} z H(z) \,dz \nonumber \\
    & = \int_{z=0}^{z=x} z \, dz \nonumber \\
    & = \frac{x^2}{2}
\end{align}

\subsection{Derivative of the $\relu$ function}

The need for the derivative of the $\relu$ function arises during the optimization process.  While many references specify an arbitrary value for the derivative at the point of discontinuity, this is somewhat unnecessary.  Again expressing the $\relu$ function in the form

\begin{equation}
    f(x) = \reluf{ax + b} H_0(ax +b)
    \label{reluH}
\end{equation}
we can consider the derivative of the $relu$ function in the distributional sense (cf. Chp.~\ref{deltachap}).  Taking the derivative of both sides of Eq.~\eqref{reluH}, we obtain the distributional derivative

\begin{align}
    \frac{df}{dx} &= \frac{d}{dx}\left[\reluf{ax + b} H_0(ax +b) \right]\nonumber \\
    &=a H_0(ax+b)+\reluf{ax+b}\delta(ax+b)
    \label{dreluH}
\end{align}
Now note that the second term in this expression involves the delta function.  Recall that in a distributional sense, the delta function is only nonzero at the point where its argument is zero, i.e., when $ax+b=0$.  But, it is multiplied on the outside by a function that will return zero when $ax+b=0$.  Thus, this term is always zero.  More correctly, we might establish this by integration of both sides of this result from $(-\infty, x)$ with $x>0$

\begin{align}
 \int_{-\infty}^x  \frac{df}{dx}\, dx&=   \int_{-\infty}^x a H_0(ax+b)+\reluf{ax+b}\delta(ax+b)\, dx \nonumber \\
 &= \int_{-\infty}^x a H_0(ax+b)\, dx+\reluf{0}\delta(0)\,  \nonumber \\
 &= \int_{-b/a}^x a \,  dx \nonumber \\
 &= a x +b ~~\textrm{for $x>-b/a$} \\
 &= \reluf{ax + b} H_0(ax +b)
\end{align}
showing that we recover the original function after integrating the derivative (as we must). The practical result of this analysis is that the derivative of the $\relu$ function can safely be set to zero at the point of discontinuity, and this result is at least consistent with the distributional interpretation of the derivative.

\section{~${ReLU}$ Functions in Multiple Dimensions}\indexme{ReLU function!multiple dimensions}

One of the interesting properties of $\relu$ functions is how complex simple compositions of such functions can become.  This becomes more evident as the number of independent variables (features) grows.  Here, the purpose of discussion is primarily to build some intuition as to how such functions behave and appear.  Restricting ourselves to 2-dimensions, we can consider $\relu$ functions that depend on the pair $(x,y)$.  For concreteness, we will select specific values of the weighting (slope and bias) parameters.  Suppose we define the following  functions

\begin{align}
    \sigma_1(x,y) &= \reluf{x + 2 y -4}+2\\
    \sigma_2(x,y)&= \reluf{-x +1/2 y -1}
\end{align}
and
\begin{align}
    \sigma_3(x,y) &= \reluf{ \sigma_1(x,y)- \sigma_2(x,y)}
\end{align}
While the composition of difference of the two $\relu$ functions does not necessarily look complicated, its plot illustrates the ability of compositions of $\relu$ functions to generate an increasing number of regions.  This increase is roughly polynomial (i.e., grows at some factor that is described by a polynomial) as the number of independent variables or network width increases, and grows exponentially as the number of layers increase \citep{pascanu2013number}.  

\begin{figure}[!ht]
            \sidecaption[t]
           \centering
            \includegraphics[scale=0.55]{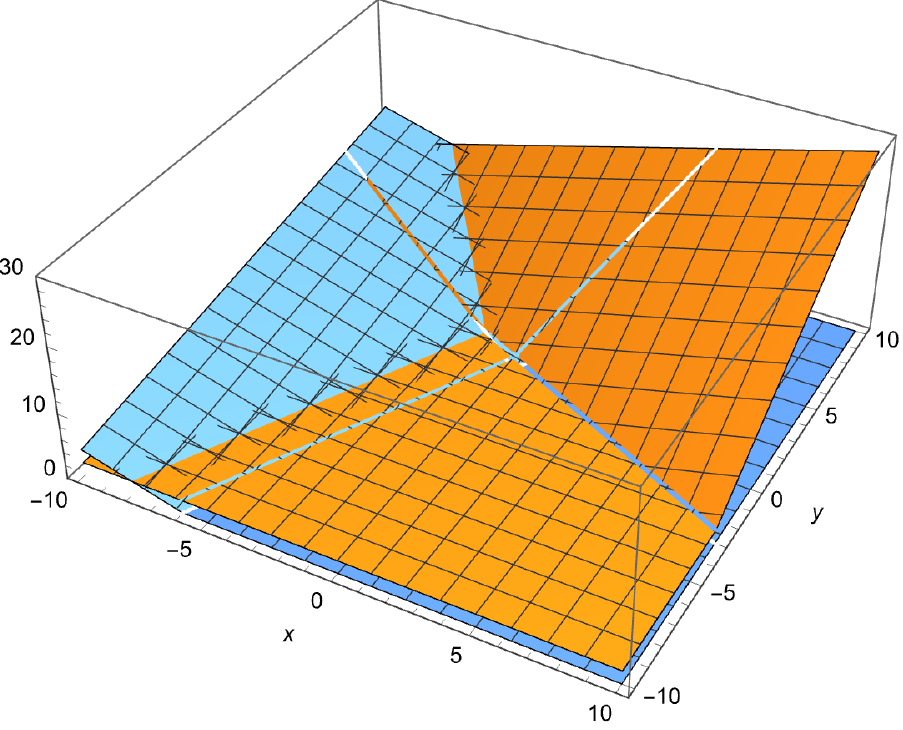}
            \includegraphics[scale=0.55]{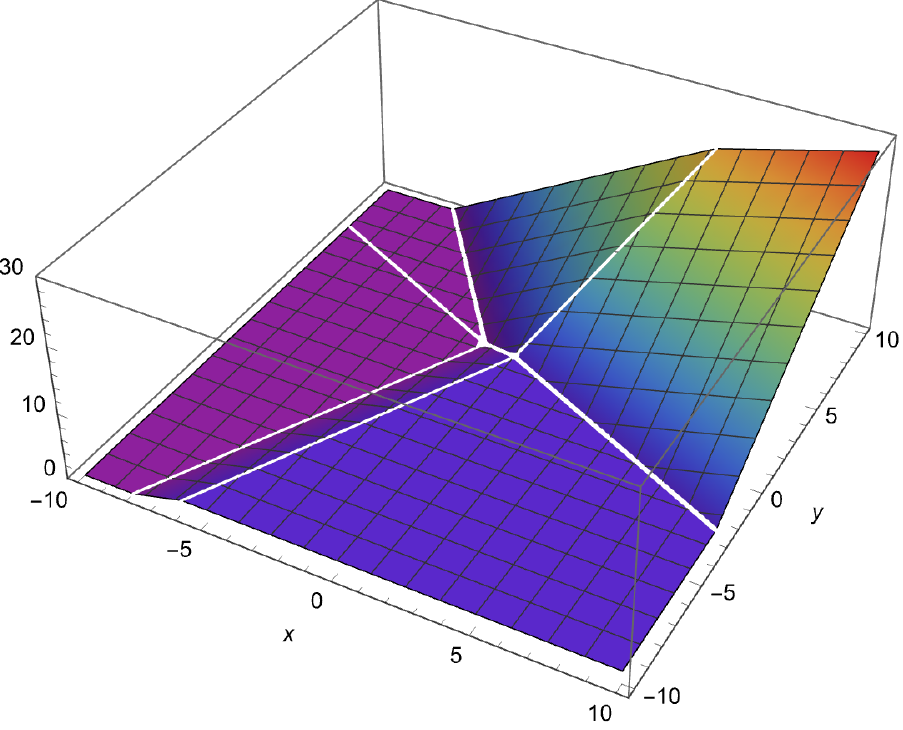}
            \vspace{-5mm}
            \caption{(Left) $\relu$ functions, $\sigma_1$ (orange) and $\sigma_2$ (blue) in two independent variables.  (Right) The function $\sigma_3=\reluf{\sigma_1-\sigma_2}$.  Color scale for the right figures denotes the value of $\sigma_3$, and is added to improve function visualization.}
            \label{fig:relucomp2}       
            \end{figure}

As an example of this behavior, the functions $\sigma_1$, $\sigma_2$, and $\sigma_3$ are plotted in Fig.~\ref{fig:relucomp2}.  

Various researchers have studied the number of regions in the independent variables formed by the $\relu$ functions for various configurations of FNNs. Estimates (and lower bounds) for the number of regions formed by $\relu$ FNNs is detailed in \citet{montufar2014number} and \citet{pascanu2013number}.

\section{~Brief Comments on $\relu$ Networks and Approximation Theory}\indexme{feedforward network (FNN)!approximation theory}
%
The subject of approximation theory is concerned with how functions can best be approximated with simpler functions, and how well such approximations represent the original function.  Viewed this way, the topic of Fourier series can be considered part of approximation theory.  However, the kinds of functions that can be used in approximation theory can be quite different than the smooth trigonometric functions with which we are already familiar.  While the topic of approximation theory as a whole is far to extensive to summarize here, concrete examples can be read about in the literature.  For example, collections of Gaussian functions can be combined to approximate other functions as described by the theory of radial basis functions.  Wavelets, which can be thought of as a generalization of Fourier transforms that emphasise the local structure of the functions they represent, are another example of functions used in developing approximations of other functions.  The primary questions that arise when discussing these methods are (1) whether or not the proposed functions are sufficiently \emph{dense} (by some measure) to represent the kinds of functions that are being approximated, and (2) whether or not the basis functions being investigated are orthogonal (like Fourier series).  The answer to the first of these has been the topic of many papers on approximation theory by ANNs (see \citep{pinkus1999approximation} for a summary of that literature); the general consensus is that $\relu$ functions are indeed sufficiently dense to represent $L_2$ functions defined on compact intervals (or volumes or for dimensions greater than 1).  The second question is an interesting one.  It has been known for some time that functions like the Gaussian and even compact \emph{bump} functions can form a dense basis for $L_2$ functions under very broad conditions.  However, constructive methods for doing so have been lacking, primarily because the methods are so under-constrained. 

Neural networks escape this problem by not attempting to construct necessarily unique solutions.  In other words, when solved numerically, the objective is often simply to find \emph{some} solution that meets a desired error condition.  Because of the unique way that gradient descent (and a related methods called stochastic gradient descent) work, these methods are able to find solutions even when the landscape of potential solutions is enormous.

\begin{figure}[!ht]
            \sidecaption[t]
           \centering
            \includegraphics[scale=0.55]{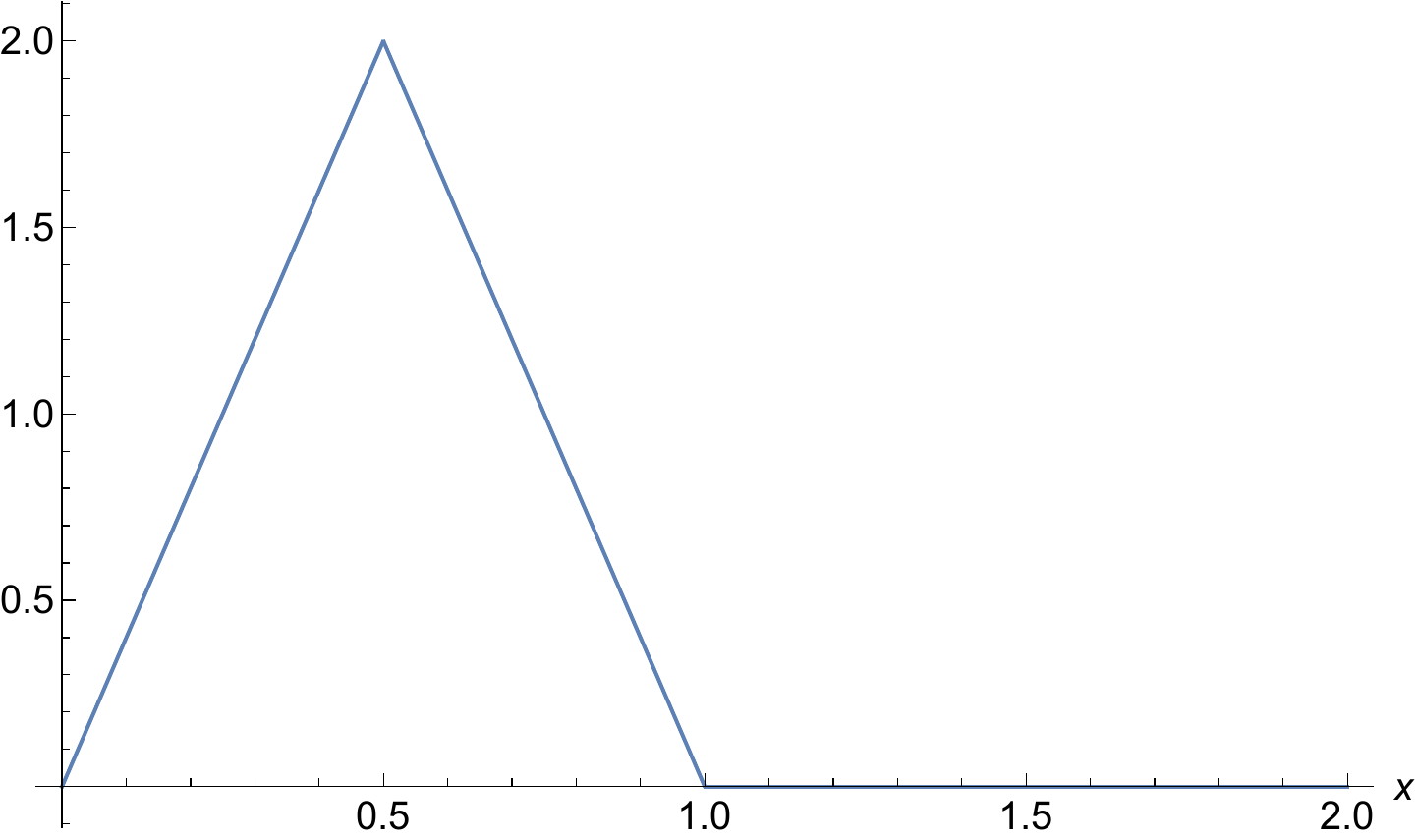}
            \vspace{-5mm}
            \caption{A symmetric compact function generated by the sum of ReLU functions $f[x]= 4 \reluf{x} - 8\reluf{x-1/2}+4 \reluf{x-1}$.}
            \label{fig:relutriangle}       
            \end{figure}

One recent method that has been adopted to examine the abilities of FNNs in particular is to compare them with existing methods for approximation.  For example, one can construct triangle functions using a combination of three $\relu$ functions as follows

\begin{equation}
    f[x]= 4 \reluf{x} - 8\reluf{x-1/2}+4 \reluf{x-1}
    \label{compactrelu}
\end{equation}
This function is plotted in Fig.~\ref{fig:relutriangle}.

One correspondence one can make is with the conventional finite-element constructions that have been well studied in the literature.  Although many functions are used in practice, the simple triangle function is the most common.  In terms of $\relu$ functions, this is

This function is useful because it indicates that, in some senses, $\relu$ is equivalent to a compact basis function for linear finite elements.  In fact, it has been shown that there are equivalences between finite element basis functions and $\relu$ functions \citep{he2022relu}.  A paper by \citet{telgarsky2015representation} provides an analysis of error bounds for compositions of sawtooth functions that correspond to deep $\relu$ networks.  

\begin{figure}[!ht]
            \sidecaption[t]
           \centering
            \includegraphics[scale=0.75]{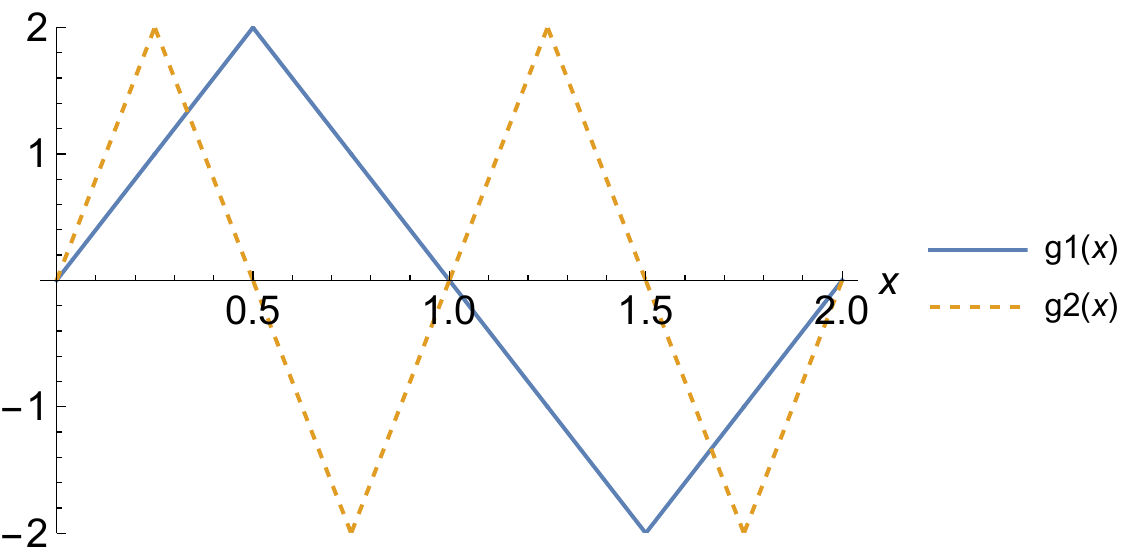}
            \vspace{-5mm}
            \caption{Orthogonal functions constructed from $relu$ functions.}
            \label{fig:reluorth}       
            \end{figure}

Another approach that has been by researchers is to illustrate how such functions can be used to create \emph{partitions of unity} (which is a finite collections of functions on $x\in[0,1]$ such that the sum of the functions at every point in the domain is 1).  There are powerful analytical frameworks for such constructions, and these represent one possible way to illustrate how well networks of $\relu$ functions can represent high-dimensional data sets \citep{yarotsky2017error,poggio2017and}.  Yet another approach that has been pursued is to construct a Fourier-like orthogonal functions from $\relu$ functions, and then show that these have good expressive properties \citep{daubechies2022nonlinear}.  For example, consider the compact function $f(x)$ given by \eqref{compactrelu}.  Now consider

\begin{align}
    g1(x) &= f(x)-f(x-1)\\
    g2(x) & =g1(2x)+g1(2(x-1))
\end{align}
These functions are plotted in Fig.~\ref{fig:reluorth}.  It is clear that these functions are orthogonal.  Functions constructed in this manner a generalization of Fourier series, and are known as Riesz bases.  More information on such constructions can be found in the literature \citep{daubechies2022nonlinear}.

%
While these examples are far from ``proofs" that $\relu$ functions can effectively reproduce many continuous functions, it does provide some intuition about how such a goal might be pursued.  There is still much that remains unknown about the expressive power of deep $\relu$ networks, but this remains an active area of investigation.  A paper by \citet{daubechies2022nonlinear} provides a particularly insightful analysis along the lines discussed above.

\subsection{~Ridge Functions}\label{ridgefunct}\indexme{ridge functions}\indexme{ReLU function!ridge functions}

Before finishing this section on the $\relu$ function and its role in approximation, it will be useful to look at a particular class of functions that can be constructed from higher-dimensional $\relu$ functions: These are known as \emph{ridge functions}.  Ridge functions are linear weighted sums of functions that form a higher dimensional result analogous to the 1-dimensional triangle function.  An example of two ridge functions in 2-dimensions is given below, with their plot appearing in Fig.~\ref{fig:ridge}.  

The characteristic feature of ridge functions is that, in an appropriate coordinate system, they become a function of fewer independent variables than define the space itself.  To be clear, examine Fig.~\ref{fig:ridge}.  If we were to rotate the coordinate system by $\pi/2$ (with a plus or minus sign, depending on the case), then one of the coordinate axes would align with the folds along $x_1+x_2-1=0$ for $r_1(x-1,x_2)$, and $x_1-x_2=0$ for $r_2(x_1,x_2)$.  In that case, these functions would cease dependence upon the variable aligned with the fold, and would depend only upon the variable perpendicular to the fold.  In the literature, this is sometimes described as finding a lower-dimensional \emph{manifold} (here, ``manifold" is just a generalization of a subset of Euclidian space) embedded in the higher dimensional one.  In other words, we would then have a one-variable representation of our function, even though it is embedded in a two-dimensional space.  Some of the current theories about how ANNs work as well as they do involve identification of such sub-manifolds.  Our purpose for illustrating these functions, however, has nothing in particular to do with manifolds.  Later, we will find a specific use for these functions in our efforts to determine a neural network for learning the exclusive or function.

            \begin{figure}[!ht]
            \sidecaption[t]
           \centering
            \includegraphics[scale=0.70]{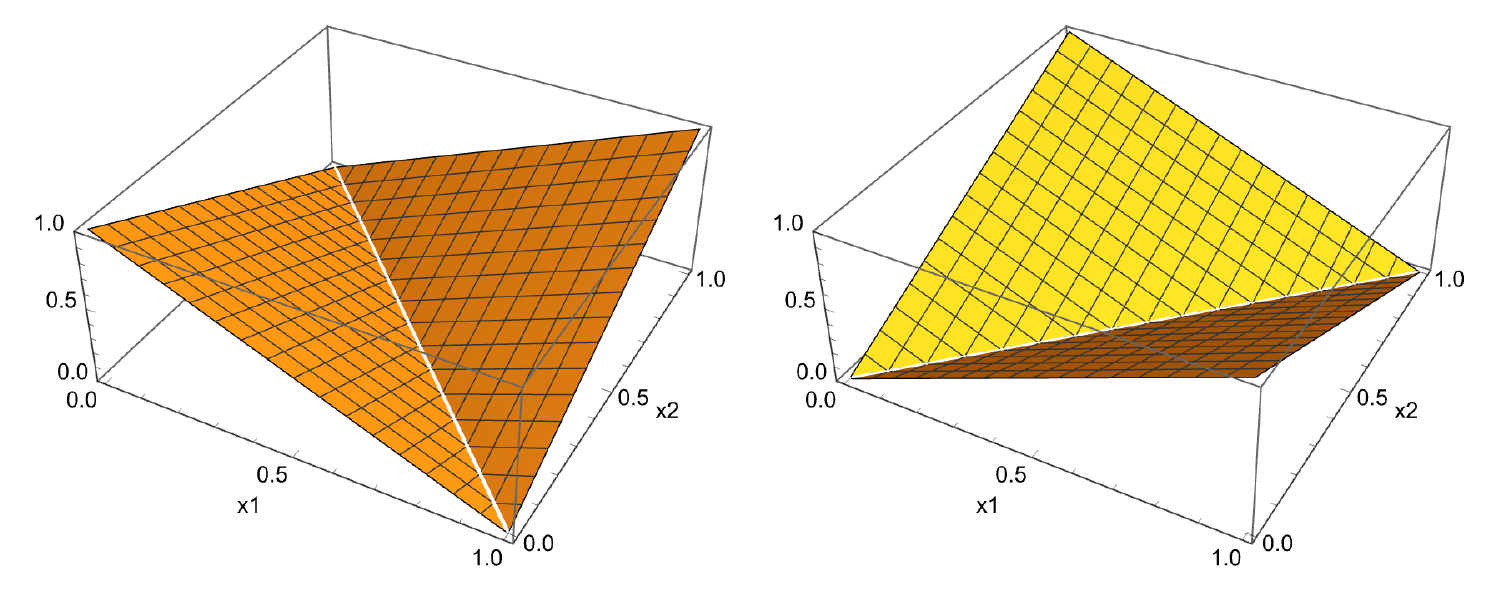}
            \vspace{-5mm}
            \caption{Two ridge functions.  The expressions for these functions appears in the text.  The function $r_1(x_1,x_2)$ is on the left, and the function $r_2(x_1,x_2)$ is on the right.}
            \label{fig:ridge}       
            \end{figure}
The expression for these two ridge functions is given as follows

\begin{align}
    r_1(x_1,x_2) &= \reluf{x_1+ x_2-1}+\reluf{-x_1-x_2 +1}\\
    r_2(x_1,x_2) &= \reluf{-x_1+x_2}+\reluf{x_1-x_2}
\end{align}
These examples illustrate two shapes that can be useful in applications.  Later, we will have an opportunity to construct a neural network for the exclusive or, $XOR$, function.  This function has been offered as an example that cannot be constructed with a one-hidden-layer neural network model; we will have the opportunity construct a one-layer model that illustrates that this is untrue.

\section{~Feedforward Networks for \uline{Nonlinear} Regression}\indexme{feedforward network (FNN)!nonlinear regression}

Now that we have outlined an appropriate description of the $\relu$ function, it is possible to consider constructing some example networks that perform specified functions.  In this section, we examine one- and two-hidden layer FNNs as they are applied to data sets.

\subsection{~Single-Hidden-Layer FNNs Using $\relu$}

The single hidden layer FNN has many analogues to the method of linear splines with free knots \citep{devore1998nonlinear}.  In fact, for single-layer $\relu$ FNNs, the two are essentially the same \citep{plonka2022spline}.  While single layer $\relu$ networks do not have the same efficiency as deeper networks, they are a useful model to begin investigations because of their relative transparency.  In other words, it is possible to readily \emph{interpret} the results of a single-layer $\relu$ network.  As a start, the following is an example of determining the binary $AND$ function with a single-layer $\relu$ network.

\begin{svgraybox}
\begin{example}[Learning the AND function]

The binary AND function is an interesting example to examine using a single-layer $\relu$ network.  The $AND$ function is a binary function whose functional relations are specified by the following table

{   \centering 
\begin{tabular}{|c|c|c|}
     $~~x_1~~$    & $~~x_2~~$ & $~~y~~$ \\
     \hline
      0  & 0 & 0 \\
      0 & 1 & 0 \\
      1 & 0 & 0  \\
      1 & 1 & 1
    \end{tabular}
    \captionof{table}{The binary AND function.}
    \label{tab:ANDfunction}
    }
    \vspace{4mm}
%
It is not difficult to determine that the solution $y(x_1,x_2)=\relu(x_1+x_2-1)$ is a continuous function that matches the binary AND function at the points $(1,1), (0,1),(1,0)$ and $(1,1)$.  The question now is, can we determine this function by constructing an appropriate neural network and then optimizing the parameters.  To do so, we first draw out an appropriate network: This is done in Fig.~\ref{fig:AND}.

To start the process, we first convert our figure into the equivalent mathematical representation.  It is not difficult to work out the representation for this simple network, it is

\begin{equation}
    \hat{y}(x_1,x_2) = \relu\left(w_{1,4} x_1 + w_{2,4}x_2 + b_{3,4}\right)
\end{equation}
Now, we have a number of choices for the \emph{loss} function for this problem.  However, the $L_2$ loss function (least-squares) is convenient.  The loss function can be represented by the difference between our \emph{data} and the predicted function, $\hat{y}$ at the same values of the independent variables for which we have data.  Thus, our loss function is given by the sum

\begin{align}
    J(x_1,x_2,w_{1,4}, w_{2,4},b_{3,4}) = [0-\hat{y}(0,0)]^2+[0-\hat{y}(1,0)]^2+[0-\hat{y}(0,1)]^2+[0-\hat{y}(1,1)]^2
\end{align}
This problem is simple, and it can be solved uniquely without imposing any additional constraints.  However, the $\relu$ function creates a number of branching options in the loss function without additional constraints.  Or, in other words, the loss function is described by four piecewise function statements, each of which has a zero portion an a positive portion arising from the application of the $\relu$ function.  To be concrete about this, the following results from the expansion of the loss function

\begin{align}
   J(x_1,x_2,w_{1,4}, w_{2,4},b_{3,4}) &=  \Bigg[1-\Bigg(
\begin{cases}
 0 & b_{3,4}+w_{1,4}+w_{2,4}<0 \\
 b_{3,4}+w_{1,4}+w_{2,4} & b_{3,4}+w_{1,4}+w_{2,4}\geq 0 
\end{cases}
\Bigg)\Bigg]^2
\nonumber\\
&+\Bigg(
\begin{cases}
 0 &b_{3,4}+w_{1,4}<0 \\
 b_{3,4}+w_{1,4} & b_{3,4}+w_{1,4}\geq 0 \\
 \end{cases}
\Bigg)^2
+\Bigg(
\begin{cases}
 0 & b_{3,4}+w_{2,4}<0 \\
 b_{3,4}+w_{2,4} & b_{3,4}+w_{2,4}\geq 0 \\
\end{cases}
\Bigg)^2 
\nonumber \\
&+\Bigg(
\begin{cases}
 0 & b_{3,4}<0 \\
 b_{3,4} & b_{3,4}\geq 0 \\
\end{cases}
\Bigg)^2
\end{align}
This equation has 4 quadratic terms, and each such term involves 2 domains each; this yields a potential 16 different combinations of cases (counting any redundant cases) that would need to be checked!  However, judicious implementation of constraints can reduce this to essentially zero.  Consider the following constraints.  First, we know at the point $(0,0)$, the only contribution to the function is $b_{34}$.  Because the function should be zero at this point, clearly $b_{3,4}\le 0$ (so that when the $\relu$ function is applied to it, the result is zero).  Also, the AND function is symmetric around the line $x_1=x_2$, that is to say $y(0,1)=y(1,0)$. The network defining it is also symmetric.  Therefore we must have $w_{1,4}=w_{2,4}$ (there is no mechanism to break the symmetry, therefore the symmetric input data must lead to a symmetric solution).  Finally, we know that $y(1,1)=1$, so we must have that $w_{1,4}=w_{2,4}>0$.  

{\hspace{14mm}
\centering\includegraphics[scale=.4]{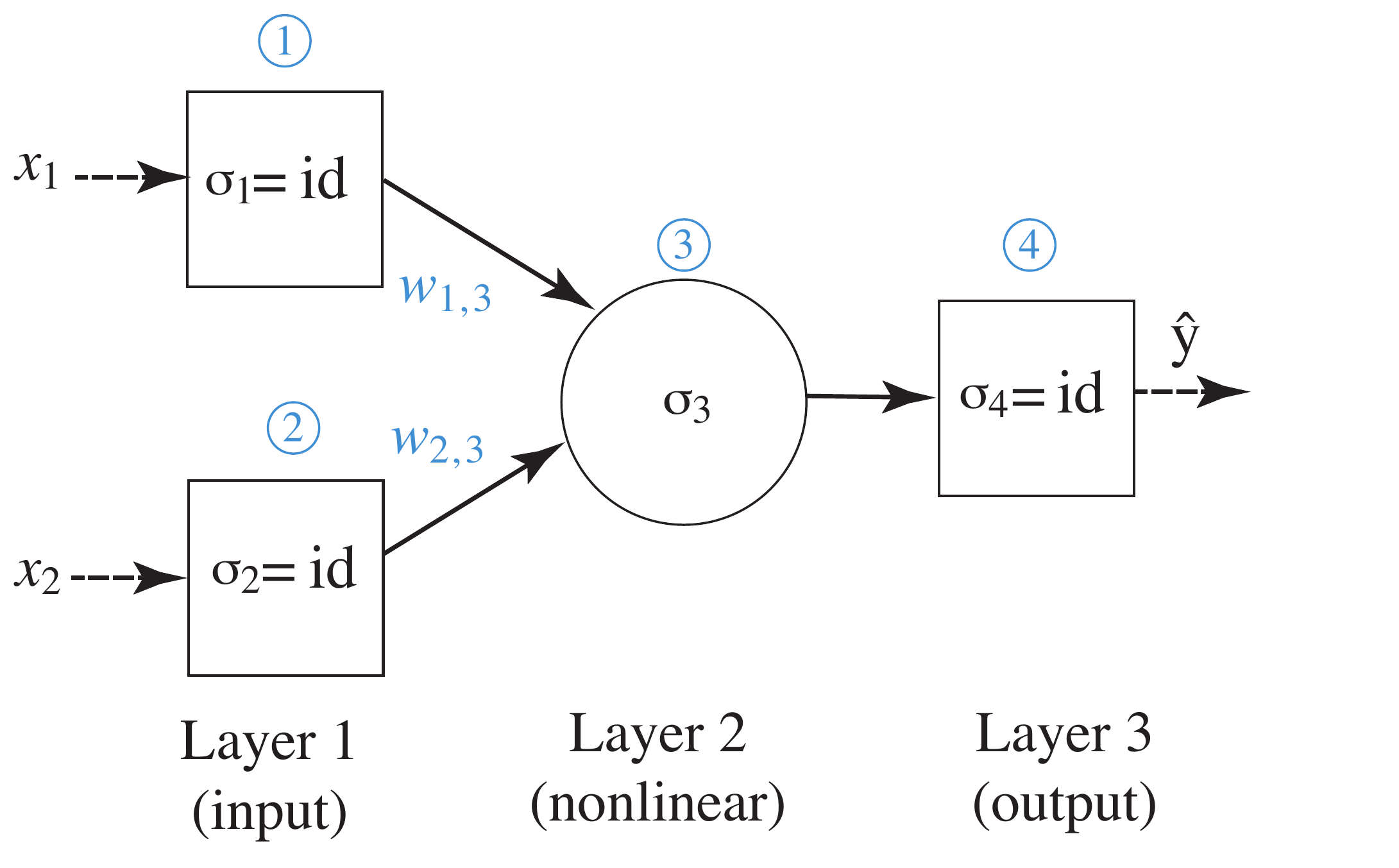}
\vspace{-2mm}
\label{fig:AND}}
\captionof{figure}{Network for determining the AND function.\vspace{4mm}}

To optimize the system, we need only find the point in the parameter space where the derivative with respect to each of the independent variables $(w_{1,4},w_{2,4}, b_{3,4})$, each give zero.  With the additional constraints in place, we have only a single expression for each derivative.  The results are as follows

\begin{align}
    a_1(w_{1,4},w_{2,4}, b_{3,4})&=\frac{\partial J}{\partial w_{1,4}} =2 (-1 + b34 + w14 + w24) \\
    a_2(w_{1,4},w_{2,4}, b_{3,4})&=\frac{\partial J}{\partial w_{2,4}} =-2 + 6 b34 + 8 w24 \\
    a_3(w_{1,4},w_{2,4}, b_{3,4})&=\frac{\partial J}{\partial w_{2,4}} =-2 + 6 b34 + 8 w24 
\end{align}
Because $a_2$ and $a_3$ are identical, we must use as our third equation the constraint $w_{1,4}-w_{2,4}=0$.  This then provides us three linear equations in three unknowns.  The solution is straightforward, and will not be detailed. The solution is 

\begin{align}
    w_{1,4} & = 1\\
    w_{2,4} &= 1\\
    b_{3,4} &= -1
\end{align}

Our final solution for the AND function is given by the following expression

\begin{equation}
    \hat{y}(x_1,x_2) = \relu\left(x_1+x_2-1 \right)
\end{equation}
It is easy to determine that this solution is correct for the four possible input values for $(x_1,x_2)$.

As a final note, again it is stressed that this problem was solvable \emph{without} the imposition of the additional constraints.  That route would have led to a solution, but each derivative function would have been defined by four separate intervals defining each piecewise derivative.  If none of the resulting intervals overlapped, this could have meant up to 64 individual cases to check.  If multiple solutions were found that minimized the loss function, then we would select the result that was the \emph{global} minimum.  Thus, the use of constraints helped us solve this problem in a simple manner; however, the brute-force approach would have also been successful!

\end{example}
\end{svgraybox}
\vspace{4mm}
As a second example, we illustrate the problem of fitting  a set of nonlinear data that has no a priori \emph{known} functional behavior other than the data that it contains.  This example is somewhat of an artifice-- three data points are fit with two $\relu$ functions, yielding a solution that simply connects the three points via straight lines.  Obviously we do not need an FNN to accomplish this task.  However, it again gives us a tractable example that helps further illustrate the process of generating and optimizing FNNs.

\begin{svgraybox}
    \begin{example}[Interpolating data points linearly with $\relu$ functions]
        
Suppose we have the following, exceptionally simple, data set, ${\bf R}$, with one independent (feature) and one dependent (target) dimension specified by the ordered pair $(x,y)$ as follows

\begin{equation}
X = \{ (0,0), (\tfrac{1}{2}, 1), (1,\tfrac{3}{4}) \}
\end{equation}
The \emph{domain} of this data set is $x\in[0,1]$.   Our task is to fit these data using two $\relu$ functions; this is an \emph{interpolation} rather than a regression (i.e., each two adjacent points are connected by a line segment), and thus should give a fit with zero error.   

Given any two points, there is some $\relu$ function that will intersect both points.  However, because of the existence of general translations and vertical displacements, this solution is not unique.  This can be seen in Fig.~\ref{fig:nonunique}.  Thus, right from the start we should be aware that we may need to add \emph{some} constraints if we are to obtain a unique solution.  When optimizing the problem numerically, the initial state for the slopes and translations is usually done as random numbers.  Thus, even without constraints, the solution may converge to \emph{some} combination that minimizes the error.  The fact that it is not the \emph{only} solution is often immaterial.  However, when approaching the problem analytically, we do not have the same luxury!

Given that any two points can be fit with a single $\relu$ function, then for our case (three points) we need two $\relu$ functions to exactly fit the data. Our first step is to write out an appropriate network that can be used to solve the problem.  In Fig.~\ref{fig:ann1} we have drawn out such a network.  In drawing out this figure we have already imposed one constraint.  Note that in layer 2 there is no addition of bias terms.  If you take a moment to think about what this does, it forces the result to be a sum of two $\relu$ functions with no vertical translations; in other words, we have implicitly set $b$ for the final linear composition at node 5 to be zero. This dispenses with the non-uniqueness problem that was illustrated in Fig.~\ref{fig:nonunique}.  As we will see, the resulting solution is now a unique one.


            {\centering
            \includegraphics[scale=0.75]{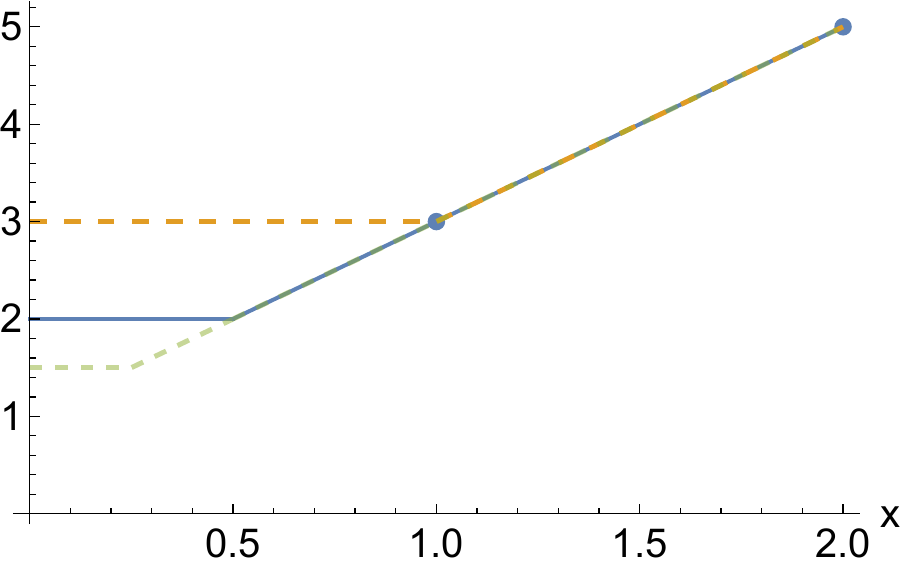}
            }
            \vspace{4mm}
            \captionof{figure}{Nonuniqueness in solutions of $\relu$ functions through two points.  Additional constraints are required to obtain a unique solution.
            }
            \label{fig:nonunique}       
            \vspace{4mm}

Translating the network graph into a mathematical statement, we obtain the following

\begin{equation}
    F(x; w_{1,3},w_{1,4},b_{2,3},b_{2,4},w_{3,5},w_{4,5}) = w_{3,5} \reluf{w_{1,3} x-b_{2,3}}+ w_{4,5} \reluf{w_{1,4} x-b_{2,4}}
    \label{ANN1}
\end{equation}
We can use the properties of $\relu$ functions and what we know about the data to impose an additional constraint that makes the system simpler.  Recall that the linear combination of two $\relu$ function with unspecified coefficients can always be re-expressed with a weighting of unity, and and unknown sign function.  In other words, we can write the function above as 

\begin{equation}
    F(x; w_{1,3},w_{1,4},b_{2,3},b_{2,4},\eta_1,\eta_2) = {\eta_1} \reluf{w_{1,3} x-b_{2,3}}+ {\eta_2} \reluf{w_{1,4} x-b_{2,4}}
\end{equation}
where $\eta_1$ and $\eta_2$ are either $+1$ or $-1$.  We can say a bit more.  A quick survey of the data indicates that the desired curve is not monotonic; that is, there must be a sign change in the slopes of the lines at the midpoint.  This requires that $\eta_1$ and $\eta_2$ be of opposite signs.  It does not matter which one is assigned which value.  If we impose the constraint $w_{1,3},w_{1,4}>0$, then we can rewrite our solution by

\begin{equation}
    F(x; w_{1,3},w_{1,4},b_{2,3},b_{2,4}) =  \reluf{w_{1,3} x-b_{2,3}}- \reluf{w_{1,4} x-b_{2,4}}
\end{equation}

Finally, for this simple problem there is a final constraint that must can be applied.  We know that our two $relu$ functions must have their origins shifted by $\beta=0$ and $\beta = 1/2$, respectively.
We can infer this from the data itself.  We desire to find two line segments that fit three points identically.  From the data itself, we know the \emph{domain for the two piecewise-linear functions} (let's call them $PW_1$ and $PW_2$). We have the following information from the data.

\begin{align}
    PW_1(x): \{x: x\in[0,\tfrac{1}{2}]\} \\
    PW_2(x): \{x: x\in[\tfrac{1}{2},1]\}
\end{align}
Note that the middle point, $x=\tfrac{1}{2}$ is in both domains, as it should be (this enforces continuity).  From this information, we note by inspection that  $\beta_{2,3} = 0$ and $b_{2,4} = \tfrac{1}{2}$.  
Substituting this into the function $F$ that we are attempting to parameterize, we find

\begin{equation}
    F(x; w_{1,3},w_{1,4}) =  \reluf{w_{1,3} x}- \reluf{w_{1,4} (x-\tfrac{1}{2})}
    \label{ANN3}
\end{equation}
Take note that we have rewritten the second $\relu$ function in the form where the translation appears explicitly.  Also recall, we have made two sets of constraints: (1) we have explicitly put in the translation terms (because our problem allowed no other choice), and (2) we have assumed $w_{1,3},w_{1,4}>0$ to get the form given by Eq.~\ref{ANN3}.

            {\centering
            \includegraphics[scale=0.45]{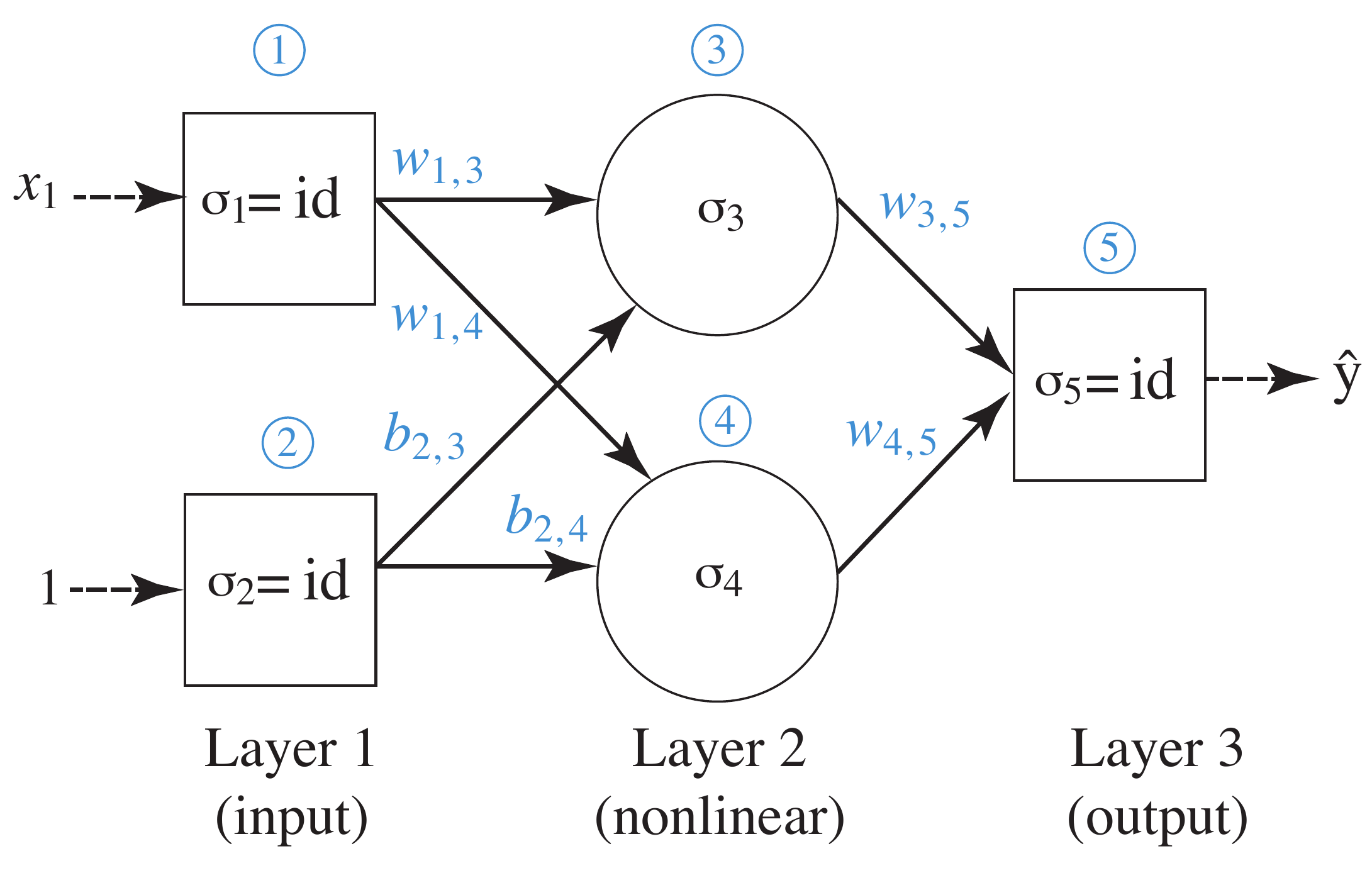}
            }\vspace{4mm}
            \captionof{figure}{The FNN associated with Eq.~\eqref{ANN1}.  Here, $\relu _1=\reluf{w_{1,3} (x - b_{2,3})}$ and $\relu _2=\reluf{\alpha_2 (x - b_2)}$
            \vspace{4mm}}
            \label{fig:ann1}       

To fit these data, as always, we need \emph{some metric} for defining the error.  In this particular case, the data form a $C^0$  (piecewise linear) function, and the $\relu$ basis functions are \emph{complete}, and are able to represent the piecewise linear function exactly, provided we can determining the values of $w_{1,3}$ and $w_{1,4}$ explicitly.  

Now, we adopt the $\ell_2$ metric as a measure the sum of the errors involved in our estimate at the three points of interest, i.e., for $N$ examples, the $\ell_2$ error is 

\begin{equation}
    J({\bf R}) =  \sum_{i=1}^{i=N} \left[y_i -\hat{y}({\bf X}_i;\boldsymbol{\theta}\right]^2
\end{equation}

\begin{align}
J & = (0-[\reluf{w_{1,3}\cdot 0}- \reluf{ w_{1,4}\cdot(0-\tfrac{1}{2})}])^2 \nonumber \\
&+(1-  [\reluf{w_{1,3}\cdot\tfrac{1}{2}}- \reluf{ w_{1,4}\cdot(\tfrac{1}{2}-\tfrac{1}{2})}])^2\nonumber \\
&+(\tfrac{3}{4}- [\reluf{w_{1,3}\cdot 1}- \reluf{w_{1,4}\cdot(1-\tfrac{1}{2}}])^2
\end{align}
This expression can be simplified using the properties of $\relu$ to give

\begin{align}
J & = (1-\tfrac{1}{2} w_{1,3})^2+(\tfrac{3}{4}-w_{1,3} + \tfrac{1}{2}w_{1,4})^2
\end{align}
To minimize the residual, we take the derivatives with respect to $w_{1,3}$ and $w_{1,4}$, and set them equal to zero.  This gives

\begin{align}
    \frac{\partial J}{\partial w_{1,3}}&= \tfrac{5}{2} w_{1,3}-w_{1,4}-\tfrac{5}{2} \\
    \frac{\partial J}{\partial w_{1,4}}&=   
    w_{1,3}-\tfrac{1}{2} w_{1,4}-\tfrac{3}{4} 
\end{align}
            {\centering
            \includegraphics[scale=0.35]{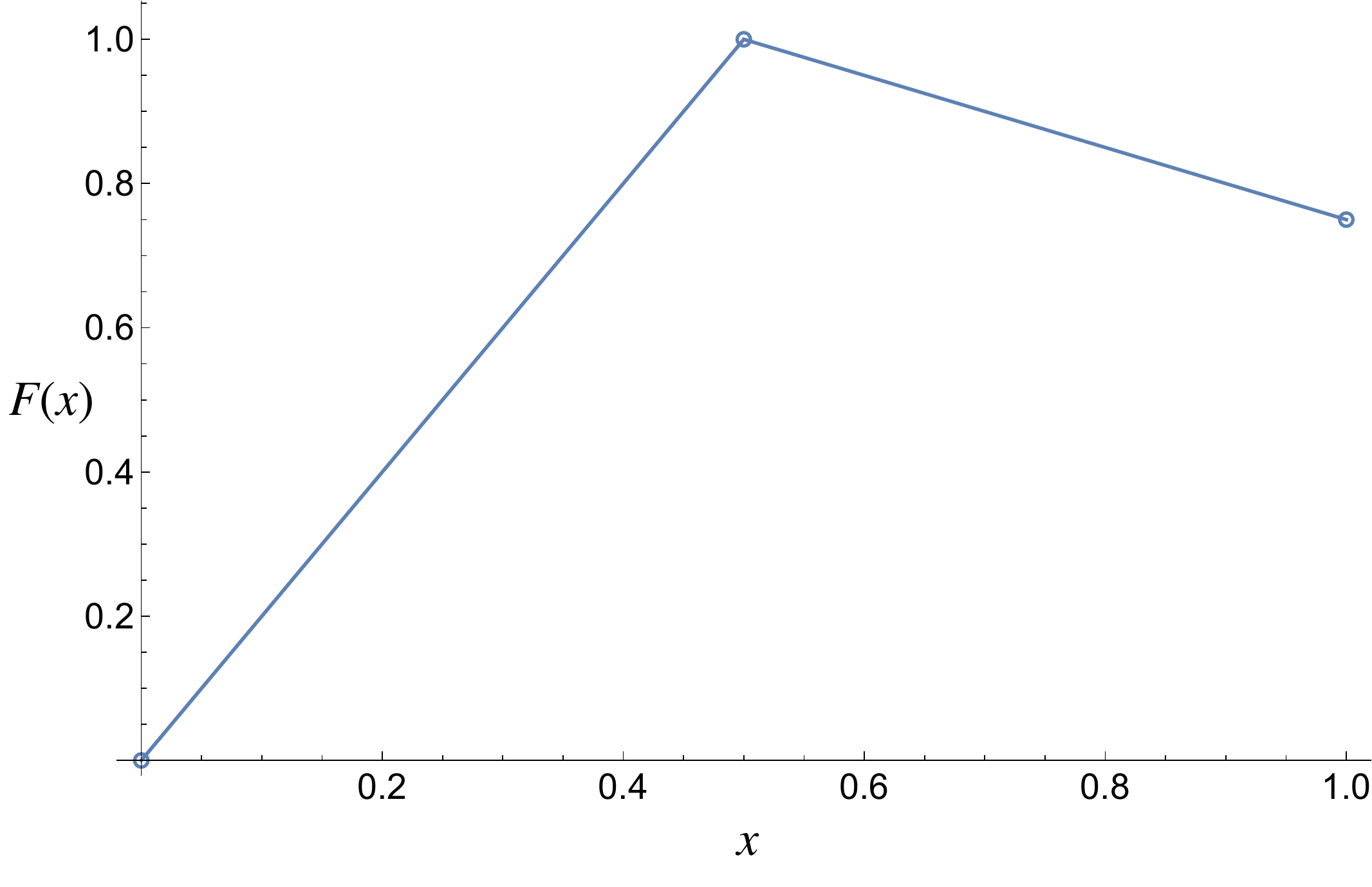}
            }
            \vspace{4mm}
            \captionof{figure}{An \emph{extremely simple} one-dimensional data set containing three points. \vspace{4mm}}
            \label{fig:datapoints2}       
%
At this point, we need only solve the following set of linear equations

\begin{align}
    \tfrac{5}{2} w_{1,3} + w_{1,4} -\tfrac{5}{2} &= 0 \\
   w_{1,3} + \tfrac{1}{2} w_{1,4} -\tfrac{3}{4} &= 0
\end{align}
The solution to this set of equations is given by $w_{1,3}= 2$ and $w_{1,4} = \tfrac{5}{2}$.  This gives the solution

\begin{equation}
 \hat{y}(x) = 2 \reluf{x}+ \tfrac{5}{2} \reluf{x-\tfrac{1}{2}}
\end{equation}
It is easy to verify that this is the has a residual error of $J=0$; the results are plotted in Fig.~\ref{fig:datapoints2}.  As a final note, the process outlined above is essentially the same if we have many \emph{more} data points than we do $\relu$ functions.  In that event, the optimization process results in the \emph{best possible} fit as measured by the least-squares ($L_2$) error metric.  However, this would be an example of \emph{regression} rather than \emph{interpolation}.  So, while it would be a best fit (meaning that the derivatives of the loss function would be identically zero for the minimizing set of parameters), the value of the error itself, $J$ would not be zero.  For the case of interpolation, the error, $J$, is identically zero.

    \end{example}
\end{svgraybox}

\begin{svgraybox}
\begin{example}[The XOR function via an inverted analysis.]

The $XOR$ function is a classic example that is often used as an example for machine learning.  This function is given by the following binary relation

{   \centering 
\begin{tabular}{|c|c|c|}
     $~~x_1~~$    & $~~x_2~~$ & $~~y~~$ \\
     \hline
      0  & 0 & 0 \\
      0 & 1 & 1 \\
      0 & 0 & 1  \\
      1 & 1 & 0
    \end{tabular}
    \captionof{table}{The binary AND function.}
    \label{tab:XORfunction}
    }
    \vspace{4mm}

In the context of \emph{classification} (which is usually how it is solved), a 2-hidden layer network is required in order to obtain a solution.  Even then, the problem is quite under-constrained, and it requires careful analysis to arrive at a correct result.  There appear to be existing analytical solutions to optimizing a network describing $XOR$. \\

In this example, we will find a method to approach a solution in a few steps. We will use our knowledge of \emph{ridge functions} that was introduced in section \ref{ridgefunct} to help inspire a solution. Previously, we examined the following ridge function.

\begin{equation}
    r_2(x_1,x_2) = \reluf{x_2-x_1} + \reluf{x_1-x_2}
\end{equation}
As a reminder, the plot of this function appears in Fig.~\ref{fig:XORplot}.  If one compares this function to the requirements of the XOR function, we find that this function represents \emph{one possible solution} (out of an infinite number of potential solutions) that is composed of only $\relu$ functions.  In fact, this solution is composed of exactly two $\relu$ functions, which is minimum possible number.

The XOR function has mirror-image symmetry, which is something that we can take advantage of in determining a solution.  This kind of symmetry means only that we can always exchange the variables $x_1$ and $x_2$ in our solution and leave the result unchanged.  This recognition builds in a significant number of constraints.  To see this, we start by writing out the expression for a generic one-hidden layer network with three input values $(x_1,x_2,1)$, two $\relu$ functions in the hidden layer, and one output, $\hat{y}$. A generic version of this network appears as Fig.~\ref{fig:genericANN}.

            {\centering
            \includegraphics[scale=0.85]{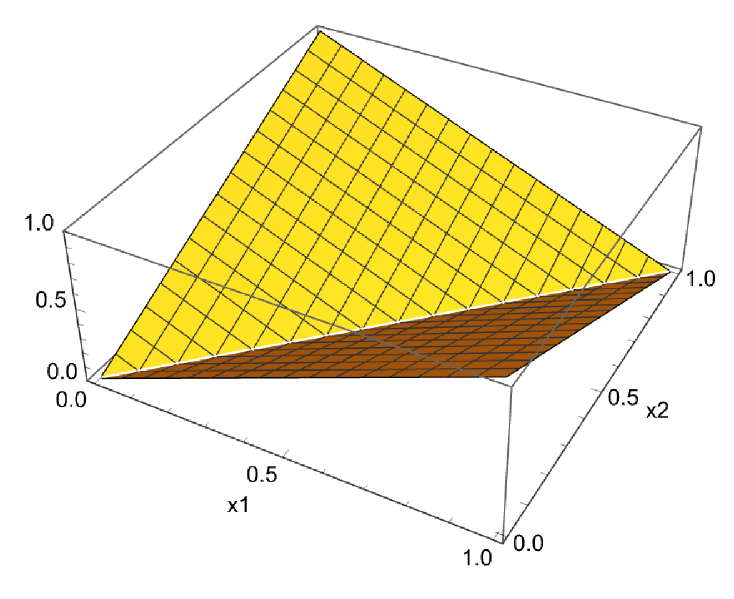}
            }
            \vspace{4mm}
            \captionof{figure}{A ridge function. \vspace{4mm}}
            \label{fig:XORplot}       
%
Using what we know about the XOR function, we can immediately determine that $b_{6,7}=0$ (this is from the point $x_1=0,~ x_2=0$ which has a solution $y=0$.  From the requirement of mirror-symmetry, we also have the condition that $w_{1,4}=w_{2,5}$ and $w_{2,4}=w_{2,5}$ (this is required so that we can always interchange $x_1$ and $x_2$, and the requirement that $w_{4,7}=w_{5,7}$. Because our solution is positive definite, and using the property that we can always replace the weights of the last layer of a $\relu$ network with a sign function, we can state that $w_{4,7}=w_{5,7}=1$.  Symmetry alone has helped us quite a lot!

This leaves us with the following expression for the network estimator, $\hat{y}$

\begin{equation}
    \hat{y}(x_1,x_2)= \reluf{w_{1,4} x_1+w_{2,4} x_2+b_{3,4}} + \reluf{w_{2,4} x_1 + w_{1,4} x_2+b_{3,5}}
\end{equation}
Again, using symmetry, and the solution for the point $x_1=0, x_2=0$ with $y=0$, we have that the two bias terms must equal, the value must be negative or zero.  To constrain these further, we can evaluate the function $\hat{y}$ at $(x_1,x_2)=(0,1),~ (1,0)$ and $(1,1)$.  We find the following options

\begin{align}
\textrm{either} & & w_{1,4}+b_{3,4}&>0, 
\textrm{~~or~~}  w_{2,4}+b_{3,4}>0 \nonumber \\
\textrm{and} && w_{1,4}+w_{2,4}+b_{3,4}&=0 \nonumber
\end{align}
Any solution that meets these constraints is valid.  The constraints are met by $b_{3,4}=0$, which is the simplest option.

            {\centering
            \includegraphics[scale=0.35]{chapter17_intro_to_the_mathematics_of_deep_neural_networks/figs/generic_network_simple.pdf}
            }
            \vspace{4mm}
            \captionof{figure}{A generic single-hidden-layer network. \vspace{4mm}}
            \label{fig:genericANN}       

Our estimator is now relatively simple.

\begin{equation}
    \hat{y}(x_1,x_2)= \reluf{w_{1,4} x_1+w_{2,4} x_2} + \reluf{w_{2,4} x_1 + w_{1,4} x_2}
\end{equation}
which contains only two unknown parameters.  The loss function is given by

\begin{align}
    J(w_{1,4},w_{2,4}) &= [0-\hat{y}(0,0)]^2+[1-\hat{y}(1,0)]^2 
    +[1-\hat{y}(0,1)]^2 + [0-\hat{y}(1,1)]^2
\end{align}
where
\begin{align}
    \hat{y}(0,0)&=0 \\
    \hat{y}(1,0)&= \reluf{w_{1,4}}+\reluf{w_{2,4}} = 1\\
    \hat{y}(0,1)&= \reluf{w_{2,4}}+\reluf{w_{1,4}} = 1\\
    \hat{y}(1,1)&=\reluf{w_{1,4}+w_{2,4}}+\reluf{w_{2,4}+w_{1,4}} = 0
\end{align}
In this case, we do not necessarily need to take derivatives of the objective function because the solution is so highly constrained (although that process \emph{will} lead to the correct result!).  The expression for $hat{y}(1,1)$ clearly indicate that $w_{1,4}=-w_{2,4}$.   In addition, either of $\hat{y}(1,0)$ or $\hat{y}(0,1)$ require that $w_{1,4}=1$ and $w_{2,4}=-1$ (or $w_{1,4}=-1$ and $w_{2,4}=1$, either choice is a solution).  The resulting network can be drawn out as given in Fig.~\ref{fig:XORnetwork}. \\

As a final note, we could have conducted a \emph{classification} scheme with this network by adding a layer before the output.  This layer would contain a single Heaviside function with one additional unknown parameter (the bias).  This would have led to a network with two hidden layers.   However, given that we have a solution, a little thought will indicate that filtering our solution \emph{a posteriori} through a Heaviside with a bias of $-1/2$ will accomplish the same result.  This is also plotted in Fig.~\ref{fig:XORnetwork}

            {\centering
            \includegraphics[scale=0.30]{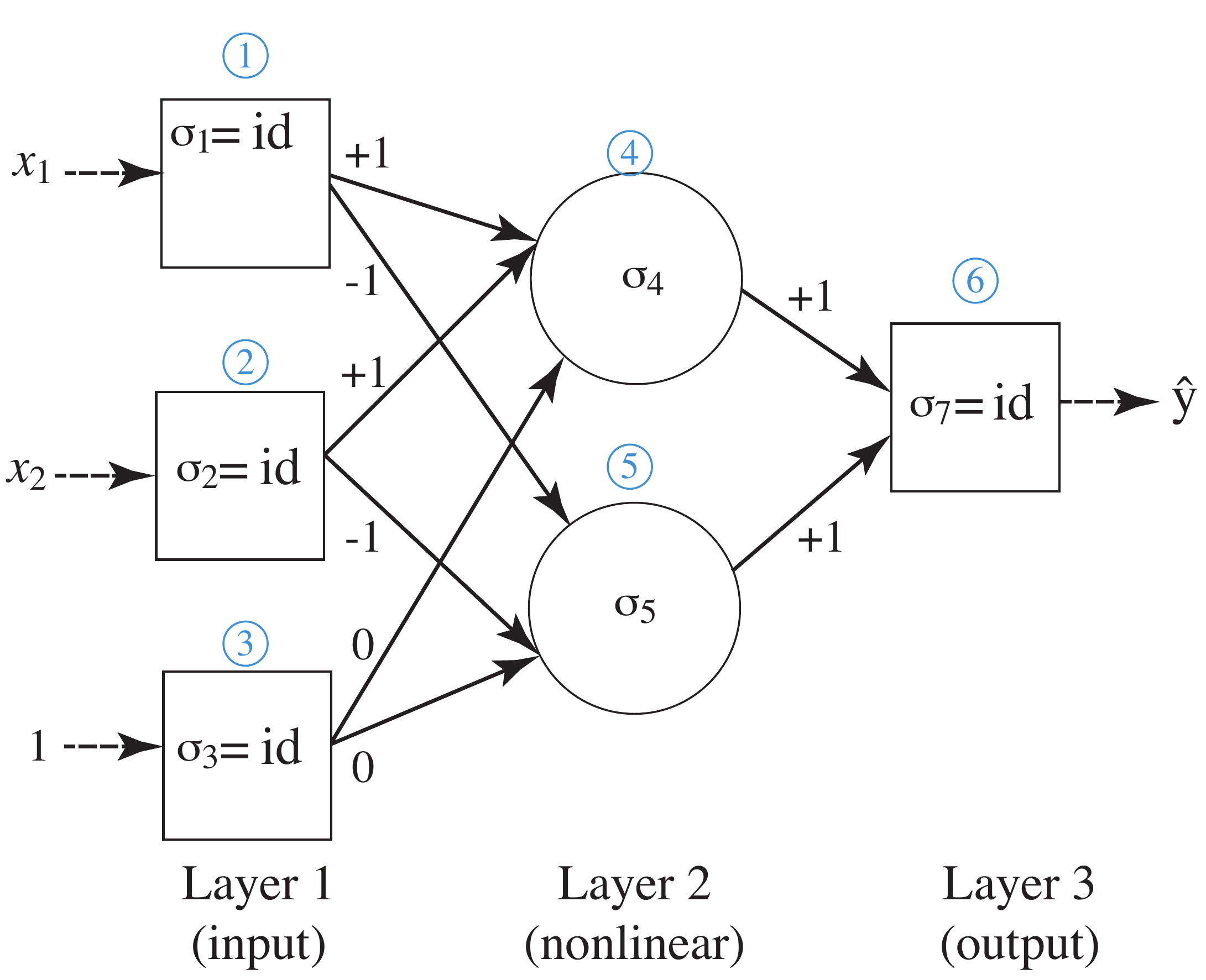}
            \includegraphics[scale=0.65]{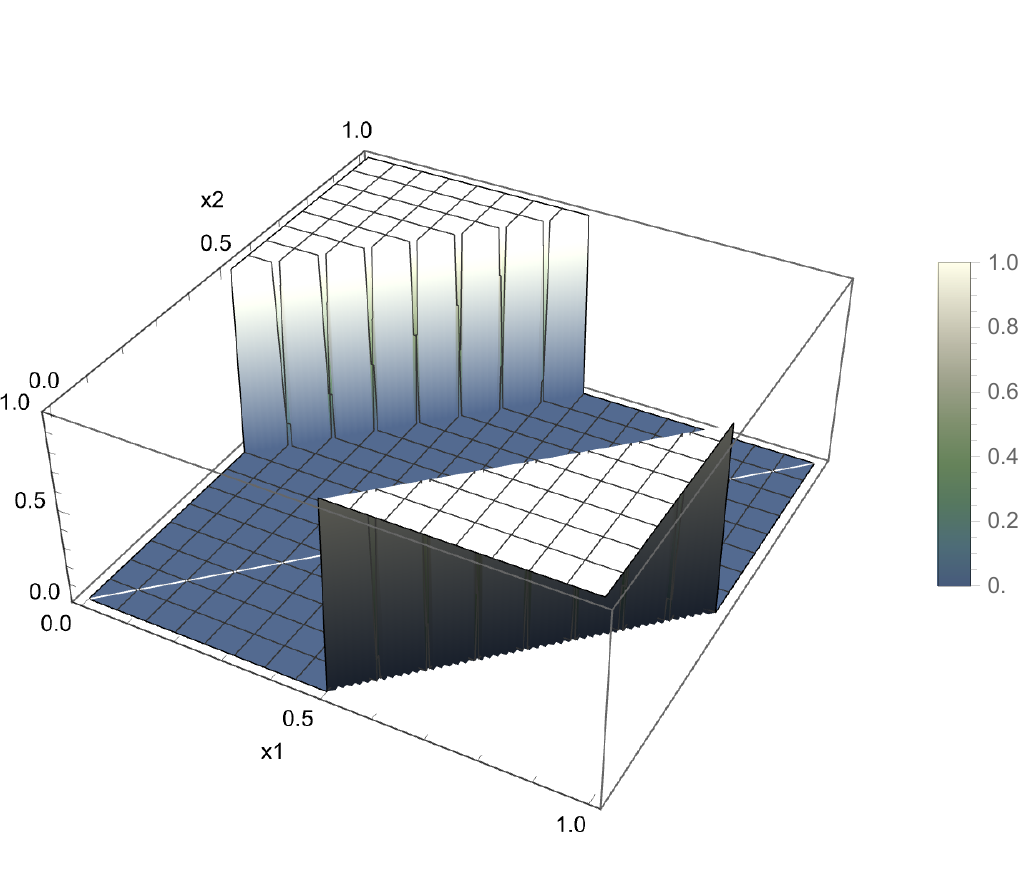}
            }
            \vspace{4mm}
            \captionof{figure}{(Left).  The XOR  network solution.  (Right). The solution filtered through the Heaviside function $H(\hat{y}-1/2)$. \vspace{4mm}}
            \label{fig:XORnetwork}       

\end{example}
\end{svgraybox}

\section{ Approximation of Known Functions with Nonlinear FNNs}
\indexme{feedforward network (FNN)!approximation of functions}

The process of using a weighted sum of simple functions to approximate a given function is part of approximation theory.  This has been mentioned in passing in the material presented previously. Here, we illustrate the process of approximating a continuous function with a piecewise-continuous $\relu$ network.  There are many reasons that one might want to approximate a given function by an approximating function; for example, it is often easier to conduct certain transformations on a piecewise-linear approximation than on the original function.   In this example, we focus on fitting polynomials with $\relu$ functions.  In particular, we examine fitting quadratics.

\begin{figure}[t]
            \sidecaption[t]
            \centering
            \includegraphics[scale=0.45]{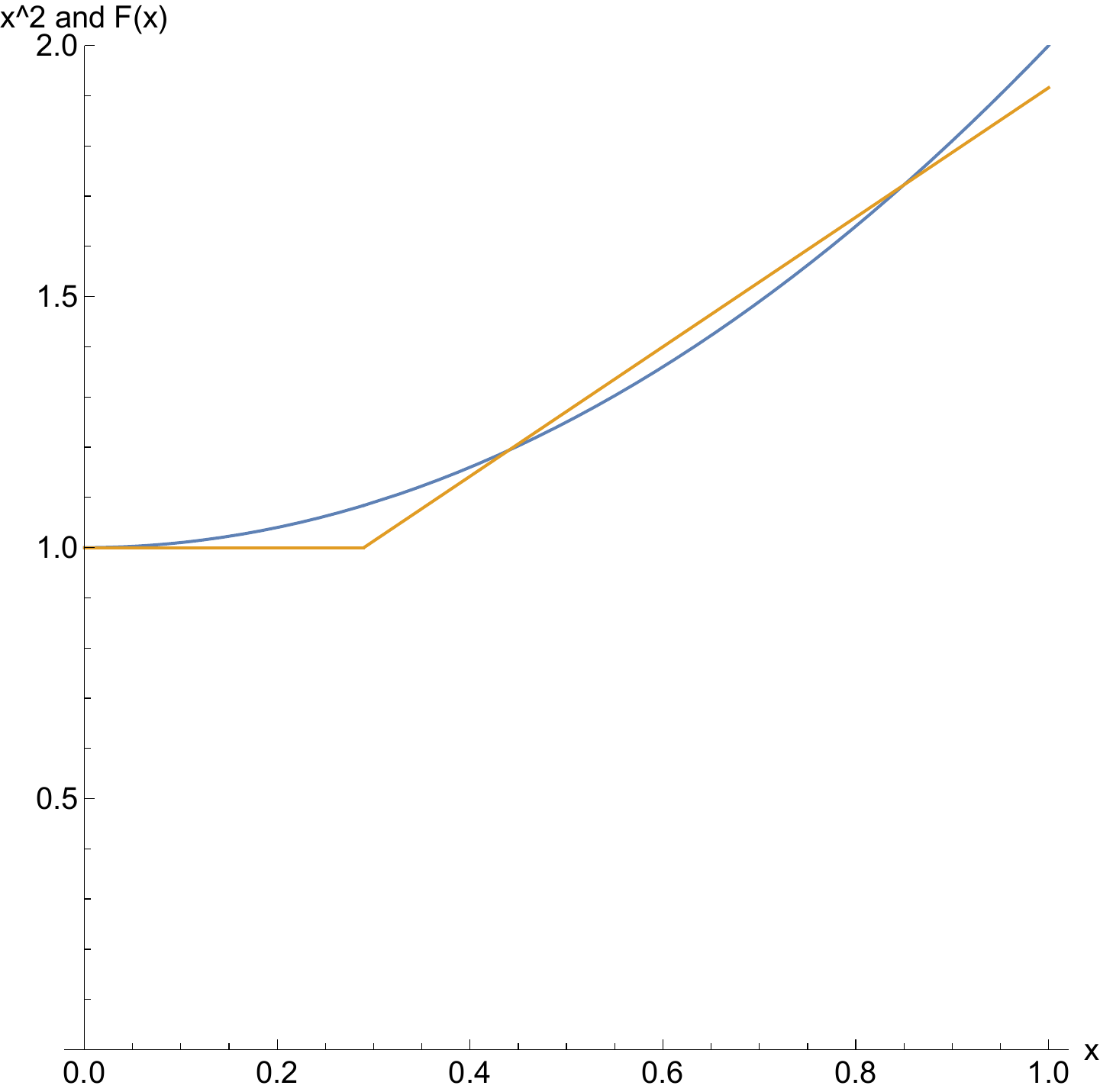}
            \caption{The solution for a single $\relu$ unit fitting the expression $f(x)=x^2+1$.  The solution represents the global minimum residual.}
            \label{singlrelu}       
            \end{figure}

As the quadratic of interest, take the function

\begin{equation}
    y_0(x) = x^2 +1, ~ x\in[0,1]
\end{equation}
This example is interesting because it also includes a global bias term, $B=1$, that shifts the function vertically by 1 unit.  This network is illustrated graphically as a network in Fig.~\ref{singlrelu}.

            \begin{figure}[!ht]
            \sidecaption[t]
           \centering
            \includegraphics[scale=0.45]{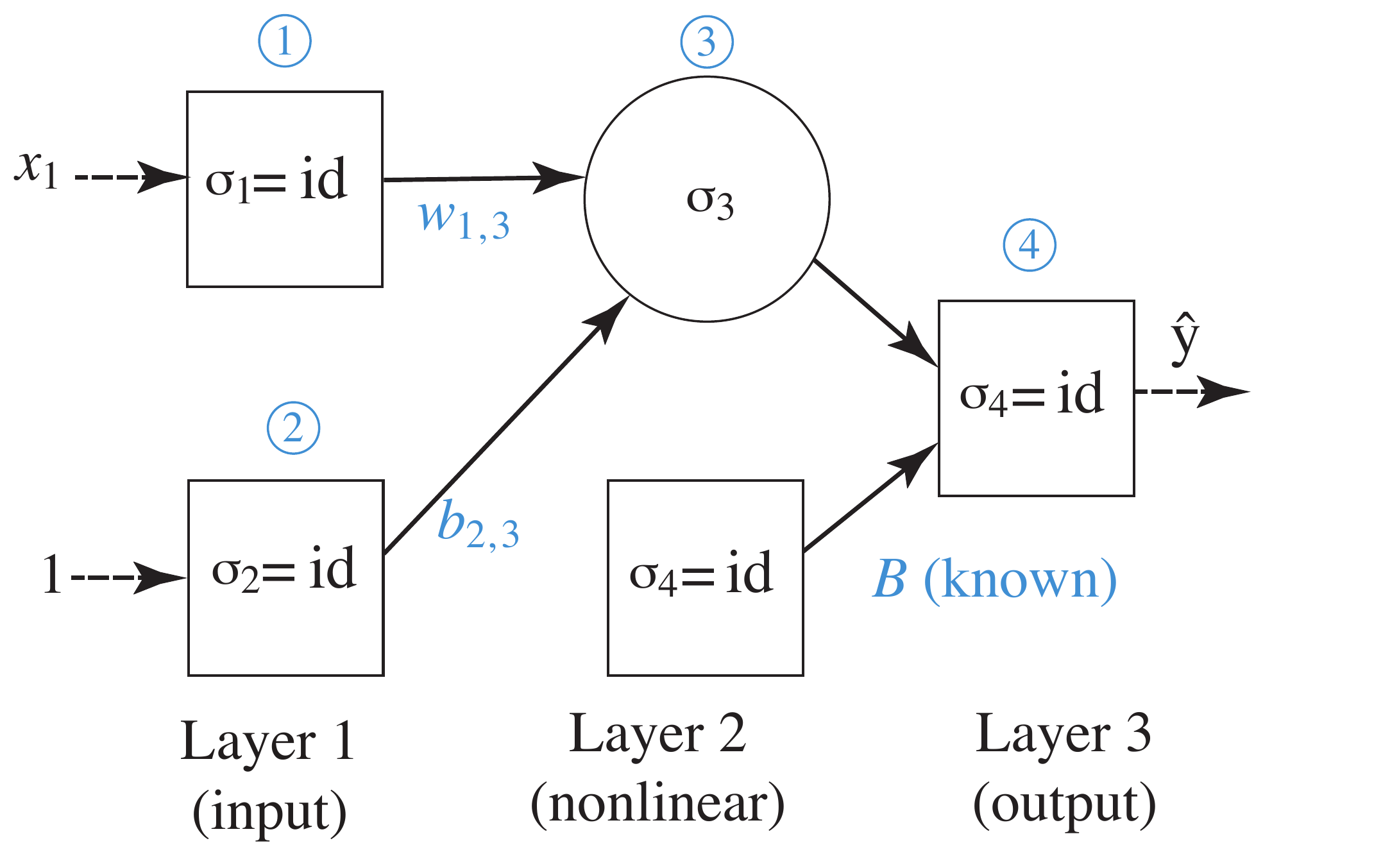}
            \vspace{-5mm}
            \caption{A single-hidden-layer, one-node network for polynomial approximation.}
            \label{fig:quadapprox}       
            \end{figure}

To begin, note that for a single-hidden-layer $\relu$ network, the bias term can always be handled by adding a bias term to the output node. In other words, we set

\begin{equation}
    B(x)= 1, ~ x\in[0,1]
\end{equation}
for the output node.  Here, we have given this constant a special symbol so that it does not depend on network numbering schemes.  From here forward, we \emph{assume} that this can be done for any constants added to a polynomial.  

The remaining problem is to fit the homogeneous polynomial
\begin{equation}
    y(x) = x^2
\end{equation}
Using $\relu$ functions.  To outline the process, consider the case where we use a single $\relu$ function to represent this polynomial.  Thus, we are approximating $y(x) = x^2$ by $\hat{y}(x)$, where

\begin{equation}
   \hat{y}(x) =  \reluf{w_{1,3}x+b_{2,3}}
\end{equation}
Because we are now in the context of continuous functions rather than discrete data points, our error function (or residual), $J$ is now an \emph{integral} rather than a sum

\begin{equation}
   J(w_{1,3}, b_{2,3}) = \int_{0}^1 \left[ x^2 -\reluf{w_{1,3} x-b_{2,3}} \right]^2 \, dx
\end{equation}
To help solve this problem, we note the following constraint, imposed by the domain width.

\begin{equation}
b_{2,3} <1
\end{equation}
This constraint prevents the contribution of the $\relu$ function from being identically zero.  Using the properties of the $\relu$ function, the residual, $J$ can be put in the form

\begin{align}
J(w_{1,3},b_{2,3}) &= \int_0^1 x^4\, dx -\int_{b_{2,3}}^1 2 w_{1,3} x^2 (x-b_{2,3}) \, dx +\int_{b_{2,3}}^1 w_{1,3}^2(x-b_{2,3})^2\, dx
\end{align}
Here, the dependence of the integral upon $b_{2,3}$ lower bound of the integration domain.  As is always helpful, we consider additional constraints on the function $\hat{y}$ (and, hence, the loss function) before attempting to optimize.  Again, it is \emph{possible} to continue without adding constraints a priori.  However, one will be faces with a multiplicity of possible solutions depending on linear combinations of parameter, and the number of such potential solutions can quickly become overwhelming.  If we can impose sensible constraints ahead of time, we avoid this problem.  However, we must assure that our proposed solutions do not unintentionally truncate the ``best" solution.  We impose the following constraints by inspection of the problem.  First, we assume that because the function $y(x)=x^2$ is strictly positive on $x\in[0,1]$, then we must have a positive solution.  This demands that $w_{1,3}>0$.  Second, we assume that the bias $b_{23}$ must be less than or equal to zero, otherwise the resulting prediction for $\hat{y}$ would be linear.  The problem is now stated

\begin{align}
J(w_{1,3},b_{2,3}) &= \int_0^1 x^4\, dx -\int_{b_{2,3}}^1 2 w_{1,3} x^2 (x-b_{2,3}) \, dx +\int_{b_{2,3}}^1 w_{1,3}^2(x-b_{2,3})^2\, dx \nonumber \\
&\textrm{Constraints}:~w_{1,3}>0,~~b_{2,3}\le 0
\end{align}
Incorporating these two constraints, the expression for the loss function be integrated to give the following

\begin{equation}
    J(w_{1,3},b_{2,3}) =
\frac{1}{5}-\frac{(b_{2,3}+w_{1,3})^2 \left(b_{2,3}^2-2 b_{2,3} w_{1,3} (w_{1,3}+1)+(3-2 w_{1,3})
   w_{1,3}^2\right)}{6 w_{1,3}^3}
\end{equation}

This is a nice result; our constraints have given us a single expression for the error function rather than a list of cases, and this dramatically simplifies the analysis.  Continuing forward, we need to determine the derivatives of the loss function.

\begin{equation}
    \frac{\partial J}{\partial w_{1,3}} =
\frac{1}{6} \left(\frac{3 b_{2,3}^4}{w_{1,3}^4}-\frac{2 b_{2,3}^3}{w_{1,3}^2}+6 b_{2,3}+4
   w_{1,3}-3\right)
\end{equation}

\begin{equation}
    \frac{\partial J}{\partial b_{2,3}} =
-\frac{2 b_{2,3}^3}{3 w_{1,3}^3}+\frac{b_{2,3}^2}{w_{1,3}}+2 b_{2,3}+w_{1,3}-\frac{2}{3}
\end{equation}
Setting the two derivatives equal to zero, and solving simultaneously gives a single solution that satisfies the constraints

\begin{align}
   \frac{3 b_{2,3}^4}{w_{1,3}^4}-\frac{2 b_{2,3}^3}{w_{1,3}^2}+6 b_{2,3}+4
   w_{1,3}-3 &=0\\
    -\frac{2 b_{2,3}^3}{3 w_{1,3}^3}+\frac{b_{2,3}^2}{w_{1,3}}+2 b_{2,3}+w_{1,3}-\frac{2}{3} &=0
\end{align}
Solving a set of polynomial equations like this is a nontrivial task.  One option is to conduct multivariate root finding for the system using numerical methods.  There are also structured approaches, such as the method of Gr\"obner bases \citep{adams2022introduction}, that can facilitate such analyses.  Here, analytical solutions have been determined using the symbolic mathematics package Mathematica.  The solutions are

\begin{align}
    w_{1,3} = \frac{1}{5} (4 + \sqrt{6}) && b_{2,3} = \frac{1}{25} (2 + 3\sqrt{6}) && J&=\frac{73-28\sqrt{6}}{3125}\approx  0.00141
\end{align}


\begin{svgraybox}
    \begin{example}[Two-hidden layer approximation]\label{ex:sine_approx}
The deeper that networks become, the more challenging it is to solve them analytically.  Generally, as $\relu$ FNNs get deeper, the number of individual regions goes up significantly.  Even for networks with a single independent variable and a constant number of nodes in the hidden layers, the rate of increase is is still a power of the number of hidden layers. 

In this example, we will approximate a compact function using a 2-hidden-layer network.  The function of interest is given by

\begin{equation}
    y(x) =
    \begin{cases}
        0 & x<0 \\
        \sin[\pi x/2] & 0\le x \le 2 \\
        0 & x > 2
    \end{cases}
\end{equation}
This is a positive function defined by a single lobe of the sine function spanning the interval $x\in[0,2]$; note, however, we do not require our solution to have the same domain. The optimal approximation might occur on a slightly larger domain, with the increasing component of the approximation starting slightly to the left of $x=0$ (this is the purpose of the bias, $b_{2,3}$ appearing in the model below).  One thing that makes this problem interesting is that it is symmetric about a vertical axis placed at $x=1$.  As before, establishing reasonable constraints ahead of time will be part of how we can resolve this problem successfully.  To begin, we will insist that our solution be symmetric.  Another recognition is that, whatever our solution, it must be some form of triangle function.  One $\relu$ function will capture the increasing part of the function to the left of $x=1$, and the other $\relu$ function will capture the decreasing part of the function to the right of $x=1$. This also, then, sets the translation for one of the two $\relu$ functions-- it must be set at $x=1$.  

With these constraints in place, we find that there is another constraint that is needed for consistency with those already defined.  The slope of the approximation $\hat{y}(x)$ must be equal in magnitude, but with opposite signs (otherwise the approximator would not be symmetric).  A little though will indicate that if the $\relu$ function approximates $\hat{y}(x)$ using a slope of $w$, then the second $\relu$ function must have a slope of $-w$.  Finally, there is no bias in the final solution because the both the triangle function approximation and the original function are both compact functions that go to zero at large values of $|x|$.  

            {\centering
            \includegraphics[scale=0.35]{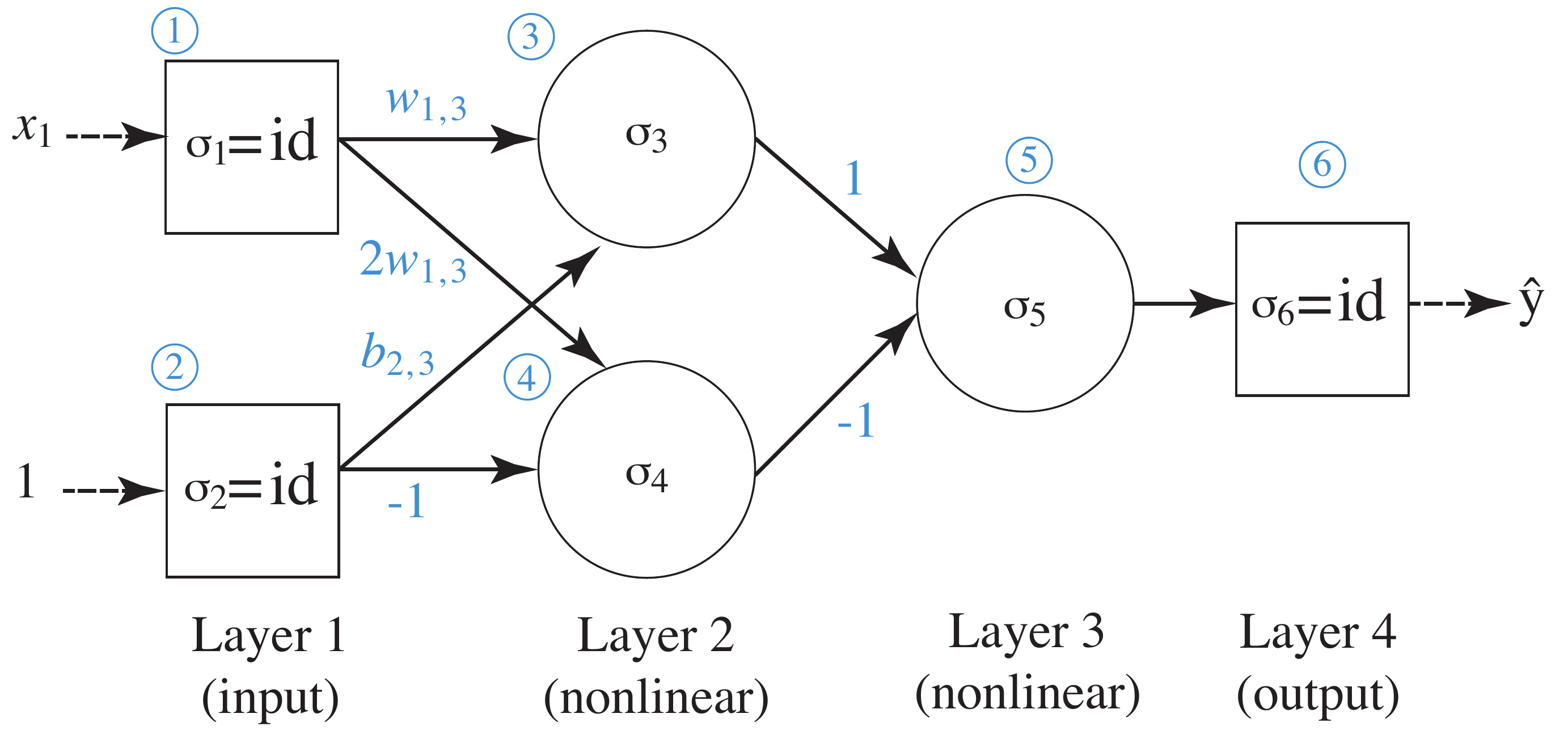}
            }
            \vspace{4mm}
            \captionof{figure}{A two-hidden-layer network graph.  Constraints have helped fix several of the parameters, and those are displayed over the appropriate link. \vspace{4mm}}

The mathematical model associated with this graph is given by

\begin{equation}
    \hat{y}(x) = \relu\left[
\reluf{w_{1,3} x+b_{2,3}}-\reluf{2w_{1,3}(x-1)}  
   \right]
\end{equation}
Note that this expression assumes that $w_{1,3}>0$ and that $b_{2,3}>0$; the minus sign before the second term is to enforce the slope equal to $-w_{1,3}$ for the portion of $\hat{y}(x)$ where $x>1$.  The factor of 2 on the $(1,3)$ link is necessary to attain this result.

            {\centering
            \includegraphics[scale=0.65]{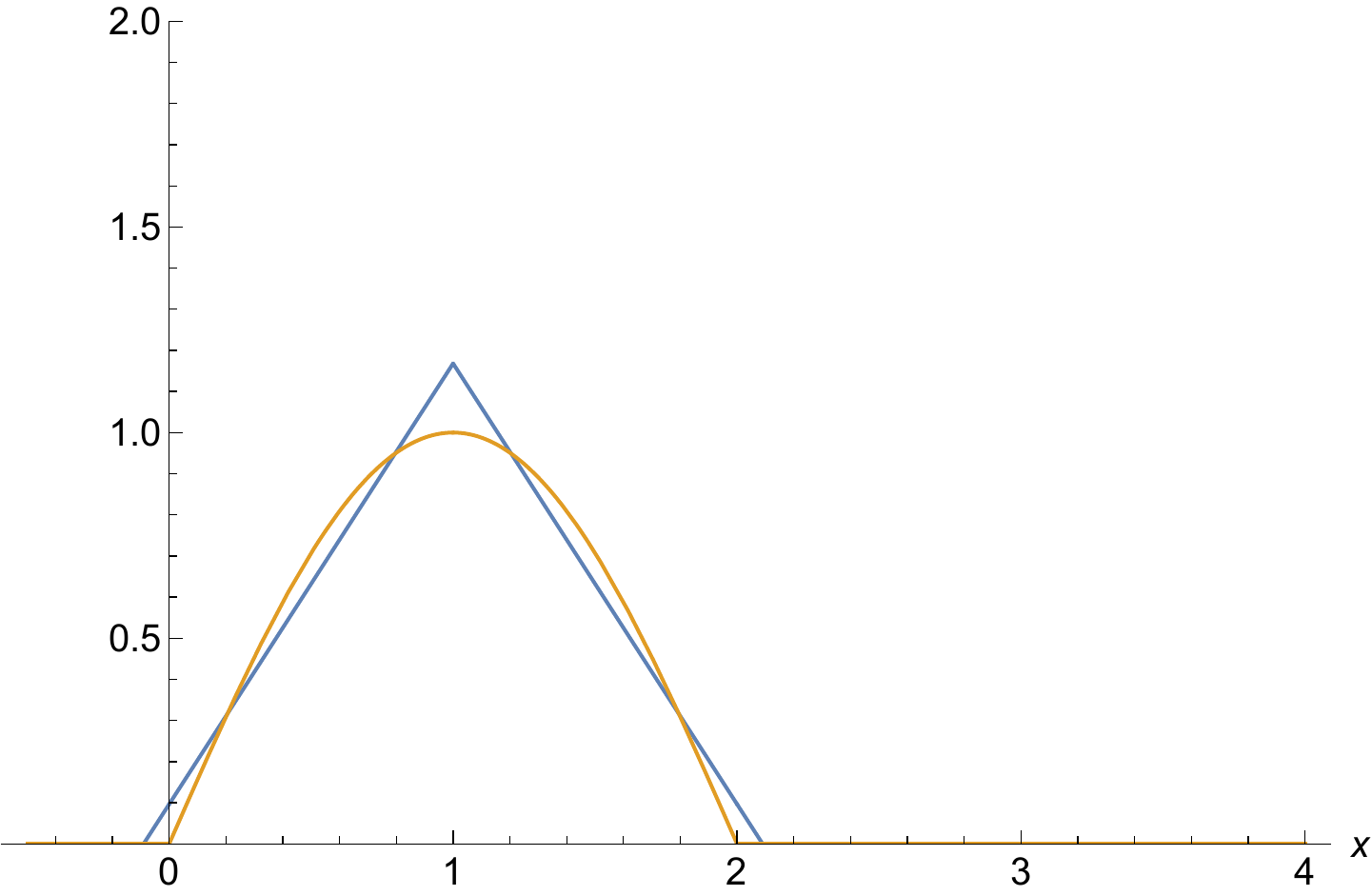}
            }
            \vspace{4mm}
            \captionof{figure}{The $relu$ FNN fit to the sine function described in the example problem statement. \vspace{4mm}} 

As is usual, we now form the loss function.  This is given by

\begin{align}
    J(w_{1,3},b_{2,3}) &= \int_{-\infty}^\infty [\sin(\pi x/2)-\hat{y}(x)]^2 \, dx\nonumber \\
    &= \int_{-\infty}^\infty\Big[\sin(\pi x/2)-\relu\Big(
\reluf{w_{1,3} x+b_{2,3}}-\reluf{2w_{1,3}[x-1]}  
   \Big)\Big]^2 \, dx\nonumber \\
   &\textrm{Constraints}:~~w_{1,3}>0, ~b_{2,3}>0
\end{align}
The integral bounds are set to $x\in(-\infty,\infty)$ here because both the function and its approximation are compact, but we do not know the exact domain of the approximating function yet (this is part of the problem we are solving!)  This domain accomplishes the integration without needing to know the exact bounds.  Another tip-- differentiation commutes with integration here.  Thus, when finding the minimum, it is easier to differentiate the integral form and then integrate.  Using that approach, we find the following result for the two derivatives.

\begin{align}
    \frac{\partial J}{\partial w_{1,3}}  & =  \frac{4 \pi ^2 b_{2,3}^3 w_{1,2}+12 \pi ^2 b_{2,3}^2 w_{1,2}+12 \pi ^2
   b_{2,3} w_{1,2}-24 \pi  b_{2,3}+4 \pi ^2 w_{1,2}-48}{3 \pi ^2}\\
    \frac{\partial J}{\partial b_{2,3}}  & = \frac{6 \pi ^2 b_{2,3}^2 w_{1,2}^2+12 \pi ^2 b_{2,3} w_{1,2}^2+6 \pi ^2
   w_{1,2}^2-24 \pi  w_{1,2}}{3 \pi ^2}
\end{align}
As before, we solve this set of equations with the derivative equal to zero, and solve them analytically using symbolic mathematics software.  The result is

\begin{align}
    w_{1,3}&=\frac{4 \pi}{36 - 36 \pi + 9 \pi^2} \\
    b_{2,3}&=\frac{-6 + 2 \pi}{\pi} \\
    J(w_{1,3},b_{2,3})
    &=\frac{1}{3 \pi ^2}
    \Bigg(
    3 \pi ^2+\frac{32 \pi ^4}{\left(36-36 \pi +9 \pi ^2\right)^2}
    +\frac{96 \pi ^3 (2
   \pi -6)}{\left(36-36 \pi +9 \pi ^2\right)^2}
   +\frac{96 \pi ^2 (2 \pi -6)^2}{\left(36-36
   \pi +9 \pi ^2\right)^2}\nonumber \\
   &+\frac{32 \pi  (2 \pi -6)^3}{\left(36-36 \pi +9 \pi
   ^2\right)^2}-\frac{192 \pi }{36-36 \pi +9 \pi ^2}-\frac{96 \pi  (2 \pi -6)}{36-36 \pi
   +9 \pi ^2}\Bigg)\nonumber \\
   &\approx 0.00861
\end{align}
\end{example}
\end{svgraybox}


\newpage

\textbf{Problems}

\begin{enumerate}
    \item Prove by constructing one representative  example that the sum of three ReLU functions 
    \begin{equation}
        f(x) = \sum_{i=1}^3 \relu_i(x+b_i)
    \end{equation}
    defined on $x\in[0,1]$ with biases of $b_i$ defined by $-1<b_3<b_2<b_1<0$ yield a continuous function with three segments.

    \item Prove by constructing a representative  example that two $\relu$ functions in two independent variables of the form

    \[ f(x,y)= ReLU(a_1 x +b_1)+ ReLU(a_2 y +b_2)\]
    Results in a continuous function.  How many separate regions in the $(x,y)$ plane are defined for the example
    
    \[ f(x,y)= ReLU( x -1)+ ReLU(2 y -2)\]

\item Define (but do not solve) an example of a 1-hidden-layer $\relu$ FNN that can be used to estimate the curve $f(x)=\sin(\pi x)$  on $x\in[0,1]$.  Assume there is only \emph{one node} in the hidden layer.  First write out the graph form of the network, and label the connections using the simplified conventions presented in the text.  After the graph is complete, write out the appropriate mathematical expression for the graph.

\item Define (but do not solve) an example of a 1-hidden-layer $\relu$ FNN that can be used to estimate the curve $f(x)=\sin(\pi x)$  on $x\in[0,1]$.  Assume there are \emph{two nodes} in the hidden layer.  First write out the graph form of the network, and label the connections using the simplified conventions presented in the text.  After the graph is complete, write out the appropriate mathematical expression for the graph.

\item Draw the graph for a $\relu$ FNN that contains 2 hidden layers with two independent variables ($x_1,x_2$) and two nodes in each hidden layer. Assume that the output is a single (weighted) scalar.  Assume that there is a bias term that enters at the input, and at the first hidden layer.   First write out the graph form of the network, and label the connections using the simplified conventions presented in the text.  After the graph is complete, write out the appropriate mathematical expression for the graph.

\item For the network graph illustrated below (Fig.~\ref{fig:probXOR}, write out the mathematical expression for the network (note that here, $H$ indicates a \emph{Heaviside} function rather than a $\relu$ function).  Once you have that written out, verify that the network produces behavior equivalent to the XOR function.\label{prob:XOR}

            \begin{figure}[!ht]
            \sidecaption[t]
           \centering
            \includegraphics[scale=0.25]{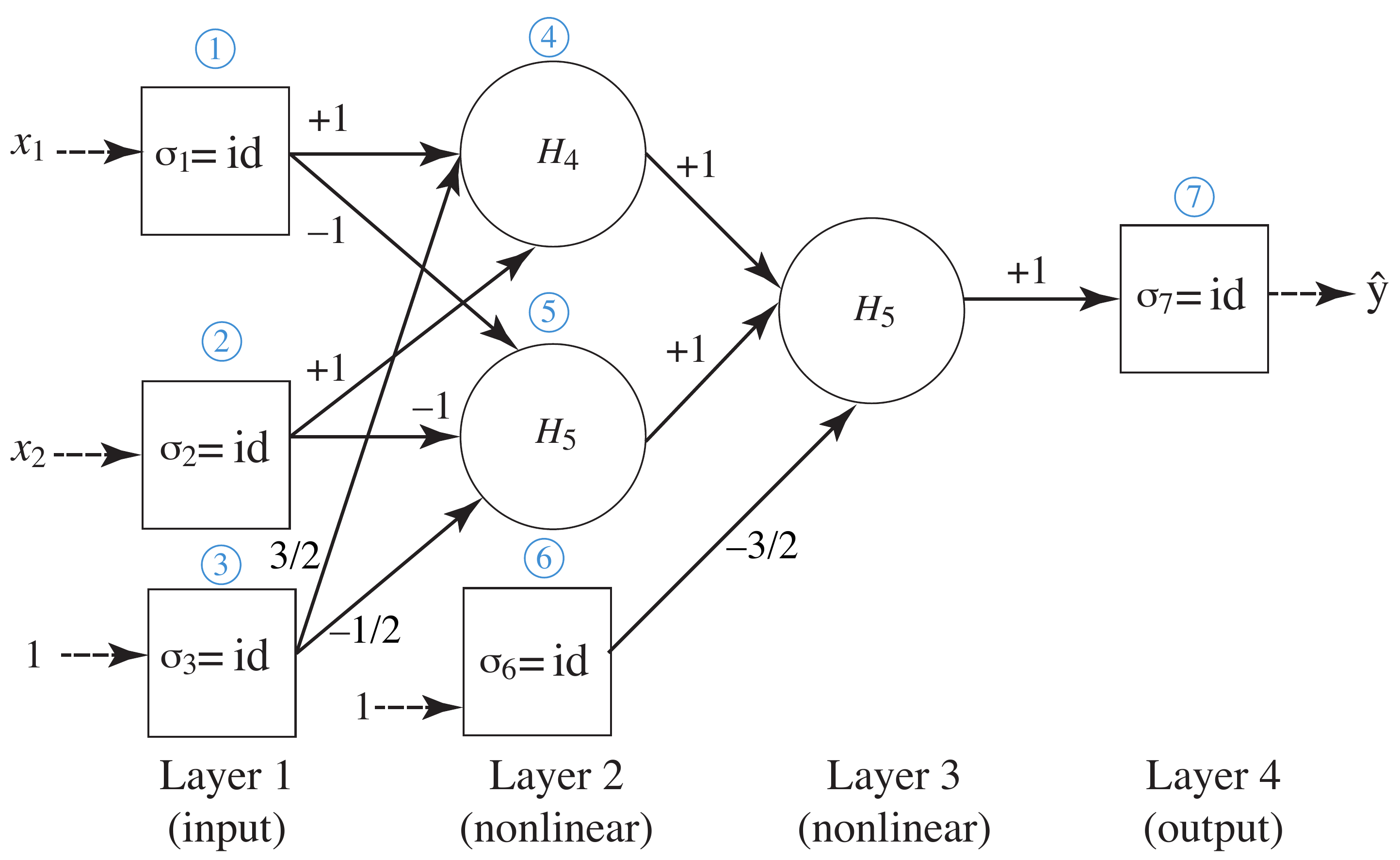}
            \vspace{-5mm}
            \caption{Problem \ref{prob:XOR}.}
            \label{fig:probXOR}       
            \end{figure}

\item Skip layer networks are still feedforward networks, but they are not simply connected (this term indicates that every node is connected to every other node only in the previous and following layers) but \emph{multiply connected}.  In the literature, these are sometimes called skip layers.   In applications, these kinds of networks are sometimes used to help avoid the ``vanishing gradients" problem in numerical backpropagation.\\

The graph illustrated in Fig.~\ref{fig:probXOR2} is an example of a skip connection FNN.  For the network graph illustrated, write out the mathematical expression for the network (note that here, $H$ indicates a \emph{Heaviside} function rather than a $\relu$ function).  Once you have that written out, verify that the network produces behavior equivalent to the XOR function.  Explain in one or two sentences why you think that the network illustrated here and the one in problem \ref{prob:XOR} both represent the same function. \label{prob:XOR2}

            \begin{figure}[!ht]
            \sidecaption[t]
           \centering
            \includegraphics[scale=0.3]{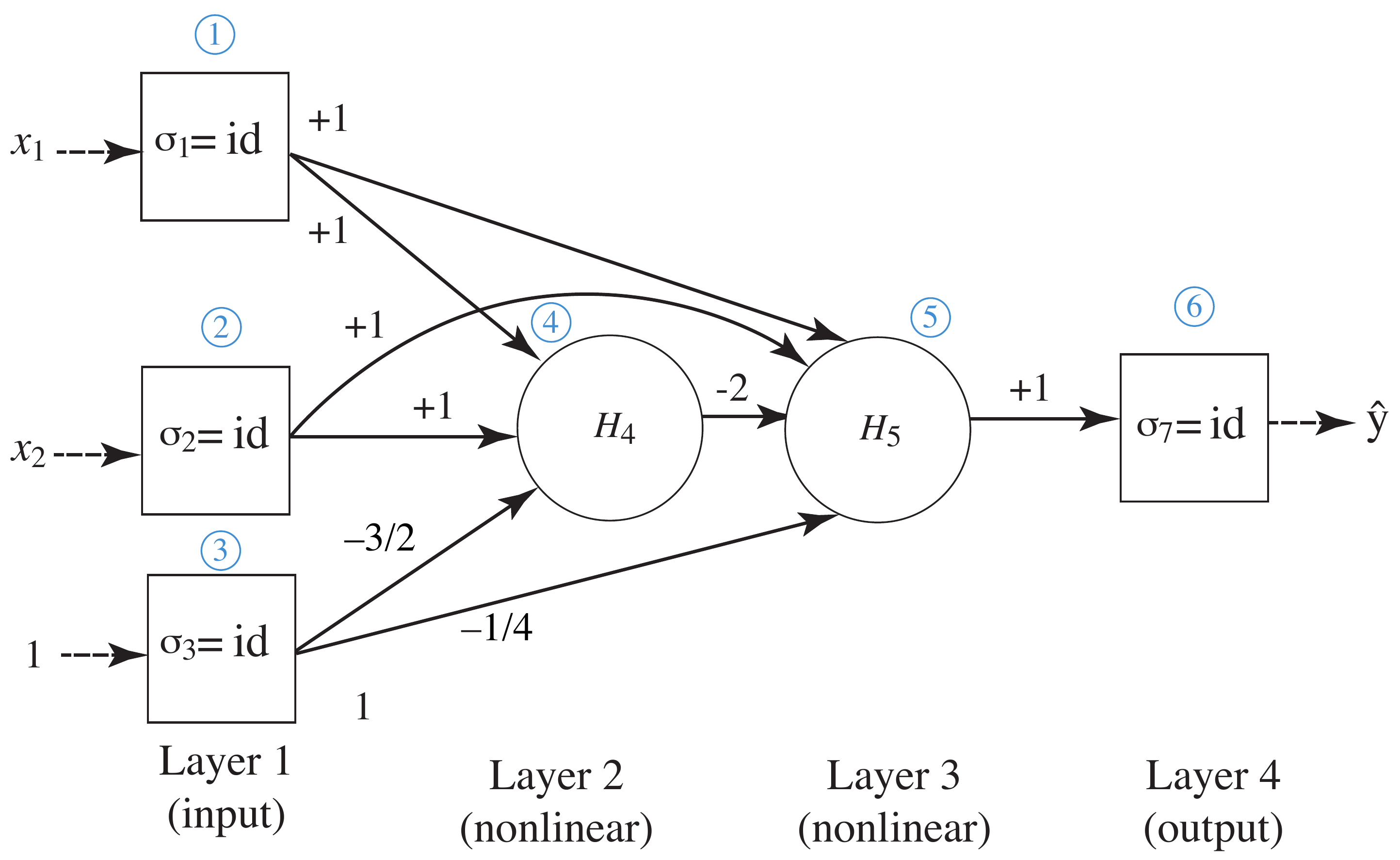}
            \vspace{-5mm}
            \caption{Problem \ref{prob:XOR2}.}
            \label{fig:probXOR2}       
            \end{figure}

\item As mentioned in the introduction, there are many different activation functions that can be used in ANNs, especially if they are not very deep.  The \emph{logistics} function is a nicely-behaved function (in that it has derivatives of all orders), and it behaves like a smooth version of the Heaviside function.  In general, the function is given by \label{prob:logistics}

\begin{equation}
    \textrm{Lo}(x) = \frac{1}{1+\exp[-k(x+b)]}
\end{equation}
where $b$ translates the function, and the parameter $k$ changes the slope of the function.  To get a Heaviside like function, we can set $k$ to a large value like $k=20$.

\begin{equation}
    \textrm{Lo}(x) = \frac{1}{1+\exp[-20(x+b)]}
\end{equation}
A plot of Lo($x$) appears in Fig.~\ref{fig:logistics}.
\newpage

            \begin{figure}[!ht]
            \sidecaption[t]
           \centering
            \includegraphics[scale=0.8]{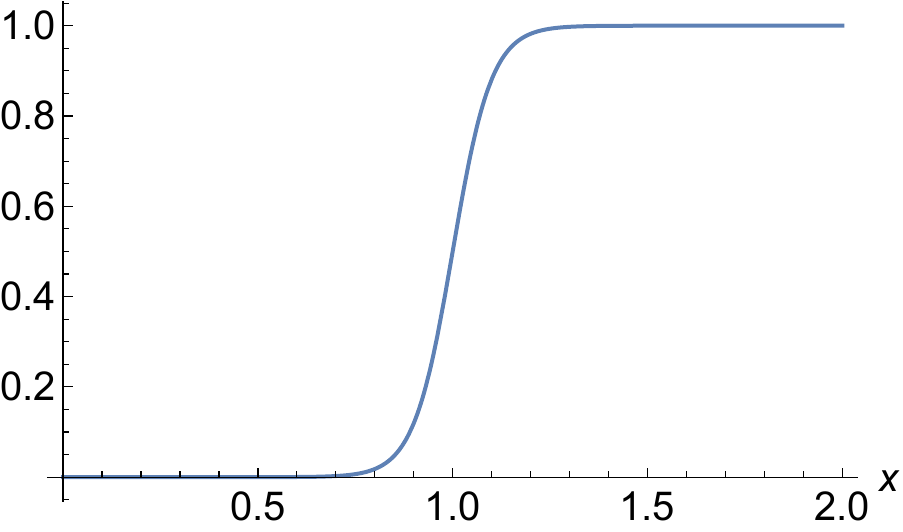}
            \vspace{-5mm}
            \caption{The logistics activation function used in Problem \ref{prob:logistics}.}
            \label{fig:logistics}       
            \end{figure}

Using this function, write out the mathematical expression for $\hat{y}$ for the network illustrated in Fig.~\ref{fig:logistics2}.  Then, determine the $L_2$ loss function using the four discrete values associated with the AND function.  As a reminder, the AND function is defined as given in Table~\ref{tab:ANDfunction2}.

\vspace{4mm}
{   \centering 
\begin{tabular}{|c|c|c|}
     $~~x_1~~$    & $~~x_2~~$ & $~~y~~$ \\
     \hline
      0  & 0 & 0 \\
      0 & 1 & 0 \\
      1 & 0 & 0  \\
      1 & 1 & 1
    \end{tabular}
    \captionof{table}{The binary AND function.}
    \label{tab:ANDfunction2}
    }
    \vspace{4mm}

Determine the value of $b_{3,4}$ by optimizing the loss function for the variable $b_{3,4}$.

            \begin{figure}[!ht]
            \sidecaption[t]
           \centering
            \includegraphics[scale=0.35]{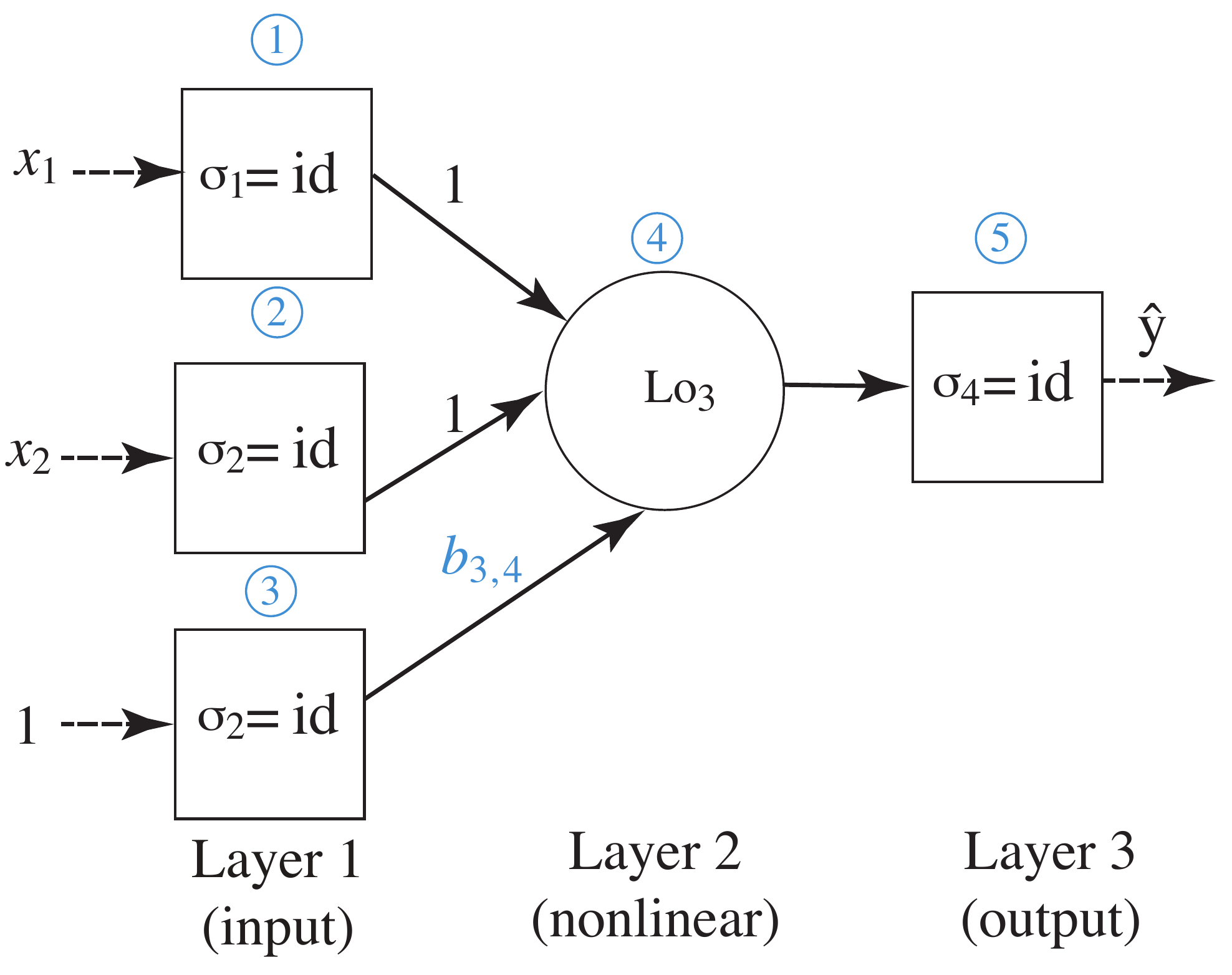}
            \vspace{-5mm}
            \caption{Problem \ref{prob:logistics}  A network for the AND function.}
            \label{fig:logistics2}       
            \end{figure}

\item Following the analysis in Example \ref{ex:sine_approx}, determine an approximation for the function

\begin{equation}
    y(x) =
    \begin{cases}
        0 & x<0 \\
        \sin[\pi x/2] & 0\le x \le 2 \\
        0 & x > 2
    \end{cases}
\end{equation}

Impose the following additional constraint on the problem: $b_{2,3}=0$.  Recall that we have demanded that the solution retain the symmetry of the original problem, so we still have that the slope to the left of $x=1$ and to the right of $x=1$ the same magnitude, but opposite signs.  This will result in a \emph{single} unknown parameter.  You are welcome to use a symbolic computation software to complete the problem, but it is also possible to optimize without such tools.\\

How does this solution differ compared to the solution obtained in the example problem?  How much difference in percent is there in the error term?  Use the average of the two errors as the denominator for your percent calculation.

\item \label{prob:constrainsine} Using symmetries can be a very powerful tool when attempting to find analytical solutions to complex problems.  As an example, consider the following function
\begin{equation}
    y(x) =
    \begin{cases}
        0 & x<0 \\
        \sin[\pi x/2] & 0\le x \le 4 \\
        0 & x > 4
    \end{cases}
\end{equation}
This problem is complicated by the fact that at least three $\relu$ units would be required to represent this function.  Thinking about using symmetries to reduce the computational burden of this network, do the following.

\begin{enumerate}
    \item Without making \emph{any} computations, sketch out the function $y(x)$ by hand.  Break the function up into three domains as follows: (1) $x < 1$, (2) $1 \le x \le 3$, and (3) $x>3$.  Note that the first and second intervals are unbounded!  Now, sketch out a linear piecewise continuous approximation of the function over these three intervals that maintains the symmetry around $x=2$.  Carefully identify the constraints on the slopes and on bias terms.  Draw the graph of a FNN that is capable of approximating this function. 

    \item Now, repeat this process, but now for a four-interval solution, where the intervals are: (1) $x < 1$, (2) $1 \le x \le 2$, (3) $2 \le x \le 3$, and (3) $x>3$.\\
    
   The solutions to Example\ref{ex:sine_approx} and to problem \ref{prob:constrainsine} contain all the information that we need to construct an optimal solution with a four-unit $\relu$ network if we employ the ideas of symmetry.   Carefully identify the constraints on the slopes and on bias terms.  Draw the graph of FNN that is capable of approximating this function. 
\end{enumerate}

\end{enumerate}

\backmatter
\bibliographystyle{splncsnat}
\bibliography{wood_book}

\printindex
\addcontentsline{toc}{chapter}{Index}
\end{document}